\documentclass[12pt]{book}

\makeindex

\usepackage{bussproofs}
\usepackage{amsmath}
\usepackage{amssymb}
\usepackage{wrapfig}
\usepackage{url}
\usepackage[normalem]{ulem}
\usepackage{stmaryrd}
\usepackage{eufrak}
\usepackage{graphicx}

\EnableBpAbbreviations

\def\HDS{\vrule width0pt height2.3ex depth1.05ex\displaystyle}
\def\f#1#2{{{\HDS #1}\over{\HDS #2}}}
\def\afrac#1{{\phantom{\HDS #1}\atop{\HDS #1}}}

\def\ul{$\backslash\;$}
\def\ud{$\slash\;$}

\def\zl{$\lfloor$}
\def\zd{$\rfloor$}

\def\riota{\rotatebox[origin=c]{180}{$\iota$}}

\newenvironment{rcases}
  {\left.\begin{aligned}}
  {\end{aligned}\right\rbrace}

\begin{document}

\pagenumbering{roman}

\title{\LARGE Logic Lectures\\[1ex] \Large G\" odel's Basic Logic\\ Course at Notre Dame}
\author{Edited by {\sc Milo\v s Ad\v zi\' c} and {\sc Kosta Do\v sen}
\\[1ex]
{\small Faculty of Philosophy, University of Belgrade}\\[-.5ex]
{\small \v Cika Ljubina 18-20, 11000 Belgrade, Serbia, and}\\
{\small Mathematical Institute, Serbian Academy of Sciences and Arts}\\[-.5ex]
{\small Knez Mihailova 36, p.f.\ 367, 11001 Belgrade, Serbia}\\[.5ex]
{\small email: milos.adzic@gmail.com, kosta@mi.sanu.ac.rs}}
\date{\small April 2017}
\maketitle

\clearpage

\pagestyle{empty}\makebox[1em]{} \clearpage

\pagestyle{myheadings}\markboth{}{LOGIC LECTURES}

\begin{small}

\noindent \emph{Abstract.} An edited version is given of the text of G\" odel's unpublished manuscript of the notes for a course in basic logic he delivered at the University of Notre Dame in 1939. G\" odel's notes deal with what is today considered as important logical problems par excellence, completeness, decidability, independence of axioms, and with natural deduction too, which was all still a novelty at the time the course was delivered. Full of regards towards beginners, the notes are not excessively formalistic. G\" odel presumably intended them just for himself, and they are full of abbreviations. This together with some other matters (like two versions of the same topic, and guessing the right order of the pages) required additional effort to obtain a readable edited version. Because of the quality of the material provided by G\" odel, including also important philosophical points, this effort should however be worthwhile. The edited version of the text is accompanied by another version, called the source version, which is quite close to G\" odel's manuscript. It is meant to be a record of the editorial interventions involved in producing the edited version (in particular, how the abbreviations were disabridged), and a justification of that later version.

\vspace{1ex}

\noindent {\it Keywords:} propositional logic, predicate logic

\vspace{1ex}

\noindent {\it Mathematics Subject Classification
(2010):} 01A60 (History of mathematics and mathematicians, 20th century), 03-01 (Mathematical logic and foundations, instructional exposition), 03-03 (Mathematical logic and foundations, historical)

\vspace{1ex}

\noindent \emph{Acknowledgements.} Our work was
supported by the Ministry of Education, Science and Technological Development of Serbia.
We are very grateful indeed to Gabriella Crocco for making G\" odel's
unpublished notes available to us. The decipherment and publication of some of these notes
is part of the project {\it Kurt G\" odel Philosopher: From
Logic to Cosmology}, which is directed by her and funded by the French
National Research Agency (project ANR-09-BLAN-0313). We are also very grateful to Peter Schroeder-Heister for encouraging us in our project. We are grateful to the Institute for Advanced Study in Princeton for granting us, through the office of its archivist Mr Casey Westerman, the permission to
publish G\" odel's unpublished notes; we were asked to give credit for that with the following text:
``Works of Kurt G\" odel used with permission of Institute for Advanced Study. Unpublished Copyright
(1939) Institute for Advanced Study. All rights reserved.''

\end{small}

\clearpage

\pagestyle{empty}\makebox[1em]{} \clearpage

\noindent\textbf{\Large Abbrev.\zl iated\zd\ editorial introduction}\label{Editorial}
\pagestyle{myheadings}\markboth{LOGIC LECTURES}{EDITORIAL INTRODUCTION}
\\*[1.5ex]
G\" odel taught a one-semester course in basic logic at the University of Notre Dame in the spring of 1939, when he turned 33. Among his unpublished writings in the Princeton University Library one can find notebooks with the manuscript of his notes for that course. The title \emph{Logic Lectures}, which we gave to these notes, is suggested by the German ``Log.\zl ik\zd\ Vorl.\zl esungen\zd '', or a variant of that, written on the front covers of the notebooks.

Besides the Notre Dame course G\" odel taught a basic logic course in Vienna in the summer of 1935, notes for which, on 43 notebook pages (27 of which are numbered), made mainly of formulae and very little accompanying text in ordinary language, have been preserved in a manuscript at the same place. The notes for the Notre Dame course, which with their 427 notebook pages are ten times bigger, are more detailed and we think more important. Propositional logic is not much present in the Vienna notes.

We have published recently in \cite{AD16} a brief, and hence not complete, summary with comments of the Notre Dame notes, and an assessment of their importance. This preceding short paper is a natural introduction to this introduction, which is more oriented towards details concerning G\" odel's text. We deal however here occasionally, in the paragraph on definite descriptions below and in the last few pages of this introduction, with some matters of logic and philosophy, partly in the sphere of the preceding paper, but not to be found there. Anyway, that paper enables us to abbreviate this introduction (which explains up to a point its title; the rest will be explained in a moment).

We will not repeat ourselves, and we will not give again all the references we gave in the preceding paper, but we want to mention however John Dawson, who in \cite{Daw} supplies biographical data on G\" odel's stay at Notre Dame, John and Cheryl Dawson who in \cite{Daw05} set what we did with the Notre Dame notes as a task for G\" odel scholars,\footnote{We are grateful to John Dawson for encouraging us to get into this publishing project.} and Pierre Cassou-Nogu\` es, who has published in \cite{CN09} a dozen printed pages extracted and edited from G\" odel's manuscript of the Notre Dame course (this concerns pp.\ \textbf{1}.-\textbf{26}.\ of Notebook~I, including small bits of Notebook~0, pp.\ \textbf{73}.\textbf{1}-\textbf{73}.\textbf{7} of Notebook~V, pp.\ \textbf{122}.-\textbf{125}., \textbf{134}.-\textbf{136}.\ of Notebook VI and pp.\ \textbf{137}.-\textbf{157}.\ of Notebook VII; altogether 60 notebook pages).\footnote{We have found sometimes useful Cassou-Nogu\` es' reading of G\" odel's manuscript, and we wish to acknowledge our debt. Our decipherment of the manuscript does not however accords always with his, and we have not followed his editorial interventions.}

Besides the edited version of G\" odel's text we have prepared another version of it, which we call the source version, and the present introduction should serve for both of them. This other, source, version is quite close to the original manuscript, and is meant to be a record of the additions and other interventions made in the manuscript to arrive at the edited version, and a justification of that later version.

G\" odel used abbreviations in the manuscript of the notes quite a lot. For example, the second sentence and the beginning of the third of Notebook~0 of the manuscript are: ``Accord.\ to this def the centr.\ part of log.\ must be the theory of inf and the theory of logically true prop. By a log true prop. I mean
a prop. which is true for merely log reasons\ldots'' In the source version this is rendered as: ``Accord.\zl ing\zd\ to this def\zl inition\zd\ the centr\zl al\zd\ part of log.\zl ic\zd\ must be the theory of inf\zl erence\zd\ and the theory of logically true prop\zl osi\-tions\zd . By a log\zl ically\zd\ true prop.\zl osition\zd\ I mean a prop.\zl osition\zd\ which is true for merely log\zl ical\zd\ reasons\ldots'' All the abbreviated words are typed in the source version as they occur in the manuscript, with a full stop after the abbreviation or without, together with their prolongation or decipherment within the parenthetical signs \zl\ and \zd\ to obtain the non-abbreviated, disabridged, word they are supposed to stand for, which one finds in the edited version. Sometimes whole words are omitted and they are restored in the source version within \zl\ and \zd .

Using abbreviations may produce problems, which are however surmountable. For example, log., with or without full stop, stands for ``logic'', ``logically'' and ``logical''. Singular or plural has to be inferred from the context; ``form.'', with or without full stop, stands for ``formula'' or ``formulas'' (G\" odel has the plural ``formulas'' while we here and in our comments use ``formulae''; he says often ``expression'' for ``formula''). Sometimes, but not very often, it is not obvious, and even not certain, what is the abbreviated word; for example, both ``proposition'' and ``property'' are abbreviated by ``prop.''. This involvement with abbreviations in the manuscript goes so far that one finds even ``probl.'' for ``problem'' and ``symb.'' for ``symbol''.  Because of their number, and some particular problems they produced occasionally, taking care of the abbreviations made our editing task considerably harder, but this number tells that they cannot be neglected if one wants to leave a more precise record of G\" odel's style (see the end of this introduction).

In the source version one may also find all the parts of the text crossed out in the manuscript, with the indication that they were found crossed out, either by being really crossed out in the source version, or if they are too long, the crossing out is mentioned as an editorial comment within \zl\ and \zd . We use \zl\ and \zd\ in the source version in connection with the abbreviations as we said above, and in general for other editorial comments too. (For example, we will have \zl unreadable text\zd .)

In a few cases we have estimated that a crossed out part of the text is worth reproducing even in the edited version. (G\" odel's crossing out a text need not mean dissatisfaction with it, but it may mean perhaps lack of time to use it in the lectures.) In one place it may compensate a little bit for a lost part of the text (see the footnote on p.\ \textbf{7}. of Notebook IV), in another (see the footnote on pp.\ \textbf{114}.-\textbf{115}.\ of Notebook VI), it completes what is needed for establishing that binary relations with composition and the identity relation make a monoid. (Composition of relations is called by G\" odel ``relative product'', and his examples for it are with relations between \emph{relatives}, nephew, son, brother, sister, uncle, father, grandfather, grandchild, child,\ldots, which is etymologically inspirative.)  A third such place, which is tied to Russell's understanding of definite descriptions (see pp.\ \textbf{123}.-\textbf{125}. of Notebook VI), is philosophically important.

Let us dwell for a moment at this third place, to justify our choice of reproducing the crossed out text. G\" odel's says there that taking ``The present king of France is bald'' as meaningless is undesirable because whether the present king of France exists is an empirical question. He then continues: ``Therefore it would depend on an empirical fact whether or not this sequence of words is a meaningful statement or nonsense, whereas one should expect that it can depend only on the grammar of the language concerned whether something makes sense.'' So G\" odel asserts the primacy and independence of the understanding of language over empirical, i.e.\ epistemological, matters. The primacy of the linguistic over the epistemological (and presumably other philosophical concerns, like the ontological, or axiological) should be one of the main, if not the main, mark of the linguistic turn in twentieth century philosophy. G\" odel's single sentence quoted above is more significant and more explanatory than thousands and thousands of others in the sea of ink spilled over the king's baldness.

The notes are written by hand in English in eight notebooks bound by a spiral, with however some loose leafs (four leaves on a different paper, not torn out from the notebook, without holes for the spiral, at the end of Notebook III with pp.\ \textbf{new page x-xiii}, nine torn out leafs towards the end of Notebook~V including pp.\ \textbf{73}.\textbf{1}-\textbf{73}.\textbf{7}, and nine torn out leafs at the end of Notebook VII with pp.\ \textbf{new page iii-iv} and \textbf{1}.-\textbf{7}.). G\" odel writes usually on the left pages, the back sides of the leafs, and he uses the right pages, the front side of the leafs, most often for inserted additions, or simply continuations of the text from the left pages. As insertion signs, one finds most often $\forall$ (which is not used in the manuscript for the universal quantifier), but also $\times$, and a few others. Insertions tied to these signs, as well as other insertions, often tied to \raisebox{1ex}{$\underbrace{}$}, but not continuations on the right pages, are marked in the source version with \ul at the beginning of the insertion and \ud at its end. Sometimes one finds remarks and examples not possible to insert simply in the main text, and they are not to be found in the edited version. Since usually only the left pages are numbered, and the right page is usually associated with the left, we do not speak of left and right pages, but say, for example, that something occurs on the right of a certain page, or use similar forms of speaking.

There are no footnotes in the source version, because G\" odel does not have them. (We do not interpret his insertions as footnotes.) All the editorial comments there are within \zl\ and \zd . All the footnotes in the edited text are ours, and they are made of editorial comments.

In general we have strived to stay as close to G\"odel's text as possible, at the cost of failing to follow standard usage. G\" odel's manners in writing are sometimes strange, according to the contemporary standards, but they always make sense. (On pp.\ \textbf{47}.-\textbf{49}.\ of Notebook II he says, for example, ``then and only then'' for ``if and only if'', which one finds later. Instead of three dots as a punctuation mark he uses two---perhaps because he wants to abbreviate---but we have rendered that both in the source and the edited version in the usual triple way.)

We have corrected G\" odel's not very numerous spelling mistakes, and did not keep in the edited text peculiar or foreign spelling (like ``tautologie'' and ``geometrie''). If however an unusual spelling (like, for example, \emph{caracter} instead of \emph{character}) is permitted by the Oxford English Dictionary, then we kept it. We have not corrected G\" odel's style in the notes, and we are aware that it is often on the edge of the grammatically correct, and perhaps even sometimes on the other side of the edge. In cases of doubt we opted for keeping his words. We made this choice because thereby the reader should be able to hear better G\" odel lecturing, to hear his voice and not the voice of somebody else. G\" odel had at that time no doubt his own foreign accent, which, since we ourselves are not native speakers of English, we did not want to replace with ours.

G\" odel omitted in the notes many punctuation marks, in particular commas and quotation marks, but also full stops, presumably for the sake of abbreviating. We have added them, in the source version with \zl\ and \zd\ and in the edited text, together with some colons, only when we considered they are absolutely indispensable, but we did not want to add all of them that would usually be written. For example, G\" odel
practically never wrote commas before ``then'', and we did not add those.

G\" odel was very sparing in using quotation marks. (Initial quotation marks he wrote in the German way ,, and not ``.) He did not use them systematically for naming words and sentences. We did put them at many places where we were afraid understanding would be endangered, but at the cost of looking unsystematic, as G\" odel, we did not restore them everywhere. We felt that in doing that, analogously to what we said in the preceding paragraph, we would be too intrusive, and get too estranged from G\" odel's customs and intentions. Perhaps he did not omit quotation marks just for the sake of abbreviating, but wanted to use words autonymously, which might be related to his involvement with self-reference (see the end of this introduction). Once one becomes accustomed to this autonymous use, it hardly leads to confusion.

To make easier comparison with the scanned manuscript (which is the only one we have seen), we have standardized only slightly the numbers of the pages G\" odel assigned to them there. These numbers are rendered in both the source and edited version with boldface Arabic figures, followed by a full stop, which is to be found in the manuscript, but not always, and also further figures, Arabic, Roman, or letters found in the manuscript; examples will come in a moment. We found five successive, not very systematic, numberings of pages in the manuscript starting from pages numbered \textbf{1}.\ in various notebooks. Some pages were left unnumbered by the numberings, and we introduced our own way of naming them, usually with the label \textbf{new page}.

We believe the first numbering is made of pp.\ \textbf{1}.-\textbf{26}.~\textbf{I} of Notebook~I (where a break occurs in that notebook). We will explain below why we think these pages of Notebook~I should precede Notebook~0.

The second numbering starts with pp.\ \textbf{1}-\textbf{38}.\ of Notebook~0 (i.e.\ the whole of that notebook),  followed by pp.\ \textbf{38}.\textbf{1}~\textbf{II}-\textbf{44}.\ \textbf{II} of Notebook~I, followed by pp.\ \textbf{33}.-\textbf{55}.\textbf{2}\ of Notebook II, followed by pp.\ \textbf{56}.-\textbf{60}.\ of Notebook~I, followed finally by pp.\ \textbf{61}.-\textbf{76}.\ of Notebook II. Our reasons for this complicated arrangement are in the sense of the text. For example, the involvement of Notebook II in this numbering has to do with the presentation of the axiom system for propositional logic (see Section 1.1.9 in the edited text below). We must warn however that though in this numbering the page numbers from different notebooks sometimes fit perfectly, and follow the sense, sometimes the fitting is somewhat less than perfect.

We have rearranged the page order in our edited version as the first and second numberings require. In the source version the original order from the scanned manuscript is kept in general, and also for the pages involved in these numberings. The order of pages required by the remaining three numberings are the same in the edited and source version and in the scanned manuscript, with a small exception which we will mention in a moment.

The third numbering is from the initial, first, p.\ \textbf{1}.\ of Notebook III up to p.\ \textbf{53}. of that notebook.

The fourth, longest, numbering is from the second p.\ \textbf{1}.\ of Notebook III, which is close to the end of the notebook, up to p.\ \textbf{157}.\ of Notebook VII, following more or less regularly the order of the notebooks and the numbering in them.

A small rearrangement guided by subject matter is made in the edited version in the last part of Section 1.1.10, where guided by subject matter four pages from Notebook IV not numbered in the manuscript have been inserted, which has made possible a perfect fitting in Section 1.1.14 \emph{Sequents and natural deduction system}.

The fifth, last and shortest, numbering is made of pp.\ \textbf{1}.-\textbf{7}.\ of Notebook VII, at the very end.

Zero precedes one, and presumably because of that, in the scanned manu\-script Notebook~0 precedes Notebook~I, while in \S 1.II of \cite{Daw05} one finds that Notebook~I ``appears to be a rewritten, somewhat condensed version'' of Notebook~0. It is however not clear in relevant cases that condensation from 0 to I is made, and sometimes the opposite, addition, from I to 0 seems to be at work. Sometimes even the text in Notebook~0 is shorter than the corresponding text in Notebook~I, from which it seems to have been obtained by tidying up (cf.\ in the source version the text pp.\ \textbf{20}.-\textbf{21}.\ of Notebook~0 with the approximately twice longer corresponding text on pp.\ \textbf{15}.-\textbf{16}.\ of Notebook~I). We want to present now additional reasons for believing that Notebook~I precedes Notebook~0, and that Notebook~0 together with the parts mentioned in the second numbering above is written later and may be considered to supersede the pages of Notebook~I in the first numbering.

From p.\ \textbf{4}.\ until the end of p.\ \textbf{21}.\ of Notebook~I propositional
variables are written first mostly as capital $P$, $Q$ and $R$, which are
later on alternated with the lower-case $p$, $q$ and $r$. In the edited version they are
all written uniformly as lower-case, because when they alternate they might be confusing, while in the source version they are as in the manuscript. After p.\ \textbf{21}.\ of Notebook~I and in Notebook~0 the lower-case letters only are used for propositional variables. This usage is kept in Notebook II and later, and the capital letters starting from  p.\ \textbf{58}.\ of Notebook~I, which belongs to our second numbering, and later, are used as schematic letters for formulae. The notation in Notebook~0 seems to be a correction of that in Notebook~I.

Before p.\ \textbf{42}.\ \textbf{II} of Notebook~I the signs + and $-$, which were used in the notes for the 1935 Vienna course, are used instead of T and F for naming truth values. The letters T and F are to be found in Notebook~0, on pages of Notebook~I that belong to our second numbering, and they are used regularly in Notebook~II and later. In the edited text we did not try to replace + and $-$ by T and F, because no confusion is likely.

The pages numbered in the manuscript with the suffix I in Notebook~I, which belong to our first numbering, could be superseded by pages after p.\ \textbf{23}.\ of Notebook~0, which leave a better impression and belong to our second numbering. The suffix II added in the manuscript to some later pages in Notebook~I would indicate that these pages belong to the second numbering.

In Notebook~I decidability is considered with tautologies on pages that make Section 1.1.7 \emph{Decidability for propositional logic} of the edited text. In Notebook~0 decidability is not considered, but it is considered more thoroughly on pp.\  \textbf{41}.\ \textbf{II}-\textbf{44}.\ \textbf{II} of Notebook~I, which belong to our second numbering.

The axioms of the system for propositional logic would appear for the first time on p.\ \textbf{53}.\ of Notebook II, which until the end Notebook II is followed by a preliminary discussion of the role of primitive rules of inference in logic (we consider this matter below in a more philosophical spirit), but no such rule is given. The primitive inference rules are to be found on pp.\ \textbf{56}.-\textbf{59}.\ of Notebook~I, and after them the four axioms are given again on p.\ \textbf{60}.\ of Notebook~I. This induced part of the order in our second numbering.

On pp.\ \textbf{11}.-\textbf{12}.\ of Notebook~I G\" odel writes something like handwritten $o$, which we put (or perhaps $\sigma$), for exclusive disjunction, while on pp.\ \textbf{16}.\ and \textbf{18}.\ of Notebook~0 he has for it $\circ$, which is then again to be found on p.\ \textbf{44}.\ of Notebook II.

On the same pages pp.\ \textbf{11}.-\textbf{12}.\ of Notebook~I, and also on p.\ \textbf{7}.\  of the same notebook, one finds a number of times a crossed out word ``wrong'' replaced by ``false''. In Notebook~0 ``wrong'' is not to be found and ``false'' is used regularly, while later ``wrong'' occurs here and there, but ``false'' predominates.

At the very beginning of the notes, the programme of the course is stated together with a reprobation of traditional logic (which we will consider below in this introduction). Citing the source version, a sentence in that part starts with: ``What the textbooks \sout{give} and also what Arist.\zl otle\zd\ gives is a more or less arbitrary selection of the \ul infinity of \ud \zl the\zd\ laws of logic'' on p.\ \textbf{1}.\ of Notebook~I, and with: ``What the trad\zl itional\zd\ logic gives is a more or less arbitrary selection from the infinity of the laws of logic'' on pp.\ \textbf{1}-\textbf{2}.\ of Notebook~0. We have not gone over the matter systematically, but it seems to us that this is an indicative sample of what happens when one passes from Notebook~I to Notebook~0. In Notebook~I we have ``selection of the infinity of laws of logic'', where ``infinity of'' has been inserted (``\zl the\zd '' means that the article has been added by us in the edited version), while in Notebook~0 we have ``selection from the infinity of the laws of logic'', which is less ambiguous and better English. Note, by the way, that Aristotle and textbooks are not mentioned here in Notebook~0 (on p.\ \textbf{1}.\ of Notebook~0 a mention of textbooks a little bit earlier has been crossed out, as marked by a footnote in the edited version).

We conclude our discussion about Notebook~I preceding Notebook~0 with a detail that sets Notebook~0 apart, and that together with the number of that notebook may point in the other direction. On the front cover of Notebook~0 one finds ``Vorl.\ Log.'', while on the front covers of all the remaining notebooks one finds ``Log.\ Vorl.'', except for Notebook VII, where ``Logik Vorl.'' is written (see the source version).

G\" odel's text has neither chapters nor sections, nor an explicit division into lectures. The edited version and the source version make two chapters in this book. We have divided the edited version into two parts, the first about propositional and the second about predicate logic, and we have further divided these parts into sections which, as the parts, we have named with our own words. Our titles of the parts and sections are not mentioned in the source version. For them we use standard modern terminology and not G\" odel's. We put ``connectives'' instead of ``connections''. G\" odel did not use the expressions ``functional completeness'', ``disjunctive normal form'', ``conjunctive normal form'', ``sequents'', ``natural deduction'', ``first-order languages'', ``valid formulas'' (he uses ``tautology'' also for these formulae, or he says that they are universally true). He uses the term ``class'' rather than ``set'', and we have kept it for naming Sections 1.2.7 and 1.2.8 in the edited text.  Our table of contents below is not exactly the same as that given in \cite{AD16}. The present one is more detailed and follows more closely the manuscript, including repetitions in it. We have added moreover to the edited text an index for it.

G\" odel did not pay very much attention in the notes to the division of the text into paragraphs, and where we found it very desirable, following either the sense of the text or rather the excessive length of the paragraphs in the manuscript, we introduced new paragraphs, with due notice, using \zl new paragraph\zd , in the source version. We did not introduce them however at all places where this might have been done, following a policy similar to the one we had with punctuation marks.

Some, but not much, of G\" odel's text is unreadable and a very small part of it is in shorthand. Sometimes it is not clear whether one has to do with shorthand or unreadable text. We have not tried to decipher the shorthand in the source version, because practically everywhere it occurs in parts omitted in the edited version, which do not belong at all to the main text, and sometimes are not directly about logic (as, for example, in the theological remarks at the beginning of Notebook VII). We did not find we need this decipherment. The unreadable portions of the text are marked with the words ``unreadable text'', ``unreadable symbol'', or something related.

Pages written not very systematically, not numbered, with lists of formulae, jottings, and some unreadable text, crossed out to a great extent, have been rendered as far as possible in the source version but not in the edited one. We did not want to be too intrusive by making a selection in this text, which we estimate should not all belong to the edited version. There are thirteen such pages at the end of Notebook III. Notebook VII starts with nine, not numbered, pages of remarks and questions mostly theological, partly unreadable, partly in shorthand, and all seemingly not closely related to the remaining notes for the course. They are rendered as far as possible in the source version but not in the edited one. The text crossed out in the manuscript is not in the edited version.

The underlined parts of the manuscript have in principle been rendered in the edited version by italics. The underlining has however been kept in derivations where it can play a special role.

As we said in \cite{AD16} (see the section \emph{Major problems and branches of logic}), \cite{HA28} influenced G\" odel in general, and that influence is to be found in the Notre Dame course too. (This influence might be seen in details like the remarks on the Latin \emph{aut} and \emph{vel} on p.\ \textbf{9}.\ of Notebook~0, which follow \cite{HA28}, Section I.\S1, but G\" odel also mentions \emph{sive}\ldots\ \emph{sive} on p.\ \textbf{7}. of Notebook~I.) In the notes G\" odel does not use the expressions ``formal language'' and ``inductive definition'', and does not have a proper inductive definition of the formal language, i.e.\ of the formulae, of propositional logic (he comes nearest to that on pp.\ \textbf{11}.\ and \textbf{15}.\ of Notebook~0 and p.\ \textbf{8}.\ of Notebook~I). The formal language of propositional logic is not defined more precisely in \cite{HA28}, though a formal language of first-order predicate logic is defined by a regular inductive definition in Section III.\S 4. In the Notre Dame notes however, the formulae of predicate logic are not defined more precisely than those of propositional logic (see pp.\ \textbf{32}.ff of Notebook IV). It seems that in many textbooks of logic, at that time and later, and even today, clear inductive definitions of formal languages might be lacking, the matter being taken for granted.

In the precise inductive definition of formulae in \cite{Go31} (Section~2, pp.\ 52-53 in the \emph{Collected Works}), his most famous paper, G\" odel has the clauses that if $a$ is a formula, then ${\sim(a)}$ is a formula, and that if $a$ and $b$ are formulae, then $(a)\vee(b)$ is a formula. This definition excludes outermost parentheses, but in complex formulae it puts parentheses around propositional letters and  negations, where they might be deemed unnecessary. This way of dealing with parentheses should explain why on pp.\ \textbf{14}.-\textbf{15}.\ of Notebook~0 (and occasionally also elsewhere, as on pp.\ \textbf{23}.ff of Notebook III) it is taken that there are parentheses around negations, as in $(\sim p)$, which are not customary, and that there should be a convention that permits to omit them.

To prefix the universal and existential quantifiers $(x)$ and $(\exists x)$ square brackets are put in the notes around formulae before which they are prefixed, which is also neither customary nor necessary, as noted on p.\ \textbf{41}.\ of Notebook IV, where in some cases it is permitted, but not required, to omit these brackets. As in some other matters of logical notation, neither the convention to write the brackets nor the permission to omit them are followed systematically (see pp.\ \textbf{32}.\textbf{a}\,ff of Notebook IV). We have not tried to mend always this and similar matters in the edited text. Besides corrections of slips of the pen, found in formulae as well as in English, but not very numerous, we have made changes of what is in the manuscript in cases where we estimated that understanding would be hampered.

G\" odel's usage in the notes is not very systematic and consistent, neither concerning formalities of logical notation, nor concerning matters of ordinary English, including punctuation marks (which he does not use as much as it is usual). One should however always bear in mind that the notes were presumably meant only for himself, and he could correct in the lectures whatever irregularity they contain. This matter concerns also sometimes the meaning of his text, which taken literally is not correct. He speaks, for instance, nearly always of substitution of \emph{objects} and not of their \emph{names} for individual variables (on p.\ \textbf{42}.\ of Notebook IV one finds, for example, ``the free variables are replaced by individual objects''). On p.\ \textbf{138}.\ of Notebook VI he says ``for any arbitrary object which you substitute for $x$'', but three lines below he says ``if you substitute for $x$ the name of an arbitrary object''. On  p.\ \textbf{139}.\ of Notebook VII he has ``if you substitute
for $x$ the name of an arbitrary object'', with ``the name of'' inserted later (which in our source version is rendered with \ul and \ud ). So one may take that G\" odel had always in mind the correct statements mentioning names, which at most places he omitted for the sake of abbreviating, which he relied on very much. (It is also possible that sometimes, except where names are mentioned, by \emph{substituting} an object for a variable G\" odel meant \emph{interpreting} the variable by the object.)

G\" odel's definition of tautology for propositional logic (see pp.\ \textbf{33}.\ of Notebook~0 and \textbf{25.~I}.\ of Notebook~I) and valid formula, i.e.\ tautology or universally true formula in his terminology, for predicate logic (see p.\ \textbf{45}.\ of Notebook IV) are not very formal. His definitions could be taken as defining syntactical notions based on substitution, if this substitution is not understood as model-theoretical interpretation (cf.\ the parenthetical remark at the end of the paragraph before the preceding one). The word ``model'' does not however occur in the notes, and the notion, which is somehow taken for granted, is not introduced with much detail.

Concerning tautologies of predicate logic, one finds on p.\ \textbf{54}.\ of Notebook IV and p.\ \textbf{55}.\ of Notebook V: ``An expression is a tautology if it is true in a
world with infinitely many individuals, i.e.\ one can prove that whenever an expression is universally true in a world with infinitely many objects it is true in any world no matter how many individuals there
may be and of course also vice versa.'' G\" odel says that he cannot enter into the proof of that. (For this matter one may consult Section III.\S 12 of \cite{HA28}.)

G\" odel seems parsimonious by relying a lot on abbreviations, but he does not spare his energy and time in explaining quite simple matters in great detail, and in repeating himself. He addresses beginners, and does not forget that they are that. This might be a reason to add to those mentioned in the following concluding remark in \S 1.II of \cite{Daw05} concerning the Notre Dame notes: ``Although the material is standard, the choice and ordering of topics, as well as some of the examples that are discussed, may well be of pedagogical interest.'' In the remainder of this introduction, we will give reasons that should be added to those given in \cite{AD16}, \cite{DA16} and \cite{DA16a} to justify our belief that the interest of these notes is not just pedagogical.

Our involvement with G\" odel's notes from Notre Dame started with an interest in G\" odel's views concerning deduction, about which we wrote in \cite{DA16} and \cite{DA16a}. This was the main reason for our getting into the project, which, as can be gathered from \cite{AD16}, led to other matters concerning the course  that we found interesting. (Also, one of us taught a logic course as a visiting professor at Notre Dame when he turned 33.) Concerning deduction, we would like to add here that on pp.\ \textbf{69}.-\textbf{70}.\ of Notebook II G\" odel commends derived rules and says ``in our system we cannot only derive formulas but also new rules of inference''. We believe this short remark is in accordance with our discussion in \cite{DA16a} and \cite{DA16} of G\" odel's natural deduction system of Notebook IV and his recommendation of it in Notebook III. G\" odel's remarks about rules of inference on pp.\ \textbf{52}.-\textbf{55}.\textbf{2}\ at the end of Notebook II, which in the edited text are at the beginning of Section 1.1.9 \emph{Axiom system for propositional logic}, are relevant too for G\" odel's opinions about deduction. G\" odel says there that if rules are not formulated explicitly and derivability is understood as, for example, in geometry, where it means ``follows by logical inference'', then ``every logical law would be derivable from any other'' (p.\ \textbf{55}.\textbf{1}\ of Notebook II; cf.\ the second p.\ \textbf{4}.\ towards the end of Notebook III).

In the edited text we entitled Section 1.1.4 of Notebook~0 and the corresponding Section 1.1.1 of Notebook~I \emph{Failure of traditional logic---the two gaps}. Before dealing with the two gaps, let us survey other aspects of this failure in connection with matters in the notes. There is first the arbitrariness and narrowness of the selection of the type of logical form to be investigated. The logical words selected are not completely pure (quantifiers are meshed with the connectives in the Aristotelian a, e, i, o forms), and they do not cover completely the propositional connectives, as G\" odel points out towards the end of Section 1.2.8 of the edited text (this is a matter in the sphere of functional completeness, treated by G\" odel in Section 1.1.8 of the edited text).

These words are also incomplete because they do not cover the quantifiers, as it is clearly shown by the envisaged axiomatization of Aristotelian syllogistic as a formal theory of propositional logic in Section 1.2.8 \emph{Classes and Aristotelian moods} of the edited text. (We have said in \S~16 of \cite{AD16} that {\L}ukasiewicz was working on such a presentation of Aristotelian syllogistic not much later than G\" odel in the Notre Dame course, if not at the same time, and they approached the subject in very much the same manner. This was a short while before the invasion of Poland and the outbreak of the Second World War, when G\" odel was back in Vienna.)

Relations of arity greater than the arity one, which properties have, are also left out in the Aristotelian approach, and this is another crucial incompleteness, as G\" odel says in the third paragraph of Section 1.2.1 \emph{First-order languages and valid formulas} of the edited text, because these relations are more important than properties ``for the applications of logic in mathematics and other sciences''. He also notes in the following paragraph of that section: ``Most of the predicates of everyday language are relations and not properties.''

Traditional logic deals exclusively with unary predicates, tied to properties, but it is incomplete also because it does not take all of them into account. Those which have an empty extension are left out, and this is detrimental for the use of logic, as G\" odel says in Section 1.2.6 \emph{Existential presuppositions} of the edited text. First, logic becomes dependent on empirical matters, and it also becomes impossible to use logic for answering in mathematics or elsewhere the question whether there is something that satisfies a property. Like leaving out zero in mathematics, it makes also the theory unnecessarily more complicated and uglier, if it does not end up in confusion and outright mistakes with the four wrong moods among the 19 moods, or with the conviction that no conclusion can be drawn where this is not the case (see the end of Section 1.2.8).

In Sections 1.1.4 and 1.1.1 G\" odel speaks about traditional logic failing to present logical laws as theorems of a deductive system. Occasionally in the past one heard boasts concerning this matter, which were based neither on a proof nor even a clear conception of the completeness in question. With a slight knowledge concerning classes and a few operations on them, which is based on a small, simple and intuitive fragment of propositional logic, of which Aristotelian logic is not aware, all the correct 15 Aristotelian moods are contained in a single formula (see Section 1.2.8). Decidability, which G\" odel calls completeness (see the remarks about the first gap below), is beyond the narrow horizon of traditional logic.

So taking into account several kinds of completeness, traditional logic failed to reach any of them. It is a complete failure. Traditional logic seems at first glance to be much present in G\" odel's course, but only in the Stoic's anticipatory discovery of connectives and propositional logical form there is something mentioned with approval---in the Aristotelian heritage nothing.

This complete failure of traditional logic in matters of completeness should certainly be taken into account in the explanation of the waste of the realm of traditional logic, which Greek mathematicians and most of the later ones ignored in their work, while some, like Descartes, condemned severely, centuries ago. G\" odel's measured but thorough condemnation is made in the light of various aspects of completeness, a modern theme developed by him with success in logic and mathematics.

G\" odel says that his chief aim in the propositional part of the course is to fill two gaps, solve two problems, which traditional logic failed to deal with, let alone solve (see the bottom of p.\ \textbf{3}.\ of Notebook ~0 and the bottom of p.\ \textbf{2}.\ and the top of p.\ \textbf{3}.\ of Notebook~I). The first is he says the problem of completeness of logical inference and logically true propositions, which he explicates as decidability, and the second is the problem of showing how all of them can be deduced from a small---he says ``minimum''---number of primitive laws. He considers the first problem solved by showing that the notion of tautology is decidable (see the bottom of p.\ \textbf{43}.~II of Notebook~I), and the second is solved by proving a deductive system for propositional logic complete (i.e.\ the sets of provable formulae and tautologies coincide; see the second p.\ \textbf{2}.\ towards the end of Notebook III). The two analogous problems for predicate logic are considered on p.\ \textbf{47}.\ of Notebook IV. G\" odel mentions that the second completeness problem was solved positively, and he gives indications concerning the negative solution of the first completeness problem, i.e.\ decidability, without entering into the proofs. He mentions the decidability of the monadic fragment.

For propositional logic G\" odel considers (at the end of p.\ \textbf{43}.~II of Notebook~I) that providing a decision procedure is even more than what is required for solving the first problem, as if he thought that providing concretely such a procedure (which is moreover easy to understand) is more than showing decidability nonconstructively. Usually today, completeness is understood in such a way that showing just the recursive enumerability of the set of tautologies is enough for it, and showing the recursiveness of that set is not compulsory. Decidability, i.e.\ the recursiveness of the set of tautologies, amounts to showing that both this set and its complement with respect to the set of formulae, are recursively enumerable, and so it makes sense to call decidability too \emph{completeness}; it is completeness in a stronger sense. G\" odel in any case distinguished the first problem, and the completeness involved in this problem, from the second problem of showing completeness with respect to a deductive system. From a positive solution of the first problem one can deduce the recursive enuberability of the set of logical laws, but that is not enough for the  second problem, which awaits to be solved. By not reducing proof theory to recursion theory, G\" odel took deduction as a separate important matter.

In that context, speaking of rules of inference G\" odel says: ``And of course we shall try to work with as few as possible.'' (p.\ \textbf{54}.\ of Notebook II) The ``of course'' in this sentence reflects something still in the air at the time the course was given, about which we spoke in Section~5 of \cite{DA16a}. G\" odel's advocacy of minimality is also related to the problem of independence of the axioms, with which he dealt in Section 1.1.12 of the edited text concerning his axiom system for propositional logic. This is besides completeness and decidability one of the main problems of logic, to which many investigations in set theory, in which G\" odel was involved too, were devoted.  We believe that his advocacy of minimality has however also to do with the following.

We said above several times that G\" odel used abbreviations very much. The economy brought by them is not only, so to speak, physical---with them less paper is needed, less ink, the reading is quicker. This economy is also of a conceptual kind. The Chinese way of writing need not have evolved from abbreviations, but it is as if it did. By moving away from the phonetic way of writing we do not represent concepts indirectly through the mediation of spoken words, which are represented in our writing. We represent the concepts directly. The written word ``two'' represents the number two indirectly through the mediation of the spoken word, while the figure 2 represents it directly. The written word ``prop.'' moves away from the representation of the spoken word ``proposition'' (and the context is practically always sufficient not to confuse it with the ``prop.'' of ``property''). The abbreviation ``log.'' in our example above stands for different words of different grammatical categories, as a Chinese character does. The Chinese way and the similar mathematical one are eminently reasonable, and bring advantages once one becomes accustomed to them.

Mathematical notation is far from phonetic. If something phonetic is still present in it, it is through abbreviations, or traces of abbreviations, often initial letters, as with functions being usually called $f$. There might be something mathematical in G\" odel's inclination towards abbreviations.

G\" odel's lectures end in the notes with Section 1.2.10 \emph{Type theory and paradoxes} of the edited text (pp.\ \textbf{127}.-\textbf{140}.\ of Notebook VI and \textbf{137}.-\textbf{157}.\ of Notebook VII, which precedes Section 1.2.11 \emph{Examples and samples of previous subjects}, which does not seem to be a lecture), where he presents Russell's paradox not explicitly as a set-theoretical matter, but through the predicate $\Phi$, read ``impredicable'', such that $\Phi(x)$ is equivalent with $\sim x(x)$ (see p.\ \textbf{142}.\ of Notebook VII; he follows there \cite{HA28}, Section IV.\S 4). Then on pp.\ \textbf{149}.-\textbf{156}.\ of Notebook VII he argues forcibly that self-reference (his term is ``self-reflexivity'') should not be blamed for the contradiction. He says that rejecting self-reference, which inspired Russell's theory of types, both in its ramified and in its simplified form, excludes many legitimate arguments based on self-reference, which do not lead to contradiction and are necessary for building set theory (pp.\ \textbf{155}.-\textbf{156}.\ of Notebook VII). The contradiction in the paradoxes is due to the illegitimacy of taking that there is a complete, achieved, totality of all objects---or to put it in other words, the impossibility to achieve completeness in the extensional realm.

It would be in G\" odel's style to write: ``Abbr.\ is an abbr''. The turn towards the conceptual here need not however be simply mathematical, because the self-reference involved could be akin not only to that made famous by \cite{Go31} but also to the intensional logic of the future (about which we said something in Section~5 of \cite{DA16}), where with legitimate self-reference the achievement of completeness is expected.

\clearpage

\pagestyle{empty}\makebox[1em]{} \clearpage

\noindent\textbf{\Large Contents}\label{Contents}
\pagestyle{myheadings}\markboth{LOGIC LECTURES}{CONTENTS}
\\*[1.5ex]
The page numbers below within square brackets, for example, [0: 1.] and [V: 73.1], are to be read as ``p.\ \textbf{1}.\ of Notebook~0'' and ``p.\ \textbf{73}.\textbf{1} of Notebook~V'' respectively. They come from the page numbers in the manuscript and the source version, which are rendered in boldface in the text below with the brackets $\mathbf{\llbracket \;\; \rrbracket}$. These page numbers refer to the first page of the section in question in the source version.

\begin{tabbing}
Editorial introduction\`\pageref{Editorial}
\\
Contents\`\pageref{Contents}
\\[2ex]
1\hspace{.4em} EDITED TEXT
\\[1ex]
1.1\hspace{.4em} Propositional logic\`\pageref{1}
\\*[.5ex]
1.1.1 \hspace{.4em} \=  Failure of traditional logic---the two gaps\quad [I: 1.]\`\pageref{1.1}
\\
1.1.2 \> Connectives\quad [I: 5.]\`\pageref{1.2}
\\
1.1.3 \> Tautologies\quad [I: 25. I]\`\pageref{1.3}
\\
1.1.4 \> Failure of traditional logic---the two gaps\quad [0: 1.]\`\pageref{1.4}
\\
1.1.5 \> Connectives\quad [0: 7.]\`\pageref{1.5}
\\
1.1.6 \> Tautologies\quad [0: 33.]\`\pageref{1.6}
\\
1.1.7 \> Decidability for propositional logic\quad [I: 38.1 II]\`\pageref{1.7}
\\
1.1.8 \> Functional completeness\quad [II: 33.]\`\pageref{1.8}
\\
1.1.9 \> Axiom system for propositional logic\quad [II: 52.]\`\pageref{1.9}
\\
1.1.10 \> Theorems and derived rules of the system for propositional logic\\*
 \`[II: 64.]\hspace{.25em}\pageref{1.10}
\\
1.1.11 \> Completeness of the axiom system for propositional logic [III: 11.]\`\pageref{1.11}
\\
1.1.12 \> Independence of the axioms\quad [III: 35.]\`\pageref{1.12}
\\
1.1.13 \> Remark on disjunctive and conjunctive normal forms\quad [III: 49.]\`\pageref{1.13}
\\
1.1.14 \> Sequents and natural deduction system\quad [III: second 1.]\`\pageref{1.14}
\\[2ex]
1.2\hspace{.4em} Predicate logic\`\pageref{2}
\\*[.5ex]
1.2.1 \> First-order languages and valid formulas\quad [IV: 24.]\`\pageref{2.1}
\\
1.2.2 \> Decidability and completeness in predicate logic\quad [IV: 46.]\`\pageref{2.2}
\\
1.2.3 \> Axiom system for predicate logic\quad [V: 58.]\`\pageref{2.3}
\\
1.2.4 \> Remarks on the term ``tautology'' and ``thinking machines''\\*
\`[V: 73.1]\hspace{.25em}\pageref{2.4}
\\
1.2.5 \> Theorems and derived rules of the system for predicate logic\\*
\`[V: 74.]\hspace{.25em}\pageref{2.5}
\\
1.2.6 \> Existential presuppositions\quad [V: 82.]\`\pageref{2.6}
\\
1.2.7 \> Classes\quad [V: 88.]\`\pageref{2.7}
\\
1.2.8 \> Classes and Aristotelian moods\quad [VI: 97.]\`\pageref{2.8}
\\
1.2.9 \> Relations\quad [VI: 107.]\`\pageref{2.9}
\\
1.2.10 \> Type theory and paradoxes\quad [VI: 126.]\`\pageref{2.10}
\\
1.2.11 \> Examples and samples of previous subjects\quad [VII: 1.]\`\pageref{2.11}
\\*[2ex]
Index\`\pageref{Index}
\\[2ex]
2\hspace{.4em} SOURCE TEXT
\\[1ex]
2.0\hspace{1em}\= Notebook 0\`\pageref{00}
\\
2.1 \> Notebook I\`\pageref{0I}
\\
2.2 \> Notebook II\`\pageref{0II}
\\
2.3 \> Notebook III\`\pageref{0III}
\\
2.4 \> Notebook IV\`\pageref{0IV}
\\
2.5 \>Notebook V\`\pageref{0V}
\\
2.6 \> Notebook VI\`\pageref{0VI}
\\
2.7 \> Notebook VII\`\pageref{0VII}
\end{tabbing}

\clearpage

\setcounter{page}{0}
\pagenumbering{arabic}

\chapter{EDITED TEXT}

\section{Propositional logic}\label{1}
\subsection{Failure of traditional logic---the two gaps}\label{1.1}

\pagestyle{myheadings}\markboth{EDITED TEXT}{NOTEBOOK I \;\;---\;\; 1.1.1\;\; Failure of traditional logic---the two gaps}

$\mathbf{\llbracket Notebook\; I \rrbracket}$ $\mathbf{\llbracket 1. \rrbracket}$ Logic is usually defined as the science whose object are the laws of correct thinking. According to this definition the central part of logic must be the theory of inference and the theory of logically true propositions [as e.g.\ the law of excluded middle] and in order to get acquainted with mathematical logic it is perhaps best to go \emph{in medias res} and begin with this central part.

Professor Menger\index{Menger, Karl} has pointed out in his introductory lecture that the treatment of these things in traditional logic\index{traditional logic} and in the current textbooks is unsatisfactory.\index{failure of traditional logic} Unsatisfactory from several standpoints. First from the standpoint of completeness. What the textbooks and also what Aristotle\index{Aristotle} gives is a more or less arbitrary selection of the infinity of the laws of logic, whereas in a systematic treatment as is given in mathematical logic we shall have to develop methods which allow $\mathbf{\llbracket 2. \rrbracket}$ us to obtain all possible logically true propositions and to decide of any given proposition whether or not it is logically true or of an inference whether it is correct or not. But secondly the classical treatment is also unsatisfactory as to the question of reducing the\footnote{If the crossed out ``inf.'', which appears at this place in the manuscript, is interpreted instead as underlined, which is possible, this might be taken as an abbreviation for ``infinity''. Above in this paragraph and at the beginning of p.\ \textbf{2}.\ of Notebook~0 one finds the phrase ``the infinity of the laws of logic''.} laws of logic to a certain number of primitive laws from which they can be deduced. Although it is sometimes claimed that everything can be deduced from the three fundamental laws of contradiction, excluded middle and identity or from the modus Barbara this claim has never been proved or even clearly formulated in traditional logic.

The chief aim in the first part of these lectures will be to fill those two gaps [solve those two problems in a satisfactory way], i.e.\ to give as far as possible a complete theory of logical inference and logically true propositions, $\mathbf{\llbracket 3. \rrbracket}$ complete at least for a certain very wide domain of propositions, and to show how they can be reduced to a certain number of primitive laws.

The theory of syllogisms\footnote{or syllogistic} as presented in the current textbooks is usually divided into two parts:
\begin{itemize}
\item[1.] The Aristotelian figures and moods\index{Aristotelian figures and moods} of inference including the inferences with one premise (e.g.\ contradiction),
\item[2.] inferences of an entirely different kind which are treated under the heading of hypothetical disjunctive conjunctive inferences and which seem to be a Stoic addition\index{Stoic addition to traditional logic} to the Aristotelian figures.
\end{itemize}
Let us begin with the syllogisms of the second kind which turn out to be much more fundamental. We have for instance the modus ponendo ponens\index{ponendo ponens}.\index{modus ponens}
\vspace{1ex}

$\mathbf{\llbracket 4. \rrbracket}$ From the two premises
\begin{tabbing}
\hspace{1.7em}\=1. \= If Leibnitz has invented the infinitesimal calculus he was a great\\
\` mathematician,\\
\>2. \> Leibnitz has invented the infinitesimal calculus,\\[.5ex]
we conclude\\[.5ex]
\>\>Leibnitz was a great mathematician.
\end{tabbing}

Generally, if $p$ and $q$ are arbitrary propositions and if we have the two premises
\begin{tabbing}
\hspace{1.7em}\= 1. \= If $p$ so $q$,\\*[.3ex]
\>2. \> $p$,\\[.5ex]
we can conclude\\[.5ex]
\>\> $q$.
\end{tabbing}
Or take a disjunctive inference tollendo ponens.\index{disjunctive syllogism}\index{tollendo ponens} If we have the two premises
\begin{tabbing}
\hspace{1.7em}\= 1. \= Either $p$ or $q$,\\*[.3ex]
\>2. \>Not $p$,\\[.5ex]
we can conclude\\[.5ex]
\>\> $q$.
\end{tabbing}

It is possible to express this is syllogism by one logically true proposition as follows:
\begin{tabbing}
\hspace{1.7em}If either $p$ or $q$ and if not-$p$ then $q$.
\end{tabbing}
This whole statement will be true whatever $p,q$ may be.

Now what is the most striking caracter of these inferences which distinguishes them from the Aristotelian syllogistic figures? It is this: $\mathbf{\llbracket 5. \rrbracket}$ that in order to make those inferences it is not necessary to know anything about the structure of $p$ and $q$. $p$ or $q$ (may themselves be disjunctive or hypothetical propositions), they may be affirmative or negative propositions, or they may be simple or as complicated as you want; all this is indifferent for this syllogism, i.e.\ only propositions as a whole occur in it and it is this fact that makes this kind of syllogism simpler and more fundamental than the Aristotelian. The law of contradiction and excluded middle would be other examples of logical laws of this kind. Because e.g.\ the law of excluded middle say for any proposition $p$ either $p$ or $\sim p$ is true and this quite independently of the structure of $p$. With the Aristotelian logical syllogism it is of course quite different; they depend on the structure of the propositions involved, e.g.\ in order to apply the mood Barbara you must know e.g.\ that the two premises are general affirmative propositions.

\subsection{Connectives}\label{1.2}
\pagestyle{myheadings}\markboth{EDITED TEXT}{NOTEBOOK I \;\;---\;\; 1.1.2\;\; Connectives}

Now the theory $\mathbf{\llbracket 6. \rrbracket}$ of logically true propositions and logical inferences in which only propositions as a whole occur is called calculus of propositions.\index{calculus of propositions} \index{propositional calculus} In order to subject it to a systematic treatment we have first to examine more in detail the connections\index{connections (connectives)}\index{connectives}\footnote{``connective'' would be more suitable than ``connection'', but G\" odel does not seem to have used that word at that time (see the last footnote on p.\ \textbf{10}. of Notebook~0); ``connection'' is put at analogous places below.} between propositions which can occur in there inferences, i.e.\ the or, and, if\ldots\ so, and the not. One has introduced special symbols to denote them, in fact there are two different symbolisms for them, the Russell\index{Russell, Bertrand} and the Hilbert\index{Hilbert, David} symbolism. I shall use in these lectures Russell's\index{Russell, Bertrand} symbolism. In this not is denoted by $\sim$, and by a dot $.\:$, or by $\vee$ and the if\ldots\ so by $\supset$, $\mathbf{\llbracket 7. \rrbracket}$ i.e.\ if $p,q$ are arbitrary propositions then $\sim p$ means $p$ is false,\index{negation} $p\: .\: q$ means both $p$ and $q$ are true,\index{conjunction} $p \vee q$ means at least one of the propositions $p,q$ is true, either both are true or one is true and the other one false.\index{disjunction} This is different from the meaning that is given to the or in traditional logic. There we have to do with the exclusive or, in Latin \emph{aut}\ldots\ \emph{aut}\index{aut aut@\emph{aut}\ldots\ \emph{aut}}, which means that exactly one of the two propositions $p,q$ is true and the other one is false, whereas this logical symbol for or has the meaning of the Latin \emph{sive}\ldots\ \emph{sive}\index{sive sive@\emph{sive}\ldots\ \emph{sive}}, i.e.\ one of the two propositions is true where it is not excluded that both are true. The exclusive or as we shall see later can be expressed by a combination of the other logistic symbols, but one has not introduced a proper symbol for it because it turns out not to be as fundamental as the or in the sense of \emph{sive}\ldots\ \emph{sive}; $\mathbf{\llbracket 8. \rrbracket}$ it is not very often used. The next symbol is the $\supset$. If $p,q$ are two propositions $p \supset q$ means if $p$ so $q$, i.e.\ $p$ implies $q$.\index{implication} Finally we introduce a fifth connection $p \equiv q$ ($p$ equivalent to $q$) which means both $p \supset q$ and $q \supset p$.\index{equivalence}

The five connections introduced so far are called respectively negation,\index{negation} conjunction,\index{conjunction} disjunction,\index{disjunction} implication,\index{implication} equivalence,\index{equivalence} and all of them are called connections\index{connections (connectives)}\index{connectives} or operations of the calculus of propositions.\index{operation of the calculus of propositions (connectives)} Conjunction and disjunction are also called logical product\index{product, logical (conjunction)}\index{logical product (conjunction)} and logical sum\index{sum, logical (disjunction)}\index{logical sum (disjunction)} respectively.
All of the mentioned logical operations\ excluding negation are operations with two arguments, i.e.\ they form a new proposition out of two given ones, for example, $p\vee q$.
Only the negation is an operation with one argument forming a new proposition $\sim p$ out of any single given proposition.

Not only the operations $\supset$, $\vee$ and $.$ are called implication, disjunction and conjunction, but also an expression of the form $p\supset q$, $p\vee q$ is called an implication etc., where $p,q$ may again be expressions involving again $\supset$, $\vee$ etc. and $p,q$ are called respectively first and second member. Of course if $p$ and $q$ are propositions then $\sim p$, $\sim q$, $p\vee q$, ${p\: .\:q}$ and $p\supset q$ are also propositions and hence to them the operations of the calculus of propositions can again be applied, so as to get more complex expressions, e.g.\ ${p \vee (q\: .\:r)}$, either $p$ is true or $q$ and $r$ are both true.

The disjunctive inference\index{disjunctive syllogism} I mentioned before would read in this symbolism as follows: $[(p\vee q)\: .\:\sim p]\supset q$. You see in more complex expressions as this one brackets\index{brackets (parentheses)}\index{parentheses} have to be used exactly as in algebra in order to indicate the order in which the operations have to be applied. E.g.\ if I put the round brackets in this expression like this ${p \vee (q\: .\:\sim p)}$, it would have a different meaning, namely either $p$ is true or $q$ and $\sim p$ are both true.

There is an interesting $\mathbf{\llbracket 9. \rrbracket}$ remark due to {\L}ukasiewicz\index{Lukasiewicz@{\L}ukasiewicz, Jan} that one can dispense with the brackets if one writes the operational symbols $\vee$, $\supset$ etc.\ always in front of the propositions\index{Polish notation} to which they are applied, e.g.\ $\supset\! p\, q$ instead of ${p\supset q}$. Then e.g.\ the two different possibilities for the expression in square brackets would be distinguished automatically because the first would be written as follows $.\vee p\, q\sim p$; the second would read $\vee p\: .\: q\sim p$, so that they differ from each other without the use of brackets as you see and it can be proved that it is quite generally so. But since the formulas in the bracket notation are more easily readable I shall stick to this notation and put the operational symbols in between the propositions.

You know in algebra one can spare many brackets by the convention that the $\mathbf{\llbracket 10. \rrbracket}$ multiplication connects stronger than addition; e.g.\ $a\cdot b+c$ means $(a\cdot b)+c$ and not $a\cdot (b+c)$. We can do something similar here by stipulating an order of the strength\index{strength, of binding for connectives} in which the logical symbols bind, so that:
\begin{itemize}
\item[1.] the $\sim$ (and similarly any operation with just one proposition as argument) connects stronger than any operation with two arguments, as $\vee$, $\supset$ and $.$, so that $\sim p\vee q$ means $(\sim p) \vee q$ and not $\sim (p\vee q)$;
\vspace{-1ex}
\item[2.] the disjunction and conjunction bind stronger than implication and equivalence, so that e.g.\ ${p\vee q \supset r\: .\:s}$ means ${(p\vee q) \supset (r\: .\:s)}$ and not perhaps $p\vee [(q \supset r)\: .\:s]$.
\end{itemize}
A third convention consists in leaving out brackets in such expressions as $(p\vee q)\vee r$ exactly as in $(a+b)+c$. A similar convention is made for $.\;$.

After those merely symbolic conventions the next thing we have to do is to examine in more detail the meaning of the operations of the calculus of propositions. $\mathbf{\llbracket 11. \rrbracket}$ Take e.g.\ disjunction $\vee$. If any two propositions $p,q$ are given $p\vee q$ will again be a proposition. Hence the disjunction is an operation which applied to any two propositions gives again a proposition. But now (and this is the decisive point) this operation is such that the truth or falsehood of the composite proposition $p\vee q$ depends in a definite way on the truth or falsehood of the constituents $p,q$. This dependence can be expressed most clearly in the form of a table\index{truth table} as follows: let us form three columns, one headed by $p$, one by $q$, one by $p\vee q$, and let us write $+$ for true and $-$ for false. Then for the proposition $p\vee q$ we have the following four possibilities:
\begin{center}
\begin{tabular}{ c|c|c|c }
 $p$ & $q$ & $p\vee q$ & $p \: o\: q$\\[.5ex]
 + & + & + & $-$ \\
 + & $-$ & + & + \\
 $-$ & + & + & + \\
 $-$ & $-$ & $-$ & $-$
 \end{tabular}
\end{center}
Now $\mathbf{\llbracket 12. \rrbracket}$ for each of these four cases we can determine whether ${p\vee q}$ will be true or false, namely since ${p\vee q}$ means that one or both of the propositions $p,q$ are true it will be true in the first, second and third case, and false in the last case. And we can consider this table as the most precise definition of what $\vee$ means.

It is usual to call truth\index{Truth} and falsehood\index{Falsehood} the truth values,\index{truth values} so there are exactly two truth values, and say that a true proposition has the truth value ``truth'' (denoted by $+$) and a false proposition has the truth value ``false'' (denoted by $-$), so that any proposition has a uniquely determined truth value. The truth table then shows how the truth value of the composite expressions depends on the truth value of the constituents. The exclusive or would have another truth table; namely if we denote it by $o$ for the moment we have that $p\: o\: q$ is false if both $p$ and $q$ are true, and it is false if both are false but true in the two other cases. The operation $\sim$ $\mathbf{\llbracket 13. \rrbracket}$ has of course the following truth table:\begin{center}
\begin{tabular}{ c|c}
 $p$ & $\sim p$ \\[.5ex]
 + & $-$ \\
 $-$ & +
 \end{tabular}
\end{center}
Here we have only two possibilities: $p$ true or $p$ wrong, and in the first case we have that not-$p$ is wrong while in the second it is true. Also the truth table for $.$ can easily be determined:
\begin{center}
\begin{tabular}{ c|c|c }
 $p$ & $q$ & $p\: .\: q$\\[.5ex]
 + & + & + \\
 + & $-$ & $-$ \\
 $-$ & + & $-$ \\
 $-$ & $-$ & $-$
 \end{tabular}
\end{center}
(I think I will leave that to you.)

A little more difficult is the question of the truth table for $\supset$.  $\mathbf{\llbracket 14. \rrbracket}$ $p\supset q$ was defined to mean ``If $p$ is true $q$ is also true''. So let us assume that for two given propositions $p,q$ we know that $p\supset q$ is true, i.e.\ assume that we know ``If $p$ then $q$'' but nothing else. What can we conclude then about the possible truth values of $p$ and $q$?
\begin{center}
\text{Assumption $p \supset q$}
\end{center}
\begin{center}
\begin{tabular}{@{}l}
 $\left.
  \begin{tabular}{@{}*2{p{0.3cm}}}
  $p$ & $q$ \\[.5ex]
  $-$ & + \\
  $-$ & $-$ \\
  + & +
  \end{tabular}
 \right\}\text{possible truth values for $p,q$}$\\[5ex]
 $\left.
 \begin{tabular}{@{}*2{p{0.3cm}}}
  + & $-$
 \end{tabular}
 \right\}\text{impossible}$
\end{tabular}
\end{center}

\noindent First it may certainly happen that $p$ is false because the assumption statement ``If $p$ then $q$'' says nothing about the truth or falsehood of $p$. And in this case where $p$ is false $q$ may be true as well as false because the assumption ``If $p$ then $q$'' says nothing about what happens to $q$ if $p$ is false but only if $p$ is true. So we have both possibilities: $p$ false $q$ true, $p$ false $q$ false. Next we have the possibility that $p$ is true. $\mathbf{\llbracket 15. \rrbracket}$ But in this case owing to the assumption $q$ must also be true. So that the possibility $p$ true $q$ false is excluded and we have only this third possibility $p$ true $q$ true, and this possibility may of course really happen. So from the assumption $p\supset q$ it follows that either one of the first three cases happens. But we have also vice versa: If one of the first three possibilities of the truth values is realized then $(p\supset q)$ is true. Because let us assume we know that one of the three cases written down is realized. I claim then we know also: ``If $p$ is true then $q$ is true''. If $p$ is true only the third of the three possibilities can be realized (in all the others $p$ is false), but in this third possibility $q$ is true. $\mathbf{\llbracket 16. \rrbracket}$ So we see that the statement $p \supset q$ is exactly equivalent with the statement that one of the three marked cases for the distribution of truth values is realized, i.e.\ $p\supset q$ is true in each of the three marked cases and false in the last case. So we have obtained a truth table for implication.
However there are two important remarks about it namely:

1. Exactly the same truth table can also be obtained by a combination of operations introduced previously, namely $\sim p \vee q$ has the same truth table
\begin{center}
\begin{tabular}{ c|c|c|c }
 $p$ & $q$ & $\sim p$ & $\sim p\vee q$\\[.5ex]
 $-$ & $-$& + & + \\
 $-$& +& + & + \\
 +& $-$& $-$ & $-$ \\
 +& +& $-$ & +
 \end{tabular}
\end{center}
$\mathbf{\llbracket 17. \rrbracket}$ Since $p\supset q$ and $\sim p \vee q$ have the same truth table they will be equivalent, i.e.\ whenever the one expression is true the other one will also be true and vice versa. This makes it possible to define $p\supset q$ by $\sim p \vee q$ and this is the standard way of introducing implication in mathematical logic.

2. The second remark about implication is this. We must be careful not to forget that $p \supset q$ was understood to mean simply ``If $p$ then $q$'' and only this made the construction of the truth table possible. We have deduced the truth table for implication from the assumption that $p \supset q$ means ``If $p$ then $q$'' and nothing else. There are other meanings $\mathbf{\llbracket 18. \rrbracket}$ perhaps even more suggested by the term implication for which our truth table would be completely inadequate. E.g.\ $p\supset q$ could be given the meaning: $q$ is a logical consequence of $p$, i.e.\ $q$ can be derived from $p$ by means of a chain of syllogisms.

This kind of implication is usually called strict implication\index{strict implication} and denoted in this way $\prec$ and the implication $p\supset q$ defined before is called material implication\index{material implication} if it is to be distinguished. Now it is easy to see that our truth table is false for strict implication. In order to prove that consider the first line of a supposed such table
\begin{center}
\begin{tabular}{ c|c|c }
 $p$ & $q$ & $p \prec q$\\[.5ex]
 + & + & \\
 & & \\
 & & \\
 & &
 \end{tabular}
\end{center}
where $p$ and $\mathbf{\llbracket 19. \rrbracket}$ $q$ are both true and ask what will be the truth value of $p \prec$ strictly $q$. It is clear that this truth value will not be uniquely determined. For take e.g.\ for $p$ the proposition ``The earth is a sphere'' and for $q$ ``The earth is not a disk''. Then $p$ and $q$ are both true and $p\prec q$ is also true because from the proposition that the earth is a sphere it follows by logical inference that it is not a disk; on the other hand if you take for $p$ again the same proposition and for $q$ ``France is a republic'' then again both $p$ and $q$ are true but $p \prec q$ is wrong. $\mathbf{\llbracket 20. \rrbracket}$ So we see the truth value of $p \prec q$ is not uniquely determined by the truth values of $p$ and $q$, and therefore no truth table exists. Such connections\footnote{The plural of ``connective'' would be more suitable (see the footnote on p.\ \textbf{6}.\ of the present Notebook~I).} for which no truth table exists are called intensional\index{intensional connections (connectives)} as opposed to extensional\index{extensional connections (connectives)} ones for which they do exist. The extensional connections are called also truth functions.\index{truth functions}

So we see the implication which we introduced does not mean logical consequence. Its meaning is best given by the simple ``if then''\index{if\ldots\ then} which has much wider significance than just logical consequ\-ence.\index{logical consequence} E.g.\ if I say ``If he cannot come he will telephone to you'', that has nothing to do with logical relations between $\mathbf{\llbracket 21. \rrbracket}$ his coming and his telephoning, but it simply means he will either come or telephone which is exactly the meaning expressed by the truth table. Now the decisive point is that we don't need any other kind of implication besides material in order to develop the theory of inference because in order to make the conclusion from a proposition $p$ to a proposition $q$ it is not necessary to know that $q$ is a logical consequence of $p$. It is quite sufficient to know ``If $p$ is true $q$ is true''.
Therefore I shall use only material implication, at least in the first half of my lectures, and use the terms ``implies''\index{implies} and ``it follows''\index{follows} only in this sense.

$\mathbf{\llbracket 22. \rrbracket}$ This simplifies very much the whole theory of inference because material implica\-tion defined by the truth table is a much simpler notion. I do not want to say by this that a theory of strict implication may not be interesting and important for certain purposes; in fact I hope to speak about it later on in my lectures. But its theory belongs to an entirely different part of logic than that with which we are dealing at present, namely it belongs to the logic of modalities.\index{logic of modalities}

Now I come to some apparently paradoxical\index{paradoxes, of material implication} consequences of our definition of implication whose paradoxicality however disappears if we remember that implication does not mean logical consequence. Namely if we look at the truth table for $p\supset q$ we see at once that $p\supset q$ is always true if $q$ is true whatever $p$ may be. So that means a true proposition is implied by any proposition. Secondly we see that $p\supset q$ is always true if $p$ is false whatever $q$ $\mathbf{\llbracket 23. \:I \rrbracket}$ may be; i.e.\ a false proposition implies any proposition whatsoever. In other words: ``An implication with true second member is true (whatever the first member may be) and an implication with a false first member is always true (whatever the second member may be).''
Or written in formulas this means $q \supset (p\supset q)$, $\sim p \supset (p\supset q)$. Both of these formulas are also immediate consequences of the fact that $p\supset q$ is equivalent with $\sim p \vee q$ because $\sim p \vee q$ says exactly either $p$ is false or $q$ is true, so it will always be true if $p$ is false and if $q$ is true whatever the other proposition may be. These formulas are rather unexpected and if we apply them to special cases we get strange consequences. E.g.\ $\mathbf{\llbracket 24. \rrbracket}$ ``The earth is not a sphere'' implies that France is a republic, but it also implies that France is not a republic because a false proposition implies any proposition whatsoever. Similarly the proposition ``France is a republic'' is implied by any other proposition whatsoever, true or false. But these consequences are only paradoxical if we understand implication to mean logical consequence. For the ``if\ldots\ so'' meaning they are quite natural, e.g.\  $q\supset(p\supset q)$ means: If $q$ is true then $q$ is true also if $p$ is true, and $\sim p\supset(p\supset q)$ If we have a false proposition $p$ then if $p$ is true anything is true. $\mathbf{\llbracket 25. \: I \rrbracket}$ Another of these so called paradoxical consequences is this $(p\supset q)\vee (q\supset p)$, i.e.\ of any two arbitrary propositions one must imply the other one. That it must be so is proved as follows: $q$ must be either true or false; if $q$ is true the first member of the disjunction is true and if $q$ is false the second member is true because a false proposition implies any other. So (one of the two members of the implication is true) either $p\supset q$ or $q\supset p$ in any case.

\subsection{Tautologies}\label{1.3}
\pagestyle{myheadings}\markboth{EDITED TEXT}{NOTEBOOK I \;\;---\;\; 1.1.3\;\; Tautologies}

We have here three examples of logically true formulas,\index{logically true formulas, of propositional logic (tautology)}\footnote{see pp.\ \textbf{23}.\,\textbf{I}-\textbf{25}.\,\textbf{I} of the present Notebook~I} i.e.\ formulas which are true whatever the propositions $p,q$ may be. Such formulas are called tautological\index{tautological formula, of propositional logic}\index{tautology, of propositional logic} and it is exactly the chief aim of the calculus of propositions to investigate those tautological formulas.

I shall begin with discussing a few more examples of such logically true propositions before going over to general considerations. $\mathbf{\llbracket 26. \: I \rrbracket}$ We have at first the traditional hypothetical and disjunctive inferences which in our notation read as follows:
\begin{itemize}
\item[1.] $p\, .\:(p\supset q)\supset q$ \quad ponendo ponens\index{ponendo ponens}\index{modus ponens} \quad
\vspace{-1ex}
\item[[2.] $\sim q\, .\:(p\supset q)\supset \;\sim p$ \quad tollendo tollens\index{tollendo tollens}]
\vspace{-1ex}
\item[3.] $(p\vee q)\, .\: \sim q\supset p$ \quad tollendo ponens\index{tollendo ponens}\\
disjunctive ponendo tollens\index{ponendo tollens} does not hold for the not exclusive $\vee$ which we have
\vspace{-1ex}
\item[4.] The inference which is called dilemma\index{dilemma}\\
  $(p\supset q)\, .\:(r\supset q)\supset(p\vee r\supset q)$
\end{itemize}

\subsection{Failure of traditional logic---the two gaps}\label{1.4}
\pagestyle{myheadings}\markboth{EDITED TEXT}{NOTEBOOK 0 \;\;---\;\; 1.1.4 \;\; Failure of traditional logic---the two gaps}

$\mathbf{\llbracket Notebook\; 0 \rrbracket}$ $\mathbf{\llbracket 1. \rrbracket}$ Logic is usually defined as the science of the laws of correct thinking. According to this definition the central part of logic must be the theory of inference and the theory of logically true propositions. By a logically true proposition I mean a proposition which is true for merely logical reasons as e.g.\ the law of excluded middle, which says that for any proposition $p$ either $p$ or $\sim p$ is true. I intend to go in medias res right away and to begin with this central part.

As Professor Menger\index{Menger, Karl} has pointed out in his introductory lecture the treatment of these things, inferences and logically true propositions, in traditional logic\index{traditional logic}\footnote{The following text is here crossed out in the manuscript: ``and in most of the current textbooks''.} is unsatisfactory in some respect.\index{failure of traditional logic} First with respect to completeness. What the $\mathbf{\llbracket 2. \rrbracket}$ traditional logic gives is a more or less arbitrary selection from the infinity of the laws of logic, whereas in a systematic treatment we shall have to develop methods which allow us to obtain as far as possible all logically true propositions and methods which allow to decide of arbitrary given propositions whether or not they are logically true. But the classical treatment is unsatisfactory also in another respect; namely as to the question of reducing the laws of logic to a certain number of primitive laws from which $\mathbf{\llbracket 3. \rrbracket}$ all the others can be deduced. Although it is sometimes claimed that everything can be deduced from the law of contradiction or from the first Aristotelian figure, this claim has never been proved or even clearly formulated in traditional logic.

The chief aim in the first part of this seminary will be to fill these two gaps of traditional logic, i.e.\ 1.\ to give as far as possible a complete theory of logical inference and of logically true propositions and 2.\ to show how all of them can be deduced from a minimum number of primitive laws.

$\mathbf{\llbracket 4. \rrbracket}$ The theory of inference as presented in the current textbooks is usually divided into two parts:
\begin{itemize}
\item[1.] The Aristotelian figures and moods\index{Aristotelian figures and moods} including the inferences with one pre\-mise, i.e.\ conversion, contraposition etc.
\item[2.] Inferences of an entirely different kind, which are treated under the heading of hypothetical disjunctive conjunctive inference, and which are a Stoic addition\index{Stoic addition to traditional logic} to the Aristotelian figures.
\end{itemize}

Let us begin with these inferences of the second kind, which turn out to be more fundamental than the Aristotelian figures.

Take the following examples of the disjunctive inference tollendo ponens:\index{disjunctive syllogism}\index{tollendo ponens}

\vspace{1ex}

\noindent $\mathbf{\llbracket 5. \rrbracket}$ From the two premises
\begin{tabbing}
\hspace{1.7em}\=1. \= Nero was either insane or a criminal,\\[.3ex]
\>2. \> Nero was not insane, \\[.5ex]
we can conclude\\[.5ex]
\>\> Nero was a criminal.
\end{tabbing}
\begin{tabbing}
\hspace{1.7em}\=1. \= Today is either Sunday or a holiday,\\[.3ex]
\>2. \> Today is not Sunday,\\[.5ex]
\>\> Today is a holiday.
\end{tabbing}

Generally, if $p,q$ are two arbitrary propositions and we have the two premises
\begin{tabbing}
\hspace{1.7em}\=1. \= Either $p$ or $q$, \\[.3ex]
\>2. \> not-$p$, \\[.5ex]
we can conclude\\[.5ex]
\>\> $q$.
\end{tabbing}
It is possible to express this syllogism by one logically true proposition as follows:
\begin{tabbing}
\hspace{1.7em}``(If either $p$ or $q$ and not-$p$) then $q$''
\end{tabbing}
This whole proposition under quotation marks will be true whatever the proposi\-tions $p$ and $q$ may be.

$\mathbf{\llbracket 6. \rrbracket}$ Now what is the caracter of this inference which distinguishes it from the Aristotelian figures? It is this that in order to make this inference it is not necessary to know anything about the structure of the propositions $p$ and $q$. $p$ and $q$ may be affirmative or negative propositions, they may be simple or complicated, they may themselves be disjunctive or hypothetical propositions; all this is indifferent for this syllogism, i.e.\ only propositions as a whole occur in it, and it is this caracter that makes this kind of syllogism simpler and more fundamental than e.g.\ the Aristotelian $\mathbf{\llbracket 7. \rrbracket}$ figures, which depend on the structure of the propositions involved. E.g.\ in order to make an inference by mood Barbara you must know that the two premises are universal affirmative. Another example of a logical law in which only propositions as a whole occur would be the law of excluded middle, which says: For any proposition $p$ either $p$ or not-$p$ is true.

\subsection{Connectives}\label{1.5}
\pagestyle{myheadings}\markboth{EDITED TEXT}{NOTEBOOK 0 \;\;---\;\; 1.1.5\;\; Connectives}

Now the theory of those laws of logic in which only propositions as a whole occur is called calculus of propositions,\index{calculus of propositions}\index{propositional calculus} and it is exclusively with this part of mathematical logic that we shall have $\mathbf{\llbracket 8. \rrbracket}$ to do in the next few lectures. We have to begin with examining in more detail the connections\index{connections (connectives)}\index{connectives} between propositions which occur in the inferences concerned, i.e.\ the or, and, if, not. One has introduced special symbols to denote them. ``Not'' is denoted by a circumflex,``and'' by a dot, ``or'' by a kind of abbreviated v (derived from vel), ``if then'' is denoted by this symbol similar to a horseshoe:
\begin{tabbing}
\hspace{1.7em}\= not\index{not}\index{negation} \hspace{5em} \=$\sim$ \quad which is an abbreviated N \quad \= $\sim p$\\[.5ex]
\>and\index{and}\index{conjunction}\> $\, .$ \> $p\: .\:q$\\[.5ex]
\>or\index{or}\index{disjunction}\> $\vee$\> $p\vee q$\\[.5ex]
\>if\ldots\ then\index{if\ldots\ then}\index{implication} \> $\supset$\> $p\supset q$\\[.5ex]
\>equivalent\index{equivalence} \> $\equiv$\> $p\equiv q$
\end{tabbing}
i.e.\ if $p$ and $q$ are arbitrary propositions $\sim p$ means $p$ is false, $p\: .\: q$ means both $p$ and $q$ is true, $p\vee q$ means either $p$ or $q$, $p\supset q$ means if $p$ then $q$, or in other words $p$ implies $q$.\footnote{Here one finds in the manuscript a broken sentence beginning with: ``So if e.g.\ $p$ is the proposition today it will rain and $q$ is the proposition tomorrow it will snow then'', of which the words after $q$ are on p.\ \textbf{9}.\ of the present Notebook~0.}

$\mathbf{\llbracket 9. \rrbracket}$ About the ``or'' namely, this logical symbol means that at least one of the two propositions $p,q$ is true but does not exclude the case where both are true, i.e.\ it means one or both of them are true, whereas the ``or'' in traditional logic is the exclusive ``or'' which means that exactly one of the two propositions $p,q$ is true and the other one false. Take e.g.\ the sentence ``Anybody who has a salary or interests from capital is liable to income tax''. Here the ``or'' is meant in the sense of the logical ``or'', because someone who has both is also liable to income tax. On the other hand, in the proposition ``Any number except $1$ is either greater or smaller than 1'' we mean the exclusive ``or''.\index{exclusive ``or''} This exclusive ``or'' corresponds to the Latin \emph{aut},\index{aut@\emph{aut}} the logical ``or'' to the Latin \emph{vel}\index{vel@\emph{vel}}.\footnote{Here one finds in the manuscript an apparently broken sentence beginning with: ``As we shall see later''.}

The exclusive ``or'' can be expressed by a combination $\mathbf{\llbracket 10. \rrbracket}$ of the other logical symbols, but no special symbol has been introduced for it, because it is not very often used. Finally, I introduce a fifth connection, the so called ``equivalence'' denoted by three horizontal lines. $p\equiv q$ means that both $p$ implies $q$ and $q$ implies $p$. This relation of equivalence would hold e.g.\ between the two propositions: ``Tomorrow is a weekday'' and ``Tomorrow is not a holiday''.\footnote{Here one finds in the manuscript an incomplete sentence: ``because we have If\ldots\ but also vice versa''.}

The five notions which we have introduced so far are called respectively operation of negation,\index{negation} conjunction,\index{conjunction} disjunction,\index{disjunction} implication,\index{implication} equivalence.\index{equivalence} By a common name they are called functions of the calculus of propositions\index{functions of the calculus of propositions (connectives)} or\footnote{text missing in the manuscript; ``connective'' would be suitable, but G\" odel does not seem to have used that word at that time. In the preceding paragraph and at the beginning of the present Section 1.1.5 he has ``connection'' instead.} Disjunc\-tion is also called $\mathbf{\llbracket 11. \rrbracket}$ logical sum\index{sum, logical (disjunction)}\index{logical sum (disjunction)} and conjunction logical product\index{product, logical (conjunction)}\index{logical product (conjunction)} because of certain analogies with the arithmetic sum and the arithmetic product. A proposition of the form $p\vee q$ is called a disjunction and $p,q$ its first and second member; similarly a proposition of the form $p\supset q$ is called an implication and $p,q$ its first and second member, and similarly for the other operations. Of course, if $p,q$ are proposi\-tions, then $\sim p$, $\sim q$, $p\vee q$, $p\: .\:q$, $p\supset q$ are also proposi\-tions and therefore to them the functions of the calculus of propositions can again be applied so as to get more complicated expres\-sions; e.g.\ $p\vee(q\, .\:r)$, which would mean: Either $p$ is true or $q$ and $r$ are both true.

The disjunctive syllogism\index{disjunctive syllogism} $\mathbf{\llbracket 12. \rrbracket}$ I mentioned before can be expressed in our symbolism as follows: $[(p\vee q)\, .\:\sim q]\supset p$. You see in more complicated expressions as e.g.\ this one brackets\index{brackets (parentheses)}\index{parentheses} have to be used exactly as in algebra to indicate in what order the operations have to be carried out. If e.g.\ I put the brackets in a different way in this expression, namely like this $(p\vee q)\, .\:r$, it would mean something entirely different, namely it would mean either $p$ or $ q$ is true and in addition $r$ is true.

There is an interesting remark due to the Polish logician {\L}ukasiewicz,\index{Lukasiewicz@{\L}ukasiewicz, Jan} namely that one can dispense entirely with brackets\index{Polish notation} if one writes the $\mathbf{\llbracket 13. \rrbracket}$ operational symbols $\vee$, $\supset$ etc.\ always in front of the proposition to which they are applied, e.g.\ $\supset p\,q$ instead of $p\supset q$. Incidentally, the word ``if'' of ordinary language is used in exactly this way. We say e.g.\ ``If it is possible I shall do it'' putting the ``if'' in front of the two propositions to which we apply it. Now in this notation where the operations are put in front the two different possibilities of this expression $p\vee q \: .\: r$ would be distinguished automatically without the use of brackets because the second would read $.\vee p\,q\,r$, with ``or'' applied to $p,q$ and the ``and'' applied to this formula and $r$, whereas the first would read ``and'' applied to $q$, $r$ and the $\vee$ applied to $p$ and this formula $\vee p \, .\: qr$. As you see, these two formulas differ from each other without the use of brackets and it can be shown that $\mathbf{\llbracket 14. \rrbracket}$ it is quite generally so. Since however the formulas in the bracket notation are more easily readable I shall keep the brackets and put the operation symbol between the propositions to which they are applied.

You know in algebra one can save many brackets by the convention that multiplication is of greater force than addition, and one can do something similar here by stipulating an order of force\index{force, of binding for connectives} between the operations of the calculus of propositions, and this order is to be exactly the same in which I introduced them, namely
\[
\sim\, .\:\vee\;{\supset \atop \equiv}
\]
No order of force is defined for $\supset\, \equiv$, they are to have equal force. Hence
\begin{tabbing}
$\mathbf{\llbracket 15. \rrbracket}$\\*[1ex]
\hspace{1.7em}\= $\sim p\vee q$ \quad me\=ans \quad \=$(\sim p)\vee q$ \hspace{1.1em} n\=ot \quad \=$\sim(p\vee q)$\\*[.5ex]
\>$p\: .\: q\vee r$ \>$''$ \>$(p\: .\: q)\vee r$ \>$''$ \>$p\: .\: (q\vee r)$\\*[.5ex]
exactly as for arithmetical sum and product\\[.5ex]
\>$p\vee q\supset r$ \>$''$ \>$(p\vee q)\supset r$ \>$''$ \>$p\vee (q\supset r)$\\[.5ex]
\>$\sim p\supset q$ \>$''$ \>$(\sim p)\supset q$ \>$''$ \>$\sim(p\supset q)$\\[.5ex]
\>$\sim p\: .\: q$ \>$''$ \>$(\sim p)\, .\: q$ \>$''$ \>$\sim(p\: .\: q)$\\[.5ex]
\>$\sim p\equiv q$ \>$''$ \>$(\sim p)\equiv q$ \>$''$ \>$\sim(p\equiv q)$
\end{tabbing}
In all these cases the expression written without brackets has the meaning of the proposition in the second column. If we have the formula of the third column in mind we have to write the brackets. Another convention used in arithmetic for saving brackets is this that instead of $(a+b)+c$ we can write $a+b+c$. We make the same conventions for logical addition and multiplication, i.e.\ $p\vee q\vee r$ means $(p\vee q)\vee r$, $p\: .\: q\, .\: r$ means $(p\: .\: q)\, .\: r$.

The letters $p,q,r$ which denote arbitrary propositions are called propositional variables,\index{propositional variables}\index{variables, propositional} and any expression composed of propositional variables and the operations $\sim$, $\vee$, $\, .\:$, $\supset$, $\equiv$\ is called meaningful expression\index{expression (formula), of propositional logic}\index{formula, of propositional logic} or formula of the calculus of propositions, where also the letters $p,q$ themselves are considered as the simplest kind of expressions.

After those merely symbolic conventions the next thing we have to do is to examine in more detail the meaning of the operations of the calculus of propositions. Take e.g.\ the disjunction $\vee$. If $\mathbf{\llbracket 16. \rrbracket}$ any two proposi\-tions $p,q$ are given $p\vee q$ will again be a proposition. But now (and this is the decisive point) this operation of ``or'' is such that the truth or falsehood of the composite proposition $p\vee q$ depends in a definite way on the truth or falsehood of the constituents $p,q$. This dependence can be expressed most clearly in the form of a table\index{truth table} as follows: Let us form three columns, one headed by $p$, one by by $q$, one by $p\vee q$, and let us write T for true and F for false. Then for the propositions $p,q$ we have the following four possibilities
\begin{center}
\begin{tabular}{ c|c|c c|c }
 $p$ & $q$ & $p\vee q$ & $p\circ q$ & $p\: .\: q$\\[.5ex]
 T & T & T & F & T \\
 T & F & T & T & F\\
 F & T & T & T & F\\
 F & F & F & F & F
 \end{tabular}
\end{center}
Now for each of these four cases we can easily determine $\mathbf{\llbracket 17. \rrbracket}$ whether $p\vee q$ will be true or false; namely, since $p\vee q$ means that one or both of the propositions $p$, $q$ are true it will be true in the first, second and third case, and false only in the fourth case. We can consider this table (called the truth table for $\vee$) as the most precise definition of what $\vee$ means.

It is usual to call truth and falsehood the truth values\index{truth values} and to say of a true proposition that it has the truth value ``Truth'',\index{Truth} and of a false proposition that it has the truth value ``Falsehood''.\index{Falsehood} T\index{T, Truth} and F\index{F, Falsehood} then denote the truth values and the truth table for $\vee$ shows how the truth value of the composite expression ${p\vee q}$ depends on the truth values of the constituents. The exclusive ``or'' would have another truth $\mathbf{\llbracket 18. \rrbracket}$ table; namely if I denote it by $\circ$ for the moment, we have $p\circ q$ is false in the case when both $p$ and $q$ are true and in the case when both $p$ and $q$ are false, and it is true in the other cases, where one of the two propositions $p,q$ is true and the other one is false. The operation $\sim$ has the following truth table
\begin{center}
\begin{tabular}{ c|c}
 $p$ & $\sim p$ \\[.5ex]
 T & F \\
 F & T
 \end{tabular}
\end{center}
Here we have only two possibilities: $p$ is true and $p$ is false, and if $p$ is true not-$p$ is false and if $p$ is false not-$p$ is true. The truth table for ``and'' can also easily be be determined: $p\: .\: q$ is true only in the case where both $p$ and $q$ are true and false in all the other three cases.

A little more $\mathbf{\llbracket 19. \rrbracket}$ difficult is the question of the truth table for $\supset$. $p\supset q$ was defined to mean: If $p$ is true then $q$ is also true. So in order to determine the truth table let us assume that for two given propositions\ $p,q$ $p\supset q$ holds, i.e.\ let us assume \emph{we know} ``If $p$ then $q$'' \emph{but nothing else}, and let us ask what can we conclude about the truth values of $p$ and $q$ from this assumption.

\vspace{1.5ex}

\hspace{2em}\begin{tabular}{c c|c|c|c|c }
Assumption & $p\supset q$ & $p$ & $q$ & $\sim p$ & $\sim p\vee q$\\[.5ex]
& T & F & T & T & T\\
& T & F & F & T & T\\
& T & T & T & F & T\\
& F & T & F & F & F
\end{tabular}

\vspace{1.5ex}

\noindent First it may certainly happen that $p$ is false, because the assumption ``If $p$ then $q$'' says nothing about the truth or falsehood of $p$, and in this case when $p$ is false $q$ may be true as well as false, because the assumption says nothing about what happens to $q$ if $p$ is false, but only if $p$ is true. $\mathbf{\llbracket 20. \rrbracket}$ So we have both these possibilities: $p\;$F$\;\;q\;$T,
$p\;$F$\;\;q\;$F. Next we have the possibility that $p$ is true, but in this case $q$ must also be true owing to the assumption so that the possibility $p$ true $q$ false is excluded and it is the only of the four possibilities that is excluded by the assumption $p\supset q$. It follows that either one of those three possibilities, which I frame in
\begin{center}
\begin{tabular}{ |c|c| }
$p$ & $q$\\[.5ex]
\hline
F & T\\
F & F\\
T & T\\
\hline
\end{tabular}
\end{center}
occurs. But we have also vice versa: If one of these three possibilities for the truth value of $p$ and $q$ is realized then $p\supset q$ holds. For let us assume we know that one of the three marked $\mathbf{\llbracket 21. \rrbracket}$ cases occurs; then we know also ``If $p$ is true $q$ is true'', because if $p$ is true only the third of the three marked cases can be realized and in this case $q$ is true. So we see that the statement ``If $p$ then $q$'' is exactly equivalent with the statement that one of the three marked cases for the truth values of $p$ and $q$ is realized, i.e.\ $p\supset q$ will be true in each of the three marked cases and false in the last case. And this gives the desired truth table for implication. However there are two important remarks about it, namely:

1. Exactly the same truth table can also be $\mathbf{\llbracket 22. \rrbracket}$ obtained by a combination of operations introduced previously, namely $\sim p\vee q$, i.e.\ either $p$ is false or $q$ is true has the same truth table. For $\sim p$ is true whenever $p$ is false, i.e.\ in the first two cases and $\sim p\vee q$ is then true if either $\sim p$ or $q$ is true, and as you see that happens in exactly the cases where $p\supset q$ is true. So we see $p\supset q$ and $\sim p\vee q$ are equivalent, i.e.\ whenever $p\supset q$ holds then also $\sim p\vee q$ holds and vice versa. This makes possible to define $p\supset q$ by $\sim p\vee q$ and this is the usual way of introducing the implication in mathematical logic.

2. The second remark about the truth table for implication is this. We must $\mathbf{\llbracket 23. \rrbracket}$ not forget that $p\supset q$ was understood to mean simply ``If $p$ then $q$'' and nothing else, and only this made the construction of the truth table possible. There are other interpretations of the term ``implication'' for which our truth table would be completely inadequate. E.g.\ $p\supset q$ could be given the meaning: $q$ is a logical consequence of $p$, i.e.\ $q$ can be derived from $p$ by means of a chain of syllogisms. In this sense e.g.\ the proposition ``Jupiter is a planet'' would imply the proposition ``Jupiter is not a fixed star'' because no planet can be a fixed star by definition, i.e.\ $\mathbf{\llbracket 24. \rrbracket}$ by merely logical reasons.

This kind and also some other similar kinds of implication are called strict implication\index{strict implication} and denoted by this symbol $\prec$ and the implication defined by the truth table is called material implication\index{material implication} if it is to be distinguished from $\prec$. Now it is easy to see that our truth table would be false for strict implication and even more, namely that there exists no truth table at all for strict implication. In order to prove this consider the first line of our truth table, where $p$ and $q$ are both true and let us ask what will the truth value of $p\prec q$ be in this case. $\mathbf{\llbracket 25. \rrbracket}$ It turns out that this truth value is not uniquely determined. For take e.g.\ for $p$ the proposition ``Jupiter is a planet'' and for $q$ ``Jupiter is not a fixed star'', then $p,q$ are both true and $p\prec q$ is also true. On the other hand if you take for $p$ again ``Jupiter is a planet'' and for $q$ ``France is a republic'' then again both $p$ and $q$ are true, but $p\prec q$ is false because ``France is a republic'' is not a logical consequence of ``Jupiter is a planet''. So we see the truth value of $p\prec q$ is not uniquely determined by the truth values of $p$ and $q$ and therefore no truth$\,$ table exists. $\mathbf{\llbracket 26. \rrbracket}$ Such functions of propositions for which no truth table exists are called intensional\index{intensional functions (connectives)} as opposed to extensional\index{extensional functions (connectives)} ones for which a truth table does exist. The extensional functions are also called truth functions,\index{truth functions} because they depend only on the truth or falsehood of the propositions involved.

So we see logical consequence\index{logical consequence} is an intensional relation between propositions and the material implication introduced by our truth table cannot mean logical consequence. Its meaning is best given by the word ``if'' of ordinary language which has a much wider signification than just logical consequence; e.g.\ if someone says: ``If I don't come I $\mathbf{\llbracket 27. \rrbracket}$ shall call you'' that does not indicate that this telephoning is a logical consequence of his not coming, but it means simply he will either come or telephone, which is exactly the meaning expressed by the truth table. Hence material implication introduced by the truth tables corresponds as closely to ``if then''\index{if\ldots\ then} as a precise notion can correspond to a not precise notion of ordinary language.

If we are now confronted with the question which one of the two kinds of implication we shall use in developing the theory of inference we have to consider two things: 1.\ material implication is the much simpler and clearer notion and 2.\ it is quite sufficient for developing the theory of inference because in order to conclude $q$ from $p$ it is quite sufficient $\mathbf{\llbracket 28. \rrbracket}$ to know $p$ implies materially $q$ and not necessary to know that $p$ implies strictly $q$. For if we know $p\supset q$ we know that either $p$ is false or $q$ is true. Hence if we know in addition that $p$ is true the first of the two possibilities that $p$ is false is not realized. Hence the second must be realized, namely $q$ is true. For these two reasons that material implication is simpler and sufficient I shall use only material implication at least in the first introductory part of my lectures, and shall use the terms ``implies''\index{implies} and ``follows''\index{follows} only in the sense of material implication. I do not want to say by this that a theory of strict implication may not be interesting and important for certain purposes. In fact I hope it will be discussed in the second half of this seminary. But this theory belongs to an entirely different part of logic than the one I am dealing with now, $\mathbf{\llbracket 29. \rrbracket}$ namely to the logic of modalities.\index{logic of modalities}

I come now to some apparently paradoxical consequences of our definition of material implication whose paradoxicality\index{paradoxes, of material implication} however disappears if we remember that it does not mean logical consequence. The first of these consequences is that a true proposition is implied by any proposition whatsoever. We see this at once from the truth table which shows that $p\supset q$ is always true if $q$ is true whatever $p$ may be. You see there are only two cases where $q$ is true and in both of them $p\supset q$ is true. But secondly we see also that $p\supset q$ is always true if $p$ is false whatever $q$ may be. So that means a false proposition implies any proposition whatsoever, which is the second of the paradoxical consequences. These properties of implication $\mathbf{\llbracket 30. \rrbracket}$ can also be expressed by saying: ``An implication with true second member is always true whatever the first member may be and an implication with false first member is always true whatever the second member may be''; we can express that also by formulas like this ${q\supset(p\supset q)}$, ${\sim p\supset(p\supset q)}$. Both of these formulas are also immediate consequences of the fact that $p\supset q$ is equivalent with ${\sim p\vee q}$ because what ${\sim p\vee q}$ says is exactly that either $p$ is false or $q$ is true; so ${\sim p\vee q}$ will always be true if $p$ is false and will be also true if $q$ is true whatever the other proposition may be. If we apply $\mathbf{\llbracket 31. \rrbracket}$ these formulas to special cases we get strange consequences; e.g.\ ``Jupiter is a fixed star'' implies ``France is a republic'', but it also implies ``France is not a republic'' because a false proposition implies any proposition whatsoever. Similarly ``France is a republic'' is implied by ``Jupiter is a planet'' but also by ``Jupiter is a fixed star''. But as I mentioned before these consequences are paradoxical only for strict implication. They are in pretty good agreement with the meaning which the word ``if'' has in ordinary language. Because the first formula then says if $q$ is true $q$ is also true if $p$ is true which is not paradoxical but trivial and the second says if $p$ is false then if $p$ is true anything $\mathbf{\llbracket 32. \rrbracket}$ is true. That this is in good agreement with the meaning which the word ``if'' has can be seen from many colloquialisms; e.g.\ if something is obviously false one says sometimes ``If this is true I am a Chinaman'', which is another way of saying ``If this is true anything is true''. Another of these so called paradoxical consequences is e.g.\ that for any two arbitrary propositions one must imply the other, i.e.\ for any $p,q$ $(p\supset q)\vee(q\supset p)$; in fact $q$ must be either true or false---if it is true the first member of the disjunction is true because it is an implication with true second member, if it is false the second member of the disjunction is $\mathbf{\llbracket 33. \rrbracket}$ true. So this disjunction is always true.

\subsection{Tautologies}\label{1.6}
\pagestyle{myheadings}\markboth{EDITED TEXT}{NOTEBOOK 0 \;\;---\;\; 1.1.6\;\; Tautologies}

Those three formulas, as well as the formula of disjunctive inference we had before,\footnote{see pp.\ \textbf{30}., \textbf{32}.\ and \textbf{5}.\ of the present Notebook~0} are examples of so called universally true formulas,\index{universally true formula, of propositional logic (tautology)} i.e.\ formulas which are true whatever the pro\-positions $p,q,r$ occurring in them may be. Such formulas are also called logically true\index{logically true formulas, of propositional logic (tautology)} or tautological,\index{tautological formula, of propositional logic}\index{tautology, of propositional logic} and it is exactly the chief aim of the calculus of propositions to investigate these tautological formulas.

I shall begin with discussing a few more examples before going over to more general considerations. I mention at first some of the traditional hypothetical and $\mathbf{\llbracket 34. \rrbracket}$ disjunctive inferences which in our notation read as follows:
\begin{itemize}
\item[1.] $(p\supset q)\, .\: p\supset q$ \quad ponendo ponens\index{ponendo ponens}\index{modus ponens} \quad (Assertion)
\vspace{-1ex}
\item[2.] $(p\supset q)\, .\: \sim q\supset\;\sim p$ \quad tollendo tollens\index{tollendo tollens}
\vspace{-1ex}
\item[3.] $(p\vee q)\, .\: \sim q\supset p$ \quad tollendo ponens\index{tollendo ponens} \quad as we had before\\
(the modus ponendo tollens\index{ponendo tollens} holds only for the exclusive $\vee$)
\vspace{-1ex}
\item[4.] An inference which is also treated in many of the textbooks under the heading of ``dilemma''\index{dilemma} is this\\[1ex]
  $(p\supset r)\, .\:(q\supset r)\supset(p\vee q\supset r)$\\[1ex]
   If both $p\supset r$ and $q\supset r$ then from $p\vee q$ follows $r$. It is usually written as an inference with three premises, $\mathbf{\llbracket 35. \rrbracket}$ namely from the three premises $(p\supset r)\, .\:(q\supset r)\, .\:(p\vee q)$ one can conclude $r$. \end{itemize}
This is nothing else but the principle of proof by cases,\index{proof by cases, dilemma} namely the premises say: one of the two cases $p,q$ must occur and from both of them follows $r$. That this formula with three premises means the same thing as the formula under consideration is clear because this earlier formula says: ``If the first two premises are true then if the third is true $r$ is true'', which means exactly the same thing as ``If all the three premises are true $r$ is true. The possibility of going over from one of these two formulas to the other is due to another important logical principle which is called importation and reads like this
\begin{tabbing}
\hspace{1.7em}$[p\supset(q\supset r)]\supset(p\: .\: q\supset r)$ \quad importation\index{importation}
\end{tabbing}
and its inverse which is called exportation and reads like this
\begin{tabbing}
\hspace{1.7em}$(p\: .\: q\supset r)\supset[p\supset(q\supset r)]$ \quad exportation.\index{exportation}
\end{tabbing}
So owing to these two implications we have also an equivalence between the left and right-hand side.
Next we have the three laws of identity,\index{law of identity} excluded middle\index{law of excluded middle}\index{excluded middle} and contradiction\index{law of contradiction} which read as follows in our notation
\begin{tabbing}
\hspace{1.7em}1. $p\supset p$\quad\quad 2. $p\:\vee\sim p$ \quad\quad 3. $\sim\!(p\: .\sim p)$
\end{tabbing}
We can add another similar law, the law of double negation\index{law of double negation} which says ${\sim\!(\sim p)\equiv p}$.
Next we have the very important formulas of transposition:\index{transposition (contraposition)}\index{contraposition}
\begin{tabbing}
\hspace{1.7em}$(p\supset q)\supset(\sim q\supset\: \sim p)$
\end{tabbing}
Other forms of this formula of transposition would be
\begin{tabbing}
\hspace{1.7em}\=$(p\supset\: \sim q)\supset(q\supset\: \sim p)$ \\
\>$(\sim p\supset q)\supset(\sim q\supset p)$ \quad proved in the same way. \end{tabbing}
In all those formulas of transposition we can write equivalence instead of the main implication,\footnote{Instead of ``the main implication'' in the manuscript one finds ``identity''.} i.e.\ $\mathbf{\llbracket 36. \rrbracket}$ we have also $(p\supset q)\equiv(\sim q\supset\: \sim p)$. Another form of transposition, namely with two premises is this $(p\: .\: q\supset r)\supset ({p\: .\sim r}\supset\: \sim q)$ because under the assumption $p\: .\: q\supset r$ if we know $p\: .\sim r$, then $q$ cannot be true because $r$ would be true in this case.

Next we have different so called reductio ad absurdum,\index{reductio ad absurdum} e.g.\ \begin{tabbing}
\hspace{1.7em}$(p\supset q)\, .\: (p\supset\:\sim q)\supset\:\sim p$
\end{tabbing}
A particularly interesting form of reductio ad absurdum is the one which Professor Menger\index{Menger, Karl} mentioned in his introductory talk and which reads as follows
\begin{tabbing}
\hspace{1.7em}$(\sim p\supset p)\supset p$
\end{tabbing}

Other examples of logically true formulas are the commutative\index{commutative law} and associative law\index{associative law} for disjunction and conjunction
\begin{itemize}
\item[1.] $p\vee q\equiv q\vee p$
\vspace{-1ex}
\item[2.] $(p\vee q)\vee r\equiv p\vee(q\vee r)$
\vspace{-1ex}
\item[3.] similar formulas hold for addition\\[.5ex]
 $p\: .\: q\equiv q\, .\: p$,\quad $(p\: .\: q)\, .\: r\equiv p\: .\:(q\, .\: r)$
\end{itemize}

$\mathbf{\llbracket 37. \rrbracket}$ Next we have some formulas connecting $\vee$ and $.$ namely at first the famous so called De Morgan formulas:\index{De Morgan formulas (laws)}
\begin{tabbing}
\hspace{1.7em}\=$\sim(p\: .\: q)\equiv \;\sim p\;\vee\sim q$\\[.5ex]
\>$\sim(p\vee q)\equiv \;\sim p\: .\sim q$
\end{tabbing}
The left-hand side of the first means not both $p,q$ are true, the right-hand side at least one is false. The left-hand side of the second means not at least one is true, the right-hand side both are false.

These formulas give a means to distribute so to speak the negation of a product on the two factors and also the negation of a sum on the two terms, where however sum has to be changed into product and product into sum in this distribution process. Another tautology connecting sum and product is $\mathbf{\llbracket 38. \rrbracket}$ the distributive law\index{distributive laws} which reads exactly analogously as in arithmetic
\begin{tabbing}
\hspace{1.9em}\=1. \= $p\: .\: (q\vee r)\equiv p\: .\: q\vee p\: .\: r$
\end{tabbing}
because let us assume left is true then we have $p$ and two cases $q$, $r$; in the first case $p\: .\:q$, in the second $p\: .\:r$ is true, hence in any case right is true
\begin{tabbing}
and \=2.\hspace{.3em}\= $p\vee q\, .\: r\equiv (p\vee q)\, .\: (p\vee r)$\\[.5ex]
\>3.\> $(p\supset q)\, .\:(q\supset r)\supset(p\supset r)$\quad Syllogism, Transitivity of $\supset$\\[.5ex]
\>4.\> $(p\supset q)\supset[(q\supset r)\supset(p\supset r)]$\\[.5ex]
\end{tabbing}

\vspace{-1.5ex}

\hspace{1.5em}$(p\: .\: q\supset r)\supset[p\supset(q\supset r)]$\quad Export

\hspace{4.9em}inverse \hspace{5.2em} Import\footnote{This line and the preceding one are crossed out in the manuscript.}

\vspace{-1.5ex}

\begin{tabbing}
\hspace{1.9em}\=5.\hspace{.3em}\= $(p\supset q)\, .\: (r\supset s)\supset(p\: .\: r\supset q\, .\: s)$ Leibnitz\index{Leibnitz theorema praeclarum} theorema praeclarum\index{theorema praeclarum}\\[.5ex]
\>\>$(p\supset q)\supset(p\: .\: r\supset q\, .\: r)$ \hspace{1.6em} factor\\[.5ex]
\>6.\> $(p\supset q)\: .\: (r\supset s)\supset(p\vee r\supset q\vee s)$\\[.5ex]
\>\>$(p\supset q)\supset(p\vee r\supset q\vee r)$ \quad Sum\\[.5ex]
\>7.\> $p\supset p\vee q$ \quad \=7$'$. $p\: .\: q\supset p$\\[.5ex]
\>8.\> $p\vee p \supset p$ \quad \>8$'$. $p\;$\=$\supset p \: .\: p$\\
\>\>\hspace{2.4em}$\equiv$\>\>$\equiv$\\[.5ex]
\>9.\> $p\supset(q\supset p\: .\:q)$
\end{tabbing}

\subsection{Decidability for propositional logic}\label{1.7}
\pagestyle{myheadings}\markboth{EDITED TEXT}{NOTEBOOK I \;\;---\;\; 1.1.7\;\; Decidability for propositional logic}

$\mathbf{\llbracket Notebook\; I \rrbracket}$ $\mathbf{\llbracket 38.1 \: II \rrbracket}$\footnote{Notebook~0 ends with p.\ \textbf{38}., and hence, judging by how it is numbered, the present page should be a continuation of Notebook~0. The content of this page does not make obvious this supposition, but does not exclude it.}  Last time and also today in the classes we set up the truth tables for some of the functions which occur in the calculus of propositi\-ons. Their purpose is to give a precise definition of the functions concerned because they state exactly the conditions under which the proposit\-ion to be defined, e.g.\ $p\vee q$, is true and under which conditions it is not true. In ordinary language we have also the notions and, or, if etc. which have very approximately the same meaning, but for setting up a mathematical theory it is necessary that the notions involved have a higher degree of preciseness than the notions of ordinary language. It is exactly this what is accomplished by the truth tables.

$\mathbf{\llbracket 40. \: II \rrbracket}$\footnote{In the scanned manuscript there is no page numbered with \textbf{39}.\ in the present Notebook~I.} Take e.g.\ the formula $p\, .\: (p\supset q)\supset q$, the modus ponendo ponens. Here we have two propositional variables $p,q$ and therefore \emph{four possibilities} for these truth values, namely
\begin{center}
\begin{tabular}{ c|c|c|c|c }
 $p$ & $q$ & $p\supset q$ & $p\, .\:(p\supset q)$ & $p\, .\:(p\supset q)\supset q$ \\[.5ex]
 T & T & T & T & T\\
 T & F & F & F & T\\
 F & T & T & F & T\\
 F & F & T & F & T
 \end{tabular}
\end{center}

\noindent $\mathbf{\llbracket 41. \: II \rrbracket}$ and what we have to do is simply to check that the truth value of the whole expression is true in each of these four cases, i.e.\ we have to ascertain that the truth table of the whole expression consists of T's only. That's very simple. let us write down all the parts of which this expression is built up. We have first $p\supset q$ is a part, then $p\, .\:(p\supset q)$ and finally the whole expression. So we see actually in all four cases the whole formula is true. Hence it is universally true. It is clear that this purely mechanical method of checking all possibilities will always give a decision whether a given formula is or is not a $\mathbf{\llbracket 42. \: II \rrbracket}$ tautology. Only if the number of variables $p,q$ occurring in the expression is large this method is very cumbersome, because the number of cases which we have to deal with is $2^n$ if the number of variables is $n$ and the number of cases is the same as the number of lines in the truth table. Here we had 2 variables $p,q$ and therefore $2^2=4$ cases. With 3 variables we would have $2^3=8$ cases and in general if the number of variables is increased by one the number of cases to be considered is doubled, because each of the previous cases is split into two new cases according as the truth value of the new variable is truth or falsehood. Therefore we have $\mathbf{\llbracket 43. \: II \rrbracket}$ $2^n$ cases for $n$ variables. In the applications however usually the number of cases actually to be considered is much smaller because mostly several cases can be combined into one, e.g.\ in our example case 1 and 2 can be treated together because if $q$ is true the whole expression is certainly true whatever $p$ may be because it is then an implication with true second member.

So we see that for the calculus of propositions we have a very simple procedure to decide\index{decision procedure}\index{decidability} for any given formula whether or not it is logically true. This solves the first of the two general problems which I mentioned in the beginning for the calculus of propositions, namely the problem to give a complete theory of logically true formulas. We have even more, namely a procedure to decide of any formula whether or not it is logically true. That this problem $\mathbf{\llbracket 44. \: II \rrbracket}$ could be solved in such a simple way is chiefly due to the fact that we introduced only extensional operations (only truth functions of propositions). If we had introduced strict implication the question would have been much more complicated. It is only very recently that one has discovered general procedures for deciding whether a formula involving strict implication is logically true under certain assumptions about strict implication.

Now after having solved this so called decision problem\index{decision problem} for the calculus of propositions I can go over to the second problem I have announced in the beginning.

\subsection{Functional completeness}\index{functional completeness}\label{1.8}
\pagestyle{myheadings}\markboth{EDITED TEXT}{NOTEBOOK II \;\;---\;\; 1.1.8\;\; Functional completeness}

$\mathbf{\llbracket Notebook\; II \rrbracket}$ $\mathbf{\llbracket 33. \rrbracket}$\footnote{In the source version, as in the manuscript, one finds in the present Notebook II first pages numbered \textbf{61}.-\textbf{76}., which is followed by pages numbered \textbf{33}.-\textbf{55}.\textbf{2}. In this edited version, the order of these two blocks of pages is permuted, which puts them in the right arithmetical order, and in between pp.\ \textbf{56}.-\textbf{60}.\ of Notebook~I fit well.} After having solved last time the first of the two problems I announced in the beginning, namely the problem of deciding of a given expression whether or not it is a tautology, I come now to the second, namely to reduce the infinite number of tautologies to a finite number of axioms from which they can be derived. So this problem consists in setting up what is called a deductive system for the calculus of propositions. Now if you think of other examples of deductive systems as e.g.\ geometry you will see that their aim is not truly to derive the theorems of the science concerned from a minimal number of axioms, but also to define the notions of the discipline concerned in terms of a minimal number of undefined or $\mathbf{\llbracket 34. \rrbracket}$ primitive notions. So we shall do the same thing for the calculus of propositions.

We know already that some of the notions introduced $\sim$, $\vee$, $ .\:$, $\supset$, $\equiv$, $\mid$ can be defined in terms of others, namely e.g.\ $p\supset q\equiv\: \sim p\vee q$, $p\equiv q\equiv p\supset q\: .\: q\supset p$, but now we want to choose some of them in terms of which all others can be defined. And I claim that e.g.\ $\sim$ and $\vee$ are sufficient for this purpose because
\begin{tabbing}
\hspace{1.7em}\=1.\quad \=$p\: .\: q\;\;$ \=$\equiv\;\sim(\sim p\;\vee \sim q)$\\[.5ex]
\>2.\> $p\supset q$ \>$\equiv\;\sim p\vee q$\\[.5ex]
\>3.\> $p\equiv q$ \>$\equiv (p\supset q)\: .\: (q\supset p)$\\[.5ex]
\>4.\> $p\mid q$ \>$\equiv\;\sim p\;\vee\sim q$
\end{tabbing}
So it is possible to take $\sim$ and $\vee$ as $\mathbf{\llbracket 35. \rrbracket}$ primitive terms for our deductive system and we shall actually make this choice. But it is important to remark that this choice is fairly arbitrary. There would be other possibilities, e.g.\ to take $\sim$, $ .\:$ because $\vee$ can be expressed in terms of $\sim$ and $\: .\:$ by $p\vee q\equiv\;\sim(\sim p\: .\:\sim q)$ and by $\vee$ and $\sim$ the others can be expressed as we have just seen. This fact that the choice of primitive terms is arbitrary to a certain extent is not surprising. The same situation prevails in any theory, e.g.\ in geometry we can take either the notion of movement of the space or the notion of congruence between figures as primitive because it is possible $\mathbf{\llbracket 36. \rrbracket}$ to define congruence of figures in terms of movement of space and vice versa. The same situation we have here. We can define $\vee$ in terms of ``and'' and ``not'' but also vice versa ``and'' in terms of ``or'' and ``not''. And there are still further possibilities for the primitive terms, e.g.\ it would be possible to take $\sim$ and $\supset$ as the only primitive terms because $\vee$ can be defined by
\begin{tabbing}
\hspace{1.7em}\=$p\vee q\equiv\;\sim p\supset q$ \quad since\\[.5ex]
\>$\sim p\supset q\equiv\;\sim\sim p\vee q\equiv p\vee q$\quad by the law of double negation
\end{tabbing}
In the three cases discussed so far we had always two primitive notions in terms of $\mathbf{\llbracket 37. \rrbracket}$ which the others could be defined. It is an interesting question whether there might not be a single operation in terms of which all the others can be defined. This is actually the case as was first discovered by the logician Sheffer.\index{Sheffer, Henry Maurice} Namely the $\mid$ function\index{Sheffer stroke} suffices to define all the others because $\sim p\equiv p\mid p$ means at least one of the propositions $p, p$ is false, but since they are both $p$ that means $p$ is false, i.e.\ $\sim p$, so $\sim$ can be defined in terms of $\mid$ and now the ``and'' can be defined in terms of $\sim$ and $\mid$ since $p\: .\: q\equiv\;\sim(p\mid q)$ for $p\mid q$ means at least $\mathbf{\llbracket 38. \rrbracket}$ one of the two propositions is false; hence the negation means both are true. But in terms of $\sim$ and the $\: .\:$ others can be defined as we saw before. It is easy to see that we have now exhausted all possibilities of choosing the primitive terms from the six operations written down here. In particular we can prove e.g.: $\sim, \equiv$ are not sufficient to define the others in terms of them. We can e.g.\ show that $p\vee q$ cannot be defined in terms of them.

Now what could it mean that $p\: .\: q$ or $p\vee q$ can be defined in terms of $\sim , \equiv$? It would mean that we can find an expression $f(p,q)$ in two variables containing only the symbols $\sim , \equiv$ besides $p,q$ and such that $p {\vee \atop .} q\equiv f(p,q)$, i.e.\ such that this expression would have the same truth table as $p {\vee \atop .} q$. But we shall prove now that such an expression does not exist.

$\Big\lceil$Let's write down the truth functions in two variables $p,q$ which we certainly can define in terms of $\sim , \equiv$; we get the following eight: 1. ${p\equiv p}$, 2. $\sim (p\equiv p)$, 3. $p$, 4. $q$, 5. $\sim p$, 6. $\sim q$, $\mathbf{\llbracket new\, page \rrbracket}$ 7. $p\equiv q$, 8. $\sim(p\equiv q)$, and now it can be shown that no others can be defined except those eight because we can show the following two things: 1. If we take one of those eight functions and negate it we get again one of those eight functions, 2. If we take any two of those eight functions and form a new one by connecting them by an equivalence symbol we get again one of the eight. I.e.\ by application of the operation of negation and of the operation of equivalence we never get outside of the set of eight functions written down. So let's see at first that by negating them $\mathbf{\llbracket new\, page \rrbracket}$ we don't get anything new. Now if we negate the first\ldots\ Now let's connect any two of them by $\equiv$. If we connect the first with any formula $P$ we get $P$ again, i.e.\ $(\top\equiv P)\equiv P$ because\ldots\ and if connect a contradiction $C$ with any formula $P$ by $\equiv$ we get the negation of $P$, i.e.\ $(C\equiv P)\equiv\;\sim P$ because\ldots\ So by combining the first two formulas with any other we get certainly nothing new. For the other cases it is very helpful that $(p\equiv\;\sim q)\equiv\;\sim(p\equiv q)$; this makes possible to factor out the negation so to speak. Now in order to apply that to the other formulas we divide them in two groups\ldots$\Big\rfloor$

$\mathbf{\llbracket 39. \rrbracket}$ For this purpose we divide the 16 truth functions of two variables which we wrote down last time into two classes according as the number of letters T occurring in their truth table is even or odd, or to be more exact according as the number of T's occurring in the last column. So e.g.\ $p\: .\: q$ is odd, $p\equiv q$ is even and an arbitrary expression in two variables will be called even if its truth function is even. And now what we can show is this: Any expression in two variables containing only $\sim$ and $\equiv$ is even (i.e.\ its truth table contains an even number of T's, i.e.\ either 0 or 2 or 4 T's).
In order to show that we prove the following three lemmas.
\begin{itemize}
\item[1.] The letter expressions, namely the letters $p,q$ are even.
\item[2.] If an expression $f(p,q)$ is even then also the expression $\sim f(p,q)$ is even.
\item[3.] If two expressions $f(p,q)$, $g(p,q)$ are even then also the expression $f(p,q)\equiv g(p,q)$ obtained by connecting them with an equivalence sign is even.
\end{itemize}
$\mathbf{\llbracket 40. \rrbracket}$ So propositions 2, 3 have the consequence:
\begin{itemize}
\item[]
By applying the operations $\sim$ and $\equiv$ to even expressions as many times as we wish we always get again an even expression.
\end{itemize}
But any expression containing only $\sim$ and $\equiv$ is obtained from the single letters $p,q$ by an iterated application of the operations $\sim$ and $\equiv$; hence since $p,q$ are even the expression thus obtained will also be even. So our theorem that every expression containing only $\sim$ and $\equiv$ is even will be proved when we shall have proved the three lemmas.

1. is clear because of the truth table for $p$ (and for $q$ the same thing). 2. also is clear because $\sim f(p,q)$ has T's when $f(p,q)$ had F's, i.e.\ the number of T's in the new expression is the same as the number of F's in the $\mathbf{\llbracket 41. \rrbracket}$ old one. But the number of F's in the old one is even because the number of T's is even and the number of F's is equal to the number of T's.

Now to the third. Call the number of T's of the first $t_1$, the number of T's of the second $t_2$ and call the number of cases in the truth table where both $f$ and $g$ have the truth value T $r$. We have that $t_1$ is even and $t_2$ is even, but we do not know anything about $r$; it may be odd or even. We shall try to find out in how many cases $f(p,q)\equiv g(p,q)$, i.e.\ $f\equiv g$, will be true and to show that this number of cases will be even. I prefer to find out in how many cases it will be false. If we know that this number is even we know also that the number of cases in which it is true will be even. Now this whole expression is false if $g$ and $f$ have different truth values, i.e.\ if $\mathbf{\llbracket 42. \rrbracket}$ either we have $g$ false and $f$ true or we have $g$ true and $f$ false. The cases where $f$ is true and $g$ false make $t_1-r$ cases because from $t_1$ cases where $f$ is true we should subtract cases when $g$ is also true, and because $r$ was the number of cases in which both are true. Hence in $t_1-r$ cases $f$ is T and $g$ is F, and similarly in $t_2-r$ cases $g$ is T and $f$ is F; hence altogether in $t_1-r + t_2-r$ cases $f$ and $g$ have different truth values, i.e.\ in $t_1+t_2-2r$ cases $f(p,q)\equiv g(p,q)$ is false, and this is an even number because $t_1$, $t_2$ and $2r$ are even, and if you add an even number to an even number, after subtracting an even number from the sum you get again an even number. Hence the number of cases in which the whole expression is false is an even number and such is also the number of cases in which it is true, i.e.\ $f(p,q)\equiv g(p,q)$ is an even expression. q.e.d.

\vspace{1ex}

So this shows that only even truth functions $\mathbf{\llbracket 43. \rrbracket}$ can be expressed in terms of $\sim$ and $\equiv$. Hence e.g.\ $\vee$ and $.$ cannot be expressed because three T's occur in their truth tables. It is easy to see that of the 16 truth functions exactly half the number is even and also that all even truth functions really can be expressed in terms of $\sim$ and $\equiv$ alone. Expressions for these eight truth functions in terms of $\sim$ and $\equiv$ are given in the notes that were distributed.\footnote{see p.\ \textbf{38}.\ in the present Notebook II above} The general theorem on even functions I proved then has the consequence that these eight truth functions must reproduce themselves by negating them or by connecting any two of them by $\sim$; i.e.\ if you negate one of those expressions the resulting expression will be equivalent to one of the eight and if you form a new expression by connecting any two of them the resulting expression will again be equivalent to one of the eight. I recommend $\mathbf{\llbracket 44. \rrbracket}$ as an exercise to show that in detail.

It is an easy corollary of this result about the undefinability of $.$ and $\vee$ in terms of $\equiv$ that also $\sim$ and the exclusive \emph{or} are not sufficient as primitive terms because as we saw last time the exclusive or can be expressed in terms of $\sim$ and $\equiv$, namely by $\sim(p\equiv q)$; hence if e.g.\ $\vee$ could be defined in terms of $\sim$ and $\circ$ (exclusive or) it could also be defined in terms of $\sim$ and $\equiv$ because the $\circ$ can be expressed in terms of $\sim$ and $\equiv$. The reason for that is of course that $\circ$ is also an even function and therefor only even functions can be defined in terms of it. So we see that whereas $\sim$ and $\vee$ are sufficient as $\mathbf{\llbracket 45. \rrbracket}$ primitive terms $\sim$ \emph{and exclusive or} are not, which is one of the reasons why the not exclusive or is used in logic. Another of those negative results about the possibility of expressing some of the truth functions by others would be that $\sim$ cannot be defined in terms of $.\:, \vee, \supset$; even in terms of all three of them it is impossible to express $\sim$. I will give that as a problem to prove.

As I announced before we shall choose from the different possibilities of primitive terms for our deductive system the one in which $\sim$ and $\vee$ are taken as primitive and therefore it is of importance to make sure that not only the particular functions $\equiv$, $.\:$, $\supset$, $\mid$ for which $\mathbf{\llbracket 46. \rrbracket}$ we introduced special symbols but that any truth function whatsoever in any number of variables can be expressed by $\sim$ and $\vee$. For truth functions with two variables that follows from the considerations of last time since we have expressed all 16 truth functions by our logistic symbols and today we have seen that all of them can be expressed by\index{functional completeness} $\sim$ and $\vee$. Now I shall prove the same thing also for truth functions with three variables and you will see that the method of proof can be applied to functions of any number of variables. For the three variables $p.q,r$ we have eight $\mathbf{\llbracket 47. \rrbracket}$ possibilities for the distribution of truth values over them, namely
\begin{center}
\begin{tabular}{ l c c c|l l}
 & $p$ & $q$ & $r$ & $f(p,q,r)$\\[.5ex]
 \hline
1.& T & T & T & \hspace{.1em} $p\: .\: q\: .\: r$ & \hspace{1em}$P_1$\\
2.& T & T & F & \hspace{.1em} $p\: .\: q\: . \sim r$ & \hspace{1em}$P_2$\\
3. & T & F & T & \hspace{.1em} $p\: . \sim q\: .\: r$\\
4.& T & F & F\\
5.& F & T & T\\
6.& F & T & F\\
7.& F & F & T\\
8.& F & F & F & & \hspace{1em}$P_8$
\end{tabular}
\end{center}

Now to define a truth function in three variables means to stipulate a truth value T or F for $f(p,q,r)$ for each of these eight cases. Now to each of these eight cases we can associate a certain expression in the following way: to 1. we associate $p\: .\: q\: .\: r$, to 2. we associate $p\: .\: q\: . \sim r$, to 3. we associate $p\: . \sim q\: .\: r$,\ldots\ So each of these expressions will have a $\sim$ before those letters which have an F in the corresponding case. Denote the expressions associated with these eight lines by $P_1$,\ldots ,$P_8$. Then the expression $P_2$ e.g.\ will be true then and only $\mathbf{\llbracket 48. \rrbracket}$ then if the \emph{second} case is realized for the truth values of $p,q,r$ ($p\: .\: q\: . \sim r$ will be true then and only then if $p$ is T, $q$ is T and $r$ is false, which is exactly the case for the truth values $p,q,r$ represented in the second line. And generally $P_i$ will be true if the $i^{\,\rm th}$ case for the truth values of $p,q,r$ is realized. Now the truth function which we want to express by $\sim$ and $\vee$ will be true for certain of those eight cases and false for the others. Assume it is true for case number $i_1$,$i_2$,\ldots,$i_n$ and false for the others. Then form the disjunction $P_{i_1}\vee P_{i_2}\ldots\vee P_{i_n}$, i.e.\ the disjunction of those $P_i$ which correspond to the cases in which the given function is true. This disjunction is an expression in the variables $p,q,r$ containing only the operations $.$, $\sim$ and $\vee$, and I claim its truth table $\mathbf{\llbracket 49. \rrbracket}$ will coincide with the truth table of the given expression $f(p, q, r)$. For $f(p, q, r)$ had the symbol T in the $i_1$,$i_2$,\ldots,$i_n^{\;\rm th}$ line but in no others and I claim the same thing is true for the expression $P_{i_1}\vee\ldots\vee P_{i_n}$.

You see at last a disjunction of an arbitrary number of members will be true then and only then if at least one of its members is true and it will be false only if all of its members are false (I proved that in my last lecture for the case of three members and the same proof holds generally). Hence this disjunction will certainly be true in the $i_1$,\ldots , $i_n^{\;\rm th}$ case because $P_{i_1}$ e.g.\ is true in the $i_1^{\;\rm th}$ case as we saw before. Therefore the $\mathbf{\llbracket 50. \rrbracket}$ disjunction is also true for the $i_1^{\;\rm th}$ case because then one of its members is true. The same holds for $i_2$\ldots\ etc. So the truth table for the disjunction will certainly have the letter T in the $i_1$,\ldots,$i_n$ line. But it will have F's in all the other lines. Because $P_{i_1}$ was true only in the $i_1^{\;\rm th}$ case and false in all the others. Hence in a case different from the $i_1$,\ldots,$i_n^{\;\rm th}$ $P_{i_1}$,\ldots ,$P_{i_n}$ will all be false and hence the disjunction will be false, i.e.\ $P_{i_1}\vee\ldots\vee P_{i_n}$ will have the letter F in all lines other than the $i_1$,\ldots,$i_n^{\;\rm th}$, i,e.\ it has T in the $i_1$,\ldots,$i_n$ line and only in those. But the same thing was true for the truth table of the given $f(p, q , r)$ by assumption. So they coincide, i.e.\ $f(p,q,r)\equiv P_{i_1}\vee \ldots\vee P_{i_n}$.

$\mathbf{\llbracket 51. \rrbracket}$ So we have proved that an arbitrary truth function of three variables can be expressed by $\sim$, $\vee$ and $.$, but $.$ can be expressed by $\sim$ and $\vee$, hence every truth$\,$ function of three variables can be expressed by $\sim$ and $\vee$, and I think it is perfectly clear that exactly the same proof applies to truth functions of any number of variables.

\vspace{1ex}

Now after having seen that \emph{two primitive notions} $\sim, \vee$ really suffice to define any truth function we can begin to set up the deductive system.

I begin with writing three \emph{definitions in terms of our primitive notions}:\index{defined connectives}
\begin{tabbing}
\hspace{1.7em}\= $P\supset Q$ \=$=_{\rm{Df}}\;\:$\=$\sim P \vee Q$\\[.5ex]
\> $P\: .\: Q$ \>$=_{\rm{Df}}$\>$\sim(\sim P\:\vee \sim Q)$\\[.5ex]
\> $P\equiv Q$ \>$=_{\rm{Df}}$\>$P\supset Q\: .\: Q\supset P$
\end{tabbing}
$\mathbf{\llbracket 52. \rrbracket}$ I am writing capital letters\index{schematic letters} because these definitions are to apply also if $P$ and $Q$ are formulas, not only if they are single letters, i.e.\ e.g.\ $p\supset p\vee q$ means $\sim p\vee(p\vee q)$ and so on.

\subsection{Axiom system for propositional logic}\label{1.9}
\pagestyle{myheadings}\markboth{EDITED TEXT}{NOTEBOOK II \;\;---\;\; 1.1.9\;\; Axiom system for propositional logic}

The\footnote{This section is made of the following blocks of pages in the following order: pp.\ \textbf{52}.-\textbf{55}.\textbf{2} of Notebook II, pp.\ \textbf{56}.-\textbf{60}.\ of Notebook~I and pp.\ \textbf{61}.-\textbf{64}.\ of Notebook II}\index{axiom system for propositional logic} next thing to do in order to have a deductive system is to set up the axioms. Again in the axioms one has a freedom of choice as in the primitive terms, exactly as also in other deductive theories, e.g.\ in geometry, many different systems of axioms have been set up each of which is sufficient to derive the whole geometry. The system of axioms for the calculus of propositions which I use is essentially the one set up by first by Russell\index{Russell, Bertrand} and then also adopted by Hilbert.\index{Hilbert, David} It has the following four axioms:\index{axioms for propositional logic}\index{four axioms for propositional logic}
\begin{tabbing}
$\mathbf{\llbracket 53. \rrbracket}$\\*[.5ex]
\hspace{1.7em}\= (1)\quad \=$p\supset p\vee q$\index{axiom (1) for propositional logic}\index{(1), axiom for propositional logic}\\[.5ex]
\> (2)\>$p\vee p\supset p$\index{axiom (2) for propositional logic}\index{(2), axiom for propositional logic}\\[.5ex]
\> (3)\>$p\vee q\supset q\vee p$\index{axiom (3) for propositional logic}\index{(3), axiom for propositional logic}\\[.5ex]
\> (4)\>$(p\supset q)\supset(r\vee p\supset r\vee q)$\index{axiom (4) for propositional logic}\index{(4), axiom for propositional logic}
\end{tabbing}

I shall discuss the meaning of these axioms later. At present I want only to say that an expression written down in our theory as an axiom or as a theorem always means that it is true for any propositions $p,q,r$ etc., e.g.\ $p\supset p\vee q$.

Now in geometry and any other theory \emph{except logic} the deductive system is completely given by stating what the primitive terms and what the axioms are. It is important to remark that it is different here for the following reason: in geometry and other theories it is clear how the theorems are to be derived from the axioms; they are to be derived by the rules of logic which are assumed to be known. In our case however we cannot assume the rules of logic to be known $\mathbf{\llbracket 54. \rrbracket}$ because we are just about to formulate the rules of logic\index{rules in logic}\index{rules of inference in logic} and to reduce them to a minimum. So this will naturally have to apply to the rules of inference as well as to the axioms with which we start. We shall have to formulate the rules of inference explicitly and with greatest possible precision, that is in such a way there can never be a doubt whether a certain rule can be applied for any formula or not. And of course we shall try to work with as few as possible. I have to warn here against an error.

One might think that an explicit formulation of the rules of inference besides the axioms is superfluous because the axioms themselves seem to express rules of inference, e.g.\ $p\supset p\vee q$ the rule that from a proposition $p$ one can conclude $p\vee q$, and one might think that the axioms themselves contain at the same time the rules by which the theorems are to be derived. But this way out of the difficulty would be entirely wrong $\mathbf{\llbracket 55. \rrbracket}$ because e.g.\ $p\supset p\vee q$ does not say that it is permitted to conclude $p\vee q$ from $p$ because those terms ``allowable to conclude'' do not occur in it. The notions in it are only $p$, $\supset$, $\vee$ and $q$. According to our definition of $\supset$ it does not mean that, but it simply says $p$ is false or $p\vee q$ is true. It is true that the axioms suggest or make possible certain rules of inference, e.g.\ the just stated one, but it is not even uniquely determined what rules of inference it suggests; e.g.\ $\sim p\vee(p\vee q)$ says either $p$ is false or $p\vee q$ is true, which suggests the rule of inference $p$ : $p\vee q$, but it also suggests $\sim(p\vee q)$ : $\sim p$. So we need written specifications, i.e.\ we have to formulate rules of inference in addition to formulas.\footnote{Here a note in a box in the manuscript mentions pp.\ 56-60 of Notebook~I.}

It is only because the ``if then'' in ordinary language is ambivalent and has besides the meaning given by the truth table also the meaning ``the second member can be inferred from the first'' that the axioms seem to express uniquely rules of inference.

$\mathbf{\llbracket 55.1 \rrbracket}$ This remark applies generally to any question whether or not certain laws of logic can be derived from others (e.g.\ whether the law of excluded middle is sufficient). Such questions have only a precise meaning if you state the rules of inference which are to be accepted in the derivation. It is different e.g.\ in geometry; there it has a precise meaning whether it follows, namely it means whether it follows by logical inference, but it cannot have this meaning in logic because then every logical law would be derivable from any other. So it could $\mathbf{\llbracket 55.2 \rrbracket}$ only mean derivable by the inferences made possible by the axioms. But as we have seen that has no precise meaning because an axiom may make possible or suggest many inferences.

\pagestyle{myheadings}\markboth{EDITED TEXT}{NOTEBOOK I \;\;---\;\; 1.1.9\;\; Axiom system for propositional logic}
$\mathbf{\llbracket Notebook\; I \rrbracket}$
$\mathbf{\llbracket 56. \rrbracket}$\footnote{In the scanned manuscript, pages numbered from \textbf{45}., with or without II, up to \textbf{55}.\ are missing in the present Notebook~I.} Now it has turned out that three rules of inference\index{rules of inference for propositional logic}\index{three rules of inference for propositional logic} are sufficient for our purposes, namely for deriving all tautologies from these formulas. Namely first the so called \emph{rule of substitution}\index{rule of substitution for propositional logic} which says:
\begin{itemize}
\item[] \emph{If we have a formula $F$} (of the calculus of propositions) \emph{which involves the propositional variables say $p_1,\ldots,p_n$ then it is permissible to conclude from it any formula obtained by sub\-stituting in $F$ for all or some of the propositional vari\-ables $p_1,\ldots,p_n$\- any arbitrary expressions, but in such a way that if a letter $p_i$ occurs in several places in $F$ we have to substitute the same formula in all places where it occurs.}
\end{itemize}
E.g.\ take the formula $(p\: .\: q\supset r)\supset [p\supset (q\supset r)]$ which is called exportation. According to the rule of substitution we can conclude from it the formula obtained by substituting $p\: .\: q$ for $r$, i.e.\ $(p\: .\: q\supset p\: .\: q)\supset [p\supset (q\supset p\: .\: q)]$. The expression which we substitute, in our case $p\: .\: q$, is quite arbitrary $\mathbf{\llbracket 57. \rrbracket}$ and it need not be a tautology or a proved formula. The only requirement is that if the same letter occurs on several places in the formula in which we substitute (as in out case the $r$) then we have to substitute the same expression in all the places where $r$ occurs as we did here. But it is perfectly allowable to substitute for different letters the same formula, e.g.\ for $q$ and $r$ and it is also allowable to substitute expressions containing variables which occur on some other places in the formula, as e.g.\ here $p\: .\: q$. It is clear that by such a substitution we get always a tautology if the expression in which we substitute is a tautology, because e.g.\ that this formula of exportation is a tautology says exactly that it is true whatever $p,q,r$ may be. So it will in particular be true if we take for $r$ the proposition $p\: .\: q$, whatever $p$ and $q$ may be $\mathbf{\llbracket 58. \rrbracket}$ and that means that the formula obtained by the substitution is a tautology.

The second rule of inference we need is the so called \emph{rule of implication}\index{rule of implication}\index{modus ponens}\index{rule 1} which reads as follows:
\begin{itemize}
\item[] \emph{If $P$ and $Q$ are arbitrary expressions then from the two premises $P,P\supset Q$ it is allowable to conclude $Q$.}
\end{itemize}
An example: take for $P$ the formula $p\: .\: q\supset p\: .\: q$ and for $Q$ the formula $p\supset (q\supset p\: .\: q))$ so that $P\supset Q$ will be the formula $(p\: .\: q\supset p\: .\: q)\supset [p\supset (q\supset p\: .\: q)]$. Then from those two premises we can conclude $p\supset (q\supset p\: .\: q)$. Again we can prove that this rule of inference is correct, i.e.\ if the two premises are tautologies then the conclusion is. Because if we assign any particular truth values to the propositional variables occurring in $P$ and $Q$, $P$ and $P\supset Q$ will both get the truth value truth because they are tautologies. Hence $Q$ will also get the truth value true if any particular truth values are assigned to its variables. Because if $P$ and $P\supset Q$ both have the truth value truth, $Q$ has also the truth. So $Q$ will have the truth value T whatever truth values are assigned to the variables occurring in it which means that it is a tautology.

Finally as the third rule of inference we have the \emph{rule of defined symbol}\index{rule of defined symbol for propositonal logic} which says (roughly speaking) that within any formula the definiens can be replaced by the definiendum and vice versa, or formulated $\mathbf{\llbracket 59. \rrbracket}$ more precisely for a particular definiens say $p\supset q$ it says:
\begin{itemize}
\item[] \emph{From a formula $F$ we can conclude any formula $G$ obtained from $F$ by replacing a part of $F$ which has the form $P\supset Q$ by the expression $\sim P\vee Q$ and vice versa}. (Similarly for the other definitions we had.)
\end{itemize}
As an example:
\begin{tabbing}
\hspace{1.7em}\= 1. \= $ \sim p\vee (p\vee q)$ \quad from the first axiom by replacing $p \supset Q$ by ${\sim p\vee Q}$\\[.5ex]
\>2.\> $\sim p\supset (\sim p \vee q)$\quad (\emph{Again clear that tautology of tautology}.)\\[.5ex]
\>\> $\sim p\supset ( p \supset q)$
\end{tabbing}
This last rule is sometimes not explicitly formulated because it is only necessary if one introduces definitions and it is superfluous in principle to introduce them because whatever can be expressed by a defined symbol can be done without (only it would sometimes be very long and cumbersome). If however one introduces definitions as we did this third rule of inference is indispensable.

Now what we shall prove is that any tautology can be derived from these four axioms\index{axioms for propositional logic}\index{four axioms for propositional logic} by means of the mentioned three rules of inference:\index{theorem, system of propositional logic}
\begin{tabbing}
$\mathbf{\llbracket 60. \rrbracket}$\\*[1ex]
\hspace{1.7em}\= (1)\quad\=$p\supset p\vee q$\index{axiom (1) for propositional logic}\index{(1), axiom for propositional logic}\\*[.5ex]
\>(2)\>$p \vee p \supset p$\index{axiom (2) for propositional logic}\index{(2), axiom for propositional logic}\\[.5ex]
\>(3)\>$p \vee q \supset q\vee p$\index{axiom (3) for propositional logic}\index{(3), axiom for propositional logic}\\[.5ex]
\>(4)\>$(p\supset q)\supset (r\vee p\supset r\vee q)$\index{axiom (4) for propositional logic}\index{(4), axiom for propositional logic}
\end{tabbing}

Let us first ascertain that all of these formulas are tautologies and let us ascertain that fact first by their meaning and then by their truth table.

The first means: If $p$ is true $p\vee q$ is true. That is clear because $p\vee q$ means at least one of the propositions $p,q$ is true, but if $p$ is true then the expression $p\vee q$ is true. The second means: If the disjunction $p\vee p$ is true $p$ is true, i.e.\ we know that the disjunction $p\vee p$ is true means that one of the two members is true, but since both members are $p$ that means that $p$ is true. The third says if $p\vee q$ is true $q\vee p$ is also true.

$\mathbf{\llbracket Notebook\; II \rrbracket}$
\pagestyle{myheadings}\markboth{EDITED TEXT}{NOTEBOOK II \;\;---\;\; 1.1.9\;\; Axiom system for propositional logic}
$\mathbf{\llbracket 61. \rrbracket}$
This does not need further explanation because the ``or'' is evidently symmetric in the two members. Finally the fourth means this: ``If $p\supset q$ then if $r\vee  p$ is true then $r\vee q$ is also true'', i.e.\ ``If you have a correct implication $p\supset q$ then you can get again a correct implication by adding a third proposition $r$ to both sides of it getting $r\vee p\supset r\vee q$''.

That this is so can be seen like this: it means ``If $p\supset q$ then if one of the propositions $r,p$ is true then also one of the propositions $r,q$ is true'', which is clear because if $r$ is true $r$ is true and if $p$ is true $q$ is true by assumption. So whichever of the two propositions $r,p$ is true always it has the consequence that one of the propositions $r,q$ is true.

$\mathbf{\llbracket 62. \rrbracket}$ Now let us ascertain the truth of these formulas by the truth-table method, combining always as many cases as possible into one case.
\begin{itemize}
\item[1.] If $p$ is F this is an implication with a false first member, hence true owing to the truth table of $\supset$; if $p$ is true then $p\vee q$ is also true according to the truth table of ``or'', hence the formula is an implication with true second member, hence again true.
\vspace{-1ex}
\item[2.] If $p$ is true this will be an implication with true second member, hence true. If $p$ is false then $p\vee p$ is a disjunction both of whose members are false, hence false according to the truth table for $\vee$. Hence in this case we have an implication with $\mathbf{\llbracket 63. \rrbracket}$ a false first member, which is true by the truth table of $\supset$.
\vspace{-1ex}
\item[3.] Since the truth table for $\vee$ is symmetric in $p,q$ it is clear that whenever the left-hand side has the truth value true also the right-hand side has it, and if the left-hand side is false the right-hand side will also be false; but an implication both of whose members are true or both of whose members are false is true by the truth table of implication, because $p\supset q$ is false only in the case when $p$ is true and $q$ false. \vspace{-1ex}
\item[4.] Here we have to consider only the following three cases:
\end{itemize}
\vspace{-4ex}
\begin{tabbing}
\hspace{5em}\=1. one of $r,q$ has the \=truth value T\\[.5ex]
\>2. both $r,q$ are F and $p$ true\\[.5ex]
\>3. both $r,q$ are F and $p$ false
\end{tabbing}
\vspace{-1ex}
$\mathbf{\llbracket 64. \rrbracket}$ These three cases evidently exhaust all possibilities.
\begin{itemize}
\item[\emph{1.}] In the first case $r\vee q$ is true, hence also $(r\vee p)\supset(r\vee q)$ is true because it is an implication with second member true; $(p\supset q)\supset(r\vee p\supset r\vee q)$ is true for the same reason.
\vspace{-1ex}
\item[\emph{2.}] In the second case $p$ is true and $q$ false, hence $p\supset q$ false, hence the whole expression is an implication with false first member, hence true.
\vspace{-1ex}
\item[\emph{3.}] In the third case all of $r$, $q$ and $p$ are false; then $r\vee p$ and $r\vee q$ are false, hence the implication $r\vee p\supset r\vee q$ is true, hence the whole formula is true because it is an implication with true second member.
\item[] So we see that the whole formula is always true.
\end{itemize}

\subsection{Theorems and derived rules of the system for\\ propositional logic}\label{1.10}
\pagestyle{myheadings}\markboth{EDITED TEXT}{NOTEBOOK II \;\;---\;\; 1.1.10\;\; Theorems and derived rules of the system\ldots}

Now I can begin with deriving other tautologies from these three axioms by means of the three rules of inference, namely the rule of \emph{substitution and implication and defin\-ed symbol}, in order to prove later on that all logically true formulas can be derived from them.

Let us first substitute $\f{\sim r}{r}$ in (4) to get $(p\supset q)\supset(\sim r\vee p\supset \; \sim r\vee q)$, but for $\sim r\vee p$ we can substitute $r\supset p$ and likewise for $\sim r\vee q,$ $\mathbf{\llbracket 65. \rrbracket}$ getting:
\begin{tabbing}
5. $(p\supset q)\supset[(r\supset p)\supset (r\supset q)]$ \quad Syllogism\index{formula of syllogism}\index{Syllogism, formula of}
\end{tabbing}
This is the so called formula of syllogism, which has a certain similari\-ty to the mood Barbara in so far as it says: If from $p$ follows $q$ then if from $r$ follows $p$ from $r$ follows $q$.

\noindent 6. Now substitute $\f{p}{q}$ in (1) $p\supset p\vee p$ and now make the following subst\-itution:
\vspace{-2ex}
\begin{tabbing}
\hspace{1.7em}$\f{\f{p\vee p}{p}\quad\f{p}{q}\quad\f{p}{r} \quad \rm{in \;\; Syllogism}}{(p\vee p\supset p)\supset[(p\supset p\vee p)\supset (p\supset p)]}$
\end{tabbing}
This is an implication and the first member of it reads $p\vee p\supset p$, which is nothing else but the second axiom. Hence we can apply the rule of implication to the $\mathbf{\llbracket 66. \rrbracket}$ two premises and get
\begin{tabbing}
\hspace{1.7em}$(p\supset p\vee p)\supset (p\supset p)$
\end{tabbing}
This is again an implication and the first member of it was proved before; hence we can again apply the rule of implication and get
\begin{tabbing}
7. $p\supset p$\quad law of identity\index{law of identity}
\end{tabbing}
Using the third rule
\begin{tabbing}
8*. we have $\sim p\vee p$ \quad the law of excluded middle\index{law of excluded middle}
\end{tabbing}
Now let us substitute $\f{\sim p}{p}$ in this formula to get {$\sim\sim p \;\vee \sim p$} and now apply to it the commutative law for $\vee$, i.e.\ substitute $\f{\sim \sim p}{p}$\quad $\f{\sim p}{q}\;\;$ in (3) to get
\begin{tabbing}
\hspace{1.7em}\=$\sim\sim p \;\vee \sim p\supset \;\sim p\;\vee \sim\sim p$ \quad rule of implication\\*[.5ex]
\>$\sim p\;\vee \sim\sim p$
\end{tabbing}
$\mathbf{\llbracket 67. \rrbracket}$ 9.* $p\supset\;\sim\sim p$

\vspace{1ex}

I have to make an important remark on how we deduced $p\supset p$ from the axioms. We had at first the two formulas $p\supset p\vee p$ and $p\vee p\; \supset p$. Now substitute them in a certain way in the formula of Syllogism $\f{p}{r}\;\;$ $\f{p\vee p}{p}\;\;$ $\f{p}{q}\;\;$ and then by applying twice the rule of implication we get $p\supset p$. If $P,Q,R$ are any arbitrary expressions and \emph{if we have succeeded in deriving $P\supset Q$ and $Q\supset R$ from the four axioms by means of the three rules of inference then we can also derive} $\mathbf{\llbracket 68. \rrbracket}$ $P\supset R$. Because we can simply substitute $\f{Q}{p}$\quad $\f{R}{q}$\quad $\f{P}{r}$ in Syllogism getting $(Q\supset R)\supset[(P\supset Q)\supset (P\supset R)]$. Then we apply the rule of implication to this formula and $Q\supset R$ getting ${(P\supset Q)}\supset (P\supset R)$ and then we apply again the rule of implication to this formula and $P\supset Q$ getting $P\supset R$.

\emph{So we know quite generally if $P\supset Q$ and $Q\supset R$ are both demonstrable then $P\supset R$ is also demonstrable whatever formula $P,Q,R$ may be} because we can obtain $P\supset R$ always in the manner just described. This fact allows us to save the trouble of repeating the whole argument by which we derived the conclusion from the two premises in each particular case, but we can state it once for all as a new $\mathbf{\llbracket 69. \rrbracket}$ rule of inference as follows:
\begin{itemize}
\item[] \emph{From the two premises $P\supset Q$, $Q\supset R$ we can conclude $P\supset R$ whatever the formulas $P,Q,R$ may be.} \quad 4.R.
\end{itemize}
So this is a fourth rule of inference, which I call Rule of syllogism.\index{rule of syllogism 4.R.} But note that this rule of syllogism is not a new independent rule, but can be derived from the other three rules and the four axioms. Therefore it is called a derived rule of inference.\index{derived rules}\index{derived rules of inference}\index{derived rules for propositional logic} So we see that in our system we cannot only derive formulas but also new rules of inference and the latter is very helpful for shortening the proofs. In principle it is superfluous to introduce such derived rules of inference because whatever can be proved with their help can also be proved without them. It is exactly this what we have shown before introducing this new rule of inference, namely we have shown that the conclusion of it can be obtained also by the former axioms and rules of inference and this was the justification for introducing it.

$\mathbf{\llbracket 70. \rrbracket}$ But although these derived rules of inference are superfluous in principle they are very helpful for shortening the proofs and therefore we shall introduce a great many of them. We now apply this rule immediately to the (1) and (3) axioms because they have this form $P \supset Q$, $Q\supset R$ for $\f{p}{P}$\quad $\f{p\vee q}{Q}$ \quad$\f{q\vee p}{R}$, and get because (1), (3)
\begin{tabbing}
10.* $p\supset q\vee p$\\[.5ex]
paradox: \= 11. \= $p\supset (q\supset p)$\quad\quad $p\supset (\sim q\vee p)$\\[.5ex]
\>\> Add$\,$* $\f{\sim q}{q}$ in last formula 10.* \\[1ex]
\> 12. [$\sim p\supset (p\supset q)$\quad\quad $\sim p\supset (\sim p\vee q)$\\[.5ex]
\>\> Add $\f{\sim p}{p}$ \quad $\f{q}{q}$] \, in\, (1)
\end{tabbing}

Other derived rules:
\begin{tabbing}
4$\cdot 1'$R $\;$\=\underline{$P_1\supset P_2\;\; P_2\supset P_3\;\; P_3\supset P_4 \;\; : \;\; P_1\supset P_4$}\quad \emph{generalized rule of syllogism}\index{generalized rule of syllogism}\index{rule of syllogism, generalized}\\[.5ex]
\hspace{6em} $P_1\supset P_3$\\[1ex]
5.R* \>\underline{$P\supset Q\;\; : \;\; R\vee P\supset R \vee Q$} \hspace{3em} \= addition from the left\index{rule of addition from the left}\index{addition from the left, rule of}\\[.5ex]
This rule is similar to the rules by which one calculates with inequalities\\[1ex]
\hspace{7em} $a<b\quad :\quad c+a<c+b$\\[1ex]
[6R \> \underline{$P\supset Q\;\; : \;\; (R\supset P)\supset (R \supset Q)$}\;]\\[1ex]
5$\cdot 1$R*\> \underline{$P\supset Q\;\; : \;\; P\vee R\supset Q \vee R$}\> addition from the right\index{rule of addition from the right}\index{addition from the right, rule of}\\[1ex]
$\mathbf{\llbracket 71. \rrbracket}$ \> 1. \hspace{.4em}\= $P\vee R\supset R\vee P$\hspace{3em}\= $\f{P}{p}$\quad $\f{R}{q}$ \quad in (3)\\[.5ex]
\> 2.\>$R\vee P\supset R\vee Q$\>by rule addition from the left\\[.5ex]
\> \underline{3.\;\, $R\vee Q\supset Q\vee R$} \>\> $\f{R}{p}$\quad $\f{Q}{q}$ \quad in (3)\\[.5ex]
\>\>$P\vee R\supset Q\vee R$\>by rule Syllogism\\[1ex]
7R* \>\underline{$P\supset Q\quad\quad R\supset S\;\; : \;\; P\vee R\supset Q \vee S$}\\[.5ex]
 \`Rule of addition of two implications\\[.5ex]
\>\>$P\vee \mbox{\underline{\emph{R}}}\supset Q\vee \mbox{\underline{\emph{R}}}$\> addition from the \=right to the\\
\` first premise ($R$)\\[.5ex]
\>\>\underline{$\mbox{\underline{\emph{Q}}}\vee R\supset \mbox{\underline{\emph{Q}}}\vee S$}\> addition from the left $\; ''\; $ second $\; ''\; $ ($Q$)\\[.5ex]
\>\>$P\vee R\supset Q\vee S$\> Syllogism, but this is the conclusion to be\\
\` proved\\[1ex]
8R* \>\underline{$P\supset Q\quad\quad R\supset Q\;\; : \;\; P\vee R\supset Q$}\quad\quad Dilemma\index{rule of dilemma}\index{dilemma, rule of}\\[.5ex]
\>\>$P\vee R\supset Q\vee Q$\\[.5ex]
\>\>\underline{$Q\vee Q\supset Q$}\> $\f{Q}{p}$ in (2)\\[.5ex]
\>\>$P\vee R\supset Q$\> Syllogism\\[1.5ex]
$\mathbf{\llbracket 72. \rrbracket}$ Application to derive formulas\\*
\>\> $p\supset \;\sim\sim p$ \> proved before, substitute $\f{\sim p}{p}$\\
\>\> $\sim p\supset \;\sim\sim\sim p$ \> addition from the right\\[.5ex]
\>\> $\sim p\vee p\supset \;\sim\sim\sim p\vee p$\> \quad rule of implication\\[.5ex]
\>\> $\sim\sim\sim p \vee p$\> rule of defined symbol\\[1ex]
13.\> \underline{$\sim\sim p\supset p$}\\[1ex]
14.\>Transposition \quad\= \underline{$(p\supset\;\sim q)\supset(q\supset\;\sim p)$}\\[.5ex]
\>Proof\> $(\sim p\;\vee\sim q)\supset(\sim q\;\vee\sim p)$ \quad\quad substitution in (3)\\
\` rule of defined symbol\\[1ex]

14$\cdot$1\>\underline{$(p\supset q)\supset(\sim q\supset\;\sim p)$}\\[.5ex]
\>$(\sim p\vee q)\supset(\sim\sim q\;\vee\sim p)$\\[.5ex]
\>Proof\hspace{2em}\= $q\supset \;\sim\sim q$\\[.5ex]
\>\>$\sim p\vee q\supset \;\sim p\;\vee \sim\sim q$\\[.5ex]
\>\>\underline{$\sim p\;\vee \sim\sim q\supset \;\sim\sim q\;\vee\sim p$}\quad Permutation (3)\\[.5ex]
\>\>$\sim p\vee q\supset\;\sim\sim q\;\vee\sim p$\quad\quad\quad rule of defined symbol\\[1ex]

14$\cdot$1\>$(p\supset q)\supset(\sim q\supset\;\sim p)$\>\>\quad$\sim p\vee q\supset\;\sim\sim q\;\vee\sim p$\\[.5ex]
14$\cdot$2\>$(\sim p\supset\;\sim q)\supset(q\supset p)$\>\>\quad$\sim\sim q\;\vee \sim p\supset\;\sim p\vee q$\footnotemark\\[.5ex]
14$\cdot$3*\>$(p\supset\;\sim q)\supset(q\supset\;\sim p)$\>\>\quad$\sim p\;\vee\sim q\supset\;\sim q\;\vee\sim p$\\[.5ex]
14$\cdot$4*\>$(\sim p\supset q)\supset(\sim q\supset p)$\>\>\quad$\sim\sim p\vee q\supset\;\sim\sim q\;\vee p$\\[1.5ex]
14$\cdot$2\>\underline{$(\sim p\supset\;\sim q)\supset(q\supset p)$}\\[1ex]
\>Proof\>$\sim\sim p\supset p$\\[.5ex]
\>\>$\sim\sim p\;\vee\sim q\supset p\;\vee\sim q$\\[.5ex]
\>\>\underline{$p\;\vee\sim q\supset\;\sim q\vee p$}\\[.5ex]
\>\>$\sim\sim p\;\vee\sim q\supset\;\sim q\vee p$\\[.5ex]
\>\>$(\sim p\supset\;\sim q)\supset(q\supset p)$
\end{tabbing}
\footnotetext{$\sim\sim p\;\vee\sim q\supset\;\sim q\vee p$, as in the proof below}
\begin{tabbing}
4$\cdot 1'$R $\;$\=Proof\hspace{2em}\= $q\supset \;\sim\sim q$\kill

$\mathbf{\llbracket 73. \rrbracket}$\\*[1.5ex]

14$\cdot$4* \>\underline{$(\sim p\supset q)\supset(\sim q\supset p)$}\\[.5ex]
\>$\sim\sim p\vee q\supset\;\sim\sim q\vee p$\\[1ex]
\>Proof\>$\sim\sim p\supset p$\\[.5ex]
\>\>\underline{$q\supset\;\sim\sim q$}\\[.5ex]
\>\>$\sim\sim p\vee q\supset p\;\vee\sim\sim q$\quad Addition of two implications\\[.5ex]
\>\>\underline{$p\;\vee\sim\sim q\supset\;\sim\sim q\vee p$}\quad Permutation\\[.5ex]
\>\>$\sim\sim p\vee q\supset\;\sim\sim q\vee p$\quad q.e.d.\quad rule of defined symbol \\[1ex]
Four transposition rules of inference:\index{rules of transposition}\index{transposition rules of inference}\\[.5ex]
9R\>\underline{$P\supset\; \sim Q\quad : \quad Q\supset\; \sim P$}\hspace{3em}\=9$\cdot$1R \quad\=\underline{$P\supset Q\quad : \quad \sim Q\supset\; \sim P$}\\[.5ex]
9$\cdot$2R \>\underline{$\sim P\supset Q\quad : \quad \sim Q\supset P$}\hspace{3em}\>9$\cdot$3R \>\underline{$\sim P\supset\;\sim Q\quad : \quad Q\supset P$}
\end{tabbing}
By them the laws for . correspond to laws for $\vee$ or can be derived, e.g.\
\begin{tabbing}
15.*\hspace{1.7em}\=\underline{$p\: .\: q\supset p$}\hspace{7em}\=\underline{$p\: .\: q\supset q$}\\[.5ex]
\>$\sim(\sim p\;\vee\sim q)\supset p$\>$\sim(\sim p\;\vee\sim q)\supset q$\qquad\=Formula 10.*\\[1ex]
Proof\>$\sim p\;\supset\;\sim p\;\vee\sim q$\>$\sim q\;\supset\;\sim p\;\vee\sim q$\>Transposition  9$\cdot$2R\\[.5ex]
\>$\sim(\sim p\;\vee\sim q)\supset p$\>$\sim(\sim p\;\vee\sim q)\supset q$
\end{tabbing}
15.2\quad Similarly for products of any number of factors we can prove that the product implies any factor, e.g.\
\begin{tabbing}
\hspace{1.7em}\=\underline{$p\: .\: q\, .\: r\supset p$}\hspace{2em}\= because $\;\;$\=$(p\: .\: q)\, .\: r\supset p\: .\: q$\\[.5ex]
\>\underline{$p\: .\: q\, .\: r\supset q$}\>\>$p\: .\: q\supset p$, \qquad $p\: .\: q\supset q$\\[.5ex]
\>\underline{$p\: .\: q\, .\: r\supset r$}\>\>$(p\: .\: q)\, .\: r\supset r$
\end{tabbing}
and for any number of factors.

From this one has the following rules of inference:
\begin{tabbing}
10R\hspace{1.7em}\=\underline{$\;P\supset Q\quad : \quad P\: .\:R\supset Q$}\hspace{1em} adjoining a new hypothesis\index{rule of adjoining a new hypothesis}\index{adjoining a new hypothesis, rule of}\\[.5ex]
10$\cdot$1R\>\underline{$\;P\supset Q\quad : \quad R\: .\:P\supset Q$}\\[.5ex]
\>because \hspace{2em}\=$P\: .\: R\supset P$\hspace{1.7em}\= by substitution\\[.5ex]
\>\>\underline{$P\supset Q$}\>by assumption\\[.5ex]
\>\>$P\: .\: R\supset Q$\>Syllogism\\[1ex]
10$\cdot$2R\>$\;\;\;$\underline{$Q$\quad : \quad $P\supset Q$}\>\> from paradox
\end{tabbing}

\begin{tabbing}
\noindent $\mathbf{\llbracket 74. \rrbracket}$ Associativity: Recall \; (1) \, $p\supset p\vee q$, $\;$ II \, $p\supset q\vee p$\\*[1ex]
15.*\quad\=\underline{$(p\vee q)\vee r\supset p\vee(q\vee r)$}\\[.5ex]
\>1.\quad\=$p\supset p\vee(q\vee r)$\qquad Addition (1)\quad $\f{q\vee r}{q}$\\[.5ex]
\>\>$q\supset q\vee r$\qquad $q\vee r\supset p\vee(q\vee r)$\qquad Formula 10.*\\
\` Addition*\quad $\f{q\vee r}{p}$\quad$\f{p}{q}$ \qquad($p\supset q\vee p$\quad $\f{q\vee r}{p}$\quad$\f{p}{q}$)\\[.5ex]
\>2.\>$q\supset p\vee(q\vee r)$\qquad Syllogism\\[.5ex]
\>a.)\> $p\vee q\supset p\vee(q\vee r)$\qquad Dilemma\\[.5ex]
\>\>$r\supset q\vee r\;$ \quad (II $\f{r}{p}$)\qquad $q\vee r\supset p\vee(q\vee r)$\quad see \emph{before}\\[.5ex]
\>b.)\> $r\supset p\vee(q\vee r)$\\[.5ex]
\>\>$(p\vee q)\vee r\supset p\vee(q\vee r)$\qquad inverse similar\\[1ex]
15$\cdot$1\>\underline{$p\vee(q\vee r)\supset (p\vee q)\vee r$}\\
\>$p\supset p\vee q$\qquad $p\vee q\supset (p\vee q)\vee r$\qquad ($p\supset p\vee q$\quad $\f{p\vee q}{p}$\quad $\f{r}{q}$)\\[.5ex]
\>$p\supset (p\vee q)\vee r$\\[.5ex]
\>$q\supset (p\vee q)\vee r$\\
\>$r\supset (p\vee q)\vee r$\qquad[II\quad $p\supset q\vee p$\quad $\f{r}{p}$\quad $\f{p\vee q}{q}$]\\[.5ex]
\>$q\vee r\supset (p\vee q)\vee r$\\[.5ex]
\>$p\vee(q\vee r)\supset (p\vee q)\vee r$\quad \\[1ex]
Exportation and importation\\[.5ex]
16.*\>\underline{$(p\: .\: q\supset r)\supset[p\supset(q\supset r)]$}\qquad Exportation\\[1ex]

$\mathbf{\llbracket 75. \rrbracket}$\>$(\sim(p\: .\: q)\vee r)\supset\:\sim p\vee(\sim q\vee r)$\\[.5ex]
\>$\sim\sim(\sim p\;\vee\sim q)\vee r\supset\;\sim p\vee(\sim q\vee r)$\\[1ex]
Proof\>$\sim\sim(\sim p\;\vee\sim q)\supset\;\sim p\;\vee\sim q$\qquad double negation\\[-1ex]
\` substitute\quad $\f{\sim p\;\vee\sim q}{p}$\\[.5ex]
\>$\sim\sim(\sim p\;\vee\sim q)\vee r\supset (\sim p\;\vee\sim q)\vee r$\quad addition from the right\\[.5ex]
\>$(\sim p\;\vee\sim q)\vee r\supset\;\sim p\vee(\sim q\vee r)$\qquad associative law\\[.5ex]
\emph{Syllogism}$\;\;\;$\underline{$\sim\sim(\sim p\;\vee\sim q)\vee r\supset\;\sim p\vee(\sim q\vee r)$\quad q.e.d.}\\[1ex]
\>\underline{$[p\supset(q\supset r)]\supset(p\: .\: q\supset r)$}\qquad Importation\\[1ex]
\>$\sim p\vee(\sim q\vee r)\supset\;\sim\sim(\sim p\;\vee\sim q)\vee r$\\[1ex]
Proof$\;\times$\quad$\sim p\vee(\sim q\vee r)\supset(\sim p\;\vee\sim q)\vee r$\hspace{2em} Associativity\\[.5ex]
\>$\sim p\;\vee\sim q\supset\;\sim\sim(\sim p\;\vee\sim q)$\\[.5ex]
$\times$\>$(\sim p\;\vee\sim q)\vee r\supset\;\sim\sim(\sim p\;\vee\sim q)\vee r$\qquad Addition right\\[.5ex]
\>$\sim p\vee(\sim q\vee r)\supset\;\sim\sim(\sim p\;\vee\sim q)\vee r$\hspace{1.8em} Syllogism $\;\times\times$\\[1ex]
\>\underline{$[p\supset(q\supset r)]\supset[q\supset(p\supset r)]$}\\[1ex]
$\times$\>$\sim p\vee(\sim q\vee r)\supset(\sim p\;\vee\;\sim q)\vee r$\\[.5ex]

$\mathbf{\llbracket 76. \rrbracket}$\>$\sim p\;\vee\sim q\supset\;\sim q\;\vee\sim p$\\*[.5ex]
$\times$\>$(\sim p\;\vee\;\sim q)\vee r\supset(\sim q\;\vee\;\sim p)\vee r$\\[.5ex]
$\times$\>$(\sim q\;\vee\;\sim p)\vee r\supset\;\sim q\vee(\sim p\vee r)$\\[.5ex]
\>\underline{$\sim p\vee(\sim q\vee r)\supset\;\sim q\vee(\sim p\vee r)$\qquad Syllogism $\;\times\times\times$}
\end{tabbing}
Rule of exportation or importation or commutativity\begin{tabbing}
11\qquad\=\underline{$P\: .\: Q\supset R\quad :\quad P\supset(Q\supset R)$}\qquad Exportation\index{rule of exportation}\index{exportation, rule of} \\[.5ex]
11$\cdot$1\>\underline{$P\supset(Q\supset R)\quad :\quad P\: .\: Q\supset R$}\qquad Importation\index{rule of importation}\index{importation, rule of} \\[.5ex]
11$\cdot$2\>\underline{$P\supset(Q\supset R)\quad :\quad Q\supset(P\supset R)$}\qquad Commutativity\index{rule of commutativity}\index{commutativity, rule of}
\end{tabbing}

\pagestyle{myheadings}\markboth{EDITED TEXT}{NOTEBOOK III \;\;---\;\; 1.1.10\;\; Theorems and derived rules of the system\ldots}
\noindent $\mathbf{\llbracket Notebook\; III \rrbracket}$
\begin{tabbing}
$\mathbf{\llbracket 1. \rrbracket}$\hspace{1em}\=$(p\supset q)\supset [(r\supset p)\supset(r\supset q)]$\\[.5ex]
\>\underline{$(q\supset r)\supset[(p\supset q)\supset(p\supset r)]$}\\[.5ex]
\>\underline{$(p\supset q)\supset [(q\supset r)\supset (p\supset r)]$}\hspace{.2em} \= Commutativity \hspace{.2em} \=$\f{q \supset r}{P}$ $\f{p \supset q}{Q}$ $\f{p \supset r}{R}$ \\[.5ex]
\>\underline{$(p\supset q)\: .\: (q\supset r)\supset (p\supset r)$} \> Importation \> $\f{p \supset q}{P}$ $\f{q \supset r}{Q}$ $\f{p \supset r}{R}$\\[.5ex]
\>\underline{$(q\supset r)\: .\: (p\supset q)\supset (p\supset r)$}\\[3ex]
\>\underline{$(p\supset q)\: .\: p \supset q$}\\[.5ex]
\> $(p\supset q)\supset (p\supset q)$\qquad $\f{p \supset q}{P}$ $\f{p}{Q}$ $\f{q}{R}$ \\[.5ex]
\>$(p\supset q)\: .\: p \supset q$\hspace{3em} Importation\end{tabbing}

\begin{tabbing}
$\mathbf{\llbracket 2. \rrbracket}$\\*[1ex]
17 \hspace{3em}\=\underline{$p\: .\: q\supset q\: .\: p$}\\[.5ex]
Proof \>$\sim q \:\vee \sim p \supset \:\sim p \:\vee \sim q$\qquad (3) $\f{\sim q}{p}$ $\f{\sim p}{q}$\\[.5ex]
\>$\sim(\sim p\;\vee\sim q)\supset \:\sim (\sim q\;\vee\sim p)$ \qquad Transposition\\[.5ex]
\>$p\: .\: q\supset q\: .\: p$ \qquad rule of defined symbol
\end{tabbing}

\begin{tabbing}
18. \hspace{3em}\=\underline{$p\supset p\: .\: p$}\\[.5ex]
Proof \>$\sim p\; \vee \sim p \supset \; \sim p$\\[.5ex]
\>$p \supset \; \sim (\sim p \;\vee \sim p)$ \qquad Transposition\\[.5ex]
\>$p\supset p\: .\: p$ \qquad defined symbol
\end{tabbing}

\begin{tabbing}
19. \hspace{3em}\=\underline{$p\supset (q\supset p\: .\: q)$}\\[-1ex]
\>$(p\: .\: q \supset p\: .\: q)\supset (p \supset (q \supset p\: .\: q))$\quad\=exportation\quad\= $\f{p\: .\: q}{r}$ \\[-1ex]
\>$p \supset (q \supset p\: .\: q)$\\[2ex]
19.1 \>\underline{$p\supset (q\supset q\: .\: p)$}\\[-1ex]
\>$(p\: .\: q \supset q\: .\: p)\supset (p \supset (q \supset q\: .\: p))$\> exportation \> $\f{q\: .\: p}{r}$
\end{tabbing}

\begin{tabbing}
$\mathbf{\llbracket 3. \rrbracket}$\\*[1ex]
12R \hspace{3em}\=\underline{$P\, ,\, Q\quad : \quad P\: .\: Q$}\qquad rule of product\index{rule of product}\index{conjunction introduction} \\[.5ex]
\>$P \supset (Q \supset P\: .\: Q)$ \\[.5ex]
\>$Q \supset P\: .\: Q$ \\[.5ex]
\>$P\: .\: Q$
\end{tabbing}

\begin{tabbing}
Inversion \hspace{3em}\=\underline{$P\: .\: Q\quad : \quad P\, ,\, Q$}\qquad rule of product\index{rule of product, inversion}\index{conjunction elimination} \\[.5ex]
\>$P\: .\: Q \supset P \quad P\: .\: Q \supset Q$
\end{tabbing}

\begin{tabbing}
13R \hspace{3em}\=\underline{$P\supset Q \qquad R\supset S\quad : \quad P\: .\: R\supset Q\: .\: S$} \qquad \emph{Rule of multiplication}\index{rule of multiplication}\\[.5ex]
\>$\sim Q\supset \: \sim P \quad \sim S\supset \: \sim R$ \\[.5ex]
\>$\sim Q \: \vee \: \sim S \supset \: \sim P \: \vee \: \sim R$ \\[.5ex]
\>$\sim(\sim P \: \vee \: \sim R) \supset \: \sim (\sim Q \: \vee \: \sim S)$
\end{tabbing}

\begin{tabbing}
$\mathbf{\llbracket 4. \rrbracket}$\\*[1ex]
\uline{13.1R} \hspace{3em}\=\underline{$P\supset Q \quad : \quad R\: .\: P\supset R\: .\: Q$} \\[.5ex]
\> because $R\supset R$ and other side
\end{tabbing}

\begin{tabbing}
\uline{13.2R} \hspace{3em}\=\underline{$P\supset Q \qquad P\supset S \quad : \quad P \supset Q\: .\: S$} \\[.5ex]
\> $P\: .\: P\supset Q\: .\: S$ \\[.5ex]
\> \underline{$P\supset P\: .\: P$} \\[.5ex]
\> $P\supset Q\: .\: S$ \qquad rule of composition
\end{tabbing}

\begin{tabbing}
F 22. \hspace{3em}\=\underline{$p\: .\:(q\vee r) \equiv p\: .\: q \vee p\: .\: r$}\\[.5ex]
I. \>$q \supset q \vee r$ \\[.5ex]
\>$p\: .\:q \supset p\: .\:(q \vee r)$ \\[.5ex]
\>$r \supset q \vee r$\\[.5ex]
\>$p\: .\:r \supset p\: .\:(q \vee r)$ \\[.5ex]
\>$p\: .\: q \vee p\: .\: r \supset p\: .\:(q\vee r)$
\end{tabbing}

\vspace{-2ex}

\noindent\rule{12cm}{0.4pt}

\vspace{-2ex}

\begin{tabbing}
II. \\[.5ex]
$\times$\hspace{.6em}\=$q \supset (p \supset p\: .\: q)$ \hspace{8.1em}\= $q\supset (p \supset p\: .\: q \vee p\: .\: r)$ \\[.5ex]
+\>$r \supset (p \supset p\: .\: r)$ \hspace{6.5em} + \>$(p\supset p\: .\: r)\supset (p \supset p\: .\: q \vee p\: .\: r)$\\[.5ex]
\>$p\: .\: q \supset p\: .\: q \vee p\: .\: r$ \> $r\supset (p \supset p\: .\: q \vee p\: .\: r)$ \\[.5ex]
\>$p\: .\: r \supset p\: .\: q \vee p\: .\: r$ \> $q\vee r \supset (p \supset p\: .\: q \vee p\: .\: r)$\\[.5ex]
$\times$\>$(p \supset p\: .\: q)\supset (p \supset p\: .\: q \vee p\: .\: r)$ \> $(q \vee r)\: .\: p \supset p\: .\: q \vee p\: .\: r$
\end{tabbing}

\noindent $\mathbf{\llbracket 5. \rrbracket}$ Equivalences

\begin{tabbing}
~~~~\hspace{3em}\=\underline{$P\supset Q\: .\: Q\supset P \quad : \quad P\equiv Q$} \\[.5ex]
\> because $(P\supset Q)\: .\: (Q\supset P)$\quad rule of defined symbol
\end{tabbing}

\begin{tabbing}
~~~~\hspace{3em}\=\underline{$P\equiv Q\quad : \quad P\supset Q\: . \: Q\supset P$}
\end{tabbing}

Transposition:

\begin{tabbing}
~~~~\hspace{3em}\=\underline{$P \equiv Q \qquad : \qquad \sim P \equiv \: \sim Q$} \\[.5ex]
\> \underline{$P \equiv \: \sim Q \qquad : \qquad \sim P \equiv Q$} \\[.5ex]
Proof \> $P \equiv Q$ \hspace{1em} $P \supset Q$ \hspace{2.3em} $Q \supset P$ \\[.5ex]
\hspace{7.5em} $\sim Q \supset \: \sim P$ \hspace{.2em} $\sim P \supset \: \sim Q$ \hspace{1em} $\sim P \equiv \: \sim Q$
\end{tabbing}

Addition and Multiplication
\begin{tabbing}
~~~~\hspace{3em}\= \underline{$P\equiv Q\qquad R\equiv S$}\hspace{2em}
$
\begin{cases}
\;\;\;\underline{P\vee R \equiv Q\vee S}\\
\;\;\;\underline{P\: .\: R \equiv Q\: .\: S}
\end{cases}
$
\end{tabbing}

\begin{tabbing}
$\mathbf{\llbracket 6. \rrbracket}$ \hspace{1em}\= $P\supset Q$ \hspace{0.5em} $R\supset S$ \hspace{1em} $Q\supset P$ \hspace{0.5em} $S\supset R$\\[.5ex]
\>$P\vee R \supset Q\vee S$ \hspace{1.2em} $Q\vee S \supset P\vee R$\\[.5ex]
\>\hspace{4.1em} $P\vee R \equiv Q\vee S$.\\[2ex]
\> Syllogism\\[.5ex]
\> \underline{$P \equiv Q\;\;,\;\;Q \equiv S \quad : \quad P \equiv S$} \\[.5ex]
\> \underline{$P \equiv Q \quad : \quad Q \equiv P$} \\[2ex]
\> \underline{$p \equiv p$} ~~~~~\qquad $p \supset p$ \quad $p \supset p$ \quad $(\f{P}{p}\: \f{Q}{p})$\\[.5ex]
\> \underline{$p \equiv \: \sim \sim p$} \qquad $p \supset \: \sim \sim p$ \quad $ \sim \sim p \supset p$\\[.5ex]
\> \underline{$\sim (p\: .\: q)\equiv \: \sim p \: \vee \sim q$}\\[0.5ex]
\> $\sim \sim (\sim p\: \vee \sim q)\equiv \: \sim p \: \vee \sim q$\\[0.5ex]
\> \underline{$\sim (p \vee q)\equiv \: \sim p\:.\sim q$}\\[0.5ex]
\> ~~~~~~~~~~~~~$\equiv \: \sim (\sim \sim p \: \vee \sim \sim q)$\\[0.5ex]
\> $p \equiv \: \sim \sim p$\\[0.5ex]
\> \underline{$q \equiv \: \sim \sim q$}\\[0.5ex]
\> $p \vee q \equiv \: \sim \sim p \: \vee \sim \sim q$ \quad $\mid$ \quad $\sim (p \vee q) \equiv \: \sim (\sim \sim p \: \vee \sim \sim q)$
\end{tabbing}

\begin{tabbing}
$\mathbf{\llbracket 6a. \rrbracket}$\\*[1ex]
23. \hspace{3em}\=\underline{$p\vee(q\: .\:r) \equiv (p\vee q)\: .\: (p \vee r)$}\\[.5ex]
~~~~~~1.) \>$p \supset p \vee q$ \\[.5ex]
\>$p \supset p \vee r$ \\[.5ex]
~~~~~~~~~~~~$\llcorner$\> $p\supset (p\vee q)\: .\: (p \vee r)$ \\[.5ex]
\>$q\: .\:r \supset p \vee q$ \quad because $q\: .\:r \supset q$\\[.5ex]
\> \underline{$q\: .\:r \supset p \vee r$}\\[.5ex]
~~~~~~~~~~~~$\llcorner$\> $q\: .\:r \supset (p\vee q)\: .\: (p \vee r)$\\[.5ex]
~~~~~~2.) ~$\llcorner$\> $p \supset [(p \vee q) \supset (p \vee q\: .\:r)]$ $\times$\\[.5ex]
\> \Big[$r \supset [(p \vee q) \supset (p \vee q\: .\:r)]$\Big]\\[.5ex]
\> $r \supset [q \supset q\: .\:r]$\\[.5ex]
\> $q \supset q\: .\:r \supset [(p \vee q)\supset (p \vee q\: .\:r)]$ \quad Summation\\[.5ex]
~~~~~~~~~~~~$\llcorner$\> $r\supset [(p \vee q) \supset (p \vee q\: .\:r)]$\\[.5ex]
\> $(p \vee r)\supset [(p \vee q)\supset (p \vee q\: .\:r)]$\\[.5ex]
\> $( p\vee r)\: .\: (p \vee q)\supset (p \vee q\: .\:r)$\\[2ex]
$\times$ because \quad \=$p\supset p \vee q\: .\:r$\\[.5ex]
\>$p \vee q\: .\:r \supset [(p \vee q)\supset(p \vee q\: .\:r)]$\\[.5ex]
\>$p \supset [(p \vee q)\supset(p \vee q\: .\:r)]$
\end{tabbing}

\noindent $\mathbf{\llbracket 7. \rrbracket}$ Syllogism under an assumption\index{syllogism under an assumption}

\begin{tabbing}
14R \hspace{3em}\=\underline{$P\supset (Q \supset R)\, ,\, P\supset (R \supset S)\quad : \quad P\supset (Q \supset S)$}\\[.5ex]
\> and similarly for any number of premises\\[.5ex]
\>$P\supset (Q \supset R)\: .\:(R \supset S)$\\[.5ex]
\> \underline{$(Q \supset R)\: .\:(R \supset S)\supset Q \supset S$} \qquad exportation syllogism \\[.5ex]
\>$P\supset (Q \supset S)$ \qquad \emph{also generalized}
\end{tabbing}

\vspace{0.5ex}
{\setlength{\parindent}{0cm}
$\left[
\begin{tabular}{@{}l@{}}
  14.1R \hspace{3em} \underline{$P\supset Q \qquad P \supset (Q \supset R) \quad : \quad P\supset R$} \\[.5ex]
\hspace{5.7em} $P\supset (Q \supset R)\: .\:Q$ \\[.5ex]
  \hspace{5.7em} \underline{$(Q \supset R)\: .\:Q \supset R$} \\[.5ex]
  \hspace{5.7em} $P\supset R$ \qquad Syllogism\end{tabular}
\right]$
}
\vspace{0.5ex}

\begin{tabbing}
$\mathbf{\llbracket 8. \rrbracket}$ \hspace{.6em}\=\underline{$(p\supset q)\: .\:(r\supset s) \supset (p \vee r \supset q \vee s)$}\\[.5ex]
1. \> $p \vee r \supset r \vee p$\\[.5ex]
2. \> $(p \supset q)\supset (r \vee p \supset r \vee q)$\\[.5ex]
3. \> $r \vee q \supset q \vee r$\\[.5ex]
4. \> $(r \supset s)\supset (q \vee r \supset q \vee s)$\\[.5ex]
5. \> $(p \supset q)\: .\: (r \supset s)\supset (p \vee r \supset r \vee p)$\\[.5ex]
6. \> $(p \supset q)\: .\: (r \supset s)\supset (r \vee p \supset r \vee q)$\\[.5ex]
7. \> $(p \supset q)\: .\: (r \supset s)\supset (r \vee q \supset q \vee r)$\\[.5ex]
8. \> \underline{$(p \supset q)\: .\: (r \supset s)\supset (q \vee r \supset q \vee r)$}\\[.5ex]
9. \> \underline{$(p \supset q)\: .\: (r \supset s)\supset (p \vee r \supset q \vee s)$}\\[.5ex]

\>\underline{$(p\supset q)\: .\:(r\supset q) \supset (p \vee r \supset q)$}\\[.5ex]
\> $(p\supset q)\: .\:(r\supset q) \supset (p \vee r \supset q \vee q)$ \qquad $\f{q}{s}$ \\[.5ex]
\> $(p\supset q)\: .\:(r\supset q) \supset (q \vee q \supset q)$\\[.5ex]
\> $(p\supset q)\: .\:(r\supset q) \supset (p \vee r \supset q)$\\[1ex]
$\mathbf{\llbracket 9. \rrbracket}$\>\underline{$(p\supset q)\: .\:(r\supset s) \supset (p \: .\: r \supset q \: .\: s)$}\\[.5ex]
\> $(p\supset q)\supset (\sim q \supset \: \sim p)$\\[.5ex]
\> $(r\supset s)\supset (\sim s \supset \: \sim r)$\\[.5ex]
$A.$ \> $(p\supset q)\: .\:(r\supset s) \supset (\sim q \supset \: \sim p) \: .\: (\sim s \supset \: \sim r)$\\[.5ex]
$B.$ \> $(\sim q \supset \: \sim p) \: .\: (\sim s \supset \: \sim r) \supset (\sim q \: \vee \sim s \supset \: \sim p \: \vee \sim r)$\\[.5ex]
$C.$ \> $(\sim q \: \vee \sim s \supset \: \sim p \: \vee \sim r) \supset (p \: .\: r \supset q \: .\: s)$\\[.5ex]
\> $(p\supset q)\: .\:(r\supset s) \supset (p \: .\: r \supset q \: .\: s)$ \qquad $A,B,C$\\[.5ex]
\>{$(p\supset q)\: .\:(p\supset s) \supset (p \supset q \: .\: s)$}\\[.5ex]
\> $(p\supset q)\: .\:(p\supset s) \supset ( p\: .\: p \supset q \: .\: s)$\\[.5ex]
\> $(p\supset q)\: .\:(p\supset s) \supset (p \supset p\: .\: p)$\\[.5ex]
\> $(p\supset q)\: .\:(p\supset s) \supset (p \supset q \: .\: s)$\\[1ex]
\>\underline{$(p\supset \: \sim p)\supset \: \sim p$}\\[.5ex]
\> $\sim p \: \vee \sim p \supset \: \sim p$\\[1ex]
$\mathbf{\llbracket 10. \rrbracket}$ \>\underline{$(\sim p\supset p)\supset p$}\\[.5ex]
\> $(\sim \sim p \vee p) \supset p$\\[.5ex]
\> $\sim \sim p \supset p$\\[.5ex]
\> \underline{$p \supset p$}\\[.5ex]
\> $(\sim \sim p \vee p) \supset p$\\[.5ex]
\>\underline{$\sim (p \: . \sim p)$} \qquad see below$^\ast$\\[.5ex]
\>\underline{$(p \supset q).(p \supset \: \sim q)\supset \: \sim p$}\\[.5ex]
\> $(p \supset q)\: .\: (p \supset \: \sim q)\supset[p \supset (q \: . \sim q)]$\\[.5ex]
\> $p \supset (q \: . \sim q)\supset (\sim(q \: . \sim q)\supset \: \sim p)$\\[.5ex]
\> $(p \supset q)\: .\: (p \supset \: \sim q)\supset (\sim(q \: . \sim q)\supset \: \sim p)$ ~ $\diagdown$\\
\` Principle of Commutativity\\
\> $\sim(q \: . \sim q)\supset [(p \supset q).(p \supset \: \sim q)\supset \:\sim p]$~~~~ $\diagup$\\[.5ex]
\> $(p \supset q).(p \supset \: \sim q)\supset \: \sim p$\\[1ex]
\>\underline{$\sim (p \: . \sim p)$}\\*[.5ex]
\hspace{1.7em}$^\ast$ \> $\sim \sim(\sim p \: \vee \sim \sim p)$
\end{tabbing}

\pagestyle{myheadings}\markboth{EDITED TEXT}{NOTEBOOK IV \;\;---\;\; 1.1.10\;\; Theorems and derived rules of the system\ldots}

\begin{tabbing}
$\mathbf{\llbracket Notebook\; IV\rrbracket}$\\*[.5ex]
$\mathbf{\llbracket new\, page\; i\rrbracket}$\footnotemark
\hspace{2.5em}\= $(p\supset q)\supset [(r \supset p)\supset(r \supset q)]$\quad\= 1. \\[.5ex]
 \> $p\supset \;\sim\sim p$ \>2.\\[1ex]
\hspace{26em}$R$:\,  $\sim p$\\[.5ex]
\hspace{26em}$S$:\, $\sim\sim\sim p$\\[.5ex]
\hspace{26em}$T$:\, $p$\\[.5ex]
Su 2. \> $\sim p\supset \;\sim\sim\sim p$ \` $R\supset S$ \\[.5ex]
 \> $\sim p\vee p\supset \;\sim\sim\sim p\vee p$ \` $R\vee T\supset S\vee T$ \\[.5ex]
Su (3) \> $\sim p\vee p\supset p\;\vee\sim p$ \>3. \` $R\vee T\supset T\vee R$ \\[.5ex]
Su (4) \> $(\sim p\supset\;\sim\sim\sim p)\supset [p\;\vee \sim p\supset p\;\vee\sim\sim\sim p]$\quad 4. \\[.5ex]
Imp 2., 4. \> $p\;\vee \sim p\supset p\;\vee\sim\sim\sim p$\> 5. \` $T\vee R\supset T\vee S$ \\[.5ex]
Su (3) \> $p\;\vee\sim\sim\sim p\supset\;\sim\sim\sim p\vee p$\> 6. \` $T\vee S\supset S\vee T$ \\[.5ex]
Su 1. \>$(p\;\vee\sim p\supset p\;\vee\sim\sim\sim p)\supset$\\ \` $[(\sim p\vee p \supset p\;\vee\sim p)\supset(\sim p\vee p \supset p\;\vee\sim\sim\sim p)]$\quad 7.\\[.5ex]
Imp twice 5., 7.; 3. \>$\sim p\vee p \supset p\;\vee\sim\sim\sim p$\> 8.\\[.5ex]
Su 1. \>$(p\;\vee\sim\sim\sim p\supset \;\sim\sim\sim p\vee p)\supset$\\ \` $[(\sim p\vee p \supset p\;\vee\sim\sim\sim p)\supset(\sim p\vee p \supset \;\sim\sim\sim p\vee p)]$\quad 10.\\[.5ex]
Imp twice 6., 10.; 8. \>$\sim p\vee p \supset \;\sim\sim\sim p\vee p$
\end{tabbing}
\footnotetext{Before p.\ {\bf 7}., the first numbered page in Notebook IV, there are in the manuscript four not numbered pages with theorems of the axiom system for propositional logic. These pages are here numbered with the prefix $\mathbf{new\, page}$ and inserted within Notebook III, at the end of the present Section 1.1.10, to which they belong by their subject matter.}

\begin{tabbing}
$\mathbf{\llbracket new\, page\; ii\rrbracket}$
\hspace{2.2em}\=\underline{$p\supset q\vee p$}\\*[.5ex]
\> $p\supset p\vee q$\hspace{3em} \=(1) \\[.5ex]
\> \underline{$p\vee q\supset q\vee p$}\> (3) \\[1ex]
Su 1.\hspace{5em}\> $(p\vee q\supset q\vee p)\supset [(p \supset p\vee q)\supset(p \supset q\vee p)]$\hspace{1.3em}\=2. \\
\` Su\quad $\f{p\vee q}{p}$\quad $\f{q\vee p}{q}$\quad$\f{p}{r}$\\[.5ex]
Imp (2., (3))\> $(p \supset p\vee q)\supset(p \supset q\vee p)$\> 3. \\[-2ex]
\underline{\hspace{29.7em}}\\
Imp (3., (1))\> $p \supset q\vee p$\> 4.
\end{tabbing}

\begin{tabbing}
$\mathbf{\llbracket new\, page\; iii\rrbracket}$\\*[1ex]
\emph{1.} \hspace{1em}\=$(\sim p\supset p)\supset p$\hspace{3em}\= $(\sim\sim p\vee p)\supset p$\\*[1ex]
\> A. \hspace{.27em} $p\supset p$\\[.5ex]
\> \underline{$\sim\sim p\supset p$}\\[.5ex]
\> $\sim\sim p\vee p\supset p$\> Dilemma\\[1ex]
\emph{2.}\>$(p\: .\: q\supset r)\supset (p\: .\sim r\supset
\;\sim q)$\\[1ex]
1.\> $(p\: .\: q\supset r)\supset[p\supset(q\supset r)]$\hspace{4.5em}\=Exportation\\[.5ex]
\>$(q\supset r)\supset(\sim r\supset \:\sim q)$\> Transposition\\[.5ex]
2.\> $[p\supset(q\supset r)]\supset[p\supset(\sim r\supset \:\sim q)]$\> Addition from the left\\[.5ex]
3. \> $[p\supset(\sim r\supset \:\sim q)]\supset [p\: .\sim r\supset \:\sim q]$\> Importation\\[-2ex]
\underline{\hspace{17em}}\\
\> $(p\: .\: q\supset r)\supset (p\: .\sim r\supset\;\sim q)$\> 1., 2., 3. Syllogism\\[1ex]
\emph{3.1}\> $(p\supset q)\supset(p\supset(p\supset q))$\\[.5ex]
\> $r\supset(p\supset r)$\hspace{9em}\>$\f{p\supset q}{r}$\\[.5ex]
\emph{3.2}\> $[p\supset(p\supset q)]\supset(p\supset q)$\>$\sim p\vee(\sim p\vee q)\supset \;\sim p\vee q$\\[2ex]

$\mathbf{\llbracket new\, page\; iv\rrbracket}$\\*[1ex]

1.\> $\sim p\vee(\sim p\vee q)\supset(\sim p\;\vee\sim p)\vee q$\\[.5ex]
\> $\sim p\;\vee\sim p\supset\; \sim p$\\[.5ex]
2.\> $(\sim p\;\vee\sim p)\vee q\supset\;\sim p \vee q$\>Addition from the right\\[.5ex]
\> $\sim p\vee(\sim p\vee q)\supset\;\sim p \vee q$\> Syllogism 1., 2.\\[.5ex]
\>$[p\supset(p\supset q)]\supset(p\supset q)$\>Rule of defined symbol
\end{tabbing}

\subsection{Completeness of the axiom system for\\ propositional logic}\label{1.11}
\pagestyle{myheadings}\markboth{EDITED TEXT}{NOTEBOOK III \;\;---\;\; 1.1.11\;\; Completeness of the axiom system for\ldots}

$\mathbf{\llbracket Notebook\; III\rrbracket}$ $\mathbf{\llbracket 11. \rrbracket}$ Now I can proceed to the proof of the completeness theorem\index{completeness theorem, propositional logic} announced in the beginning which says that any tautology whatsoever can actually be derived in a finite number of steps from our four axioms by application of the three primitive rules of inference (substitution, implication, defined symbol) or shortly ``Every tautology is demonstrable''. We have already proved the inverse theorem which says: ``Every demonstrable expression is a tautology''.

$\mathbf{\llbracket 12. \rrbracket}$ But the proposition which we are interested in now is the inverse one, which says ``Any tautology is demonstrable''. In order to prove it we have to use again the formulas $P_i$ which we used for proving that any truth table function can be expressed by $\sim$ and $\vee$. If we have say $n$ propositional variables $p_1, p_2, p_3,\ldots, p_n$ then consider the conjunction of them $p_1\: .\:p_2\: .\:p_3\: .\:\:\ldots\:\: .\:p_n$ and call a ``fundamental conjunction''\index{fundamental conjunction} of these $\mathbf{\llbracket 13. \rrbracket}$ letters $p_1, \ldots, p_n$ any expression obtained from this conjunction by negating some or all of the variables $p_1, \ldots, p_n$. So e.g.\ ${p_1\: . \sim p_2\:.\: p_3\:.\:\:\ldots\:\:.\: p_n}$ would be a fundamental conjunction, another one ${\sim p_1\:.\: p_2\:. \sim p_3\:.\: p_4\:.\:\:\ldots\:\: .\:p_n}$ etc.; in particular we count also ${p_1\:.\:\:\ldots\:\: .\:p_n}$ itself and ${\sim p_1\:. \sim p_2\:.\:\:\ldots\:\: . \sim p_n}$ (in which all variables are negated) as fundamental conjunctions.

\begin{tabbing}
\hspace{1.7em}\=$2$ for one \qquad \=$p_1$, \quad $\sim p_1$\\[0.5ex]
$2^2$\>$4$ for two\>$p_1\:.\:p_2$,\quad $p_1\:. \sim p_2$,\quad $\sim p_1\:.\:p_2$,\quad $\sim p_1\:. \sim p_2$\\[0.5ex]
$2^3$\>$8$ for three\>$p_1\:.\:p_2\:.\:p_3$,\quad $p_1\:.\:p_2\:.\sim p_3$,\quad $p_1\:.\sim p_2\:.\:p_3$,\quad  $p_1\:.\sim p_2\:.\sim p_3$ \\[0.5ex]
\>$\sim p_1\:.\:p_2\:.\:p_3$,\quad $\sim p_1\:.\:p_2\:. \sim p_3$,\quad $\sim p_1\:. \sim p_2\:.\:p_3$,\quad $\sim p_1\:. \sim p_2\:. \sim p_3$
\end{tabbing}

So for the $n$ variables $p_1,\ldots, p_n$ there are exactly $2^n$ fundamental conjunc\-tions in general; $2^n$ because you see by adding a new variable $p_{n+1}$ the number of fundamental conjunctions is doubled, because we can combine $p_{n+1}$ and $\sim p_{n+1}$ with any of the previous $\mathbf{\llbracket 14. \rrbracket}$ fundamental conjunctions (as e.g.\ here $p_3$ with any of the previous four and $\sim p_3$ getting eight). I denote those $2^n$ fundamental con\-junctions for the variables $p_1, \ldots, p_n$ by $P^{(n)}_1, P^{(n)}_2, \ldots, P^{(n)}_i, \ldots, P^{(n)}_{2^n}$. I am using $^{(n)}$ as an upper index to indicate that we mean the fundamental con\-junction of the $n$ variables $p_1, \ldots, p_n$. The order in which they are enumerated is arbitrary. [We may stick e.g.\ to the order which we used in the truth tables.] From our formulas considered for $n=3$ we know $\mathbf{\llbracket 14.1 \rrbracket}$ that to each of these fundamental conjunctions $P^{(n)}_i$ corresponds exactly one line in a truth table for a function of the $n$ variables $p_1, \ldots, p_n$ in such a way that $P^{(n)}_i$ will be true in this line and false in all the others. So if we numerate the lines correspondingly we can say $P^{(n)}_i$ will be true in the $i^{\text{th}}$ line and false in all other lines.

$\mathbf{\llbracket 15. \rrbracket}$ Now in order to prove the completeness theorem I prove first the following auxiliary theorem.\index{completeness proof, auxiliary theorem}
\begin{itemize}
\item[] \emph{Let $E$ be any expression which contains no other propositional variables but $p_1, \ldots, p_n$ and $P^{(n)}_i$ any funda\-mental conjunction of the variables $p_1, \ldots, p_n$. Then either $P^{(n)}_i \supset E$ or $P^{(n)}_i \supset \: \sim E$ is demonstrable}
\end{itemize}
\emph{where by} either or \emph{I mean at least one}.\footnote{perhaps ``at most one'', or ``exactly one''}
\begin{tabbing}
\hspace{11.7em}$E$\\
Example\: \=$p_1\:.\:p_2\:.\:p_3 \supset [p\: .\:q \supset r]$ \hspace{7em} $p_1\:. \sim p_2\:.\:p_3$ \\[.5ex]
\> \underline{$p_1\:. \sim p_2\:.\:p_3 \supset (p_1\:.\:p_2 \supset p_3)$} \qquad or\\[.5ex]
\> $p_1\:. \sim p_2\:.\:p_3 \supset\: \sim (p_1\:.\:p_2 \supset p_3)$\\[.5ex]
\> \underline{$\sim p \:. \sim q \: .\: r \supset\: \sim (p \: . \: q \supset r)$}
\end{tabbing}

It is to be noted that $E$ need not actually contain all the variables $p_1, \ldots, p_n$; it is only required that it contains no other variables but $p_1, \ldots, p_n$. So e.g.\ $p_1 \vee p_2$ would be an expression for which the theorem applies, i.e.\
\begin{equation*}
\begin{rcases}
P^{(n)}_i\supset(p_1 \vee p_2)\;\;\;\;\\
\;\;\;\;\;\supset \; \sim (p_1 \vee p_2)
\end{rcases}
\text{demonstrable}
\end{equation*}

$\mathbf{\llbracket 19. \rrbracket}$\footnote{The text that follows should be a continuation of p.\ \textbf{15}.\ of the present Notebook III, according to a note at the bottom of that page. Page \textbf{16}.\ is crossed out in the manuscript and pages \textbf{17}.-\textbf{18}.\ are missing in the scanned manuscript.} I shall prove the auxiliary theorem only for such expressions as contain only the primitive symbols $\sim , \vee$ (but do not contain $\supset , \equiv$) because that is sufficient for our purpose, and I prove it by a kind of complete induction, which we used already once in order to show that $\vee$ cannot be defined in terms of $\sim , \equiv$ . $\mathbf{\llbracket 20. \rrbracket}$ Namely I shall prove the following three lemmas:

\begin{itemize}
\item [1.] The theorem is true for the simplest kind of expression $E$, namely the variables $p_1, \ldots, p_n$ themselves, i.e.\ for any variable $p_k$ of the above series $p_1, \ldots, p_k$ and any fundamental conjunction $P^{(n)}_i$, $P^{(n)}_i \supset p_k$ or $P^{(n)}_i \supset \: \sim p_k$ is demonstrable.\vspace{-1ex}
\item[2.] If the theorem is true for an expression $E$, then it is also true for the negation $\sim E$.\vspace{-1ex}
\item[3.] If it is true for two expressions $G,H$ then it is also true for the expression $G \vee H$.\end{itemize}

After having proved these three lemmas we are finished. Because any expres\-sion $\mathbf{\llbracket 21. \rrbracket}$ $E$ containing only the variables $p_1, \ldots, p_n$ and the operations $\sim , \vee$ is formed by iterated application of the operations $\sim , \vee$ beginning with the variables $p_1, \ldots, p_n$. Now by $(1.)$ we know that the theorem is true for the variables $p_1, \ldots, p_n$ and by $(2.)$ and $(3.)$ we know that it remains true if we form new expressions by application of $\sim$ and $\vee$ to expressions for which it is true. Hence it will be true for any expression of the considered type. So it remains only to prove these three auxiliary propositions.

$\mathbf{\llbracket 22. \rrbracket}$ $(1.)$ means: For any variable $p_k$ (of the series $p_1, \ldots, p_n$) and any fundamental conjunction $P^{(n)}_i$ either $P^{(n)}_i \supset p_k$ or $P^{(n)}_i \supset \: \sim p_k$ is demonstra\-ble. But now the letter $p_k$ or the negation $\sim p_k$ must occur among the members of this fundamental conjunction $P^{(n)}_i$ by definition of a fundamental conjunction. On the other hand we know that for any conjunction it is demonstra\-ble that the conjunction implies any of its members. (I proved that explicitly for conjunctions of two and three members and remarked that the same method will prove it for conjunctions of any $\mathbf{\llbracket 23. \rrbracket}$ number of members. The exact proof would have to go by an induction on the number of members. For two, proved. Assume $P^{(n)}$ has $n$ members and $p$ is a variable among them. Then $P^{(n)}$ is $P^{({n-1})} \: .\: r$:

$1.$ $p$ occurs in $P^{({n-1})}$; then $P^{({n-1})} \supset p$, hence $P^{({n-1})} \: .\: r \supset p$.

$2.$ $r$ is $p$; then $P^{({n-1})} \: .\: p \supset p$ is demonstrable.) Hence if $p_k$ occurs among the members of $P^{(n)}_i$ then $P^{(n)}_i \supset p_k$ is demonstrable and if $\sim p_k$ occurs among them then $P^{(n)}_i \supset \: \sim p_k$ is demonstrable. So one of these two formulas is demonstrable in any case and that is exactly the assertion of lemma $(1.)$.

Now to $(2.)$, i.e.\ let us assume the theorem is true for $E$, i.e.\ for any fundamental conjunction $P^{(n)}_i$ either $P^{(n)}_i \supset E$ or $P^{(n)}_i \supset \: \sim E$ is demonstrable and let us show that the theorem is true also for the expression $\sim E$, i.e.\ for any $P^{(n)}_i$either $P^{(n)}_i \supset \: \sim E$ or $P^{(n)}_i \supset \: \sim (\sim E)$ is demonstrable for any $P^{(n)}_i$
\begin{center}
\begin{tabular}{l|l}
$P^{(n)}_i \supset E$ & $P^{(n)}_i \supset \: \sim E$ \\
$P^{(n)}_i \supset \: \sim E$ & $P^{(n)}_i \supset \: \sim (\sim E)$ \\
\end{tabular}
\end{center}
(because it is $\mathbf{\llbracket 24. \rrbracket}$ this what the theorem says if applied to $\sim E$).
But now in the first case if $P^{(n)}_i \supset E$ is demonstrable then $P^{(n)}_i \supset \: \sim (\sim E)$ is also demonstrable because $E \supset \: \sim (\sim E)$ is demonstrable by substitution in the law of double negation, and if both $P^{(n)}_i \supset E$ and $E \supset \: \sim (\sim E)$ are demonstrable then also $P^{(n)}_i \supset \: \sim (\sim E)$ by the rule of syllogism. So we see if the first case is realized for $E$ then the second case is realized for $\sim E$ and of course if the second case is realized for $E$ the first case is realized for $\sim E$ (because they say the same thing). $\mathbf{\llbracket 25. \rrbracket}$ So if one of the two cases is realized for $E$ then also one of the two cases is realized for $\sim E$, i.e.\ if the theorem is true for $E$ it is also true for $\sim E$ which was to be proved.

Now to $(3.)$. Assume the theorem true for $G , H$ and let $P^{(n)}_i$ be any arbi\-trary fundamental conjunction of $p_1, \ldots, p_n$. Then $P^{(n)}_i \supset G$ is demonstra\-ble or $P^{(n)}_i \supset \: \sim G$ is demonstrable and $P^{(n)}_i \supset H$ is demonstrable or $P^{(n)}_i \supset \: \sim H$ is demonstrable by assumption and we have to prove from these assump\-tions that also:

\begin{tabbing}
	\hspace{1.7em}\=$P^{(n)}_i \supset G\vee H$ \hspace{3em} \=or\\[0.5ex]
	\>$P^{(n)}_i \supset \: \sim(G\vee H)$ \> is demonstrable.
\end{tabbing}
In order to do that distinguish three cases:

\vspace{1ex}

\noindent $\mathbf{\llbracket 26. \rrbracket}$

\begin{itemize}
\item [1.] [For $G$ first case realized, i.e.] $P^{(n)}_i \supset G$ is demonstrable; then we have $G \supset G\vee H$ also by substitution in axiom, hence $P^{(n)}_i \supset G\vee H$ demonstrable by rule of syllogism [hence first case realized for $G \vee H$].\vspace{-1ex}
\item [2.] case [For $H$ first case realized] $P^{(n)}_i \supset H$ is demonstrable; then $H \supset G\vee H$ by substitution in formula 10.*, hence $P^{(n)}_i \supset G\vee H$ is demonstrable by rule of syllogism [hence first case realized for $G \vee H$].\vspace{-1ex}
\item [3.] case Neither for $G$ is $P^{(n)}_i \supset G$ nor for $H$ is $P^{(n)}_i \supset H$ the first case realized. Thus for both of them second case happens, i.e.\ $P^{(n)}_i \supset \: \sim G$ and $P^{(n)}_i \supset \: \sim H$ are both demonstrable by assumption, but then by rule of transposition $G \supset \: \sim P^{(n)}_i$ and $H \supset \: \sim P^{(n)}_i$ are demon\-strable. Hence $G \vee H \supset \: \sim P^{(n)}_i$ by rule of Dilemma. Hence $P^{(n)}_i \supset \: \sim (G \vee H)$ by transposition [i.e.\ second case realized for $G \vee H$].\end{itemize}

$\mathbf{\llbracket 27. \rrbracket}$ So we see in each of the three cases which exhaust all possibilities ei\-ther $P^{(n)}_i \supset G\vee H$ or $P^{(n)}_i \supset \: \sim(G\vee H)$ is demonstrable, namely the first happens in case $1$ and $2$, the second in case $3$. But that means that the theorem is true for $G \vee H$ since $P^{(n)}_i$ was any arbitrary fundamental conjunction. So we have proved the three lemmas and therefore the auxiliary theorem for all expressions $E$ containing only $\sim , \vee$.

Now\index{completeness proof, main theorem} let us assume in particular that $E$ is a tautology of this kind (i.e.\ containing only the letters $p_1, \ldots, p_n$ and only $\sim , \vee$); then I maintain $\mathbf{\llbracket 28. \rrbracket}$ that $P^{(n)}_i \supset E$ is demonstrable for any fundamental conjunction $P^{(n)}_i$. Now we know from the preceding theorem that certainly either $P^{(n)}_i \supset E$ or $P^{(n)}_i \supset \: \sim E$ is demonstrable. So it remains only to be shown that the second case, that $P^{(n)}_i \supset \: \sim E$ is demonstrable, can never occur if $E$ is a tautology and that can be shown as follows: As I mentioned before any demonstrable proposition is a tautology. But on the other hand we can easily see that $P^{(n)}_i \supset \: \sim E$ is certainly not a tautology if $E$ is a tautology because the truth value of $P^{(n)}_i \supset \: \sim E$ will be false $\mathbf{\llbracket 29. \rrbracket}$ in the $i^{th}$ line of its truth table. For in the $i^{th}$ line $P^{(n)}_i$ is true as we saw before and $E$ is also true in the $i^{th}$ line because it is assumed to be a tautology, hence true in any line. Therefore $\sim E$ will be false in the $i^{th}$ line, therefore $P_i \supset \: \sim E$ will be false in the $i^{th}$ line because $P_i$ is true and $\sim E$ false and therefore $P_i \supset \: \sim E$ false by the truth table of $\supset$. So this expression $P_i \supset \: \sim E$ has F in the $i^{th}$ line of its truth table, hence is not a tautology, hence cannot be demonstrable and therefore $P^{(n)}_i \supset E$ is demonstrable for any fundamental conjunction $P^{(n)}_i$, if $E$ $\mathbf{\llbracket 30. \rrbracket}$ is a tautology containing only $\sim , \vee , p_1, \ldots, p_n$.

But from the fact that $P^{(n)}_i \supset E$ is demonstrable for any $P^{(n)}_i$ it follows that $E$ is demonstrable in the following way: We can show \emph{first} that also for any fundamental conjunction $P^{({n-1})}_i$ of the ${n-1}$ variables $p_1, \ldots, p_{{n-1}}$, $P^{({n-1})}_i \supset E$ is demonstrable because if $P^{({n-1})}_i$ is a fundamental conjunction of the ${n-1}$ variables $p_1, \ldots, p_{{n-1}}$ then $P^{({n-1})}_i \: .\: p_n$ is a fundamental conjunction of the $n$ variables $p_1, \ldots, p_n$ and likewise $P^{({n-1})}_i \: . \sim p_n$ is a fundamental conjunction of the $n$ variables $p_1, \ldots, p_n$; therefore by our previous theorem $\mathbf{\llbracket 31. \rrbracket}$ $P^{({n-1})}_i \: .\: p_n \supset E$ and $P^{({n-1})}_i \: . \sim p_n \supset E$ are both demonstrable. Applying the rule of exportation and commu\-tativity to those two expressions we get $p_n \supset (P^{({n-1})}_i \supset E)$ and $\sim p_n \supset (P^{({n-1})}_i \supset E)$ are both demonstrable. To be more exact we have to apply first the rule of exportation and then the rule of commutativity because the rule of exportation gives $P^{({n-1})}_i \supset (p_n \supset E)$. But now we can apply the rule of dilemma to these two formulas ($P \supset R,Q \supset R \: : \: P \vee Q \supset R$) and obtain $\sim p_n \vee p_n \supset (P^{({n-1})}_i \supset E)$ is demonstrable; and now since $\sim p_n \vee p_n$ is demonstrable we can apply the rule of implication again and obtain $P^{({n-1})}_i \supset E$ is demonstrable which was to be shown. Now since this holds $\mathbf{\llbracket 32. \rrbracket}$ for any fundamental conjunction $P^{({n-1})}_i$ of the ${n-1}$ variables $p_1, \ldots, p_{n-1}$ it is clear that we can apply the same argument again and prove that also for any fundamental conjunction $P^{(n-2)}_i$of the $n-2$ variables $p_1, \ldots, p_{n-2}$, $P^{(n-2)}_i \supset E$ is demonstrable. So by repeating this argument ${n-1}$ times we can finally show that for any fundamental conjunction of the one variable $p_1$ this implication is demonstrable, but that means $p_1 \supset E$ is demonstrable and $\sim p_1 \supset E$ is demonstrable (because $p_1$ and $\sim p_1$ are the fundamental conjunc\-tion of the one variable $\mathbf{\llbracket 33. \rrbracket}$ $p_1$), but then $\sim p_1 \vee p_1 \supset E$ is demonstrable by rule of dilemma and therefore $E$ is demonstrable by rule of implication.

Incidentally so we have shown that any tautology containing only $\sim$ and $\vee$ is demonstrable, but from this it follows that any tautology whatsoever is demonstrable because: let $P$ be one containing the defined symbols $. \: , \supset , \equiv$ . I then denote by $P'$ the expression obtained from $P$ by replacing $. \: , \supset , \equiv$ by their definiens, i.e.\ $R\: . \: S$ by $\sim (\sim R \: \vee \sim S)$ wherever it occurs in $P$ etc. Then $P'$ will also be a tautology. But $P'$ is a tautology containing only $\sim, \vee$ hence $P'$ is demonstrable, but then also $P$ is demonstrable because it is obtained from $P'$ by one or several applications of the rule of defined symbol, namely since $P'$ was obtained from $P$ by replacing $p\: .\: q$ by $\sim (\sim p\: \vee \: \sim q)$ etc.\ $P$ is obtained from $P'$ by the inverse substitution, but each such substitution is an application of rule of defined symbol, hence: If $P'$ is demonstrable then also $P$ is demonstrable.

As an example take the formula $(p \supset q) \vee (q \supset p)$ which is a tautology.
\begin{tabbing}
\hspace{1.7em}\= 1. Without defined symbols \quad $(\sim p \vee q) \vee (\sim q \vee p)=E$\\[.5ex]
\> 2. Fundamental conjunctions in $p,q$\\
\` $p\: .\: q$, \hspace{.5em} $p\: . \sim q$, \hspace{.5em} $\sim p\: .\: q$, \hspace{.5em} $\sim p\: . \sim q$
\end{tabbing}
To prove that $p\: .\: q \supset E$ etc.\ are all demonstrable we have to verify our auxiliary theorem successively for all particular formulas, i.e.\ for $p$, $q$, $\sim p$, $\sim q$, $\sim p \vee q$, $\sim q \vee p$, $E$.

\vspace{2ex}

\noindent $\mathbf{\llbracket 34. \rrbracket}$
\begin{center}
\begin{tabular}{ l|c|c|c|c|c|c}
 & $p$ & $q$ & $\sim p$ & $\sim q$ & $\sim p \vee q$ & $\sim q \vee p$\\[.5ex] \hline
 \hspace{1.1em}$p \: .\: q \supset$ & $p$ & $q$ & $\sim (\sim p)$ & $\sim (\sim q)$ & $\sim p \vee q$ & $\sim q \vee p$ \\
 \hspace{1.1em}$p \: .\:\! \sim q \supset$ & $p$ & $\sim q$ & $\sim (\sim p)$ & $\sim q$ & $\sim (\sim p \vee q)$ & $\sim q \vee p$ \\
 $\sim p \: .\: q \supset$ & $\sim p$ & $q$ & $\sim p$ & $\sim (\sim q)$ & $\sim p \vee q$ & $\sim (\sim q \vee p)$ \\
 $\sim p \: .\:\! \sim q \supset$ & $\sim p$ & $\sim q$ & $\sim p$ & $\sim q$ & $\sim p \vee q$ & $\sim q \vee p$
 \end{tabular}
\end{center}
\begin{flushright}
\begin{tabular}{ l|c }
  & $(\sim p \vee q) \vee (\sim q \vee p)$\\[.5ex] \hline
  & $E$\\
  & $E$\\
  & $E$\\
  & $E$
 \end{tabular}
\end{flushright}
\begin{tabbing}
\hspace{1.7em}\=$p\: . \sim q \supset \: \sim (\sim p)$ \qquad \=$\sim p \supset\: \sim (p \: . \sim q)$ \\[.5ex]
\> $p\: . \sim q \supset \: \sim q$ \qquad \> $q \supset\: \sim (p \: . \sim q)$ \\[.5ex]
\>\> $\sim p \vee q \supset \: \sim (p\: . \sim q) $\\[.5ex]
\>\> \underline{$p \:. \sim q \supset \: \sim (\sim p \vee q)$}\\[2ex]

\>$p\: .\: q \supset E$ \>$p \supset (q \supset E)$\\[.5ex]
\> $\sim p\: .\: q \supset E$ \>\underline{$\sim p \supset (q \supset E)$} \\
\>\> $\sim p \vee p \supset (q \supset E)$\\[.5ex]
\>\> $q \supset E$\\

\>\rule{8cm}{0.4pt}\\[1ex]

\>$p\: . \sim q \supset E$ \>$p \supset (\sim q \supset E)$\\[.5ex]
\> $\sim p\: . \sim q \supset E$ \>\underline{$\sim p \supset (\sim q \supset E)$} \\[.5ex]
\>\> $\sim p \vee p \supset (\sim q \supset E)$\\[.5ex]
\>\> $\sim q \supset E$\\[2ex]
$\mathbf{\llbracket 35. \rrbracket}$ \>\>$\sim q \vee q \supset E$ \hspace{5em} $E$
\end{tabbing}

\subsection{Independence of the axioms}\label{1.12}
\pagestyle{myheadings}\markboth{EDITED TEXT}{NOTEBOOK III \;\;---\;\; 1.1.12\;\; Independence of the axioms}

Now after having proved that any tautology can be derived from the four axioms, the next question which arises is whether all of those four axioms are really necessary to derive them or whether perhaps one or the other of them is superfluous. That would mean one of them could be left out and nevertheless the remaining three would allow to derive all tautologies. If this were the case then in particular also the superfluous axiom (since it is a tautology) could be derived from the three other, $\mathbf{\llbracket 36. \rrbracket}$ i.e.\ it would not be independent from the other. So the question comes down to investigating the independence of the four axioms from each other.\index{independence of the axioms} That such an investigation is really necessary is shown very strikingly by the last development. Namely when Russell\index{Russell, Bertrand} first set up this system of axioms for the calculus of propositions he assumed a fifth axiom, namely the associative law for disjunction and only many years later it was proved by P.\ Bernays\index{Bernays, Paul} that this associative law was superfluous, i.e.\ could $\mathbf{\llbracket 37. \rrbracket}$ be derived from the others. You have seen in one of the previous lectures how this derivation can be accomplished. But Bernays\index{Bernays, Paul} has shown at the same time that a similar thing cannot happen for the four remaining axioms, i.e.\ that they are really independent from each other.

Again here as in the completeness proof the interest does not lie so much in proving that these particular four axioms are independent but in the method to prove it, because so far we have only had an opportunity to prove that certain propositions follow from other propositions. But now we are confronted with the opposite problem to show that certain propositions do not follow from certain others and this problem requires evidently an entirely new method for its solution. This method is very interesting and somewhat connected with the questions of many-valued logics.

You know the calculus of propositions can be interpreted as an algebra in which $\mathbf{\llbracket 38. \rrbracket}$ we have the two operations of logical addition and multiplication as in usual algebra but in addition to them a third opera\-tion, the negation and besides some operations defined in terms of them ($\supset, \equiv$ etc.). The objects to which those operations are applied are the proposi\-tions. So the propositions can be made to correspond to the numbers of ordin\-ary algebra. But as you know all the operations $.\:, \vee$ etc. which we introduced are ``truth functions'' and therefore it is only the truth value of the proposi\-tions that really matters in this algebra, $\mathbf{\llbracket 39. \rrbracket}$ i.e.\ we can consider them as the numbers of our algebra instead of the propositions (simply the two ``truth values'' T and F). And this is what we shall do, i.e.\ our algebra (as opposed to usual algebra) has only two numbers T, F and the result of the operations $.\:,\vee ,\sim$ applied to these two numbers is given by the truth table, i.e.\ T $\vee$ F = T (i.e.\ the sum of the two numbers T and F is T) T $\vee$ T = T, F $\vee$ T = T, F $\vee$ F = F, $\sim$T = F, $\sim$F = T. In order to stress $\mathbf{\llbracket 40. \rrbracket}$ more the analogy to algebra I shall also write $1$ instead of T and $0$ instead of F. Then in this notation the rules for logical multiplication would look like this: $1\:.\:1=1$ , $0\:.\:1=0$, $1\:.\:0=0$, $0\:.\:0=0$.
If you look at this table you see that logical and arithmetical multiplication exactly coincide in this notation. Now what are the tautologies considered from this algebraic standpoint? They are expressions $f(p , q , r , \ldots)$ which have always the value $1$ whatever numbers $p,q,r$ may be, $\mathbf{\llbracket 41. \rrbracket}$ i.e.\ in algebraic language expressions identically equal to one $f(p , q , \ldots)=1$ and the contradictions expressions identically zero $f(p , q , \ldots)=0$. So an expression of usual algebra which would correspond to a contradiction would be e.g.\ $x^2-y^2-(x+y)(x-y)$; this is equal to~0.

But now from this algebraic standpoint nothing can prevent us to consider also other similar algebras with say three numbers $0,1,2$ instead of two and with the operations $\vee , \:.\: , \sim$ defined in some different manner. For any such algebra we shall have tautologies, $\mathbf{\llbracket 42. \rrbracket}$ i.e.\ formulas equal to 1 and contradictions equal to 0, but they will of course be different formulas for different algebras. Now such algebra with three and more num\-bers were used by Bernays\index{Bernays, Paul} for the proof of independence, e.g.\ in order to prove the independence of the second axiom Bernays\index{Bernays, Paul} considers the following algebra:
\begin{tabbing}
\hspace{1.7em}\=$3$ numbers \hspace{1em}\=$0,1,2$\\[.5ex]
\>negation \>$\sim 0=1$ \qquad $\sim 1=0$ \qquad \=$\sim 2=2$ \\[.5ex]
\>addition \>$1 \vee x=x \vee 1=1$ \> $2 \vee 2=1$\\[.5ex]
\>\> $0 \vee 0=0$ \> $2 \vee 0=0 \vee 2=2$
\end{tabbing}
Implication and other operations need not be defined separately because $p \supset q=\; \sim p \vee q$.

$\mathbf{\llbracket 43. \rrbracket}$ A tautology is a formula equal to 1, e.g.\ $\sim p \vee p$ because for $p$ equal to 0 or 1 it is equal to 1, because the operations for $0,1$ as arguments coincide with the operations of the usual calculus of propositions; if $p=2$ then $\sim p=2$ and $2 \vee 2=1$ is also true. Also $p \supset p$ is a tautology because by definition it is the same as $\sim p \vee p$.

Now for this algebra one can prove the following proposition:\begin{itemize}
\item [1.] Axioms $(1),(3),(4)$ are tautologies in this algebra.\vspace{-1ex}
\item [2.] For each of the three rules of inference we have: If the premises are tautologies in this algebra then so is the conclusion.
\item[$\mathbf{\llbracket 44. \rrbracket}$] I.e.\
\begin{itemize}
\item [1.] If $P$ and $P \supset Q$ are tautologies then $Q$ is a tautology.
\item [2.] If $Q'$ by substitution from $Q$ and $Q$ is a tautology then also $Q'$ is a tautology.
\item [3.] If $Q'$ is obtained from $Q$ by replacing $P \supset Q$ by $\sim P \vee Q$ etc. and $Q$ is a tautology then also $Q'$ is a tautology.\end{itemize}
\item[3.] The axiom $(2)$ is not a tautology in this algebra.\end{itemize}

After having shown these three lemmas we are finished because by $1,2$: Any formula demonstrable from axioms $(1),(3),(4)$ by the three rules of inference is a tautology for our algebra but axiom $(2)$ is not a tautolo\-gy for our $\mathbf{\llbracket 45. \rrbracket}$ algebra. Hence it cannot be demonstrable from $(1),(3),(4)$.

Now to the proof of the lemmas $1,2,3$. First some auxiliary theorems (for $1$ I say true and for $0$ false because for $1$ and $0$ the tables of our algebra coincide with those for T and F):
\begin{itemize}
\item [1.] $p \supset p$ \qquad (we had that before, because $\sim p \vee p=1$ also $\sim 2 \vee 2=1$)\vspace{-1ex}
\item [2.] $1 \vee p=p \vee 1=1$ \qquad $0 \vee p=p \vee 0=p$\vspace{-1ex}
\item [3.] $p \vee q=q \vee p$\vspace{-1ex}
\item [4.] Also in our three-valued algebra we have: An implication whose first member is $0$ is $1$ and an implication whose second member is $1$ is also $1$ whatever the other member may be, i.e.\ $0\supset p=1$ and $p\supset 1=1$ because:

\end{itemize}

\begin{tabbing}
\hspace{1.7em}\=1.) $0\supset p= \: \sim 0 \vee p=1 \vee p=1$\\[1ex]
$\mathbf{\llbracket 46. \rrbracket}$\\*[1ex]
\>2.) $p \supset 1= \: \sim p \vee 1=1$\\[1ex]
Now (1) \quad $p \supset p \vee q=1$\\[1ex]
\>1. $p=0$ \quad $\rightarrow$ \quad $p \supset p \vee q=1$\\[.5ex]
\>2. $p=1$ \quad $\rightarrow$ \quad $1 \supset 1 \vee q=1 \supset 1=1$\\[1ex]
(3) $p \vee q=q \vee p$ \quad $\rightarrow$ \quad $p \vee q=q \vee p=1$\\[1ex]
(4) $(p \supset q)\supset (r \vee p \supset r \vee q)$ \quad $E$\\[1ex]
\>\underline{1.} \,$r=0$ \quad $r \vee p=p$ \quad $r \vee q=q$ \quad $E=(p \supset q)\supset (p \supset q)=1$\\[5ex]
\>\underline{2.} \,$r=1$ \quad $r \vee p=r\vee q=1$ \quad\!\! $E=(p \supset q)\supset (1 \supset 1)=(p \supset q)\supset 1=1$
\end{tabbing}
\begin{tabbing}
$\mathbf{\llbracket 47. \rrbracket}$\\*[1ex]
\hspace{1.7em}\=\underline{3.}\, $r=2$\\

\>\hspace{1.7em}\=\underline{$\alpha .)$}\,\,\=$q=1,2$\quad $r \vee q$\!\!\!

\begin{tabular}{ l }
 =\;\underline{$2 \vee 1$}\;=\;1 \\
 =\;\underline{$2 \vee 2$}\;=\;1 \\
 \end{tabular} \\[.5ex]
\>\>\>$r \vee p \supset r \vee q=1$\\[.5ex]
\>\>\>$(p \supset q)\supset(r \vee p\supset r \vee q)=1$\\[1ex]

\>\>\underline{$\beta .)$} \>$q=0$\\[1ex]
\>\>\hspace{1.7em}\>\underline{1.}\;\:\= $p=0$ \qquad $r \vee p=r \vee q$\\[.5ex]
\>\>\>\>$(r \vee p)\supset(r \vee q)=1$\\[.5ex]
\>\>\>\>$(p \supset q)\supset (r\supset p)\supset(r \vee q)=1$\\[1ex]
\>\>\>\underline{2.}\> $p=1$ \qquad $p \supset q=0$\\[.5ex]
\>\>\>\>$E=1$\\[1ex]
\>\>\>\underline{3.}\> $p=2$ \\[.5ex]
\>\>\>\>$(2\supset 0)\supset(2 \vee 2 \supset 2 \vee 0)=2\supset(1 \supset 2)=2\supset 2=1$
\end{tabbing}

\begin{tabbing}
$\mathbf{\llbracket 48. \rrbracket}$ Lemma 2.\ A. \quad\= $p=1$ \quad $p \supset q=1$ \quad $\rightarrow$ \quad $q=1$\\[.5ex]
\>$1=\:\sim p \vee q=0 \vee q=q$
\end{tabbing}

\hspace{1em}Hence if $f(p, q , \ldots)=1$ then
$$ \f{f(p, q , \ldots) \supset g(p, q , \ldots)=1}{g(p, q , \ldots)=1}$$

\begin{itemize}
\item[B.] Rule of substitution holds for any truth-value algebra, i.e.\ if $f(p, q , \ldots)=1$ then $f(g(p, q , \ldots), q , \ldots)=1$.
\item[C.] Rule of defined symbol likewise holds because $p \supset q$ and $\sim p \vee q$ have the same truth table.
\end{itemize}
\begin{tabbing}
$\mathbf{\llbracket 49. \rrbracket}$ Lemma 3.\ (2) $p \vee p \supset p$ is not a tautology,  i.e.\\[1.5ex]
\hspace{6em} $2 \vee 2 \supset 2$ = $1\supset 2=\:\sim 1 \vee 2=0 \vee 2=2 \neq 1$
\end{tabbing}
So the lemmas are proved and therefore also the theorem about the independence of Axiom $(2)$.

\subsection{Remark on disjunctive and conjunctive normal\\ forms}\label{1.13}
\pagestyle{myheadings}\markboth{EDITED TEXT}{NOTEBOOK III \;\;---\;\; 1.1.13\;\; Remark on conjunctive and disjunctive\ldots}

We have already developed a method for deciding of any given expression whether or not it is a tautology, namely the truth-table method. I want to develop another method which uses the analogy of the rules of the $\mathbf{\llbracket 50. \rrbracket}$ calculus of proposition with the rules of algebra. We have the two distributive laws:

\begin{tabbing}
\hspace{1.7em}\=$\underline{p}\:.\:(q\vee r)\equiv (\underline{p}\:.\:q)\vee(\underline{p}\:.\:r)$ \qquad \=$\underline{p}\:.\:q\equiv q$\\[0.5ex]
\>$\underline{p}\vee (q\:.\: r)\equiv (\underline{p}\vee q)\:.\:(\underline{p}\vee r)$\> $\underline{p}\vee q\equiv q$
\end{tabbing}
In order to prove them by the shortened truth-table method I use the following facts which I mentioned already once at the occasion of one of the exercises:
\begin{tabbing}
	\hspace{1.7em}\=if $p$ is true \qquad \=$p\:.\:q\equiv q$\\[.5ex]
	\>if $p$ is false \>$p \vee q\equiv q$
\end{tabbing}
In order to prove those equivalences I distinguish two cases: $1.$ $p$ true and $2.$ $p$ false.\footnote{There seems to be a gap in the text here.}

$\mathbf{\llbracket 51. \rrbracket}$ Now the distributive laws in algebra make it possible to decide of any given expression containing only letters and $+, -, \cdot$ whether or not it is identically zero, namely by factorizing out all products of sums, e.g.\ $x^2-y^2-(x+y)(x-y)=0$. A similar thing is to be expected in the algebra of logic. Only two differences: $1.$ In logic we have the negation which has no analogue in algebra. But for negation we have also a kind of distributive law given by the De Morgan formulas\index{De Morgan formulas (laws)} $\sim (p \vee q) \equiv \: \sim p \:. \sim q$ $\mathbf{\llbracket 52. \rrbracket}$ and $\sim (p\:.\:q)\equiv\: \sim p \:\vee \sim q$. (Proved very easily by the truth-table method.) These formulas allow us to get rid of the negations by shifting them inwards to the letters occurring in the expression. The \emph{second} difference is that we have two distributive laws and therefore two possible ways of factorizing. If we use the first law we shall get as the final result a sum of products\index{disjunctive normal form} of single letters as in algebra. By using the other law of distribution we get a product of sums\index{conjunctive normal form} unlike in algebra. I think it is best to explain that on an $\mathbf{\llbracket 53. \rrbracket}$ example:
\begin{tabbing}
\hspace{1.7em}\=$\times$ 1. \hspace{.1em}\=$(p \supset q)\supset (\sim q \supset \:\sim p)$\hspace{14em}\=\\[.5ex]
\>\>$\sim (\sim p \vee q)\vee (q \:\vee \sim p)$ \\[.5ex]
\>\> $(p\:. \sim q)\vee q \:\vee \sim p$ \> disjunctive\\[.5ex]
\>\> $(p \vee q \:\vee \sim p)\:.\:(\sim q \vee q \vee \sim p)$ \> conjunctive\\[2ex]

\>$\times$ 2. \>$(p \supset q)\:.\:(p \supset \:\sim q)\:.\:p$\\[.5ex]
\>\>$(\sim p \vee q)\:.\:(\sim p \:\vee \sim q)\:.\:p$ \> conjunctive\\[.5ex]
\>\> $(\sim p\:. \sim p \vee q\:. \sim p \:\vee \sim p\:. \sim q \vee q\:. \sim q)\:.\: p$ \\[.5ex]
\>\> $ (\sim p \: .\: p)\vee (q \:. \sim p\:.\: p)\vee (\sim p\:. \sim q \:.\: p)\vee (q\:. \sim q\:.\: p)$\> disjunctive\\[2ex]

\>\;\;\; 3. \>$( p \supset q)\supset (r \vee p \supset r \vee q)$\\[.5ex]
\>\>$\sim (\sim p \vee q) \vee [\sim(r \vee p)\vee r \vee q]$\\[.5ex]
\>\> $(p\:. \sim q)\vee (\sim r \:. \sim p) \vee r \vee q$ \> disjunctive\\[.5ex]
\>\> $(p \:\vee \sim r \vee r \vee q)\:.\:(p \:\vee \sim p \vee r \vee q)\: .\:$ \\[.5ex]
\>\> \hspace{4em}$(\sim q \:\vee \sim r \vee r \vee q)\:.\:(\sim q \:\vee \sim p \vee r \vee q)$ \> conjunctive
\end{tabbing}

\subsection{Sequents and natural deduction system}\label{1.14}
\pagestyle{myheadings}\markboth{EDITED TEXT}{NOTEBOOK III \;\;---\;\; 1.1.14\;\; Sequents and natural deduction system}

$\mathbf{\llbracket 1. \rrbracket}$\footnote{Here the numbering of pages in the present Notebook III starts anew with \textbf{1}.} In the last two lectures a proof for the completeness of our system of axioms for the calculus of propositions was given, i.e.\ it was shown that any tautology is demonstrable from these axioms. Now a tautology is exactly what in traditional logic would be called a law of logic or a logically true proposition. $\mathbf{\llbracket 2. \rrbracket}$ Therefore this completeness proof solves for the calculus of propositions the second of the two problems which I announced in the beginning of my lectures, namely it shows how all laws of a certain part of logic namely of the calculus of propositions can be deduced from a finite number of logical axioms and rules of inference.

I wish to stress that the interest of this result\index{completeness, interest of} does not lie so much in this that our particular four axioms and three rules are sufficient to deduce everything, $\mathbf{\llbracket 3. \rrbracket}$ but the real interest consists in this that here for the first time in the history of logic it has really been \emph{proved} that one \emph{can} reduce all laws of a certain part of logic to a few logical axioms. You know it has often been claimed that this can be done and sometimes the laws of identity, contradiction, excluded middle have been considered as the logical axioms. But not even the shadow of a proof was given that every logical inference can be derived from them. Moreover the assertion to be proved was not even clearly formulated, because $\mathbf{\llbracket 4. \rrbracket}$ it means nothing to say that something can be derived e.g.\ from the law of contradiction unless you specify in addition the rules of inference which are to be used in the derivation.

As I said before it is not so very important that just our four axioms are sufficient. After the method has once been developed, it is possible to give many other sets of axioms which are also sufficient to derive all tautologies of the calculus $\mathbf{\llbracket 5. \rrbracket}$ of propositions, e.g.\
\begin{tabbing}
	\hspace{1.7em}\=$p \supset (\sim p \supset q)$\\[0.5ex]
	\>$(\sim p \supset p)\supset p$\\[0.5ex]
	\>$(p \supset q)\supset [(q \supset r)\supset(p \supset r)]$	
\end{tabbing}

I have chosen the above four axioms because they are used in the standard textbooks of logistics. But I do not at all want to say that this choice was particularly fortunate. On the contrary our system of axioms is open to some objections from the aesthetic point of view; e.g.\ one of the aesthetic requirements for a set of axioms is that the axioms should be as simple and evident as possible, in any case simpler than the theorems to be proved, whereas in our system $\mathbf{\llbracket 6. \rrbracket}$ e.g.\ the last axiom is pretty complicated and on the other hand the very simple law of identity $p \supset p$ appears as a theorem. So in our system it happens sometimes that simpler propositions are proved from more complicated axioms, which is to be avoided if possible. Recently by the mathematician G.\ Gentzen\index{Gentzen, Gerhard} a system was set up which avoids these disadvantages. I want to reference briefly about this system\footnote{At the end of Notebook III there are in the manuscript thirteen not numbered pages with formulae, sometimes significant, and jottings. Since it would be too intrusive to make a selection of what would be appropriate for the edited text, they are not given here.}
\pagestyle{myheadings}\markboth{EDITED TEXT}{NOTEBOOK IV \;\;---\;\; 1.1.14\;\; Sequents and natural deduction system}
$\mathbf{\llbracket Notebook\; IV \rrbracket}$ $\mathbf{\llbracket 7. \rrbracket}$\footnote{The present p.\ \textbf{7}., is in the manuscript the first numbered page of Notebook IV. It is there preceded by four pages, which have been fitted in this edited text at the end of Section 1.1.10 \emph{Theorems and derived rules of the system for propositional logic}.} or to be more exact on a system which is based on Gentzen's\index{Gentzen, Gerhard} idea, but simpler than his. The idea consists in introducing another kind of implication (denoted by an arrow $\rightarrow$).\footnote{The remainder of p.\ \textbf{7}.\ is crossed out in the manuscript, but since pp.\ \textbf{8}.-\textbf{9}.\ in the present Notebook IV are missing in the scanned manuscript, and because of the interest of this part of the text, this crossed out remainder is cited here: ``such that $P\rightarrow Q$ means $Q$ is true under the assumption $P$. The difference of this implication as opposed to our former one is

1. There can be any number of premises, e.g.\ $P,Q\rightarrow R$ means $R$ holds under the assumptions $P,Q$ (i.e.\ the same thing which would be denoted by $P\: .\: Q\supset R$. In particular the number of premises\dots''

The next page, p.\ \textbf{10}., begins with the second part of a broken sentence.}

$\mathbf{\llbracket 10. \rrbracket}$ system with altogether three primitive terms $\rightarrow$, $\sim$, $\supset$. We have now to distinguish between expressions in the former sense, i.e.\ containing only $\sim$, $\supset$ and variables, e.g.\ $p\supset q$, $\sim p\supset q$, $q\supset p\vee r$, etc., and secondary formulas\index{sequent}\index{secondary formula (sequent)} containing the arrow, e.g.\ $p,p\supset q\rightarrow q$. I shall use capital Latin letters $P,Q$ only to denote expressions of the first kind, i.e.\ expressions in our former sense, and I use capital Greek letters $\Delta,\Gamma$ to denote sequences of an arbitrary number of assumptions $\underbrace{P,Q,R\ldots}_\text{$\Delta$}$

$\mathbf{\llbracket 11. \rrbracket}$ Hence a formula of Gentzen's\index{Gentzen, Gerhard} system will always have the form $\Delta\rightarrow S$, a certain sequence of expressions of the first kind implies an expression of the first kind. Now to the axioms and rules of inference.

\vspace{1ex}

I\quad Any formula $P\rightarrow P$ where $P$ is an arbitrary expression of the first kind is an axiom and only those formulas are axioms.

\vspace{1ex}

\noindent $\mathbf{\llbracket 12. \rrbracket}$ So that is the law of identity which appears here as an axiom and as the only axiom.

\vspace{1ex}

As to the rules of inference we have four, namely

\vspace{1ex}

\noindent 1.\quad The rule of addition of premises,\index{rule of addition of premises for sequents} i.e.\ from $\Delta\rightarrow A$ one can conclude $\Delta,P \rightarrow A$ and $P,\Delta\rightarrow A$, i.e.\ if $A$ is true under the assumptions $\Delta$ then it is a fortiori true under the assumptions $\Delta$ and the further assumption $P$.

\begin{tabbing}
$\mathbf{\llbracket 13. \rrbracket}$\\*[1ex]
2.\quad \= The Rule of exportation:\index{rule of exportation for sequents}\\*[1ex]
\> $\Delta,P\rightarrow Q$\hspace{2em}\= : \hspace{1.2em} $\Delta\rightarrow (P\supset Q)$
\end{tabbing}
If the propositions $\Delta$ and $P$ imply $Q$ then the propositions $\Delta$ imply that $P$ implies $Q$.
\begin{tabbing}
3.\quad The Rule of implication:\index{rule of implication for sequents}\index{implication rule for sequents}\index{modus ponens for sequents}\\*[1ex]
\hspace{1.2em}\begin{tabular}{ l|l }
$\Delta\rightarrow P$ & \\[-1ex]
& \hspace{1.5em}$\Delta\rightarrow Q$ \\[-1ex]
$\Delta\rightarrow(P\supset Q)$ &
\end{tabular}
\end{tabbing}
So that is so to speak the rule of implication under some assumptions: If $A$ and $A\supset B$ both hold under the assumptions $\Delta$ then $B$ also holds under the assumptions $\Delta$. \begin{tabbing}
4.\quad Rule of Reductio ad absurdum or rule of indirect proof:\index{rule of reductio ad absurdum for sequents}\index{rule of indirect proof for sequents}\\[1ex]
\hspace{1.2em}\begin{tabular}{ l|l }
$\Delta,\sim P\rightarrow Q$ & \\[-1ex]
& \hspace{1.2em}$\Delta\rightarrow P$ \\[-1ex]
$\Delta,\sim P\rightarrow\;\sim Q$ &
\end{tabular}
\end{tabbing}
Here the premises mean that from the assumptions $\Delta$ and $\sim P$ a contradiction follows, i.e.\ $\sim P$ is incompatible $\mathbf{\llbracket 14. \rrbracket}$ with the assumptions $\Delta$, i.e.\ from $\Delta$ follows $P$.

Again it can be proved that every tautology follows from the axioms and rules of inference. Of course only the tautologies which can be expressed in terms of the symbols introduced, i.e.\ $\sim$, $\supset$ and $\rightarrow$. If we want to introduce also $\vee$,\;$.$ etc.\ we have to add the rule of the defined symbol $.$ or other rules concerning $\vee$,\;$.$ etc.

Now you see that in this system the aforementioned disadvantages have been avoided. All the axioms are really very simple and $\mathbf{\llbracket 15. \rrbracket}$ evident. It is particularly interesting that also the pseudo-paradoxical propositions about the implication follow from our system of axioms although nobody will have any objections against the axioms themselves, i.e.\ everybody would admit them if we interpret both the $\rightarrow$ and the $\supset$ to mean ``if\ldots\ then''. Perhaps I shall derive these pseudo-paradoxes as examples for derivations from this system. The first reads:
\begin{tabbing}
By \= I\hspace{1.5em}\= $q\rightarrow q$\kill

\>\>$q\rightarrow p\supset q$\hspace{2em} Proof:\\[1ex]

$\mathbf{\llbracket 16. \rrbracket}$\\[1ex]

By \> I\> $q\rightarrow q$\\[.5ex]
$\; ''$\> 1\> $q,p\rightarrow q$\\[.5ex]
$\; ''$\> 2\> $q\rightarrow (p\supset q)$
\end{tabbing}
Incidentally, again applying 2 we get $\rightarrow q\supset(p\supset q)$ which is another form for the same theorem. The second paradox reads like this:
\begin{tabbing}
By \= I\hspace{1.5em}\= $q\rightarrow q$\kill

\>\>$\sim p\rightarrow p\supset q$\hspace{2em} Proof:\\[1ex]
\> I\> $p\rightarrow p$\\[.5ex]
\> 1\> $\sim p,p,\sim q\rightarrow p$\\[.5ex]
\> I\> $\sim p\rightarrow \:\sim p$\\[.5ex]
\> 1\> $\sim p,p,\sim q\rightarrow \:\sim p$\\[.5ex]
\> 4\> $\sim p,p\rightarrow q$\\[.5ex]
\> 2\> $\sim p\rightarrow (p\supset q)$
\end{tabbing}

$\mathbf{\llbracket 17. \rrbracket}$ Incidentally this formula $\sim p,p\rightarrow q$ which we derived as an intermediate step of the proof is interesting also on its own account; it says: From a contradictory assumption everything follows since the formula is true whatever the proposition $q$ may be. I am sorry I have no time left to go into more details about this Gentzen\index{Gentzen, Gerhard} system. I want to conclude now this chapter about the calculus of proposition.\footnote{Here p.\ \textbf{17}. ends and pp.\ \textbf{18}.-\textbf{23}. are missing in the scanned manuscript from the present Notebook IV.}

\section{Predicate logic}\label{2}

\subsection{First-order languages and valid formulas}\label{2.1}
\pagestyle{myheadings}\markboth{EDITED TEXT}{NOTEBOOK IV \;\;---\;\; 1.2.1\;\; First-order languages and valid formulas}

$\mathbf{\llbracket 24. \rrbracket}$ I am concluding now the chapter about the calculus of propositions and begin with the next chapter which is to deal with the so called calculus of functions or predicates. As I explained formerly the calculus of propositions is characterized by this that only propositions as a whole occur in it. The letters $p,q,r$ etc. denoted arbitrary propositions and all the formulas and rules which we proved are valid whatever propositions $p,q,r$ may be, i.e.\ they are independent of the structure of the propositions involved. Therefore we could use single letters $p,q\ldots$ to denote whole propositions.

$\mathbf{\llbracket 25. \rrbracket}$ But now we shall be concerned with inferences which depend on the structure of the propositions involved and therefore we shall have to study at first how propositions are built up of their constituents. To this end we ask at first what do the simplest propositions which one can imagine look like. Now evidently the simplest kind of propositions are those in which simply some predicate is asserted of some subject, e.g.\ Socrates is mortal. Here the predicate mortal is asserted to belong to the subject Socrates. Thus far we are in agree- $\mathbf{\llbracket 26. \rrbracket}$ ment with classical logic.

But there is another type of simple proposition which was very much neglected in classical logic,\index{failure of traditional logic} although this second type is more important for the applications of logic in mathematics and other sciences. This second type of simple proposition consists in this that a predicate is asserted of several subjects, e.g.\ New York is larger than Washington. Here you have two subjects, New York and Washington, and the predicate larger says that a certain relation\index{relations} subsists between those two subjects. Another example is ``Socrates is the teacher of Plato''. So you see there are two different kinds $\mathbf{\llbracket 27. \rrbracket}$ of predicates, namely predicates with one subject\index{properties} as e.g.\ \emph{mortal} and predicates with several subjects\index{relations} as e.g.\ \emph{greater}.

The predicates of the first kind may be called properties,\index{properties} those of the second kind are called relations.\index{relations} So e.g.\ ``mortal'' is a property, ``greater'' is a relation. Most of the predicates of everyday language are relations and not properties. The relation ``greater'' as you see requires two subjects and therefore is called a dyadic relation.\index{dyadic relation (binary relation)}\index{relation, dyadic (binary)}\index{binary relation} There are also relations which require three or more subjects, e.g.\ \emph{betweenness} is a relation with three subjects, i.e.\ triadic relation.\index{triadic relation}\index{relation, triadic} If I say e.g.\ New York $\mathbf{\llbracket 28. \rrbracket}$ lies between Washington and Boston., the relation of betweenness is asserted to subsist for the three subjects New York, Washington and Boston, and always if I form a meaningful proposition involving the word between I must mention three objects of which one is to be in between the others. Therefore ``betweenness'' is called a triadic relation and similarly there are tetradic, pentadic relations etc. Properties may be called monadic predicates\index{monadic predicate}\index{predicate, monadic}\index{properties} in this order of ideas.

I don't want to go into any discussions of what predicates are (that could lead $\mathbf{\llbracket 29. \rrbracket}$ to a discussion of nominalism and realism). I want to say about the essence of a predicate only this. In order that a predicate\index{predicate} be well-defined it must be (uniquely and) unambiguously determined of any objects (whatsoever) whether the predicate belongs to them or not. So e.g.\ a property is given if it is uniquely determined of any object whether or not the predicate belongs to it and a dyadic relation is given if it is \dots uniquely determined of any two objects whether or not the relation subsists between them. I shall use capital letters $M,P,$ to denote individual predicates---as e.g.\ mortal, greater etc.\ $\mathbf{\llbracket 30. \rrbracket}$ and small letters $a,b,c$ to denote individual objects as e.g.\ Socrates, New York etc. (of which the predicates $M,P\ldots$ are asserted). Those objects are usually called \emph{individuals} in mathematical logic.

Now let $M$ be a monadic predicate, e.g.\ ``mor\-tal'', and $a$ an individual, e.g.\ Socrates. Then the proposi\-tion that $M$ belongs to $a$ is denoted by $M(a)$. So $M(a)$ means ``Socrates is mortal'' and similarly if $G$ is a dyadic relation, e.g.\ larger, and $b,c$ two individuals, e.g.\ New York and Washington, then $G(b,c)$ means ``The relation $G$ subsists between $b$ and $c$'', i.e.\ in our case ``New York is larger than Washington''. So in this notation there is no copula, but e.g.\ the proposition ``Socrates is mortal'' $\mathbf{\llbracket 31. \rrbracket}$ has to be expressed like this Mortality(Socrates), and that New York is greater than Washington by Larger(New York, Washington).

That much I have to say about the simplest type of propositions which say that some definite predicate belongs to some definite subject or subjects. These propositions are sometimes called atomic propositions\index{atomic formulas in predicate logic}\index{atomic propositions in predicate logic} because they constitute so to speak the atoms of which the more complex propositions are built up. But now how are they built up? We know already one way of forming $\mathbf{\llbracket 32. \rrbracket}$ compound propositions namely by means of the operations of the propositional calculus $.\, , \vee, \supset$ etc., e.g.\ from the two atomic propositions ``Socrates is a man'' and ``Socrates is mortal'' we can form the composit proposition ``If Socrates is a man Socrates is mortal''; in symbols, if $T$ denotes the predicate of mortality it would read $M(a)\supset T(a)$, or e.g.\ $M(a)\:\vee\sim M(a)$ would mean ``Either Socrates is a man or Socrates is not a man''. $M(a)\: .\: T(a)$ would mean ``Socrates is a man and Socrates is mortal'', and so on. The propositions which we can obtain in this way, i.e.\ by combining atomic propositions by means $\mathbf{\llbracket 31.a \rrbracket}$ of the truth functions $\vee, .$ etc. are sometimes called molecular propositions.\index{molecular propositions}

But there is still another way of forming compound propositions which we have not yet taken account of in our symbolism, namely by means of the particles ``every'' and ``some''.\index{universal quantifier} These are expressed in logistics by the use of variables as follows: Take e.g.\ the proposition ``Every man is mortal''. We can express that in other words like this: ``Every object which is a man is mortal'' or ``For every object $x$ it is true that $M(x)\supset T(x)$''. Now in order to indicate that this implication $\mathbf{\llbracket 32.a \rrbracket}$ is asserted of any object $x$ one puts $x$ in brackets in front of the proposition and includes the whole proposition in brackets to indicate that the whole proposition is asserted to be true for every $x$. And generally if we have an arbitrary expression, say $\Phi(x)$ which involves a variable $x$, then $(x)[\Phi(x)]$ means ``For every object $x$, $\Phi(x)$ is true'', i.e.\ if you take an arbitrary individual $a$ and substitute it for $x$ then the resulting proposition $\Phi(a)$ is true. As in our example $(x)[M(x)\supset T(x)]$, $\mathbf{\llbracket 33. \rrbracket}$ if you substitute Socrates for $x$ you get the true proposition. And generally if you substitute for $x$ something which is a man you get a true proposition because then the first and second term of the implication are true. If however you substitute something which is not a man you also get a true proposition because\ldots\ So for any arbitrary object which you substitute for $x$ you get a true proposition and this is indicated by writing $(x)$ in front of the proposition. $(x)$ is called the universal quantifier.

$\mathbf{\llbracket 34. \rrbracket}$ As to the particle ``some'' or ``there exists''\index{existential quantifier} it is expressed by a reversed $\exists$ put in brackets together with a variable $(\exists x)$. So that means: there is an object $x$; e.g.\ if we want to express that some men are not mortal we have to write $(\exists x)[M(x)\: .\: \sim T(x)]$ and generally if $\Phi(x)$ is a propositional function\index{propositional function} with the variable $x$, $(\exists x)[\Phi(x)]$ means $\mathbf{\llbracket 35. \rrbracket}$ ``There exits some object $a$ such that $\Phi(a)$ is true''. Nothing is said about the number of objects for which $\Phi(a)$ is true; there may be one or several. $(\exists x)\Phi(x)$ only means there is at least one object $x$ such that $\Phi(x)$. $(\exists x)$ is called the existential quantifier. From this definition you see at once that we have the following equivalences:
\begin{tabbing}
\hspace{10em}$(\exists x)\Phi(x)\;$\=$\equiv\; \sim(x)[\sim\Phi(x)]$\\[.5ex]
\hspace{10em}$(x)\Phi(x)$\>$\equiv\; \sim(\exists x)[\sim\Phi(x)]$
\end{tabbing}

Generally $(x)[\sim\Phi(x)]$ means $\Phi(x)$ holds for no object and $\sim (\exists x)[\Phi(x)]$ means there is no object $x$ such that $\Phi(x)$. Again you see that these two statements are equivalent with each other. It is easy e.g.\ to express the traditional four $\mathbf{\llbracket 36. \rrbracket}$ types of propositions a, e, i, o in our notation.\index{a, e, i, o propositions} In each case we have two predicates, say $P$, $S$ and
\begin{tabbing}
\hspace{4em}\=$S$a$P$\index{a proposition} \hspace{.3em}\=means\hspace{.4em} \=every $S$ is a $P$ \hspace{1.3em}\=i.e.\ \hspace{.2em}\=$(x)[S(x)\supset P(x)]$\\[.5ex]
\>$S$i$P$\index{i proposition} \>means \>some $S$ are $P$ \>i.e.\ \>$(\exists x)[S(x)\: .\: P(x)]$\\[.5ex]
\>$S$e$P$\index{e proposition} \>means \>no $S$ is a $P$ \>i.e.\ \>$(x)[S(x)\supset\;\sim P(x)]$\\[.5ex]
\>$S$o$P$\index{o proposition} \>means \>some $S$ are $\sim P$ \>i.e.\ \>$(\exists x)[S(x)\: .\: \sim P(x)]$
\end{tabbing}
You see the universal propositions have the universal quantifier in front of them and the particular propositions the existential quantifier. I want to mention that in classical logic two entirely different types of proposi\-tions are counted as universal affirmative, namely propositions of the type ``Socrates is mortal'' expressed by $P(a)$ and ``Every man is mortal'' $(x)[S(x)\supset P(x)]$.

$\mathbf{\llbracket 37. \rrbracket}$ Now the existential and universal quantifier can be combined with each other and with the truth functions $\sim$,$\,$\ldots\ in many ways so as to express more complicated propositions.

$\mathbf{\llbracket 37.1 \rrbracket}$\footnote{This page followed by the \textbf{new page} below is inserted within p.\ \textbf{37}, which continues with the paragraph after the next starting with ``I want now to give''.} Thereby one uses some abbreviations, namely: Let $\Phi(xy)$ be an expression containing two variables; then we may form: $(x)[(y)[\Phi(xy)]]$. That means ``For any object $x$ it is true that for any object $y$ $\Phi(xy)$'' that evidently means ``$\Phi(xy)$ is true whatever objects you take for $x,y$'' and this is denoted by $(x,y)\Phi(xy)$. Evidently the order of the variables is arbitrary here, i.e.\ $(x, y)\Phi(xy)\equiv (y, x)\Phi(xy)$. Similarly $(\exists x)[(\exists y)[\Phi(xy)]]$ means ``There are some objects $x,y$ such that $\Phi(xy)$'' and this is abbreviated by $(\exists x,y)\Phi(xy)$ and means:\ldots\ But it has to be noted that this does not mean that there are really two different objects $x,y$ satisfying $\Phi(xy)$. This formula is also true if there is one object $a$ such that $\Phi(aa)$ because then there exists an $x$, namely $a$, such that there exists a $y$, namely again $a$, such that etc. Again $(\exists x,y)\Phi(xy)\equiv (\exists y,x)\Phi(xy)$.

But it is to be noted that this interchangeability holds $\mathbf{\llbracket new\, page \rrbracket}$ only for two universal or two existential quantifiers. It does not hold for an universal and an existential quantifier, i.e.\ $(x)[(\exists y)[\Phi(yx)]]\not\equiv (\exists y)[(x)[\Phi(yx)]]$. Take e.g.\ for $\Phi(yx)$ the proposition ``$y$ greater than $x$''; then the first means ``For any object $x$ it is true that there exists an object $y$ greater than $x$''; in other words ``For any object there exists something greater''. The right-hand side however means ``There exists an object $y$ such that for any $x$ $y$ is greater than $x$'', there exists a greatest object. So that in our case the right side says just the opposite of what the left side says. The above abbreviation is also used for more than two variables, i.e.\ $(x,y,z)[\Phi(xyz)]$ $(\exists x,y,z)[\Phi(xyz)]$.

I want now to give some examples for the notation introduced. Take e.g.\ the proposition ``For any integer there exists a greater one''. The predicates occurring in this proposition are: 1.\ integer and 2.\ greater. Let us denote them by $I$ and $>$ so $I(x)$ is to be read ``$x$ is an integer'' and $>\!(xy)$ is to be read ``$x$ greater $y$'' or ``$y$ smaller $x$''. Then the proposition is expressed in logistic symbolism as follows:
\begin{tabbing}
\hspace{1.7em}$(x)[I(x)\supset (\exists y)[I(y)\: .\: >\!(yx)]]$. \end{tabbing}
We can express the same fact by saying $\mathbf{\llbracket 38. \rrbracket}$ there is no greatest integer. What would that look like in logistic symbolism: \begin{tabbing}
\hspace{1.7em}$\sim\!(\exists x)[I(x)\: .\:$ such that no integer is greater i.e.\
$(y)[I(y)\supset \;\sim \;>\!(yx)]$].
\end{tabbing}
As another example take the proposition ``There is a smallest integer'' that would read:
\begin{tabbing}
\hspace{1.7em}$(\exists x)[I(x)\: .\:$ such that no integer is smaller i.e.\
$(y)[I(y)\supset \;\sim \; >\!(xy) ]]$.
\end{tabbing}
I wish to call your attention to a near at hand mistake. It would be wrong to express this last proposition like this:
\begin{tabbing}
\hspace{1.7em}$(\exists x)[I(x)\: .\:(y)[I(y)\supset\; >\!(yx)]]$
\end{tabbing}
because that would mean there is an integer smaller than every integer. But such an integer does not exist $\mathbf{\llbracket 39. \rrbracket}$ since it would have to be smaller than itself. An integer smaller than every integer would have to be smaller than itself---that is clear. So the second proposition is false whereas the first is true, because it says only there exists an integer $x$ which is not greater than any integer.

Another example for our notation may be taken from Geometry. Consider the proposition ``Through any two different points there is exactly one straight line''. The predicates which occur in this proposition are 1.\ point $P(x)$, $\mathbf{\llbracket 40. \rrbracket}$ 2.\ straight line $L(x)$, 3.\ \emph{different} that is the negation of identity. Identity\index{identity} is denoted by $=$ and difference by $\neq$. =$\,(xy)$ means $x$ and $y$ are the same thing, e.g.\ =$\,$(Shakespeare, author of Hamlet), and $\neq(xy)$ means $x$ and $y$ are different from each other. There is still another relation that occurs in our geometric proposition, namely the one expressed by the word ``through''. That is the relation which holds between a point $x$ and a line $y$ if ``$y$ passes through $x$'' or in other words if ``$x$ lies on $y$''. Let us denote that relation by $J(xy)$. Then the geometric proposition mentioned, in order to be expressed in logistic symbolism, has to be splitted into two parts, namely there is at least one line and there is at most one line. The first reads: $(x,y)[P(x)\: .\: P(y)\: .\: \neq(xy)\supset$ $\mathbf{\llbracket 41. \rrbracket}$ $(\exists u)[L(u)\: .\: J(xu)\: .\: J(yu)]]$. So that means that through any two different points there is\ldots\ But it is not excluded by that statement that there are two or three different lines passing through two points. That there are no two different lines could be expressed like this
\begin{tabbing}
\hspace{1.7em}$(x,y)[P(x)\: .\: P(y)\: .\: \neq(xy)\supset\;\sim(\exists u,v)[L(u)\: .\: L(v)\: .\: \neq(uv)\: .\:$ \\
\` $J(xu)\: .\: J(yu)\: .\: J(xv)\: .\: J(yv)]]$
\end{tabbing}

I hope these examples will suffice to make clear how the quantifiers are to be used. For any quantifier occurring in an expression there is a definite portion of the expression to which it relates (called the scope of the expression), e.g.\ scope of $x$ whole expression, of $y$ only this portion\ldots\ So the scope\index{scope}\index{scope of a quantifier} is the proposition of which it is asserted that it holds for all or every object. It is indicated by the brackets which begin immediately behind the quantifier. There are some conventions about leaving out these brackets, namely they may be left out 1.\ if the scope is atomic, e.g.\ $(x)\varphi(x)\supset p\;$: $(x)[\varphi(x)]\supset p$, not $(x)[\varphi(x)\supset p]$, 2.\ if the scope begins with $\sim$ or a quantifier, e.g.\
\begin{tabbing}
\hspace{1.7em}\=$(x)\sim[\varphi(x)\: .\: \psi(x)]\vee p\;$ \=: $\; (x)[\sim[\varphi(x)\: .\: \psi(x)]]\vee p$\\[.5ex]
\>$(x)(\exists y)\varphi(x)\vee p$\>: $\; (x)[(\exists y)[\varphi(x)]]\vee p$\end{tabbing}
But these rules are only facultative, i.e.\ we may also write all the brackets if it is expedient for the sake of clarity.

A variable to which a quantifier $(x)$, $(y)$, $(\exists x)$, $(\exists y)$ refers is called a ``bound variable''.\index{bound variable}\index{variable, bound} In the examples which I gave, all variables $\mathbf{\llbracket 42. \rrbracket}$ are bound (e.g.\ to this $x$ relates this quantifier etc.) and similarly to any variable occurring in those expressions you can associate a quantifier which refers to it. If however you take e.g.\ the expression: $I(y)\: .\: (\exists x)[I(x)\: .\: >(yx)]$, which means: there is an integer $x$ smaller than $y$, then here $x$ is a bound variable because the quantifier $(\exists x)$ refers to it. But $y$ is not bound because the expression contains no quantifier referring to it. Therefore $y$ is called a free variable\index{free variable}\index{variable, free} of this expression. An expression containing free variables is not a proposition, but it only becomes a proposition if the free variables are replaced by individual objects, e.g.\ this expression here means $\mathbf{\llbracket 43. \rrbracket}$ ``There is an integer smaller than the integer $y$''. That evidently is not a definite assertion which is either true or wrong. But if you substitute for the free variable $y$ a definite object, e.g.\ 7, then you obtain a definite proposition, namely: ``There is an integer smaller than~7''.

The bound variables have the property that it is entirely irrelevant by which letters they are denoted; e.g.\ $(x)(\exists y)[\Phi(xy)]$ means exactly the same thing as $(u)(\exists v)[\Phi(uv)]$. The only requirement is that you must use different letters for different bound variables. But even that is only necessary for variables $\mathbf{\llbracket 44. \rrbracket}$ one of whom is contained in the scope of the other as e.g.\ in $(x)[(\exists y)\Phi(xy)]$, where $y$ is in the scope of $x$ which is the whole expression, and therefor it has to be denoted by a letter different from $x$; $(x)[(\exists x)\Phi(xx)]$ would be ambiguous. Bound variables whose scopes lie outside of each other however can be denoted by the same letter without any ambiguity, e.g.\ $(x)\varphi(x)\supset(x)\psi(x)$. For the sake of clarity we also require that the free variables in a propositional function should always be denoted by letters different from the bound variables; so e.g.\ $\varphi(x)\: .\: (x)\psi(x)$ is not a correctly formed propositional function, but $\varphi(x)\: .\: (y)\psi(y)$ is one.

The examples of formulas which I gave last time and also the problems to be solved were propositions concerning certain definite predicates $I$, $<$, =, etc. They are true only for those particular predicates occurring in them. But now exactly as in the calculus of propositions there are certain formulas which are true whatever propositions the letters $p,q,r$ may be so also in the calculus of predicates $\mathbf{\llbracket 45. \rrbracket}$ there will be certain formulas which are true for any arbitrary predicates. I denote arbitrary predicates by small Greek letters $\varphi,\psi$. So these are supposed to be variables for predicates\index{predicate variables}\index{variables, predicate} exactly as $p,q\ldots$ are variables for propositions and $x,y,z$ are variables for objects.

Now take e.g.\ the proposition $(x)\varphi(x)\vee(\exists x)\sim\varphi(x)$, i.e.\ ``Either every individual has the property $\varphi$ or there is an individual which has not the property $\varphi$''. That will be true for any arbitrary monadic predicate $\varphi$. We had other examples before, e.g.\ ${(x)\varphi(x)\equiv\:\sim(\exists x)\sim\varphi(x)}$ that again is true for any arbitrary monadic predicate $\varphi$. Now exactly as in the calculus of propositions such expressions which are true for all predicates are called tautologies or universally true.\index{valid formula}\index{tautology, of predicate logic (valid formula)}\index{universally true formula (valid formula)} Among them are e.g.\ all the formulas which express the Aristotelian $\mathbf{\llbracket 46. \rrbracket}$ moods of inference, e.g.\ the mood Barbara\index{Barbara} is expressed like this:
\begin{tabbing}
\hspace{1.7em}$(x)[\varphi(x)\supset\psi(x)]\: .\: (x)[\psi(x)\supset\chi(x)]\supset (x)[\varphi(x)\supset\chi(x)]$
\end{tabbing}
The mood Darii\index{Darii} like this
\begin{tabbing}
\hspace{11em}\=$\varphi$\hspace{1em}\= $M$a$P$\hspace{1em}\=$\psi$\\[.5ex]
\>$\chi$\>\uline{$S$i$M$}\>$\varphi$\\[.5ex]
\>\>$S$i$P$\\[1ex]
\hspace{1.7em}$(x)[\varphi(x)\supset\psi(x)]\: .\: (\exists x)[\chi(x)\: .\:\varphi(x)]\supset (\exists x)[\chi(x)\: .\:\psi(x)]$
\end{tabbing}

\subsection{Decidability and completeness in predicate logic}\label{2.2}
\pagestyle{myheadings}\markboth{EDITED TEXT}{NOTEBOOK IV \;\;---\;\; 1.2.2\;\; Decidability and completeness in predicate\ldots}

It is of course the chief aim of logic to investigate the tautologies and exactly as in the calculus of propositions there are again two chief problems which arise. Namely: 1.\ To develop methods for finding out about a given expression whether or not it is a tautology, 2.\ To reduce all tautologies to a finite number of logical axioms and rules of inference from which they can be derived. I wish to mention right now that only $\mathbf{\llbracket 47. \rrbracket}$ the second problem can be solved for the calculus of predicates. One has actually succeeded in setting up a system of axioms for it and in proving its completeness\index{completeness in predicate logic} (i.e.\ that every tautology can be derived from it).

As to the first problem, the so called decision problem,\index{decision problem in predicate logic} it has also been solved in a sense but in the negative, i.e.\ one has succeeded in proving that there does not exist any mechanical procedure to decide of any given expression whether or not it is a tautology of the calculus of predicates. That does not mean that there are any individual formulas of which one could not decide whether or not they are $\mathbf{\llbracket 48. \rrbracket}$ tautologies. It only means that it is not possible to decide that by a purely mechanical procedure. For the calculus of propositions this was possible, e.g.\ the truth-table method is a purely mechanical procedure which allows to decide of any given expression whether or not it is a tautology. So what has been proved is only that a similar thing cannot exist for the calculus of predicates. However for certain special kinds of formulas such methods of decision have been developed, e.g.\ for all formulas with only monadic predicates (i.e.\ formulas without relations in it); $\mathbf{\llbracket 49. \rrbracket}$ e.g.\ all formulas expressing the Aristotelian moods are of this type because no relations occur in the Aristotelian moods.

Before going into more detail about that I must say a few more words about the notion of a tautology of the calculus of predicates.

There are also tautologies which involve variables both for propositions and for predicates, e.g.\
\begin{tabbing}
\hspace{1.7em}$p\: .\:(x)\varphi(x)\equiv (x)[p\: .\: \varphi(x)]$\index{laws of passage}
\end{tabbing}
i.e.\ if $p$ is an arbitrary proposition and $\varphi$ an arbitrary predicate then the assertion on the left, i.e.\ ``$p$ is true and for every $x$, $\varphi(x)$ is true'' is equivalent with the assertion on the right, i.e.\ ``for every object $\mathbf{\llbracket 50. \rrbracket}$ $x$ the conjunction $p\: .\: \varphi(x)$ is true''. Let us prove that, i.e.\ let us prove that the left side implies the right side and vice versa the right side implies the left side. If the left side is true that means: $p$ is true and for every $x$, $\varphi(x)$ is true, but then the right side is also true because then for every $x$, $p\: .\: \varphi(x)$ is evidently true. But also vice versa: If for every $x$, $p\: .\: \varphi(x)$ is true then 1. $p$ must be true because otherwise $p\: .\: \varphi(x)$ would be true for no $x$ and 2. $\varphi(x)$ must be true for every $x$ since by assumption even $p\: .\: \varphi(x)$ is true for every $x$. So you see this equivalence holds for any predicate $\varphi$, $\mathbf{\llbracket 51. \rrbracket}$ i.e.\ it is a tautology.

There are four analogous tautologies obtained by replacing $.$ by $\vee$ and the universal quantifier by the existential quantifier, namely
\begin{tabbing}
\hspace{1.7em}\=2.\quad\=$p\vee (x)\varphi(x)\equiv (x)[p\vee\varphi(x)]$\index{laws of passage}\\[.5ex]
\>3.\>$p\: .\:(\exists x)\varphi(x)\equiv (\exists x)[p\: .\: \varphi(x)]$\index{laws of passage}\\[.5ex]
\>4.\>$p\vee(\exists x)\varphi(x)\equiv (\exists x)[p\vee\varphi(x)]$\index{laws of passage}
\end{tabbing}
I shall give the proof for them later on. These four formulas are of a great importance because they allow to shift a quantifier over a symbol of conjunction or disjunction. If you write $\sim p$ instead of $p$ in the first you get $[p\supset(x)\varphi(x)]\equiv (x)[p\supset\varphi(x)]$. This law of logic is used particularly frequently in proofs as you will see later. Other examples of tautologies are e.g.\
\begin{tabbing}
\hspace{1.7em}\=$(x)\varphi(x)\: .\: (x)\psi(x)\equiv (x)[\varphi(x)\: .\:\psi(x)]$\\[.5ex]
\>$(\exists x)\varphi(x)\vee (\exists x)\psi(x)\equiv (\exists x)[\varphi(x)\vee\psi(x)]$
\end{tabbing}
or e.g.\
\begin{tabbing}
\hspace{1.7em}$\sim(x)(\exists y)\varphi(xy)\equiv (\exists x)(y)\sim\varphi(xy)$
\end{tabbing}
$\mathbf{\llbracket 52. \rrbracket}$ That means:

Proof. $\sim(x)(\exists y)\varphi(xy)\,{\equiv \atop \mbox{\scriptsize\rm means}}\,(\exists x)\sim(\exists y)\varphi(xy)$, but $\sim(\exists y)\varphi(xy)\equiv(y)\sim\varphi(xy)$ as we saw before. Hence the whole expression is equivalent with $\equiv(\exists x)(y)\sim\varphi(xy)$ which was to be proved.

\vspace{1ex}

Another example: $(x)\varphi(x)\supset(\exists x)\varphi(x)$, i.e.\ If every individual has the property $\varphi$ then a fortiori there are individuals which have the property $\varphi$. The inverse of this proposition is no tautology, i.e.\ \begin{tabbing}
\hspace{1.7em}$(\exists x)\varphi(x)\supset(x)\varphi(x)$ is not a tautology
\end{tabbing}
because if there is an object $x$ which has the property $\varphi$ that does not imply that every individual has the property $\varphi$.

But here there is an important remark $\mathbf{\llbracket 53. \rrbracket}$ to be made. Namely: In order to prove that this formula here is not a tautology we must know that there exists more than one object in the world. For if we assume that there exists only one object in the world then this formula would be true for every predicate $\varphi$, hence would be universally true because if there is only one object, say $a$, in the world then if there is an object $x$ for which $\varphi(x)$ is true this object must be $a$ (since by assumption there is no other object), hence $\varphi(a)$ is true; but then $\varphi$ is true for every object because by assumption there exists only this object $a$. I.e.\ in a world with only one $\mathbf{\llbracket 54. \rrbracket}$ object $(\exists x)\varphi(x)\supset(x)\varphi(x)$ is a tautology. It is easy to find some expressions which are universally true if there are only two individuals in the world etc., e.g.\
\begin{tabbing}
\hspace{1.7em}$(\exists x,y)[\psi(x)\: .\:\psi(y)\: .\:\varphi(x)\: .\:\sim\varphi(y)]\supset (x)[\psi(x)]$
\end{tabbing}

At present I only wanted to point out that the notion of a tautology of the calculus of predicates\index{valid formula}\index{tautology, of predicate logic (valid formula)} needs a further specification in order to be precise. This specification consists in this that an expression is called a tautology only if it is universally true no matter how many individuals are in the world assuming only that there is at least one (otherwise the meaning of the quantifiers is not definite). So e.g.\ $(x)\varphi(x)\supset(\exists y)\varphi(y)$; this is a tautology because it is true\ldots\ but this inverse is not because\ldots\ It can be proved that this means the same thing as if I said: An expression is a tautology if it is true in a world with infinitely many individuals, i.e.\ one can prove that whenever an expression is universally true in a world $\mathbf{\llbracket Notebook\; V \rrbracket}$ $\mathbf{\llbracket 55. \rrbracket}$
\pagestyle{myheadings}\markboth{EDITED TEXT}{NOTEBOOK V \;\;---\;\; 1.2.2\;\; Decidability and completeness in predicate\ldots}
with infinitely many objects it is true in any world no matter how many individuals there may be and of course also vice versa. I shall not prove this equivalence but shall stick to the first definition.

The formulas by which we expressed the tautologies contain free variables (not for individuals) but for predicates and for propositions, e.g.\ $\varphi$ here is a free variable in this expression (no quantifier related to it, i.e.\ no $(\varphi)$ $(\exists \varphi)$ occurs); similarly here, so these formulas are really propositional functions since they contain free variables.

And the definition of a tautology\index{valid formula}\index{tautology, of predicate logic (valid formula)} was that whatever particular proposition or predicate you substitute for those free variables of predicates or propositions you get a true proposition. The variables for individuals were all bound.

We can extend the notion of a $\mathbf{\llbracket 56. \rrbracket}$ tautology\index{valid formula}\index{tautology, of predicate logic (valid formula)} also to such expressions as contain free variables for individuals, e.g.\
\begin{tabbing}
\hspace{1.7em} $\varphi(x) \:\vee \sim \varphi(x)$
\end{tabbing}
This is a propositional function containing one free functional variable and one free individual variable $x$ and whatever object and predicate you substitute for $\varphi, x$ you get a true proposition. Formula
\begin{tabbing}
\hspace{1.7em}$(x)\varphi(x) \supset \varphi(y)$
\end{tabbing}
contains $\varphi ,y$ and is universally true because if $M$ is an arbitrary predicate and $a$ an arbitrary individual then
\begin{tabbing}
\hspace{1.7em}$(x)M(x)\supset M(a)$
\end{tabbing}
So in general a tautological logical formula\index{valid formula}\index{tautology, of predicate logic (valid formula)} of the calculus of functions is a propositional function composed of the above mentioned symbols and which is true whatever particular $\mathbf{\llbracket 57. \rrbracket}$ objects and predicates and propositions you substitute for free variables no matter how many individuals there exist. We can of course express this fact, namely that a certain formula is a universally true, by writing quantifiers in front, e.g.\
\begin{tabbing}
\hspace{1.7em}$(\varphi ,x)[\varphi(x) \:\vee \sim \varphi(x)]$
\end{tabbing}
or
\begin{tabbing}
\hspace{1.7em}$(\varphi, y)[(x)\varphi(x) \supset \varphi(y)]$
\end{tabbing}
For the tautology of the calculus of propositions
\begin{tabbing}
\hspace{1.7em}$(p,q)[p \supset p \vee q]$
\end{tabbing}
But it is more convenient to make the convention that universal quantifiers whose scope is the whole expression may be left out. So if a formula containing free variables is written down as an assertion, e.g.\ as an axiom or theorem, it means that it holds for everything substituted for the free variables, i.e.\ it means the same thing as if all variables were bound by quantifiers whose scope is the whole expression. This convention is in agreement with the way in which theorems are expressed in mathematics, e.g.\ the law of raising a sum to the square is written $(x+y)^2=x^2+2xy+y^2$, i.e.\ with free variables $x,y$ which express that this holds for any numbers. $\mathbf{\llbracket 57.1 \rrbracket}$ It is also in agreement with our use of variables for propositions in the calculus of propositions. The axioms and theorems of the propositional calculus were written with free variables, e.g.\ $p \supset p \vee q$, and a formula like this was understood to mean that it holds for any propositions $p,q$.

\subsection{Axiom system for predicate logic}\label{2.3}
\pagestyle{myheadings}\markboth{EDITED TEXT}{NOTEBOOK V \;\;---\;\; 1.2.3\;\; Axiom system for predicate logic}

$\mathbf{\llbracket 58. \rrbracket}$ I hope that these examples will be sufficient and that I can now begin with setting up the axiomatic system for the calculus of predicates which allows to derive all tautologies of the calculus of predicates. The primitive notions will be $1.$ the former $\sim ,\vee$ $2.$ the universal quantifier $(x), (y)$. The existential quantifier need not be taken as a primitive notion because it can be defined in terms of $\sim$ and $(x)$ by ${(\exists x)\varphi(x)\equiv\: \sim (x) \sim \varphi (x)}$. The formulas of the calculus of predicates\index{formula, of predicate logic} will be composed of three kinds of letters: $p,q,\ldots$ propositional variables, $\varphi, \psi,\ldots$ functional variables\index{functional variables}\index{variables, functional}\index{predicate variables}\index{variables, predicate} for predicates, $x,y,\ldots$ variables for individuals.\index{individual variables}\index{variables, invidual} Furthermore they will contain $\mathbf{\llbracket 59. \rrbracket}$ $(x), (y), \sim , \vee$ and the notions defined by those three, i.e.\ $(\exists x), (\exists y), \supset, \:.\:, \equiv, \mid$ etc. So the quantifiers apply only to individual variables, propositional and functional variables are free, i.e.\ that something holds for all $p,\varphi$ is to be expressed by free variables according to the convention mentioned before.

So all formulas given as examples before are examples for expressions of the calculus of functions but also e.g.\ $(\exists x)\psi(xy)$ and $[p\:.\:(\exists x)\psi(xy)]\vee \varphi(y)$ would be examples etc. I am using the letters $\Phi,\Psi,\Pi$ to denote arbitrary expressions of the calculus of predicates and if I wish to indicate that some variable say $x$ occurs in a formula as a free variable denote the formula by $\Phi(x) \vee \Psi(xy)$ if $x,y$ occur both free, which does not exclude that there may be other free variables besides $x$, or $x$ and $y$, in the formula.

The axioms\index{axiom system for predicate logic}\index{five axioms for predicate logic}\index{axioms for predicate logic} are like this:
\begin{tabbing}
\hspace{1.7em}\=I.\hspace{1em}\= The four axioms of the calculus of propositions\index{four axioms for propositional logic}\index{axioms for propositional logic}
\\[1ex]
\>\>\hspace{2em}\=$p \supset p \vee q$\hspace{2em}\=$p \vee q \supset q \vee p$\\[.5ex]
\>\>\>$p \vee p \supset p$\>$(p \supset q)\supset (r \vee p \supset r\vee q)$\\[1ex]
\>II. \>One specific axiom for the universal quantifier \\[1ex]
\>\>\>Ax.\ 5 \qquad $(x)\varphi(x)\supset \varphi(y)$\index{Ax.\ 5, axiom for predicate logic}\index{axiom Ax.\ 5 for predicate logic}
\end{tabbing}
This is the formula mentioned before which says: ``For any $y$, $\varphi$ it is true that if $\varphi$ holds for every $x$ then it holds for $y$''.

These are all axioms which we need. The rules of inference are the following four:\index{rules of inference for predicate logic}\index{four rules of inference for predicate logic}
\begin{itemize}
	\item[$\mathbf{\llbracket 60. \rrbracket}$]
\item[1\;\;] The rule of implication\index{rule of implication}\index{modus ponens}\index{rule 1} which reads exactly as for the calculus of propositions: If $\Phi,\Psi$ are any expressions then from $\Phi,\Phi \supset \Psi$ you can conclude $\Psi$.
\end{itemize}
The only difference is that now $\Phi,\Psi$ are expressions which may involve quantifiers and functional variables and individual variables in addition to the symbols occurring in the calculus of propositions. So e.g.\
\begin{tabbing}
 \hspace{1.7em}\= from $[p \vee (x)[\varphi(x)\supset \varphi(x)]]\supset \varphi(y) \: \vee \sim \varphi (y)$\\[.5ex]
\>\underline{and $[p \vee (x)[\varphi (x)\supset \varphi (x)]]$}\\[.5ex]
\>conclude $\varphi (y)\: \vee \sim \varphi (y)$
\end{tabbing}
\begin{itemize}
\item[2\;\;] The rule of Substitution\index{rule of substitution for predicate logic}\index{rule 2} which has now three parts (according to the three kinds of variables):
\begin{itemize}
\item[a)] For individual variables $x,y$ bound or free any other individual variable may be substituted as long as our conventions about the notion of free variables are observed, i.e.\ bound variable whose scopes do not lie outside of each other must be denoted by different letters and all free variables must be denoted by letters different from all bound variables -- [Rule of renaming the individual variables].
\end{itemize}

\vspace{-2ex}

\item[$\mathbf{\llbracket 61. \rrbracket}$]
\begin{itemize}
\item[b)] For a propositional variable any expression may be substituted with a certain restriction formulated later.
\item[c)] If you have an expression $\Pi$ and $\varphi$ a functional variable occurring in $\Pi$ perhaps on several places and with different arguments $\varphi (x)$, $\varphi(y)$,\ldots\ and if $\Phi(x)$ is an expression containing $x$ free then you may substitute $\Phi(x)$ for $\varphi(x)$, $\Phi(y)$ for $\varphi(y)$ etc. simultaneously in all places where $\varphi$ occurs. Similarly for $\varphi(xy)$ and $\Phi(xy)$.
\end{itemize}
\end{itemize}

$\mathbf{\llbracket 61.1 \rrbracket}$ It is clear that this is a correct inference, i.e.\ gives a tautology if the formula in which we substitute is a tautology, because if a formula is a tautology that means that it holds for any property or relation $\varphi, \psi$, but any propositional function with one or several free variables defines a certain property or relation; therefore the formula must hold for them. Take e.g.\ the tautology $(x)\varphi(x)\supset \varphi(y)$ and substitute for $\varphi$ the expression $(\exists z)\psi(zx)$ which has one free individual variable. Now the last formula says that for every property $\varphi$ and any individual $y$ we have: ``If for any $x$ $\varphi(x)$ then $\varphi(y)$''. But if $\psi$ is an arbitrary relation then $(\exists z)\psi(zx)$ defines a certain property because it is a propositional function with one free variable $x$. Hence the above formula must hold also for this property, i.e.\ we have: If for every object $(x)[(\exists z)\psi(zx)]$ then also for $y$ $(\exists z)\psi(zy)$ and that will be true whatever the relation $\psi$ and the object $y$ may be, i.e.\ it is again a tautology.

$\mathbf{\llbracket 62. \rrbracket}$ You see in this process of substitution we have sometimes to change the free variables, as here we have to change $x$ into $y$ because the $\varphi$ occurs with the variable $y$ here; if the $\varphi$ occurred with the variable $u$ $\varphi(u)$ we would have to substitute $(\exists z)\psi(zu)$ in this place. In this example we substituted an expression containing $x$ as the only free variable, but we can substitute for $\varphi(x)$ here also an expression which contains other free individual variables besides $x$, i.e.\ also in this case we shall obtain a tautology. Take e.g.\ the expression $(\exists z)\chi(zxu)$. This is a propositional function with the free individual variable $x$ but it has the free individual variable $u$ in addition. Now if we replace $\chi$ by a special triadic relation $R$ and $u$ by a special object $a$ then $(\exists z)R(zxa)$ is a propositional function with one free variable $x$; hence $\mathbf{\llbracket 63.1\rrbracket}$ it defines a certain property, hence the above formula holds, i.e.\
\begin{tabbing}
\hspace{1.7em}$(x)(\exists z)R(zxa)\supset (\exists z)R(zya)$
\end{tabbing}
whatever $y$ may be, but this will be true whatever $R,a$ may be; therefore if we replace them by variables $\chi$, $u$ the formula obtained:
\begin{tabbing}
\hspace{1.7em}$(x)(\exists z)\chi(zxu)\supset (\exists z)\chi(zyu)$
\end{tabbing}
will be true for any $\chi,u,y$, i.e.\ it is a tautology. So the rule of substitution is also correct for expressions containing additional free variables $u$, and therefore this $\Phi(x)$ is to mean an expression containing the free variable $x$ but perhaps some other free variables in addition.

$\mathbf{\llbracket 64.\rrbracket}$ Examples for the other two rules of substitution:
\begin{tabbing}
\hspace{1.7em} \= For propositional variable\\[1ex]
\>$p\:.\:(x)\varphi(x) \equiv (x)[p\:.\:\varphi(x)]$
\end{tabbing}
substitute $(\exists z)\psi(z)$. Since this holds for every proposition it holds also for $(\exists z)\psi(z)$ which is a proposition if $\psi$ is any arbitrary predicate. Hence we have for any predicates $\psi,\varphi$
\begin{tabbing}
\hspace{1.7em}$(\exists z)\psi(z)\:.\:(x)\varphi(x)\equiv (x)[(\exists z)\psi(z)\:.\: \varphi(x)]$
\end{tabbing}
But we are also allowed to substitute expressions containing free variables and propositional variables e.g.\ $(z)\chi(zu)$ (free variable $u$) because if you take for $u$ any individual object $a$ [and $p$ any individual proposition $\pi$] and $\chi$ any relation $R$ then $\mathbf{\llbracket 65.\rrbracket}$ this will be a proposition. And $p\:.\:(x)\varphi(x) \equiv (x)[p\:.\:\varphi(x)]$ holds for any proposition. So it will also hold for this, i.e.\
\begin{tabbing}
\hspace{1.7em}$[(z)\chi(zu)]\:.\:(x)\varphi(x) \equiv (x)[(z)\chi(zu)\:.\:\varphi(x)]$
\end{tabbing}
will be true whatever $p,\chi,\varphi,u$ may be, i.e.\ a tautology.

Finally an example for substitution of individual variables:
\begin{itemize}
\item[] For a bound $(x)\varphi(x)\supset \varphi(y)$\; :\; $(z)\varphi(z)\supset \varphi(y)$. So this inference merely brings out the fact that the notation of bound variables is arbitrary.
\item[] The rule of substitution applied for free variables is more essential; e.g.\ from $(x,y)\varphi(xy)\supset \varphi(uv)$ we can conclude $(x,y)\varphi(xy)\supset \varphi(uu)$ by substituting $u$ for $v$. This is an allowable substitution because the variable which you substitute, $u$, does not occur as a bound variable. It occurs as a free variable but that does not matter.
\end{itemize}

Of course if a variable occurs in several places it has to be replaced by the same other variable $\mathbf{\llbracket 66.\rrbracket}$ in all places where it occurs. In the rule of substitution for propositional and functional variable there is one restriction to be made as I mentioned before, namely one has to be careful about the letters which one uses for the bound variables, e.g.\
\begin{tabbing}
\hspace{1.7em}$(\exists x)[p\:.\:\varphi(x)]\:.\:(x)\varphi(x)\supset (x)[p\:.\:\varphi(x)]$
\end{tabbing}
is a tautology. Here we cannot substitute $\psi(x)$ for $p$ because
\begin{tabbing}
\hspace{1.7em}$(\exists x)[\psi(x)\:.\:\varphi(x)]\:.\:(x)\varphi(x)\supset (x)[\psi(x)\:.\:\varphi(x)]$
\end{tabbing}
is not a tautology, because here the expression which we substituted contains a variable $x$ which is bound in the expression in which we substitute. Reason: This formula holds for any proposition $p$ but not for any propositional function with the free variable $x$.

Now if we substitute for $p$ an expression $\Phi$ containing perhaps free variables $y,z,\ldots$ (but not the free variable $x$) then $y,z$ will be free in the whole expression. Therefore if $y,z,\ldots$ are replaced by definite things then $\Phi$ will become a proposition because then all free variables contained in it are replaced by definite objects.

Therefore the expression to be substituted must not contain $x$ as a free variable because it would play the role of a propositional function and not of a proposition. In order to avoid such occurrences we have to make in the rule of substitution the stipulation that the expression to be substituted should contain no variable $\mathbf{\llbracket 67.\rrbracket}$ (bound or free) which occurs in the expression in which we substitute bound or free, excluding of course the variable $x$ here. If you add this restriction you obtain the formulation of the rule of substitution which you have in your notes that were distributed.

So far I formulated two rules of inference (implication, substitution). The third is
\begin{itemize}
\item[3\;\;] the rule of defined symbol\index{rule of defined symbol for predicate logic}\index{rule 3} which reads:
\begin{itemize}
\item [1.] For any expressions $\Phi,\Psi$ , $\Phi \supset \Psi$ may be replaced by $\sim \Phi \vee \Psi$ and similarly for $.$ and $\equiv$.
\end{itemize}
\end{itemize}

\vspace{-2ex}

\noindent $\mathbf{\llbracket 68.\rrbracket}$

\vspace{-2ex}

\begin{itemize}
\item[]
\begin{itemize}
\item[2.] $(\exists x)\Phi(x)$ may be replaced by $\sim (x) \sim \Phi(x)$ and vice versa where $\Phi(x)$ is any expression containing the free variable $x$. (So that means that the existential quantifier is defined by means of the universal quantifier in our system.)
\end{itemize}
\end{itemize}

The three rules of inference mentioned so far (implication, substitution, defined symbol) correspond exactly to the three rules of inference which we had in the calculus of propositions. Now we set up a fourth one which is specific for the universal quantifier, namely:

\begin{itemize}
\item [4\;\;] Rule of the universal quantifier:\index{rule of the universal quantifier}\index{rule 4} From $\Pi\supset \Phi(x)$, if $\Pi$ does not contain $x$ as a free variable we can conclude $\mathbf{\llbracket 69.\rrbracket}$ $\Pi \supset (x)\Phi(x)$.
\end{itemize}

That this inference is correct can be seen like this: Assume $\pi$ is a definite proposition and $M(x)$ a definite propositional function with exactly one free variable $x$ and let us assume we know: $\pi \supset M(x)$ \emph{holds for every} $x$. Then I say we can conclude: $\pi\supset (x)M(x)$. For 1. if $\pi$ is false the conclusion holds, 2. if $\pi$ is true then by assumption $M(x)$ is true for every $x$, i.e.\ $(x)M(x)$ is true; hence the conclusion again holds because it is an implication both terms of which are true. So we have proved that in any case $\pi\supset (x)M(x)$ is true if $\pi\supset M(x)$ is true for every $x$. But from this consideration about a particular proposition $\pi$ and a particular propositional $\mathbf{\llbracket 70.\rrbracket}$ function with one free variable $M(x)$ it follows that the above rule of inference yields tautologies if applied to tautologies. Because assume $\Pi\supset \Phi(x)$ is a tautology. Now then $\Pi$ will contain some free variables for propositions $p,q, \ldots$ for functions $\varphi, \psi, \ldots$ and for individuals $y, z, \ldots$ ($x$ does not occur among them) and $\Phi(x)$ will also contain free variables $p, q,\ldots ,\varphi, \psi, \ldots$ and free variables for individuals $x, y, z$ ($x$ among them). Now if you substitute definite propositions for $p, q$, definite predicates for $\varphi, \psi$ and definite objects for $y, z,\ldots$ but leave $x$ where it stands then $\mathbf{\llbracket 71.\rrbracket}$ by this substitution all free variables of $\Pi$ are replaced by individual objects, hence $\Pi$ becomes a definite proposition $\pi$ and all free variables of $\Phi$ excluding $x$ are replaced by objects; hence $\Phi(x)$ becomes a propositional function with one free variable $M(x)$ and we know $\pi \supset M(x)$ is true for any object $x$ because it is obtained by substitution of individual predicates, propositions and objects in a tautology. But then as we have just seen under this assumption $\pi \supset (x)M(x)$ is true. But this argument applies whatever particular predicate, $\mathbf{\llbracket 72.\rrbracket}$ proposition etc.\ we substitute; always the result $\pi \supset (x)M(x)$ is true, i.e.\ $\Pi\supset (x)\Phi(x)$ is a tautology. This rule of course is meant to apply to any other individual variable $y, z$ instead of $x$. So these are the axioms and rules of inference of which one can prove that they are complete: i.e.\ every tautology of the calculus of functions can be derived.

Now I want to give some examples for derivations from these axioms. Again an expression will be called demonstrable or derivable\index{theorem, system of predicate logic} if it can be obtained from Axioms (1)\ldots(4) and Ax.\ 5 by rules 1--4. First of all I wish to remark that, since among our axioms and rules all axioms and rules of the calculus of propositions occur, we can derive from our axioms and rules all formulas and rules which we formerly derived in the calculus of propositions. But the rules are now formulated for all expressions of the calculus of predicates, e.g.\ if $\Phi,\Psi$ are such expressions
\begin{tabbing}
\hspace{15em}\=$\Phi \supset \Psi$\\[.5ex]
\> \underline{$\Psi \supset \Pi$}\\[.5ex]
\> $\Phi \supset \Pi$
\end{tabbing}
So we are justified to use them in the subsequent $\mathbf{\llbracket 73.\rrbracket}$ derivations. At first I mention some further rules of the calculus of propositions which I shall need:
\begin{tabbing}
\hspace{1.7em}\=1. \hspace{2em}\=$P \equiv Q \qquad : \qquad$ \=$P \supset Q,\; Q \supset P$ \;\; and vice versa \\[.5ex]
\>2.\>$P \equiv Q \qquad : \qquad$ \>$\sim P \equiv \: \sim Q$ \\[.5ex]
\>$1'$.\> $p \equiv\: \sim \sim p$ \> ($2'$. $p \equiv p$)\\[.5ex]
\>$3'$.\> $(p \supset q)\:.\: p \supset q$ \> $(p \supset q)\supset (p \supset q)$ \;\; Importation\end{tabbing}
\begin{tabbing}
\hspace{1.7em}\=1. \hspace{2em}\=\underline{$\varphi(y) \supset (\exists x)\varphi(x)$} \\[.5ex]
\>\>$(x)[\sim \varphi(x)] \supset\: \sim \varphi(y)$ \qquad \=Substitution, Ax.\ 5\\[.5ex]
\>\>$\varphi(y)\supset\: \sim (x)[\sim \varphi(x)]$ \>Transposition \quad $\f{\sim \varphi(x)}{\varphi(x)}$ \\
\>\>$\varphi(y) \supset(\exists x)\varphi(x)$ \> \emph{defined symbol}\\[2ex]
\>2. \> \underline{$(x)\varphi(x) \supset (\exists x)\varphi(x)$}\\[.5ex]
\>\>$(x)\varphi(x) \supset \varphi(y)$ \> Ax.\ 5\\[.5ex]
\>\>$\varphi(y)\supset (\exists x)\varphi(x)$ \> 1.
\end{tabbing}

\subsection{Remarks on the term ``tautology'' and ``thinking\\ machines''}\label{2.4}
\pagestyle{myheadings}\markboth{EDITED TEXT}{NOTEBOOK V \;\;---\;\; 1.2.4\;\; Remarks on the terms ``tautology'' and\ldots}

$\mathbf{\llbracket 73.1 \rrbracket}$ Last time I set up a system of axioms and rules of inference from which it is possible to derive all tautologies of the calculus of predicates. Incidentally I wish to mention that the technical term tautology is somewhat out of fashion at present, the word analytical\index{analytical expression}\index{tautological vs. analytical} (which goes back to Kant)\index{Kant, Immanuel} is used in its place, and that has certain advantages because analytical is an indifferent term whereas the term tautological suggests a certain philosophy of logic, namely the theory that the propositions of logic are in some sense void of content, that they say nothing. Of course it is by no means necessary for a $\mathbf{\llbracket 73.2 \rrbracket}$ mathematical logician to adopt this theory, because mathematical logic is a purely mathematical theory which is wholly indifferent towards any philosophical question. So if I use this term tautological I don't want to imply by that any definite standpoint as to the essence of logic, but the term tautological is only to be understood as a shorter expression for universally true. Now as to our axiomatic system the Axioms were as follows 1.\footnote{These axioms, which are omitted at this place in the manuscript, are presumably those on p.\ \textbf{59}.\ of the present Notebook~V.}
\begin{tabbing}
2. Rules of inference\index{rules of inference for predicate logic}\index{four rules of inference for predicate logic}\\[1ex]
\hspace{1.7em}\= 1 \hspace{1em}\= Implication\index{rule of implication}\index{modus ponens}\index{rule 1} \qquad\= $\Phi, \Phi \supset \Psi \quad : \quad \Psi$\\[1ex]
\>2\> Substitution\index{rule 2}\hspace{1em}\> a)\;\;\= individual variables\\[.5ex]
\>\>\>b)\> propositional variables \\[.5ex]
\>\>\>c)\> functional variables\\[1ex]
\>3\> Rule of defined symbol\index{rule 3} \\[.5ex]

\>\hspace{3em}1. For $.$ , $\supset$, $\equiv$\,\, as formerly\\[.5ex]
\>\hspace{3em}2. $(\exists x)\Phi(x)$ may be replaced by $\sim (x) \sim \Phi(x)$ and vice versa\\[1ex]

\>4\>Rule of the universal quantifier\index{rule of the universal quantifier}\index{rule 4} \qquad $\Phi\supset \Psi(x) \quad : \quad \Phi \supset (x)\Psi(x)$
\end{tabbing}

$\mathbf{\llbracket 73.3 \rrbracket}$ It may seem superfluous to formulate so carefully the stipulations about the letters which we have to use for the bound variables here in rule 2 because if you take account of the meaning of the expressions involved you will observe these rules automatically, because otherwise they would either be ambiguous or not have the intended meaning. To this it is to be answered that it is exactly the chief purpose of the axiomatization of logic to avoid this reference to the meaning of the formulas, i.e.\ we want to set up a calculus which can be handled purely mechanically (i.e.\ a calculus which makes thinking superfluous $\mathbf{\llbracket 73.4 \rrbracket}$ and which can replace thinking for certain questions).

In other words we want to put into effect as far as possible Leibnitz's\index{Leibnitz, Gottfried Wilhelm} program\index{Leibnitz's program} of a ``calculus ratiocinator''\index{calculus ratiocinator} which he characterizes by saying that he expects there will be a time in the future when there will be no discussion or reasoning necessary for deciding logical questions but when one will be able simply to say ``calculemus'',\index{calculemus} let us reckon exactly as in questions of elementary arithmetic. This program has been partly carried out by this axiomatic system for logic. For you see that the rules of inference can be applied $\mathbf{\llbracket 73.5 \rrbracket}$ purely mechanically, e.g.\ in order to apply the rule of syllogism $\Phi$, $\Phi \supset \Psi$ you don't have to know what $\Phi$ or $\Psi$ or the sign of implication means, but you have only to look at the outward structure of the two premises. All you have to know in order to apply this rule to two premises is that the second premise contains the $\supset$ and that the part preceding the $\supset$ is conform with the first premise. And similar remarks apply to the other axioms.

Therefore as I mentioned already it would actually be possible to construct a machine which would do the following thing: The supposed machine is to have a crank\index{crank} and whenever you turn the crank once around the machine would write down a tautology of the calculus of predicates and it would write down every existing tautology of the calculus of predicates $\mathbf{\llbracket 73.6 \rrbracket}$ if you turn the crank sufficiently often. So this machine would really replace thinking completely as far as deriving of formulas of the calculus of predicates is concerned. It would be a thinking machine\index{thinking machines} in the literal sense of the word.

For the calculus of propositions you can do even more. You could construct a machine in the form of a typewriter\index{typewriter} such that if you type down a formula of the calculus of propositions then the machine would ring a bell\index{bell} if it is a tautology and if it is not it would not. You could do the same thing for the calculus $\mathbf{\llbracket 73.7 \rrbracket}$ of monadic predicates. But one can prove that it is impossible to construct a machine which would do the same thing for the whole calculus of predicates. So here already one can prove that Leibnitz's\index{Leibnitz, Gottfried Wilhelm} program\index{Leibnitz's program} of the ``calculemus''\index{calculemus} cannot be carried through, i.e.\ one knows that the human mind will never be able to be replaced by a machine already for this comparatively simple question to decide whether a formula is a tautology or not.

\subsection{Theorems and derived rules of the system for\\ predicate logic}\label{2.5}
\pagestyle{myheadings}\markboth{EDITED TEXT}{NOTEBOOK V \;\;---\;\; 1.2.5\;\; Theorems and derived rules of the system\ldots}

\begin{tabbing}
$\mathbf{\llbracket 74\rrbracket}$
\hspace{0.9em}\= $(x)\varphi(x) \supset (\exists x)\varphi(x)$ \hspace{5em}\=Syllogism\\[2ex]

3.\> \underline{$\sim (\exists x)\varphi(x)\equiv (x)\sim \varphi(x)$} \\
\> $\sim \sim (x)\sim \varphi(x) \equiv (x) \sim \varphi(x)$ \hspace{2em} $p \equiv \:\sim \sim p$ \hspace{2em} $\f{(x) \sim \varphi(x)}{p}$ \\
\>$\sim (\exists x)\varphi(x) \equiv (x)\sim \varphi(x)$ \> defined symbol\\[2ex]
4.\> \underline{$p\:.\: (x)\varphi(x) \equiv (x)[p\:.\: \varphi(x)]$}\\[.5ex]
\> $(x)\varphi(x) \supset \varphi(x)$ \\[.5ex]
\> $p\:.\: (x)\varphi(x) \supset p\:.\: \varphi(y)$ \hspace{2.7em} \=Multiplication from left\\[.5ex]
\> $p \: .\: (x)\varphi(x)\supset (y)[p\:.\: \varphi(y)]$ \> \emph{Rule} 4 \; $\Phi:p\:.\: (x)\varphi(x)$ \; $\Psi(y):p\:.\: \varphi(y)$\\[.5ex]
\> $(x)[p\:.\: \varphi(x)]\supset p\:.\: \varphi(y)$ \> Ax.\ 5 \quad Substitution\;\; $\f{p\:.\: \varphi(x)}{\varphi(x)}$\\
\> $p\:.\: \varphi(y) \supset \varphi(y)$ \>$p\:.\: q \supset q$ \quad $\f{\varphi(y)}{q}$\\
\> $p\:.\: \varphi(y) \supset p$ \>$p\:.\: q \supset p$\\[.5ex]
\> $(x)[p\:.\: \varphi(x)]\supset \varphi(y)$ \> Syllogism\\[.5ex]
\> $(x)[p\:.\: \varphi(x)]\supset p$ \> Syllogism\\[.5ex]

$\mathbf{\llbracket 75\rrbracket}$\> $(x)[p\:.\: \varphi(x)]\supset (y)\varphi(y)$ \> Rule 4\\[.5ex]

\> $(x)[p\:.\: \varphi(x)]\supset p\:.\: (y)\varphi(y)$ \> Composition\\[2ex]
5.?\> \underline{$p \vee (x)\varphi(x) \equiv (x)[p \vee \varphi(x)]$}\\[.5ex]
\> $(x)\varphi(x) \supset \varphi(y)$ \> Ax.\ 5 \\[.5ex]
\> $p \vee (x)\varphi(x)\supset p \vee \varphi(y)$ \> Addition from left\\[.5ex]
\> $p \vee (x)\varphi(x)\supset (y)[p \vee \varphi(y)]$ \> Rule 4\\[.5ex]
\> $(x)[p \vee \varphi(x)]\supset p \vee \varphi(y)$ \> Ax.\ 5 \\[.5ex]
\> $p \vee \varphi(y) \supset (\sim p \supset \varphi(y))$ \> $p \vee q \supset (\sim p \supset q)$\\[.5ex]
\> $(x)[p \vee \varphi(x)]\supset (\sim p \supset \varphi(y))$ \quad\= Syllogism\\[.5ex]
\> $(x)[p \vee \varphi(x)]\:.\: \sim p \supset \varphi(y)$ \> Importation\\[.5ex]
\> $(x)[p \vee \varphi(x)]\:.\: \sim p \supset (y)\varphi(y)$ \> Rule 4\\[.5ex]
\> $(x)[p \vee \varphi(x)]$\=$\supset [\sim p \supset (y)\varphi(y)]$\: Exportation\\[.5ex]
\>\> $\supset [p \vee (y)\varphi(y)]$\\[.5ex]

$\mathbf{\llbracket 76\rrbracket}$\\*[.5ex]

6.\> \underline{$(x)[\varphi(x) \supset \psi(x)]\supset [(x)\varphi(x) \supset (x)\psi(x)]$}\\[-.5ex]
\> $(x)[\varphi(x)\supset \psi(x)]\supset [\varphi(y)\supset \psi(y)]$ \quad \= Ax. 5 \quad $\f{\varphi(x)\supset \psi(x)}{\varphi(x)}$ \\[-.5ex]
\>\hspace{4em}$(x)\varphi(x)\supset \varphi(y)$ \> Ax. 5\\[.5ex]
\> $(x)[\varphi(x)\supset \psi(x)]\:.\: (x)\varphi(x) \supset [\varphi(y)\supset \psi(y)]\:.\: \varphi(y)$ \quad Multiplication\\
\> $[\varphi(y)\supset \psi(y)]\:.\: \varphi(y) \supset \psi(y)$ \quad $(p \supset q)\:.\: p \supset q$ \quad $\f{\varphi(y)}{p}$\;\; $\f{\psi(y)}{q}$ \\
\> $(x)[\varphi(x)\supset \psi(x)]\:.\: (x)\varphi(x)$\:\=$\supset \psi(y)$ \hspace{3em}\= Syllogism \\[.5ex]
\>\> $\supset (y)\psi(y)$ \> Rule 4\\[.5ex]
\> $(x)[\varphi(x)\supset \psi(x)]\supset [(x)\varphi(x)\supset (y)\psi(y)]$ \>\> Exportation\\[2ex]

7.\> \emph{Derived Rule I}\index{derived rule I}\index{rule I}\\*[.5ex]
\> $\Phi(x) \quad : \quad (x)\Phi(x)$\\
\hspace{20.3em}$P \supset Q \quad : \quad P\:.\: R \supset Q$\\*
\> $p \:\vee\sim p \supset \Phi(x)$ \quad by addition of premises\quad $Q \quad : \quad P\supset Q$ \\[.5ex]

$\mathbf{\llbracket 77\rrbracket}$\> $p \:\vee\sim p \supset (x)\Phi(x)$ \quad Rule 4\\[.5ex]
\>\underline{$p \:\vee \sim p$}\\[.5ex]
\> $(x)\Phi(x)$ \quad Rule of implication\\[2ex]

8.\> \emph{Derived rule II}\index{derived rule II}\index{rule II}\\*[.5ex]
\> $\Phi(x) \supset \Psi(x)\quad : \quad (x)\Phi(x)\supset (x)\Psi(x)$\\[.5ex]
~~~1.\> $(x)[\Phi(x)\supset \Psi(x)]$\\[.5ex]
~~~2.\> Substitution: $(x)[\Phi(x)\supset \Psi(x)] \supset (x)\Phi(x)\supset (x)\Psi(x)$\\[.5ex]
~~~3.\> Implication\\[2ex]

?9.\> \emph{Derived rule III}\index{derived rule III}\index{rule III}\\*[.5ex]
\>$\Phi(x) \equiv \Psi(x)\quad : \quad$\=$(x)\Phi(x)\equiv (x)\Psi(x)$\\[.5ex]
\>$\Phi(x) \supset \Psi(x)$ \> $(x)\Phi(x)\supset (x)\Psi(x)$\\[.5ex]
\>$\Psi(x) \supset \Phi(x)$ \> \underline{$(x)\Psi(x)\supset (x)\Phi(x)$}\\[.5ex]
\>\> \hspace{3.2em} \ldots \\[.5ex]

$\mathbf{\llbracket 78\rrbracket}$\\*[.5ex]

?10.\> \underline{$\sim (x)\varphi(x) \equiv (\exists x)\sim \varphi(x)$}\\[.5ex]
\> $\varphi(x) \equiv\: \sim \sim \varphi(x)$ \qquad double negation\\[.5ex]
\> $(x)\varphi(x) \equiv (x)\sim \sim \varphi(x)$ \qquad Rule II\\[.5ex]
\> $\sim (x) \varphi(x)$\:\=$\equiv\: \sim (x) \sim \sim \varphi(x)$ \qquad Transposition\\[.5ex]
\>\> $\equiv (\exists x)\sim \varphi(x)$ \qquad defined symbol\\[2ex]
?$10'$.\> \underline{$(x)\varphi(x) \vee (\exists x) \sim \varphi(x)$}\\[.5ex]
\> $(x)\varphi(x) \:\vee \sim (x)\varphi(x)$ \qquad Excluded middle\\[.5ex]
\> $\sim (x)\varphi(x)\supset (\exists x) \sim \varphi(x)$ \qquad ?10.\\[.5ex]
\> \underline{$(x)\varphi(x) \:\vee \sim (x) \varphi(x)$} $\supset (x)\varphi(x) \vee (\exists x)\sim \varphi(x)$\qquad Implication\\[2ex]
?11.\> \underline{$(x)[\varphi(x)\:.\:\psi(x)]\equiv (x)\varphi(x)\:.\:(x)\psi(x)$}\\[.5ex]
\> $\varphi(x)\:.\:\psi(x)\supset \varphi(x)$\\[.5ex]
\> $(x)[\varphi(x)\:.\:\psi(x)]\supset (x)\varphi(x)$ \qquad Rule II\\*[.5ex]
\> $(x)[\varphi(x)\:.\:\psi(x)]\supset (x)\psi(x)$ \qquad \quad $''$\\[.5ex]
\> $(x)[\varphi(x)\:.\:\psi(x)]\supset (x)\varphi(x)\:.\:(x)\psi(x)$ \qquad Composition
\end{tabbing}

\vspace{-6ex}

\begin{equation*}
\begin{rcases}
(x)\varphi(x)\supset \varphi(y)\;\;\\
(x)\psi(x)\supset \psi(y)\;\;
\end{rcases}
\text{\hspace{.2em} Ax. 5 \hspace{15.4em}}
\end{equation*}

\vspace{-4ex}

\begin{tabbing}
$\mathbf{\llbracket 74\rrbracket}$ \hspace{0.9em}\= \kill

\> $(x)\varphi(x)\:.\:(x)\psi(x) \supset \varphi(x)\:.\:\psi(x)$ \qquad Composition \\[.5ex]

$\mathbf{\llbracket 79\rrbracket}$\> $(x)\varphi(x)\:.\:(x)\psi(x) \supset (x)[\varphi(x)\:.\:\psi(x)]$ \qquad Rule 4\\[2ex]

?12.\>\underline{$(x)[\varphi(x)\supset \psi(x)]\:.\:(x)[\psi(x)\supset \chi(x)]\supset (x)[\varphi(x)\supset \chi(x)]$}\\[.5ex]
*\> $(x)[\varphi(x)\supset \psi(x)]\:.\:(x)[\psi(x)\supset \chi(x)]\supset (x)\{ [\varphi(x)\supset \psi(x)]\:.$\\[.5ex]
\` $[\psi(x)\supset \chi(x)]\}$\quad Substitution ?11.\\[.5ex]
\> $[\varphi(x)\supset \psi(x)]\:.\:[\psi(x)\supset \chi(x)]\supset [\varphi(x)\supset \chi(x)]$  \` Substitution Syllogism \\[.5ex]
**\>$(x)\{ [\varphi(x)\supset \psi(x)]\:.\:[\psi(x)\supset \chi(x)]\} \supset (x)[\varphi(x)\supset \chi(x)]$\` Rule II \\[1ex]
\> * and ** with Syllogism give the result.\\[.5ex]

$\mathbf{\llbracket 80.\rrbracket}$\\*[.5ex]

13.\> \underline{Rule \quad $\Psi(x)\supset \Phi \quad : \quad (\exists x)\Psi(x) \supset \Psi$}\\[.5ex]
\> $\sim \Phi \supset\: \sim \Psi(x)$\\[.5ex]
\> $\sim \Phi \supset (x)\sim \Psi(x)$\\[.5ex]
\> $\sim (x)\sim \Psi (x) \supset \Phi$\\[.5ex]
\> $(\exists x)\Psi(x)\supset \Phi$\\[2ex]

$13'$.\> \underline{$\varphi(y)\supset (\exists x)\varphi(x)$}\\[.5ex]
\> $(x)\sim \varphi(x)\supset\: \sim \varphi(y)$\\[.5ex]
\> $\varphi(y)\supset\: \sim (x)\sim \varphi(x)$ \qquad defined symbol\\[2ex]

14.\> \underline{$(x)[\varphi(x) \supset \psi(x)]\supset[(\exists x)\varphi(x)\supset (\exists x)\psi(x)]$}\\[.5ex]
\> $[\varphi(x) \supset \psi(x)]\supset$\:\=$[\sim \psi(x) \supset\: \sim \varphi (x)]$\\*[.5ex]
$\times$\> $(x)$ \hspace{1em}$''$ \>$(x)$ \hspace{2.1em}$''$\\[.5ex]
$\times$\> $(x)[\sim \psi(x) \supset\: \sim \varphi (x)]\supset (x)\sim \psi(x)\supset (x) \sim \varphi (x) $ \\[.5ex]
$\times$\> $[(x)\sim \psi(x)\supset (x) \sim \varphi (x)] \supset \:\sim (x) \sim \varphi (x) \supset \:\sim (x) \sim \psi(x)$\\
\` $(p \supset q)\supset (\sim q \supset\: \sim p)$ \quad $\f{(x)\sim \psi(x)}{p}$ \;\; $\f{(x)\sim \varphi(x)}{q}$\\
\> $(x)[\varphi(x)\supset \psi(x)]\supset[\sim (x)\sim \varphi(x) \supset \:\sim (x) \sim \psi(x)]$\` Rule of defined \\*[.5ex]
\` symbol\\[.5ex]

$\mathbf{\llbracket 81.\rrbracket}$\\*[.5ex]

15.\> Rule corresponding to 14.\\[2ex]
16. \> \underline{$(\exists x)[\varphi(x) \vee \psi (x)] \equiv (\exists x)\varphi(x) \vee (\exists x)\psi(x)$}\\[.5ex]
\> $\varphi(x) \supset \varphi(x) \vee \psi(x)$\\[.5ex]
\> $(\exists x)\varphi(x) \supset (\exists x)[\varphi(x) \vee \psi(x)]$\\
\> \ldots\\[.5ex]
\> Dilemma\\[-1.5ex]
\> \underline{\hspace{12.5em}}\\[.5ex]
\> $\varphi(y) \supset (\exists x)\varphi(x)$\\[.5ex]
\> $\psi(y) \supset (\exists x)\psi(x)$\\[.5ex]
\> $\varphi(y) \vee \psi(y)$\:\=$\supset (\exists x)\varphi(x) \vee (\exists x) \psi(x)$\\*[.5ex]
\> $(\exists y)[$\quad $''$\quad $]$ \>$\supset$ \qquad $''$ \qquad \qquad $''$
\end{tabbing}

An example where we have to substitute for $\varphi(x)$ something containing other free variables besides $x$:
\begin{tabbing}
\hspace{3em}\=\underline{$(y)(x)\psi(xy)\equiv (x)(y)\psi(xy)$}\\[.5ex]
\> \underline{$(x)\varphi(x) \supset \varphi(y)$} \\[.5ex]
\> $(x)\varphi(x) \supset \varphi(u)$ \qquad $\f{\psi(xy)}{\varphi(x)}$\\[.5ex]
\hspace{1.7em}*\> $(x)\psi(xy) \supset \psi(uy)$\\[.5ex]
\> $(z)\varphi(z) \supset \varphi(y)$ \qquad $\f{(x)\psi(xz)}{\varphi(z)}$\\[.5ex]
\hspace{1.7em}*\> $(z)(x)\psi(xz) \supset (x)\psi(xy)$ \qquad\= *\:* Syllogism \\[.5ex]
\> $(z)(x)\psi(xz) \supset \psi(uy)$ \> Rule 4 \quad\= $y$\\*[.5ex]
\> $(z)(x)\psi(xz) \supset (y)\psi(uy)$ \>\hspace{.9em} $''$ \> $u$\\[.5ex]
\> $($\=$z)(x)\psi(xz) \supset ($\=$u)(y)\psi(uy)$ \\[.5ex]
\>\> $y$ \> $x$\\[.5ex]
\> $(y)(x)\psi(xz) \supset (x)(y)\psi(uy)$
\end{tabbing}

\subsection{Existential presuppositions}\label{2.6}
\pagestyle{myheadings}\markboth{EDITED TEXT}{NOTEBOOK V \;\;---\;\; 1.2.6\;\; Existential presuppositions}

$\mathbf{\llbracket 82.\rrbracket}$ I have mentioned already that among the tautological formulas of the calculus of predicates are in particular those which express the Aristotelian moods of inference, but that not all of the 19 Aristotelian moods are really valid in the calculus of propositions. Some of them require an additional third premise in order to be valid, namely that the predicates involved be not vacuous;\index{existential presuppositions}\index{vacuous predicates} e.g.\ the mood Darapti\index{Darapti} is one of those not valid, it says
\begin{tabbing}
\hspace{1.7em}\=$M$a$S$, $M$a$P$ \; : \; $S$i$P$, \quad in symbols:\\[.5ex]
\>$(x)[M(x)\supset S(x)]\:.\:(x)[M(x)\supset P(x)] \supset (\exists x)[S(x)\:.\:P(x)]$
\end{tabbing}
But this is not a tautological formula because that would mean it holds for any monadic predicates $M,S,P$ whatsoever. But $\mathbf{\llbracket 83.\rrbracket}$ we can easily name predicates for which it is wrong; namely if you take for $M$ a vacuous predicate which belongs to no object,\index{vacuous predicates} say e.g.\ the predicate president of America born in South Bend and for $S$ and $P$ any two mutually exclusive predicates, i.e.\ such that no $S$ is $P$, then the above formula will be wrong because the two premises are both true. Since $\mathbf{\llbracket 84.\rrbracket}$ $M(x)$ is false for every $x$ we have $M(x)\supset S(x)$ is true for every $x$ (because it is an implication with false first term); likewise $M(x)\supset P(x)$ is true for every $x$. I.e.\ the premises are both true but the conclusion is false because $S,P$ are supposed to be two predicates such that there is no $S$ which is a $P$. Hence for the particular predicate we chose the first term of this whole implication is true and the second is false, i.e.\ the whole formula is false. So there are predicates which substituted in this formula yield a false proposition, hence this formula is not a tautology. If we want to transform that expression into a real tautology we have to add the further premise that $M$ is not $\mathbf{\llbracket 85.\rrbracket}$ vacuous, i.e.\
\begin{tabbing}
\hspace{1.7em}$(\exists x)M(x)\:.\:(x)[M(x)\supset S(x)]\:.\:(x)[M(x)\supset P(x)]\supset (\exists x)[S(x)\:.\:P(x)]$
\end{tabbing}
would really be a tautology. Altogether there are four of the $19$ Aristotelian moods which require this additional premise. Furthermore $S\text{a}P \supset S\text{i}P$, $P$i$S$ (conversion) as I mentioned last time also requires that $S$ is non-vacuous. Also $S\text{a}P\supset\: \sim (S\text{e}P)$, i.e.\ $S$a$P$ and $S$e$P$ cannot both be true, does not hold in the logical calculus because if $S$ is vacuous both $S$a$P$ and $S$e$P$ are true $(x)[S(x)\supset P(x)]\:.\:(x)[S(x)\supset\: \sim P(x)]$; $S(x)=$ $x$ is a president of the States born in South Bend, $P(x) =$ $x$ is bald, then both
\begin{tabbing}
\hspace{1.7em}\= Every president\ldots \;\;\= is bald\\[.5ex]
\> No president\ldots \> is bald
\end{tabbing}

So we see Aristotle\index{Aristotle} makes the implicit assumption that all predicates which he speaks of are non-vacuous; in the logistic calculus however we do not make this assumption, i.e.\ all tautologies and all formulas derivable from our axioms hold for any predicates whatsoever they may be, vacuous or not. $\mathbf{\llbracket 86.\rrbracket}$ Now one may ask: which procedure is preferable, to formulate the laws of logic in such a way that they hold for all predicates vacuous and non-vacuous or in such a way that they hold only for non-vacuous. I think there can be no doubt that the logistic way is preferable for many reasons:\index{failure of traditional logic}

1. As we saw it may depend on purely empirical facts whether or not a predicate is vacuous (as we saw in the example of a president of America born in South Bend). Therefore if we don't admit vacuous predicates it will depend on empirical facts which predicates are to be admitted in logical reasonings or which inferences are valid, but that $\mathbf{\llbracket 87.\rrbracket}$ is very undesirable. Whether a predicate can be used in reasoning (drawing inferences) should depend only on mere logical considerations and not on empirical facts.

But a second and still more important argument is this: that to exclude vacuous predicates would be a very serious hampering, e.g.\ in mathematical reasoning, because it happens frequently that we have to form predicates of which we don't know in the beginning of an argument whether or not they are vacuous, e.g.\ in indirect proofs. If we want to prove that there does not exist an algebraic equation whose root is $\pi$ we operate $\mathbf{\llbracket 88.\rrbracket}$ with the predicate ``algebraic equation with root $\pi$'' and use it in conclusions, and later on it turns out that this predicate is vacuous. But also in everyday life it happens frequently that we have to make general assertions about predicates of which we don't know whether they are vacuous. E.g.\ assume that in a university there is the rule that examinations may be repeated arbitrarily often; then we can make the statement: A student which has\ldots\ ten times is allowed to\ldots\ for an eleventh time. But if we want to exclude vacuous predicates we cannot express this true proposition if we don't know whether there exists a student who has\ldots\ But of course this proposition (rule) has nothing to do with the existence of a student\ldots\ Or e.g.\ excluding vacuous predicates has the consequence that we cannot always form the conjunction of two predicates, e.g.\ president of U.S.A. is an admissible predicate, born in South Bend is admissible, but president of America born in South Bend is not admissible. So if we want to avoid absolutely unnecessary complications we must not exclude the vacuous predicates and have to formulate the laws of logic in such a way that they apply both to vacuous and non-vacuous predicates. I don't say that it is false to exclude them, but it leads to absolutely unnecessary complications.

\subsection{Classes}\label{2.7}
\pagestyle{myheadings}\markboth{EDITED TEXT}{NOTEBOOK V \;\;---\;\; 1.2.7\;\; Classes}

As to the 15 valid moods of Aristotle\index{Aristotle} they can all be expressed by one logistic formula. However in order to do that I have first to embody the calculus of monadic predicates in a different form, namely in the form of the calculus of classes.\index{calculus of classes}\index{classes} $\mathbf{\llbracket 89.\rrbracket}$ The calculus of classes also yields the solution of the decision problem for formulas with only monadic predicates.

If we have an arbitrary monadic predicate, say $P$, then we can consider the extension of this predicate,\index{extension of a predicate} i.e.\ the totality of all objects satisfying $P$; it is denoted by $\hat{x}[P(x)]$. These extensions of monadic predicates are all called classes.\index{classes} So this symbol $\hat{x}$ means: the class of objects $x$ such that the subsequent is true. It is applied also to propositional functions, e.g.\ $\hat{x}[I(x)\:.\: x>7]$ means ``the class of integers greater than seven''. $\mathbf{\llbracket 90.\rrbracket}$ So to any monadic predicate belongs a uniquely determined class of objects as its ``extension'', but of course there may be different predicates with the same extension, as e.g.\ the two predicates: heat conducting, elasticity conducting. These are two entirely different predicates, but every object which has the first property also has the second one and vice versa; therefore their extension is the same, i.e.\ if $H,E$ denotes them, $\hat{x}[H(x)]=\hat{x}[E(x)]$ although $H \neq E$. I am writing the symbol of identity and distinctness in between the two identical objects as is usual in mathematics. I shall speak about this way of writing in more detail later. In general we have if $\varphi,\psi$ are two monadic predicates then\index{extensionality}
\begin{tabbing}
\hspace{1.7em}$\hat{x}[\varphi(x)]=\hat{x}[\psi(x)] \equiv (x)[\varphi(x)\equiv \psi(x)]$
\end{tabbing}
This equivalence expresses the essential property of extensions of predicates. It is to be noted that we have not defined what classes are because we explained it by the term extension, and extensions we explained by the term totality, and a totality is the same thing as a class. So this definition would be circular. The real state of affairs is this: that we consider $\hat{x}$ as a new primitive (undefined) term, which satisfies this axiom here. Russell\index{Russell, Bertrand} however has shown that one can dispense with this $\hat{x}$ as a primitive term by introducing it by a kind of implicit definition, but that would take too much time to explain it; so we simply can consider it as a primitive.

The letters $\alpha, \beta, \gamma,\ldots$ are used as variables for classes and the statement that
$\mathbf{\llbracket Notebook\; VI \rrbracket}$
\pagestyle{myheadings}\markboth{EDITED TEXT}{NOTEBOOK VI \;\;---\;\; 1.2.7\;\; Classes}
$\mathbf{\llbracket 91. \rrbracket}$ an object $a$ belongs to $\alpha$ (or is an element of $\alpha$) by $a\,\varepsilon\,\alpha$. Hence\index{extensionality}\index{comprehension}
\[
\fbox{\begin{minipage}{7.8em}
$y\,\varepsilon\,\hat{x}[\varphi(x)]\equiv\varphi(y)$
\end{minipage}} \quad\quad {\rm Furthermore}
  \begin{cases}
   \alpha = \hat{x}[x\,\varepsilon\,\alpha]\\
   (x)[x\,\varepsilon\,\alpha\equiv x\,\varepsilon\,\beta] \supset \alpha =\beta
  \end{cases}
\]
So far we spoke only of extensions of monadic predicates; we can also introduce extensions of dyadic (and polyadic) predicates. If e.g.\ $Q$ is a dyadic predicate then $\hat{x}\hat{y}[Q(xy)]$ (called the extension of $Q$) will be something that satisfies the condition:
\begin{tabbing}
\hspace{1.7em}$\hat{x}\hat{y}[\psi(xy)]=\hat{x}\hat{y}[\chi(xy)] .\equiv . (x,y)[\psi(xy)\equiv \chi(xy)]$
\end{tabbing}
e.g.\ the class of pairs $(x,y)$ such that $Q(xy)$ would $\mathbf{\llbracket 92. \rrbracket}$ be something which satisfies this condition, but the extension of a relation is not defined as the class of ordered pairs, but is considered as an undefined term because ordered pair is defined in terms of extension of relations. An example for this formula, i.e.\ an example of two different dyadic predicates which have the same extension would be $x<y$, $x>y\vee x=y$, $x$ exerts an electrostatic attraction on $y$, $x$ and $y$ are loaded by electricities of different sign.

Extensions of monadic predicates are called classes,\index{classes} extensions of polyadic predicates are called relations\index{relations} in logistic. So in logistic the term relation is used not for the polyadic predicates themselves but for their extensions, that conflicts with the meaning of the term relation in everyday life and also with the meaning in which I introduced this term a few lectures ago, but since it is usual to use this term relation in this extensional sense I shall stick to this use and the trouble is that there is no better term. If $R$ is a relation, the statement that $x$ bears $R$ $\mathbf{\llbracket 93. \rrbracket}$ to $y$ is denoted by $xRy$. This way of writing, namely to write the symbol denoting the relation between the symbols denoting the objects for which the relation is asserted to hold, is adapted to the notation of mathematics, e.g.\ $<$, $x<y$, $=$, $x=y$. Of course we have:
\begin{tabbing}
\hspace{1.7em}$(x,y)[xRy\equiv xSy]\supset R=S$
\end{tabbing}
for any two relations $R,S$, exactly as before $(x)[x\,\varepsilon\,\alpha\equiv x\,\varepsilon\,\beta ] \supset \alpha =\beta$. So a relation is uniquely determined if you know all the pairs which have this relation because by this formula there cannot exist two different relations which subsist between the same pairs (although there can exist many different dyadic predicates).

Therefore a relation can be represented e.g.\ by a figure of arrows
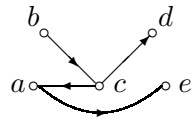
\begin{center}
\begin{picture}(100,40)(-5,10)

\put(45,20){\circle{3}}
\put(20,20){\circle{3}}
\put(24,40){\circle{3}}
\put(65,40){\circle{3}}
\put(70,20){\circle{3}}

\put(43.5,20){\line(-1,0){22}}
\put(43.5,21){\line(-1,1){18}}
\put(46,21){\line(1,1){18}}

\put(17,20){\small\makebox(0,0)[r]{$a$}}
\put(20,43){\small\makebox(0,0)[b]{$b$}}
\put(55,20){\small\makebox(0,0)[r]{$c$}}
\put(70,43){\small\makebox(0,0)[b]{$d$}}
\put(80,20){\small\makebox(0,0)[r]{$e$}}

\put(32,20){\vector(-1,0){2}}
\put(35,29.3){\vector(1,-1){2}}
\put(62,37){\vector(1,1){2}}
\put(46,10){\vector(1,0){2}}

\qbezier(22,20)(45,0)(68.5,20)

\end{picture}
\end{center}
or by a quadratic scheme e.g.\ \begin{center}
\begin{tabular}{ c|c c c c c }
& $a$ & $b$ & $c$ & $d$ & $e$\\
\hline
$a$ &&&&& $\bullet$\\
$b$ &&& $\bullet$ \\
$c$ & $\bullet$ &&& $\bullet$\\
$d$ & \\
$e$ &
\end{tabular}
\end{center}
Such a figure determines a unique relation; in general it will be infinite.

The letters $R,S,T$ are mostly used as variables for relations. But now let us return to the extensions of monadic predicates, i.e.\ the classes for which we want to set up a calculus.

First we have two particular classes $\bigwedge$ (vacuous class),\index{vacuous class}\index{zero class} $\bigvee$ (the universal class)\index{universal class} which are defined as the extension $\mathbf{\llbracket 94. \rrbracket}$ of a vacuous predicate and of a predicate that belongs to everything. So
\begin{tabbing}
\hspace{1.7em}\=$\bigwedge\;$\=$=\; \hat{x}[\varphi(x)\: .\, \sim\varphi(x)]$\\[.5ex]
\>$\bigvee$\>$=\; \hat{x}[\varphi(x)\:\vee \sim\varphi(x)]$
\end{tabbing}
It makes no difference which vacuous predicate I take for defining $\bigwedge$. If $A$, $B$ are two different vacuous predicates then $\hat{x}[A(x)]=\hat{x}[B(x)]$ because $(x)[A(x)\equiv B(x)]$. And similarly if $C,D$ are two different predicates belonging to everything $\hat{x}[C(x)]=\hat{x}[D(x)]$ because $(x)[C(x)\equiv D(x)]$, i.e.\ there exists exactly one 0-class and exactly one $\mathbf{\llbracket 95. \rrbracket}$ universal class, although of course there exist many different vacuous predicates. But they all have the same extension, namely nothing which is denoted by $\bigwedge$. So the zero class\index{zero class} is the class with no elements $(x)[\sim x\,\varepsilon\, \bigwedge]$, the universal class\index{universal class} is the class of which every object is an element $(x)(x\,\varepsilon\, \bigvee )$; $\bigwedge$ and $\bigvee$ are sometimes denoted by 0 and 1 because of certain analogies with arithmetic.

Next we can introduce certain operations for classes which are analogous to the arithmetical operations: namely
\begin{tabbing}
\hspace{.5em}\= Addition\index{addition of classes} or sum\index{sum of classes}\hspace{7.5em}\= $\alpha +\beta$ \= $=\;\; \hat{x}[x\,\varepsilon\,\alpha\vee x\,\varepsilon\,\beta]$\\[.5ex]
\` $y\,\varepsilon\, \alpha +\beta\equiv y\,\varepsilon\,\hat{x}[x\,\varepsilon\,\alpha\vee x\,\varepsilon\,\beta]\equiv y\,\varepsilon\,\alpha\vee y\,\varepsilon\,\beta$\\[.5ex]
\hspace{15.8em} mathematician or democrat \\[.5ex]
\> Multiplication\index{multiplication of classes} or intersection\index{intersection of classes} \> $\alpha\cdot\beta$ \> $=\;\; \hat{x}[x\,\varepsilon\,\alpha\: .\: x\,\varepsilon\,\beta]$\\[.5ex]
\hspace{15.8em} mathematician democrat \\[.5ex]
\> Opposite\index{opposite of a class} or complement\index{complement of a class} \> $-\alpha$ \> $=\;\; \hat{x}[\sim x\,\varepsilon\,\alpha]$\quad or \quad$\overline{\alpha}$\\[.5ex]
\hspace{15.8em} non mathematician \\[.5ex]
\>Difference\index{difference of classes}\> $\alpha -\beta$\> $=\;\;\alpha\cdot(-\beta)\;
=\;\hat{x}[x\,\varepsilon\,\alpha\: .\, \sim x\,\varepsilon\,\beta]$\\[.5ex]
\hspace{15.8em} mathematician not democrat \\
\hspace{15.8em} (New Yorker not sick)\footnotemark
\end{tabbing}
\footnotetext{On the right of this table, two intersecting circles, as in Euler or Venn diagrams, are drawn in the manuscript.}

Furthermore we have a relation classes which corresponds to the arithmetic relation of $<$, namely the relation of subclass\index{subclass, relation of}
\begin{tabbing}
\hspace{1.7em}$\alpha\subseteq\beta\equiv(x)[x\,\varepsilon\,\alpha\supset x\,\varepsilon\,\beta]$ \qquad Man $\subseteq$ Mortal
\end{tabbing}
All these operations obey laws very similar $\mathbf{\llbracket 96. \rrbracket}$ to the corresponding arithmetical laws: e.g.\
\begin{tabbing}
\hspace{1.7em}\= $\alpha +\beta = \beta +\alpha$\hspace{2em}$\alpha\cdot\beta = \beta\cdot\alpha$\\[.5ex]
\>$(\alpha +\beta)+\gamma = \alpha +(\beta +\gamma)$\hspace{2em}$(\alpha\cdot\beta)\cdot\gamma = \alpha\cdot(\beta\cdot\gamma)$\\[.5ex]
\>$(\alpha +\beta)\cdot\gamma = \alpha\cdot\gamma +\beta\cdot\gamma$\\[.5ex]
\>$(\alpha\cdot\beta)+\gamma = (\alpha +\gamma)\cdot(\beta +\gamma)$
\end{tabbing}
They follow from the corresponding laws of the calculus of propositions: e.g.\
\begin{tabbing}
\hspace{1.7em}\= $x\,\varepsilon(\alpha +\beta)\equiv x\,\varepsilon\,\alpha\vee x\,\varepsilon\,\beta\equiv x\,\varepsilon\,\beta\vee x\,\varepsilon\,\alpha\equiv x\,\varepsilon(\beta +\alpha)$\\[1ex]
\>$x\,\varepsilon(\alpha +\beta)\cdot\gamma\equiv x\,\varepsilon(\alpha +\beta)\: .\: x\,\varepsilon\,\gamma\equiv(x\,\varepsilon\,\alpha\vee x\,\varepsilon\,\beta)\: .\: x\,\varepsilon\,\gamma$\\[.5ex]
\`$\equiv(x\,\varepsilon\,\alpha\: .\: x\,\varepsilon\,\gamma)\vee (x\,\varepsilon\,\beta
 \: .\: x\,\varepsilon\,\gamma)\equiv x\,\varepsilon\,\alpha\cdot\gamma\vee x\,\varepsilon\,\beta\cdot\gamma \equiv x\,\varepsilon(\alpha\cdot\gamma +\beta\cdot\gamma)$\\[1ex]
\>$\alpha + 0\;$\=$=\alpha$\hspace{2em}\=$\alpha\cdot 0$\hspace{.7em}\=$=0$\\[.5ex]
\>$\alpha\cdot 1$\>$=\alpha$\>$\alpha +1$\>$=1$\\[1ex]
\>$(x)\sim(x\,\varepsilon\, 0)$\hspace{3em}$x\,\varepsilon(\alpha +0)\equiv x\,\varepsilon\,\alpha\vee x\,\varepsilon\, 0\equiv x\,\varepsilon\,\alpha$\\[.5ex]
\>$(x)(x\,\varepsilon\, 1)$\footnotemark
\end{tabbing}
\footnotetext{On the right of this table, three intersecting circles, as in Euler or Venn diagrams, with $\alpha$, $\beta$ and perhaps $\gamma$ marked in them, and some areas shaded, are drawn in the manuscript.}
\begin{tabbing}
\hspace{3.5em}$\alpha\subseteq\beta$\hspace{7em}\=$\alpha\subseteq\beta\:.\:\beta\subseteq\gamma\supset\alpha\subseteq\gamma$\\[.5ex]
\hspace{3.5em}$\gamma\subseteq\delta$\hspace{8em} Law of transitivity\\[-2ex]
\hspace{1.7em}\uline{\hspace{6em}}\\
\hspace{1.7em}$\alpha +\gamma\subseteq\beta +\delta$\\
\hspace{2.2em}$\alpha\cdot\gamma\subseteq\beta\cdot\delta$\>$\alpha\subseteq\beta\: .\:\beta\subseteq\alpha\supset\alpha =\beta$.
\end{tabbing}
Laws different from arithmetical:
\begin{tabbing}
\hspace{1.7em}\=$\alpha +\alpha =\alpha\cdot\alpha=\alpha$\hspace{3.5em}$x\,\varepsilon\,\alpha  +\alpha\equiv x\,\varepsilon\,\alpha\vee x\,\varepsilon\,\alpha\equiv x\,\varepsilon\,\alpha$\\[.5ex]
\>$\alpha\subseteq\beta\supset[\alpha +\beta =\beta\;\; .\;\; \alpha\cdot\beta=\alpha]$\hspace{2em}$\beta\subseteq\alpha +\beta$\hspace{2em}\=$\alpha\subseteq\beta$\\
\>\>$\beta\subseteq\beta$\\[-2ex]
\hspace{22.7em}\uline{\hspace{8em}}\\
\hspace{22.7em}$\alpha +\beta\subseteq\beta +\beta=\beta$
\end{tabbing}
\begin{tabbing}
$\mathbf{\llbracket 97. \rrbracket}$\\*[1ex]
\hspace{1.7em}\=$-(\alpha +\beta)=(-\alpha)\cdot(-\beta)$\hspace{2em} De Morgan\index{De Morgan law for classes}\\[.5ex]
\` $x\,\varepsilon\, -(\alpha +\beta)\equiv\;\sim x\,\varepsilon\, (\alpha +\beta)\equiv\;\sim(x\,\varepsilon\,\alpha\vee x\,\varepsilon\,\beta)\equiv\; \sim(x\,\varepsilon\,\alpha)\: .\:\sim(x\,\varepsilon\,\beta)\equiv$\\*
\` $x\,\varepsilon\, -\alpha\: .\: x\,\varepsilon\, -\beta\equiv\; x\,\varepsilon\, (-\alpha)\cdot(-\beta)$\\
\>$-(\alpha\cdot\beta)=(-\alpha)+(-\beta)$\\[.5ex]
\>$\alpha\cdot(-\alpha)=0$\hspace{2em}$\alpha +(-\alpha)=1$\\[.5ex]
\>$-(-\alpha)=\alpha$
\end{tabbing}
The complement of $\alpha$ is sometimes also denoted by $\overline{\alpha}$ (so that $\overline{\alpha}= -\alpha$).

Exercise Law for difference:
\begin{tabbing}
\hspace{1.7em}\=$\alpha\cdot(\beta -\gamma)=\alpha\cdot\beta -\alpha\cdot\gamma$\\[.5ex]
\>$\alpha\cdot\beta = \alpha -(\alpha -\beta)$\\[.5ex]
\>$\alpha\subseteq\beta \supset \overline{\beta}\subseteq\overline{\alpha}$
\end{tabbing}

\subsection{Classes and Aristotelian moods}\label{2.8}
\pagestyle{myheadings}\markboth{EDITED TEXT}{NOTEBOOK VI \;\;---\;\; 1.2.8\;\; Classes and Aristotelian moods}

If $\alpha\cdot\beta=0$, that means the classes $\alpha$ and $\beta$ have no common element, then $\alpha$ and $\beta$ are called mutually exclusive. We can now formulate the four Aristotelian types of judgement\index{a, e, i, o propositions} a, e, i, o also in the symbolism of the calculus of classes as follows:
\begin{tabbing}
\hspace{1.7em}\=$\alpha\, {\rm a}\, \beta$\hspace{1em}\=$\equiv$\hspace{1.7em}$\alpha\subseteq\beta$\hspace{1.7em}\=
$\equiv$\hspace{1.5em}\=$\alpha\cdot\overline{\beta}=0$\index{a proposition}\\[.5ex]
$\mathbf{\llbracket 98. \rrbracket}$\\*[.5ex]
\>$\alpha\, {\rm e}\, \beta$\>$\equiv$\hspace{1.6em}$\alpha\subseteq\overline{\beta}$\>$\equiv$\>$\alpha\cdot\beta=0$\index{e proposition}\\[.5ex]
\>$\alpha\, {\rm i}\, \beta$\>$\equiv$ $\;\sim(\alpha\subseteq\overline{\beta})$\>$\equiv$\>$\alpha\cdot\beta\neq 0$\index{i proposition}\\[.5ex]
\>$\alpha\, {\rm o}\, \beta$\>$\equiv$ $\;\sim(\alpha\subseteq\beta)$\>$\equiv$\>$\alpha\cdot\overline{\beta}\neq 0$\index{o proposition}
\end{tabbing}
So all of these four types of judgements can be expressed by the vanishing, respectively not vanishing, of certain intersections.

Now the formula which compresses all of the 15 valid Aristotelian inferences reads like this
\[
\sim(\alpha\cdot\beta =0\;\; .\;\; \overline{\alpha}\cdot\gamma =0\;\; .\;\; \beta\cdot\gamma\neq 0)
\]
So this is a universally true formula because $\alpha\cdot\beta=0$ means $\beta$ outside of $\alpha$, $\overline{\alpha}\cdot\gamma =0$ means $\gamma$ inside of $\alpha$. If $\beta$ outside $\gamma$ inside they can have no element in $\mathbf{\llbracket 99. \rrbracket}$ common, i.e.\ the two first propositions imply $\beta\cdot\gamma =0$, i.e.\ it cannot be that all three of them are true. Now since this says that all three of them cannot be true you can always conclude the negation of the third from the two others; e.g.\ \begin{tabbing}
\hspace{1.7em}\=$\alpha\cdot\beta =0\;\; .\;\; \overline{\alpha}\cdot\gamma =0\;\supset\;\beta\cdot\gamma =0$\\[.5ex]
\>$\alpha\cdot\beta =0\;\; .\;\; \beta\cdot\gamma\neq 0\;\supset\;\overline{\alpha}\cdot\gamma\neq 0$\quad etc. \end{tabbing}
and in this way you obtain all valid 15 moods if you substitute for $\alpha,\beta,\gamma$ in an appropriate way the minor term, the major term and the middle term or their negation, e.g.\
\begin{tabbing}
$\mathbf{\llbracket 100. \rrbracket}$\\*[1ex]
\hspace{1.7em} I \quad \emph{Barbara}\index{Barbara}
\hspace{2em}\begin{tabular}{ l|l }
$M$a$P$ & \\[-1ex]
& \hspace{.5em}$S$a$P$ \\[-1ex]
$S$a$M$&
\end{tabular}\\[1ex]
\hspace{2.7em}$M\cdot\overline{P} =0\;\; .\;\; S\cdot\overline{M} =0\;\supset \;S\cdot \overline{P}=0$\\[.5ex]
\hspace{1.7em}$\sim(M\cdot\overline{P} =0\;\; .\;\; S\cdot\overline{M} =0\,\;\; .\,\;\; S\cdot\overline{P}\neq 0)$\\[.5ex]
\hspace{1.7em}$\alpha=M$\hspace{2em}$\beta=\overline{P}$\hspace{2em}$\gamma =S$\\[2ex]

\hspace{1.7em} III \quad \emph{Feriso}\index{Feriso}
\hspace{2em}\begin{tabular}{ l|l }
$M$e$P$ & \\[-1ex]
& \hspace{.5em}$S$o$P$ \\[-1ex]
$M$i$S$&
\end{tabular}\\[1ex]
\hspace{3.6em}$M\cdot P =0\;\; .\;\; M\cdot S \neq 0\;\supset\; S\cdot\overline{P}\neq 0$\\[.5ex]
\hspace{1.7em}$\sim(M\cdot P =0\;\; .\;\; M\cdot S \neq 0 \;\; .\;\; S\cdot\overline{P}= 0)$\\[.5ex]
\hspace{1.7em}$\alpha=P$\hspace{2em}$\beta=M$\hspace{2em}$\gamma =S$.
\end{tabbing}

The four moods which require an additional premise can also be expressed by one formula, namely:
\[
\sim(\alpha\neq 0\;\; .\;\;\alpha\cdot\beta =0\;\; .\;\; \alpha\cdot\gamma =0\;\; .\;\; \overline{\beta}\cdot\overline{\gamma}= 0)
\]

$\mathbf{\llbracket 101. \rrbracket}$ Darapti\index{Darapti}
\begin{tabbing}
\hspace{9.8em}\=$M$a$P$\\[.5ex]
\>$M$a$S$\\[-2ex]
\>\uline{\hspace{2.8em}}\\
\>$\;S$i$P$
\end{tabbing}
e.g.\ is obtained by taking
\begin{tabbing}
\hspace{1.7em}\=$\alpha =M$\hspace{2em}$\beta=\overline{P}$\hspace{2em}$\gamma =\overline{S}$\\[.5ex]
\>$M$a$P\: .\: M$a$S\supset S$i$P$\\[.5ex]
\>$M\cdot\overline{P}=0\:.\: M\cdot\overline{S}=0\:\supset\: S\cdot P\neq 0$
\end{tabbing}
However, this second formula is an easy consequence of the first, i.e.\ we can derive it by two applications of the first. To this end we have only to note that $\alpha\neq 0$ can be expressed by $\alpha\, {\rm i}\,\alpha$ because
\begin{tabbing}
\hspace{1.7em}\=$\varphi\,{\rm i}\,\psi\equiv(\exists x)[\varphi(x).\psi(x)]$\\[.5ex]
\>$\varphi\,{\rm i}\,\varphi\equiv(\exists x)[\varphi(x)\: .\:\varphi(x)]\equiv(\exists x)\varphi(x)$\\[1ex]
\>$\sim(\alpha\cdot\beta =0\;\; .\;\; \overline{\alpha}\cdot\gamma =0\;\; .\;\; \beta\cdot\gamma\neq 0)$\\[.5ex]
\>$\alpha\cdot\alpha\neq 0$\hspace{2em}$\alpha\cdot\beta=0$\hspace{2em}
$\alpha\cdot\overline{\beta}=0$\\[.5ex]
\>$\alpha:\beta$\hspace{2em}$\beta :\alpha$\hspace{2em}$\gamma :\alpha$\\[1ex]
\>$\alpha\cdot\overline{\beta}\neq 0$\hspace{2em}$\alpha\cdot\gamma=0$\hspace{2em}
$\overline{\beta}\cdot\overline{\gamma}= 0$\\[.5ex]
\>$\alpha:\gamma$\hspace{2em}$\beta :\alpha$\hspace{2em}$\gamma :\overline{\beta}$\\[1ex]

III \;Feriso\quad $\alpha\cdot\alpha\neq 0\;\; .\;\;\alpha\cdot\beta=0\;\; .\;\;\alpha\cdot\gamma=0 \;\supset\; \overline{\beta}\cdot\overline{\gamma}\neq 0$\\[.5ex]

\hspace{5.5em}$\f{\f{\alpha\cdot\alpha\neq 0\hspace{1.6em}\alpha\cdot\beta=0}{\alpha\cdot\overline{\beta}\neq 0}\hspace{1.2em}\afrac{\alpha\cdot\gamma=0}}
{\overline{\beta}\cdot\overline{\gamma}\neq 0}$
\end{tabbing}

$\mathbf{\llbracket 102. \rrbracket}$ In general it can be shown that every correct formula expressed by the Aristotelian terms a, e, i, o and operations of the calculus of proposi\-tions can be derived from this principle; to be more exact, fundamental notions a,~i
\begin{tabbing}
def\hspace{2em}\=$\alpha\,{\rm e}\,\beta\equiv\;\sim(\alpha\,{\rm i}\,\beta)$\\[.5ex]
\>$\alpha\,{\rm o}\,\beta\equiv\;\sim(\alpha\,{\rm a}\,\beta)$\\[1ex]
\>1. \hspace{.5em}\=$\alpha\,{\rm a}\,\alpha$\hspace{2em}Identity\\[.5ex]
\>2. \>$\alpha\,{\rm a}\,\beta\: .\: \beta\,{\rm a}\,\gamma\supset\alpha\,{\rm a}\,\gamma$\hspace{2em}I\hspace{.5em}Barbara\index{Barbara} \\[.5ex]
\>3. \> $\alpha\,{\rm i}\,\beta\: .\: \beta\,{\rm a}\,\gamma\supset\gamma\,{\rm i}\,\alpha$\hspace{2em} IV\hspace{.5em}Dimatis\index{Dimatis}
\end{tabbing}
and all axioms of the propositional calculus; then if we have a formula composed only of such expressions $\alpha\,{\rm a}\,\beta$, $\alpha\,{\rm i}\,\gamma$ and $\sim,\vee\ldots$ and which is universally true, i.e.\ holds for all classes $\alpha,\beta,\gamma$ involved, then it is derivable from these axioms\ by rule of substitution and implication and defined symbol. $\mathbf{\llbracket 103. \rrbracket}$ I am sorry I have no time to give the proof.

So we can say that the Aristotelian theory of syllogisms for expressions of this particular type a, e, i, o is complete, i.e.\ every true formula follows from the Aristotelian moods. But those Aristotelian moods are even abundant because those two moods alone are already sufficient to obtain everything else. The incompleteness of the Aristotelian theory\index{failure of traditional logic} lies in this that there are many $\mathbf{\llbracket 104. \rrbracket}$ propositions which cannot be expressed in terms of the Aristotelian primitive terms. E.g.\ all formulas which I wrote down for $+,\cdot, -$ (distributive law, De Morgan law etc.)\ because those symbols $+,\cdot, -$ do not occur in Aristotle. But there are even simpler things not expressible in Aristotelian terms; e.g.\ $\overline{a}\cdot\overline{c}=0$ (some not $a$ are not $c$), e.g.\ ${\alpha\,{\rm e}\,\beta \atop \underline{\beta\, {\rm o}\,\gamma}}$ according to Aristotle\index{Aristotle} there is no conclusion from that (there is a principle that from two negative premises no conclusion can be drawn)
$\mathbf{\llbracket 105. \rrbracket}$ and that is true if we take account only of propositions expressible by the a, e, i, o. But there is a conclusion to be drawn from that, namely ``Some not $\alpha$ are not $\gamma$'' $\overline{\alpha}\cdot\overline{\gamma}\neq 0$. Since some $\beta$ are not $\gamma$ and every $\beta$ is not $\alpha$ we have some not $\alpha$ (namely the $\beta$) are not $\gamma$. The relation which holds between two classes $\alpha,\gamma$ if $\overline{\alpha}\cdot\overline{\gamma}\neq 0$ cannot be expressed by a, e, i, o, but it is arbitrary to exclude that relation. Another example
\begin{tabbing}
\hspace{1.7em}\= $\alpha\, {\rm i}\,\beta$\\[.5ex]
\>$\alpha\,{\rm o}\,\beta$\\[-2ex]
\>\underline{\hspace{2.5em}}\\
\> $\alpha$ contains at least two elements
\end{tabbing}
$\mathbf{\llbracket 106. \rrbracket}$ Such propositions: ``There are two different objects to which the predicate $\alpha$ belongs'' can of course not be expressed by a, e, i, o, but they can in the logistic calculus by
\begin{tabbing}
\hspace{1.7em}$(\exists x, y)[x\neq y\: .\: x\,\varepsilon\,\alpha \: .\: y\,\varepsilon\,\alpha]$. \end{tabbing}

\subsection{Relations}\label{2.9}
\pagestyle{myheadings}\markboth{EDITED TEXT}{NOTEBOOK VI \;\;---\;\; 1.2.9\;\; Relations}

$\mathbf{\llbracket 107. \rrbracket}$ Last time I developed in outline the calculus of classes in which we introduced certain operations +, $\cdot$, $-$ which obey laws similar to those of arithmetic. One can develop a similar calculus for relations. First of all we can introduce for relations operations +, $\cdot$, $-$ in a manner perfectly analogous to the calculus of classes.

$\mathbf{\llbracket 108. \rrbracket}$ If $R$ and $S$ are any two dyadic relations I put\index{relations, operations on}\index{operations on relations}
\begin{tabbing}
\hspace{1.7em}\= $R +S\;$ \= $=\;\; \hat{x}\hat{y}[xRy\vee xSy]$\\[.5ex]
\>$R\cdot S$ \> $=\;\; \hat{x}\hat{y}[xRy\; .\; xSy]$\\[.5ex]
\>$-R$\> $=\;\; \hat{x}\hat{y}[\sim xRy]$\quad p.\ 110\\[.5ex]
\>$R-S$\> $=\;\; \hat{x}\hat{y}[xRy\; .\: \sim xSy]$
\end{tabbing}

So e.g.\ if $R$ is the relation of father, $S$ the relation of mother one has for the relation of parent:
\begin{tabbing}
\hspace{1.7em}\= parent \= = \hspace{.2em} \= father + mother\\[.5ex]
\> $x$ is a parent of $y$ $\equiv$ $x$ is a father of $y$ $\vee$ $x$ is a mother of $y$\\[.5ex]
\> $\leq$\> = \>$(< + =)$\\[.5ex]
\> child\> = \>son + daughter
\end{tabbing}

$\mathbf{\llbracket 109. \rrbracket}$ Or consider similarity for polygons and the relation of same size and the relation of congruence, then Congruence = Similarity $\cdot$ Same size, or consider the four relations parallelism, without common points, coplanar, and skew, then we have
\begin{tabbing}
\hspace{1.7em}\=Parallelism = without common point $\cdot$ coplanar,\\[.5ex]
\>or Parallelism = without common point $\cdot$ $-$ skew
\end{tabbing}
or $-$brother will subsist between two objects $x,y$ if 1. $x,y$ are two human beings and $x$ is not a brother of $y$ or 2. if $x$ or $y$ is not a human being because $x$ brother $y$ is true only if $x$ and $y$ are human beings and in addition $x$ is a brother of $y$. So if $x$ or $y$ are not human beings the relation \emph{eo ipso} will not $\mathbf{\llbracket 110. \rrbracket}$ hold, i.e.\ the relation $-$brother will hold. Exactly as for classes there will exist also a vacuous and a universal relation denoted by $\dot{\Lambda}$ and $\dot{\rm{V}}$. $\dot{\Lambda}$ is the relation which subsists between no objects $(x,y)\sim x \dot{\Lambda} y$, and $(x,y)x \dot{\rm{V}} y$, e.g.\ \begin{tabbing}
\hspace{1.7em}\=greater $\cdot$ smaller = $\dot{\Lambda}$\\[.5ex]
\>greater + (not greater) = $\dot{\rm{V}}$
\end{tabbing}
Also there exists an analogon to the notion of subclass, namely $R\subseteq S$ if $xRy \supset xSy$, e.g.\
\begin{tabbing}
\hspace{1.7em}\=father \hspace{1.2ex}\=$\subseteq$ ancestor\\[.5ex]
\> brother \>$\subseteq$ relative\\[.5ex]
\> smaller \>$\subseteq$ not greater
\end{tabbing}

These operations for relations (i.e.\ +, $\cdot$, $-$) are exactly analogous to the corresponding for classes and therefore will obey the same laws, e.g.\ $(R+S)\cdot T = R\cdot T + S\cdot T$. But in addition to them there are certain operations specific for relations and therefore more interesting, e.g.\ for any relation $R$ we can form what is called the \emph{inverse of $R$}\index{inverse of a relation} (denoted by $\breve{R}$ or $R^{-1}$ ) where $\breve{R}=\hat{x}\hat{y}[yRx]$, hence $x\breve{R}y \equiv yRx$, i.e.\ if $y$ has the relation $R$ to $x$ then $x$ has the relation $\breve{R}$ $\mathbf{\llbracket 111. \rrbracket}$ to $y$, e.g.\
\begin{tabbing}
\hspace{1.7em}\=child = (parent)$^{-1}$\\[.5ex]
\> $x$ child $y$ $\equiv$ $y$ parent $x$\\[.5ex]
\> $< \;= (>)^{-1}$ \\[.5ex]
\> smaller = (greater)$^{-1}$\\[.5ex]
\> (nephew + niece) = (uncle + aunt)$^{-1}$
\end{tabbing}
There are also relations which are identical with their inverse, i.e.\ $xRy\equiv yRx$. Such relations are called symmetric.\index{symmetric relation}\index{relation, symmetric} Other example (brother + sister) is symmetric because --\ldots; brother is not symmetric, sister isn't either.

$\mathbf{\llbracket 112. \rrbracket}$ Another operation specific for relations and particularly important is the so called relative product\index{relative product of relations}\index{composition of relations} of two relations rendered by $R|S$ and defined by
\[
R|S=\hat{x}\hat{y}[(\exists z)(xRz\; .\; zSy)]
\]
i.e.\ $R|S$ subsists between $x$ and $y$ if there is some object $z$ to which $x$ has the relation $R$ and which has the relation $S$ to $y$, e.g.\
\begin{tabbing}
\hspace{1.7em}\=nephew = son$|$(brother or sister)
\end{tabbing}
$\mathbf{\llbracket 113. \rrbracket}$ $x$ is a nephew to $y$ if $x$ is son of some person $z$ which is brother or sister of $y$. In everyday language the proposition $xRy$ is usually expressed by $x$ is an $R$ of $y$ or $x$ is the $R$ of $y$. Using this we can say $xR|Sy$ means $x$ is an $R$ of an $S$ of $y$, e.g.\ $x$ is a nephew of $y$ means $x$ is a son of a brother or sister of $y$. Other example:\footnote{A note inserted in the manuscript at this example mentions a continuation on p.\ \textbf{119}.}
\begin{tabbing}
\hspace{1.7em}\=paternal uncle = brother$|$father
\end{tabbing}
The relative product can also be applied to a relation and the same relation again, i.e.\ we can form $R|R$ (by def= $R^2$) square of a relation, $\mathbf{\llbracket 114. \rrbracket}$ e.g.\
\begin{tabbing}
\hspace{1.7em} \= paternal grandfather = (father)$^2$\\[.5ex]
\> grandchild = (child)$^2$
\end{tabbing}
Similarly we can form $(R|R)|R= R^3$, e.g.\footnote{A note inserted in the manuscript at this example mentions a continuation on p.\ \textbf{117}.}
\begin{tabbing}
\hspace{1.7em} \= great grandchild = (child)$^3$
\end{tabbing}

The relative product again follows laws very similar to the arithmetic ones, e.g.\
\begin{tabbing}
\hspace{1.7em} \= Associativity: \hspace{.3em} \=$(R|S)|T = R|(S|T)$\\[.5ex]
\> Distributivity:\>$R|(S+T)=R|S +R|T$\\[.5ex]
\> also\> $R|(S\cdot T)\subseteq R|S \cdot R|T$\\[1ex]
but not commutativity\\[1ex]
\> $R|S=S|R$ is false\\[.5ex]
\> brother$|$father $\neq$ father$|$brother
\end{tabbing}

$\mathbf{\llbracket 115. \rrbracket}$\footnote{The whole of pages \textbf{115}.\ and \textbf{116}.\ are crossed out, but the beginning of the present page, p.\ \textbf{115}., is given here because it completes naturally what was said before about the relative product, i.e.\ composition, of relations.} Identity\index{identity relation} $I$ is a unity for this product, i.e.\ $R|I=I|R= R$ because
\begin{tabbing}
\hspace{1.7em} \= $xR|Iy\;$\=$\equiv xIz\; .\; zRy$ for some $z$\\[.5ex]
\>\>$\equiv xRy$\\[1ex]
\> Monotonicity: $R\subseteq S, P\subseteq R\supset R|P\subseteq S|Q$
\end{tabbing}

$\mathbf{\llbracket 117. \rrbracket}$\footnote{see the preceding footnote, at the beginning of p.\ \textbf{115}.} A relation $R$ is called transitive\index{transitive relation}\index{relation, transitive} if
\begin{tabbing}
\hspace{1.7em} $(x,y,z)[xRy\; .\; yRz\supset xRz]\equiv \;  R$ is transitive
\end{tabbing}
In other words if an $R$ of an $R$ of $z$ is an $R$ of $z$; e.g.\ brother is transitive, a brother of a brother of a person is a brother of this person, in other words
\begin{tabbing}
\hspace{1.7em} $x$ brother $y\; .\; y$ brother $z \supset x$ brother $z$
\end{tabbing}
Smaller is also transitive, i.e.\
\begin{tabbing}
\hspace{1.7em} $x<y\; .\; y<z \supset x<z$
\end{tabbing}
Very many relations in mathematics are transitive: congruence, parallelism, isomorphism, ancestor. Son is not transitive, a son of a son of a person is not a son of a person.

$\mathbf{\llbracket 118. \rrbracket}$ Therefore called intransitive; friend is an example of a relation which is neither transitive nor intransitive. A friend of a friend of $x$ is not always a friend of $x$, but is sometimes a friend of $x$. By means of the previously introduced operation transitivity can be expressed by
\begin{tabbing}
\hspace{1.7em}\= $R^2 \subseteq R$ \quad because \\[.5ex]
\> $ xR^2 y \:.\!\supset (\exists z)(xRz\; .\; zRy) \supset xRy$
\end{tabbing}
if $R$ is transitive, but also vice versa if $R$ satisfies the condition $R^2 \subseteq R$ then $R$ is transitive
\begin{tabbing}
\hspace{1.7em}\= $xRy\; .\; yRz \supset xR^2z \supset xRz$
\end{tabbing}

$\mathbf{\llbracket 119. \rrbracket}$ A very important property of relations is the following one: A binary relation $R$ is called one-many\index{one-many relations} if for any object $y$ there exists at most one object $x$ such that $xRy$:
\begin{tabbing}
\hspace{1.7em}$(x,y,z)[xRy\; .\; zRy\supset x=z]\equiv \; R$ is one-many
\end{tabbing}
and many-one if $R^{-1}$ is one-many; e.g.\ father is one-many, every object $x$ can have at most one father, it can have no father if it is no man, but it never has two or more fathers. The relation $<$ is not one-many: for any number there are many different numbers smaller than it.

The \footnote{This sentence and the beginning of the next one, until the end of p.\ \textbf{119}., are crossed out in the manuscript, though the remainder of the paragraph on p.\ \textbf{120}.\ is not.} relation $x$ is the reciprocal of number $y$ is one-many. Every number has at most $\mathbf{\llbracket 120. \rrbracket}$ one reciprocal. Some numbers have no reciprocal, namely 0 (but that makes no difference). The relation of reciprocal is at the same time many-one;\index{many-one relations} such relations are called one-one.\index{one-one relations}

The relation of husband in Christian countries e.g.\ is an example of a one-one relation. The relation smaller is neither one-many nor many-one; for any number there exist many different numbers smaller than it and many different numbers greater than it.

One-many-ness can also be defined for polyadic relations. $\mathbf{\llbracket 121. \rrbracket}$ A triadic relation $M$ is called one-many if
\begin{tabbing}
\hspace{1.7em}$(x,y,z,u)[xM(zu)\; .\; yM(zu)\supset x=y]$
\end{tabbing}
e.g.\ $\hat{x}\hat{y}\hat{z}(x=y+z)$, $\hat{x}\hat{y}\hat{z}[x = \frac{y}{z}]$ have this property. For any two numbers $y$ and $z$ there exists at most one $x$ which is the sum or difference. $\hat{x}\hat{y}(x$ is a square root of $y$) is not one-many because there are in general two different numbers which are square roots of $y$. You see the one-many relations\index{one-many relations} are exactly the same thing which is called ``functions''\index{functions} in mathematics. The dyadic one-many relations are the functions with one argument as e.g.\ $x^2$, the $\mathbf{\llbracket 122. \rrbracket}$ triadic one-many relations are the functions with two arguments as e.g.\ $x+y$.

In order to make statements about functions, i.e.\ one-many relations it is very convenient to introduce a notation usual in mathematics and also in everyday language; namely $R`x$ denotes the $y$ which has the relation $R$ to $x$, i.e.\ the $y$ such that $yRx$ provided that this $y$ exists and is unique. Similarly for a triadic relation $M`(yz)$ denotes the $x$ such that\ldots\ Instead of this also $yMz$ is written, e.g.\ + denotes a triadic relation between $\mathbf{\llbracket 123. \rrbracket}$ numbers (sum) and $y+z$ denotes the number which has this triadic relation to $y$ and $z$ provided that it exists.  In everyday language the $`$ is expressed by the words The\ldots\ of, e.g.\ The sum of $x$ and $y$, The father of $y$.\index{definite descriptions}

There\footnote{The text in this paragraph, until the end of p.\ \textbf{125}., is crossed out in the manuscript, but because of its interest it is given here.} is only one tricky point in this notation. Namely what meaning are we to assign to propositions containing this symbol $R`x$ if there does not exist a unique $y$ such that $yRx$ (i.e.\ none or several), e.g.\ The present king of $\mathbf{\llbracket 124. \rrbracket}$ France is bald. We may convene that such propositions are meaningless (neither true nor false). But that has certain undesirable consequences, namely whether or not the present king of France exists or not is an empirical question. Therefore it would depend on an empirical fact whether or not this sequence of words is a meaningful statement or nonsense whereas one should expect that it can depend only on the grammar of the language concerned whether something makes sense. $\mathbf{\llbracket 125. \rrbracket}$ Russell\index{Russell, Bertrand} makes the convention that such statements are false and not meaningless. The convention is: That every \emph{atomic proposition} in which such an $R`x$ (describing something nonexistent) occurs is false, i.e.\
\begin{tabbing}
\hspace{1.7em}$\varphi(R`x)\equiv (\exists y)[(z)[zRx\equiv\ z=y]\; .\; \varphi(y)]$
\end{tabbing}

\subsection{Type theory and paradoxes}\label{2.10}
\pagestyle{myheadings}\markboth{EDITED TEXT}{NOTEBOOK VI \;\;---\;\; 1.2.10\;\; Type theory and paradoxes}

$\mathbf{\llbracket 126. \rrbracket}$ All aforementioned notions of the calculus of classes and relations are themselves relations; e.g.\ $\alpha\subseteq\beta$ is a binary relation between classes, $\alpha + \beta$ is a dyadic function, i.e.\ a triadic relation between classes (which subsists between $\alpha,\beta,\gamma$ if $\gamma=\alpha+\beta$). The operation of inverse is a relation between relations subsisting between $R$ and $S$ if $R=S^{-1}$ or the relative product is a triadic relation between relations subsisting between $R,S,T$ if $R=S|T$. Symmetry defines a certain class of relations (the class of symmetric relations). So we see that we have obtained a $\mathbf{\llbracket 127. \rrbracket}$ new kind of concepts (called concepts of second type or second order)\index{concepts of second order}\index{concepts of second type} which refer to the concepts of first order, i.e.\ which expresses properties of concepts of first order or relations between concepts of first order or to be more exact properties and relations of extensions of concepts of first order. But this is not very essential since we can define corresponding concepts which express properties and relations of the predicates themselves, e.g.\ $\chi$ sum of $\varphi,\psi$ if $\chi(x)\equiv \varphi(x)\vee\psi(x)$ etc.

And it is possible to (go on) continue in this way, i.e.\ we can define concepts of third type or order,\index{concepts of third order}\index{concepts of third type} which refer to the concepts of second order. An example would be: ``mutually exclusive''; a class of classes $\cal U$, i.e.\ a class whose elements are themselves classes, is called a mutually exclusive class of classes if $\alpha,\beta\:\varepsilon\:{\cal U} \supset \alpha\cdot\beta = \Lambda$. This concept of ``mutually exclusive class of classes'' expresses a property of classes of classes, i.e.\ of an object of third order, therefore is of third order. So you see in this way we get a whole hierarchy of concepts $\mathbf{\llbracket 128. \rrbracket}$ which is called the hierarchy of types.\index{hierarchy of types} In fact there are two different hierarchies of types, namely the hierarchy of extensions and the hierarchy of predicates.\index{hierarchy of extensions}\index{hierarchy of predicated}

An interesting example of predicates of higher type are the natural numbers. According to Russell\index{Russell, Bertrand} and Frege\index{Frege, Gottlob} the natural numbers are properties of predicates. If I say e.g.: There are eight planets, this expresses a property of the predicate $\mathbf{\llbracket 129. \rrbracket}$ ``planet''. So the number 8 can be defined to be a property of predicates which belongs to a predicate $\varphi$ if there are exactly 8 objects falling under this predicate. If this definition is followed up it turns out that all notions of arithmetic can be defined in terms of logical notions and that the laws of arithmetic can be derived from the laws of logic except for one thing, namely for building up arithmetic one needs the proposition that there are infinitely many objects, which cannot be proved from the axioms of logic.

$\mathbf{\llbracket 130. \rrbracket}$ The lowest layer in the hierarchy of types\index{hierarchy of types} described are the individuals or objects of the world; what these individuals are is an extralogical question which depends on the theory of the world which we assume; in a materialist theory it would be the atoms or the points of space and time, in a spiritualist theory it would be the spirits and so on. As to the higher types (classes, classes of classes, predicates of predicates etc.) each type must be distinguished very carefully from any other as can be shown e.g.\ by the following $\mathbf{\llbracket 131. \rrbracket}$ example. If $a$ is an object one can form the class whose only element is $a$ (denoted by $\riota\, ` a$). So this $\riota\, ` a$ would be the extension of a predicate, which belongs to $a$ and only to $a$. It is near at hand to identify this $a$ and $\riota\, ` a$, i.e.\ to assume that the object $a$ and the class whose only element is $a$ are the same. However it can be shown that this is not admissible, i.e.\ it would lead to contradictions to $\mathbf{\llbracket 132. \rrbracket}$ assume this identity $\riota\, ` a=a$ to be generally true because if we take for $x$ a class (which has several elements) then certainly $\riota\, ` \alpha$ and $\alpha$ are distinct from each other; since $\riota\, ` \alpha$ is a class which has only one element, namely $\alpha$, whereas $\alpha$ is a class which has several elements, so they are certainly distinct from each other. But although we have to distinguish very carefully between the different types there is on the other hand a very close analogy between the different types. E.g.\ classes of individuals $\mathbf{\llbracket 133. \rrbracket}$ and classes of classes of individuals will obey exactly the same laws. For both of them we can define an addition and a multiplication and the same laws of calculus will hold for them. Therefore it is desirable not to formulate these laws separately for classes of classes and classes of individuals, but to introduce a general notion of a class comprising in it all those particular cases: classes of individuals, classes of relations, classes of classes etc. And it was actually in $\mathbf{\llbracket 134. \rrbracket}$ this way that the logistic calculus was first set up (with such a general notion of a class and of a predicate and of a relation and so on embracing under it all types) and this way also corresponds more to natural thinking. In ordinary language e.g.\ we have such a general notion of a class without a distinction of different types.

The more detailed working out of logic on this typeless base has led to the discovery of the most interesting $\mathbf{\llbracket 135. \rrbracket}$ facts in modern logic. Namely to the fact that the evidences of natural thinking are not consistent with themselves, i.e.\ lead to contradictions which are called ``logical paradoxes''. The first of these contradictions was discovered by the mathematician Burali-Forti\index{Burali-Forti, Cesare} in 1897. A few years later Russell\index{Russell, Bertrand} produced a similar contradiction which however avoided the unessential mathematical by-work of Burali-Forti's\index{Burali-Forti, Cesare} contradiction and showed the real logical structure of the contradiction much clearer. This so $\mathbf{\llbracket 136. \rrbracket}$ called \emph{Russell paradox}\index{Russell paradox}\index{Russell paradox} has remained up to now the classical example of a logical paradox and I want to explain it now in detail. I shall first enumerate some apparently evident propositions from which the paradox follows in a few steps.

The paradox under consideration involves only the following notions:
\begin{itemize}
\item[1.] object in the most general sense, which embraces everything that can be made an object of thinking; in particular it embraces the individuals, classes, predicates of all types
\end{itemize}

$\mathbf{\llbracket 137. \rrbracket}$\footnote{A note in the manuscript at the bottom of the preceding page, p.\ \textbf{136}., and at the top of this page, seem to suggest that pp.\ \textbf{137}.-\textbf{140}. of the present Notebook VI are to be superseded by pages in Notebook VII starting with p.\ \textbf{137}., the first numbered page in Notebook VII. These four pages of Notebook VI are nevertheless given here.}

\begin{itemize}
\item[2.] monadic predicate (briefly predicate), also in the most general sense comprising predicates of individuals as well as predicates of predicates etc. And this term predicate is to be so understood that it is an essential requirement of a predicate that it is well-defined for any object whatsoever whether the given predicate belongs to it or not
\end{itemize}

Now of these two notions ``object'' and ``predicate'' we have the following apparently evident propositions:
\begin{itemize}
\item[1.] \emph{If $\varphi$ is a predicate and $x$ an object then it is uniquely determined whether $\varphi$ belongs to $x$ or not}.
\end{itemize}
Let us denote the proposition $\varphi$ belongs to $x$ by $\varphi(x)$. So we have if $\varphi$ is a well-defined predicate and $x$ an object then $\varphi(x)$ is always a meaningful proposition $\mathbf{\llbracket 138. \rrbracket}$ which is either true or false. \begin{itemize}
\item[2.] Vice versa: If we have a combination of words or symbols $A(x)$ which contains the letter $x$ and is such that it becomes a meaningful proposition for any arbitrary object which you substitute for $x$ then $A(x)$ defines a certain predicate $\varphi$ which belongs to an object $x$ if and only if $A(x)$ is true. \end{itemize}
So the assumption means that if you substitute for $x$ the name of an arbitrary object then it is always uniquely determined whether the resulting proposition is true or false.
\begin{itemize}
\item[3.] It is uniquely determined of any object whether or not it is a predi\-cate.
\end{itemize}
Let us denote by $P(x)$ the proposition ``$x$ is a predicate'' so that $P$(red), $\sim P$(smaller), $\sim P$(New York); then by 3 $P(x)$ is always a meaningful proposition whatever $x$ $\mathbf{\llbracket 139. \rrbracket}$ may be.
\begin{itemize}
\item[4.] \emph{Any predicate is an object}.
\end{itemize}

I think these four propositions are all evident to natural thinking. 1 and 2 can be considered as a definition of the term predicate and 3 says that the notion of predicate thus defined is well-defined.

And now let us consider the following statement ${P(x)\: .\: \sim x(x)}$ that means $x$ is a predicate and it belongs to $x$ (i.e.\ to itself). According to our four assumptions that is a meaningful proposition which is either true or false whatever you substitute for $x$. Namely, at first by 3 it is uniquely defined: if you $\mathbf{\llbracket 140. \rrbracket}$ substitute for $x$ something which is not a predicate it becomes false, if you substitute for $x$ a predicate then $P(x)$ is true but $x(y)$ is either true or false for any object $y$ written over $x$ by 1. But $x$ is a predicate, hence an object by assumption 4, hence $x(x)$ is either true or false, hence the whole statement is always meaningful, i.e.\ either true or false. Therefore by 2 it defines a certain predicate $\Phi$ such that $\Phi(x)\,{\equiv \atop \mbox{\scriptsize\rm means}}\, P(x)\: .\: \sim x(x)$.

\vspace{2ex}

\noindent $\mathbf{\llbracket Notebook\; VII \rrbracket}$
\pagestyle{myheadings}\markboth{EDITED TEXT}{NOTEBOOK VII \;\;---\;\; 1.2.10\;\; Type theory and paradoxes}
$\mathbf{\llbracket 137. \rrbracket}$\footnote{see the footnote at the top of p.\ \textbf{137}. of the preceding Notebook VI. Notebook VII starts with nine, not numbered, pages of remarks and questions mostly theological, partly unreadable, partly in shorthand, and all seemingly not closely related to the remaining notes for the course. They are rendered as far as possible in the source version and deleted here.}

\begin{itemize}
\item[2.] The notion of a ``well-defined monadic predicate''.
\end{itemize}
That is a monadic predicate $\varphi$ such that for any object $x$ whatsoever it is uniquely determined by the definition of $\varphi$ whether or not $\varphi$ belongs to $x$, so that for any arbitrary object $x$ $\varphi(x)$ is a meaningful proposition which is either true or false. Since I need no other kind of predicate in the subsequent considerations but only well-defined monadic predicates, I shall use the term ``\emph{predicate}'' in the sense of monadic well-defined predicate.
\begin{itemize}
\item[3.] The concept which is expressed by the word ``is'' or ``belongs'' in ordinary language and which we expressed by $\varphi(x)$, which means the predicate $\varphi$ belongs to $x$. \end{itemize}

Now for these notions (of object and predicate) we have the following apparently evident propositions:

\vspace{2ex}

\noindent $\mathbf{\llbracket 138. \rrbracket}$
\begin{itemize}
\item[1.] For any object $x$ it is uniquely determined whether or not it is a predicate; in other words well-defined predicate is itself a well-defined predicate.

\item[2.] If $y$ is a predicate and $x$ an object then it is well-defined whether the predicate $y$ belongs to $x$. This is an immediate consequence of the definition of a well-defined predicate.
\end{itemize}

Let us denote for any two objects $y, x$ by $y(x)$ the proposition $y$ is a predicate and belongs to $x$. So for any two objects $y,x$ $y(x)$ will be a meaningful proposition of which it is uniquely determined whether it is true or false, namely if $y$ is no predicate it is false whatever $x$ may be, if it is a predicate then it is true or false according as the predicate $y$ belongs to $x$ or does not belong to $x$, which is uniquely determined.

\vspace{2ex}

\noindent $\mathbf{\llbracket 139. \rrbracket}$
\begin{itemize}
\item[3.] If we have a combination of symbols or words containing the letter $x$ (denote it by $A(x)$) and if this combination is such that it becomes a meaningful proposition whatever object you substitute for $x$ then $A(x)$ defines a certain well-defined predicate $\varphi$ which belongs to an object $x$ if and only if $A(x)$ is true.
\end{itemize}
(I repeat the hypothesis of this statement: It is as follows, that if you substitute for $x$ the name of an arbitrary object then the resulting expression is always a meaningful proposition of which it is uniquely determined whether it is true or false.) Now this statement too could be considered as a consequence of the definition of a well-defined predicate.
\begin{itemize}
\item[4.] Any predicate is an object. That $\mathbf{\llbracket 140. \rrbracket}$ follows because we took the term object in the most general sense according to which anything one can think of is an object.
\end{itemize}

I think these four propositions are all evident to natural thinking. But nevertheless they lead to contradictions, namely in the following way. Consider the expression $\sim x(x)$ that is an expression involving the variable $x$ and such that for any object substituted for this variable $x$ you do obtain a meaningful proposition of which it is uniquely determined whether it is true or false. $\mathbf{\llbracket 141. \rrbracket}$ Namely if $x$ is not a predicate this becomes false by the above definition of $y(x)$; if $x$ is a predicate then by 1 for any object $y$ it is uniquely determined whether $x$ belongs to $y$, hence also for $x$ it is uniquely determined because $x$ is a predicate, hence an object (by 4). $\sim x(x)$ means $x$ is a predicate not belonging to itself. It is easy to name predicates which do belong to themselves, e.g.\ the predicate ``predicate''; we have the concept ``predicate'' is a predicate. Most of the predicates of course do not belong to themselves. Say e.g.\ the predicate man is not a man, $\mathbf{\llbracket 142. \rrbracket}$ so it does not belong to itself. But e.g.\ the predicate not man does belong to itself since the predicate not man is certainly not a man, so it is a not man, i.e.\ belongs to itself.

Now since $\sim x(x)$ is either true or false for any object $x$ it defines a certain predicate by 3. Call this well-defined predicate $\Phi$, so that $\Phi(x)\equiv\;\sim x(x)$. For $\Phi$ even a term in ordinary language was introduced, namely the word ``impredicable'', and for the negation of it the word ``predicable''; so an object is called predicable if it $\mathbf{\llbracket 143. \rrbracket}$ is a predicate belonging to itself and impredicable in the opposite case, i.e.\ if it is either not a predicate or is a predicate and does not belong to itself. So predicate is predicable, not man is predicable, man is impredicable, Socrates is impredicable.

And now we ask is the predicate ``impredicable'' predi\-cable or impredicable. Now we know this equivalence holds for any object $x$ (it is the definition of impredicable); $\Phi$ is a predicate, hence an object, hence this equivalence holds for $\Phi$, i.e.\ $\Phi(\Phi)\equiv\;\sim \Phi(\Phi)$. What does $\Phi(\Phi)$ say? Since $\Phi$ means impredicable it says impredicable is impredicable. So we see that this proposition is equivalent with its own negation.

$\mathbf{\llbracket 144. \rrbracket}$ But from that it follows that it must be both true and false, because we can conclude from this equivalence:
\begin{tabbing}
\hspace{13.5em}\=$\Phi(\Phi)\supset\; \sim \Phi(\Phi)$\\[.5ex]
\>$\sim \Phi(\Phi) \supset \Phi(\Phi)$
\end{tabbing}
By the first implication, $\Phi(\Phi)$ cannot be true, because the assumption that it is true leads to the conclusion that it is false, i.e.\ it leads to a contradiction; but $\Phi(\Phi)$ cannot be false either because by the second implication the assumption that it is false leads to the conclusion that it is true., i.e.\ again to a contradiction. So this $\Phi(\Phi)$ would be a proposition which is neither true nor false, hence it would be both true and false $\mathbf{\llbracket 145. \rrbracket}$ because that it is not true implies that it is false and that it is not false implies that it is true. So we apparently have discovered a proposition which is both true and false, which is impossible by the law of contradiction.

The same argument can be given without logical symbols in the following form. The question is: Is the predicate ``impredicable'' predicable or impredicable. 1.\ If impredicable were predicable that would mean that it belongs to itself, i.e.\ then impredicable is impredicable. So from the assumption that impredicable is predi\-cable we derived that it is impredicable; so it is not predicable. 2.\ On the other hand assume impredicable is impredicable; then it belongs to itself, hence is predicable. So from the assumption that it is impredicable we derived that it is predicable. So it is certainly not impredicable. So it is neither predicable nor impredicable. But then it must be both predicable and impredicable because since it is not predicable it is impredicable and since it is not impredicable it is predicable. So again we have a proposition which is both true and false.

Now what are we to do about this situation? One may first try to say: Well, the law of contradiction is an error. There do exist such strange things as propositions which are both true and false. But this way out of the difficulty is evidently not possible $\mathbf{\llbracket 146. \rrbracket}$ because that would imply that every proposition whatsoever is both true and false. We had in the calculus of propositions the formula $p\; .\sim p\supset q$ for any $p,q$, hence also $p\; .\sim p\supset \; \sim q$ where $p$ and $q$ are arbitrary propositions. So if we have one proposition $p$ which is both true and false then any proposition $q$ has the undesirable property of being both true and false, which would make any thinking completely meaningless. So we have to conclude that we arrived at this contradictory conclusion
\[
\Phi(\Phi) \quad \rm{and} \quad \sim \Phi(\Phi)
\]
$\mathbf{\llbracket 147. \rrbracket}$ by some error or fallacy, and the question is what does this error consist in [i.e.\ which one of our evident propositions is wrong].

The nearest at hand conjecture about this error is that there is some circular fallacy hidden in this argument, because we are speaking of predicates belonging to themselves or not belonging to themselves. One may say that it is meaningless to apply a predicate to itself. I don't think that this is the correct solution. For the following reasons:
\begin{itemize}
\item[1.] It is not possible to except for any predicate $P$ $\mathbf{\llbracket 148. \rrbracket}$ just this predicate $P$ itself from the things to which it can be applied
\end{itemize}
i.e.\ we cannot modify the assumption 1.\ by saying the property $\varphi(x)$ is well-defined for any $x$ except $\varphi$ itself because if you define e.g.\ a predicate $\mu$ by two predicates $\varphi,\psi$ by $\mu(x)\equiv \varphi(x)\; .\; \psi(x)$ then we would have already three predicates $\mu$, $\varphi$ and $\psi$ to which $\mu$ cannot be applied:
\begin{tabbing}
\hspace{1.7em} $\mu(\varphi)\;{\equiv\atop\rm{Df}}\;\varphi(\varphi)\; .\; \psi(\varphi)$\quad where this makes no sense. \end{tabbing}
$\mathbf{\llbracket 149. \rrbracket}$ So it is certainly not sufficient to exclude just self-reflexivity\index{self-reflexivity}\index{self-reference} of a predicate because that entails automatically that we have to exclude also other things and it is very difficult and leads to very undesirable results if one tries to formulate what is to be excluded on the basis of this idea to avoid self-reflexivities. That was done by Russell\index{Russell, Bertrand} in his so called ramified theory of types which since has been abandoned by practically all logicians. On the other hand it is not even justified to exclude self-reflexivities of every formula because self-reflexivity does not always lead to contradiction but is perfectly legitimate in many cases. If e.g.\ I say: ``Any sentence of the English language contains a verb'' then it is perfectly alright to apply this proposition to itself and to conclude from it that also this proposition under consideration contains a verb.

Therefore the real fallacy seems to lie $\mathbf{\llbracket 150. \rrbracket}$ in something else than the self-reflexivity, namely in these notions of object and predicate in the most general sense embracing objects of all logical types. The Russell paradox seems to show that there does not exist such a concept of everything. As we saw the logical objects form a hierarchy of types and however far you may proceed in the construction of these types you will always be able to continue the process still farther and therefore it is illegitimate and makes no sense to speak of the totality of all objects.

$\mathbf{\llbracket 151. \rrbracket}$ One might think that one could obtain the totality of all objects in the following way: take first the individuals and call them objects of type 0, then take the concepts of type 1, then the concepts of type 2, 3 etc. for any natural number. But it is by no means true that we obtain in this manner the totality of all concepts, because e.g.\ the concept of the totality of concepts thus obtained for all integers $n$ as types is itself a concept not occurring in this totality, i.e.\ it is a concept of a type higher than $\mathbf{\llbracket 152. \rrbracket}$ any finite number, i.e.\
of an infinite type. It is denoted as a concept of type $\omega$. But even with this type $\omega$ we are by no means at an end, because we can define e.g.\ relations between concepts of
type $\omega$ and they would be of a still higher type $\omega + 1$. So we see there are in a sense much more than infinitely many logical types; there are so many that it is not possible to form a concept of the totality of all of them, because whichever concept we form we can define a concept of a higher type, hence not falling under $\mathbf{\llbracket 153. \rrbracket}$ the given concept.

So if we want to take account of this fundamental fact of logic that there does not exist a concept of the totality of all objects whatsoever, we must drop the words ``object'', ``predicate'', ``everything'' from our language and replace them by the words: object of a given type, predicate of a given type, everything which belongs to a given type. In particular, proposition 4 has now to be formulated like this. If $A(x)$ is an expression which becomes a
meaningful proposition for any object $x$ of a given type $\alpha$ then it defines a concept of type $\alpha +1$. We cannot even
formulate the proposition in its previous form, because we don't have such words as object, predicate etc. in our
language. Then the Russell paradox disappears immediately because
we can form the concept $\Phi$ defined by $\Phi(x)\equiv\;\sim x(x)$ only for $x$'s of a given type $\alpha$, i.e.\ $\mathbf{\llbracket 154. \rrbracket}$ we can define a concept $\Phi$ such that this equivalence holds for every $x$ of type $\alpha$. (We cannot even formulate that it holds for every object because we have dropped these words from our language). But then $\Phi$ will be a concept of next higher type because it is a property of objects of type $\alpha$. Therefore we
cannot substitute $\Phi$ here for $x$ because this equivalence holds only for objects of type $\alpha$.

So this seems to me to be the (satisfactory) true solution of the $\mathbf{\llbracket 155. \rrbracket}$ Russell paradox. I only wish to mention that the hierarchy of types as I sketched it here is
considerably more general than it was when it was first presented by its inventor B.\
Russell.\index{Russell, Bertrand} Russell's\index{Russell, Bertrand} theory of types was given in two different forms, the so called
simplified and the ramified theory of types, both of which are much more
restrictive then the one I explained here; e.g.\ in both of them it would be
impossible to form concepts of type $\omega$, also the statement $x(x)$ would always be
meaningless. Russell's theory\index{Russell, Bertrand} of $\mathbf{\llbracket 156. \rrbracket}$ types is more based on the first idea of solving the paradoxes (namely to exclude self-reflexivities) and the totality of all objects is only excluded because it would be self-reflexive (since it would itself be
an object). However the development of axioms of set theory has shown that
Russell's\index{Russell, Bertrand} system is too restrictive, i.e.\ it excludes many arguments which (as far as one can see) do not lead to contradictions and which are necessary for building up
abstract set theory.

There are other logical paradoxes which are solved by the theory of types, i.e.\ by
excluding the terms object, every etc. But there are others in which the fallacy is of
an entirely different nature. They are the so called epistemological paradoxes. $\mathbf{\llbracket 157. \rrbracket}$ The oldest of them is the Epimenides.\index{Epimenides} In the form it is usually
presented, it is no paradox. But if a man says ``I am lying now''\index{liar paradox} and says nothing else,
or if he says: The proposition which I am pronouncing right now is false, then this
statement can be proved to be both true and false, because this proposition $p$ says
that $p$ is false; so we have $p\equiv (p$ is false), $p\equiv\;\sim p$, from which it follows that $p$ is both true and false as we saw before. The same paradox can be brought to a much more
conclusive form as follows:\footnote{Here the text in the manuscript is interrupted and the subsequent numbering of pages in the present Notebook VII starts anew from \textbf{1}.}

\subsection{Examples and samples of previous subjects}\label{2.11}
\pagestyle{myheadings}\markboth{EDITED TEXT}{NOTEBOOK VII \;\;---\;\; 1.2.11\;\; Examples and samples of previous\ldots}

$\mathbf{\llbracket 1. \rrbracket}$ All four rules are purely formal, i.e.\ for applying them it is not necessary to know the meaning of the expressions. Examples of derivations from the axioms. Since all axioms and rules of the calculus of propositions are also axioms and rules of the calculus of functions we are justified in assuming all formulas and rules formerly derived in the calculus of propositions.

\vspace{1ex}

1. Example\footnote{see 13$'$.\ on p.\ \textbf{80}.\ of Notebook~V} $\varphi(y)\supset (\exists x)\varphi(x)$
\begin{tabbing}
Derivation:\\[.5ex]
\hspace{1.7em}\= (1)\hspace{1em}\= $(x)[\sim\varphi(x)]\supset\;\sim\varphi(y)$\quad obtained by substituting $\sim\varphi(x)$ for $\varphi(x)$\\
\` in Ax.\ 5\\
\>(2)\>$\varphi(y)\supset\;\sim(x)[\sim\varphi(x)]$\quad by rule of transposition applied to (1)\\[.5ex]
\>(3)\>$\varphi(y)\supset (\exists x)\varphi(x)$\quad by rule of defined symbol from (2)
\end{tabbing}
$\mathbf{\llbracket 2. \rrbracket}$ 2. Example\footnote{see 6.\ on p.\ \textbf{76}.\ of Notebook~V} $(x)[\varphi(x)\supset\psi(x)]\supset[(x)\varphi(x)\supset(x)\psi(x)]$
\begin{tabbing}
\hspace{1.7em}\= (1)\hspace{1em}\= $(x)[\varphi(x)\supset\psi(x)]\supset[\varphi(y)\supset\psi(y)]$\quad by substituting $\varphi(x)\supset\psi(x)$\\
\` for $\varphi(x)$ in Ax.\ 5\\
\>(2)\>$(x)\varphi(x)\supset\varphi(y)$\quad Ax.\ 5\\[.5ex]
\>(3)\>$(x)[\varphi(x)\supset\psi(x)]\: .\: (x)\varphi(x)\supset[\varphi(y)\supset\psi(y)]\: .\: \varphi(y)$\quad by rule of\\
\` multiplication of implications applied to (1) and (2)\\[.5ex]
\>(4)\> $[\varphi(y)\supset\psi(y)]\: .\: \varphi(y)\supset\psi(y)$\quad by substituting $\varphi(y)$ for $p$ and $\psi(y)$\\
\` for $q$ in the demonstrable formula $(p\supset q)\: .\: p\supset q$\\[.5ex]
$\mathbf{\llbracket 3. \rrbracket}$\>(5)\>$(x)[\varphi(x)\supset\psi(x)]\: .\: (x)\varphi(x)\supset\psi(y)$\quad by rule of syllogism applied to\\
\` (3) and (4)\\[.5ex]
\>(6)\>$(x)[\varphi(x)\supset\psi(x)]\: .\: (x)\varphi(x)\supset(y)\psi(y)$\` by rule of quantifier from (5)\\[.5ex]
\>(7)\>$(x)[\varphi(x)\supset\psi(x)]\supset[(x)\varphi(x)\supset (y)\psi(y)]$\quad by rule of exportation\\
\` from (6) \\[.5ex]
\>(8)\>$(x)[\varphi(x)\supset\psi(x)]\supset[(x)\varphi(x)\supset (x)\psi(x)]$\` by rule of substitution for\\
\` individual variables
\end{tabbing}

Predicates which belong to no object are called vacuous\index{vacuous predicates} (e.g.\ president of U.S.A. born in South Bend). $S$a$P$ and $S$e$P$ are both true if $S$ is vacuous whatever $P$ may be. $\mathbf{\llbracket 4. \rrbracket}$ All tautologies are true also for vacuous predicates but some of the Aristotelian inferences are not, e.g.\
\begin{tabbing}
\hspace{1.7em}\= $S$a$P\supset S$i$P$\hspace{3em}\= (false if $S$ is vacuous)\\[.5ex]
\> $S$a$P\supset \;\sim(S$e$P)$\> (false $\, ''$ $\, ''$ $\,''$ \hspace{1.55em}$''$\hspace{1.55em}),
\end{tabbing}
the mood Darapti\index{Darapti} $M$a$P\: .\: M$a$S\supset S$i$P$ is false if $M$ is vacuous and if $S,P$ are any two predicates such that $\sim(S$i$P)$.

The totality of all objects to which a monadic predicate $P$ belongs is called the extension\index{extension of a predicate} of $P$ and denoted by $\hat{x}[P(x)]$, so that the characteristic $\mathbf{\llbracket 5. \rrbracket}$ property of the symbol $\hat{x}$ is:
\[
\hat{x}\varphi(x)=\hat{x}\psi(x)\equiv (x)[\varphi(x)\equiv\psi(x)]
\]

Extensions of monadic predicates are called classes\index{classes} (denoted by $\alpha,\beta,\gamma\ldots$). That $y$ belongs to the class $\alpha$ is expressed by $y\varepsilon\alpha$ so that $y\varepsilon\hat{x}\varphi(x)\equiv \varphi(y)$. $\hat{x}$ is applied to arbitrary propositional functions $\Phi(x)$, i.e.\ $\hat{x}\Phi(x)$ means the class of objects satisfying $\Phi(x)$, e.g.\ $\hat{x}[I(x)\: .\: x>7]=$ class of integers greater than seven.
Also for dyadic predicates $Q(xy)$ extensions denoted by $\hat{x}\hat{y}[Q(xy)]$ are introduced, which satisfy the equivalence
\[
\hat{x}\hat{y}[\psi(xy)]=\hat{x}\hat{y}[\chi(xy)]\equiv (x,y)[\psi(xy)\equiv\chi(xy)]
\]

$\mathbf{\llbracket 6. \rrbracket}$ It is usual to call these extensions (not the dyadic predicates themselves) relations.\index{relations} If $\Phi(xy)$ is a propositional function with two variables $\hat{x}\hat{y}\Phi(xy)$ denotes the relation which is defined by $\Phi(xy)$. If $R$ is a relation $xRy$ means that $x$ bears the relation $R$ to $y$ so that
\[
u\{\hat{x}\hat{y}[\varphi(xy)]\}v\equiv\varphi(uv)
\]

The extension of a vacuous predicate is called zero class\index{vacuous class}\index{zero class} and denoted by 0 (or $\Lambda$); the extension of a predicate belonging to every object is called universal class\index{universal class} and denoted by 1 (or V).

$\mathbf{\llbracket 7. \rrbracket}$ For classes operation of $+$, $\cdot$~, $-$ which obey laws similar to the arithmetic laws are introduced by the following definitions:\footnote{see p.\ \textbf{95}.\ of Notebook VI}
\begin{tabbing}
\hspace{1.7em}\= $\alpha +\beta$ \= $=\;\; \hat{x}[x\,\varepsilon\,\alpha\vee x\,\varepsilon\,\beta]$\quad\=(sum)\index{sum of classes}\\[.5ex]
\> $\alpha\cdot\beta$ \> $=\;\; \hat{x}[x\,\varepsilon\,\alpha\: .\: x\,\varepsilon\,\beta]$\> (intersection)\index{multiplication of classes}\index{intersection of classes}\\[.5ex]
\> $-\alpha$ \> $=\;\; \hat{x}[\sim x\,\varepsilon\,\alpha]$\> (complement)\index{opposite of a class}\index{complement of a class}\\[.5ex]
\>$\alpha -\beta$\> $=\;\;\alpha\cdot(-\beta)$\> (difference)\index{difference of classes}
\end{tabbing}

\begin{theindex}\label{Index}
\pagestyle{myheadings}\markboth{EDITED TEXT}{INDEX}

  \item (1), axiom for propositional logic, 31, 34
  \item (2), axiom for propositional logic, 31, 34
  \item (3), axiom for propositional logic, 31, 35
  \item (4), axiom for propositional logic, 31, 35

  \indexspace

  \item a proposition, 70, 97
  \item a, e, i, o propositions, 70, 97
  \item addition from the left, rule of, 39
  \item addition from the right, rule of, 39
  \item addition of classes, 95
  \item adjoining a new hypothesis, rule of, 41
  \item analytical expression, 84
  \item and, 12
  \item Aristotelian figures and moods, 1, 11
  \item Aristotle, 1, 91, 93, 100
  \item associative law, 22
  \item atomic formulas in predicate logic, 68
  \item atomic propositions in predicate logic, 68
  \item \emph{aut}, 13
  \item \emph{aut}\ldots\ \emph{aut}, 3
  \item Ax.\ 5, axiom for predicate logic, 78
  \item axiom (1) for propositional logic, 31, 34
  \item axiom (2) for propositional logic, 31, 34
  \item axiom (3) for propositional logic, 31, 35
  \item axiom (4) for propositional logic, 31, 35
  \item axiom Ax.\ 5 for predicate logic, 78
  \item axiom system for predicate logic, 78
  \item axiom system for propositional logic, 31
  \item axioms for predicate logic, 78
  \item axioms for propositional logic, 31, 34, 78

  \indexspace

  \item Barbara, 73, 98, 100
  \item bell, 86
  \item Bernays, Paul, 58, 59
  \item binary relation, 67
  \item bound variable, 72
  \item brackets (parentheses), 4, 14
  \item Burali-Forti, Cesare, 108

  \indexspace

  \item calculemus, 85, 86
  \item calculus of classes, 93
  \item calculus of propositions, 3, 12
  \item calculus ratiocinator, 85
  \item classes, 93, 94, 116
  \item commutative law, 22
  \item commutativity, rule of, 43
  \item complement of a class, 96, 117
  \item completeness in predicate logic, 74
  \item completeness proof, auxiliary theorem, 51
  \item completeness proof, main theorem, 55
  \item completeness theorem, propositional logic, 50
  \item completeness, interest of, 63
  \item composition of relations, 102
  \item comprehension, 94
  \item concepts of second order, 106
  \item concepts of second type, 106
  \item concepts of third order, 106
  \item concepts of third type, 106
  \item conjunction, 3, 4, 12, 13
  \item conjunction elimination, 44
  \item conjunction introduction, 44
  \item conjunctive normal form, 62
  \item connections (connectives), 3, 4, 12
  \item connectives, 3, 4, 12
  \item contraposition, 21
  \item crank, 86

  \indexspace

  \item Darapti, 91, 99, 116
  \item Darii, 73
  \item De Morgan formulas (laws), 22, 62
  \item De Morgan law for classes, 97
  \item decidability, 24
  \item decision problem, 24
  \item decision problem in predicate logic, 74
  \item decision procedure, 24
  \item defined connectives, 31
  \item definite descriptions, 105
  \item derived rule I, 88
  \item derived rule II, 88
  \item derived rule III, 88
  \item derived rules, 38
  \item derived rules for propositional logic, 38
  \item derived rules of inference, 38
  \item difference of classes, 96, 117
  \item dilemma, 10, 20
  \item dilemma, rule of, 39
  \item Dimatis, 100
  \item disjunction, 3, 4, 13
  \item disjunctive normal form, 62
  \item disjunctive syllogism, 2, 4, 11, 14
  \item distributive laws, 22
  \item dyadic relation (binary relation), 67

  \indexspace

  \item e proposition, 70, 97
  \item Epimenides, 115
  \item equivalence, 4, 13
  \item excluded middle, 21
  \item exclusive ``or'', 13
  \item existential presuppositions, 91
  \item existential quantifier, 69
  \item exportation, 21
  \item exportation, rule of, 43
  \item expression (formula), of propositional logic, 15
  \item extension of a predicate, 93, 116
  \item extensional connections (connectives), 8
  \item extensional functions (connectives), 18
  \item extensionality, 93, 94

  \indexspace

  \item F, Falsehood, 16
  \item failure of traditional logic, 1, 10, 67, 92, 100
  \item Falsehood, 5, 16
  \item Feriso, 98
  \item five axioms for predicate logic, 78
  \item follows, 8, 19
  \item force, of binding for connectives, 14
  \item formula of syllogism, 36
  \item formula, of predicate logic, 78
  \item formula, of propositional logic, 15
  \item four axioms for propositional logic, 31, 34, 78
  \item four rules of inference for predicate logic, 79, 85
  \item free variable, 72
  \item Frege, Gottlob, 106
  \item functional completeness, 25, 29
  \item functional variables, 78
  \item functions, 105
  \item functions of the calculus of propositions (connectives), 13
  \item fundamental conjunction, 51

  \indexspace

  \item generalized rule of syllogism, 39
  \item Gentzen, Gerhard, 64--66

  \indexspace

  \item hierarchy of extensions, 106
  \item hierarchy of predicated, 106
  \item hierarchy of types, 106, 107
  \item Hilbert, David, 3, 31

  \indexspace

  \item i proposition, 70, 98
  \item identity, 71
  \item identity relation, 103
  \item if\ldots\ then, 8, 13, 18
  \item implication, 4, 13
  \item implication rule for sequents, 65
  \item implies, 8, 19
  \item importation, 21
  \item importation, rule of, 43
  \item independence of the axioms, 57
  \item individual variables, 78
  \item intensional connections (connectives), 8
  \item intensional functions (connectives), 18
  \item intersection of classes, 96, 117
  \item inverse of a relation, 102

  \indexspace

  \item Kant, Immanuel, 84

  \indexspace

  \item law of contradiction, 21
  \item law of double negation, 21
  \item law of excluded middle, 21, 37
  \item law of identity, 21, 37
  \item laws of passage, 74, 75
  \item Leibnitz theorema praeclarum, 23
  \item Leibnitz's program, 85, 86
  \item Leibnitz, Gottfried Wilhelm, 85, 86
  \item liar paradox, 115
  \item logic of modalities, 9, 19
  \item logical consequence, 8, 18
  \item logical product (conjunction), 4, 13
  \item logical sum (disjunction), 4, 13
  \item logically true formulas, of propositional logic (tautology), 10,
		20
  \item {\L}ukasiewicz, Jan, 4, 14

  \indexspace

  \item many-one relations, 105
  \item material implication, 8, 18
  \item Menger, Karl, 1, 10, 22
  \item modus ponens, 2, 10, 20, 33, 79, 85
  \item modus ponens for sequents, 65
  \item molecular propositions, 69
  \item monadic predicate, 68
  \item multiplication of classes, 96, 117

  \indexspace

  \item negation, 3, 4, 12, 13
  \item not, 12

  \indexspace

  \item o proposition, 70, 98
  \item one-many relations, 104, 105
  \item one-one relations, 105
  \item operation of the calculus of propositions (connectives), 4
  \item operations on relations, 101
  \item opposite of a class, 96, 117
  \item or, 13

  \indexspace

  \item paradoxes, of material implication, 9, 19
  \item parentheses, 4, 14
  \item Polish notation, 4, 14
  \item ponendo ponens, 2, 10, 20
  \item ponendo tollens, 10, 20
  \item predicate, 68
  \item predicate variables, 73, 78
  \item predicate, monadic, 68
  \item product, logical (conjunction), 4, 13
  \item proof by cases, dilemma, 21
  \item properties, 67, 68
  \item propositional calculus, 3, 12
  \item propositional function, 69
  \item propositional variables, 15

  \indexspace

  \item reductio ad absurdum, 21
  \item relation, dyadic (binary), 67
  \item relation, symmetric, 102
  \item relation, transitive, 103
  \item relation, triadic, 67
  \item relations, 67, 94, 117
  \item relations, operations on, 101
  \item relative product of relations, 102
  \item rule 1, 33, 79, 85
  \item rule 2, 79, 85
  \item rule 3, 82, 85
  \item rule 4, 82, 85
  \item rule I, 88
  \item rule II, 88
  \item rule III, 88
  \item rule of addition from the left, 39
  \item rule of addition from the right, 39
  \item rule of addition of premises for sequents, 65
  \item rule of adjoining a new hypothesis, 41
  \item rule of commutativity, 43
  \item rule of defined symbol for predicate logic, 82
  \item rule of defined symbol for propositonal logic, 34
  \item rule of dilemma, 39
  \item rule of exportation, 43
  \item rule of exportation for sequents, 65
  \item rule of implication, 33, 79, 85
  \item rule of implication for sequents, 65
  \item rule of importation, 43
  \item rule of indirect proof for sequents, 65
  \item rule of multiplication, 44
  \item rule of product, 44
  \item rule of product, inversion, 44
  \item rule of reductio ad absurdum for sequents, 65
  \item rule of substitution for predicate logic, 79
  \item rule of substitution for propositional logic, 33
  \item rule of syllogism 4.R., 38
  \item rule of syllogism, generalized, 39
  \item rule of the universal quantifier, 82, 85
  \item rules in logic, 32
  \item rules of inference for predicate logic, 79, 85
  \item rules of inference for propositional logic, 33
  \item rules of inference in logic, 32
  \item rules of transposition, 41
  \item Russell paradox, 108
  \item Russell, Bertrand, 3, 31, 57, 93, 106, 108, 113--115

  \indexspace

  \item schematic letters, 31
  \item scope, 72
  \item scope of a quantifier, 72
  \item secondary formula (sequent), 64
  \item self-reference, 113
  \item self-reflexivity, 113
  \item sequent, 64
  \item Sheffer stroke, 26
  \item Sheffer, Henry Maurice, 26
  \item \emph{sive}\ldots\ \emph{sive}, 3
  \item Stoic addition to traditional logic, 2, 11
  \item strength, of binding for connectives, 5
  \item strict implication, 8, 18
  \item subclass, relation of, 96
  \item sum of classes, 95, 117
  \item sum, logical (disjunction), 4, 13
  \item syllogism under an assumption, 47
  \item Syllogism, formula of, 36
  \item symmetric relation, 102

  \indexspace

  \item T, Truth, 16
  \item tautological formula, of propositional logic, 10, 20
  \item tautological vs. analytical, 84
  \item tautology, of predicate logic (valid formula), 73, 76, 77
  \item tautology, of propositional logic, 10, 20
  \item theorem, system of predicate logic, 83
  \item theorem, system of propositional logic, 34
  \item theorema praeclarum, 23
  \item thinking machines, 86
  \item three rules of inference for propositional logic, 33
  \item tollendo ponens, 2, 10, 11, 20
  \item tollendo tollens, 10, 20
  \item traditional logic, 1, 10
  \item transitive relation, 103
  \item transposition (contraposition), 21
  \item transposition rules of inference, 41
  \item triadic relation, 67
  \item Truth, 5, 16
  \item truth functions, 8, 18
  \item truth table, 5, 15
  \item truth values, 5, 16
  \item typewriter, 86

  \indexspace

  \item universal class, 95, 117
  \item universal quantifier, 69
  \item universally true formula (valid formula), 73
  \item universally true formula, of propositional logic (tautology),
		20

  \indexspace

  \item vacuous class, 95, 117
  \item vacuous predicates, 91, 116
  \item valid formula, 73, 76, 77
  \item variable, bound, 72
  \item variable, free, 72
  \item variables, functional, 78
  \item variables, invidual, 78
  \item variables, predicate, 73, 78
  \item variables, propositional, 15
  \item \emph{vel}, 13

  \indexspace

  \item zero class, 95, 117

\end{theindex}

\clearpage

\chapter{SOURCE TEXT}

\setcounter{section}{-1}

\section{Notebook 0}\label{00}
\pagestyle{myheadings}\markboth{SOURCE TEXT}{NOTEBOOK 0}
\zl Folder 58, on the front cover of the notebook ``Vorl.\zl esungen\zd\ Log.\zl ik\zd\ \zl German: Lectures Logic\zd\ N.D.\zl Notre Dame\zd\ 0'' together with some crossed out practically unreadable text in which one can recognize what is presumably: Arb, Beg., Res, Vol, N.D.\zd\

\vspace{1ex}

\zl Before p.\ {\bf 1}.\ one finds on a page not numbered the following apparently incomplete note, which does not seem directly related to the text that follows:

\vspace{1ex}

\noindent $x$ is called D-pair (resp\zl ectively\zd\ D-trip\zl le\zd ) if $z=\langle\;\rangle$ (resp\zl ectively\zd\ $z=\langle\;\rangle$) where the $x,y,z$ are \sout{then} evid.\zl ently\zd\ uniquely det.\zl ermined\zd\ by $z$ \sout{\zl unreadable text\zd }\zd\
\vspace{2ex}

$\mathbf{\llbracket 1. \rrbracket}$ Log\zl ic\zd\ is usually def.\zl ined\zd\ a\zl s\zd\ the science of the laws of \zl presuma\-bly ``corr'', which abbreviates ``correct''; if ``corr'' is read instead as ``con'', then this would abbreviate ``consistent''\zd\ thinking. Accord.\zl ing\zd\ to this def\zl i\-ni\-tion\zd\ the centr.\zl al\zd\ part of log.\zl ic\zd\ must be the theory of inf\zl erence\zd\ and the theory of logically true prop\zl osi\-tions\zd . By a log\zl ically\zd\ true prop.\zl osition\zd\ I mean a prop.\zl osition\zd\ which is true for merely log\zl ical\zd\ reasons as e.g.\ the law of excluded middle\zl ,\zd\ which says that for any prop\zl osition\zd\ $p$ either $p$ or $\sim p$ is true. \ul I intend to go in med\zl ias\zd\ res right away an\zl d\zd\ to begin with this centr.\zl al\zd\ part. \ud

\zl new paragraph\zd\ As Prof\zl essor\zd\ M\zl enger\zd\ has pointed out in his introductory lecture the treatment of these things\zl ,\zd\ \ul inferences and log.\zl ically\zd\ true prop.\zl ositions,\zd\ \ud in traditional logic \sout{and in most of the current textbooks} is unsatisfactory in some resp\zl ect\zd . \zl \sout{1.}\zd\ First with resp.\zl ect\zd\ to completeness. What the $\mathbf{\llbracket 2. \rrbracket}$ trad\zl itional\zd\ logic gives is a more or less arbitrary selection from the infinity of the laws of logic\zl ,\zd\ whereas in a systematic treatment we shall have to develop methods which allow us to obtain \ul as far as possible \ud all logically true prop.\zl ositions\zd\ and \ul \sout{least for cert.\ domains of logic, and furthermore} \ud methods \ul which allow \ud to decide of arbitrary given prop.\zl ositions\zd\ \zl \sout{of}\zd\ \ul \sout{these domains} \ud whether or not they are logically true. But the classical treatment is unsatisfactory also \sout{from} in another respect.\zl ;\zd\ namely as to the question of reducing the laws of logic to a cert.\zl ain\zd\ number of prim.\zl itive\zd\ laws from which $\mathbf{\llbracket 3. \rrbracket}$ all the others can be deduced. Although it is sometimes claimed that everything can be deduced from the law of contradiction or from the first Aristotelian figure\zl ,\zd\ this claim has never been proved or even clearly formulated in traditional logic. \zl dash from the manuscript deleted\zd\

\zl new paragraph\zd\ The chief aim in the first part of this seminary will be to fill these two gaps \ul of trad.\zl itional\zd\ log\zl ic\zd\ \ud\zl ,\zd\ i\zl .\zd e.\ 1.\ to give as far as possible \sout{to give} a complete theory of log.\zl ical\zd\ inf\zl erence\zd\ and of log.\zl ically\zd\ true prop.\zl ositions\zd\ and 2.\ to show how \ul all of them \ud can be deduced from a minimum number of prim.\zl itive\zd\ laws.

$\mathbf{\llbracket 4. \rrbracket}$ The theory of inf\zl erence\zd\ as present.\zl ed\zd\ in the current textbooks is usually divided into two parts\zl :\zd\
\begin{itemize}
\item[1.] The Arist\zl otelian\zd\ figures and moods including the inf.\zl erences\zd\ with one pre\-m.\zl ise,\zd\ i\zl .\zd e.\ conv.\zl ersion,\zd\ contr.\zl aposition\zd\ etc.
\item[2.] Inferences of an entirely different kind\zl ,\zd\ which are treated under the heading of hyp.\zl othetical\zd\ disj.\zl unctive\zd\ conj\zl unctive\zd\ inf.\zl erence,\zd\ and which are a Stoic addition to the Arist.\zl otelian\zd\ figures\zl .\zd\
\end{itemize}

Let us begin with these inf\zl erences\zd\ of the sec.\zl ond\zd\ kind\zl ,\zd\ which turn out to be \sout{much} more fundamental than the Arist\zl otelian\zd\ figures.

Take the following example\zl s\zd\ of the disj.\zl unctive\zd\ inf.\zl erence\zd\ tollendo ponens:

\vspace{1ex}

\noindent $\mathbf{\llbracket 5. \rrbracket}$ From the two premis\zl \sout{s}\zd es
\begin{tabbing}
\hspace{1.7em}\=1. \= Nero was either insane or a criminal\zl ,\zd \\[.3ex]
\>2. \> Nero was not insane\zl ,\zd \\[.5ex]
we can conclude\\[.5ex]
\>\> Nero was a criminal\zl .\zd\
\end{tabbing}
\zl ``Nero'' above, in all three instances, is written almost as ``New''.\zd\
\begin{tabbing}
\hspace{1.7em}\=\ul \= \zl 1.\zd\ \= Today is either Sunday or a holiday\zl ,\zd \\[.3ex]
\>\> \zl 2.\zd\ \> Today is not Sunday\zl ,\zd \\[.5ex]
\>\>\> Today is a holiday\zl .\zd\ \ud
\end{tabbing}

Generally\zl , i\zd f $p,q$ are \ul two \ud arbitrary prop\zl ositions\zd \zl inserted \underline{!!} from the manu\-script deleted\zd\ and we have the two premis\zl \sout{s}\zd es
\begin{tabbing}
\hspace{1.7em}\=1. \= Either $p$ or $q$\zl ,\zd \\[.3ex]
\>2. \> not-$q$\zl not-$p$,\zd \\[.5ex]
we can conclude\\[.5ex]
\>\> $p$\zl $q$.\zd\
\end{tabbing}
It is possible to express this syll\zl ogism\zd\ by one log.\zl ically\zd\ true prop.\zl osition\zd\ as follows:
\begin{tabbing}
\hspace{1.7em},,\zl ``\zd (If either $p$ or $q$ and \sout{if} not-$p$) then $q$'' \zl !! from the manuscript dele-\\
\` ted\zd\
\end{tabbing}
This whole prop.\zl osition\zd\ under quotation marks will be true whatever the prop.\zl o\-sitions\zd\ $p$ and $q$ may be\zl .\zd\

$\mathbf{\llbracket 6. \rrbracket}$ Now what is the caract.\zl er\zd\ of this inf.\zl erence\zd\ which distinguishes them \zl it\zd\ from the Arist.\zl otelian\zd\ figures? It is this that in order to make this inf.\zl erence\zd\ it is not necessary to know anything about the structure of \ul the prop\zl ositions\zd\ \ud $p$ and $q$. $p$ and $q$ may be \zl \sout{may be}\zd\ aff.\zl irmative\zd\ or neg.\zl ative\zd\ prop.\zl ositions,\zd\ they may be simple or complicated\zl ,\zd\ they may themselves be disj.\zl unctive\zd\ or hyp.\zl othetical\zd\ prop.\zl ositions;\zd\ all this is indifferent for this syllogism\zl ,\zd\ i.e.\ only prop\zl ositions\zd\ as a whole occur in it\zl ,\zd\ and it is this \ul caract.\zl er\zd\ \ud that makes this kind of syl\zl logism\zd\ simpler and more fund.\zl amental\zd\ than \ul e.g.\ \ud the Arist\zl otelian\zd\ $\mathbf{\llbracket 7. \rrbracket}$ figures\zl ,\zd\ which depend on the structure of the prop.\zl ositions\zd\ involved. \zl E\zd .g.\ in order to make an inf\zl erence\zd\ by mood Barbara you must know that the two prem.\zl ises\zd\ are universal affirmative. Another example of a log.\zl ical\zd\ law in which only prop\zl ositions\zd\ as a whole occur would be the law of excl.\zl uded\zd\ middle\zl ,\zd\ which says: For any prop\zl osition\zd\ $p$ either $p$ or not-$p$ is true.

\zl dash from the manuscript deleted, and new paragraph introduced\zd\ Now the theory of those laws of logic in which only prop.\zl ositions\zd\ as a whole occur is called calculus of proposition\zl s,\zd\ and it is exclusively with this part of math.\zl ematical\zd\ logic that we shall have $\mathbf{\llbracket 8. \rrbracket}$ to do in the next \ul few \ud lectures\zl .\zd\ \zl dash from the manuscript deleted\zd\ We have to begin with examining in more detail the connections between prop.\zl o\-sitions\zd\ which occur in the inf.\zl erences\zd\ concerned\zl ,\zd\ i\zl .\zd e.\ the or, and, if, not. One has introduced special symbols to denote them. \zl ``N\zd ot\zl ''\zd\ is denoted by a circumflex\zl ,\zd\ \zl ``\zd and\zl ''\zd\ by a dot\zl ,\zd\ ,,\zl ``\zd or'' by a kind of \ul abbrev.\zl ated\zd\ \ud v (derived from vel)\zl ,\zd\ \zl ``\zd if then\zl ''\zd\ is denoted by this symbol similar to a horseshoe \ul \zl !! from the manuscript deleted; it indicated presumably where the following table should be inserted.\zd\ \ud \zl :\zd\
\begin{tabbing}
\ul \zl $p$, $2>1$, $q$ and $3>2$, which are presumably given as examples in the\\
\` manuscript, are here deleted\zd\ \\[.5ex]
\hspace{1.7em}\=not \hspace{5em} \=$\sim$ \quad \ul which is an abbrev\zl iated\zd\ N \ud\quad \= $\sim p$\\[.5ex]
\>and\> $\, .$ \> $p\: .\:q$\\[.5ex]
\>or\> $\vee$\> $p\vee q$\\[.5ex]
\>if\ldots\ then \> $\supset$\> $p\supset q$\\[.5ex]
\>equivalent \> $\equiv$\> $p\equiv q$ \ud
 \end{tabbing}
i.e\zl .\zd\ if $p$ and $q$ are arbitrary prop.\zl ositions\zd\ $\sim p$ m.\zl eans\zd\ $p$ is false\zl ,\zd\ $p\: .\:q$ means both $p$ and $q$ is true\zl ,\zd\ $p\vee q$ means either $p$ or $q$\zl ,\zd\ $p\supset q$ means \zl i\zd f $p$ then $q$\zl ,\zd\ or in other words $p$ implies $q$\zl .\zd\ So if e.g.\ $p$ is the prop\zl osition\zd today it will rain and $q$ \zl is\zd\ $\mathbf{\llbracket 9. \rrbracket}$ the prop.\zl osition\zd tomorrow it will snow then \zl text in the manuscript broken\zd\

\zl A\zd bout the \zl ``\zd or\zl '':\zd\ namely\zl ,\zd\ this log\zl ical\zd\ symb.\zl ol\zd\ means that at least one of the two prop.\zl ositions\zd\ $p,q$ is true but does not exclude the case where both are true\zl ,\zd\ \ul i\zl .\zd e\zl .\zd\ it means one or both of them are true\zl ,\zd\ \ud \sout{i e it corresponds to the latin vel} whereas the \zl ``\zd or\zl ''\zd\ in trad.\zl itional\zd\ logic is the exclusive \zl ``\zd or\zl ''\zd\ which \sout{corresp.\ to the latin aut and} means that exactly one of the two prop\zl ositions\zd\ $p,q$ is true and the other one false. \ul Take e.g.\ the sentence \zl ``\zd Anybody who has a salary or interests \ul from cap\zl ital\zd\ \ud is liable to income tax\zl ''\zd . Here the \zl ``\zd or\zl ''\zd\ is meant in the sense of the log\zl ical\zd\ \zl ``\zd or\zl '',\zd\ because someone who has both is also liable to income tax\zl .\zd\ On the other hand\zl ,\zd\ in the prop.\zl osition\zd \zl ``A\zd ny number \zl minus written over another sign; should be: except\zd\ $1$ is either greater or smaller \zl than\zd\ 1\zl ''\zd\ we mean the excl.\zl usive\zd \zl ``\zd or\zl ''\zd . This excl.\zl usive\zd\ \zl ``\zd or\zl ''\zd\ corresp.\zl onds\zd\ to the \zl L\zd at.\zl in\zd\ aut\zl \emph{aut},\zd\ the log.\zl ical\zd\ \zl ``\zd or\zl ''\zd\ to the \zl L\zd at.\zl in\zd\ vel\zl \emph{vel}\zd . As we shall see later\zl .\zd\ \ud

The excl.\zl usive\zd\ ,,\zl ``\zd or'' can be expressed by a comb.\zl ination\zd\ $\mathbf{\llbracket 10. \rrbracket}$ of the other logical symb.\zl ols,\zd\ but no special symbol \ul has been \ud introduced for it, because it is not very often used. Finally\zl ,\zd\ I introduce a fifth connection\zl ,\zd\ \ul the so\zl $\:$\zd called \ud ,,\zl ``\zd equivalence'' denoted by three horiz.\zl on\-tal\zd\ lines. $p\equiv q$ means that both $p$ implies $q$ and $q$ implies $p$. This relation of equivalence would hold e.g.\ between the two prop\zl ositions: \zl ``\zd T\zl \sout{w}\zd omorrow is a weekday\zl ''\zd\ and \zl ``T\sout{w}\zd omorrow is not \ul a \ud holiday\zl ''\zd \zl full stop added here, which in the manuscript is followed by the words: ``because we have --- \ul If\ldots\ but also vice versa \ud ''\zd\

The five notions which we have introduced so far are called resp.\zl ectively\zd\ \ul operation of \ud neg\zl ation\zd , conj\zl unction\zd , disj\zl unction\zd , implic.\zl ation,\zd\ equivalence. By a common name they are called f\zl u\zd nct\zl ions\zd\ of the calc.\zl u\-lus\zd\ of prop.\zl ositions\zd\ \ul or \zl missing text, full stop from the manuscript deleted\zd\ Disj\zl unc\-tion\zd\ is also called $\mathbf{\llbracket 11. \rrbracket}$ log.\zl ical\zd\ sum and conj.\zl unction\zd\ log.\zl ical\zd\ prod.\zl uct\zd\ because of cert\zl ain\zd\ analogies with the arithmetic sum and the ar.\zl ithmetic\zd\ prod\zl uct\zd . A prop\zl osition\zd\ of the form ${p\vee q}$ is called a dis\-j.\zl unction\zd\ \sout{or a logical sum} and $p,q$ its first and sec.\zl ond\zd\ member\zl ;\zd\ similarly a prop\zl osition\zd\ of the form $p\supset q$ is called an impl\zl ication\zd\ and $p,q$ its first and sec.\zl ond\zd\ member\zl ,\zd\ and similarly for the other op\zl erations\zd . Of course\zl ,\zd\ if $p,q$ are prop.\zl osi\-tions,\zd\ then \underline{$\sim p$\zl ,\zd\ $\sim q$\zl ,\zd\ $p\vee q$\zl ,\zd\ $p\: .\:q$\zl ,\zd\ $p\supset q$} \zl underlining omitted in the edited version\zd\ are also prop.\zl osi\-tions\zd\ and therefore to them the functions of the calc.\zl ulus\zd\ of prop\zl ositions\zd\ can again be applied so as to get more complicated expr\zl es\-sions;\zd\ e.g.\ \underline{$p\vee(q\, .\:r)$}\zl underlining omitted in the edited version\zd \zl ,\zd\ which would mean: Either $p$ is true or $q$ and $r$ are both true.

\zl new paragraph\zd\ The disj.\zl unctive\zd\ syllogism $\mathbf{\llbracket 12. \rrbracket}$ I mentioned before can be expressed in our symbolism as follows: \underline{$[(p\vee q)\, .\:\sim q]\supset p$}\zl underlining omitted in the edited version\zd \zl .\zd\ You see in more complicated expressions as e.g.\ this one brackets have to be used exactly as in algebra to indicate in what order the op.\zl erations\zd\ have to be carried out. If e.g.\ I put the brackets in a diff.\zl erent\zd\ way in this expr.\zl ession,\zd\ namely like this \underline{$(p\vee q)\, .\:r$}\zl underlining omitted in the edited version\zd \zl ,\zd\ it would mean something entirely diff.\zl erent,\zd\ namely \ul it would mean \ud either $p$ or $ q$ is true and in addition $r$ is true.

\zl new paragraph\zd\ There is an interesting remark due to the Polish log.\zl i\-cian\zd\ L\zl {\L}\zd ukasiewicz\zl ,\zd\ namely that one can dispense entirely with brackets if one writes the $\mathbf{\llbracket 13. \rrbracket}$ \zl \sout{the}\zd\ operational symb.\zl ols\zd\ $\vee$, $\supset$ etc\zl .\zd\ always in front of the prop\zl osition\zd\ to which they are applied\zl ,\zd\ e.g.\ \underline{$\supset p\,q$}\zl un\-der\-lining omitted in the edited version\zd\ instead of \underline{$p\supset q$}\zl underlining omitted in the edited version\zd . \ul Inc.\zl identally,\zd\ the word \zl ``\zd if\zl ''\zd\ \ul of ordinary lang\zl uage\zd\ \ud is used in exactly this way. We say e.g.\ \zl ``\zd If it is possible I shall \zl do it\zd \zl ''\zd\ putting the \zl ``\zd if\zl ''\zd\ in front of the \ul two \ud prop\zl ositions\zd\ to which we apply it. \ud \ul Now \ud in this notation \ul where the op.\zl erations\zd\ are put in front \ud the two diff.\zl erent\zd\ possibilities of this expression $p\vee q \lfloor .\rfloor r$ would be dist\zl inguished\zd\ automatically without the use of brackets because the sec.\zl ond\zd\ would read \underline{$.\vee p\,q\,r$}\zl underlining omitted in the edited version\zd \zl , with\zd\ ,,\zl ``\zd or'' appl\zl ied\zd\ to $p,q$ and the ,,\zl ``\zd and'' applied to this form.\zl ula\zd\ and $r$\zl ,\zd\ whereas the first would read ,,\zl ``\zd and'' applied to $q$, $r$ and the $\vee$ applied to $p$ and this form\zl ula\zd\ \underline{$\vee p \, .\: qr$}\zl underlining omitted in the edited version\zd \zl .\zd\ \ul As you see\zl ,\zd\ \ud t\zl written over T\zd hese two form\zl ulas\zd\ differ from each other without the use of brackets and it can be shown that $\mathbf{\llbracket 14. \rrbracket}$ it is quite generally so. Since however the formulas in the bracket notation are more easily readable I shall keep the brackets and put the operat.\zl ion\zd\ symb.\zl ol\zd\ \sout{in} between the prop.\zl ositions\zd\ to which they are applied.

\zl new paragraph\zd\ You know in algebra one can save many brackets by the conv.\zl ention\zd\ that multipl\zl ication\zd\ is of greater force than addition\zl ,\zd\ and one can do something similar here by stipulating an order of force between the op.\zl erations\zd\ of the calc.\zl ulus\zd\ of prop.\zl ositions,\zd\ and this order is to be exactly the same in which I introduced them\zl ,\zd\ namely
\[
\sim\, .\:\vee\;{\supset \atop \equiv}
\]
\zl N\zd o order of force is def.\zl ined\zd\ for $\supset\,\equiv$\zl ,\zd\ they are to have equal force. Hence
\begin{tabbing}
$\mathbf{\llbracket 15. \rrbracket}$\\*[1ex]
\hspace{1.7em}\= $\sim p\vee q$ \quad me\=ans \quad \=$(\sim p)\vee q$ \hspace{1.1em} n\=ot \quad \=$\sim(p\vee q)$\\*[.5ex]
\>$p\: .\: q\vee r$ \>$''$ \>$(p\: .\: q)\vee r$ \>$''$ \>$p\: .\: (q\vee r)$\\*[.5ex]
\ul exactly as for arith.\zl metical\zd\ sum and prod.\zl uct\zd\ \ud \\[.5ex]
\>$p\vee q\supset r$ \>$''$ \>$(p\vee q)\supset r$ \>$''$ \>$p\vee (q\supset r)$\\[.5ex]
\>$\sim p\supset q$ \>$''$ \>$(\sim p)\supset q$ \>$''$ \>$\sim(p\supset q)$\\[.5ex]
\>$\sim p\: .\: q$ \>$''$ \>$(\sim p)\, .\: q$ \>$''$ \>$\sim(p\: .\: q)$\\[.5ex]
\>$\sim p\equiv q$ \>$''$ \>$(\sim p)\equiv q$ \>$''$ \>$\sim(p\equiv q)$
\end{tabbing}
\ul \zl I\zd n all these cases the expr\zl ession\zd\ written without brackets has the meaning of the prop\zl osition\zd\ in the sec.\zl ond\zd\ col\zl umn\zd . If we have the form\zl ula\zd\ of the 3\zl third\zd\ col\zl umn\zd\ in mind we have to write the brackets. \ud

Another conv\zl ention\zd\ used in arithm.\zl etic\zd\ for saving brack.\zl ets\zd\ is this that inst\zl ead\zd\ of $(a+b)+c$ we can write $a+b+c$. We make the same conventions for log.\zl ical\zd\ addition and mult.\zl iplication,\zd\ i\zl .\zd e\zl .\zd\ $p\vee q\vee r$ mean\zl s\zd\ $(p\vee q)\vee r$\zl ,\zd\ $p\: .\: q\, .\: r$ \zl means\zd\ $(p\: .\: q)\, .\: r$\zl .\zd\ \ul \zl

\zl new paragraph\zd\ T\zd he letters $p,q,r$ which den.\zl ote\zd\ arb.\zl itrary\zd\ prop.\zl o\-si\-tions\zd\ are called prop.\zl ositional\zd\ variables\zl ,\zd\ and any expression composed of prop.\zl ositional\zd\ var.\zl iables\zd\ and the oper.\zl ations\zd\ $\sim$\zl ,\zd\ $\vee$\zl ,\zd\ $\, .\:$\zl ,\zd\ $\supset$\zl ,\zd\ $\equiv$\ is called meaningful expression or formula of the calc.\zl ulus\zd\ of prop.\zl ositions,\zd\ where also the letters $p,q$ themselves are considered as the simplest kind of expressions\zl .\zd\

After those merely symbolic conventions the next thing we have to do is to examine in more detail the meaning of the op.\zl erations\zd\ of the calc.\zl ulus\zd\ of prop\zl ositions\zd . Take e.g.\ the disj.\zl unction\zd\ $\vee$\zl .\zd\ If $\mathbf{\llbracket 16. \rrbracket}$ any two prop.\zl osi\-tions\zd\ $p,q$ are given $p\vee q$ will again be a prop\zl osition\zd . But now (and this is the decisive point) this op.\zl eration\zd\ of \zl ``\zd or\zl ''\zd\ is such that the truth or falsehood of the composit\zl e\zd\ prop.\zl osition\zd\ $p\vee q$ depends in a def.\zl inite\zd\ way on the truth or falsehood of the const.\zl ituents\zd\ $p,q$. This dependence can be expressed most clearly in \zl the\zd\ form of a table as follows: Let us form three col.\zl umns,\zd\ one headed by $p$\zl ,\zd\ one by by $q$\zl ,\zd\ one by $p\vee q$\zl ,\zd\ and let us write T for true and F for false. Then for the prop\zl ositions\zd\ $p,q$ we have the foll.\zl owing\zd\ four possibilities \zl dots pointing in the manuscript to the following tables deleted\zd\
\begin{center}
\begin{tabular}{ c|c|c c|c }
 $p$ & $q$ & $p\vee q$ & $p\circ q$ & $p\: .\: q$\\[.5ex]
 T & T & T & F & T \\
 T & F & T & T & F\\
 F & T & T & T & F\\
 F & F & F & F & F
 \end{tabular}
\end{center}
Now for each of these 4\zl four\zd\ cases we can easily determine $\mathbf{\llbracket 17. \rrbracket}$ \ul whether \ud $p\vee q$ will be true or false\zl ;\zd\ namely\zl ,\zd\ since $p\vee q$ means that one or both of the prop\zl ositions\zd\ $p$ \zl \sout{$\vee$} ,\zd\ $q$ are true it will be true in the first\zl ,\zd\ sec.\zl ond\zd\ and third case\zl ,\zd\ and false only in the fourth case.\ul We can consider this table (called the truth\zl $\,$\zd table \ul for $\vee$ \ud) as the most precise def.\zl inition\zd\ of what $\vee$ means. \ud

\zl new paragraph\zd\ It is usual to call truth and falsehood the truth values and to say of a true prop.\zl osition\zd\ that it has the truth value ,,\zl ``\zd Truth''\zl ,\zd\ and of a false prop.\zl osition\zd\ that it has the truth value ,,\zl ``\zd Falsehood''\zl .\zd\ T and F then denote the truth values and the \sout{this table called the} truth table \ul for $\vee$ \ud shows how the truth value of the composit\zl e\zd\ expr\zl ession\zd\ \ul $p\vee q$ \ud depends on the truth values of the constituents. The exclusive \zl ``\zd or\zl ''\zd\ would have another truth $\mathbf{\llbracket 18. \rrbracket}$ table\zl ;\zd\ namely if I denote it by $\circ$ for the moment, we have $p\circ q$ is false in the case when both $p$ and $q$ are true and in the case when both $''$ $''$ $''$ \zl $p$ and $q$\zd\ are false\zl ,\zd\ and it is true in the other cases, where one of the two prop.\zl ositions\zd\ $p,q$ is true and the other one is false. The op.\zl eration\zd\ $\sim$ has the following truth\zl $\,$\zd table
\begin{center}
\begin{tabular}{ c|c}
 $p$ & $\sim p$ \\[.5ex]
 T & F \\
 F & T
 \end{tabular}
\end{center}
Here we have only two poss.\zl ibilities:\zd\ $p$ is true and $p$ is false\zl ,\zd\ and if $p$ is true not-$p$ is false and if $p$ is false not-$p$ is true. The truth\zl $\,$\zd table for ,,\zl ``\zd and'' can also easily be determined\zl :\zd\ $p\: .\: q$ is true only in the case where \sout{$p$} both $p$ and $q$ are true and false in all the other three cases.

\zl new paragraph\zd\ A little more $\mathbf{\llbracket 19. \rrbracket}$ difficult is the question of the truth\zl $\,$\zd\ table for $\supset$. $p\supset q$ was defined to mean: If $p$ is true then $q$ is also true. So \ul in order to determine the truth\zl $\,$\zd table \ud let us assume that for two given prop\zl ositions\zd\ $p,q$ $p\supset q$ holds\zl ,\zd\ i.e\zl .\zd\ let us assume \uline{we know \zl ``\zd If $p$ then $q$\zl ''\zd\ but nothing else} \zl underlining replaced partially in the edited version by italics\zd \zl ,\zd\ and let us ask what can \ul we conclude about \ud the truth values of $p$ and $q$ from this assumption. \zl It is not indicated in the manuscript where the following table should be inserted. The text in the manuscript that follows it is a comment upon it. In this table the first three lines in the columns beneath $p$ and $q$ are put in a box, which in the edited text is printed separately in the next display, further down.\zd\

\vspace{1.5ex}

\hspace{1em}\begin{tabular}{c c|c|c|c|c }
Ass\zl umption\zd\ & $p\supset q$ & $p$ & $q$ & $\sim p$ & $\sim p\vee q$\\[.5ex]
& T & F & T & T & T\\
& T & F & F & T & T\\
& T & T & T & F & T\\
& F & T & F & F & F
\end{tabular}

\vspace{1.5ex}

\noindent First it may certainly happen that $p$ is false\zl ,\zd\ bec.\zl ause\zd\ the as\-s.\zl ump\-tion\zd\ ,,\zl ``\zd If $p$ then $q$'' says nothing about the truth or falsehood of $p$\zl ,\zd\ and in this case when $p$ is false $q$ may be true as well as false\zl ,\zd\ because the ass.\zl umption\zd\ says nothing about what happens to $q$ if $p$ is false\zl ,\zd\ but only if $p$ is true\zl .\zd\ $\mathbf{\llbracket 20. \rrbracket}$ So we have both these poss.\zl ibilities:\zd\ $p\;$F$\;\;q\;$T\zl ,\zd\
$p\;$F$\;\;q\;$F. Next we have the poss.\zl ibility\zd\ that $p$ is true\zl ,\zd\ but in this case $q$ must also be true owing to the ass.\zl umption;\zd\ so that the poss.\zl ibility\zd\ $p$ true $q$ false is excluded and it is the only of the four possibilities that is excluded by the ass.\zl umption\zd\ $p\supset q$. It follows that either one of those three possib.\zl ilities,\zd\ \zl \sout{(}\zd\ which I frame in \zl
\begin{center}
\begin{tabular}{ |c|c| }
$p$ & $q$\\[.5ex]
\hline
F & T\\
F & F\\
T & T\\
\hline
\end{tabular}
\end{center}
\zd \zl \sout{)}\zd\ occurs. But we have also vice versa: If one of these three possib\zl ilities\zd\ for the truth\zl $\,$\zd val.\zl ue\zd\ of $p$ and $q$ is realized then $p\supset q$ holds. For let us assume we know that one of the three marked $\mathbf{\llbracket 21. \rrbracket}$ cases occurs\zl ;\zd\ then we know also ,,\zl ``\zd If $p$ is true $q$ is true''\zl ,\zd\ because if $p$ is true only the third of the three marked cases can be realized and in this case $q$ is true. So we see that the statement \zl ``\zd If $p$ then $q$\zl ''\zd\ is exactly equivalent with the statement that one of the three marked cases for the truth values of $p$ and $q$ is realized\zl ,\zd\ i.e.\ $p\supset q$ will be true in each of the three marked cases and false in the last case. And this gives the desired truth\zl $\,$\zd table for implication. However there are two important remarks about it\zl ,\zd\ namely\zl :\zd\

1. Exactly the same truth\zl $\,$\zd table can also be $\mathbf{\llbracket 22. \rrbracket}$ obtained by a combination of operations introduced previously\zl ,\zd\ namely $\sim p\vee q$\zl ,\zd\ i\zl .\zd e.\ either $p$ is false or $q$ is true has the same truth table. For $\sim p$ is true whenever $p$ is false\zl ,\zd\ i.e\zl .\zd\ in the first two cases and $\sim p\vee q$ is then true if either $\sim p$ or $q$ is true\zl ,\zd\ and as you see that happens in exactly the cases where $p\supset q$ is true\zl .\zd\ So we see $p\supset q$ and $\sim p\vee q$ are equivalent\zl ,\zd\ i\zl .\zd e\zl .\zd\ whenever $p\supset q$ holds then also $\sim p\vee q$ holds and vice versa. This makes possible to define $p\supset q$ by $\sim p\vee q$ and \ul this \ud is the usual way of introducing the impl.\zl ication \ul in math.\zl ematical\zd\ log\zl ic\zd\ \ud\zl .\zd

\zl new paragraph, 2.\zd\ The sec.\zl ond\zd\ remark about the truth\zl $\,$\zd table for impl.\zl ica\-tion\zd\ is this. We must $\mathbf{\llbracket 23. \rrbracket}$ not forget that $p\supset q$ was understood to mean simply \zl ``\zd If $p$ then $q$\zl ''\zd\ and nothing else\zl ,\zd\ and only this made the constr.\zl uction\zd\ of the truth\zl $\,$\zd table possible. There are other interpretations of the term ,,\zl ``\zd implic.\zl a\-tion\zd '' for which our truth\zl $\,$\zd table would be completely inadequate\zl .\zd\ E.g.\ $p\supset q$ could be given the meaning: $q$ is a log.\zl ical\zd\ consequence of $p$\zl ,\zd\ i\zl .\zd e.\ $q$ can be derived from $p$ by means of a chain of syllogisms. In this sense e.g.\ the prop.\zl osition\zd \zl ``\zd Jup.\zl iter\zd\ is a planet\zl ''\zd\ would imply the prop\zl osition\zd \zl ``\zd Jup.\zl i\-ter\zd\ is not a fix\zl ed\zd\ star\zl ''\zd\ because no planet can be a fix\zl ed\zd\ star by def.\zl inition,\zd\ i\zl .\zd e.\ $\mathbf{\llbracket 24. \rrbracket}$ by merely log\zl ical\zd\ reasons.

\zl new paragraph\zd\ This kind \ul and also some other similar kinds \ud of impl.\zl ica\-tion\zd\ is\zl are\zd\ \sout{usually} called strict impl.\zl ication\zd\ and denoted by this symbol \ul $\prec$ \ud and the implication defined \sout{before} \ul by the truth\zl $\,$\zd table \ud is called material impl.\zl ication\zd\ if it is to be distinguished from $\prec$. Now it is easy to see \sout{not only} that our truth\zl $\,$\zd table would be false for strict impl.\zl ication\zd\ and even more\zl ,\zd\ namely that there exists no truth\zl $\,$\zd table at all for strict implication. In order to prove this consider the first line of our truth table, where $p$ and $q$ are both true and let us ask what will the truth\zl $\,$\zd value of $p\prec q$ be in this case\zl .\zd\ $\mathbf{\llbracket 25. \rrbracket}$ It turns out that this truth\zl $\,$\zd value is not \sout{be} uniquely det\zl ermined\zd . For take e.g.\ for $p$ the prop\zl osi\-tion\zd \zl ``\zd Jup\zl iter\zd\ is a planet\zl ''\zd\ and for $q$ \zl ``\zd Ju.\zl piter\zd\ is not a fix\zl ed\zd\ star\zl '',\zd\ then $p,q$ are both true \ul and \ud $p\prec q$ is also true\zl .\zd\ On the other hand if you take for $p$ again \zl ``\zd Ju.\zl piter\zd\ is a planet\zl ''\zd\ and for $q$ \zl ``\zd France is a republic\zl ''\zd\ then again both $p$ and $q$ are true\zl ,\zd\ but $p\prec q$ is false because \zl ``\zd France is a republic\zl ''\zd\ is not a log.\zl ical\zd\ consequ.\zl ence\zd\ of \zl ``\zd Ju.\zl piter\zd\ is a planet\zl ''\zd . So we see the truth value of $p\prec q$ is not uniquely det.\zl ermined\zd\ by the truth values of $p$ and $q$ and therefore no truth\zl $\,$\zd table exists\zl .\zd\ $\mathbf{\llbracket 26. \rrbracket}$ Such functions of prop.\zl ositions\zd\ for which no truth\zl $\,$\zd table exists are called intensional as opposed to extensional ones for which a truth\zl $\,$\zd table does exist. The ext.\zl ensional\zd\ f\zl u\zd nct\zl ions\zd\ are also called truth\zl $\,$\zd functions, because they depend only on the truth or falsehood of the prop.\zl ositions\zd\ involved\zl .\zd\

So we see logical consequ\zl ence\zd\ is an intensional rel.\zl ation\zd\ \ul betw.\zl een\zd\ prop.\zl ositions\zd\ \ud and \zl \sout{there are}\zd \zl the\zd\ mat\zl erial\zd\ impl.\zl ication\zd\ introd\zl uced\zd\ by our \sout{a} truth\zl $\,$\zd table cannot mean logical consequence\zl .\zd\ Its meaning is best given by the word \zl ``\zd if\zl ''\zd\ of ordinary language which has a much wider sign.\zl ification\zd\ than just log.\zl ical\zd\ cons.\zl equence;\zd\ e.g.\ \zl if\zd\ \ul someone \ud says: \zl ``\zd If I don't come I $\mathbf{\llbracket 27. \rrbracket}$ shall call you\zl ''\zd\ that does not indicate that this telephoning is a log.\zl ical\zd\ consequ.\zl ence\zd\ of \ul his not \ud coming\zl ,\zd\ but it means simply he will either come or telephone\zl ,\zd\ which is exactly the meaning expressed by the truth\zl $\,$\zd table. \ul Hence mat.\zl erial\zd\ implication introduced by the truth table\zl s\zd\ corresponds as closely to \zl ``\zd if then\zl ''\zd\ as a precise notion can correspond to a not precise notion of ordinary language\zl .\zd\:~\ud

\zl dash from the manuscript deleted, and new paragraph introduced\zd\ If we are now confronted with the question which one of the two kinds of impl.\zl ication\zd\ we shall use in developing the theory of inf.\zl erence\zd\ we have to consider two things\zl :\zd\ 1.\ mat.\zl erial\zd\ implication is the much simpler and clearer notion and 2.\ it is quite sufficient for developing the theory of inf.\zl erence\zd\ because in order to conclude $q$ from $p$ it is quite sufficient $\mathbf{\llbracket 28. \rrbracket}$ to know $p$ implies mat\zl erially\zd\ $q$ and not nec.\zl essary\zd\ to know that $p$ impl.\zl ies\zd\ strictly $q$\zl .\zd\ \ul For if we know $p\supset q$ we know that either $p$ is false or $q$ is true. Hence if we know in add.\zl ition\zd\ that $p$ is true the first of the two poss.\zl ibilities\zd\ that $p$ is false is not realized\zl .\zd\ Hence the sec.\zl ond\zd\ must be realized\zl ,\zd\ namely $q$ is true\zl .\zd\ \ud For these two reasons \ul that mat.\zl erial\zd\ impl.\zl ication\zd\ is simpler and sufficient \ud I shall use only mat.\zl erial\zd\ impl.\zl ication\zd\ at least in th\zl e\zd\ \ul first \ud introductory part of my lectures\zl ,\zd\ and shall use the terms ,,\zl ``\zd implies'' and ,,\zl ``\zd follows'' only in the sense \ul of mat\zl erial\zd\ imp.\zl lication\zd\ \ud . I do not want to say by this that a theory of strict impl\zl ication\zd\ may not be interesting and important for cert.\zl ain\zd\ purposes. In fact I hope it will be discussed in the sec\zl ond\zd\ half of this seminary. But this theory bel\zl ongs\zd\ to an entirely diff.\zl erent\zd\ part of logic than the one I am dealing with now\zl ,\zd\ $\mathbf{\llbracket 29. \rrbracket}$ namely to the logic of modalities.

I come now to some apparently parad.\zl oxical\zd\ consequences of our def\zl ini\-tion\zd\ of mat\zl erial\zd\ impl.\zl ication\zd\ whose parad\zl oxicality; one finds however ``paradoxity'' on p.\ \textbf{22}.\ of Notebook~I\zd\ however disappears if we remember that it does not mean log.\zl ical\zd\ consequ\zl ence\zd . The first of these con\-se\-qu.\zl ences\zd\ is that a true prop.\zl osition\zd\ is implied by any prop.\zl osition\zd\ whatsoever. We see this at once from the truth\zl $\,$\zd table which shows that $p\supset q$ is always true if $q$ is true whatever $p$ may be. \ul You see there are only two cases where $q$ is true \zl \sout{namely}\zd\ and in both of them $p\supset q$ is true. \ud But sec.\zl ondly\zd\ we see also that $p\supset q$ is always true if $p$ is false whatever $q$ may be. \zl \sout{\ul bec. you see \ud}\zd\ So that means \sout{that} \zl a\zd\ false propo\zl osition\zd\ implies any prop.\zl osition\zd\ whatsoever\zl ,\zd\ which is the sec\zl ond\zd\ of the paradoxical consequences. These properties of impl.\zl ication\zd\ $\mathbf{\llbracket 30. \rrbracket}$ can also be expressed by saying\zl : ``\zd An implication with true sec.\zl ond\zd\ member is always true whatever the first member may be and an impl.\zl ication\zd\ with false first member is always true whatever the second member may be\zl '';\zd\ we can express that also by formulas like this $q\supset(p\supset q)$\zl ,\zd\ $\sim p\supset(p\supset q)$\zl .\zd\ Both of these form\zl ulas\zd\ are also immediate consequences of the fact that $p\supset q$ is equiv\zl alent\zd\ with $\sim p\vee q$ because what $\sim p\vee q$ says is exactly that either $p$ is false or $q$ is true\zl ;\zd\ so $\sim p\vee q$ will always be true if $p$ is false and \ul will be also true \ud if $q$ is true whatever the other prop\zl osition\zd\ may be. If we apply $\mathbf{\llbracket 31. \rrbracket}$ these formulas to special cases we get strange cons.\zl equences;\zd\ e.g.\ \zl ``\zd J.\zl upiter\zd\ is a fix\zl ed\zd\ star\zl ''\zd\ implies \zl ``\zd France is a republic\zl '',\zd\ but it also implies \zl ``\zd France is not a republic\zl ''\zd\ because a false prop\zl osition\zd\ implies any prop\zl osition\zd\ whatsoever. Similarly \zl ``\zd France is a republic\zl ''\zd\ is implied by \zl ``\zd Ju.\zl piter\zd\ is a planet\zl ''\zd\ but also by \zl ``\zd Ju.\zl piter\zd\ is a fix\zl ed\zd\ star\zl ''\zd . But as I mentioned before these consequ\zl ences\zd\ are \sout{only} paradoxical only for strict impl\zl ication\zd . They are in pretty good agreement with the meaning which the word \zl ``\zd if\zl ''\zd\ has in ord.\zl inary\zd\ langu\zl age\zd . \sout{if the} Because the first formula then says if $q$ is true $q$ is also true if $p$ is true \ul which is not paradoxical but trivial \ud and the sec.\zl ond\zd\ says if $p$ is false then if $p$ is true anything $\mathbf{\llbracket 32. \rrbracket}$ is true. That this is in \ul good \ud agreement with the meaning which the word ,,\zl ``\zd if'' has can be seen from many colloquialisms\zl ;\zd\ e.g\zl .\zd\ if something is obviously false one says sometimes \zl ``I\zd f this is true I am a Chinaman\zl '',\zd\ which is another way of saying \zl ``I\zd f this is true anything is true\zl ''.\zd\ Another of these so called parad.\zl oxical\zd\ cons.\zl equences\zd\ is e.g\zl .\zd\ that for any two arbitrary prop\zl ositions\zd\ one must imply the other\zl ,\zd\ i\zl .\zd e.\ for any $p,q$ $(p\supset q)\vee(q\supset p)$\zl ;\zd\ in fact $q$ must be either true or false\zl ---\zd if it is true the first member of the disj.\zl unction\zd\ is true bec.\zl ause\zd\ it is an impl.\zl ication\zd\ with true sec\zl ond\zd\ member\zl ,\zd\ if it is false the second member of the disj\zl unction\zd\ is $\mathbf{\llbracket 33. \rrbracket}$ true. \ul So this disjunction is always true\zl .\zd\ \ud

\zl new paragraph\zd\ Those three formulas\zl ,\zd\ as well as the form\zl ula\zd\ of disj.\zl unc\-tive\zd\ inf\zl erence\zd\ we had before\zl ,\zd\ are examples of \ul so called \ud universally true formulas\zl ,\zd\ i\zl .\zd e.\ formulas which are true whatever the pro\-p\zl ositions\zd\ $p,q,r$ occurring in them may be. Such form.\zl ulas\zd\ are also called logically true or tautological\zl ,\zd\ and it is exactly the chief aim of the calc.\zl ulus\zd\ of prop.\zl ositions\zd\ to investigate these tautol\zl ogical\zd\ formulas.

\zl new paragraph\zd\ I shall begin with discussing a few more examples before going \ul over \ud to more general considerations\zl .\zd\ I mention at first \ul some of \ud the trad\zl itional\zd\ hyp.\zl othetical\zd\ and $\mathbf{\llbracket 34. \rrbracket}$ disj.\zl unctive\zd\ inferences which in our notation read as follows:
\begin{itemize}
\item[1.] $(p\supset q)\, .\: p\supset q$ \quad pon\zl endo\zd\ pon.\zl ens\zd\ \quad (Assertion)
\vspace{-1ex}
\item[2.] $(p\supset q)\, .\: \sim q\supset\;\sim p$ \quad toll\zl endo\zd\ toll\zl ens\zd\
\vspace{-1ex}
\item[3.] $(p\vee q)\, .\: \sim q\supset p$ \quad toll.\zl endo\zd\ pon.\zl ens\zd\ \quad as we had bef.\zl ore\zd \\
(the mod.\zl us\zd\ pon.\zl endo\zd\ toll\zl ens\zd\ holds only for the exc\zl lusive\zd\ $\vee$)
\vspace{-1ex}
\item[4.] An inf\zl erence\zd\ which is also treated in many of the textbooks under the heading of ,,\zl ``\zd dilemma'' is this\\[1ex]
 $(p\supset r)\, .\:(q\supset r)\supset(p\vee q\supset r)$\\[1ex]
 If both $p\supset r$ and $q\supset r$ then from $p\vee q$ follows $r$. It is usually written as an inf\zl erence\zd\ with three prem.\zl ises,\zd\ $\mathbf{\llbracket 35. \rrbracket}$ namely from the three premis\zl \sout{s}\zd es $(p\supset r)\, .\:(q\supset r)\, .\:(p\vee q)$ one \ul can \ud conclude\zl \sout{s}\zd\ $r$\zl .\zd\
\end{itemize}
\ul This is nothing else but the principle of proof by cases\zl ,\zd\ namely the prem.\zl ises\zd\ say: one of the two cases $p,q$ must occur and from both of them follows $r$\zl .\zd\ That this \ul form\zl ula\zd\ with 3\zl three\zd\ prem\zl ises\zd\ \ud means the same thing as the form\zl ula\zd\ under cons\zl ideration\zd\ is clear because this earlier form\zl ula\zd\ \ul says: \ud \zl ``\zd If the first two prem\zl ises\zd\ are true then if the third is true $r$ is true\zl '',\zd\ which means exactly the same thing as \zl ``\zd If all the three premis\zl \sout{s}\zd es are true $r$ is true\zl .\zd\ The possibility of going over from one of these two form\zl ulas\zd\ to the other is due to another \ul import\zl ant\zd\ \ud log.\zl ical\zd\ principle which is called importation and reads like this
\begin{tabbing}
\hspace{1.7em}$[p\supset(q\supset r)]\supset(p\: .\: q\supset r)$ \quad imp.\zl ortation\zd\
\end{tabbing}
and its inverse which is called exp.\zl ortation\zd\ and reads like this
\begin{tabbing}
\hspace{1.7em}$(p\: .\: q\supset r)\supset[p\supset(q\supset r)]$ \quad exp\zl ortation\zd .
\end{tabbing}
So owing to these two impl\zl ications\zd\ we have also an equiv.\zl alence\zd\ between the left and right\zl -\zd h.\zl and\zd\ side\zl .\zd\ \ud

Next we have the \ul three \ud law\zl s\zd\ of identity\zl ,\zd\ excl\zl uded\zd\ middle and contr.\zl a\-diction\zd\ which read as follows in our not.\zl ation\zd\
\begin{tabbing}
\hspace{1.7em}1. $p\supset p$\quad\quad 2. $p\:\vee\sim p$ \quad\quad 3. $\sim\!(p\: .\sim p)$
\end{tabbing}
\zl W\zd e can add another sim\zl ilar\zd\ law\zl ,\zd\ the law of double neg\zl ation\zd\ which says $\sim\!(\sim p)\equiv p$\zl .\zd\

Next we have the very important formulas of transpos\zl ition\zd :
\begin{tabbing}
\hspace{1.7em}$(p\supset q)\supset(\sim q\supset\: \sim p)$ \quad \zl \sout{if from $p$ foll\zl ows\zd\ $q$ then \dots}\zd\
\end{tabbing}
\zl O\zd ther forms of this form\zl ula\zd\ of trans\zl position\zd\ would be
\begin{tabbing}
\hspace{1.7em}\=$(p\supset\: \sim q)\supset(q\supset\: \sim p)$ \quad \zl \sout{if.}\zd \\
\>$(\sim p\supset q)\supset(\sim q\supset p)$ \quad proved in the same way\zl .\zd\
\end{tabbing}
\zl I\zd n all those formulas of transp\zl osition\zd\ we can write equ.\zl ivalence\zd\ in\-st.\zl ead\zd\ of \zl \sout{identity} the main implication,\zd\ i\zl .\zd e.\ $\mathbf{\llbracket 36. \rrbracket}$ we have also ${(p\supset q)}\equiv(\sim q\supset\: \sim p)$\zl .\zd\ \zl A\zd nother form \ul of transpos\zl ition,\zd\ namely with two prem\-\zl ises,\zd\ is this \ud $(p\: .\: q\supset r)\supset (p\: .\sim r\supset\: \sim q)$ because under the ass.\zl umption\zd\ $p\: .\: q\supset r$ if we know $p\: .\sim r$\zl , then\zd\ $q$ cannot be \ul true \ud because $r$ would be true in this case\zl .\zd\

Next we have diff.\zl erent\zd\ so called red.\zl uctio\zd\ ad abs\zl urdum,\zd\ e.g\zl .\zd\
\begin{tabbing}
\hspace{1.7em}$(p\supset q)\, .\: (p\supset\:\sim q)\supset\:\sim p$
\end{tabbing}
\zl A\zd\ part.\zl icularly\zd\ interest\zl ing\zd \zl \sout{the}\zd\ form of red\zl uctio\zd\ ad abs.\zl urdum\zd\ is the one which Prof.\zl essor\zd\ M.\zl enger\zd\ mentioned in his intr.\zl oductory\zd\ talk and which reads as foll.\zl ows\zd\
\begin{tabbing}
\hspace{1.7em}$(\sim p\supset p)\supset p$
\end{tabbing}

Other ex\zl amples of log\zl ically\zd\ true form\zl ulas\zd\ are the commut\zl ative\zd\ and associative law for disj\zl unction\zd\ and conj\zl unction\zd\
\begin{itemize}
\item[1.] $p\vee q\equiv q\vee p$
\vspace{-1ex}
\item[2.] $(p\vee q)\vee r\equiv p\vee(q\vee r)$\quad \ul \zl \sout{If either the disj\zl unction\zd\ of $p$ and $q$ is true or $r$ is true then}\zd\ \ud
\vspace{-1ex}
\item[3.] similar formulas hold for add.\zl ition\zd \\[.5ex]
 $p\: .\: q\equiv q\, .\: p$\zl ,\zd \quad $(p\: .\: q)\, .\: r\equiv p\: .\:(q\, .\: r)$
\end{itemize}

$\mathbf{\llbracket 37. \rrbracket}$ Next we have some form\zl ulas\zd\ connecting $\vee$ and $.$ namely at first the famous so called De Morg.\zl an\zd\ formulas:
\begin{tabbing}
\hspace{1.7em}\=$\sim(p\: .\: q)\equiv \;\sim p\;\vee\sim q$\\[.5ex]
\>$\sim(p\vee q)\equiv \;\sim p\: .\sim q$
\end{tabbing}
The left\zl -\zd h.\zl and\zd\ side of the first means not both $p,q$ are true\zl ,\zd\ the right\zl -\zd h\zl and\zd\ side at least one is false \zl \sout{which is..}\zd . $''$ $''$ $''$ $''$ $''$ \zl The left\zl -\zd h.\zl and\zd\ side of the\zd\ sec\zl ond\zd\ $''$ \zl means\zd\ not at least one \zl is\zd\ true\zl ,\zd\ $''$ $''$ $''$ $''$ \zl the right\zl -\zd h\zl and\zd\ side\zd\ both are false\zl .\zd\

These formulas give a means to distribute \ul so to speak \ud the neg\zl ation\zd\ of a product on the two fact\zl ors\zd\ and also the neg\zl ation\zd\ of a sum on the two terms\zl ,\zd\ where however sum has to be changed into prod\zl uct\zd\ and prod\zl uct\zd\ into sum in this distrib.\zl ution\zd\ process\zl .\zd\ Another tautologie\zl y\zd\ conn\zl ecting\zd\ sum and prod\zl uct\zd\ is $\mathbf{\llbracket 38. \rrbracket}$ the distr\zl ibutive\zd\ law which reads exactly analogously as in arith.\zl metic\zd\
\begin{tabbing}
\hspace{1.9em}\=1. \= $p\: .\: (q\vee r)\equiv p\: .\: q\vee p\: .\: r$
\end{tabbing}
\ul bec.\zl ause\zd\ let us ass\zl ume\zd\ left is true then we have \zl \sout{then}\zd\ $p$ \zl full stop deleted\zd\ and two cases $q$\zl ,\zd\ $r$\zl ;\zd\ in the first case $p\: .\:q$\zl ,\zd\ in the sec\zl ond\zd\ $p\: .\:r$ is true\zl ,\zd\ hence in any case \zl right is true\zd\ \ud
\begin{tabbing}
and \=2.\hspace{.3em}\= $p\vee q\, .\: r\equiv (p\vee q)\, .\: (p\vee r)$\\[.5ex]
\>3.\> $(p\supset q)\, .\:(q\supset r)\supset(p\supset r)$\quad Syllog\zl ism,\zd\ \ul Transitivity of $\supset$ \ud \\[.5ex]
\>4.\> $(p\supset q)\supset[(q\supset r)\supset(p\supset r)]$\\[.5ex]
\>\>\sout{$(p\: .\: q\supset r)\supset[p\supset(q\supset r)]$\quad Export}\\*[.5ex]
\>\>\hspace{3.4em}inverse \hspace{5.2em} \sout{Import}\\[.5ex]
\>5.\> $(p\supset q)\, .\: (r\supset s)\supset(p\: .\: r\supset q\, .\: s)$ \sout{\underline{factor}} Leibnitz theorema praeclarum\\[.5ex]
\> \ul \> $(p\supset q)\supset(p\: .\: r\supset q\, .\: r)$ \quad factor \ud \\[.5ex]
\>6.\> $(p\supset q)\: .\: (r\supset s)\supset(p\vee r\supset q\vee s)$\\[.5ex]
\> \ul \> $(p\supset q)\supset(p\vee r\supset q\vee r)$ \quad Sum \ud \\[.5ex]
\>7.\> $p\supset p\vee q$ \zl unreadable word\zd\ \=7$'$. $p\: .\: q\supset p$\\[.5ex]
\>8.\> $p\vee p\supset p$ \zl \sout{taut}\zd\ \>8$'$. $p\;$\= $\supset p \: .\: p$\\
\>\>\hspace{2.4em}$\equiv$\>\>$\equiv$\\[.5ex]
\>9.\> $p\supset(q\supset p\: .\:q)$
\end{tabbing}

\zl On a page after p.\ {\bf 38}., which is not numbered, one finds the following short text containing perhaps exercises or examination questions, which does not seem directly related to the preceding and succeeding pages of the course:
\begin{tabbing}
Log\zl ic\zd\ Notre Dame\\[.3ex]
1. \= \zl a text in shorthand\zd \\[.3ex]
2. \> Trans\zl itivity and\zd\ irrefl\zl exivity\zd\ $\supset$ As\zl \sout{s}\zd ym.\zl metry\zd
\end{tabbing}
On the last page of the notebook, which is also not numbered, there are just two letters ``aq'' or ``ag''.\zd\

\section{Notebook I}\label{0I}
\pagestyle{myheadings}\markboth{SOURCE TEXT}{NOTEBOOK I}
\zl Folder 59, on the front cover of the notebook ``Log.\zl ik\zd\ Vorl.\zl esungen\zd\ \zl German: Logic Lectures\zd\ Notre Dame I''\zd\

\vspace{1ex}

$\mathbf{\llbracket 1. \rrbracket}$ Log\zl ic\zd\ is usually def\zl ined\zd\ as the science whose object are the laws of \zl presumably ``corr.'', which abbreviates ``correct''; if ``corr.'' is read instead as ``con.'', then this would abbreviate ``consistent''\zd\ thinking. According to this def\zl inition\zd\ the cent.\zl ral\zd\ part of log.\zl ic\zd\ must be the theory of inference and the theory of logically true prop.\zl ositions\zd\ [as e.g.\ the law of excl.\zl uded\zd\ middle \zl right square bracket put before the inserted text which follows\zd\ \ul and in order to get acqu.\zl ainted\zd\ with math.\zl ematical\zd\ log\zl ic\zd\ it is perhaps best to go in medias res\zl \emph{in medias res}\zd\ and begin with this centr.\zl al\zd\ part. \ud \zl full stop and right square bracket deleted\zd\

\zl new paragraph\zd\ Prof\zl essor\zd\ Men.\zl ger\zd\ has pointed out in his introduct\zl o\-ry\zd\ lecture that the treatment of these things in trad.\zl itional\zd\ logic and in the current textbooks is \sout{very} unsatisfactory\zl .\zd\ Unsatisfactory \ul from several standp\zl oints\zd . \ud 1.\zl First\zd\ from the standpoint of completeness\zl .\zd\ What the textbooks \sout{give} and also what Arist.\zl otle\zd\ gives is a more or less arbitrary selection of the \ul infinity of \ud \zl the\zd\ laws of logic\zl ,\zd\ whereas in \ul a \ud systematic treatment as is given in math.\zl ematical\zd\ log.\zl ic\zd\ we shall have to develop methods which allow $\mathbf{\llbracket 2. \rrbracket}$ us to obtain all possible logically true prop.\zl ositions\zd\ and to decide of any given prop.\zl osition\zd\ whether or not they are\zl it is\zd\ logically \zl true\zd\ or of an inf.\zl erence\zd\ whether it is correct or not. But 2.\zl secondly\zd\ the class.\zl ical\zd\ treatment is also unsatisf.\zl actory\zd\ as to the question of reducing the \sout{inf.} \ul laws \ud of logic \sout{true prop.} to a cert.\zl ain\zd\ number of primitive laws \sout{to} \ul from \ud which they can be deduced. Although it is sometimes claimed that everything can be deduced from the three fund\zl amental\zd\ laws of contr.\zl adiction,\zd\ excl.\zl uded\zd\ middle and identity or \ul from \ud the modus \zl B\zd arbara this claim has never been \sout{\zl unreadable symbol\zd }\, proved \sout{in trad\zl itional\zd }\, or even clearly formul.\zl ated\zd\ in trad.\zl itional\zd\ logic.

\zl new paragraph\zd\ The chief aim in the first part of these lectures will be to \ul fill those two gaps \zl unreadable word\zd\ [solve those two probl.\zl ems\zd\ in a satisf.\zl actory\zd\ way]\zl ,\zd\ i\zl .\zd e.\ to give \ul as far as possible \ud a complete theory of log\zl ical\zd\ \ul inf\zl erence\zd\ and log\zl ically\zd\ \ud true prop.\zl ositions,\zd\ $\mathbf{\llbracket 3. \rrbracket}$ \ul complete at least for a cert.\zl ain\zd\ very wide domain of prop.\zl ositions,\zd\ \ud and \ul \zl 2 followed by unreadable symbols, perhaps ``.1''\zd\ \ud to show how they can be reduced to a cert\zl ain\zd\ number of primitive laws.

\zl dash from the manuscript deleted, and new paragraph introduced\zd\ The theory of syl.\zl ab\-bre\-vi\-ation for ``syllogisms'' or ``syllogistic''\zd\ as presented in the current textbook\zl s\zd\ is \ul usually \ud divided into two parts\zl :\zd\

\zl display\zd\ 1. The Arist.\zl otelian\zd\ figures and moods of inf.\zl erence\zd\ incl.\zl ud\-ing\zd\ the inf.\zl erences\zd\ with one premise (e.g.\ contrad.\zl iction\zd)\zl ,\zd\

\zl display\zd\ 2. inf.\zl erences\zd\ of \zl unreadable word, should be ``an''\zd\ \ul entirely \ud diff.\zl erent\zd\ kind which are treated under the heading of hypoth.\zl etical\zd\ disj.\zl unc-tive\zd\ conj.\zl unctive\zd\ inferences \ul \zl unreadable text\zd\ \zl \sout{they}\zd\ \ud and which seem to be a Stoic add.\zl ition\zd\ to the Arist.\zl otelian\zd\ figures.

\vspace{1ex}

\noindent Let us begin with the syl.\zl logisms\zd\ of the sec\zl ond\zd\ kind which turn out to be much more fundamental. We have for inst.\zl ance\zd\ the modus ponendo ponens\zl .\zd\

\vspace{1ex}

$\mathbf{\llbracket 4. \rrbracket}$ From the two premises
\begin{tabbing}
\hspace{1.7em}\= 1. \= If Leibn\zl itz\zd\ has inv\zl ented\zd\ the inf.\zl initesimal\zd\ calc\zl ulus\zd\ he was a\\
\` great math.\zl ematician,\zd\ \\
\> 2. \> Leibn\zl itz\zd\ has \zl invented the infinitesimal calculus,\zd\ \\[.5ex]
we conclude\\[.5ex]
\>\>Leibn.\zl itz\zd\ was a great math.\zl ematician.\zd\
\end{tabbing}

\zl From the next paragraph until the end of p.\ \textbf{21}.\ the lower-case propositional letters $p$, $q$ and $r$ are written first as capital $P$, $Q$ and $R$, which are later on alternated with the lower-case letters. In the edited text they are all uniformly lower-case, while in the present source text they are as in the manuscript.\zd\

Generally\zl ,\zd\ if $p$ \zl and\zd\ $q$ are arbitr.\zl ary\zd\ prop.\zl ositions\zd\ and if we have the two premises
\begin{tabbing}
\hspace{1.7em}\= 1. \= If $P$ so $Q$\zl ,\zd \\*[.3ex]
\>2. \> $P$\zl ,\zd \\[.5ex]
we \ul can \ud conclude\\[.5ex]
\>\> $Q$\zl .\zd\
\end{tabbing}
\zl O\zd r \ul take \ud a disjunctive inf.\zl erence\zd\ \ul tollendo ponens\zl .\zd\ \ud If we have the two premises
\begin{tabbing}
\hspace{1.7em}\= 1. \= Either $P$ or $Q$\zl ,\zd \\*[.3ex]
\>2. \>Not $P$\zl ,\zd \\[.5ex]
we \ul can \ud conclude\\[.5ex]
\>\> $Q$\zl .\zd\
\end{tabbing}

It is possible to \sout{write} \ul express \ud those \ul \zl th\zd is \ud syllogism\zl \sout{s}\zd\ \zl ``as'' or ``is'' and a superscripted minus from the manuscript deleted\zd\ by one logically true prop.\zl osition\zd\ as follows:
\begin{tabbing}
\hspace{1.7em}If either $P$ or $Q$ and if not-$P$ then $Q$.
\end{tabbing}
\ul \sout{Other examples of log.\zl ically\zd\ true prop\zl ositions\zd\ of this kind would be} This whole statement will be true whatever $P,Q$ may be\zl .\zd\ \ud

\zl new paragraph\zd\ Now what is the most striking caract\zl er\zd\ of these inf.\zl er\-ences\zd\ which distinguishes them from the Arist.\zl otelian\zd\ syll.\zl ogistic\zd\ \ul fig\-ures \ud ? It is this\zl :\zd\ $\mathbf{\llbracket 5. \rrbracket}$ that in order to make those inf.\zl erences\zd\ it is not nec.\zl essary\zd\ to know anything about the structure of $P$ and $Q$. $P$ or $Q$ (may themselves be disju.\zl nctive\zd\ or hyp.\zl othetical\zd\ prop.\zl ositions\zd )\zl ,\zd\ they may be aff\zl irmative\zd\ or neg.\zl ative\zd\ prop.\zl ositions,\zd\ or they may be s\zl i\zd mple or as compl.\zl icated\zd\ as you want\zl ;\zd\ \zl \ul (\ldots) \ud from the manuscript deleted\zd\ all this is indiff.\zl erent\zd\ for this syl.\zl logism,\zd\ i\zl .\zd e\zl .\zd\ only prop\zl ositions\zd\ as a whole occur in it and it is this fact that makes this kind of syl\zl logism\zd\ simpler and more fundamental than the Arist\zl otelian\zd .
\ul \zl T\zd he law of contrad.\zl iction\zd\ and excl.\zl uded\zd\ middle would be \ul \sout{an} \ud other ex.\zl amples\zd\ of log.\zl ical\zd\ \sout{true prop.\zl ositions\zd} \ul laws \ud of this kind. Bec.\zl ause\zd\ ause e.g.\ the l.\zl aw\zd\ of e.\zl xcluded\zd\ m.\zl iddle\zd\ say\zl s\zd\ for any prop\zl o\-si\-tion\zd\ $P$ either $P$ or $\sim P$ is true and this quite indep.\zl endently\zd\ of the struct.\zl ure\zd\ of $P$. \ud \ul With \ud t\zl written over T\zd he\sout{se} Arist.\zl otelian\zd\ \ul log.\zl ical\zd\ syl.\zl logisms\zd\ it is of course quite diff.\zl erent;\zd\ \sout{moods of course} they \ud depend on the struct.\zl ure\zd\ of the prop.\zl ositions\zd\ \ul involved\zl ,\zd\ \ud  e.g.\ in order to apply the mood Barbara you must know \ul e.g.\ \ud that the two premises are gen.\zl eral\zd\ affirmat.\zl ive\zd\ prop\zl ositions.\zd  \zl insertion sign crossed out in the manuscript\zd\

\zl new paragraph\zd\ Now the theory $\mathbf{\llbracket 6. \rrbracket}$ of log.\zl ically\zd\ true prop.\zl ositions\zd\ and log \zl ical\zd\ inferences in which only prop.\zl ositions\zd\ as a whole occur is called calcul.\zl us\zd\ of prop\zl ositions\zd . In order to \zl unreadable word, ``subject'' or perhaps ``bring''\zd\ it to a syst.\zl ematic\zd\ treatment we have first to examine more in detail the \zl unreadable word, presumably ``connection''\zd\ between prop.\zl ositions\zd\ which \ul can \ud occur in there inf.\zl erences,\zd\ i.e.\ the or, and, if\ldots\ so, and the not. One has introduced special symbols to denote them\zl ,\zd\ in fact there are two diff.\zl erent\zd\ symbol.\zl isms\zd\ for them\zl ,\zd\ the Russell and the Hilb.\zl ert\zd\ symb\zl olism\zd . I shall use in these lect.\zl ures\zd\ Russell's symb\zl olism\zd . In this not is den.\zl oted\zd\ by $\sim$\zl ,\zd\ and by a \sout{point} \ul dot \ud $.$\zl ,\zd\ or by $\vee$ and \ul the \ud if\ldots\ so \sout{i.e.\ the} \zl crossed out unreadable word, presumably ``connection''\zd\ \sout{of impl.\zl ication\zd }\, by $\supset$\zl ,\zd\ $\mathbf{\llbracket 7. \rrbracket}$ i.e.\ if $P,Q$ are arbitrary prop\zl ositions\zd\ then $\sim P$ means $P$ is \sout{wrong} \ul false\zl ,\zd\ \ud $P\: .\: Q$ means \ul both \ud $P$ \ul and \ud $Q$ \sout{are both} \ul are \ud true\zl ,\zd\ $P \vee Q$ means at least one of the prop\zl ositions\zd\ $P,Q$ is true\zl ,\zd\ \ul either both are true or one is true and the other one \sout{wrong} \ul false \ud \ud . This is \zl \sout{a}\zd\ diff\zl erent\zd\ from the mean.\zl ing\zd\ that is given to \ul the \ud or in trad.\zl itional\zd\ logic. There we have to do with the exclusive or\zl ,\zd\ \ul in \zl L\zd at\zl in\zd\ aut\zl \ldots\zd\ aut\zl ,\zd\ \ud which means \ul that exactly \ud one of the two prop.\zl ositions\zd\ $P,Q$ is true and the other one is  \sout{wrong} \ul false\zl ,\zd\ \ud whereas this log.\zl ical\zd\ symb\zl ol\zd\ for or has the meaning of the \zl L\zd at\zl in\zd\ sive\ldots\ sive\zl ,\zd\ \zl a right parenthesis in the manuscript over the second \emph{sive} deleted\zd\ i.e\zl .\zd\ one of the two prop\zl ositions\zd\ is true where it is not excl.\zl uded\zd\ that both are true. \sout{Of course} \zl T\zd he excl.\zl usive\zd\ or \ul as we shall see later \ud can be expressed by a comb.\zl ination\zd\ of the other logistic symb.\zl ols,\zd\ but one has not introduced a proper symb.\zl ol\zd\ for it because it turns out not to be a\zl s\zd\ fund.\zl amental\zd\ as the \sout{sive sive} or in the sense of sive-\zl \ldots\zd\ sive\zl ;\zd\ $\mathbf{\llbracket 8. \rrbracket}$ it \zl is\zd\ not very often used. The \zl n\zd ext symb\zl ol\zd\ is the $\supset$\zl .\zd\ If $P,Q$ are two prop\zl ositions\zd\ \zl dash from the manuscript deleted\zd\ $P \supset Q$ \ul \sout{read as $P$ implies $Q$} \ud means I\zl i\zd f $P$ so $Q$\zl ,\zd\ i\zl .\zd e\zl .\zd\ $P$ implies $Q$\zl .\zd\ \ul \sout{So this $\supset$ is the symb\zl ol\zd\ of implication} Finally we introduce a fifth \zl unreadable word, presumably ``connection''\zd\ $p \equiv q$ ($p$ equiv.\zl alent to\zd\ $q$\zl )\zd\ which means both $p \supset q$ and $q \supset p$\zl .\zd\

\zl new paragraph\zd\ The 5\zl five\zd \zl written over presumably 4\zd\ \zl unreadable word, presumably ``connections''\zd\ introd\zl uced\zd\ so far are called resp.\zl ective\-ly\zd\ negation, conj\zl unction\zd , disj\zl unction\zd , implic\zl ation,\zd\ \ul equivalence\zl ,\zd\ \ud and all of them are called \zl unreadable word, presumably ``connections''\zd\ or operations of the calc.\zl ulus\zd\ of prop\zl ositions\zd . \sout{inst.\zl ead\zd\ of} \zl C\zd onj\zl unction\zd\ and disj\zl unction\zd\ \sout{one s} are also called logical prod\zl uct\zd\ and log\zl ical\zd\ sum respectively\zl .\zd\ \zl A\zd ll of the ment.\zl ioned\zd\ log\zl ical\zd\ op.\zl erations\zd\ \ul exc.\zl lud\-ing\zd\ neg.\zl ation\zd\ \ud are op\zl erations\zd\ with two arg.\zl uments,\zd\ i\zl .\zd e.\ they form a new prop\zl osition\zd\ out of two given ones\zl ,\zd\ exp\zl should be ``for example,''\zd\ $P\vee Q$\zl .\zd\ Only the neg\zl ation\zd\ is an op.\zl era-tion\zd\ with one arg\zl ument\zd\ forming a new prop\zl osition\zd\ $\sim P$ out of any\zl written over something else\zd\ single given prop\zl osition.\zd\ \sout{$P,Q$}

$\Big\lceil$Not only the \sout{symb.} op\zl erations\zd\ $\supset$\zl ,\zd\ $\vee$ \zl and\zd\ $.$ are called impl.\zl ication,\zd\ disj\zl unction and\zd\ conj\zl unction,\zd\ but also an expr.\zl ession\zd\ of the form ${p\supset q}$ \zl ,\zd\ $p\vee q$ is\zl written over ``are''\zd\ called \sout{\zl unreadable word\zd }\, an impl.\zl ication\zd\ etc\zl .,\zd\ where $p,q$ may again be expressions inv.\zl olving\zd\ again $\supset$\zl ,\zd\ $\vee$ etc\zl .\zd\ and $p$\zl ,\zd\ $q$ are called resp.\zl ectively\zd\ first and sec.\zl ond\zd\ member \zl \sout{of}\zd .

Of course if $P$ \zl and\zd\ $Q$ are prop\zl ositions\zd\ then $\sim P$\zl ,\zd\ $\sim Q$\zl ,\zd\ $P\vee Q$\zl ,\zd\ ${P\: .\:Q}$ \zl and\zd\ $P\supset Q$ are \sout{again} \ul also \ud prop.\zl ositions\zd\ and hence to them the op.\zl erations\zd\ of the calc.\zl ulus\zd\ of prop\zl ositions\zd\ can again be applied\zl ,\zd\ so as to get more compl.\zl ex\zd\ expr.\zl essions,\zd\ e.g\zl .\zd\ ${P \vee (Q\: .\:R)}$ \zl in this formula a left parenthesis before $P$ and a right parenthesis after $Q$ have been crossed out in the manuscript\zd \zl ,\zd\ eithe\zl r\zd\ $P$ is true or $Q$ and $R$ are both true\zl .\zd\ \ud

The disj.\zl unctive\zd\ inf.\zl erence\zd\ I mentioned before would read in this symbolism as follows: $[(P\vee Q)\: .\:\sim P]\supset Q$\zl .\zd\ \zl The text to be inserted to which the sign $\otimes$ in the manuscript at this place should refer to is missing.\zd\ \sout{Of course} \ul You see \ud in more comp\zl icated\zd\ expressions as this one brack.\zl ets\zd\ have to be used exactly as in algebra in order to indicate the order in which the operations have to be applied. E.g.\ if \zl unreadable symbol, perhaps ``I''\zd\ put the \ul round \ud brackets in this expr.\zl ession\zd\ like this ${P\vee (Q\: .\:\sim P)}$\zl ,\zd\ it would have a diff\zl erent\zd\ mean\zl ing,\zd\ namely either $P$ is true or $Q$ and $\sim P$ are both true.

\zl new paragraph\zd\ There is an interest.\zl ing\zd\ $\mathbf{\llbracket 9. \rrbracket}$ remark \zl a small superscribed $\vee$ from the manuscript deleted\zd\ due to L\zl {\L}\zd uk.\zl asiewicz\zd\ that one can dispense with the brackets if one writes the operational symb\zl ols\zd\ $\vee$\zl ,\zd\ $\supset$ \zl a dot in the manuscript under $\supset$ deleted\zd\ etc\zl .\zd\ always in front of the prop\zl ositions\zd\ to which they are applied\zl ,\zd\ e.g.\ ${\supset\! p\,q}$ inst\zl ead of\zd\ ${p\supset q}$. Then \ul e.g.\ \ud the two diff.\zl erent\zd\ possibilities for the expr.\zl ession\zd\ in squ.\zl are\zd\ brackets would be \sout{aut.} distinguished \ul aut\zl omatically\zd\ \ud bec.\zl ause\zd\ the first would be written as foll\zl ows\zd\ $.\vee PQ\sim P$\zl ;\zd the sec.\zl ond\zd\ would read $\vee P\: .\: Q\sim P$\zl ,\zd\ so \sout{its} \ul that \ud \zl they\zd\ diff\zl er\zd\ul from each other without the use of brack.\zl ets\zd\ \ud as you see and it can be proved that it is quite generally so.
\ul But since the formulas in the bracket notation are more easily readable I shall stick to this \ul not.\zl ation\zd\ \ud and put the op.\zl erational\zd\ symb.\zl ols\zd\ in betw.\zl een\zd\ the prop.\zl ositions\zd . \ud

\zl new paragraph\zd\ You know in algebra one can spare \sout{the} many brackets by the convention that the $\mathbf{\llbracket 10. \rrbracket}$ mult.\zl iplication\zd\ connects stronger than add\zl ition\zd ; \ul e.g\zl .\zd\ $a\cdot b+c$ means $(a\cdot b)+c$ and not $a\cdot (b+c)$\zl .\zd\ \zl The $\cdot$ in the manuscript, here and in later notebooks, in particular Notebook VI, where $\cdot$ is meant to stand for set intersection, is often indistinguishable from $.$, but the meaning of the text makes it possible to make the distinction, and this will be done without notice.\zd\ \ud We can do something similar here by stipulating an order of the strength in which the log\zl ical\zd\ symb.\zl ols\zd \zl bind\zd\ \zl , so\zd\ that\zl :\zd

\begin{itemize}
\item[1.] the $\sim$ (and similarly any op.\zl eration\zd\ with just one prop.\zl osition\zd\ as ar\-g\zl u\-ment\zd ) connects stronger than any op.\zl eration\zd\ with two ar\-g\zl u\-ments,\zd\ as $\vee$\zl ,\zd\ $\supset$ \zl and\zd\ $.$\zl ,\zd\ so that $\sim p\vee q$ means $(\sim p) \vee q$ and not $\sim (p\vee q)$;
\vspace{-1ex}
\item[2.] the disj.\zl unction\zd\ and conj.\zl unction\zd\ bind stronge.\zl r\zd\ than implic.\zl a\-tion\zd\ and equiv\zl alence,\zd\ so that e.g.\ $p\vee q \supset r\: .\:s$ means $(p\vee q) \supset (r\: .\:s)$ and not \ul \zl unreadable word, perhaps: perh\zl aps\zd \zd\ \ud $p\vee [(q \supset r)\: .\:s]$\zl .\zd
\end{itemize}
\ul \zl A\zd\ third conv.\zl ention\zd\ consists in leaving out brack.\zl ets\zd\ in such \zl ex\-pres\-sions\zd\ as $(p\vee q)\vee r$ exactly as in $\lfloor ( \rfloor a+b)+c$ \zl here, after ``whereas'', the sentence in the manuscript is broken\zd \zl .\zd\ A similar convention is made for $.\;$. \ud

After those merely symb.\zl olic\zd\ conventions the next thing we have to do is to examine in more detail the meaning of the op.\zl erations\zd\ of the calc\zl ulus\zd\ of prop\zl ositions\zd . $\mathbf{\llbracket 11. \rrbracket}$ Take e.g.\ disj.\zl unction\zd\ $\vee$. If any two prop\zl ositions\zd\ $P,Q$ are given $P\vee Q$ will again be a prop\zl osition\zd . Hence the disj.\zl unction\zd\ is an operation which applied to any two prop.\zl ositions\zd\ gives again a prop\zl osition.\zd\ But now (and this is the dec.\zl isive\zd\ point) \ul this op.\zl eration\zd\ is such that \ud the truth or falsehood of the composite prop\zl osition\zd\ $P\vee Q$ depends in a def.\zl inite\zd\ way on the truth or falsehood of the const.\zl ituents\zd\ $P,Q$\zl .\zd\ \ul \sout{and depends only on the truth or falsehood of the const.} \ud This dependence can be expressed most clearly in \zl the\zd\ form of a table as follows\zl : l\zd et us form three col\zl umns,\zd\ one headed by \ul $p$\zl ,\zd\ \ud one by $q$\zl ,\zd\ one by $p\vee q$\zl ,\zd\ and let us write $+$ for true and $-$ for \sout{wrong} \ul false\zl .\zd\ \ud Then for the prop\zl osition\zd\ $p\vee q$ we have the foll.\zl owing\zd\ $4$\zl four\zd\ possibilities\zl :\zd
\begin{center}
\begin{tabular}{ c|c|c|c }
 $p$ & $q$ & $p\vee q$ & $p\: o\: q$\\[.5ex]
 + & + & + & $-$ \\
 + & $-$ & + & + \\
 $-$ & + & + & + \\
 $-$ & $-$ & $-$ & $-$
 \end{tabular}
\end{center}
Now $\mathbf{\llbracket 12. \rrbracket}$ for each of these fo\zl u\zd r cases we can det.\zl ermine\zd\ wheth.\zl er\zd\ ${p\vee q}$ will be true or false\zl ,\zd\ namely since ${p\vee q}$ means that one or both of the prop\zl ositions\zd\ $p,q$ are true it will be true in the first\zl ,\zd\ sec.\zl ond and\zd\ third case\zl ,\zd\ and \sout{wrong} \ul false \ud in the last case. And we can consider this table as the most precise def\zl inition\zd\ of what $\vee$ means.

\zl new paragraph\zd\ \ul \sout{One also} \ul It is usual to \ud call truth and falsehood the truth values\zl ,\zd\ \ul so there are exactly two truth values\zl ,\zd\ \ud and say\zl \sout{s}\zd\ that a true prop\zl osition\zd\ has the truth value ,,\zl ``\zd truth'' \zl (\zd den\zl oted\zd\ by $+$) and a false prop.\zl osition\zd\ has the truth value \zl ``\zd false\zl ''\zd\ (den\zl oted\zd\ by $-$)\zl ,\zd\ \ul so that any prop.\zl osition\zd\ has a un.\zl iquely\zd\ det.\zl ermined\zd\ truth value\zl .\zd\ \ud  The truth table then shows how the truth value of the comp\zl osite\zd\ expr\zl essions\zd\ depends on the truth value of the constituents. \ud The excl.\zl usive\zd\ or would have another truth table\zl ;\zd\ namely if we denote it by $o$ for the moment we have \zl that\zd\ $p\: o\: q$ is \sout{wrong} \ul false \ud if both $p$ and $q$ are true\zl ,\zd\ and it is \sout{wrong} \ul false \ud if both are \sout{wrong} \zl false\zd\ but true in the two other cases. \zl It is not clear what the words ``where exactly'' inserted in the manuscript at this place refer to, and they will be deleted.\zd\ The op.\zl eration\zd\ $\sim$ $\mathbf{\llbracket 13. \rrbracket}$ has of course the foll.\zl owing\zd\ truth table\zl :\zd
\begin{center}
\begin{tabular}{ c|c}
 $p$ & $\sim p$ \\[.5ex]
 + & $-$ \\
 $-$ & +
 \end{tabular}
\end{center}
Here we have only two poss.\zl ibilities:\zd\ $p$ true or $p$ wrong\zl ,\zd\ and in the first case we have \zl that\zd\ not-$p$ is wrong \zl while\zd\ in the sec\zl ond\zd\ it \zl is\zd\ true. Also the truth table for $.$ can easily be det.\zl ermined:\zd
\begin{center}
\begin{tabular}{ c|c|c }
 $p$ & $q$ & $p\: .\: q$\\[.5ex]
 + & + & + \\
 + & $-$ & $-$ \\
 $-$ & + & $-$ \\
 $-$ & $-$ & $-$
 \end{tabular}
\end{center}
(I think \zl unreadable text, possibly ``I will'' or ``I can''\zd\ leave that to you).\zl .)\zd\

\zl new paragraph\zd\ A little more diff.\zl icult\zd\ is the question of the truth table for $\supset$. \zl The following text, until the end of p.\ \textbf{13}., is crossed out in the manuscript: In fact $\subset$\zl $\supset$\zd\ can be interpreted in different ways\zl ,\zd\ and for cert.\zl ain\zd\ interpret.\zl ations\zd\ there exist no truth table\zl ,\zd\ e.g\zl .\zd\ if we define $P\supset Q$ to mean ,,\zl ``\zd From $P$ $Q$ follows logically'' then the truth value of $P\supset Q$ is not determined at all by the\zd\ $\mathbf{\llbracket 14. \rrbracket}$ $p\supset q$ was defined to mean \zl ``\zd If $p$ is true $q$ is \ul also \ud true\zl ``\zd . So let us assume that \sout{we know} for two given prop.\zl ositions\zd\ $P,Q$ we know that $P\supset Q$ is true\zl ,\zd\ i\zl .\zd e.\ \ul assume that \ud we know \zl ``\zd If $P$ then $Q$\zl ''\zd\ \ul but nothing else \ud \zl . W\zd hat can we conclude then about the \ul possible \ud truth values of $P$ and $Q$.\zl ?\zd\ \zl As for the analogous table on p.\ \textbf{19}.\ of Notebook 0, it is not indicated in the manuscript where the following table should be inserted. The text that follows is a comment upon it.\zd\
\begin{center}
\text{Ass.\zl umption\zd\ $p \supset Q$}
\end{center}
\begin{center}
\begin{tabular}{@{}l}
 $\left.
 \begin{tabular}{@{}*2{p{0.3cm}}}
 $P$ & $Q$ \\[.5ex]
 $-$ & + \\
 $-$ & $-$ \\
 + & +
 \end{tabular}
 \right\}\text{possible truth val.\zl ues\zd\ for $P,Q$\zl full stop deleted\zd }$\\[5ex]
$\left.
 \begin{tabular}{@{}*2{p{0.3cm}}}
 + & $-$
 \end{tabular}
 \right\}\text{impossible}$
\end{tabular}
\end{center}

\noindent First it may \ul cert.\zl ainly\zd\ \ud happen that $P$ is false because the \ul as\-sump\-tion \ud stat\zl written over another letter\zd ement \zl ``\zd If $P$ then $Q$\zl ''\zd\ says nothing about the truth \ul or falsehood \ud of $P$. \ul And\zl written over another word, perhaps ``Now''\zd\ \ud \zl i\zd n this case where $P$ is false \ul $Q$ may be true as well as false because \ud the assumption \zl ``\zd If $P$ then $Q$\zl ''\zd\ \sout{determines nothing about the truth of $Q$ because it only says \zl ``\zd if $P$ is true $Q$ is true\zl ''\zd\ but it says nothing about the case where $P$ is false} \ul says nothing about what happens to $Q$ if $P$ is false but only if $P$ is true\zl .\zd\ \zl crossed out unreadable word\zd\ \ud So we have both possib\zl ilities:\zd\ $P$ false $Q$ true\zl ,\zd\ $P$ f\zl alse\zd\ $Q$ f\zl alse.\zd\ Next we have the \sout{case in which} \ul possibility that \ud $P$ is true\zl .\zd\ $\mathbf{\llbracket 15. \rrbracket}$ \ul But \ud I\zl i\zd n this case \sout{it follows from the assumption $p\supset q$ that} \ul owing to the ass.\zl umption\zd\ \ud $q$ \sout{is} \ul must \ud also \zl be\zd\ true. So that the poss\zl ibility\zd\ $P$ true $Q$ false is excluded and we have only this third possibility $p$ true $q$ true\zl ,\zd\ \ul and this poss.\zl ibility\zd\ may of course really happen \ud . So from the ass.\zl umption\zd\ $P\supset Q$ it follows that either one of the first three cases happens\zl .\zd\ \ul \sout{i\zl .\zd e\zl .\zd\ if $P\supset Q$ then} \ud \sout{$(\sim P\: .\: Q) \vee (\sim P\: .\:\sim Q)\vee (P\: .\:Q)$} But \sout{als} \ul we have \ud also vice versa\zl :\zd\ If \sout{$(\sim P\: .\: Q) \vee (\sim P\: .\:\sim Q)\vee (P\: .\:Q)$} \ul one of the first three poss.\zl ibilities\zd\ of the truth values is realis\zl z\zd ed \ud then $(p\supset q)$ \ul is true \ud . Because let us assume we know that one \zl of\zd\ the three \ul cases \ud written down \sout{happens} \ul is realis\zl z\zd ed\zl .\zd\ \ud \ul I claim \ud then we know also: ,,\zl ``\zd If $p$ \ul is true \ud then $q$ \ul is true\zl ''\zd\ \ud . \sout{Because} \ul \sout{That's easy bec.\zl ause\zd } \ud If $p$ is true only the third of the three poss.\zl ibilities\zd\ can be realis\zl z\zd ed (in \ul all \ud the other\zl s\zd\ $p$ is false)\zl ,\zd\ but in this third possib.\zl ility\zd\ $Q$ is true\zl .\zd

\zl Here begins a page with its number \textbf{16}.\ crossed out, and the following crossed out text: this third \ul poss.\zl ibil\-ity\zd\ \ud once $Q$ is also true so we have If $P$ then $Q$ \zl new paragraph in the manuscript\zd\ it is the only one of the first three in which $P$ is true of the three possible cases can hold and in this case $Q$ is also true\ldots\zd\

\zl The following text on the remainder of this page with number \textbf{16}. crossed out is not very clearly crossed out, but a text with the same content can be found on the next page numbered \textbf{16}.: So we have proved a complete equiv.\zl alence\zd\ between $p\supset q$ on the one hand and the disj.\zl unction\zd\ $(\sim p\: .\: q) \vee (\sim p\: .\:\sim q)\vee (p\: \lfloor .\rfloor \:q)$ on the other hand\zl ,\zd\ so that we can define impl\zl ication\zd\ by this disj\zl unction\zd . But this disj\zl unction\zd\ can be written in a simpler form as follows $\sim p \vee q$\zl .\zd\ \zl I\zd t is easy to see that this disj.\zl unction\zd\ of three cases is equivalent with the \zl here the text in the manuscript breaks\zd \zd\ $\mathbf{\llbracket 16. \rrbracket}$ So we see that the statement $p \supset q$ is exactly equivalent with the statement that one of the three marked cases \ul for the \ud distr.\zl ibution\zd\ of truth values is realis\zl z\zd ed\zl ,\zd\ \sout{and not the fourth one} i\zl .\zd e.\ $p\supset q$ is true in each of the three marked cases \sout{realis\zl z\zd ed and only then i\zl .\zd e\zl .\zd\ it is} and false in the last case (\sout{where none of these three poss.\zl ibilities\zd\ is realis\zl z\zd ed})\zl .\zd\ So we have obtained \sout{the} a truth table for implication. However there a\zl re\zd\ two imp.\zl ortant\zd\ remarks about it\zl ,\zd\ namely\zl :\zd\

\zl new paragraph\zd\ 1. Exactly the same \ul truth \ud table can also be obtained by a combin.\zl ation\zd\ of op.\zl erations\zd\ introd\zl uced\zd\ previously\zl ,\zd\ namely $\sim p \vee q$ has the same truth table \zl deleted ``for\ldots'' from the manuscript\zd\
\begin{center}
\begin{tabular}{ c|c|c|c }
 $p$ & $q$ & $\sim p$ & $\sim p\vee q$\\[.5ex]
 $-$ & \zl $-$\zd & + & + \\
 \zl $-$\zd & \zl +\zd & + & + \\
 \zl +\zd & \zl $-$\zd & $-$ & $-$ \\
 \zl +\zd & \zl +\zd & $-$ & +
 \end{tabular}
\end{center}
$\mathbf{\llbracket 17. \rrbracket}$ \sout{Hence} Since $p\supset q$ and $\sim p \vee q$ have the same truth table they will be equ.\zl ivalent,\zd\ i\zl .\zd e.\ whenever the one expr.\zl ession\zd\ is true the other one will also be true and vice versa. This makes \zl it\zd\ possible to define $p\supset q$ by $\sim p \vee q$ and this is the standard way of introd.\zl ucing\zd\ impl.\zl ication\zd\ in math\zl ematical\zd\ log\zl ic\zd .

\zl new paragraph, 2.\zd\ The sec.\zl ond\zd\ remark \ul about impl.\zl ication\zd\ \ud is this \sout{more imp\zl ortant\zd }. We must be careful not to forget that \zl $\;\mid$ from the manuscript deleted\zd\ $p \supset q$ was understood to mean simply \zl ``\zd If $p$ then $q$\zl ''\zd\ and only this made the const.\zl ruction\zd\ of the truth table possible\zl .\zd\ \zl $\;\mid$ from the manuscript deleted\zd\ \ul We have deduced the truth table for impl\zl ication\zd\ from the ass\zl umption\zd\ that $p \supset q$ means \zl ``\zd If $p$ then $q$\zl ''\zd\ and nothing else\zl .\zd\ \ud There are other meanings $\mathbf{\llbracket 18. \rrbracket}$ perhaps even more suggested by the term impl.\zl ication\zd\ for which our truth table would be completely inadequate. E.g.\ \sout{if we assume} \ul $P\supset Q$ could be \ud given \zl \sout{$P \supset Q$}\zd\ the meaning \sout{of}: $Q$ is a log\zl ical\zd\ consequence of $P$ \zl written over $Q$\zd\ \zl ,\zd\ i.e.\ $Q$ can be derived from $P$ by means of \ul a chain of \ud syllogisms\zl .\zd\ \zl \sout{then}\zd\ \sout{it is easy to see that there can exist no truth table at all for this impl\zl ication\zd\ thus defined\zl .\zd\ For consider the first line of the supp.\zl osed\zd }

This kind of impl\zl ication\zd\ is usually called strict impl.\zl ication\zd\ and denoted in this way \zl $\prec$\zd\ \ul \sout{as opposed to} \ul and the impl\zl ication\zd\ $p\supset q$ def\zl ined\zd\ before is called \ud material impl.\zl ication\zd\ \ul if it is to be distinguished\zl .\zd\ \ud \sout{$p \supset q$ def.\zl ined\zd\ by $\sim p \vee q$} Now it is easy to see that our truth table is false for strict impl\zl ication.\zd\ \zl \sout{but that $''$} $\;$I\zd n order to prove that \sout{that there exists no truth t.\zl able\zd\ for it} \ud \sout{Now \ul \zl unreadable word\zd\ \ud } consider the first line \ul of \ud a supp.\zl osed\zd\ such table
\begin{center}
\begin{tabular}{ c|c|c }
 $p$ & $q$ & $p \prec q$\\[.5ex]
 + & + & \\
 & & \\
 & & \\
 & &
 \end{tabular}
\end{center}
where $p$ and $\mathbf{\llbracket 19. \rrbracket}$ $q$ are both true and ask what will be the truth value of $p \prec$ strictly $q$\zl . I\zd t is clear that this truth value will not be uniquely det\zl ermined\zd . For take e.g.\ for $p$ the prop.\zl osition ``\zd The earth is a sphere\zl ''\zd\ and for $q$ \zl ``T\zd he earth is not a disk\zl ''.\zd\ \zl T\zd hen $p$ \zl and\zd\ $q$ are both true and $p\prec q$ is also true bec.\zl ause\zd\ \sout{If} \ul from the prop\zl osition\zd\ that \ud the earth is a sp\zl here\zd\ it foll.\zl ows\zd\ by log\zl ical\zd\ inf\zl erence\zd\ that it is not a disk\zl ;\zd\ on the other hand if you take for $p$ again\ldots\zl the same proposition\zd\ and for $q$ \zl ``\zd France is a rep.\zl ublic''\zd\ then again both $p$ and $q$ are true but $p \prec q$ is wrong\zl .\zd\ \zl ``bec from\ldots'' from the manuscript deleted\zd\ $\mathbf{\llbracket 20. \rrbracket}$ So we see the truth value of $p \prec \lfloor q\rfloor$ is not \sout{extens.\zl ionally\zd }\, uniquely det.\zl ermined\zd\ by the truth values of $P$ and $Q$\zl ,\zd\ and therefore no truth table exists. Such \zl unreadable word, presumably ``connections''\zd\ for which no truth t\zl able\zd\ exists are called intensional as opposed to extensional ones for which they \ul do \ud exist. \ul The ext.\zl ensional\zd\ conn.\zl ections\zd\ are \sout{als} called also truth\zl \,\zd functions\zd .

\zl new paragraph\zd\ So we see the impl.\zl ication\zd\ which we introd.\zl uced\zd\ does not mean log.\-\zl ical\zd\ consequence. Its meaning is best given by the simple \zl ``\zd if then\zl ''\zd\ which \sout{is used in many cases where the} \ul has much wider signif\zl icance\zd\ than just \ud logical consequ\-\zl ence\zd . E.g\zl .\zd\ if I say \zl ``\zd If \sout{I cannot} \ul he cannot \ud come \zl \sout{I}\zd\ \sout{shall} \ul he will \ud  telephone to you\zl '',\zd\ that has nothing to do with log\zl ical\zd\ rel\zl ations\zd\ betw\zl een\zd\ $\mathbf{\llbracket 21. \rrbracket}$ \sout{my} \ul his \ud coming and \sout{the} \ul his \ud telephoning\zl ,\zd\ but it simply means he will either come or telephone which is exactly the meaning expressed by the truth table. \ul Now the decisive point is that we don't need any other kind of impl\zl ication\zd\ besides material in order to develop \zl full stop and crossed out unreadable word from the manuscript deleted\zd\ the theory of inf.\zl erence\zd\ \sout{And therefore} because in order to make the concl.\zl usion\zd\ from a prop\zl osition\zd\ $P$ to a prop\zl osition\zd\ $Q$ it is not necessary to know that $Q$ is a log.\zl ical\zd\ cons.\zl equence\zd\ of $P$. It is quite sufficient to know\zl ''\zd If $P$ is true $Q$ is true\zl ''.\zd\ \zl ``e.g'' from the manuscript deleted\zd\ Therefore I shall \sout{not introduce} \sout{material} \sout{strict} \ul use only mat\zl erial\zd\ impl\zl ication,\zd\ \ud at least in the first half of my lectures\zl ,\zd\ and use \sout{als} the terms \zl ``\zd implies\zl ''\zd\ and \zl ``\zd it follows\zl ''\zd\ only in this sense. \ud

\zl The following text within square brackets until the end of p.\ \textbf{21}.\ seems to be crossed out: [Perhaps the term impl.\zl ication\zd\ is not very well chosen from this st.\zl andpoint\zd\ because it \sout{convey} suggests something like log.\zl ical\zd\ consequ\zl ence\zd\ but since it \ul has been \ud \zl in\zd\ comm\zl on\zd\ use for this notion \ul for many years \ud it is not adv.\zl antageous\zd\ to change it and it is not nec\zl essary\zd\ if we keep in \zl mind\zd\ what it means. I shall also use the term \zl ``\zd it follows\zl ''\zd\ to denote $\supset$ sometimes because in more complicated expr.\zl essions\zd\ it is desirable to have sev.\zl eral\zd\ diff.\zl erent\zd\ words for implic\zl ation\zd . So I don't want to use the word \zl ``\zd follow\zl ''\zd\ in the sense of log.\zl ical\zd\ conse\-qu.\zl ence\zd\ but of consequence in a more general sense\zl .\zd ]\zd\

$\mathbf{\llbracket 22. \rrbracket}$ \zl This page in the manuscript begins with the following sentence, which seems to be crossed out: A confusion with strict impl\zl ication\zd\ is not to be feared because I shall confine myself to mat.\zl erial\zd\ impl.\zl ication\zd\ in the whole develop.\-\zl ement\zd\ of \zl the\zd\ calc\zl ulus\zd\ of prop.\zl ositions\zd ]\zd\ \ul This simplifies very much the whole theory of inf.\zl erence\zd\ bec.\zl ause\zd\ mat\zl erial\zd\ impl.\zl ica\-tion\zd\ def\zl ined\zd\ by the truth table is a much simpler notion. I do not want to say by this that \ul a theory \ud \zl of\zd\ strict impl.\zl ication\zd\ may not be interesting and important for certain purp.\zl oses;\zd\ in fact I hope to \sout{develop it} speak about it later on in my lectures. But its theory bel\zl ongs\zd\ to an entirely diff\zl erent\zd\ part of logic than that \zl with\zd\ which we are deal.\zl ing\zd\ at present\zl ,\zd\ namely \ul it bel.\zl ongs\zd\ \ud to the log.\zl ic\zd\ of modalities\zl .\zd

\zl The following words to be deleted from the manuscript are presumably superseded by the inserted words above them immediately after: Our def\zl inition\zd\ of impli\zl cation\zd\ has some\zd\ \ul Now I come to some \ud apparently parad.\zl oxical\zd\ consequ\zl ences\zd\ \ul of our def\zl inition\zd\ of impl\zl ication\zd\ \ud whose paradoxity \zl replaced by ``paradoxicality'' in the edited text; see also the beginning of p.\ 29.\ of Notebook 0\zd\ however disappears if we remember that \sout{it} \ul implic.\zl ation\zd\ \ud does not mean log\zl ical\zd\ consequ\zl ence\zd . \sout{We have at first} \ul \sout{We have} \ud \zl colon from the manuscript deleted\zd\ Namely \sout{since $p\supset q$ is equiv\zl alent\zd\ with $\sim p \vee q$ we have} \zl i\zd f we look at the truth table for $p\supset q$ we see at once \zl that\zd\ $p\supset q$ is always true if $q$ is true what\zl e\zd ver $p$ may be. So that means a true prop\zl osition\zd\ is implied by any prop\zl osition\zd . Sec.\zl ondly\zd\ we see that $p\supset q$ is always true if $p$ is false whatever $q$ $\mathbf{\llbracket 23.\:I \rrbracket}$ may be i\zl .\zd e\zl .\zd\ \zl a\zd\ false prop.\zl osition\zd\ implies any prop.\zl osition\zd\ whatsoever\zl .\zd\ In other words:
\zl ``\zd An impl\zl ication\zd\ with true sec\zl ond\zd\ member is true (whatever\ldots\zl the first member may be\zd ) and an impl.\zl ication\zd\ with a false first member is always true (what\zl ever\zd\ the sec\zl ond\zd \ldots \zl member may be\zd ).\zl ''\zd\ Or written in formulas this means $q \supset (p\supset q)$, $\sim p \supset (p\supset q)$. \sout{Of course this is} \ul Both of these form.\zl ulas\zd\ are \ud also \sout{an} im\zl mediate\zd\ consequ\zl ences\zd\ of the fact that $p\supset q$ is equiv\zl alent\zd\ with $\sim p \vee q$ because \ul $\sim p \vee q$ \sout{that} says \sout{just} \ul exactly \ud either $p$ is false or $q$ is true\zl,\zd\ so it will \ul always \ud be true if $p$ is false and if $q$ is true whatever the other prop\zl osition\zd\ may be. \sout{But} \ul These formula\zl s\zd\ are rather unexp.\zl ected\zd\ and \ud if we apply them \sout{formulas} to spec\zl ial\zd\ cases we get strange consequences. E.g.\ $\mathbf{\llbracket 24. \rrbracket}$ \zl ``\zd The earth is not a sphere\zl ''\zd\ implies that France is a rep\zl ublic,\zd\ but it also impl.\zl ies\zd\ that France is no\zl t a\zd\ rep.\zl ublic\zd\ because a false prop.\zl osition\zd\ implies any prop\zl osition\zd\ whatsoever. Similarly the prop.\zl osition\zd \zl ``\zd France is a rep.\zl ublic''\zd\ is implied by any other prop.\zl osition\zd\ whatsoever\zl ,\zd\ \zl \sout{it may be}\zd\ true or false. But these consequences are only paradoxical if we understand implic.\zl ation\zd\ to mean logical consequence. For the \zl ``\zd if\ldots\ so\zl ''\zd\ meaning they are quite natural\zl ,\zd\ e.g.\ \zl \sout{the}\zd\ \sout{first} $q\supset(p\supset q)$ means: If $q$ is true then $q$ is \ul \sout{also} \ud true \ul also \ud if $p$ is true\zl ,\zd\ and $\sim p\supset(p\supset q)$ If we have a false prop.\zl osition\zd\ $p$ then if $p$ is true anything is true\zl .\zd\ $\mathbf{\llbracket 25.\: I \rrbracket}$ Another of this\zl these\zd\ so called parad\zl oxical\zd\ consequences is this $(p\supset q)\vee (q\supset p)$.\zl ,\zd\ i\zl .\zd e.\ of any two arbitrary prop.\zl ositions\zd\ one must imply the other one. That it must be so is proved as foll.\zl ows:\zd\ \sout{\ul is \ud} \sout{clear because} \sout{\ul of the foll\zl owing\zd\ reason \ud} $q$ must be either true \sout{\zl unreadable word\zd }\, or false\zl ;\zd\ if $q$ is true the first member of the disj\zl unction\zd\ is true \zl bec.\zl ause\zd\ inserted in the manuscript deleted\zd\  \zl and\zd\  if $q$ is false the sec.\zl ond\zd\ member is true because a false prop\zl osition\zd\ implies any other. So (one of the two members of the impl\zl ication\zd\ is true) \ul either $p\supset q$ or $q\supset p$ \ud in any case. \ul

\zl new paragraph\zd\ We have here three examples of logically true formulas\zl ,\zd\ i\zl .\zd e.\ formulas which are true whatever the prop.\zl ositions\zd\ $p,q$ may be. Such formulas are called \zl unreadable text, probably in shorthand\zd\ tautological \sout{and \zl unreadable word\zd\ It is in their} and it is exactly the chief aim of the calc.\zl ulus\zd\ of prop\zl ositions\zd\ to investigate those tautolog\zl ical\zd\ form\zl ulas\zd . \sout{In order to get more acquainted with the our symbolism} \sout{and} \sout{\ul also with the our notation \ud the fund.\zl amental op.\zl erations\zd\ of the calc.\zl ulus\zd\ of prop.\zl osi\-tions\zd}

I shall \ul begin with \ud discussing \sout{now} a few \ul more \ud examples of \ul such \ud logically true prop.\zl ositions\zd\ before going \ul over \ud to general con\ul sider \ud \zl ations.\zd\
$\mathbf{\llbracket 26. \: I \rrbracket}$ We have at first the trad.\zl itional\zd\ hyp.\zl othetical\zd\ \sout{inf.\zl er\-ences\zd }\, and disj\zl unctive\zd\ inferences which in our notation read as follows\zl :\zd

\begin{itemize}
\item[1.] $p\, \lfloor .\rfloor \:(p\supset q)\supset q$ \quad pon.\zl endo\zd\ pon\zl ens\zd\ \quad
\vspace{-1ex}
\item[[2.] $\sim q\, .\:(p\supset q)\supset \sim p$ \quad toll.\zl endo\zd\ toll.\zl ens\zd ]
\vspace{-1ex}
\item[3.] $(p\vee q)\, .\: \sim q\supset p$ \quad toll.\zl endo\zd\ pon\zl ens\zd \\
disj.\zl unctive\zd\ \sout{toll} pon.\zl endo\zd\ tollens does not hold for the not ex\-cl.\zl u\-sive\zd\ $\vee$ which we have
\vspace{-1ex}
\item[4\zl .\zd ] The inf\zl erence\zd\ which is called dilemma \sout{by}\\
 $(p\supset q)\, .\:(r\supset q)\supset(p\vee r\supset q)$
\end{itemize}

$\mathbf{\llbracket 38.1 \: II \rrbracket}$ \zl Notebook 0 ends with p.\ \textbf{38}., and hence judging by the number given to the present page, it might be a continuation of Notebook 0. The content of this page does not make obvious this supposition, but does not exclude it.\zd\ Last time and also t\zl written over another symbol\zd oday in the \zl unreadable symbol\zd\ classes we set up the truth tables for some of the funct\zl ions\zd\ which occur in the calc.\zl ulus\zd\ of prop\zl ositi\-ons\zd . \zl For an insertion sign in the manuscript at this place no text to be inserted is given on this page, but this might point to the crossed out text ``because all f\zl u\zd nct.\zl ions\zd\ involved are truth functions but'' on p.\ ${\textbf{40}.\: \textbf{II}}$ preceded by an insertion sign, for which text it is not indicated where on p.\ ${\textbf{40}.\: \textbf{II}}$ it is to be inserted. The insertion sign at this place is deleted, because anyway the text from p.\ ${\textbf{40}.\: \textbf{II}}$ is crossed out.\zd\ Their purp\zl ose\zd\ is to give a\zl \sout{n}\zd\ \sout{absolutely} precise def.\zl inition\zd\ of the funct\zl ions\zd\ concerned because they state exactly the cond\zl itions\zd\ under which \ul the prop.\zl osit\-ion\zd\ to be def.\zl ined,\zd\ \ud e.g.\ $p\vee q$\zl ,\zd\ is true and under which cond\zl itions\zd\ it is not true. In ordinary language we have also \zl the\zd\ notions and, or, I\zl i\zd f etc\zl .\zd\ which have \ul very\zl underlined or crossed out\zd\ \ud approximately the same meaning\zl ,\zd\ but for setting up a math.\zl ematical\zd\ theory it is nec.\zl essary\zd\ that the not.\zl ions\zd\ involved have a higher degree of preciseness than the notions of ordinary language \sout{and and}\zl .\zd\ It is exactly this what is accompli\zl shed\zd\ by the truth tables\zl . \sout{and in}\zd\

$\mathbf{\llbracket 40. \: II \rrbracket}$ \zl There is no page numbered \textbf{39}.\zd\ \sout{the truth tables give us almost immediately such a method} Take e.g.\ the formula $p\, .\:\lfloor (\rfloor p\supset q)\supset q$\zl ,\zd\ the mod\zl us\zd\ ponendo ponens\zl .\zd\ \sout{\Big\lceil In order to ascertain that it is logically true we have to make sure that it is true whatever the prop\zl ositions\zd\ $p$ and $q$ may be\Big\rfloor} Here we have two prop.\zl ositional\zd\ var.\zl iables\zd\ $p,q$ and therefore \ul \sout{because all f\zl u\zd nct.\zl ions\zd\ involved are truth functions but} \ud \sout{are only} \underline{four possibilities} for these truth values\zl ,\zd\ namely
\begin{center}
\begin{tabular}{ c|c|c|c|c }
 $p$ & $q$ & $p\supset q$ & $p\, .\:(p\supset q)$ & $p\, .\:(p\supset q)\supset q$ \\[.5ex]
 T & T & T & T & T\\
 T & F & F & F & T\\
 F & T & T & F & T\\
 F & F & T & F & T
 \end{tabular}
\end{center}

\noindent $\mathbf{\llbracket 41. \: II \rrbracket}$ and what we have to do is simply to check that the truth value of the whole expr\zl ession\zd\ is true \sout{if} in each of these four cases\zl ,\zd\ \ul i\zl .\zd e.\ we have to ascertain that the truth table of the whole expression consists of T's only\zl .\zd\ \ud That\zl '\zd s very simple. Let us write down all the parts of which \sout{it} \ul this expr.\zl ession\zd\ \ud is built u\zl written over another symbol\zd p\zl .\zd\ We have first $p\supset q$ \sout{$p$} is a part\zl ,\zd\ then $p\, .\:(p\supset q)$ and finally the whole expr\zl ession\zd . \zl \sout{Now in the first case \ldots}\zd\ So we see \ul actually \ud in all four cases the \ul whole \ud formula is true\zl .\zd\ Hence it is universally true. It is clear that this purely mech\zl anical\zd\ method of checking all possibilities will always give a dec.\zl ision\zd\ whether a \zl given\zd\ formula is or is not a $\mathbf{\llbracket 42. \: II \rrbracket}$ tautologie\zl y\zd . Only if the nu\zl mber\zd\ of variables $p,q$ \ul occurring in the expr\zl ession\zd\ \ud is large this method is very cumbersome\zl ,\zd\ because the number of cases which we have to deal with is $2^n$ if the number of var.\zl iables\zd\ is $n$ \ul and the nu.\zl mber\zd\ of cases is the same as the nu.\zl mber\zd\ of lines in the truth table \ud . Here we had \sout{as 4} 2 var.\zl iables\zd\ $p,q$ and therefore $2^2=4$ cases. With 3 var.\zl iables\zd\ we would have $2^3=8$ cases and in gen.\zl eral\zd\ if the number of var\zl iables\zd\ is increased by one \sout{we} the number of cases \ul to be considered \ud is doubled, because each of the previous cases is split into two \ul new \ud cases according as \sout{to} the truth value of the new var\zl iable\zd\ is truth or falsehood. Therefore we have $\mathbf{\llbracket 43. \: II \rrbracket}$ \sout{that the} \zl \sout{number}\zd\ $2^n$ \ul cases \ud for $n$ variables\zl .\zd\ In the appl.\zl ications\zd\ however \ul usually \ud the number of cases actually to be considered is much smaller bec.\zl ause\zd\ mostly several cases can be combined into one\zl ,\zd\ e.g\zl .\zd\ in our ex.\zl ample\zd\ case 1 and 2 can be treated together bec\zl ause\zd\ if $q$ is true the whole exp\zl ression\zd\ is cert\zl ainly\zd\ true whatever \ul $p$ may be because it is then an impl.\zl ication\zd\ with true second member \ud \zl repeated line: because it is then an impl.\zl ication\zd\ with a true sec\zl ond\zd\ member\zd .

So we see that for the calc\zl ulus\zd\ of prop.\zl ositions\zd\ we have a very simple procedure to decide for any given formula whether or not it is log.\zl ically\zd\ true. \ul This solves the first of the two gen.\zl eral\zd\ problems which I ment\zl ioned\zd\ in the beginning \ul for the calc.\zl ulus\zd\ of prop.\zl ositions,\zd\ \ud namely \ul the probl.\zl em\zd\ \ud to give a complete theory of logically true form\zl ulas\zd . We have even more\zl ,\zd\ namely a procedure to decide of any form\zl ula\zd\ whether\ldots\ \zl at the beginning of this paragraph we find ``or not it is logically true''\zd\ \ud . That this problem $\mathbf{\llbracket 44. \: II \rrbracket}$ could be solved in such a simple way is chiefly due to the fact that we introduced only ext\zl ensional\zd\ operations (only truth funct.\zl ions\zd\ of prop.\zl ositions\zd )\zl .\zd\ If we had introd.\zl uced\zd\ strict impl\zl ication\zd\ the question would have been \sout{very} much more compl\zl icated\zd . It is only very recently \ul that \ud one has discovered general procedures for deciding whether a formula involving \sout{the not} strict implication are \zl should be ``is''\zd\ logically true under cert.\zl ain\zd\ assumptions about strict \zl implication.\zd\

Now after having solved this so called decision probl.\zl em\zd\ \zl For an insertion sign in the manuscript at this place no text to be inserted is given on this page.\zd\ for the calc.\zl ulus\zd\ of prop.\zl ositions\zd\ I can go over to the second probl.\zl em\zd\ I have announced in the beg.\zl inning.\zd\ \zl the text in the manuscript breaks at the end of this page after: namely the probl\zl em\zd\ of reducing\Big\rfloor\zd\

$\mathbf{\llbracket 56. \rrbracket}$ \zl Pages numbered from \textbf{45}., with or without II, up to \textbf{55}.\ are missing in the present notebook.\zd\ Now it has turned out that three rules of inf\zl erence\zd\ are sufficient for our purp\zl oses,\zd\ \ul namely for deriving all tautologies from these form.\zl ulas\zd\ \sout{\zl unreadable symbol\zd }\, \ud . Namely first the so called \zl in italics in the edited text: rule of subst.\zl itution\zd \zd\ which says: \uline{If we have a formula $F$} (of the calc.\zl ulus\zd\ of prop.\zl ositions\zd ) \uline{which involves the prop\zl ositional\zd\ var.\zl iables\zd\ say $p_1,\ldots,p_n$ then it is permissible to conclude from it any form.\zl ula\zd\ obtained \sout{obtained} \sout{from it} \ul by \ud sub\-st.\zl ituting\zd\ \ul in $F$ \ud for \ul all or some of \ud the prop.\zl osi\-tional\zd\ var.\zl i\-ables\zd\ \ul $p_1,\ldots,p_n$ \ud any arb.\zl itrary\zd\ expressions\zl ,\zd\ but in such a way that if a letter $p_i$ occurs in several places \ul in $F$ \ud we have to subst\zl itute\zd\ the same form.\zl ula\zd\ in all places where it occurs\zl .\zd } \zl E\zd .g.\ take the formula $(p\: .\: q\supset r)\supset [p\supset (q\supset r)]$ which is \zl called exportation.\zd According to the rule of subst.\zl itution\zd\ we can conclude \ul from it the for\zl mula\zd\ obt\zl ained\zd\ by subst\zl ituting\zd\ \ud $p\: .\: q$ for $r$\zl ,\zd\ i.e.\ \sout{\zl unreadable word ending in ``ing''\zd }\, $(p\: .\: q\supset p\: .\: q)\supset [p\supset (q\supset p\: .\: q)]$\zl .\zd\ The expression which we \ul substitute\zl ,\zd\ in our case $p\: .\: q$\zl ,\zd\ \ud is quite arbitrary $\mathbf{\llbracket 57. \rrbracket}$ \zl unreadable text, perhaps ``and it''\zd\ need not be a tautologie\zl y\zd\ \ul or a proved formula\zl .\zd\ \ud \zl \ul \sout{must be a \zl unreadable word\zd }\, \ud \zd\ The only requirement is that if the same letter occurs on several places in the formula in which we subst\zl itute\zd\ (as in out case the $r$) then we have to subst\zl itute\zd\ the same expression in all the places where $r$ occurs as we did here\zl .\zd\ \ul But it is perfectly allowable to subst\zl itute\zd\ for different letters the same formula\zl ,\zd\ e.g.\ for $q$ and $r$ and it is also allowable to subst.\zl itute\zd\ \sout{\zl unreadable crossed out word, perhaps ``for''\zd } expr\zl essions\zd\ containing variables which occur on some other places in the form.\zl ula,\zd\ as e.g.\ here $p\: \lfloor .\rfloor \: q$\zl .\zd\ \ud \sout{in our case} It is clear that by such a subst\zl itution\zd\ we get always a tautologie\zl y\zd\ if the \ul expr.\zl ession\zd\ \ud \sout{form\zl ula\zd }\, in which we subst\zl itute\zd\ is a taut.\zl ology,\zd\ because e.g.\ that this formul.\zl a\zd\ of export.\zl ation\zd\ is a tautol.\zl ogy\zd\ says \ul exactly \ud that it is true whatever $p$\zl , $q$,\zd $r$ may be\zl .\zd\ So it will in part.\zl icular\zd\ be true if we take \sout{$p\: .\: q$} for $r$ the prop\zl osition\zd\ $p\: .\: q$, whatever $p$ and $q$ may be $\mathbf{\llbracket 58. \rrbracket}$ and that means that the form.\zl ula\zd\ obt.\zl ained\zd\ by the \ul subst.\zl itution\zd\ \ud is a tautologie\zl y\zd .

\zl new paragraph\zd\ The sec.\zl ond\zd\ rule of inf.\zl erence\zd\ we need is the so called \zl in italics in the edited text: rule of impl.\zl ication,\zd \zd\ which reads as follows:
\begin{itemize}
\item[] \uline{If $P$ and $Q$ are arbitrary expressions then from the two premises $P$, $P\supset Q$ it is allowable to conclude $Q$.}
\end{itemize}
\zl A\zd n example\zl :\zd\ take for $P$ the form\zl ula\zd\ ${p\: \lfloor .\rfloor\: q\supset p\: \lfloor .\rfloor\: q}$ and \zl for\zd\ $Q$ the form\zl ula\zd\ $p\supset (q\supset p\: .\: q))$ \ul so that $P\supset Q$ will be the for\zl mula\zd\ $(p\: \lfloor .\rfloor\: q\supset p\: \lfloor .\rfloor\: q)\supset [p\supset (q\supset p\: \lfloor .\rfloor\: q)]$\zl .\zd\ \ud Then from those two premises we can conclude ${p\supset (q\supset p\: .\: q)}$\zl .\zd\ \ul Again we can prove that this rule of inf\zl erence\zd\ is corr.\zl ect,\zd\ i\zl .\zd e\zl .\zd\ if the two premises are tautologies then the concl\zl usion\zd\ \ul is \ud . Bec\zl ause\zd \zl i\zd f we assign any partic.\zl ular\zd\ truth values to the prop\zl osi\-tional\zd\ var.\zl iables\zd\ occurring in $P$ and $Q$\zl ,\zd\ $P$ and $P\supset Q$ will both get the truth value truth because they are taut\zl ologies\zd . Hence $Q$ will also get the truth value true if any part.\zl icular\zd\ truth values are assigned to its variables.
\zl B\zd ec.\zl ause\zd\ if $P$ and $P\supset Q$ both have the truth value truth $Q$ has also the truth \ul \zl \sout{according to}\zd\ \ud . So $Q$ will have the truth v.\zl alue\zd\ T whatever truth v.\zl alues\zd\ are ass.\zl igned\zd\ to the var.\zl iables\zd\ occur\zl r\zd ing in it which means that it is a tautol\zl ogy\zd . \ud

Finally as the third rule of inf\zl erence\zd\ we have the \uline{rule of defined sym\-b.\zl ol\zd } which \zl \sout{\ul (-) \ud }\zd\ says (roughly speak\zl ing\zd ) that within any form.\zl ula\zd\ the def\zl iniens\zd\ can be replaced by the definiendum \ul and vice versa\zl ,\zd\ \ud or formulated $\mathbf{\llbracket 59. \rrbracket}$ more precisely for a part.\zl icular\zd\ def.\zl iniens\zd\ say $p\supset q$ \sout{we had} it says\zl :\zd\ \uline{From a formula $F$ we can conclude any form\zl ula\zd\ $G$ obtained from $F$ by replacing \sout{for} a part of $F$ which has the form $P\supset Q$ by the expr.\zl ession\zd\ $\sim P\vee Q$ and vice versa}. (\zl S\zd imilarly for the other def\zl initions\zd\ we had\zl .\zd )

As \zl an\zd\ ex.\zl ample:\zd
\begin{tabbing}
\hspace{1.7em}\= 1. \= $ \sim p\vee (p\vee q)$ \quad from the first\zl axiom by replacing\zd\ $p \supset q\lfloor Q\rfloor$ by\\
\` ${\sim p\vee q\lfloor Q\rfloor}$\\
\>2.\> $\sim p\supset (\sim p \vee q)$ $\;\mid\;$\zl (\zd\ \underline{Again clear that taut\zl ology\zd\ of taut\zl ology\zd }.\zl )\zd\ \\[.5ex]
\>\> $\sim p\supset ( p \supset q)$
\end{tabbing}
This \ul last \ud rule is sometimes not explicitly form.\zl ulated\zd\ because it is only nec\zl essary\zd\ if one introd.\zl uces\zd\ def.\zl initions\zd\ and it is superfluous \zl in\zd\ principle \ul to introduce \zl them\zd\ \ud because whatever can be expr.\zl essed\zd\ by a defined symbol \zl can be done without\zd\ $\mid\;$\zl (\zd\ only it would sometimes be very long and cumbersome\zl )\zd . If however one introduces def.\zl initions\zd\ \ul as we did \ud this third rule of inf.\zl erence\zd\ is indispensable\zl .\zd\

\ul Now what we shall prove is that any taut\zl ology\zd\ can be derived from these four ax\zl ioms\zd\ by means of the \ul ment\zl ioned\zd\ \ud $3$\zl three\zd\ rules of inf.\zl erence:\zd\ \sout{but before proving this we shall first make sure that all the ax.\zl ioms\zd\ really are taut\zl ologies\zd\ and then give examples of derivations of indiv\zl idual\zd\ formulas from \zl crossed out ``it''\zd\ them} \ud
\begin{tabbing}
$\mathbf{\llbracket 60. \rrbracket}$\hspace{5em} \ul \zl unreadable text, probably in shorthand\zd\ \ud \\*[1ex]
\hspace{1.7em}\= 1.\zl (1)\zd\ \quad \= $p\supset p\vee q$\\*[.5ex]
\>2\;\zl (2)\zd\ \>$p \vee p \supset p$\\[.5ex]
\>3.\zl (3)\zd\ \>$p \vee q \supset q\vee p$\\[.5ex]
\>4.\zl (4)\zd\ \>$\lfloor ( \rfloor p\supset q)\supset (r\vee p\supset r\vee q)$
\end{tabbing}

Let us first asc.\zl ertain\zd\ that all of these form\zl ulas\zd\ \sout{form} are tautologies \ul and let us ascert.\zl ain\zd\ that \zl fa\zd ct first by their mean\zl ing\zd\ and then by their truth table\zl .\zd\ \ud \Big\lceil The first means: If $p$ is true $p\vee q$ is true. That is clear bec\zl ause\zd\ $p\vee q$ means at least one of the prop\zl ositions\zd\ $p,q$ is true\zl ,\zd\ but if $p$ is true then \zl the expression $p\vee q$ is true.\zd\ The sec\zl ond\zd\ means\zl :\zd\ If \ul the disj\zl unction\zd\ \ud $p\vee p$ is true $p$ is true\zl ,\zd\ i\zl .\zd e.\ \zl we k\zd now that the disj\zl unction\zd\ $p\vee p$ is true means that one of the two members is true\zl ,\zd\ but since both members are $p$ that means that $p$ is true. The third \sout{that} says\zl : I\zd f $p\vee q$ is true $q\vee p$ is also true\zl .\zd

\zl at the bottom of this page: Forts p 60 Heft II Anfang. \zl German: continued on p.\ \textbf{60}.\ Notebook II, beginning\zd \zd\

\zl On the last page of the present notebook, which is not numbered, one finds various notes, at the beginning and at the end consisting of formulae, which do not seem directly related to the preceding and succeeding pages of the course.\zd\
\begin{tabbing}
\hspace{1.7em}\= $p\supset\;$\=$q\; .\supset\;$\=$p\; .\equiv p$\hspace{7em}\= $p\: .\: q$, T, $p\supset q$, $q\supset p$, $p\vee q$\\[.5ex]
\>+\>+\>+\>$p\: .\: q\vee (p\supset q)\equiv$ \sout{$(p\supset q)$} T\\[.5ex]
\>+\>$-$\>+\>\sout{$\sim p \vee (p\: .\: q)$}\\[.5ex]
\>$-$\>$-$\>$-$\>\sout{\zl unreadable symbol\zd }\\[.5ex]
\>$-$\>+\>$-$
\end{tabbing}

\zl Then one finds in the middle part of that page four notes numbered 1.-4. written almost entirely in what seems to be shorthand, which contain perhaps exercises or examination questions. In the fourth of them one recognizes the following words that are not in shorthand: strict impl.\zl ication\zd , taut\zl ology\zd , if then, mat\zl erial\zd\ impl.\zl ication\zd , mod\zl us\zd\ pon\zl endo\zd\ pon\zl ens\zd .\zd\

\begin{tabbing}
\hspace{1.7em}\=T, F, $p$, $q$, $\sim p$, $\sim q$, $p\: .\: q$, $p\: . \sim q$, $\sim p\: .\: q$, $\sim p\: . \sim q$\\[.5ex]
\>$p\vee q$, $p\supset q$, $q\supset p$, $p\equiv q$, $p \mid q$, $p\equiv \:\sim q$
\end{tabbing}

\section{Notebook II}\label{0II}
\pagestyle{myheadings}\markboth{SOURCE TEXT}{NOTEBOOK II}
\zl Folder 60, on the front cover of the notebook ``Log.\zl ik\zd\ Vorl.\zl esungen\zd\ \zl German: Logic Lectures\zd\ Notre\zl $\,$\zd Dame II''\zd\

\vspace{1ex}

\zl Before p.\ \textbf{61}.\ one finds on a page not numbered the formulae
\begin{tabbing}
\hspace{1.7em}\=$p\: .\: q\, .\! \vee(\sim p\: . \sim q)$\\*[.5ex]
\>$(p\vee q) \, .\:{(\sim p\:\vee \sim q)}$
\end{tabbing}
and a few scattered letters and symbols from partly missing unreadable formulae.\zd\

\vspace{2ex}

$\mathbf{\llbracket 61. \rrbracket}$ This does\zl $\,$\zd not need further explan\zl ation\zd\ because the \zl ``\zd or\zl ''\zd\ is evidently sim\zl symmetric\zd\ in the two members. Finally the fourth means \ul this\zl : ``\zd If $p\supset q$ then if $r\vee p$ is true \zl then $r\vee q$ is also\zd\ true \zl '',\zd\ i.e.\ \ud \zl ``\zd I\zl written over i\zd f you\zl written over ``we''\zd\ have a correct impl.\zl ication\zd\ $p\supset q$ then you can get again a corr.\zl ect\zd\ \ul impl.\zl ication\zd\ \ud by adding a third prop.\zl osition\zd\ $r$ to both sides of it getting $r\vee p\supset r\vee q$\zl ''.\zd\

\indent\zl The following text in big square brackets is crossed out in the manuscript: \Big[That \zl \sout{has a very analogie} is very analogous\zd\ to the laws by which one calculates with equations or inequalities in math.\zl ematics,\zd\ e.g.\ from $a<b$ you can conclude $c+a<c+b$\zl ,\zd\ \ul i.e.\ it is allowable to add an \ul arb.\zl itrary\zd\ \ud number to both sides of an inequality and likewise it is allowable to \zl add?\zd\ a prop\zl osition\zd\ to both sides of an impl\zl ication.\zd\ \ud \Big] \zd\

That this \ul is so \ud can be seen like this\zl :\zd\ it means \zl ``\zd If $p\supset q$ then if one of the prop\zl ositions\zd\ \zl\sout{(}\zd $r$\zl written over $p$\zd , $p$\zl\sout{)}\zd\ is true then also one of the prop\zl ositions\zd\ $r,q$ \zl is true'',\zd\ which is clear bec.\zl ause\zd\ if $r$ is true $r$ is true and if $p$ is true $q$ is true by ass\zl umption.\zd\ So whichever of the two prop\zl ositions\zd\ $r,p$ is true always it has the cons.\zl equence\zd\ that one of the prop\zl ositions\zd\ $r,q$ is true\zl .\zd \Big\rfloor

$\mathbf{\llbracket 62. \rrbracket}$ \Big\lceil Now let us ascertain the truth of these form\zl ulas\zd\ by the truth\zl -\zd table method, combining always as many cases as possible into one case\zl .\zd\
\begin{itemize}
\item[1.] If $p$ is F this is an impl\zl ication\zd\ with a false first member\zl ,\zd\ hence true owing to the truth\zl $\,$\zd table of $\supset$\zl ;\zd\ if $p$ is true then $p\vee q$ is also true acc\zl ording\zd\ to the truth \zl table\zd\ of \zl ``\zd or\zl '',\zd\ hence the form\zl ula\zd\ is an implic\zl ation\zd\ with true sec.\zl ond\zd\ memb.\zl er,\zd\ hence \sout{true} again true\zl .\zd\
\vspace{-1ex}
\item[2.] If $p$ is true this will be an impl\zl ication\zd\ with true sec\zl ond\zd\ mem\zl ber,\zd\ hence true\zl .\zd\ If $p$ is false then $p\vee p$ is a disj.\zl unction\zd\ both of whose memb\zl ers\zd\ are false\zl ,\zd\ hence false acc.\zl ording\zd\ to the truth \zl table\zd\ for $\vee$\zl .\zd\ Hence in this case we have an impl\zl ication\zd\ with $\mathbf{\llbracket 63. \rrbracket}$ a false first member, which is true by the truth\zl $\,$\zd table of \sout{or} $\supset$\zl .\zd\
\vspace{-1ex}
\item[3.] Since the truth\zl $\,$\zd t\zl able\zd\ for $\vee$ is sy\zl written over i\zd m\zl m\zd etric in $p,q$ it is clear that whenever the left\zl -\zd hand side has the truth value true also the right\zl -\zd h\zl and\zd\ side \zl has it,\zd\ and if the left\zl -\zd hand side is false the right\zl -\zd h\zl and\zd\ side will also be false\zl ;\zd\ but an impl.\zl ication\zd\ both of whose mem\zl bers\zd\ are true or both \zl of\zd\ whose \zl members\zd\ are false is true by the truth\zl $\,$\zd t\zl able\zd\ of impl.\zl ication,\zd\ bec.\zl ause\zd\ $p\supset q$ \zl is\zd\ false \zl only in?\zd\ the case when $p$ is true and $q$ false\zl .\zd\
\vspace{-1ex}
\item[4.] Here we have to consider only the foll.\zl owing\zd\ three cases\zl :\zd\
\end{itemize}
\vspace{-2ex}
\begin{tabbing}
\hspace{5em}\=1. one \zl of\zd\ \sout{the tw\zl o?\zd }\, $\, r$, $q$ \zl \sout{has the} has the\zd\ truth\zl $\,$\zd $\overline{\rm{v}}$.\zl value\zd \\
\`   T \zl ``th'' in ``truth'' written over ``e''\zd \\[.5ex]
\>2. both $r,q$ \zl are\zd\ F and $p$ true \\[.5ex]
\>3. \zl ditto marks interpreted as ``both $r,q$ are''\zd\ F and $p$ false
\end{tabbing}

$\mathbf{\llbracket 64. \rrbracket}$ These three cases evid.\zl ently\zd\ exhaust all poss\zl ibilities\zd .
\begin{itemize}
\item[\underline{1.}] \zl I\zd n the first case $r\vee q$ \zl unclear sign\zd\ is true, hence also $(r\vee p)\supset(r\vee q)$ is true bec.\zl ause\zd\ it is an impl\zl ication\zd\ with \zl \sout{false}\zd\ sec.\zl ond\zd\ memb\zl er true;\zd\ $(p\supset q)\supset(r\vee p\supset r\vee q)$ is true for the same reas\zl on.\zd\
\vspace{-1ex}
\item[\underline{2.}] \zl I\zd n the 2.\zl second\zd\ case $p$ \zl is true and\zd\ $q$ false\zl ,\zd\ hence $p\supset q$ false\zl ,\zd\ hence the whole expr\zl ession\zd\ is an impl\zl ication\zd\ with false first mem\-ber\zl ,\zd\ hence true\zl .\zd\
\vspace{-1ex}
\item[\underline{3.}] \zl In the\zd\ 3\zl third\zd\ case all \zl unreadable text, perhaps: of them all, should be: of $r$, $q$ and $p$\zd\ are false\zl ;\zd\ then \zl unreadable text, should be:\ $r\vee p$ and $r\vee q$\zd\ are false\zl ,\zd\ hence \zl the\zd\ \zl unreadable text, perhaps: impl\zl ication\zd , should be: $r\vee p\supset r\vee q$ is\zd\ true, hence \zl the\zd\ whole form\zl ula is\zd\ true bec\zl ause it is an\zd\ impl\zl ication\zd\ with true sec\zl ond\zd\ member\zl .\zd\
 \item[] \underline{So we see that the whole formula is always true.}
\end{itemize}

Now I can begin with deriving other taut.\zl ologies\zd\ from these 3\zl three\zd\ ax\zl ioms\zd\ by means of the two\zl three\zd\ rules of inf.\zl erence,\zd\ namely \zl the\zd\ rule of \uline{subst\zl itution\zd\ and implication \ul and def.\zl in\-ed\zd\ symb\zl ol,\zd\ \ud} in order to prove later on that all logic.\zl ally\zd\ true form\zl ulas\zd\ can be derived from them\zl .\zd\

Let us first \zl substitute\zd\ \underline{$\f{\sim r}{r}\;$} \sout{(4)} in 4\zl (4) to get\zd\ \underline{$(p\supset q)\supset(\sim r\vee p\supset$} \underline{$\; \sim r\vee q)$} \zl ,\zd\ but for $\sim r\vee p$ we can subs.\zl titute\zd\ $r\supset p$ and likewise for $\sim$ \zl $\sim r\vee q,$\zd\ $\mathbf{\llbracket 65. \rrbracket}$ \sout{so this means the same thing as} getting: \zl Some of the figures of the numbered formulae below are written over other symbols, not always recognizable, but sometimes they are, as for example with 7. and 8*, which are written over 3 and 4.\zd\
\begin{tabbing}
5. \underline{$(p\supset q)\supset[(r\supset p)\supset (r\supset q)]$} \quad Syl.\zl logism\zd\
\end{tabbing}
This is the so called form\zl ula\zd\ of syllog\zl ism,\zd\ which has a cert\zl ain\zd\ simil.\zl ari\-ty\zd\ to \ul the \ud mood \zl B\zd arbara in so far as it says: If from $p$ follows $q$ then if from $r$ foll\zl ows\zd\ $p$ from $r$ foll\zl ows\zd\ $q$.

\noindent 6. Now subst.\zl itute\zd\ \underline{$\f{p}{q}$ in (1)} \underline{$p\supset p\vee p$} and now make the foll.\zl owing\zd\ subst\-\zl itution:\zd\ \sout{in Syl.\zl logism}\zd\
\begin{tabbing}
\hspace{1.7em}$\f{\f{p\vee p}{p}\quad\f{p}{q}\quad\f{p}{r} \quad \rm{in \;\; Syl\lfloor logism\rfloor}}{(p\vee p\supset p)\supset[(p\supset p\vee p)\supset (p\supset p)]}$
\end{tabbing}
\zl T\zd his is an impl.\zl ication\zd\ and the first memb.\zl er\zd\ of it reads $p\vee p\supset p$\zl ,\zd\ which is nothing else but the \zl \sout{first} second\zd\ ax\zl iom.\zd\ Hence we can apply the rule of imp.\zl lication\zd\ to the $\mathbf{\llbracket 66. \rrbracket}$ two prem\zl ises\zd\ and get
\begin{tabbing}
\hspace{1.7em}$(p\supset p\vee p)\supset (p\supset p)$
\end{tabbing}
This is again an impl.\zl ication\zd\ and the first memb\zl er\zd\ of it was proved before\zl ;\zd\ hence we can again apply the rule of implic.\zl ation\zd\ and get
\begin{tabbing}
7. \underline{$p\supset p$}\quad law of identity
\end{tabbing}
\zl U\zd sing the third rule
\begin{tabbing}
8*\zl . we\zd\ have \underline{$\sim p\vee p$}\quad the law of excl.\zl uded\zd\ middle
\end{tabbing}
Now let us subst\zl itute\zd\ \underline{$\f{p}{\sim p}$ \zl $\f{\sim p}{p}$\zd\ in} this form.\zl ula to get\zd\ {$\sim\sim p \;\vee \sim p$} and now apply to \ul it the \ud com.\zl mutative\zd\ law for $\vee$\zl ,\zd\ i\zl .\zd e.\ subst\zl itute\zd\ \underline{$\f{\sim \sim p}{p}$\quad $\f{\sim p}{q}\;\;$ \zl in\zd } \underline{\underline{(3)}} \zl to get\zd\
\begin{tabbing}
\hspace{1.7em}\=$\sim\sim p \;\vee \sim p\supset \;\sim p\;\vee \sim\sim p$ \quad rule of impl.\zl ication\zd \\*[.5ex]
\>\underline{$\sim p\;\vee \sim\sim p$}
\end{tabbing}
$\mathbf{\llbracket 67. \rrbracket}$ 9.* \underline{$p\supset\;\sim\sim p$}

\vspace{1ex}

\zl The following inserted text is crossed out in the manuscript: \ul Here we have some ex.\zl amples\zd\ of form.\zl ulas\zd\ derived from ax.\zl ioms\zd\ by rules of inf.\zl erence;\zd\ form.\zl ulas\zd\ for which this is the case I call demonstr.\zl able\zd\ (\zl unreadable text\zd\ from the 4 ax.\zl ioms,\zd\ but I leave that expl.\zl unreadable text\zd )\ud \zd \zl The following inserted text from the manuscript is deleted: \ul So these form.\zl ulas\zd\ are dem\zl onstrable\zd\ $|$\underline{before going on \zl unreadable text\zd }$|$ \ud \zd

\vspace{1ex}

\sout{Now} I have to make an imp.\zl ortant\zd\ remark \zl unreadable text\zd\ \zl on\zd\ how we ded.\zl uced\zd\ $p\supset p$ from the ax\zl ioms\zd . We had at first the two \zl formulas $p\supset p\vee p$ and $p\vee p\; \supset p$. Now\zd\ subst.\zl itute\zd\ them in a certain way in the form.\zl ula of S\zd  yllog\zl ism\zd\ \ul $\f{p}{r}\;\;$ $\f{p\vee p}{p}\;\;$ $\f{p}{q}$\ud and then by appl\zl ying\zd\ twice the rule of impl\zl ication\zd\ we get $p\supset p$\zl .\zd\ \zl unreadable scarcely visible text in more than two lines, where one can recognize the words: this, the two, not, to these two, $p\vee p$ but\zd\ \ul If $P,Q,R$ are any arb.\zl itrary\zd\ expr.\zl essions\zd\ and \ud \uline{if we have succeeded in deriving $P\supset Q$ \zl and\zd\ $Q\supset R$ from the four ax\zl ioms\zd\ by means of the two\zl three\zd\ rules of proc\zl should be: inference\zd\ then we can also derive} $\mathbf{\llbracket 68. \rrbracket}$ \underline{$P\supset R$}\zl .\zd\ Because we can simply subst\zl itute\zd\ $\f{P}{p}$\quad $\f{Q}{q}$\quad $\f{R}{r}$ \zl $\f{Q}{p}$\quad $\f{R}{q}$\quad $\f{P}{r}$\zd\ in Syl.\zl logism\zd\ getting $(Q\supset R)\supset[(P\supset Q)\supset (P\supset R)]$. Then we apply the rule of impl\zl ication\zd\ to this form\zl ula\zd\ and $P\supset Q$ \zl $Q\supset R$\zd\ getting \ul \ldots \zl ${(P\supset Q)}\supset (P\supset R)$\zd\ \ud and then we apply \ul again the rule of \ud impl\zl ication\zd\ to this form\zl ula\zd\ \zl unreadable text; should be: and $P\supset Q$\zd\  \zl gett\zd ing $P\supset R$\zl .\zd\

\uline{So we know \ul \zl quite\zd\ \ud generally if $P\supset Q$ and $Q\supset R$ are both \sout{provable} \ul demonstrable \ud then \sout{also} $P\supset R$ is \sout{provable} \ul also demonstrable \ud whatever for\-mula $P,Q,R$ may be} \ul because we can obtain $P\supset R$ always in the manner just described\zl .\zd\ \sout{and} \zl T\zd his \sout{cogn.} \ul fact \ud allows us to save the trouble of repeating the whole arg\zl ument\zd\ by which we derived the concl.\zl usion\zd\ from the two prem.\zl ises\zd\ in each part.\zl icular\zd\ case, but we can state it once for all as a new \ul \sout{Since we know that it can \zl unreadable text\zd\ and how it can be done i.e} \ud $\mathbf{\llbracket 69. \rrbracket}$ rule of inf\zl erence\zd\ as foll\zl ows\zd :
\begin{itemize}
\item[]\uline{From the two prem\zl ises\zd\ $P\supset Q$, $Q\supset R$ we can conclude \ul $P\supset R$ \ud \ul what\-ever the form\zl ulas\zd\ $P,Q,R$ may be\zl .\zd } \quad \underline{4.R.}
\end{itemize}
So this is a $4^{\rm th}$\zl fourth\zd\ rule of inf.\zl er\-ence,\zd\ which I call Rule of syllogism\zl .\zd\ \ud But note \ul that \ud this rule of syllogism is not a new indep.\zl endent\zd\ rule, but can be derived from the other two\zl three\zd\ rules and the 4\zl four\zd\ axioms. Therefore it is called a derived rule of inf\zl erence\zd . So we see that \zl unreadable text, should be: in\zd\ our syst\zl em\zd \zl unreadable text\zd\ we cannot only derive formul.\zl as\zd\ but also new rules of inf.\zl erence\zd\ and the \ul latter \ud is very helpful for shortening the proofs\zl .\zd\ O\zl I\zd n principle it is \sout{of course} superfluous to introduce such derived rules of inf\zl erence\zd\ because whatever can be proved with their help can also be proved without them. It is exactly this what we have shown before introd\zl ucing\zd\ \sout{\zl unreadable text\zd }\, this \ul new \ud rule of inf.\zl erence\zd , namely we have shown that the conclusion of it can be obtained also by the former axioms and rules of inf.\zl erence\zd\ and this was the justification for introducing it.

\zl new paragraph\zd\ $\mathbf{\llbracket 70. \rrbracket}$ But although these \ul derived \ud rules of inf\zl erence\zd\ are superfl\zl uous\zd\ \ul o\zl i\zd n principle \ud they are very helpful for shortening the proofs and \sout{shorter} \ul therefore we shall introduce a great \zl many\zd\ of them\zl .\zd\ \zl We now\zd\ apply this rule immediately to \zl unreadable text, could be: the 1\zl (1)\zd\ and 3\zl (3)\zd\ ax\zl ioms\zd \zd\ because they have this form $P \supset Q$\zl ,\zd\ $Q\supset R$ \zl unreadable text, should be: for\zd\ $\f{p}{P}$\quad $\f{p\vee q}{Q}$ \quad$\f{q\vee p}{R}$\zl ,\zd\ and get \zl unreadable symbol or left parenthesis\zd\ bec\zl ause\zd\ (1)\zl ,\zd\ (3)
\begin{tabbing}
10.* \underline{$p\supset q\vee p$}\\[.5ex]
parad\zl ox\zd : \= 11. \= \underline{$p\supset (q\supset p)$}\quad\quad $p\supset (\sim q\vee p)$\\[.5ex]
\>\> \ul Add$\,$* \zl crossed out what seems to be: $\f{q}{p}$\zd\ $\f{\sim q}{q}$ \zl written\\
\` over something unreadable\zd\ \zl in\zd\ last formula (10)\zl 10.*\zd\ \ud \\[1ex]
\> 12. [\underline{$\sim p\supset (p\supset q)$}\quad\quad $\sim p\supset (\sim p\vee q)$\\[.5ex]
\>\> \ul Add $\f{\sim p}{p}$ \quad $\f{q}{q}$] \, \zl in\zd\ \, (1) \ud
\end{tabbing}

Other derived rules\zl :\zd\
\begin{tabbing}
4$\cdot 1'$R $\;$\=\underline{$P_1\supset P_2\;\; P_2\supset P_3\;\; P_3\supset P_4 \;\; : \;\; P_1\supset P_4$} \underline{generalized rule of syll.\zl ogism\zd }\\[.5ex]
\hspace{6em} $P_1\supset P_3$\\[1ex]
5.R* \>\underline{$P\supset Q\;\; : \;\; R\vee P\supset R \vee Q$} \quad \zl addition from the left\zd\
\end{tabbing}
This rule is sim\zl ilar\zd\ to the rules by which one calc.\zl ulates\zd\ with inequ.\zl alities\zd

\vspace{.5ex}

\hspace{7em} $a<b$\quad\zl :\zd\ \quad $c+a<c+b$\zl .\zd\ \zl \sout{ie.}\zd\

\vspace{.5ex}

\noindent \zl The name of the formula in the following deleted text is transferred next to 5.R* above: Call\zl ed\zd\ i\zl .\zd e. \underline{addition from \ul the \ud left}.\zd
\begin{tabbing}
4$\cdot 1'$R $\;$\=\underline{$P_1\supset P_2\;\; P_2\supset P_3\;\; P_3\supset P_4 \;\; : \;\; P_1\supset P_4$}\quad \underline{generalized rule of syll.\zl ogism\zd }\kill

[6R \> \underline{$P\supset Q\;\; : \;\; (R\supset P)\supset (R \supset Q)$}\;]\\[1ex]
5$\cdot 1$$_R$\zl R\zd *\> \hspace{1.8em}\underline{$P\supset Q$}\;\; : \;\; \underline{$P\vee R\supset Q \vee R$}\quad add.\zl ition\zd\ from \zl the\zd\ right\\[1ex]
$\mathbf{\llbracket 71. \rrbracket}$ \> 1. \hspace{.4em}\= $P\vee R\supset R\vee P$\hspace{3em}\= $\f{P}{p}$\quad $\f{R}{q}$ \quad in (3.)\zl (3)\zd\ \\[.5ex]
\> 2.\>$R\vee P\supset R\vee Q$\>by rule add.\zl ition from the\zd\ left\\[.5ex]
\> \underline{3.\;\, $R\vee Q\supset Q\vee R$} \>\> $\f{R}{p}$\quad $\f{Q}{q}$ \quad in (3.)\zl (3)\zd\ \\[.5ex]
\>\>$P\vee R\supset Q\vee R$\>by rule Syllog.\zl ism\zd \\[1ex]
7R* \>\underline{$P\supset Q\quad\quad R\supset S\;\; : \;\; P\vee R\supset Q \vee S$}\\[.5ex]
 \`Rule of addition of two impl\zl ications\zd \\[.5ex]
\>\>$P\vee \mbox{\underline{\emph{R}}}\supset Q\vee \mbox{\underline{\emph{R}}}$\> add\zl ition\zd\ from \zl the\zd\ \=right \zl to the\zd \\
\` first premis\zl \sout{s}\zd e ($R$)\\[.5ex]
\>\>\underline{$\mbox{\underline{\emph{Q}}}\vee R\supset \mbox{\underline{\emph{Q}}}\vee S$}\> \zl addition from the\zd\ left $\; ''\; $ sec\zl ond\zd\ $\; ''\; $ ($Q$)\\[.5ex]
\>\>$P\vee R\supset Q\vee S$\> \zl S\zd yllog\zl ism,\zd\ but this is the conclusion to\\
\` be proved\\[1ex]
8R* \>\underline{$P\supset Q\quad\quad R\supset Q\;\; : \;\; P\vee R\supset Q$}\quad\quad Dilemma\\[.5ex]
\>\>$P\vee R\supset Q\vee Q$\\[.5ex]
\>\>\underline{$Q\vee Q\supset Q$}\>$\f{Q}{p}$ \zl in\zd\ (2)\\[.5ex]
\>\>$P\vee R\supset Q$\> \zl S\zd yll\zl ogism\zd \\[1ex]
$\mathbf{\llbracket 72. \rrbracket}$ Applic.\zl ation\zd\ to derive formulas\\[.5ex]
\>\> $p\supset \;\sim\sim p$ \> proved before\zl ,\zd\ subst.\zl itute\zd\ $\f{\sim p}{p}$\\[.5ex]
\>\> $\sim p\supset \;\sim\sim\sim p$ \> ad\zl dition from the\zd\ right\\[.5ex]
\>\> $\sim p\vee p\supset \;\sim\sim\sim p\vee p$\> \quad rule of impl\zl ication\zd \\[.5ex]
\>\> $\sim\sim\sim p \vee p$\> \sout{\zl abbr\zd }\, rule of def\zl ined\zd\ symb.\zl ol\zd \\[1ex]
13.\> \underline{$\sim\sim p\supset p$}\\[1ex]
14\zl .\zd \>\zl Trans\zd pos\zl ition\zd\ \zl \sout{\underline{to d}}\zd\ \underline{$(p\supset\;\sim q)\supset(q\supset\;\sim p)$} \quad\sout{to} \\[.5ex]
\>Proof\>\hspace{3em}$(\sim p\;\vee\sim q)\supset(\sim q\;\vee\sim p)$ \hspace{1.9em} subst\zl itution\zd\ in (3)\\
\` rule \zl of defined\zd\ symb\zl ol\zd \\[1ex]
14$\cdot$1\>\underline{$(p\supset q)\supset(\sim q\supset\;\sim p)$}\\[.5ex]
\>$(\sim p\vee q)\supset(\sim\sim q\;\vee\sim p)$\\[.5ex]
Proof\>\> $q\supset \;\sim\sim q$\\[.5ex]
\>\>$\sim p\vee q\supset \;\sim p\;\vee \sim\sim q$\\[.5ex]
\>\>\underline{$\sim p\;\vee \sim\sim q\supset \;\sim\sim q\;\vee\sim p$}\quad Perm\zl utation\zd\ (3)\\[.5ex]
\>\>$\sim p\vee q\supset\;\sim\sim q\;\vee\sim p$\quad \sout{q.e.d.} \quad \zl rule of\zd\ def.\zl ined\zd\ symb.\zl ol\zd \\[1ex]
14$\cdot$1\>$(p\supset q)\supset(\sim q\supset\;\sim p)$\>\>\quad$\sim p\vee q\supset\;\sim\sim q\;\vee\sim p$\\[.5ex]
14$\cdot$2\>$(\sim p\supset\;\sim q)\supset(q\supset p)$\>\>\quad$\sim\sim q\;\vee \sim p\supset\;\sim p\vee q$\\[.5ex]
14$\cdot$3*\>$(p\supset\;\sim q)\supset(q\supset\;\sim p)$\>\>\quad\underline{$|$$\sim p\;\vee\sim q\supset\;\sim q\;\vee\sim p$$\:|$}\\[.5ex]
14$\cdot$4*\>$(\sim p\supset q)\supset(\sim q\supset p)$\>\>\quad$\sim\sim p\vee q\supset\;\sim\sim q\;\vee p$\\[1.5ex]
\zl 14$\cdot$2\zd \>\underline{$(\sim p\supset\;\sim q)\supset(q\supset p)$}\\[1ex]
\underline{Proof}\>$\sim\sim p\supset p$\\[.5ex]
\>$\sim\sim p\;\vee\sim q\supset p\;\vee\sim q$\\[.5ex]
\>\underline{$p\;\vee\sim q\supset\;\sim q\vee p$}\\[.5ex]
\>$\sim\sim p\;\vee\sim q\supset\;\sim q\vee p$\\[.5ex]
\>$\lfloor(\rfloor\sim p\supset\;\sim q\lfloor)\rfloor\supset\lfloor(\rfloor q\supset p\lfloor)\rfloor$
\end{tabbing}

\begin{tabbing}
$\mathbf{\llbracket 73. \rrbracket}$\\*[1ex]
4$\cdot 1'$R $\;$\=\underline{$P_1\supset P_2\;\; P_2\supset P_3\;\; P_3\supset P_4 \;\; : \;\; P_1\supset P_4$}\quad \underline{generalized rule of syll.\zl ogism\zd }\kill

14$\cdot$2\zl 14$\cdot$4*\zd\ \quad\underline{$(\sim p\supset q)\supset(\sim q\supset p)$}\\[.5ex]
\>$\sim\sim p\vee q\supset\;\sim\sim q\vee p$\\[1ex]
Proof\>$\sim\sim p\supset p$\\[.5ex]
\>\underline{$q\supset\;\sim\sim q$}\\[.5ex]
\>$\sim\sim p\vee q\supset p\;\vee\sim\sim q$\quad Add.\zl ition of two implications\zd \\[.5ex]
\>\underline{$p\;\vee\sim\sim q\supset\;\sim\sim q\vee p$}\quad Perm\zl utation\zd \\[.5ex]
\>$\sim\sim p\vee q\supset\;\sim\sim q\vee p$\quad q.e.d.\quad rule of def\zl ined\zd\ symb.\zl ol\zd \\[1.5ex]
4\zl Four\zd \zl unreadable text, perhaps: transposition\zd\ rules of inf\zl erence:\zd \\[1ex]
9$_R$\zl R\zd \>\hspace{1.5em}\underline{$P\supset\; \sim Q\quad : \quad Q\supset\; \sim P$}\hspace{2.5em}\=9$\cdot$1$_R$\zl R\zd  \hspace{1em}\=\underline{$P\supset Q\quad : \quad \sim Q\supset\; \sim P$}\\[.5ex]
9$\cdot$2$_R$\zl R\zd \>\hspace{1.5em}\underline{$\sim P\supset Q\quad : \quad \sim Q\supset P$}\hspace{3em}\>\zl 9$\cdot$3R\zd\ \>\underline{$\sim P\supset\;\sim Q\quad : \quad Q\supset P$}
\end{tabbing}
By \zl them\zd\ the laws for . corresp\zl ond\zd\ to laws for $\vee$ or can be derived\zl ,\zd\ e\zl .\zd g.\
\begin{tabbing}
15.*\hspace{1.7em}\=\underline{$p\: .\: q\supset p$}\hspace{7em}\=\underline{$p\: .\: q\supset q$}\\[.5ex]
\>$\sim(\sim p\;\vee\sim q)\supset p$\>$\sim(\sim p\;\vee\sim q)\supset q$\qquad\=Form\zl ula\zd\ 10\zl 10.*\zd\ \\[1ex]
Proof\>$\sim p\;\supset\;\sim p\;\vee\sim q$\>$\sim q\;\supset\;\sim p\;\vee\sim q$\>Transp\zl osition\zd\ 2.\\
\` \zl 9$\cdot$2R\zd\ \\[.5ex]
\>$\sim(\sim p\;\vee\sim q)\supset p$\>$\sim(\sim p\;\vee\sim q)\supset q$
\end{tabbing}
\ul 15.2\quad Similarly for prod.\zl ucts\zd\ of any number of fact.\zl ors\zd\ we can prove that the prod.\zl uct\zd\ implies any fact.\zl or,\zd\ e.g.\
\begin{tabbing}
\hspace{1.7em}\=\underline{$p\: .\: q\, .\: r\supset p$}\hspace{2em}\= bec.\zl ause\zd\ $\;\;$\=$(p\: .\: q)\, .\: r\supset p\: .\: q$\\[.5ex]
\>\underline{$p\: .\: q\, .\: r\supset q$}\>\>$p\: .\: q\supset p$\zl ,\zd \qquad $p\: .\: q\supset q$\\[.5ex]
\>\underline{$p\: .\: q\, .\: r\supset r$}\>\>$(p\: .\: q)\, .\: r\supset r$
\end{tabbing}
and for any numb.\zl er\zd\ of fact\zl ors\zd . \ud

\ul \Big\lceil From this \zl one has\zd\ the \zl following rules of inference:\zd\
\begin{tabbing}
10$_R$\zl R\zd \hspace{1em}\=\underline{$\;P\supset Q\quad : \quad P\: .\:R\supset Q$}\hspace{1em} adjoining \zl a\zd\ new hyp.\zl othesis\zd \\[.5ex]
10$\cdot$1$_R$\zl R\zd\>\underline{$\;P\supset Q\quad : \quad R\: .\:P\supset Q$}\\[.5ex]
\>bec\zl ause\zd\ \hspace{1em}\=$P\: .\: R\supset P$\hspace{1.7em}\= by subst\zl itution\zd \\[.5ex]
\>\>\underline{$P\supset Q$}\>by ass\zl umption\zd \\[.5ex]
\>\>$P\: .\: R\supset Q$\>\zl S\zd yll.\zl ogism\zd \\[1ex]
\Big\lfloor10$\cdot$2$_R$\zl R\zd\>$\;\;\;$\underline{$Q$}\quad : \quad \underline{$P\supset Q$}\>\> from paradox \ud
\end{tabbing}

\begin{tabbing}
$\mathbf{\llbracket 74. \rrbracket}$ Assoc\zl iativity\zd\ bes\zl ?\zd : Recall \underline{$|$I.\zl (1)\zd\ \; $p\supset p\vee q$\zl ,\zd\ $\;$ II \; $p\supset q\vee p$$\:|$}\\*[1ex]
15.*\quad\=\underline{$(p\vee q)\vee r\supset p\vee(q\vee r)$}\\[.5ex]
\>1.\quad\=$p\supset p\vee(q\vee r)$\qquad Add\zl ition\zd\ (1)\quad $\f{q\vee r}{q}$\\[.5ex]
\>\>$q\supset q\vee r$\qquad $q\vee r\supset p\vee(q\vee r)$\qquad Form\zl ula\zd\ 10\zl 10.*\zd\ \\
\` \ul Add.\zl ition\zd * \quad $\f{q\vee r}{p}$\quad$\f{p}{q}$ \qquad($p\supset q\vee p$\quad $\f{q\vee r}{p}$\quad$\f{p}{q}$) \ud \\[.5ex]
\>2.\>$q\supset p\vee(q\vee r)$\qquad \zl S\zd yll\zl ogism\zd \\[.5ex]
\>a.)\> $p\vee q\supset p\vee(q\vee r)$\qquad Dilemma\\[.5ex]
\>\>$r\supset q\vee r\;$ \quad (II $\f{r}{p}$)\qquad $q\vee r\supset p\vee(q\vee r)$\quad \zl see\zd\ \underline{before}\\
\` \zl an arrow is drawn from \underline{before} to the same formula three lines above\zd \\[.5ex]
\>b.)\> $r\supset p\vee(q\vee r)$\\[.5ex]
\>\>$(p\vee q)\vee r\supset p\vee(q\vee r)$\qquad inverse similar\\[1ex]
15$\cdot$1\>\underline{$p\vee(q\vee r)\supset (p\vee q)\vee r$}\\
\hspace{1.2em}\ul \>$p\supset p\vee q$\qquad $p\vee q\supset (p\vee q)\vee r$\qquad ($p\supset p\vee q$\quad $\f{p\vee q}{p}$\quad $\f{r}{q}$)\\[.5ex]
\>$p\supset (p\vee q)\vee r$\\[.5ex]
\>$q\supset (p\vee q)\vee r$\\
\>$r\supset (p\vee q)\vee r$\qquad[II \quad $p\supset q\vee p$\quad $\f{r}{p}$\quad $\f{p\vee q}{q}$]\\[.5ex]
\>$q\vee r\supset (p\vee q)\vee r$\\[.5ex]
\>$p\vee(q\vee r)\supset (p\vee q)\vee r$\quad \ud \\[1ex]
Export\zl ation\zd\ and import.\zl ation\zd \\[.5ex]
16.*\>\underline{$(p\: .\: q\supset r)\supset[p\supset(q\supset r)]$}\qquad Export\zl ation\zd \\[1ex]
$\mathbf{\llbracket 75. \rrbracket}$\>$(\sim(p\: .\: q)\vee r)\supset\:\sim p\vee(\sim q\vee r)$\\*[.5ex]
\>$\sim\sim(\sim p\;\vee\sim q)\vee r\supset\;\sim p\vee(\sim q\vee r)$\\[1ex]
Proof\>$\sim\sim(\sim p\;\vee\sim q)\supset\;\sim p\;\vee\sim q$\qquad double neg.\zl ation\zd \\
\` subst\zl itute\zd \quad $\f{\sim p\;\vee\sim q}{p}$\\[.5ex]
\>$\sim\sim(\sim p\;\vee\sim q)\vee r\supset (\sim p\;\vee\sim q)\vee r$\quad add.\zl ition\zd\ from \zl the\zd\ right\\[.5ex]
\>$(\sim p\;\vee\sim q)\vee r\supset\;\sim p\vee(\sim q\vee r)$\qquad associat.\zl ive\zd\ law\\[.5ex]
\underline{S\zl y\zd ll.\zl ogism}\zd $\;\;\;$\underline{$\sim\sim(\sim p\;\vee\sim q)\vee r\supset\;\sim p\vee(\sim q\vee r)$\quad q.e.d.}\\[1ex]
\>\underline{$[p\supset(q\supset r)]\supset(p\: .\: q\supset r)$}\qquad Importation\\[1ex]
\>$\sim p\vee(\sim q\vee r)\supset\;\sim\sim(\sim p\;\vee\sim q)\vee r$\\[1ex]
Pr.\zl oof\zd $\;\times$\quad$\sim p\vee(\sim q\vee r)\supset(\sim p\;\vee\sim q)\vee r$\qquad Associat\zl ivity\zd \\[.5ex]
\>$\sim p\;\vee\sim q\supset\;\sim\sim(\sim p\;\vee\sim q)$\\[.5ex]
$\times$\>$(\sim p\;\vee\sim q)\vee r\supset\;\sim\sim(\sim p\;\vee\sim q)\vee r$\qquad Add.\zl ition\zd\ right\\[.5ex]
\>$\sim p\vee(\sim q\vee r)\supset\;\sim\sim(\sim p\;\vee\sim q)\vee r$\qquad Syll\zl ogism\zd $\;\times\times$\\[1ex]
\>\underline{$[p\supset(q\supset r)]\supset[q\supset(p\supset r)]$}\\[1ex]
$\times$\>$\sim p\vee(\sim q\vee r)\supset(\sim p\;\vee\;\sim q)\vee r$\\[.5ex]
$\mathbf{\llbracket 76. \rrbracket}$\>$\sim p\;\vee\sim q\supset\;\sim q\;\vee\sim p$\\[.5ex]
$\times$\>$(\sim p\;\vee\;\sim q)\vee r\supset(\sim q\;\vee\;\sim p)\vee r$\\[.5ex]
$\times$\>$(\sim q\;\vee\;\sim p)\vee r\supset\;\sim q\vee(\sim p\vee r)$\\[.5ex]
\>\underline{$\sim p\vee(\sim q\vee r)\supset\;\sim q\vee(\sim p\vee r)$\qquad Syll\zl ogism\zd $\;\times\times\times$}
\end{tabbing}
\zl Here on p.\ \textbf{76}., which is on the right of p.\ \textbf{75}., one finds in a box on the left margin three lines that belong to that preceding page; they are inserted at appropriate places on p.\ \textbf{75}.\ in the text above.\zd\

Rule of exp.\zl ortation\zd\ or import.\zl ation\zd\ or commut.\zl ativity\zd
\begin{tabbing}
11\qquad\=\underline{$P\: .\: Q\supset R\quad :\quad P\supset(Q\supset R)$}\qquad Exp.\zl ortation\zd \\[.5ex]
11$\cdot$1\>\underline{$P\supset(Q\supset R)\quad :\quad P\: .\: Q\supset R$}\qquad Imp.\zl ortation\zd \\[.5ex]
11$\cdot$2\>\underline{$P\supset(Q\supset R)\quad :\quad Q\supset(P\supset R)$}\qquad Commut.\zl ativity\zd
\end{tabbing}

\zl After p.\ \textbf{76}.\ in this notebook comes p.\ \textbf{33}.\zd\

\vspace{1ex}

$\mathbf{\llbracket 33. \rrbracket}$ After having solved last time the first of the two probl\zl ems\zd\ I announced in the beg.\zl inning,\zd\ namely the probl.\zl em\zd\ of dec.\zl iding\zd\ of a given expr\zl ession\zd\ wheth.\zl er\zd\ or not it is a taut.\zl ology,\zd\ I come now to the sec.\zl ond,\zd\ namely to reduce the inf.\zl inite\zd\ nu.\zl mber\zd\ of taut\zl ologies\zd\ to a finite nu.\zl mber\zd\ of ax.\zl ioms\zd\ from which they can be derived. So this probl.\zl em\zd\ consists in setting up what is called a deductive syst.\zl em\zd\ for the calc.\zl ulus\zd\ of prop\zl ositions\zd . Now if you think of other ex\zl amples\zd\ of ded.\zl uctive\zd\ systems as e.g.\ geom.\zl etry\zd\ you will see that their aim is not truly to derive the theor.\zl ems\zd\ of the science concerned from a min.\zl imal\zd\ num.\zl ber\zd\ of ax.\zl ioms\zd , but also to define the notions \sout{\zl unreadable symbol\zd }\, of the disc.\zl ipline\zd\ con.\zl cerned\zd\ in terms of a min.\zl imal\zd\ nu.\zl mber\zd\ of undefined or $\mathbf{\llbracket 34. \rrbracket}$ primitive notions. So we shall do the same thing for \zl the\zd\ calc.\zl ulus\zd\ of prop\zl ositions\zd .

\zl new paragraph\zd\ We know already that some of the not.\zl ions\zd\ intro\-d\zl uced\zd\ $\sim$\zl ,\zd\ $\vee$\zl ,\zd\ $ .\:$\zl ,\zd\ $\supset$\zl ,\zd\ $\equiv$\zl ,\zd\ $\mid$ can be defined in terms of others\zl ,\zd\ namely e.g\zl .\zd\ $p\supset q\equiv\: \sim p\vee q$\zl ,\zd\ $p\equiv q\lfloor\equiv\rfloor p\supset q\: .\: q\supset p$\zl ,\zd\ but now we want to choose some of them in terms of which all others can be def\zl ined\zd . And I claim that e.g.\ $\sim$ and $\vee$ are suff\zl icient\zd\ for this purp.\zl ose\zd\ bec\zl ause\zd\
\begin{tabbing}
\hspace{1.7em}\=1.\quad \=$p\: .\: q\;\;$ \=$\equiv\;\sim(\sim p\; \vee \sim q)$\\[.5ex]
\>2.\> $p\supset q$ \>$\equiv\;\sim p\vee q$\\[.5ex]
\>3.\> $p\equiv q$ \>$\equiv (p\supset q)\: .\: (q\supset p)$\\[.5ex]
\>4.\> $p\mid q$ \>$\equiv\;\sim p\;\vee\sim q$
\end{tabbing}
So it is possible to take $\sim$ and $\vee$ as $\mathbf{\llbracket 35. \rrbracket}$ prim.\zl itive\zd\ terms for our ded.\zl uc\-tive\zd\ syst\zl em\zd\ and we shall actually \sout{do that} \ul make this choice \ud . But it is important to remark that this choice is fairly arb\zl itrary\zd . There would be other poss.\zl ibilities,\zd\ e.g.\ to take $\sim$\zl ,\zd\ $ .\:$ bec.\zl ause\zd\ $\vee$ can be expressed in terms of $\sim$ and $\: .\:$ by $p\vee q\equiv\;\sim(\sim p\: .\:\sim q)$ and by $\vee$ and $\sim$ the others can be expr\zl essed\zd\ as we have just seen. This fact that the choice of prim\zl itive\zd\ terms is arb\zl itrary\zd\ to a cert\zl ain\zd\ ext.\zl ent\zd\ is not surpr\zl ising\zd . The same situat\zl ion\zd\ prevails in any theory\zl ,\zd\ e.g\zl .\zd\ in geometrie\zl y\zd\ we can take \sout{the} either the notion of movement of the space or the notion of congr\zl uence\zd\ between \zl unreadable symbol, could stand for ``figures''\zd\ as prim\zl itive\zd\ because \zl it is\zd\ possible $\mathbf{\llbracket 36. \rrbracket}$ to define congr\zl uence\zd\ of \zl word missing, ``figures'' suggested above\zd\ in terms of movement of space and vice versa. The same situat.\zl ion\zd\ we have here. We can define $\vee$ in terms of ,,\zl ``\zd and'' and ,,\zl ``\zd not'' but also vice versa \zl\sout{or} ``and''\zd\ in terms of \zl ``\zd or\zl ''\zd\ and \zl ``\zd not\zl ''\zd . And there are still further poss.\zl ibilities\zd\ for the prim.\zl itive\zd\ terms\zl ,\zd\ e.g\zl .\zd\ it would be possible\zl p written over another letter\zd\ to take $\sim$ and $\supset$ as \zl the\zd\ only prim\zl itive\zd\ terms bec.\zl ause\zd\ $\vee$ can be defined \sout{in t} by
\begin{tabbing}
\hspace{1.7em}\=$p\vee q\equiv\;\sim p\supset q$ \quad since\\[.5ex]
\>$\sim p\supset q\equiv\;\sim\sim p\vee q\equiv p\vee q$\quad by \zl the\zd\ law of double neg\zl ation\zd
\end{tabbing}
In the three \sout{possible} cases disc.\zl ussed\zd\ so far we had always two prim\zl itive\zd\ notions in terms of $\mathbf{\llbracket 37. \rrbracket}$ which the others could be def\zl ined\zd . It is an interest\zl ing\zd\ quest.\zl ion\zd\ whether there might not be a single op.\zl eration\zd\ in terms of which all the others can be defined. This is actually the case as was first disc.\zl overed\zd\ by the log.\zl ician\zd\ Sheffer. Namely the $\mid$ f\zl u\zd nct.\zl ion\zd\ suffices to define all the others bec\zl ause\zd\ $\sim p\equiv p\mid p$ means at least one of the \sout{\zl unreadable text, perhaps: two\zd }\, prop\zl ositions\zd\ $p, p$ is false\zl ,\zd\ but since they are both $p$ that means $p$ is false\zl ,\zd\ i\zl .\zd e\zl .\zd\ $\sim p$\zl ,\zd\ \ul so $\sim$ can be def\zl ined\zd\ in terms of $\mid$ \ud and now the ,,\zl ``\zd and'' can be defined in terms of $\sim$ and $\mid$ \sout{by} since $p\: .\: q\equiv\;\sim(p\mid q)$ for $p\mid q$ means at least $\mathbf{\llbracket 38. \rrbracket}$ one of the two prop.\zl ositions\zd\ is false\zl ;\zd\ hence the neg\zl ation\zd\ means both are true. But in terms of $\sim$ and the \ul $\: .\:$ \ud others can be def.\zl ined\zd\ as we saw before. It is easy to see that we have now exhausted all possibilities of choosing \ul the \ud primit\zl ive\zd\ terms \sout{\zl unreadable symbol\zd }\, from the \ul 6\zl six\zd\ \ud operations written down here. In part\zl icular\zd\ we can prove \ul e.g.\ \ud: $\sim , \equiv$ are not suff\zl icient\zd\ to def.\zl ine\zd\ the others in terms of them. We can e.g\zl .\zd\ show that $p\vee q$ cannot be def\zl ined\zd\ in terms \ul of them \ud\zl .\zd\

\zl It is not indicated in the manuscript where exactly the following paragraph is to be inserted:\zd\ \ul \ul Now \ud \sout{W}what could it mean that \sout{\zl unreadable symbol\zd\ e.g} $p\: .\: q$ \zl or $p\vee q$\zd\ can be def\zl ined\zd\ in terms of $\sim , \equiv$\zl ?\zd\ It would mean that we can find an expr.\zl ession\zd\ $f(p,q)$ in two var.\zl iables\zd\ containing only the symb\zl ols\zd\ $\sim , \equiv$ besides $p,q$ and such that $p {\vee \atop .} q\equiv f(p\lfloor ,\rfloor q)$\zl ,\zd\ i\zl .\zd e\zl .\zd\ such that this expr.\zl ession\zd\ would have the same truth table as $p\lfloor {\vee \atop .} \rfloor q$\zl .\zd\ But we shall prove now that such an expression does not exist. \ud

$\Big\lceil$Let's write down the truth \zl functions\zd\ in two variables $p,q$ which we certainly can define in terms of $\sim , \equiv$\zl ;\zd\ we get the following eight\zl :\zd\ 1. ${p\equiv p}$\zl ,\zd\ 2. $\sim (p\equiv p)$\zl ,\zd\ 3. $p$\zl ,\zd\ 4. $q$\zl ,\zd\ 5. $\sim p$\zl ,\zd\ 6. $\sim q$\zl ,\zd\ $\mathbf{\llbracket new\, page \rrbracket}$ 7. $p\equiv q$\zl ,\zd\ 8. $\sim(p\equiv q)$\zl ,\zd\ and now it can be shown that no others can be def.\zl ined\zd\ exc.\zl ept\zd\ those eight because we can show the foll\zl owing\zd\ two things: 1. If we take one of those eight f\zl u\zd nct\zl ions\zd\ and negate it we get again one of those eight f\zl u\zd nct\zl ions,\zd 2. If we take any two of those eight f\zl u\zd nct\zl ions\zd\ and form a new one by connecting them by an equiv\zl alence\zd\ symbol we get again one of the eight. \zl I.\zd e\zl .\zd\ by appl.\zl ication\zd\ of the op.\zl eration\zd\ of neg\zl ation\zd\ and of the op\zl eration\zd\ of equiv.\zl alence\zd\ we never get outside of the set of eight f\zl u\zd nct\zl ions\zd\ written down \sout{here}\zl .\zd\ So let\zl '\zd s see at first that by negating them $\mathbf{\llbracket new\, page \rrbracket}$ we don't get anything new. Now if we neg.\zl ate\zd\ the first \zl text missing\zd . Now let\zl '\zd s connect any two of them by $\equiv$. If we connect \zl the first\zd\ with any form.\zl ula\zd\ $P$ we get $P$ again\zl ,\zd\ i\zl .\zd e.\ \sout{\zl unreadable symbols\zd }\, $(\top\equiv P)\equiv P$ \ul bec.\zl ause\zd\ \ud \zl text missing\zd\ and if connect a contrad.\zl iction\zd\ $C$ with any form\zl ula $P$\zd\ by \zl unreadable symbol, should be: $\equiv$\zd\ we get the neg\zl ation\zd\ of $P$\zl , i.e.\zd\ $(C\equiv P)\equiv\;\sim P$ bec\zl ause\zd \zl text missing\zd . So by comb.\zl ining\zd\ the first two form\zl ulas\zd\ with any other we get cert.\zl ainly\zd\ nothing new. For the other cases it is very helpful that $(p\equiv\;\sim q)\equiv\;\sim(p\equiv q)$\zl ;\zd\ this makes possible to factor out the neg.\zl ation\zd\ so to speak\zl .\zd\ Now in order to apply that to the other form.\zl ulas\zd\ we divide them in two groups\ldots\ \zl text missing\zd $\Big\rfloor$

$\mathbf{\llbracket 39. \rrbracket}$ For this purp.\zl ose\zd\ we divide the 16 \ul \zl truth functions\zd\ of two var.\zl iables\zd\ which we wrote down last time \ud into two classes according as the number of letters T occurring in their \zl\sout{their}\zd\ truth\zl $\,$\zd table is even or odd\zl ,\zd\ or to be more exact \ul accord\zl ing\zd\ as \ud the nu.\zl mber\zd\ \ul of T\zl 's\zd\ \ud occurring in the last col\zl umn\zd . So e.g.\ $p\: .\: q$ is odd\zl ,\zd\ $p\equiv q$ is even \ul and an arb\zl itrary\zd\ expr.\zl ession\zd\ in two var.\zl iables\zd\ will be called even if its truth\zl $\,$\zd f\zl unction\zd\ is even\zl .\zd\ And now what we can show is this: Any expr\zl ession\zd\ in two var.\zl iables\zd\ containing only $\sim$ \zl and\zd\ $\equiv$ is even (i\zl .\zd e.\ its truth\zl $\,$\zd table contains an even \zl\sout{)}\zd\ nu.\zl mber\zd\ of T's\zl ,\zd\ i\zl .\zd e\zl .\zd\ either 0 or 2 or 4 T's)\zl .\zd\

\sout{And} In order to show that we prove the following three lemmas.
\begin{itemize}
\item[1.] The let\zl ter\zd\ expr.\zl essions,\zd\ \sout{\zl unreadable text\zd\  form\zl ulas\zd }\, namely the letters $p,q$ are even\zl .\zd\

\item[2.] If an expr\zl ession\zd\ $f(p,q)$ is even then also the expr\zl ession\zd\ $\sim f(p,q)$ is even\zl .\zd\

\item[3.] If two expr.\zl essions\zd\ $f(p,q)$\zl ,\zd\ $g(p,q)$ are even then also the exp\zl res\-sion\zd\ $f(p,q)\equiv g(p,q)$ \ul obtained by connecting them with an equ\zl iva\-lence\zd\ sign \ud is even\zl .\zd\
\end{itemize}
$\mathbf{\llbracket 40. \rrbracket}$ So prop.\zl ositions\zd\ 2, 3 have the consequence\zl :\zd\
\begin{itemize}
\item[] By applying the op\zl erations\zd\ \ul $\sim$ \zl and\zd\ $\equiv$ to even expr\zl essions\zd\ \ud as many times as we wish we always get \ul again \ud \zl an\zd\ even expression \sout{if we start with \ul even \ud expr\zl ession\zd\ \zl unreadable word\zd\ the \zl unreadable word\zd }.
\end{itemize}
But any expr.\zl ession\zd\ cont.\zl aining\zd\ only $\sim$ \zl and\zd\ $\equiv$ is obtained from the single letters $p,q$ by an iterated appl\zl ication\zd\ of the op.\zl erations\zd\ $\sim$ \zl and\zd\ $\equiv$\zl ;\zd\ hence since $p,q$ are even the expr.\zl ession\zd\ thus obt.\zl ained\zd will also be even. So our theorem that every exp\zl ression\zd\ cont\zl aining\zd\ only $\sim$ \zl and\zd\ $\equiv$ is even will be proved \ul \zl when\zd\ we \ul shall \ud have proved the 3\zl three\zd\ lemmas\zl .\zd\ \zl \sout{(\zl unreadable symbol\zd )} \ud \zd\

\zl The following from the manuscript is deleted: \zl unreadable text\zd\ how to prove them \sout{\zl unreadable text 3\zd }\zd\ \underline{One}\zl 1.\zd\ is clear because \zl of the truth table\zd\ for $p$\ldots\ (and \zl unreadable text, perhaps: and for $q$\zd\ the same thing)\zl .\zd\ 2. also \zl is\zd\ clear because $\sim f(p,q)$ has T's when $f(p,q)$ had F's\zl ,\zd\ i.e\zl .\zd\ the nu.\zl mber\zd\ of T's in the new expr.\zl ession\zd\ is the same as the nu\zl mber\zd\ of F's in the \zl an insertion sign referring to nothing occurs in the manuscript on the right-hand side of this page\zd\ $\mathbf{\llbracket 41. \rrbracket}$ old one\zl .\zd\ But the nu\zl mber\zd\ of F's in the old one is even bec\zl ause the\zd\ number of \zl T's is\zd\ even and the nu\zl mber\zd\ of F's is \zl unreadable symbol, should be: equal to the\zd\ nu.\zl mber\zd\ of T's\zl .\zd\

\zl new paragraph\zd\ Now to the third\zl .\zd\ \zl unreadable text\zd\ \zl C\zd all the nu.\zl m\-ber\zd\ of T\zl 's\zd\ of the first $t_1$\zl ,\zd\ the nu.\zl mber\zd\ of T\zl 's\zd\ of the sec\zl ond\zd\ $t_2$ and call the nu.\zl mber\zd\ of cases \ul in the truth table \ud where both \zl $f$\zd\ and $g$ \sout{\zl unreadable text\zd }\, have \zl the\zd\ truth v\zl alue\zd\ T $r$\zl .\zd\ \zl We have\zd\ that \zl $t_1$ is\zd\ even and \zl $t_2$\zd\ is even\zl ,\zd\ but \zl unreadable text, should be: we do not\zd\ know anything about $r$\zl ;\zd\ it may be odd or even\zl .\zd\ \zl unreadable text, perhaps: We shall try\zd\ to find out in how many cases \zl unreadable text, should be: $f(p,q)\equiv g(p,q)$, i.e.\ $f\equiv g$,\zd\ will be true \ul and to show that this number of cases will be even\zl .\zd\ \ud I prefer to find out in how many cases it will be false. If we know that this nu.\zl mber\zd\ is even we know also that the nu.\zl mber\zd\ of cases in which it is true will be even\zl .\zd\ \ul Now this whole \zl expression\zd\ \ud is false \ul if \ud $g$ \zl and $f$\zd\ have diff\zl erent\zd truth v\zl alues,\zd\ i\zl .\zd e.\ if $\mathbf{\llbracket 42. \rrbracket}$ either \zl unreadable text, should be: we have $g$ false and $f$ true or we have\zd\ $g$ true \zl and $f$\zd\ false\zl .\zd\ \sout{But} The \zl unreadable text, should be: cases where $f$ is true and $g$ false make\zd\ $t_1-r$ cases bec.\zl ause\zd \zl unreadable text, should be: from $t_1$ cases where $f$ is true we should subtract cases\zd\ when $g$ is also true\zl , and\zd\ bec.\zl ause\zd\ $r$ was the nu\zl mber\zd\ of cases in which both are true\zl .\zd\ \zl H\zd\ ence in $t_1-r$ cases \zl unreadable text, should be: $f$ is T and\zd\ $g$ \zl is\zd\ F, sim\zl ilarly\zd\ in $t_2-r$ cases $g$ \zl is\zd\ T \zl and\zd\ $f$ \zl is\zd\ F\zl;\zd\ hence alt.\zl ogether\zd\ \ul in \ud $t_1-r+t_2-r$ \sout{)} \sout{$t_1+t_2-2r$} $\,$ cases $f$ \zl and\zd\ $g$ have diff\zl erent\zd\ truth values\zl ,\zd\ i\zl .\zd e.\ in \zl unreadable text, should be: $t_1+t_2-2r$\zd\ cases $f(p,q)\equiv g(p,q)$ is false\zl ,\zd\ and this is an even nu\zl mber\zd\ bec\zl ause\zd\ $t_1$\zl ,\zd\ $t_2$ \zl and\zd\ $2r$ are even \zl unreadable text\zd\ and if you add \zl unreadable text, should be: an even number to an even number, after subtracting an even number from the sum\zd you get again an even nu\zl mber\zd . Hence the number of cases in which the whole expr.\zl ession\zd\ \ul is false \ud is an even nu\zl mber\zd\ and \zl \sout{there are} such is\zd\ also the nu.\zl mber\zd\ of cases in which it is true\zl ,\zd\ i.e.\ $f(p,q)\equiv g(p,q)$ is an even expr\zl ession\zd . q.e.d.

So this shows that only even \sout{ex} truth\zl $\,$\zd f\zl u\zd nct\zl ions\zd\ $\mathbf{\llbracket 43. \rrbracket}$ \sout{truth table an even nu of T becau} can be expr\zl essed\zd\ in terms of $\sim$ \zl and\zd\ $\equiv$\zl .\zd\ Hence e.g\zl .\zd\ $\vee$ \ul and $.$ \ud cannot be expr.\zl essed\zd bec\zl ause\zd\ three T\zl 's\zd\ occur in \ldots\zl their truth tables.\zd\ It is easy to see that of the 16 truth\zl $\,$\zd f.\zl unctions\zd\ exactly half the nu.\zl mber\zd\ is even and also that all even truth\zl $\,$\zd f.\zl unctions\zd\ really can be expressed in terms of $\sim$ \zl and\zd\ $\equiv$ alone\zl .\zd\ Expr.\zl essions\zd\ for these eight \zl unreadable text, should be: truth functions\zd\ in terms of $\sim$ \zl and\zd\ $\equiv$ are given in the notes that were distributed \zl (see p.\ \textbf{38}.\ above)\zd . \sout{Our} The gen\zl eral\zd\ theor\zl em\zd\ \zl seems to be German ``\" uber'', translatable as ``on'') even f\zl u\zd nct\zl ions\zd\ I proved then has the consequ.\zl ence\zd\ that these eight \ul truth \ud f\zl u\zd nct\zl ions\zd\ must reproduce themselves by \sout{\zl unreadable text\zd }\, negating \sout{some of} them or by connecting any two of them by $\sim$\zl ;\zd\ i\zl .\zd e\zl .\ i\zd f you neg\zl ate\zd\ one of those \ul \zl unreadable word\zd\ \ud expr\zl essions\zd\ the result\zl ing\zd\ expr\zl ession\zd\ will be equiv.\zl alent\zd\ to one of the eight and if you form a new expr\zl ession\zd\ by connect.\zl ing\zd\ any two of them the resulting expression will again be equivalent to one of the eight\zl .\zd\ I recom\zl mend\zd\ \zl \sout{that}\zd\ $\mathbf{\llbracket 44. \rrbracket}$ as an exercise to show that in detail. \ul
It is an easy corol.\zl lary\zd\ of \ud this result about the undefinability of $.$ \zl and\zd\ $\vee$ in terms of $\equiv$ that also $\sim$ and the excl.\zl usive\zd\ \underline{or} are not suff.\zl icient\zd\ as primit\zl ive\zd\ terms because as we saw last time the excl.\zl usive\zd\ or can be expr.\zl essed\zd\ in terms of $\sim$ \zl and\zd\ $\equiv$\zl ,\zd\ namely by $\sim(p\equiv q)$\zl ;\zd\ hence if \ul e.g.\ \ud $\vee$ could be def.\zl ined\zd\ in terms of $\sim$ \zl and\zd\ $\circ$ \zl (exclusive or)\zd\ it could also be def\zl ined\zd in terms of $\sim$ \zl and\zd\ $\equiv$ bec.\zl ause\zd the $\circ$ can be expr\zl essed\zd\ in terms of $\sim$ \zl and\zd\ $\equiv$. The reason for that is of course that $\circ$ is also an even f.\zl u\zd nct\zl ion\zd\ and therefor only even f.\zl u\zd nct\zl ions\zd\ can be def.\zl ined\zd\ in terms of it\zl .\zd\ So we see that whereas $\sim$ \zl and\zd\ $\vee$ are suff\zl icient\zd as $\mathbf{\llbracket 45. \rrbracket}$ prim.\zl itive\zd\ terms \underline{$\sim$ \ul and \ud excl.\zl usive\zd\ or} are not\zl ,\zd\ which is one of the reasons why the not excl.\zl usive\zd\ or is used in log\zl ic\zd . Another of those neg.\zl ative\zd\ results about the poss.\zl ibility\zd\ of expressing some of the truth f.\zl unctions\zd\ by others would be \ul that \ud $\sim$ cannot be def\zl ined\zd\ in terms of $.\:, \vee, \supset$\zl ;\zd\ even in terms of all three of them it is impossible to expr.\zl ess\zd\ $\sim$\zl .\zd\ \zl I will\zd\ give that as a problem to prove.

$\Big\lceil$As I announced before we shall choose from the diff.\zl erent\zd\ possib.\zl ilities\zd\ of primitive terms for our ded.\zl uctive\zd\ syst.\zl em\zd\ the \sout{case where} \ul one in which \ud $\sim$ and $\vee$ \zl \sout{is} are\zd\ taken as prim\zl itive\zd\ and therefore it is of imp.\zl or\-tance\zd\ to make sure that not only the part.\zl icular\zd\ f\zl u\zd nct.\zl ions\zd\ $\equiv$\zl ,\zd\ $.\:$\zl ,\zd\ $\supset$\zl ,\zd\ $\mid$ for which $\mathbf{\llbracket 46. \rrbracket}$ we introduced special symbols but that any truth\zl $\,$\zd f\zl unction\zd\ whatsoever in any number of var\zl iables\zd\ can be expressed by $\sim$ \zl and\zd\ $\vee$. For truth\zl $\,$\zd f\zl unctions\zd\ with 2\zl two\zd\ variables that follows from the consid.\zl erations\zd\ of last time since we have expr.\zl essed\zd\ all 16 truth\zl $\,$\zd f\zl unctions\zd\ by our logistic symbols and today we have seen that all of them can be expr.\zl essed\zd\ by $\sim$ \zl and\zd\ $\vee$. Now I shall prove the same thing \ul also \ud for truth\zl $\,$\zd f\zl u\zd nct\zl ions\zd\ with 3\zl three\zd\ variables and you will see that the method of proof can be applied to f\zl u\zd nct\zl ions\zd\ of any number of variables. For the three var.\zl iables\zd\ $p.q,r$ we have eight $\mathbf{\llbracket 47. \rrbracket}$ possibilities for the distr\zl ibution\zd\ of truth values over them\zl ,\zd\ namely
\begin{center}
\begin{tabular}{ l c c c|l l}
 & $p$ & $q$ & $r$ & $f(p\lfloor ,\rfloor q\lfloor ,\rfloor r)$\\[.5ex]
 \hline
1\zl .\zd & T & T & T & \hspace{1em} $p\: .\: q\: .\: r$ & \hspace{1em}$P_1$\\
2\zl .\zd & T & T & F & \hspace{1em} $p\: .\: q\: . \sim r$ & \hspace{1em}$P_2$\\
3. & T & F & T & \hspace{1em} $p\: . \sim q\: .\: r$\\
\zl 4.\zd & T & F & F\\
\zl 5.\zd & F & T & T\\
6\zl .\zd & F & T & F\\
7\zl .\zd & F & F & T\\
8\zl .\zd & F & F & F & & \hspace{1em}$P_8$
\end{tabular}
\end{center}

Now to define a truth \ul fu\zl nction\zd\ \ud in three var\zl iables\zd\ means \zl comma from the manuscript deleted\zd\ to stipulate a truth value \ul T or F \ud for $f(p,q\lfloor , r\rfloor)$ for each of these eight cases. Now to each of these 8\zl eight\zd\ cases we can associate a cert.\zl ain\zd\ expr\zl ession\zd\ in the foll\zl owing\zd\ way\zl :\zd\ to 1. \ul we associate \zl colon from the manuscript deleted\zd\ \ud $p\: .\: q\: .\: r$\zl , to\zd\ 2. \zl we associate\zd\ $p\: .\: q\: . \sim r$\zl , to\zd\ 3. \zl we associate\zd\ $p\: . \sim q\: .\: r$\zl ,\zd\ \ldots\ So each of these expr\zl essions\zd\ will have a $\sim$ before those letters which have an F in the corresp.\zl onding\zd\ case. Denote the expr\zl essions\zd\ associated with these eight lines by $P_1$\zl ,\zd \ldots\zl ,\zd $P_8$. Then the expr\zl ession\zd\ $P_2$ e.g.\ will be true then and only $\mathbf{\llbracket 48. \rrbracket}$ then if the \underline{sec\zl ond\zd } case is realis\zl z\zd ed for the truth values of $p,q,r$ ($p\: .\: q\: . \sim r$ will be true then and only then \zl if\zd $p$ \zl is\zd T\zl ,\zd\ $q$ \zl is\zd T\zl and\zd\ $r$ \zl is\zd\ false\zl ,\zd\ which is exactly the case for the truth val\zl ues\zd\ $p$, $q$\zl ,\zd\ $r$ represented in the 3\zl second\zd\ line\zl .\zd\ And general\zl l\zd y $P_i$ will be true \sout{then and only then} if the $i^{\,\rm th}$ case for the truth values of $p,q,r$ is realis\zl z\zd ed\zl .\zd\ Now the \zl unreadable text\zd\ truth\zl $\,$\zd f.\zl unction\zd \zl unreadable symbol\zd which we want to expr.\zl ess\zd\ by $\sim$ \zl and\zd\ $\vee$ will be true for cert.\zl ain\zd\ of those 8\zl eight\zd\ cases and false for the others. Assume it is true for case numbe\zl r\zd\ $i_1$\zl ,\zd $i_2$\zl ,\zd \ldots\zl ,\zd $i_n$ and false for the others. Then form the disj.\zl unction\zd\ $P_{i_1}\vee P_{i_2}\ldots\vee P_{i_n}$\zl ,\zd\ i.e\zl .\zd\ the disj\zl unction of those $P_i$ which correspond to the cases in which the given f\zl u\zd nct\zl ion\zd\ is true. This \ul disj\zl unction\zd\ \ud is an expr.\zl ession\zd\ in the \zl \sout{the}\zd\ var\zl iables\zd\ $p,q,r$ containing only the op.\zl erations\zd\ $.$\zl ,\zd\ $\sim$ \zl and\zd\ $\vee$\zl ,\zd\ and I claim its truth table $\mathbf{\llbracket 49. \rrbracket}$ will coincide with the truth table of the given expr.\zl ession\zd\ $f(p, q\lfloor , \rfloor r)$\zl .\zd\ For \zl colon from the manuscript deleted\zd\ $f(p, q\lfloor , \rfloor r)$ had the symb.\zl ol\zd\ T in the $i_1$\zl ,\zd $i_2$\zl ,\zd \ldots\zl ,\zd $i_n^{\;\rm th}$ line but in no others and I claim the same thing is true for the expr\zl ession\zd\ $P_{i_1}\vee\ldots\vee P_{i_n}$\zl .\zd\ \sout{Now}

\zl new paragraph\zd\ You see \ul at \zl last?\zd\ \ud a disj\zl unction\zd\ of an arb\zl itrary\zd\ nu.\zl mber\zd\ of members will be true then and only then if at least one of its members is true and it will be false only if all of its members are false (I proved that in my last lecture for the case of 3\zl three\zd\ members and the same proof holds generally). Hence this disj\zl unction\zd\ will cert\zl ainly\zd\ be true in \zl the\zd\ $i_1$\zl ,\zd\ \ldots\zl ,\zd\ $i_n^{\;\rm th}$ case because $P_{i_1}$ \ul e.g.\ \ud is true in the $i_1^{\;\rm th}$ case as we \zl unreadable text, perhaps: saw before\zd . Therefore the $\mathbf{\llbracket 50. \rrbracket}$ disj\zl unction\zd\ is \ul also \ud true for the $i_1^{\;\rm th}$ case \ul because then one of its memb.\zl ers\zd\ is true\zl .\zd\ The same holds for $i_2$ \ldots\ etc. So the truth table for the disj\zl unction\zd\ will cert.\zl ainly\zd\ have the \ul letter T \ud in the $i_1$\zl ,\zd \ldots\zl ,\zd $i_n$ line. But it will have F's in all the other lines. Bec.\zl ause\zd\ $P_{i_1}$ was true only in the $i_1^{\;\rm th}$ case and false in all the others. Hence in a case diff\zl erent\zd\ from \ul the \ud $i_1$\zl ,\zd \ldots\zl ,\zd $i_n^{\;\rm th}$ $P_{i_1}$\zl ,\zd \ldots\zl ,\zd $P_{i_n}$ will all be false and hence the disj\zl unction\zd\ will be false, i.e.\ $P_{i_1}\vee\ldots\vee P_{i_n}$ will have the letter F in all lines other than the $i_1$\zl ,\zd \ldots\zl ,\zd $i_n^{\;\rm th}$\zl ,\zd\ i\zl ,\zd e.\ it has T in the $i_1$\zl ,\zd \ldots\zl ,\zd $i_n$ line and only in those. But the same thing was true for the truth t\zl able\zd\ of the given $f(p, q , \lfloor r \rfloor)$ \ul by ass\zl umption.\zd\ \ud So they coincide\zl ,\zd\ i\zl .\zd e.\ $f(p\lfloor ,\rfloor q\lfloor , \rfloor r)\equiv P_{i_1}\vee \ldots\vee P_{i_n}$\zl .\zd\

$\mathbf{\llbracket 51. \rrbracket}$ So we have proved that an arb.\zl itrary\zd\ truth funct\zl ion\zd\ of 3\zl three\zd\ variables can be expr.\zl essed by $\sim$\zl ,\zd\ $\vee$ \zl and\zd\ $.$\zl ,\zd\ but $.$ can be expr.\zl essed\zd\ by $\sim$ and $\vee$\zl ,\zd\ hence every truth\zl $\,$\zd f\zl unction\zd\ of three var.\zl iables\zd\ can be expr.\zl essed\zd\ by $\sim$ and $\vee$\zl ,\zd\ and I think it is perfectly clear that exactly the same proof applies to truth\zl $\,$\zd f.\zl unctions\zd\ of any number of variables. \zl wavy vertical line dividing the page\zd\

Now after having \zl unreadable text, should be: seen that\zd\ \uline{two prim\zl itive\zd\ noti}ons \ul \zl $\sim, \vee$\zd\ \ud really suffice to define any truth\zl $\,$\zd f\zl unction\zd\ we can begin to set up the ded.\zl uctive\zd\ syst\zl em.\zd\

I begin with \zl writing\zd\ three \uline{def.\zl initions\zd\ in terms of our prim\zl itive\zd } no\-tions\zl :\zd
\begin{tabbing}
\hspace{1.7em}\= $P\supset Q$ \=$=_{Df}\;\:$\=$\sim P \vee Q$\\[.5ex]
\> $P\: .\: Q$ \>$=_{Df}$\>$\sim(\sim P\:\vee \sim Q)$\\[.5ex]
\> $P\equiv Q$ \>$=_{Df}$\>$P\supset Q\: .\: Q\supset P$
\end{tabbing}
$\mathbf{\llbracket 52. \rrbracket}$ I am writing cap.\zl ital\zd\ letters because these def\zl initions\zd\ are to apply also if $P$\zl and\zd\ $Q$ are form\zl ulas\zd , not only if they are single letters\zl ,\zd\ i.e.\ e.g\zl .\zd\ $p\supset p\vee q$ means $\sim p\vee(p\vee q)$ and so on\zl .\zd\

The next thing to do in order to have a ded\zl uctive\zd\ syst\zl em\zd\ is to set up the ax\zl ioms\zd . Again in the axioms one has a freed.\zl om\zd\ of choice as in the primit.\zl ive\zd\ terms\zl ,\zd\ exactly as \ul also \ud in other ded\zl uctive\zd\ theories \sout{also}\zl ,\zd\ e.g.\ in geometr\zl y,\zd\ many diff\zl erent\zd\ syst\zl ems\zd\ of ax.\zl ioms\zd \zl \sout{for}\zd \sout{geo} have been set up each of which is suff\zl icient\zd\ to derive the whole geom\zl etry\zd . The syst.\zl em\zd\ of ax.\zl ioms\zd\ \ul for the calc.\zl ulus\zd\ of prop\zl ositions\zd\ \ud which I use is ess.\zl entially\zd\ the one set up by first by Russell and then also adopted by Hilbert. It has the foll.\zl owing\zd\ four ax\zl ioms:\zd\

\vspace{1ex}

\noindent $\mathbf{\llbracket 53. \rrbracket}$
\begin{tabbing}
\hspace{1.7em}\= \underline{1.}\zl (1)\zd\quad \=$p\supset p\vee q$ \quad\zl in the manuscript 1., 2. and 3. are in one line\zd\ \\[.5ex]
\> \underline{2.}\zl (2)\zd\>$p\vee p\supset p$\\[.5ex]
\> \underline{3.}\zl (3)\zd\>$p\vee q\supset q\vee p$\\[.5ex]
\> \underline{4.}\zl (4)\zd\>$(p\supset q)\supset(r\vee p\supset r\vee q)$
\end{tabbing}

I shall discuss the meaning of these ax\zl ioms\zd\ later \sout{discuss later}. A\zl t\zd\ present I want only to say that an expr\zl ession\zd\ written down in our theory as an axiom or as a theorem always means that it is true for any prop\zl ositions\zd\ $p,q,r$ etc\zl .,\zd\ e.g.\ $p\supset (p\vee q)$\zl $p\supset p\vee q$.\zd\

Now in geom.\zl etry\zd\ \sout{\zl unreadable word\zd }\, and any other \sout{disc} \ul theor\zl y\zd\ \ud $\mid$\uline{ex\-c\zl ept\zd\ logic}$\mid$ the ded.\zl uctive\zd\ syst.\zl em\zd\ is completely given by stating what the prim\zl itive\zd\ terms and what the ax.\zl ioms\zd\ are. It is important to remark that it is different here for the following reason: in geom\zl etry\zd\ \ul and other theor\zl ies\zd\ \ud it is clear how the theorems are to be derived from the ax.\zl ioms;\zd\ they are to be derived by the rules of logic which are assumed to be known. In our case however we cannot assume the rules of logic to be known \sout{by the rules of log.\zl ic\zd }\, $\mathbf{\llbracket 54. \rrbracket}$ because we are just about to formulate the rules of logic and to reduce them to a min\zl imum.\zd\ So this will naturally have to apply to the rules of inference as well as to the ax\zl ioms\zd\ with which we start. We shall have to formulate the\zl written over ``them''\zd\ \ul rules of inf.\zl erence\zd\ \ud explicitly and with greatest possible precision \sout{and}\zl ,\zd\ that is in such a way there can never be a doubt whether a cert\zl ain\zd\ rule can be applied for any form\zl ula\zd\ or not. And of course we shall try to\zl comma from the manuscript deleted\zd \zl \ldots\ text omitted in the manuscript, could be: work\zd\ with as few as possible. I have to warn here against a\zl n\zd\ error\zl .\zd\ \sout{one might}

$\Big\lceil$One might think that an expl\zl icit\zd\ formulation of the rules of inf.\zl erence\zd\ \ul besides the ax.\zl ioms\zd\ \ud is superfl\underline{u}ous bec.\zl ause\zd\ the ax.\zl ioms\zd\ themselves \ul seem to \ud express rules of inf.\zl erence,\zd\ e.g.\ $p\supset p\vee q$ \ul the rule \ud that from \sout{$p$} a prop\zl osition\zd\ $p$ one can conclude $p\vee q$\zl ,\zd\ and \zl unreadable symbol\zd\ one might think that the ax.\zl ioms\zd\ themselves contain at the same time the rules by which the theorems are to be derived. But this way out of the diff\zl iculty\zd\ would be entirely wrong $\mathbf{\llbracket 55. \rrbracket}$ bec.\zl ause\zd\ e.g.\ $p\supset p\vee q$ does not say that it is perm\zl itted\zd to conclude $p\vee q$ from $p$ because those terms \zl ``\zd allowable to conclude\zl ''\zd\ do\zl $\,$\zd not occur in it. The notions \zl unreadable text, should be: in it are\zd\ only $p$\zl ,\zd\ $\supset$, $\vee$ and \zl $q$.\zd\ \zl A\zd cc.\zl ording\zd\ to our def\zl inition\zd\ of $\supset$ it \ul does not mean that\zl ,\zd\ but it \ud simply says $p$ is false or $p\vee q$ is true. It is true that the axioms suggest \ul or make possible \ud cert.\zl ain\zd\ rule\zl s\zd\ of inf.\zl erence,\zd e.g\zl .\zd\ the just \zl stated one,\zd\ but it is not even uniquely det.\zl ermined\zd\ what rules of inf\zl erence\zd\ it suggests\zl ;\zd\ e.g\zl .\zd\ $\sim p\vee(p\vee q)$ says either $p$ is false or $p\vee q$ is true\zl ,\zd\ which sugg\zl ests\zd\ the rule of inf.\zl erence\zd\ $p$\zl :\zd\ $p\vee q$\zl ,\zd\ but it also sug.\zl gests\zd\ $\sim(p\vee q)$\zl :\zd\ $\sim p$\zl .\zd\ So \zl \sout{its} we\zd\ need written spec.\zl ifications,\zd\ i\zl .\zd e.\ we have to formulate rules of inf\zl erence\zd\ in add.\zl ition\zd\ to formulas\zl .\zd\ \zl Note in a box: p 56 - p 60 \} Heft\zl German: Notebook\zd\ I\zd\

It is only bec\zl ause\zd\ the \zl ``\zd if then\zl ''\zd\ in ord.\zl inary\zd\ langu\zl age\zd\ is amb.\zl i\-valent\zd\ and has besides the mean\zl ing\zd\ given by the truth\zl $\,$\zd t\zl able\zd\ also the mean\zl ing\zd\ \zl unreadable symbol, should be ``the second member\zd\ can be inferred from \zl unread\-able symbol, should be: the first''\zd\ \zl comma from the manuscript deleted\zd\ that the ax.\zl ioms\zd\ seem to express \ul uniquely \zl unread\-able text\zd\ \ud rules of inf\zl erence\zd .

$\mathbf{\llbracket 55.1 \rrbracket}$ This remark applies gen.\zl erally\zd\ to any quest\zl ion\zd\ whether \ul or not \ud cert\zl ain\zd\ laws of log.\zl ic\zd\ can be derived from others (e.g.\ whether \zl the\zd\ law of excl.\zl uded\zd\ middle are\zl is\zd\ sufficient)\zl .\zd\ Such quest.\zl ions\zd\ have only a precise mean\zl ing\zd\ if you state the rules of inf\zl erence\zd\ which are to be accept\zl ed\zd\ in the deriv\zl ation\zd . \zl The remaining text on p.\ \textbf{55.1} is in a box\zd\ It is diff.\zl rent\zd\ e.g\zl .\zd\ in geom\zl etry;\zd\ there it has a precise mean\zl ing\zd\ whether it follows\zl ,\zd\ namely it means whether it foll\zl ows\zd\ by log\zl ical\zd\ inf\zl erence,\zd\ but it cannot have this mean\zl ing\zd\ in log\zl ic\zd\ because then every log\zl ical\zd\ law would be der\zl ivable\zd\ from any other. So it could $\mathbf{\llbracket 55.2 \rrbracket}$ only mean derivable by the inf.\zl erences\zd\ made possible by the ax\zl ioms\zd . But as we have seen that has no precise mean\zl ing\zd\ bec.\zl ause\zd\ an ax.\zl iom\zd\ may make possible or sugg.\zl est\zd\ many inferences.

\vspace{2ex}

\zl On a not numbered page after p.\ \textbf{55}.\textbf{2}, which is the last page of the present notebook, one finds the following crossed out text:\zd\

\vspace{1ex}

\noindent which describe unambiguously how the mean.\zl ingful\zd\ expr.\zl essions\zd\ are to be formed from the basic symb.\zl ols\zd\ (rules of the grammar of the langu.\zl age\zd )

\section{Notebook III}\label{0III}
\pagestyle{myheadings}\markboth{SOURCE TEXT}{NOTEBOOK III}
\zl Folder 61, on the front cover of the notebook ``Log.\zl ik\zd\ Vorl.\zl esungen\zd\ \zl German: Logic Lectures\zd\ N.D\zl .\zd\ \zl Notre Dame\zd\ III''\zd\

\begin{tabbing}
$\mathbf{\llbracket 1. \rrbracket}$\hspace{1em}\= $(p\supset q)\supset [(r\supset p)\supset(r\supset q)]$\\[.5ex]
\>\underline{$(q\supset r)\supset[(p\supset q)\supset(p\supset r)\lfloor ]\rfloor$}\\
\>\underline{$(p\supset q)\supset [(q\supset r)\supset (p\supset r)]$}\hspace{.2em} \= Commut\zl ativity\zd\ \=\;\;$\f{p \supset q}{P}$ $\f{q \supset r}{Q}$ $\f{p \supset r}{R}$\\
\>\>\>\zl $\f{q \supset r}{P}$ $\f{p \supset q}{Q}$ $\f{p \supset r}{R}$\zd\ \\[.5ex]
\>\underline{$(p\supset q)\: .\: (q\supset r)\supset (p\supset r)$} \> Import.\zl ation\zd\ \> \hspace{.8em} \;$''$ \hspace{1.6em} \;$''$ \hspace{1.5em} \;$''$\\
\>\>\>\zl $\f{p \supset q}{P}$ $\f{q \supset r}{Q}$ $\f{p \supset r}{R}$\zd\ \\

\>\underline{$(q\supset r)\: .\: (p\supset q)\supset (p\supset r)$}\\[2ex]

\zl in a box: p.\ 42, 45 Examples p 53\zd\ \\[2ex]

\> \underline{$(p\supset q)\: .\: p \supset q$}\\
\> $(p\supset q)\supset (p\supset q)$\qquad $\f{p \supset q}{P}$ $\f{p}{Q}$ $\f{q}{R}$ \\[.5ex]
\>$(p\supset q)\: .\: p \supset q$\hspace{3em} Import.\zl ation\zd
\end{tabbing}

\begin{tabbing}
$\mathbf{\llbracket 2. \rrbracket}$\\*[1ex]
17 \hspace{3em}\=\underline{$p\: .\: q\supset q\: .\: p$}\\
Pr.\zl oof\zd\ \>$\sim q \:\vee \sim p \supset \:\sim p \:\vee \sim q$\qquad (3) $\f{\sim q}{p}$ $\f{\sim p}{q}$ \zl fraction bars omitted in\\[-2ex]
\` the manuscript\zd\ \\[.5ex]
\>$\sim(\sim p\;\vee\sim q)\supset \:\sim (\sim q\;\vee\sim p)$ \qquad Transp.\zl osition\zd \\[.5ex]
\>$p\: .\: q\supset q\: .\: p$ \qquad rule of def\zl ined\zd\ symb\zl ol\zd
\end{tabbing}

\begin{tabbing}
18. \hspace{3em}\=\underline{$p\supset p\: .\: p$}\\[1ex]
Pr\zl oof\zd\ \>$\sim p\; \vee \sim p \supset \; \sim p$\\[.5ex]
\>$p \supset \; \sim (\sim p \;\vee \sim p)$ \qquad Transp.\zl osition\zd \\[.5ex]
\>$p\supset p\: .\: p$ \qquad def\zl ined\zd\ symb\zl ol\zd
\end{tabbing}

\begin{tabbing}
19. \hspace{3em}\=\underline{$p\supset (q\supset p\lfloor .\rfloor q)$}\\[1ex]
\>$(p\lfloor .\rfloor q \supset p\lfloor .\rfloor q)\supset (p \supset (q \supset p\lfloor .\rfloor q))$\quad\= export\zl ation\zd\ \quad\=\\
\>\>$\f{p\lfloor .\rfloor q}{r}$\\[-3ex]
\>$p \supset (q \supset p\lfloor .\rfloor q)$\\[2ex]

19.1 \>\underline{$p\supset (q\supset q\lfloor .\rfloor p)$}\\[-1ex]
\>$(p\lfloor .\rfloor q \supset q\lfloor .\rfloor p)\supset (p \supset (q \supset q\lfloor .\rfloor p))$\> export.\zl ation\zd\ \> $\f{q\lfloor .\rfloor p}{r}$
\end{tabbing}

\begin{tabbing}
$\mathbf{\llbracket 3. \rrbracket}$\\*[1ex]
12\zl over 11\zd $_R$\zl R\zd\ \hspace{2em}\=\underline{$P\, ,\,Q\quad .\:.\lfloor :\rfloor \quad P\: .\: Q$}\qquad rule of prod\zl uct\zd \\[.5ex]
\>$P \supset (Q \supset P\: .\: Q)$ \\[.5ex]
\>$Q \supset P\: .\: Q$ \\[.5ex]
\>$P\: .\: Q$
\end{tabbing}

\begin{tabbing}
12\zl over 11\zd $_R$\zl R\zd\ \hspace{2em}\=\underline{$P,Q\quad .\:.\lfloor :\rfloor \quad P\: .\: Q$}\qquad rule of prod\zl uct\zd \kill

Inv.\zl ersion\zd\ \>\underline{$P\: .\: Q\quad .\:.\lfloor :\rfloor \quad P\, ,\,Q$}\qquad rule of prod.\zl uct\zd \\[.5ex]
\>$P\: .\: Q \supset P \quad P\: .\: Q \supset Q$
\end{tabbing}

\noindent \zl The following three lines, up to 13$_R$\zl R\zd , are crossed out in the manuscript:\zd\
\begin{tabbing}
21\zl .\zd\ \hspace{3em}\=\underline{$\sim (p\: \lfloor .\rfloor \: q)\equiv \: \sim p \: \vee \sim q$}\\[.5ex]
\>$\sim \: \sim(\sim p \: \vee \sim q) \equiv \: \sim p \: \vee \sim q$ \qquad $\f{\sim p \: \vee \sim q}{p}$ \\[.5ex]
\>$\lceil \sim (p \vee q) \equiv \: \sim p\: . \sim q \rceil$
\end{tabbing}

\begin{tabbing}
13$_R$\zl R\zd\ \hspace{1.7em}\=\underline{$P\supset Q\quad R\supset S\quad .\:.\lfloor :\rfloor \quad P\: .\: R\supset Q\: .\: S$} \` \underline{Rule of multiplic.\zl ation\zd }\\[.5ex]
\>$\sim Q\supset \: \sim P \quad \sim S\supset \: \sim R$ \\[.5ex]
\>$\sim Q \: \vee \: \sim S \supset \: \sim P \: \vee \: \sim R$ \\[.5ex]
\>$\sim(\sim P \: \vee \: \sim R) \supset \: \sim (\sim Q \: \vee \: \sim S)$
\end{tabbing}

\begin{tabbing}
$\mathbf{\llbracket 4. \rrbracket}$\\*[1ex]
\underline{13.1$_R$\zl R\zd } \hspace{1.7em}\=\underline{$P\supset Q \quad \ldots \quad R\: .\: P\supset R\: .\: Q$} \\[.5ex]
\> bec\zl ause\zd\ $R\supset R$ and other side
\end{tabbing}

\begin{tabbing}
\underline{13.2$_R$\zl R\zd} \hspace{1.7em}\=\underline{$P\supset Q$} \zl ,\zd\ \underline{$P\supset S$} \quad : \quad\underline{$P \supset Q\: .\: S$} \\[.5ex]
\> $P\: .\: P\supset Q\: .\: S$ \\[.5ex]
\> \underline{$P\supset P\: .\: P$} \\[.5ex]
\> $P\supset Q\: .\: S$ \qquad rule of composition
\end{tabbing}

\noindent \zl An insertion sign from the manuscript followed by ``p 5-6'' is deleted.\zd\

\begin{tabbing}
F 22. \hspace{3em}\=\underline{$p\: .\:(q\vee r) \equiv p\: .\: q \vee p\: .\: r$}\\[.5ex]
\quad\quad I. \>$q \supset q \vee r$ \\[.5ex]
\>$p\: .\:q \supset p\: .\:(q \vee r)$ \\[.5ex]
\>$r \supset q \vee r$\\[.5ex]
\>$p\: .\:r \supset p\: .\:(q \vee r)$ \\[.5ex]
\>$p\: .\: q \vee p\: .\: r \supset p\: .\:(q\vee r)$
\end{tabbing}

\vspace{-2ex}

\noindent\rule{12cm}{0.4pt}\\
II. \zl The following two columns of formulae are separated by a vertical line in the manuscript:\zd\

\begin{tabbing}
$\times$\hspace{.6em}\=$q \supset (p \supset p\: .\: q)$ \hspace{8.1em}\= $q\supset (p \supset p\: .\: q \vee p\: .\: r)$ \\[.5ex]
+\>$r \supset (p \supset p\: .\: r)$ \hspace{6.5em} + \>$(p\supset p\: .\: r)\supset (p \supset p\: .\: q \vee p\: .\: r)$\\[.5ex]
\>$p\: .\: q \supset p\: .\: q \vee p\: .\: r$ \> $r\supset (p \supset p\: .\: q \vee p\: .\: r)$ \\[.5ex]
\>$p\: .\: r \supset p\: .\: q \vee p\: .\: r$ \> $q\vee r \supset (p \supset p\: .\: q \vee p\: .\: r)$\\[.5ex]
$\times$\>$(p \supset p\: .\: q)\supset (p \supset p\: .\: q \vee p\: .\: r)$ \> $(q \vee r)\: .\: p \supset p\: .\: q \vee p\: .\: r$
\end{tabbing}

\noindent $\mathbf{\llbracket 5. \rrbracket}$ Ae\zl E\zd quivalences

\begin{tabbing}
~~~~\hspace{3em}\=\underline{$P\supset Q\: .\: Q\supset P \quad .\:.\lfloor :\rfloor \quad P\equiv Q$} \\[.5ex]
\> bec.\zl ause\zd\ $(P\supset Q)\: \lfloor .\rfloor \: (Q\supset P)$\quad rule of def\zl ined\zd\ symb\zl ol\zd\
\end{tabbing}

\begin{tabbing}
~~~~\hspace{3em}\=\underline{$P\equiv Q$} $\quad .\:.\lfloor :\rfloor \quad$ \underline{$P\supset Q\: \lfloor .\rfloor \: Q\supset P$}
\end{tabbing}

Transpos\zl ition\zd :

\begin{tabbing}
~~~~\hspace{3em}\=\underline{$P \equiv Q \qquad .\:.\lfloor :\rfloor \qquad \sim P \equiv \: \sim Q$} \\[.5ex]
\> \underline{$P \equiv \: \sim Q \qquad .\:.\lfloor :\rfloor \qquad \sim P \equiv Q$} \\[.5ex]
Proof \> $P \equiv Q$ \hspace{1em} \=$P \supset Q$ \hspace{2.3em} \=$Q \supset P$ \\[.5ex]
\>\> $\sim Q \supset \: \sim P$ \> $\sim P \supset \: \sim Q$ \hspace{1em} $\sim P \equiv \: \sim Q$
\end{tabbing}

Add.\zl ition\zd\ and Multipl.\zl ication\zd\
\begin{tabbing}
~~~~\hspace{3em}\= \underline{$P\equiv Q\qquad R\equiv S$}\hspace{2em}
$
\begin{cases}
\;\;\;\underline{P\vee R \equiv Q\vee S}\\
\;\;\;\underline{P\: .\: R \equiv Q\: .\: S}
\end{cases}
$
\end{tabbing}

\begin{tabbing}
$\mathbf{\llbracket 6. \rrbracket}$\hspace{3em}\= $P\supset Q$ \hspace{0.5em} $R\supset S$ \hspace{1em}\= $Q\supset P$ \hspace{0.5em} $S\supset R$ \\[0.5ex]
\> $P\vee R \supset Q\vee S$ \> $Q\vee S \supset P\vee R$ \\[0.5ex]
\>\hspace{4em}$P\vee R \equiv Q\vee S$
\end{tabbing}

Syll.\zl ogism\zd

\begin{tabbing}
~~~~\hspace{3em}\=\underline{$P \equiv Q\; ,\; Q \equiv S \quad \lfloor : \rfloor \quad P \equiv S$} \\[.5ex]
\> \underline{$P \equiv Q \quad \lfloor : \rfloor \quad Q \equiv P$} \\[2ex]
\> \underline{$p \equiv p$} ~~~~~\qquad $p \supset p$ \quad $p \supset p$ \quad $(\f{P}{p}\: \f{Q}{p})$ \\[.5ex]
\`\` \zl fraction bars omitted in manuscript\zd\ \\[.5ex]
\> \underline{$p \equiv \: \sim \sim p$} \qquad $p \supset \: \sim \sim p$ \quad $ \sim \sim p \supset p$\\[.5ex]
\> \underline{$\sim (p\: .\: q)\equiv \: \sim p \: \vee \sim q$}\\[0.5ex]
\> $\sim \sim (\sim p\: \vee \sim q)\equiv \: \sim p \: \vee \sim q$\\[0.5ex]
\> \underline{$\sim (p \vee q)\equiv \: \sim p\:.\sim q$}\\[0.5ex]
\> ~~~~~~~~~~~~~$\equiv \: \sim (\sim \sim p \: \vee \sim \sim q)$\\[0.5ex]
\> $p \equiv \: \sim \sim p$ \` Forts\zl German: continued\zd\  p 4. F\\[0.5ex]
\> \underline{$q \equiv \: \sim \sim q$}\\[0.5ex]
\> $p \vee q \equiv \: \sim \sim p \: \vee \sim \sim q$ \quad $\mid$ \quad $\sim (p \vee q) \equiv \: \sim (\sim \sim p \: \vee \sim \sim q)$
\end{tabbing}

\begin{tabbing}
$\mathbf{\llbracket 6a. \rrbracket}$\\*[1ex]
23.\zl written over unreadable figure\zd\ \\*
\hspace{4em}\=\underline{$p\vee(q\: .\:r) \equiv (p\vee q)\: .\: (p \vee r)$}\\[.5ex]
~~~~~1.) \>$p \supset p \vee q$ \\[.5ex]
\>$p \supset p \vee r$ \\[.5ex]
~~~~~~~~~~$\llcorner$\> $p\supset (p\vee q)\: .\: (p \vee r)$ \\[.5ex]
\>$q\: .\:r \supset p \vee q$ \quad bec.\zl ause\zd\ $q\: .\:r \supset q$\\[.5ex]
\> \underline{$q\: .\:r \supset p \vee r$}\\[.5ex]
~~~~~~~~~~$\llcorner$\> $q\: .\:r \supset (p\vee q)\: .\: (p \vee r)$\\[.5ex]
~~~~~2.) $\llcorner$\> $p \supset [(p \vee q) \supset (p \vee q\: .\:r)]$ $\times$ \sout{bec.\zl ause\zd\ $(p \vee q) \supset [p \supset (p \vee q\: .\:r$\zl )]\zd } \\[.5ex]
\Big[ \sout{\zl unreadable word\zd}\\
\> $r \supset [(p \vee q) \supset (p \vee q\: .\:r)]$\Big] \quad \sout{because}\\[.5ex]
\> $r \supset [q \supset q\: .\:r]$\\[.5ex]
\> $q \supset q\: .\:r \supset [(p \vee q)\supset (p \vee q\: .\:r)]$ \quad Summation\\[.5ex]
~~~~~~~~~~$\llcorner$\> $r\supset [(p \vee q) \supset (p \vee q\: .\:r)]$\\[.5ex]
\> $(p \vee r)\supset [(p \vee q)\supset (p \vee q\: .\:r)]$\\[.5ex]
\> $( p\vee r)\: .\: (p \vee q)\supset (p \vee q\: .\:r)$\\[2ex]
$\times$ bec.\zl ause\zd\ \quad \=$p\supset p \vee q\: .\:r$\\[.5ex]
\>$p \vee q\: .\:r \supset [(p \vee q)\supset(p \vee q\: .\:r)]$\\[.5ex]
\>$p \supset [(p \vee q)\supset(p \vee q\: .\:r)]$
\end{tabbing}

\noindent \zl Here ends the page numbered $\mathbf{6a}$. The following not numbered page, until p.\ \textbf{7}., is crossed out.\zd

\begin{tabbing}
~~~~\hspace{3em}\=$(p \supset p\: .\:q)\supset[(p\supset (p\: .\:q \vee p\: .\:r)]$\\[.5ex]
\> $q\supset [p\supset (p\: .\:q \vee p\: .\:r)]$\\[.5ex]
\> $r\supset [p\supset (p\: .\:q \vee p\: .\:r)]$\\[.5ex]
\> $(q\vee r)\supset [p\supset (p\: .\:q \vee p\: .\:r)]$ \qquad importation\\[.5ex]
\> $(q\vee r)\: .\:p \supset (p\: .\:q \vee p\: .\:r)$\\[.5ex]
\> $p\: .\:(q\vee r)\supset (p\: .\:q \vee p\: .\:r)$
\end{tabbing}

\begin{tabbing}
~~~~\hspace{3em}\=$(p\: \vee \sim p)\lfloor .\rfloor (q\: \vee \sim q)\lfloor .\rfloor (r\: \vee \sim r)$\\[.5ex]
\> $p\: .\:q\: .\:(r\: \vee \sim r)\,\vee \sim p \: . \sim q$\\[.5ex]
\> $p\: .\:r\lfloor .\rfloor (q\: \vee \sim q)$
\end{tabbing}

\begin{tabbing}
~~~~\hspace{3em}\=$\lfloor (\rfloor \sim p\: \vee q)\lfloor .\rfloor (\sim p\: \vee \sim q)$ \\[.5ex]
\> $\sim p \vee q\: . \sim p$ \qquad $\vee \:q\: . \sim q$\\[.5ex]
\> $p\: .\:q\supset r$
\end{tabbing}

\noindent  $\mathbf{\llbracket 7. \rrbracket}$ Syllog.\zl ism\zd\ under an assumpt.\zl ion\zd

\begin{tabbing}
14$_R$\zl R\zd\ \hspace{3em}\=\underline{$P\supset (Q \supset R)\, ,\,P\supset (R \supset S)\quad .\:.\lfloor :\rfloor \quad P\supset (Q \supset S)$}\\[.5ex]
\> and similarly for any num.\zl ber\zd\ of premises\\[.5ex]
\>$P\supset (Q \supset R)\: .\:(R \supset S)$\\[.5ex]
\> \underline{$(Q \supset R)\: .\:(R \supset S)\supset Q \supset S$} \qquad exp.\zl ortation\zd\ syll.\zl ogism\zd \\[.5ex]
\>$P\supset (Q \supset S)$ \qquad \underline{also generalized}
\end{tabbing}
\vspace{0.5ex}
{\setlength{\parindent}{0cm}
$\left[
\begin{tabular}{@{}l@{}}
 14.1$_R$\zl R\zd\ \hspace{3em} \underline{$P\supset Q \qquad P \supset (Q \supset R) \quad \lfloor :\rfloor \quad P\supset R$} \\[.5ex]
\hspace{7.3em} $P\supset (Q \supset R)\: .\:Q$ \\[.5ex]
 \hspace{7.3em} \underline{$(Q \supset R)\: .\:Q \supset R$} \\[.5ex]
 \hspace{7.3em} $P\supset R$ \qquad Syll\zl ogism\zd
\end{tabular}
\right]$
}

\vspace{3ex}
\noindent \zl The following four lines, up to p.\ \textbf{8}., are in a box in the manuscript and are crossed out:
\begin{tabbing}
\hspace{2.5em}\= $(r \vee q \supset q \vee r)$ \\[.5ex]
\> \underline{$(r\supset s)\supset (q \vee r \supset q \vee r)$} \\[.5ex]
\> Add.\zl ition\zd\ of assumpt.\zl ions\zd \\[.5ex]
\> $(p\supset q)\: .\:(r\supset s) \supset [(p \vee r)\supset (q \vee s)]$\zd\
\end{tabbing}
\begin{tabbing}
$\mathbf{\llbracket 8. \rrbracket}$ \hspace{.6em}\=\underline{$(p\supset q)\: .\:(r\supset s) \supset (p \vee r \supset q \vee s)$}\\[.5ex]
1. \> $p \vee r \supset r \vee p$\\[.5ex]
2. \> $(p \supset q)\supset (r \vee p \supset r \vee q)$\\[.5ex]
3. \> $r \vee q \supset q \vee r$\\[.5ex]
4. \> $(r \supset s)\supset (q \vee r \supset q \vee s)$\\[.5ex]
5. \> $(p \supset q)\: .\: (r \supset s)\supset (p \vee r \supset r \vee p)$\\[.5ex]
6. \> $(p \supset q)\: .\: (r \supset s)\supset (r \vee p \supset r \vee q)$\\[.5ex]
7. \> $(p \supset q)\: .\: (r \supset s)\supset (r \vee q \supset q \vee r)$\\[.5ex]
8. \> \underline{$(p \supset q)\: .\: (r \supset s)\supset (q \vee r \supset q \vee r)$}\\[.5ex]
9. \> \underline{$(p \supset q)\: .\: (r \supset s)\supset (p \vee r \supset q \vee s)$}\\[.5ex]

\>\underline{$(p\supset q)\: .\:(r\supset q) \supset (p \vee r \supset q)$}\\[.5ex]
\> $(p\supset q)\: .\:(r\supset q) \supset (p \vee r \supset q \vee q)$ \qquad $\f{q}{s}$ \\[.5ex]
\> $(p\supset q)\: .\:(r\supset q) \supset (q \vee q \supset q)$\\[.5ex]
\> $(p\supset q)\: .\:(r\supset q) \supset (p \vee r \supset q)$\\[1ex]
$\mathbf{\llbracket 9. \rrbracket}$\>\underline{$(p\supset q)\: .\:(r\supset s) \supset (p \: .\: r \supset q \: .\: s)$}\\[.5ex]
\> $(p\supset q)\supset (\sim q \supset \: \sim p)$\\[.5ex]
\> $(r\supset s)\supset (\sim s \supset \: \sim r)$\\[.5ex]
$A.$ \> $(p\supset q)\: .\:(r\supset s) \supset (\sim q \supset \: \sim p) \: .\: (\sim s \supset \: \sim r)$\\[.5ex]
$B.$ \> $(\sim q \supset \: \sim p) \: .\: (\sim s \supset \: \sim r) \supset (\sim q \: \vee \sim s \supset \: \sim p \: \vee \sim r)$\\[.5ex]
$C.$ \> $(\sim q \: \vee \sim s \supset \: \sim p \: \vee \sim r) \supset (p \: .\: r \supset q \: .\: s)$\\[.5ex]
\> $(p\supset q)\: .\:(r\supset s) \supset (p \: .\: r \supset q \: .\: s)$ \qquad $A,B,C$\\[.5ex]
\>{$(p\supset q)\: .\:(p\supset s) \supset (p \supset q \: .\: s)$}\\[.5ex]
\> $(p\supset q)\: .\:(p\supset s) \supset ( p\: .\: p \supset q \: .\: s)$\\[.5ex]
\> $(p\supset q)\: .\:(p\supset s) \supset (p \supset p\: .\: p)$\\[.5ex]
\> $(p\supset q)\: .\:(p\supset s) \supset (p \supset q \: .\: s)$\\[1ex]
\>\underline{$\lfloor (\rfloor p\supset \: \sim p)\supset \: \sim p$}\\[.5ex]
\> $\sim p \: \vee \sim p \supset \: \sim p$\\[1ex]
$\mathbf{\llbracket 10. \rrbracket}$ \>\underline{$(\sim p\supset p)\supset p$}\\[.5ex]
\> $(\sim \sim p \vee p) \supset p$\\[.5ex]
\> $\sim \sim p \supset p$\\[.5ex]
\> \underline{$p \supset p$}\\[.5ex]
\> $(\sim \sim p \vee p) \supset p$\\[.5ex]
\>\underline{$\sim (p \: . \sim p)$} \qquad siehe unten$^\ast$ \zl German: see below\zd\ \\[.5ex]
\>\underline{$(p \supset q).(p \supset \: \sim q)\supset \: \sim p$}\\[.5ex]
\> $(p \supset q)\: .\: (p \supset \: \sim q)\supset[p \supset (q \: . \sim q)]$\\[.5ex]
\> $p \supset (q \: . \sim q)\supset (\sim(q \: . \sim q)\supset \: \sim p)$\\[.5ex]
\> $(p \supset q)\: .\: (p \supset \: \sim q)\supset (\sim(q \: . \sim q)\supset \: \sim p)$ ~ $\diagdown$\\
\` Princ.\zl iple of\zd\ Com\zl mutativity\zd\ \\
\> $\sim(q \: . \sim q)\supset [(p \supset q).(p \supset \: \sim q)\supset \:\sim p]$~~~~ $\diagup$\\[.5ex]
\> $(p \supset q).(p \supset \: \sim q)\supset \: \sim p$
\end{tabbing}

\noindent \zl The caption in the right margin in the following line from the manuscript is deleted, as well as the dash connecting it to a crossed out formula.\zd\

\begin{tabbing}
\`Princ.\zl iple of\zd\ Com\zl mutativity\zd\ \\
\hspace{3em}\= \sout{$(q \vee \sim q)\supset [(p \supset q)\: .\:(p \supset \: \sim q)\supset \:\sim p]$}~~ $\diagup$\\[.5ex]
\> \sout{$(p \supset q)\: .\:(p \supset \: \sim q)\supset \sim p$}\\[1ex]
\>\underline{$\sim (p \: . \sim p)$}\\[.5ex]
$^\ast$ \> $\sim \sim(\sim p \: \vee \sim \sim p)$
\end{tabbing}

$\mathbf{\llbracket 11. \rrbracket}$ Now I can \sout{now} \ul proceed \ud to the \ul proof of the \ud completeness theorem announced in the beg.\zl inning\zd\ which says that any tautology whatsoever can actually be derived in a finite number of steps from our four axioms by application of the 3\zl three\zd\ primitive rules of inf.\zl erence\zd\ (subst\zl itution\zd , implic\zl ation\zd , defined symbol) or shortly ,,\zl ``\zd Every tautology is demonstra\-ble''\zl .\zd\ \zl the following inserted text is crossed out: \ul since ,,\zl ``\zd demonstrable'' was defined to mean derivable from the 4\zl four\zd\ ax.\zl ioms\zd\ by the 3\zl three\zd\ rules of inf\zl erence.\zd\ \ud \zd\ We have already proved the inverse theor\zl em\zd\ which says: ,,\zl ``\zd Every demonstrable \sout{prop\zl osition\zd }\, \ul expression \ud is a taut\zl olo\-gy\zd ''\zl .\zd\ \sout{because of the following facts:} \zl The following assertions numbered 1. and 2. are crossed out in the manuscript:

\begin{itemize}
\item[1.] Each of the four ax.\zl ioms\zd\ is a tautology (as can \ul easily \ud be checked up\zl ,\zd\ \ul e.g.\ \ud by the truth t\zl able\zd\ method)
\item[$\mathbf{\llbracket 12. \rrbracket}$]
\item[2.] The 3\zl three\zd\ prim.\zl itive\zd\ rules of inf\zl erence\zd\ give \sout{only} tautologies \ul as conclu\-sions \ud if the premises are tautologies\zl ,\zd\ i\zl .\zd\ e\zl .\zd\ applied to tautologies they give again tautol\zl ogies\zd .\zd\
\end{itemize}

But the prop.\zl osition\zd\ which we are interested in now \zl \sout{in}\zd\ is the \sout{other} \ul inverse \ud one, which says ,,\zl ``\zd Any tautology is dem.\zl onstrable\zd ''. In order to prove it we have to use again the formulas $P_i$ which \sout{I (needed)} \ul we used \ud for proving that any truth\zl $\,$\zd table f\zl u\zd nct\zl ion\zd\ can be expressed by $\sim$ and $\vee$. If we have say $n$ propositional var.\zl iables\zd\ $p_1\lfloor ,\rfloor p_2\lfloor ,\rfloor p_3\lfloor ,\rfloor\ldots\lfloor ,\rfloor p_n$ then consider the conj.\zl unction\zd\ of them $p_1\: .\:p_2\: .\:p_3\: .\:\ldots\: .\:p_n$ and call a ,,\zl ``\zd fund.\zl amental\zd\ conj\zl unction\zd '' of these $\mathbf{\llbracket 13. \rrbracket}$ letters \ul $p_1\lfloor ,\rfloor \ldots\lfloor ,\rfloor p_n$ \ud any expression obtained from this conj.\zl unction\zd\ by negating some or all of the variables $p_1\lfloor ,\rfloor \ldots\lfloor ,\rfloor p_n$. So e.g\zl .\zd\ $p_1\: . \sim p_2\lfloor .\rfloor p_3\lfloor .\rfloor\ldots\lfloor .\rfloor$ $p_n$ would be a fund.\zl amental conjunction,\zd\ another one $\sim p_1\lfloor .\rfloor p_2\lfloor .\rfloor \sim p_3\lfloor .\rfloor p_4 \lfloor .\rfloor$ $\ldots\lfloor .\rfloor p_n$ etc.\zl ;\zd\ in part.\zl icular\zd\ we count also $p_1\lfloor .\rfloor \ldots\lfloor .\rfloor p_n$ \ul itself \ud and $\sim p_1\lfloor .\rfloor \sim p_2 \lfloor .\rfloor \ldots\lfloor .\rfloor \sim p_n$ (\zl in\zd\ which all \ul var.\zl iables\zd\ \ud are neg\zl ated\zd ) as fund.\zl amen\-tal\zd\ conj\zl unctions\zd .

\begin{tabbing}
\ul\hspace{1.7em}\=$2$ for $1$\zl one\zd \qquad \=$p_1$\zl ,\zd \quad $\sim p_1$\\[0.5ex]
$2^2$\>$4$ for two\>$p_1\lfloor .\rfloor p_2$\zl ,\zd \quad $p_1\lfloor .\rfloor \sim p_2$\zl ,\zd \quad $\sim p_1\lfloor .\rfloor p_2$\zl ,\zd \quad $\sim p_1\lfloor .\rfloor \sim p_2$\\[0.5ex]
$2^3$\>$8$ for three\>$p_1\lfloor .\rfloor p_2\lfloor .\rfloor p_3$\zl ,\zd \hspace{3.5em}\= $p_1\lfloor .\rfloor p_2\lfloor .\rfloor \sim p_3$\zl ,\zd \\[0.5ex]
\>\>$p_1\lfloor .\rfloor \sim p_2\lfloor .\rfloor p_3$\zl ,\zd \>$p_1\lfloor .\rfloor \sim p_2\lfloor .\rfloor \sim p_3$\zl ,\zd\\[0.5ex]
\>\>$\sim p_1\lfloor .\rfloor p_2\lfloor .\rfloor p_3$\zl ,\zd \>$\sim p_1\lfloor .\rfloor p_2\lfloor .\rfloor \sim p_3$\zl ,\zd \\[0.5ex]
\>\>$\sim p_1\lfloor .\rfloor \sim p_2\lfloor .\rfloor p_3$\zl ,\zd \>$\sim p_1\lfloor .\rfloor \sim p_2\lfloor .\rfloor \sim p_3$ \ud
\end{tabbing}

So for the $n$ var.\zl iables\zd\ $p_1\lfloor,\rfloor \ldots\lfloor,\rfloor p_n$ there are exactly $2^n$ fund\zl amental\zd\ conj.\zl unc\-tions\zd\ in gen\zl eral;\zd\ $2^n$ because \ul you see \ud by adding a \ul new \ud variable \ul $p_{n+1}$ \ud the num.\zl ber\zd\ of fund\zl amental\zd\ conj.\zl unctions\zd\ is doubled \zl ,\zd\ bec\zl ause\zd\ we can combine $p_{n+1}$ and $\sim p_{n+1}$ with any of the previous $\mathbf{\llbracket 14. \rrbracket}$ fund\zl amental\zd\ conj.\zl unctions\zd\ (as e.g.\ here $p_3$ with any of the prev\zl ious\zd\ $4$\zl four\zd\ and $\sim p_3$ getting $8$\zl eight\zd ) \sout{with each of these two poss.\zl ibilities\zd\ for any other var.\zl iable\zd\ so that we have alt.\zl ogether\zd\ $2 \times 2 \times \ldots 2=2^n$ possibilities.} ,\zl ;\zd\ I denote those $2^n$ fund.\zl amental\zd\ con\-j.\zl unctions\zd\ for the var.\zl iables\zd\ $p_1\lfloor,\rfloor \ldots\lfloor,\rfloor p_n$ by $P^{(n)}_1$\zl ,\zd $P^{(n)}_2$\zl ,\zd\ \ldots \zl ,\zd $P^{(n)}_i$ \zl ,\zd \ldots \zl ,\zd $P^{(n)}_{2^n}$. I am using $^{(n)}$ as an upper ind.\zl ex\zd\ to indicate that we mean the fund\zl amental\zd\ con\-j.\zl unction\zd\ of the $n$ variables $p_1\lfloor ,\rfloor \ldots\lfloor ,\rfloor p_n$ . The order in which they are enumerated is arb.\zl itrary\zd . [We may stick e.g.\ to the order which we used in the truth\zl $\;$\zd tables\zl .\zd ] From our formu\zl las\zd\ consid.\zl ered\zd\ for $n=3$ we know \sout{also} \ul $\mathbf{\llbracket 14.1 \rrbracket}$ that to each of these fund.\zl amental\zd\ conj\zl unctions\zd\ $P^{(n)}_i$ corresp\zl onds\zd\ exactly one line in a truth\zl $\,$\zd table for \ul a funct\zl ion\zd\ of the \ud $n$ variables \ul $p_1\lfloor ,\rfloor \ldots\lfloor ,\rfloor p_n$ \ud in such a way that $P^{(n)}_i$ will be true in this line and false in all the others. So if we numerate the lines correspondingly we can say $P^{(n)}_i$ will be true in the $i^{\text{th}}$ line and false in all other lines. \ud

$\mathbf{\llbracket 15. \rrbracket}$ Now \ul in order to prove the completeness theor.\zl em\zd\ \ud I prove first the foll\zl owing\zd\ \sout{lemma} \ul aux.\zl iliary\zd\ theorem\zl .\zd\ \ud
\begin{itemize}
\item[] \uline{Let $E$ be any expr.\zl es\-sion\zd\ which contains no other prop\zl ositional\zd\ var.\zl ia\-bles\zd\ but $p_1\lfloor ,\rfloor \ldots\lfloor ,\rfloor p_n$ and $P^{(n)}_i$ any fund\zl a\-mental\zd\ conj.\zl unc\-tion\zd\ of the var.\zl iables\zd\ $p_1\lfloor ,\rfloor \ldots\lfloor ,\rfloor p_n$\zl .\zd\ \zl T\zd hen either $P^{(n)}_i \supset E$ or $P^{(n)}_i \supset \: \sim E$ is de\-monstrable}
\end{itemize}
\ul \zl  exclamation mark deleted\zd\ \uline{\zl where\zd\ by either or I mean at least one}
\begin{tabbing}
\hspace{13.7em}$E$\\
Ex.\zl ample\zd\: \=$p_1 \lfloor .\rfloor p_2 \lfloor .\rfloor p_3 \supset [p\: .\:q \supset r]$ \hspace{7em} $p_1 \lfloor .\rfloor \sim p_2 \lfloor .\rfloor p_3$ \\[.5ex]
\> $\mid$\underline{$p_1 \lfloor .\rfloor \sim p_2 \lfloor .\rfloor p_3 \supset (p_1\lfloor .\rfloor p_2 \supset p_3)$}$\mid$ \qquad or\\[.5ex]
\> $p_1 \lfloor .\rfloor \sim p_2 \lfloor .\rfloor p_3 \supset\: \sim (p_1\lfloor .\rfloor p_2 \supset p_3)$\\[.5ex]
\> $\mid$\underline{$\sim p \:. \sim q \: .\: r \supset\: \sim (p \: . \: q \supset r)$}$\mid$ \ud
\end{tabbing}

It is to be noted that $E$ need not actually contain all the var.\zl iables\zd\ $p_1\lfloor ,\rfloor \ldots\lfloor ,\rfloor$ $p_n$\zl ;\zd\ it is only required that it contains no other variables but $p_1\lfloor ,\rfloor \ldots\lfloor ,\rfloor p_n$. So e.g.\ $p_1 \vee p_2$ would be an expr.\zl ession\zd\ for which the theor.\zl em\zd\ applies\zl ,\zd\ i\zl .\zd e.\
\begin{equation*}
\begin{rcases}
P^{(n)}_i\supset(p_1 \vee p_2)\;\;\;\;\\
\;\;\;\;\;\supset \; \sim (p_1 \vee p_2)
\end{rcases}
\text{dem.\zl onstrable\zd }
\end{equation*}
\sout{Let us first consider what that means} \zl The note ``$\mid$p19$\mid$'' in the manuscript at the bottom of this page, p.\ \textbf{15}., is deleted.\zd\

$\mathbf{\llbracket 16. \rrbracket}$ \zl The number of the page as well as the following text until the second half of p.\ \textbf{19}.\ starting with ``I shall prove'' are crossed out in the manuscript, while pp.\ \textbf{17}.-\textbf{18}.\ are missing from it:\zd

It is clear \ul at first \ud that under the ass.\zl umption\zd\ ment.\zl ioned\zd\ either $P^{(n)}_i \supset E$ or $P^{(n)}_i \supset \: \sim E$ must be a tautology bec\zl ause\zd : Let us write down the truth\zl $\,$\zd t\zl a\-ble\zd\ of the expr.\zl ession\zd\ $E$ \sout{it will have} (Note that we can consider $E$ as a funct\zl ion\zd\ in $n$ var.\zl iables\zd\ \sout{which is possible} also if it should not \ul actually \ud cont.\zl ain\zd\ all of the var.\zl iables\zd\ we have\zl ;\zd\ e.g\zl .\zd\ \sout{$p$} considered \ul $p$ \ud as a f\zl u\zd nct\zl ion\zd\ of $p,q$ and written down its truth\zl $\,$\zd table and gen\zl erally\zd\ if $E$ cont.\zl ains\zd\ say \ul only \ud $p_1\lfloor ,\rfloor \ldots\lfloor ,\rfloor p_k$ then its truth\zl $\,$\zd value is det\zl ermined\zd\ by the truth values of $p_1\lfloor ,\rfloor \ldots\lfloor ,\rfloor p_k$ hence a fortiori the truth \zl $\,$\zd val.\zl ues\zd\ of $p_1\lfloor ,\rfloor \ldots\lfloor ,\rfloor p_k \lfloor , \rfloor \ldots\lfloor ,\rfloor p_n$) $\mathbf{\llbracket 19. \rrbracket}$ differ from each other only in so far as some of the def.\zl ined\zd\ symb.\zl ols\zd\ are replaced by their definiens in $E'$. Sim.\zl ilarly\zd\ $P \supset \: \sim E^*_i$ can be der.\zl ived\zd\ from $P \supset \: \sim E'_i$. Hence we have: If one of the$\Big\rfloor$ \zl The whole of the text from the beginning of this page to this point is crossed out in the manuscript.\zd

$\Big\lceil$I shall prove \sout{that} \ul the aux.\zl iliary\zd\ theor\zl em\zd\ \ud only for such expressions as contain only the primit.\zl ive\zd\ symbols $\sim , \vee$ (but \ul do \ud not \ul contain \ud $\supset \lfloor , \rfloor \equiv$) bec. \zl ause\zd\ that is suff.\zl icient\zd\ for our purpose, and I prove it by a kind of complete induction , which we used already once in order to show that $\vee$ cannot be defined in terms of $\sim , \equiv$\zl .\zd\ $\mathbf{\llbracket 20. \rrbracket}$ Namely I shall prove the foll.\zl owing\zd\ three lemmas:

\begin{itemize}
\item [1.] The theorem is true for the simplest kind of expr.\zl ession\zd\ \ul $E$ \ud \zl ,\zd\ namely the var\zl iables\zd\ $p_1\lfloor ,\rfloor \ldots\lfloor ,\rfloor p_n$ themselves\zl ,\zd\ i.e.\ for any variable $p_k$ \ul of the above series $p_1\lfloor ,\rfloor \ldots\lfloor ,\rfloor p_k$ \ud and any fund\zl amental\zd\ conj\zl unction\zd\ $P^{(n)}_i$\zl ,\zd\ $P^{(n)}_i \supset p_k$ or $P^{(n)}_i \supset \: \sim p_k$ is demonstrable\zl .\zd
\vspace{-1ex}
\item[2.] If \ul the theor\zl em\zd\ \ud is true for an expr.\zl ession\zd\ $E$\zl ,\zd\ then it is also true for the neg\zl ation\zd\ $\sim E$\zl .\zd
\vspace{-1ex}
\item[3.] If it \ul is \ud true for two expr.\zl essions\zd\ $G,H$ then it is also true for the expression $G \vee H$\zl .\zd
\end{itemize}

After having proved these three lemmas we are finished. Because any expr.\zl es\-sion\zd\ $\mathbf{\llbracket 21. \rrbracket}$ $E$ containing only the var\zl iables\zd\ $p_1\lfloor ,\rfloor \ldots\lfloor ,\rfloor p_n$ and the op.\zl erations\zd\ $\sim \lfloor ,\rfloor \vee$ is formed by iterated appl.\zl ication\zd\ of the op\zl erations\zd\ $\sim , \vee$ beginning with the var.\zl iables\zd\ $p_1\lfloor ,\rfloor \ldots\lfloor ,\rfloor p_n$. Now by (1\zl .\zd) we know that the theorem is true for the variables $p_1\lfloor ,\rfloor \ldots\lfloor ,\rfloor p_n$ and by (2\zl .\zd) \zl and\zd\ (3\zl .\zd) we know that it remains true if we form new expr.\zl essions\zd\ by appl\zl ica\-tion\zd\ of $\sim$ \zl and\zd\ $\vee$ to expr.\zl essions\zd\ for which it is true. Hence it will be true for any expr.\zl ession\zd\ of the considered \zl unreadable word, perhaps ``type'' or ``kind''\zd . So it remains only to prove these three aux.\zl iliary\zd\ propositions\zl .\zd

$\mathbf{\llbracket 22. \rrbracket}$ (1\zl .\zd) means: For any var.\zl iable\zd\ $p_k$ (of the series $p_1\lfloor ,\rfloor \ldots\lfloor ,\rfloor p_n$) and any fund.\zl amental\zd\ conj.\zl unction\zd\ $P^{(n)}_i$ either $P^{(n)}_i \supset p_k$ or $P^{(n)}_i \supset \: \sim p_k$ is dem\zl onstra\-ble\zd . But now the letter $p_k$ or the neg\zl ation\zd\ $\sim p_k$ must occur among the members of this \ul fund\zl amental\zd\ \ud conj.\zl unction\zd\ \ul $P^{(n)}_i$ \ud by def.\zl inition\zd\ of a fund\zl amental\zd\ conj\zl unction\zd . On the other hand we know that \zl f\zd or any conj.\zl unction\zd\ it is demonstr.\zl a\-ble\zd\ that the conj.\zl unction\zd\ implies any of its memb.\zl ers.\zd\ (I proved that explicitly for conj.\zl unctions\zd\ of $2$\zl two\zd and $3$\zl three\zd\ members and remarked that the same method will prove it for conj.\zl unctions\zd\ of any $\mathbf{\llbracket 23. \rrbracket}$ num.\zl ber\zd\ of members. \ul The exact proof would have to go by an ind\zl uction\zd\ on the num\zl ber\zd\ of members\zl .\zd\ For two\zl ,\zd\ proved. \zl A\zd ssume $P^{(n)}$ \zl has\zd\ $n$ members and $p$ \zl is\zd\ a var.\zl iable\zd\ among them\zl . T\zd hen $P^{(n)}$ \zl is\zd\ $P^{(n-1)} \: .\: r$\zl :\zd\

\zl new paragraph\zd\ $1.$ $p$ occurs in $P^{(n-1)}$\zl ;\zd\ then $P^{(n-1)} \supset p$\zl ,\zd\ hence $P^{(n-1)} \: .\: r$ $\supset p$\zl .\zd\

\zl new paragraph\zd\ $2.$ $r$ is $p$\zl ; then\zd\ $P^{(n-1)} \: .\: p \supset p$ \zl is\zd\ dem.\zl onstrable\zd\ \ud .) \zl H\zd ence if $p_k$ occurs among the members of $P^{(n)}_i$ then $P^{(n)}_i \supset p_k$ is demonstrable and if $\sim p_k$ occurs among them then $P^{(n)}_i \supset \: \sim p_k$ is demonstr\zl able\zd . So one of these two form.\zl ulas\zd\ is demonstr.\zl able\zd\ in any case and that is exactly the assertion of lemma (1\zl .\zd).

\zl new paragraph\zd\ Now to (2\zl .\zd)\zl ,\zd\ i.e.\ let us assume \ul the theor\zl em\zd\ is true for $E$\zl ,\zd\ i\zl .\zd e.\ for any fund\zl amental\zd\ conj\zl unction\zd\ $P^{(n)}_i$ either \ud $P^{(n)}_i \supset E$ or $P^{(n)}_i \supset \: \sim E$ is demonstrable and let us show that the theor\zl em\zd\ is true also for the expr.\zl ession\zd\ $\sim E$\zl ,\zd\ i\zl .\zd e.\ \ul for any $P^{(n)}_i$ \ud either $P^{(n)}_i \supset \: \sim E$ or $P^{(n)}_i \supset \: \sim (\sim E)$ is demonstr\zl able\zd\ \ul for any $P^{(n)}_i$ \ud \zl The following formulae mentioned in this paragraph are in the manuscript on the right of the present page:\zd\
\begin{center}
\begin{tabular}{l|l}
$P^{(n)}_i \supset E$ & $P^{(n)}_i \supset \: \sim E$ \\
$P^{(n)}_i \supset \: \sim E$ & $P^{(n)}_i \supset \: \sim (\sim E)$ \\
\end{tabular}
\end{center}
(bec.\zl ause\zd\ it is $\mathbf{\llbracket 24. \rrbracket}$ this what the theor.\zl em\zd\ says if applied to $\sim E$ \zl \sout{it says:}\zd )\zl .\zd\
But now \ul in the 1.\zl first\zd\ case \ud if $P^{(n)}_i \supset E$ \ul is dem.\zl onstrable\zd\ then $P^{(n)}_i \supset \: \sim (\sim E)$ is also dem.\zl onstrable\zd\ bec.\zl ause\zd\ $E \supset \: \sim (\sim E)$ is dem\zl onstrable\zd\ by subst\zl itution\zd\ in the law of double neg.\zl ation,\zd\ and if \ul both \ud $P^{(n)}_i \supset E$ and $E \supset \: \sim (\sim E)$ are dem\zl onstrable\zd\ \zl semicolon deleted\zd\ then also $P^{(n)}_i \supset \: \sim (\sim E)$ by the rule of syllog\zl ism\zd . So we see if the first case is real.\zl ized\zd\ for $E$ then the sec.\zl ond\zd\ case is real\zl ized\zd\ for $\sim E$ and of course if the sec.\zl ond\zd\ case is real.\zl ized\zd\ for $E$ the first case is realis\zl z\zd ed for $\sim E$ (bec\zl ause\zd\ they say the same thing)\zl .\zd\ $\mathbf{\llbracket 25. \rrbracket}$ So if one of the two cases is real.\zl ized\zd\ for $E$ then also one of the two cases is real.\zl ized\zd\ for $\sim E$\zl ,\zd\ i.e.\ if \zl the\zd\ theor.\zl em\zd\ is true for $E$ it is also true for $\sim E$ which was to be proved\zl .\zd

\zl new paragraph\zd\ Now to (3\zl .\zd)\zl .\zd\ Assume \zl the\zd\ theor\zl em\zd\ true for $G \lfloor, \rfloor H$ and let $P^{(n)}_i$ be any arb.\zl i\-trary\zd\ fund.\zl amental\zd\ conj\zl unction\zd\ of $p_1\lfloor ,\rfloor \ldots\lfloor ,\rfloor$ $p_n$. Then $P^{(n)}_i \supset G$ \zl is\zd\ dem\zl onstra\-ble\zd\ or $P^{(n)}_i \supset \: \sim G$ \zl is\zd\ dem.\zl onstrable\zd\ and $P^{(n)}_i \supset H$ \zl is\zd\ dem\zl onstrable\zd\ or $P^{(n)}_i \supset \: \sim H$ \zl is\zd\ dem\zl onstrable\zd\ by ass.\zl umption\zd\ and we have to prove from these assump.\zl \-tions\zd\ that also:

\begin{tabbing}
\hspace{1.7em}\=$P^{(n)}_i \supset G\vee H$ \hspace{3em} \=or\\[0.5ex]
\>$P^{(n)}_i \supset \: \sim(G\vee H)$ \> is demonst\zl rable\zd .
\end{tabbing}

In order to do that dist.\zl inguish\zd\ three cases\zl :\zd\

\vspace{1ex}

\noindent $\mathbf{\llbracket 26. \rrbracket}$

\begin{itemize}
\item [1.] \ul [For $G$ I \zl first\zd\ case real\zl ized,\zd\ i\zl .\zd e.] \ud $P^{(n)}_i \supset G$ \zl is\zd\ dem.\zl onstrable;\zd\ then we have $G \supset G\vee H$ also by subst\zl itution\zd\ in ax.\zl iom,\zd\ hence $P^{(n)}_i \supset G\vee H$ ''\zl ``demonstrable'' in the edited text\zd\ by rule of syll.\zl ogism\zd\ [hence I \zl first\zd\ case real\zl ized\zd\ for $G \vee H$]\zl .\zd
\vspace{-1ex}
\item [2.] case \ul [For $H$ I \zl first\zd\ case real.\zl ized\zd ] \ud $P^{(n)}_i \supset H$ \zl is\zd\ dem\zl onstrable;\zd\ then $H \supset G\vee H$ by subst\zl itution\zd\ in form\zl ula\zd\ $10.$\zl , hence\zd\ $P^{(n)}_i \supset G\vee H$ \zl is\zd\ dem.\zl onstrable\zd\ by rule of syl.\zl logism\zd\ [hence I \zl first\zd\ case real\zl ized\zd\ for $G \vee H$.]\zl ].\zd
\vspace{-1ex}
\item [3.] case Neither for $G$ \zl is\zd\ \ul $P^{(n)}_i \supset G$ \zl \sout{nor}\zd\ \ud nor for $H$ \zl is\zd\ \ul $P^{(n)}_i \supset H$ \ud the I\zl first\zd\ case \zl \sout{is}\zd\ real\zl ized.\zd\ Thus for both of them sec.\zl ond\zd\ case happens\zl ,\zd\ i\zl .\zd e.\ $P^{(n)}_i \supset \: \sim G$ and $P^{(n)}_i \supset \: \sim H$ \zl are\zd\ both dem.\zl onstrable\zd\ \zl \sout{(bec\zl ause\zd }\zd\ by ass.\zl umption\zd\ \zl \sout{)}\zd\ \zl ,\zd\ but then by rule of transpos.\zl i\-tion\zd\ $G \supset \: \sim P^{(n)}_i$ \zl and\zd\ $H \supset \: \sim P^{(n)}_i$ \zl are\zd\ dem\zl on\-strable.\zd\ Hence $G \vee H \supset \: \sim P^{(n)}_i$ by rule of \zl Di\zd lemma\zl .\zd\ Hence $P^{(n)}_i \supset \: \sim (G \vee H)$ by transpos.\zl ition\zd\ [i.e.\ sec\zl ond\zd\ case realis\zl z\zd ed for $G \vee H$]\zl .\zd
\end{itemize}

$\mathbf{\llbracket 27. \rrbracket}$ So we see in each of the 3\zl three\zd\ cases which exh\zl aust\zd\ all poss\zl ibili\-ties\zd\ ei\-ther $P^{(n)}_i \supset G\vee H$ or $P^{(n)}_i \supset \: \sim(G\vee H)$ is dem\zl onstrable,\zd\ namely the first happens in case 1 \zl and\zd\ 2\zl ,\zd\ the sec.\zl ond\zd\ in \ul case \ud $3$. But that means that the theor\zl em\zd\ is true for $G \vee H$ since $P^{(n)}_i$ was any arb.\zl itrary\zd\ fund.\zl amental\zd\ conj\zl unction\zd . So we have proved the 3\zl three\zd\ lemmas and therefore the auxil.\zl iary\zd\ theor\zl em\zd\ for all expr.\zl essions\zd\ $E$ containing only $\sim , \vee$.

\zl new paragraph\zd\ Now let us assume in part\zl icular\zd\ that $E$ is a tautolo\-gie\zl y\zd\ of this kind (\ul i.e.\ \ud containing only the letters $p_1\lfloor ,\rfloor \ldots\lfloor ,\rfloor p_n$ and only $\sim \lfloor ,\rfloor \vee$)\zl ;\zd\ then I maintain $\mathbf{\llbracket 28. \rrbracket}$ that $P^{(n)}_i \supset E$ is demonstr\zl able\zd\ for any fund.\zl amental\zd\ conj.\zl unction\zd\ $P^{(n)}_i$\zl .\zd\ Now we know from the prec\zl eding\zd\ theor\zl em\zd\ that cert.\zl ain\-ly\zd\ either $P^{(n)}_i \supset E$ or \underline{$P^{(n)}_i \supset \: \sim E$} is demonstr\zl able\zd . So it remains only to be shown that the sec.\zl ond\zd\ \ul case\zl ,\zd\ that $P^{(n)}_i \supset \: \sim E$ is dem.\zl onstrable,\zd\ \ud can never occur if $E$ is \zl a\zd\ tautology and that can be shown as foll\zl ows\zd : As I ment\zl ioned\zd\ before any dem\zl onstrable\zd\ prop.\zl osition\zd\ is a taut\zl ology\zd . But on the other hand we can easily \ul see \ud that $P^{(n)}_i \supset \: \sim E$ is certainly not a taut.\zl ology\zd\ if $E$ is \zl a\zd\ taut.\zl ology\zd\ because the truth\zl $\,$\zd v\zl alue\zd\ of $P^{(n)}_i \supset \: \sim E$ will be false $\mathbf{\llbracket 29. \rrbracket}$ in the $i^{th}$ line of its truth\zl $\,$\zd t\zl able\zd . For in the $i^{th}$ line $P^{(n)}_i$ is true as we saw before and $E$ is also true in the $i^{th}$ line bec\zl ause\zd\ it is assumed to be a taut.\zl ology,\zd\ hence true in any line. Therefore $\sim E$ will be false in the $i^{th}$ line\zl ,\zd\ \sout{and} therefore $P_i \supset \: \sim E$ \ul will be false in the $i^{th}$ line \ud because $P_i$ is true and $\sim E$ false and therefore $P_i \supset \: \sim E$ false by the truth\zl $\,$\zd t.\zl able\zd\ of $\supset$. So this expr.\zl ession\zd\ \ul $P_i \supset \: \sim E$ \ud has F in the $i^{th}$ line of its truth\zl $\,$\zd t\zl able,\zd\ hence is not a tautology, hence cannot be demonstr.\zl able\zd\ and therefore $P^{(n)}_i \supset E$ is dem.\zl onstrable\zd\ for any fund\zl amental conjunction\zd\ $P^{(n)}_i$, if $E$ $\mathbf{\llbracket 30. \rrbracket}$ is a taut.\zl ology\zd\ containing only $\sim \lfloor ,\rfloor \vee \lfloor ,\rfloor p_1\lfloor ,\rfloor \ldots\lfloor ,\rfloor p_n$\zl .\zd\

But from the fact that $P^{(n)}_i \supset E$ is demonstrable for any $P^{(n)}_i$ it follows that $E$ is demonstr.\zl able\zd\ in the following way: We can show \underline{first} that also for any fund\zl amental\zd\ conj.\zl unction\zd\ \ul $P^{(n-1)}_i$ \ud of the $n-1$ var.\zl iables\zd\ $p_1\lfloor ,\rfloor \ldots\lfloor ,\rfloor p_{n-1}$\zl ,\zd\ $P^{(n-1)}_i \supset E$ is dem.\zl onstrable\zd\ bec.\zl ause\zd\ if $P^{(n-1)}_i$ is a fund\zl amental\zd\ conj.\zl unc\-tion\zd\ of the $n-1$ variables $p_1\lfloor ,\rfloor \ldots\lfloor ,\rfloor p_{n-1}$ then $P^{(n-1)}_i \: .\: p_n$ is a fund\zl amental\zd\ conj.\zl unction\zd\ of the $n$ var.\zl iables\zd\ $p_1\lfloor ,\rfloor \ldots\lfloor ,\rfloor$ $p_n$ and likewise $P^{(n-1)}_i \: . \sim p_n$ is a fund\zl amental\zd\ conj\zl unction\zd\ of the $n$ var.\zl i\-ables\zd\ $p_1\lfloor ,\rfloor \ldots\lfloor ,\rfloor p_n$\zl ;\zd\ therefore by our previous theor.\zl em\zd\ $\mathbf{\llbracket 31. \rrbracket}$ $P^{(n-1)}_i \: .\: p_n$ $\supset E$ and $P^{(n-1)}_i \: . \sim p_n \supset E$ are both demonst\zl rable.\zd\ \zl A\zd pplying the rule of exp.\zl ortation\zd\ \ul and commu\-t\zl ativity\zd\ \ud to those two expr.\zl essions\zd\ we get $p_n \supset (P^{(n-1)}_i \supset E)$ \zl and\zd\ $\sim p_n \supset (P^{(n-1)}_i \supset E)$ are both demonstr\zl able\zd . \ul \zl t\zd o be more exact we have to apply first the rule of exp.\zl ortation\zd\ and then the rule of commut.\zl ativity\zd\ bec.\zl ause\zd\ the rule of exp.\zl ortation\zd\ gives $P^{(n-1)}_i \supset (p_n \supset E)$ \zl .\zd\ \ud But now we can apply the rule of dilemma to these two form\zl ulas\zd\ ($P \supset R,Q \supset R \: \lfloor :\rfloor \: P \vee Q \supset R$) and obt\zl ain\zd\ $\sim p_n \vee p_n \supset (P^{(n-1)}_i \supset E)$ \ul is dem\zl onstrable;\zd\ and now since $\sim p_n \vee p_n$ is dem\zl onstrable\zd\ we can apply the rule of impl.\zl ication\zd\ \ul again \ud\ and obt.\zl ain\zd\ $P^{(n-1)}_i \supset E$ is dem.\zl onstrable\zd\ which was to be shown. Now since this holds $\mathbf{\llbracket 32. \rrbracket}$ for any fund\zl amental\zd\ conj.\zl unction\zd\ $P^{(n-1)}_i$ of the $n-1$ var.\zl iables\zd\ $p_1\lfloor ,\rfloor \ldots\lfloor ,\rfloor p_{n-1}$ it is clear that we can apply the same arg.\zl ument\zd\ again and prove that also for any fund\zl amental\zd\ conj.\zl unction\zd\ \ul $P^{(n-2)}_i$ \ud of the $n-2$ var.\zl iables\zd\ $p_1\lfloor ,\rfloor \ldots\lfloor ,\rfloor p_{n-2}$\zl ,\zd\ $P^{(n-2)}_i \supset E$ is dem\zl onstrable\zd . So by repeating this arg.\zl ument\zd\ $n-1$ times we can finally show that for any fund.\zl amental\zd\ conj\zl unction\zd\ of the one var\zl iable\zd\ $p_1$ this impl.\zl ication\zd\ is dem.\zl onstrable,\zd\ but that means $p_1 \supset E$ is dem\zl onstrable\zd\ and $\sim p_1 \supset E$ is dem.\zl onstrable\zd\ (bec\zl ause\zd\ $p_1$ and $\sim p_1$ are the fund\zl amental\zd\ conj.\zl unc\-tion\zd\ of the one var.\zl iable\zd\ $\mathbf{\llbracket 33. \rrbracket}$ \zl Above the page number in the manuscript the following list of rules and tautologies is written: Syll\zl ogism\zd , Transp\zl osition\zd , Dilemma, \underline{$p \;\vee \sim p$}, Export\zl ation\zd\ Com\zl mutativity\zd , \underline{$p \supset\; \sim \sim p$}\zd\ $p_1$)\zl ,\zd\ but then $\sim p_1 \vee p_1 \supset E$ is dem\zl onstrable\zd\ by rule of dil.\zl emma\zd\ and therefore $E$ is dem\zl onstrable\zd\ by rule of impl\zl ication\zd .

\ul Incident.\zl ally s\zd o we have shown that any taut.\zl ology\zd\ cont.\zl aining\zd\ only $\sim$ \zl and\zd\ $\vee$ is demonstr.\zl able,\zd\ but from this it follows that any taut.\zl olo\-gy\zd\ whatsoever is dem.\zl onstrable\zd\ bec\zl ause\zd : let $P$ be one containing \sout{perhaps} the def.\zl ined\zd\ symbols $. \: \lfloor ,\rfloor \supset \lfloor ,\rfloor \equiv$\zl .\zd\ I then denote by $P'$ the expr.\zl ession\zd\ \sout{form\zl ula\zd }\, obt.\zl ained\zd\ from $P$ by replacing $. \: \lfloor ,\rfloor \supset \lfloor ,\rfloor \equiv$ by their def.\zl iniens,\zd\ i\zl .\zd e\zl .\zd\ $R\: . \: S$ by $\sim (\sim R \: \vee \sim S)$ wherever it occurs in $P$ etc. Then $P'$ will \ul also \ud be a taut.\zl ology\zd\ \ul \zl \sout{bec}\zd\ \ud . \ul But $P'$ is a taut\zl ology\zd\ \ud containing only $\sim, \vee$ \sout{(truth table not changed)} hence $P'$ \zl is\zd\ dem.\zl onstrable,\zd\ but then also $P$ is dem\zl onstrable\zd\ bec\zl ause\zd\ it is obtained from $P'$ by one or several applications of the rule of def.\zl ined\zd\ symbol\zl ,\zd\ namely since $P'$ was obt\zl ained\zd\ from $P$ by rep\zl lacing\zd\ $p\: .\: q$ by $\sim (\sim p \vee \: \sim q)$ etc\zl .\zd\ $P$ is obt\zl ained\zd\ from $P'$ by the inv\zl erse\zd\ subst.\zl itution,\zd\ but each such subst.\zl itution\zd\ is an applic\zl ation\zd\ of rule of def.\zl ined\zd\ symbol\zl ,\zd\ hence: If $P'$ is demonstrable then also $P$ \zl is\zd\ dem\zl onstrable\zd . \ud

As an example take \zl the\zd\ form\zl ula\zd\ $(p \supset q) \vee (q \supset p)$ which is a tautol\zl ogy\zd .
\begin{tabbing}
\hspace{1.7em}\= 1. Without def\zl ined\zd\ symb.\zl ols\zd\ \quad $(\sim p \vee q) \vee (\sim q \vee p)=E$\\[.5ex]
\>2. Fund\zl amental\zd\ conj\zl unctions\zd\ \ul in \ud $p,q$\\
\`$p \lfloor .\rfloor q$\zl ,\zd\ \hspace{.5em} $p \lfloor .\rfloor \sim q$\zl ,\zd\ \hspace{.5em} $\sim p \lfloor .\rfloor q$\zl ,\zd\ \hspace{.5em}$\sim p \lfloor .\rfloor \sim q$
\end{tabbing}
 To prove that $p\: .\: q \supset E$ etc\zl .\zd \zl are\zd\ all dem\zl onstrable\zd \zl w\zd e have to verify our aux.\zl iliary\zd\ theor\zl em\zd\ successively for all particul\zl ar\zd\ form\zl ulas,\zd\ i.e.\ for $p$, $q$, $\sim p$, $\sim q$, $\sim p \vee q$, $\sim q \vee p$, $E$\zl .\zd\

\vspace{2ex}

\noindent $\mathbf{\llbracket 34. \rrbracket}$
\begin{center}
\begin{tabular}{ l|c|c|c|c|c|c}
 & $p$ & $q$ & $\sim p$ & $\sim q$ & $\sim p \vee q$ & $\sim q \vee p$\\[.5ex] \hline
 \hspace{1.1em}$p \lfloor .\rfloor q \supset$ & $p$ & $q$ & $\sim (\sim p)$ & $\sim (\sim q)$ & $\sim p \vee q$ & $\sim q \vee p$ \\
 \hspace{1.1em}$p \lfloor .\rfloor\! \sim q \supset$ & $p$ & $\sim q$ & $\sim (\sim p)$ & $\sim q$ & $\sim (\sim p \vee q)$ & $\sim q \vee p$ \\
 $\sim p \lfloor .\rfloor q \supset$ & $\sim p$ & $q$ & $\sim p$ & $\sim (\sim q)$ & $\sim p \vee q$ & $\sim (\sim q \vee p)$ \\
 $\sim p \lfloor .\rfloor\! \sim q \supset$ & $\sim p$ & $\sim q$ & $\sim p$ & $\sim q$ & $\sim p \vee q$ & $\sim q \vee p$
 \end{tabular}
\end{center}
\begin{flushright}
\begin{tabular}{ l|c }
  & $(\sim p \vee q) \vee (\sim q \vee p)$\\[.5ex] \hline
  & $E$\\
  & $E$\\
  & $E$\\
  & $E$
 \end{tabular}
\end{flushright}

\begin{tabbing}
\hspace{1.7em}\=$p\: . \sim q \supset \: \sim (\sim p)$ \qquad \=$\sim p \supset\: \sim (p \: . \sim q)$ \\[.5ex]
\> $p\: . \sim q \supset \: \sim q$ \qquad \> $q \supset\: \sim (p \: . \sim q)$ \\[.5ex]
\>\> $\sim p \vee q \supset \: \sim (p\: . \sim q) $\\[.5ex]
\>\> \underline{$p \:. \sim q \supset \: \sim (\sim p \vee q)$}\\[2ex]

\>$p\: .\: q \supset E$ \>$p \supset (q \supset E)$\\[.5ex]
\> $\sim p\: .\: q \supset E$ \>\underline{$\sim p \supset (q \supset E)$} \\
\>\> $\sim p \vee p \supset (q \supset E)$\\[.5ex]
\>\> $q \supset E$\\

\rule{8cm}{0.4pt}\\[1ex]

\>$p\: . \sim q \supset E$ \>$p \supset (\sim q \supset E)$\\[.5ex]
\> $\sim p\: . \sim q \supset E$ \>\underline{$\sim p \supset (\sim q \supset E)$} \\[.5ex]
\>\> $\sim p \vee p \supset (\sim q \supset E)$\\[.5ex]
\>\> $\sim q \supset E$
\end{tabbing}

\zl The following formulae, which in the manuscript are on the right of this page, are deleted:

\begin{tabbing}
\hspace{1.7em}\=$P^{(n)}_i \supset \sim A$ \qquad \=$\sim (A \vee B)$\\[.5ex]
\> ~~~~~~~$\supset \sim B$ \> $P^{(n)}_i \supset A$ \\[.5ex]
\>\> $A \supset (A \vee B)$\\[.5ex]
\> $A \supset \: \sim P^{(n)}_i$\\[.5ex]
\> $B \supset \: \sim P^{(n)}_i$\\[.5ex]
\> $A \vee B \supset \: \sim P^{(n)}_i$\\[2ex]
\> $p \vee q \supset E$ \> $\sim p \: .\: q$\\[.5ex]
\>$\underline{p} \supset (\underline{q} \supset E)$ \\[.5ex]
\>$\sim p \supset (q \supset E)$\zd
\end{tabbing}

\begin{tabbing}
\hspace{1.7em}\=$p\: . \sim q \supset \: \sim (\sim p)$ \qquad \=\kill

$\mathbf{\llbracket 35. \rrbracket}$ \;\;$\sim q \vee q \supset E$ \>\> \underline{$E$}\\
\rule{8cm}{0.4pt}
\end{tabbing}

Now after having proved the\zl at\zd\ any taut.\zl ology\zd\ can be derived from the 4\zl four\zd\ ax.\zl ioms,\zd\ the next quest.\zl ion\zd\ which arises is\zl comma from the manuscript deleted\zd\ whether all of those 4\zl four\zd\ ax\zl ioms\zd\ are really necessary to derive them or whether perhaps one \ul or the other \ud of them is superfluous\zl .\zd\ That would mean one of them could be left out and nevertheless the rem.\zl aining\zd\ three would allow to derive all taut\zl ologies\zd . If this were the case then in part.\zl icular\zd\ also the superfluous ax.\zl iom\zd\ (since it is a taut.\zl ology\zd ) could be derived from the three other, $\mathbf{\llbracket 36. \rrbracket}$ i\zl .\zd e.\ it would not be independent from the other. So the question comes down to investigating the indep.\zl endence\zd\ of the 4\zl four\zd\ ax.\zl ioms\zd\ from each other. That such an invest.\zl igation\zd\ is really nec.\zl essary\zd\ is shown very strikingly by the last development. Namely when Russell first set up this sys.\zl tem\zd\ of ax\zl ioms\zd\ for the calc.\zl ulus\zd\ of prop.\zl ositions\zd\ he assumed a fifth ax.\zl iom,\zd\ namely the associat.\zl ive\zd\ law for disj.\zl unction\zd\ and only many years later it was proved by \ul P.\ \ud Bern.\zl ays\zd\ that this ass.\zl ociative\zd\ law was superfluous\zl ,\zd\ i.e.\ could $\mathbf{\llbracket 37. \rrbracket}$ be derived from the others. You have seen in one of the prev.\zl ious\zd\ lect.\zl ures\zd\ how this derivation can be accomplished. But Bern\zl ays\zd\ has shown at the same time that a similar thing cannot happen for the 4\zl four\zd\ rem.\zl aining\zd\ axioms\zl ,\zd\ i\zl .\zd e.\ that they are really ind.\zl ependent\zd\ from each other\zl .\zd\ \ul

\zl new paragraph\zd\ Again here as in the completeness proof the interest does not ly\zl ie\zd\ so much in proving that these part\zl icular\zd\ 4\zl four\zd\ ax.\zl ioms\zd\ are independent but in the method to prove it, \zl b\zd ecause so far we have only had an opport.\zl unity\zd\ to prove that \sout{\zl unreadable word\zd} cert.\zl ain\zd\ prop.\zl ositions\zd\ follow from other prop\zl ositions\zd . But now we are confronted with the \ul opposite \ud problem to show that cert.\zl ain\zd\ prop.\zl ositions\zd\ do not follow from \ul certain \ud others and this problem requires evidently an entirely new method for its solu\-tion\zl .\zd\ \ud

\sout{And I intend to give his proof \ul here \ud for at least one of the ax.\zl ioms\zd\ bec.\zl ause\zd }\, This method is very interest\zl ing\zd\ and \ul somewhat\zd\ conn.\zl ected\zd\ with the quest\zl ions\zd\ of many\zl -\zd valued logics.

You know the calc.\zl ulus\zd\ of prop.\zl ositions\zd\ can be interpret.\zl ed\zd\ as an alg\zl ebra\zd\ in which $\mathbf{\llbracket 38. \rrbracket}$ \sout{in which} we have the two op.\zl erations\zd\ of log.\zl ical\zd\ add.\zl ition\zd\ and mult.\zl iplication\zd\ as in usual alg.\zl ebra\zd\ but in add.\zl ition\zd\ to them a 3\zl third\zd\ op.\zl era\-tion,\zd\ the negation and bes.\zl ides\zd\ some op\zl erations\zd\ def.\zl ined\zd\ in terms of them ($\supset, \equiv$ etc\zl .\zd ). The objects to which those op.\zl erations\zd\ are applied are the prop\zl osi\-tions\zd . So the prop.\zl ositions\zd\ can be made to corresp.\zl ond\zd\ to the numb.\zl ers\zd\ of ord.\zl in\-ary\zd\ alg\zl ebra\zd . But as you know all the op.\zl erations\zd\ $.\:\lfloor ,\rfloor \vee$ etc\zl .\zd\ which we introd.\zl uced\zd\ are ,,\zl ``\zd truth\zl $,$\zd f\zl u\zd nct\zl ions\zd\ '' and therefore it is only the truth\zl $\,$\zd value of the prop.\zl osi\-tions\zd\ that really matters in this alg.\zl ebra,\zd\ $\mathbf{\llbracket 39. \rrbracket}$ i\zl .\zd e.\ we can consider \ul them \ud as the numbers of our alg.\zl ebra\zd\ inst.\zl ead\zd\ of the prop.\zl osi\-tions\zd\ (simply the two ,,\zl ``\zd truth \zl $\,$\zd values'' T and F)\zl .\zd\ And this is what we shall do\zl ,\zd\ i\zl .\zd e.\ our alg.\zl ebra\zd\ (as opposed to usual alg.\zl ebra\zd ) has only two numbers T, F and the result of the op.\zl erations\zd\ $.\:,\vee ,\sim$ applied to these two num.\zl bers\zd\ is given by the truth\zl $\,$\zd t\zl able,\zd\ i\zl .\zd e.\ T $\vee$ F = T (i\zl .\zd e.\ the sum of the two nu\zl mbers\zd\ T and F is T) T $\vee$ T = T\zl ,\zd\ F $\vee$ T = T\zl ,\zd\ F $\vee$ F = F\zl ,\zd\ $\sim$T = F\zl ,\zd\ $\sim$F = T\zl .\zd\ In order to stress $\mathbf{\llbracket 40. \rrbracket}$ more the anal.\zl ogy\zd\ to alg.\zl ebra\zd\ I shall \ul also \ud write $1$ inst\zl ead\zd\ of T and $0$ inst\zl ead\zd\ of F. Then in this not.\zl ation\zd\ the rules for log.\zl ical\zd\ mult.\zl iplication\zd\ would look like this\zl :\zd\ $1\:.\:1=1$\zl ,\zd\ $0\:.\:1=0$\zl ,\zd\ $1\:.\:0=0$\zl ,\zd\ $0\:.\:0=0$\zl .\zd\
If you look at this table you see that log.\zl ical\zd\ and arithm.\zl etical\zd\ mult.\zl iplication\zd\ exactly coincide in this notation. Now what are the tautologies consider.\zl ed\zd\ from this algebraic standpoint? They are expr.\zl essions\zd\ $f(p \lfloor ,\rfloor q \lfloor ,\rfloor r \lfloor ,\rfloor \ldots)$ which have always the value $1$ whatever nu\zl mbers\zd\ $p,q,r$ may be\zl ,\zd\ $\mathbf{\llbracket 41. \rrbracket}$ i\zl .\zd e.\ in alg.\zl ebraic\zd\ language expressions ident.\zl ically\zd\ equ.\zl al\zd\ to one $f(p \lfloor ,\rfloor q \lfloor ,\rfloor \ldots)=1$ and the contrad.\ul ictions \ud expr.\zl essions\zd\ id.\zl entically\zd\ zero $f(p \lfloor ,\rfloor q \lfloor ,\rfloor \ldots)=0$\zl .\zd\ So an expr.\zl ession\zd\ of usual alg\zl ebra\zd\ which would corresp.\zl ond\zd\ to a contrad.\zl iction\zd\ would be e.g.\ $x^2-y^2-(x+y)(x-y)$\zl ;\zd\ this is =\zl equal to\zd 0\zl .\zd\

\zl new paragraph\zd\ But now from this algebr.\zl aic\zd\ standp.\zl oint\zd\ nothing \ul can \ud prevent us to consider also other sim.\zl ilar\zd\ alg\zl ebras\zd\ with say three nu\zl mbers\zd\ \ul$0,1,2$ \ud inst.\zl ead\zd\ of two and with the op.\zl erations\zd\ $\vee \lfloor ,\rfloor \:.\: \lfloor ,\rfloor \sim$ defined in some diff.\zl erent\zd\ manner. For any such alg.\zl ebra\zd\ we shall have taut.\zl ologies,\zd\ $\mathbf{\llbracket 42. \rrbracket}$ i\zl .\zd e.\ form\zl ulas\zd\ =\zl equal to\zd\ 1 and contr\zl adictions\zd\ =\zl equal to\zd\ 0\zl ,\zd\ but they will of course be diff\zl erent\zd\ form\zl ulas\zd\ for diff.\zl er\-ent\zd\ alg\zl ebras.\zd\ Now such alg.\zl ebra\zd\ with 3\zl three\zd\ and more nu\zl m\-bers\zd\ were used by Bern.\zl ays\zd\ for the proof of indep\zl endence,\zd\ e.g.\ in order to prove the ind.\zl ependence\zd\ of the sec\zl ond\zd\ ax.\zl iom\zd\ Bern\zl ays\zd\ considers the foll.\zl owing\zd\ alg\zl ebra\zd :
\begin{tabbing}
\hspace{1.7em}\=$3$ N\zl n\zd umbers \hspace{1em}\=$0,1,2$\\[.5ex]
\>neg\zl ation\zd \>$\sim 0=1$ \qquad $\sim 1=0$ \qquad \=$\sim 2=2$ \\[.5ex]
\>add\zl ition\zd \>$1 \vee x=x \vee 1=1$ \> $2 \vee 2=1$\\[.5ex]
\>\> $0 \vee 0=0$ \> $2 \vee 0=0 \vee 2=2$
\end{tabbing}
\zl The equations on the right involving $2$ are in a box in the manuscript.\zd

\sout{or $0 \vee x=x \vee 0=x$} \zl I\zd mpl.\zl ication\zd\ and other op\zl erations\zd\ \zl ``not nec to'' from the manuscript rendered by ``need not be''\zd\ def.\zl ined\zd\ sep.\zl arately\zd\ because $p \supset q=\; \sim p \vee q$\zl .\zd\

$\mathbf{\llbracket 43. \rrbracket}$ \zl A t\zd aut.\zl ology\zd\ =\zl is a\zd\ formula =\zl equal to\zd\ 1\zl ,\zd\ e.g.\ $\sim p \vee p$ bec.\zl ause\zd\ for $p=0$ or $=1$ $=1$\zl $p$ equal to 0 or 1 it is equal to 1,\zd\ bec.\zl ause\zd\ the op.\zl erations\zd\ for $0,1$ as arg\zl uments\zd\ coincide with the op.\zl erations\zd\ of \zl the\zd\ usual calc.\zl ulus\zd\ of prop.\zl ositions;\zd\ if $p=2$ \zl then\zd\ $\sim p=2$ \zl and\zd\ $2 \vee 2=1$ \zl is\zd\ also true. Also $p \supset p$ \zl is a\zd\ taut.\zl ology\zd\ bec.\zl ause\zd\ by def\zl inition it is\zd\ the same as $\sim p \vee p$.

Now for this alg.\zl ebra\zd\ one can prove the foll\zl owing\zd\ prop.\zl osition:\zd
\begin{itemize}
\item [1.] Ax\zl ioms\zd\ $\lfloor(\rfloor 1 \lfloor)\rfloor,\lfloor(\rfloor 3 \lfloor)\rfloor,\lfloor(\rfloor 4 \lfloor)\rfloor$ are taut.\zl ologies\zd\ in this alg\zl ebra\zd .
\vspace{-1ex}
\item [2.] For each of the \ul three \ud rules of inf\zl erence\zd\ we have\zl :\zd\ If the premises are taut\zl ologies\zd\ in this alg\zl ebra\zd\ then \zl \sout{al}\zd so \zl is\zd\ the concl\zl usion\zd .
\item[$\mathbf{\llbracket 44. \rrbracket}$] \zl I.\zd e.\
\begin{itemize}
\item [1.] If $P$ and $P \supset Q$ \zl are tautologies\zd\ then $Q$ \zl is a tautology.\zd

\item [2.] If $Q'$ by subst\zl itution\zd\ from $Q$ and $Q$ is a taut\zl ology\zd\ then also $Q'$ \zl is a tautology.\zd

\item [3.] If $Q'$ \zl is\zd\ obt\zl ained\zd\ from $Q$ by replacing $P \supset Q$ by $\sim P \vee Q$ etc\zl .\zd\ and $Q$ \zl is a tautology\zd\ then also $Q'$ \zl is a tautology.\zd
\end{itemize}
\item[3.] The ax.\zl iom\zd\ $\lfloor(\rfloor 2 \lfloor)\rfloor$ is not a taut\zl ology\zd\ in this alg\zl ebra\zd .
\end{itemize}

After having shown these $3$\zl three\zd\ lemmas we are finished bec\zl ause\zd\ by $1,2$\zl :\zd\ Any form\zl ula\zd\ dem\zl onstrable\zd\ from ax.\zl ioms\zd\ 1, 2, 4 \zl (1), (3), (4)\zd\ by $3$\zl the three\zd\ rules of inf\zl erence\zd\ is a taut\zl ology\zd\ for our alg\zl ebra\zd\ but ax\zl iom\zd\ $\lfloor(\rfloor3\lfloor)\rfloor$\zl $(2)$\zd\ is not a taut\zl olo\-gy\zd\ for our $\mathbf{\llbracket 45. \rrbracket}$ alg\zl ebra.\zd\ Hence it cannot be dem\zl onstrable\zd\ from $\lfloor(\rfloor 1 \lfloor)\rfloor,\lfloor(\rfloor 3 \lfloor)\rfloor,\lfloor(\rfloor 4 \lfloor)\rfloor$.

Now to the proof of the lemmas $1,2,3$. First some aux\zl iliary\zd\ theor.\zl ems\zd\ \ul \zl (\zd\ for $1$ I say T\zl t\zd rue and for $0$ false bec.\zl ause\zd\ for $1$ and $0$ the tables of our alg\zl ebra\zd\ coincide with those for T,\zl and\zd\ F\zl ):\zd\ \ud

\begin{itemize}
\item [1.] $p \supset p$ \quad \zl (\zd\ we had that before\zl ,\zd\ bec\zl ause\zd\ $\sim p \vee p=1$ also $\sim 2 \vee 2=1$\zl )\zd
\vspace{-1ex}
\item [2.] $1 \vee p=p \vee 1=1$ \qquad $0 \vee p=p \vee 0=p$
\vspace{-1ex}
\item [3.] $p \vee q=q \vee p$
\vspace{-1ex}
\item [4.] Also in our three\zl -\zd val\zl ued\zd\ algebra we have: An impl.\zl ication\zd\ whose first member is $0$ is $1$ and an impl\zl ication\zd\ whose sec.\zl ond\zd\ me\zl mber\zd\ is $1$ is also $1$ whatever the other memb\zl er\zd\ may be\zl ,\zd\ i\zl .\zd e.\ $0\supset p=1$ \zl and\zd\ $p\supset 1=1$ bec\zl ause:\zd
\end{itemize}

\begin{tabbing}
\hspace{1.7em}\=1.) $0\supset p= \: \sim 0 \vee p=1 \vee p=1$\\[1ex]
$\mathbf{\llbracket 46. \rrbracket}$\\*[1ex]
\>2.) $p \supset 1= \: \sim p \vee 1=1$\zl full stop deleted\zd\ \\[1ex]
Now I\zl (1)\zd\ \quad $p \supset p \vee q=1$\\[1ex]
\>1. $p=0$ \quad $\rightarrow$ \quad $p \supset p \vee q=1$\\[.5ex]
\>2. $p=1$ \quad $\rightarrow$ \quad $1 \supset 1 \vee q=1 \supset 1=1$\\[1ex]
III\zl (3)\zd\ $p \vee q=q \vee p$ \quad $\rightarrow$ \quad $p \vee q=q \vee p=1$\\[1ex]
IV\zl (4)\zd\ $(p \supset q)\supset (r \vee p \supset r \vee q)$ \quad $E$\\[1ex]
\>\underline{1.} \,$r=0$ \hspace{.5em} $r \vee p=p$ \hspace{.5em} $r \vee q=q$ \hspace{.5em} $E=(p \supset q)\supset (p \supset q)=1$\\[5ex]
\>\underline{2.} \,$r=1$ \hspace{.5em} $r \vee p=r\vee q=1$ \hspace{.5em} $E=(p \supset q)\supset (1 \supset 1)=(p \supset q)\supset 1=1$
\end{tabbing}
\begin{tabbing}
$\mathbf{\llbracket 47. \rrbracket}$\\*[1ex]
\hspace{1.7em}\=\underline{3.}\, $r= \underline{2}$\\

\>\hspace{1.7em}\=\underline{$\alpha .)$}\,\,\=$q=2,1$\zl 1, 2\zd\ \quad $r \vee q$\!\!\!

\begin{tabular}{ l }
 =\;\underline{$2 \vee 1$}\;=\;1 \\
 =\;\underline{$2 \vee 2$}\;=\;1 \\
 \end{tabular} \\[.5ex]
\>\>\>$r \vee p \supset r \vee q=1$\\[.5ex]
\>\>\>$(p \supset q)\supset(r \vee p\supset r \vee q)=1$\\[1ex]

\>\>\underline{$\beta .)$} \>$q=\underline{0}$\\[1ex]
\>\>\hspace{1.7em}\>\underline{1.}\;\:\= $p=0$ \qquad $r \vee p=r \vee q$\\[.5ex]
\>\>\>\>$(r \vee p)\supset(r \vee q)=1$\\[.5ex]
\>\>\>\>$(p \supset q)\supset (r\supset p)\supset(r \vee q)=1$\\[1ex]
\>\>\>\underline{2.}\> $p=1$ \qquad $p \supset q=0$\\[.5ex]
\>\>\>\>$E=1$\\[1ex]
\>\>\>\underline{3.}\> $p=2$ \\[.5ex]
\>\>\>\>$(2\supset 0)\supset(2 \vee 2 \supset 2 \vee 0)=2\supset(1 \supset 2)=2\supset 2=1$
\end{tabbing}

\noindent $\mathbf{\llbracket 48. \rrbracket}$ 2.\ Lemma\zl Lemma 2.\zd\ A. \quad $p=1$ \quad $p \supset q=1$ \quad $\rightarrow$ \quad $q=1$\\[.5ex]
$1=\:\sim p \vee q=0 \vee q=q$

\vspace{1ex}

Hence if $f(p\lfloor ,\rfloor q \lfloor ,\rfloor \ldots)=1$ \zl then\zd\ $$ \frac{f(p\lfloor ,\rfloor q \lfloor ,\rfloor \ldots) \supset g(p\lfloor ,\rfloor q \lfloor ,\rfloor \ldots)=1}{g(p\lfloor ,\rfloor q \lfloor ,\rfloor \ldots)=1}$$

\begin{itemize}
\item[B.] Rule of subst.\zl itution\zd\ holds for any truth-value algebra\zl ,\zd\ i\zl .\zd e.\ if $f(p\lfloor ,\rfloor q \lfloor ,\rfloor$ $\ldots)=1$ then $f(g(p\lfloor ,\rfloor q \lfloor ,\rfloor \ldots)\lfloor ,\rfloor q \lfloor ,\rfloor \ldots)=1$\zl .\zd\
\item[C.] Rule of defined symb\zl ol\zd\ likewise \zl holds\zd\ bec\zl ause\zd\ $p \supset q$ \zl and\zd\ $\sim p \vee q$ have the same truth\zl $\,$\zd table\zl .\zd\
\end{itemize}

$\mathbf{\llbracket 49. \rrbracket}$ \zl The following on the right of this page in the manuscript is deleted: gen\zl eral\zd\ remark about the mean\zl ing\zd\ of derivability from axioms.\zd\

\begin{tabbing}
\hspace{1.7em} Lemma 3.\ II.\zl (2)\zd\ $p \vee p \supset p$ is not a taut.\zl ology\zd\ \\[1.5ex]
\hspace{6em} i.e.\ $2 \vee 2 \supset 2$ = $1\supset 2=\:\sim 1 \vee 2=0 \vee 2=2 \neq 1$
\end{tabbing}
So the lemmas are proved and therefore also the theorem about the independence of Ax\zl iom\zd\ II\zl $(2)$.\zd\

We have already developed a method for deciding of any given expr.\zl es\-sion\zd\ whether or not it is a tautology\zl ,\zd\ namely the truth\zl -\zd table method. I want to develop another method which uses the analogy of the rules of the $\mathbf{\llbracket 50. \rrbracket}$ calc.\zl ulus\zd\ of prop.\zl osition\zd\ with the rules of algebra. We have the two distrib\zl utive\zd\ laws:

\begin{tabbing}
	\hspace{1.7em}\=$\underline{p}\:.\:(q\vee r)\equiv (\underline{p}\:.\:q)\vee(\underline{p}\:.\:r)$ \qquad \=$\underline{p}\:.\:q\equiv q$\\[0.5ex]
	\>$\underline{p}\vee (q\:.\: r)\equiv (\underline{p}\vee q)\:.\:(\underline{p}\vee r)$\> $\underline{p}\vee q\equiv q$
\end{tabbing}
 In order to prove them by the shortened truth\zl -\zd table method I use the following facts \ul which I ment.\zl ioned\zd\ already once at the occasion of one of the exercises\zl :\zd\ \ud

\begin{tabbing}
	\hspace{1.7em}\=if $p$ is true \qquad \=$p\:.\:q\equiv q$\\[.5ex]
	\>if $p$ is false \>$p \vee q\equiv q$
\end{tabbing}
In order to prove those equivalences I distinguish two cases\zl :\zd\

\vspace{1ex}

$1.$ $p$ true \zl and\zd\ $2.$ $p$ false

\vspace{1ex}

\noindent \zl The text on this page breaks here with the words: in both cases.\zd\

$\mathbf{\llbracket 51. \rrbracket}$ Now the distrib.\zl utive\zd\ laws in algebra make it possible to decide of any given expr.\zl ession\zd\ cont.\zl aining\zd\ only letters and $+\lfloor ,\rfloor -\lfloor ,\rfloor \cdot$ whether or not it is identically zero, namely by factori\zl z\zd ing out all prod.\zl ucts\zd\ of sums\zl ,\zd\ e.g\zl .\zd\ $x^2-y^2-(x+y)(x-y)=0$\zl .\zd\ A similar thing \zl is\zd\ to be exp.\zl ected\zd\ in \zl the\zd\ alg\zl ebra\zd\ of logic. Only $2$\zl two\zd\ differences\zl :\zd\ $1.$ In log\zl ic\zd\ we have the neg.\zl ation\zd\ which has no analogue in algebra. But for neg\zl ation\zd\ we have also a kind of distr\zl ibutive\zd\ law given by the De Morgan form.\zl ulas\zd\ $\sim (p \vee q) \equiv \: \sim p \:. \sim q$ $\mathbf{\llbracket 52. \rrbracket}$ \zl and\zd\ $\sim (p\:.\:q)\equiv\: \sim p \:\vee \sim q$\zl .\zd\ (Proved very easily by \zl the\zd\ truth\zl -\zd table method.) These formula\zl s\zd\ allow us to get rid of the neg.\zl ations\zd\ by shifting them \ul inwards \ud to the letters occurring in the expr\zl ession\zd . The \underline{sec.\zl ond\zd } difference is that we have two distr.\zl ibutive\zd\ laws and therefore two possible ways of factorizing. If we use the first law we shall get \ul as the final result \ud a sum of products \ul of single letters \ud as in algebra. By using the other law of distr.\zl ibution\zd\ we get a product of sums unlike in algebra. I think it is best to explain that on an $\mathbf{\llbracket 53. \rrbracket}$ example\zl :\zd

\vspace{1ex}

\noindent \zl The formula $(p \supset q)\:.\:p \supset q$, written in the manuscript at the top of the page, above the page number \textbf{53}., appears also in $\times$~4. after the examples done.\zd\

\begin{tabbing}
\hspace{1.7em}\=$\times$ 1. \hspace{.1em}\=$(p \supset q)\supset (\sim q \supset \:\sim p)$\hspace{10em}\=\\[.5ex]
\>\>$\sim (\sim p \vee q)\vee (q \:\vee \sim p)$ \\[.5ex]
\>\> $(p\:. \sim q)\vee q \:\vee \sim p$ \hspace{6em}\> disj.\zl unctive\zd\ \\[.5ex]
\>\> $(p \vee q \vee \sim p)\:.\:(\sim q \vee q \vee \sim p)$ \> conj.\zl unctive\zd\ \\[2ex]

\>$\times$ 2. \>$(p \supset q)\:.\:(p \supset \:\sim q)\:.\:p$\\[.5ex]
\>\>$(\sim p \vee q)\:.\:(\sim p \:\vee \sim q)\:.\:p$ \> conj.\zl unctive\zd\ \\[.5ex]
\>\> $(\sim p\:. \sim p \vee q\:. \sim p \:\vee \sim p\:. \sim q \vee q\:. \sim q)\:.\: p$ \\[.5ex]
\>\> $ (\sim p \: .\: p)\vee (q \:. \sim p\:.\: p)\vee (\sim p\:. \sim q \:.\: p)\vee (q\:. \sim q\:.\: p)$\zl full stop\\
\` deleted\zd\ \zl disjunctive\zd\ \\

\>\;\;\; 3. \>$\lfloor( \rfloor p \supset q)\supset (r \vee p \supset r \vee q)$\\[.5ex]
\>\>$\sim (\sim p \vee q) \vee [\sim(r \vee p)\vee r \vee q]$\\[.5ex]
\>\> $(p\:. \sim q)\vee (\sim r \:. \sim p) \vee r \vee q$ \> disj\zl unctive\zd\ \\[.5ex]
\>\> $(p \:\vee \sim r \vee r \vee q)\:.\:(p \:\vee \sim p \vee r \vee q)\lfloor .\rfloor$ \> conj.\zl unctive\zd\ \\[.5ex]
\>\> $(\sim q \:\vee \sim r \vee r \vee q)\:.\:(\sim q \:\vee \sim p \vee r \vee q)$
\end{tabbing}
\zl The line ``$\times$~4. $(p \supset q)\:.\:p \supset q$'' inserted at the end, which seems to be the beginning of an example not done, is deleted.\zd\

\vspace{1ex}

$\mathbf{\llbracket 1. \rrbracket}$ \zl Here the numbering of pages in this notebook starts anew.\zd\ In the last two lectures a proof for the completeness of our system of axioms for the calc.\zl ulus\zd\ of prop.\zl ositions\zd\ \ul was given\zl ,\zd\ \ud i.e.\ it was \ul shown \ud that any tautology is demonstrable from these axioms. Now a tautology is exactly what in trad.\zl itional\zd\ logic would be called a law of logic or a logically true prop\zl osition\zd . $\mathbf{\llbracket 2. \rrbracket}$ Therefore this completeness proof solves \ul for the calc.\zl ulus\zd\ of prop\zl ositions\zd\ \ud the second of the two problems which I announced in the beginning of my lectures\zl ,\zd\ namely it shows how all laws of a certain part of logic \ul namely \ud of the calc\zl ulus\zd\ of prop\zl ositions\zd\ can be deduced from a finite nu\zl mber\zd\ of logical axioms and rules of inference.

\zl new paragraph\zd\ \sout{And} I wish to stress that the interest of this result does not ly\zl ie\zd\ \sout{in this} so much in this that our particular four ax.\zl ioms\zd\ and three rules and four ax.\zl ioms\zd \zl repeated phrase ``four axioms''\zd\ are sufficient to deduce everything\zl ,\zd\ $\mathbf{\llbracket 3. \rrbracket}$ but the real interest consists in this that here for the first time in \zl the\zd\ history of logic it has really been \underline{proved} that one \underline{can} reduce all laws of a cert.\zl ain\zd\ part of logic to a few logical axioms\zl .\zd\ You know it has often been claimed that this can be done and sometimes the law\zl s\zd\ of id\zl entity\zd , contr.\zl adiction\zd , excl.\zl uded\zd\ middle have been considered as the log.\zl ical\zd\ axioms. But not even the shadow of a proof was given \ul that every logical inference can be derived from them \ud . Moreover the assertion to be proved was not even clearly formulated, because $\mathbf{\llbracket 4. \rrbracket}$ it means nothing to say that something \sout{prop.\zl erty\zd }\, can be derived e\zl .\zd g\zl .\zd\ from the law of contradiction unless you \sout{formulate} \ul specify \ud in addition the rules of inference which are to be used in the derivation.

As I said before it is not so very important that just our four ax.\zl ioms\zd\ are sufficient. After the method has once been developed, it is possible to give many other sets of axioms which are also sufficient to derive all (\sout{logically true prop.\zl ositions\zd }\, \ul tautologies \ud ) of the calc.\zl ulus\zd\ $\mathbf{\llbracket 5. \rrbracket}$ of prop\zl ositions\zd , \ul e.g.\

\begin{tabbing}
\hspace{1.7em}\=$p \supset (\sim p \supset q)$\\[0.5ex]
\>$(\sim p \supset p)\supset p$\\[0.5ex]
\>$(p \supset q)\supset [(q \supset r)\supset(p \supset r)]$ \ud	
\end{tabbing}

I have chosen the above four axioms because they are used in the standard textbooks of logistics\zl .\zd\ But I do not \ul at all \ud want to say that this choice was particularly fortunate. On the contrary our system of axioms is open to \sout{many} \ul some \ud objections from the aesthetic point of view\zl ;\zd\ e.g.\ one of the aesthetic requirements for a \sout{good} set of axioms \ul is that \ud the axioms should be as simple \ul and evident \ud as possible\zl ,\zd\ in any case simpler than the theor\zl ems\zd\ to be proved, whereas in our system $\mathbf{\llbracket 6. \rrbracket}$ e.g.\ the last axiom is pretty complicated and on the \ul other hand \ud the very simple law of identity $p \supset p$ appears as a theorem\zl .\zd\ \ul So in our system it happens sometimes that simpler propositions are \sout{based} proved \zl from\zd\ \sout{on} more complicated \sout{ones} \ul axioms\zl ,\zd\ \ud which is to be avoided if possible. \ud Recently by the \sout{Gentzen} mathematician G.\ Gentzen a system was set up which avoids these disadvantages. \zl The sentence broken here starting with ``I want to refer\zl ence\zd\ briefly about this system \sout{but wish to remark first that what I can}'' is continued on p.\ \textbf{7}.\ of Notebook IV.\zd\

\zl At the end of the present notebook there are in the manuscript thirteen not numbered pages with formulae and jottings. These pages are numbered here with the prefix $\mathbf{new\, page}$. It seems \textbf{new page i-iii} have been filled up backwards.\zd\
\begin{tabbing}
\noindent$\mathbf{\llbracket new\, page\; i\rrbracket}$\\*[1ex]
41 \hspace{3em}\=$(x)\varphi (x)\equiv (x\lfloor ,\rfloor y)\varphi (x)\:.\:\varphi (y)$\\*[.5ex]
?$\times$! 42\>$\varphi(x) \:.\: \psi(xy)\supset_{xy} \chi(xy)\equiv \varphi(x) \supset_x [\psi(xy)\supset_y \chi(xy)]$\\[.5ex]
43 \> $(\exists z)\varphi(z)\:.\: (\exists v)\chi(v) \supset [\varphi(z)\supset_z \psi(z)\:.\: \chi(v)\supset_v \vartheta(v) \equiv$ \\
\`$\varphi(z)\:.\: \chi(v)\supset_{zv} \psi(z)\:.\: \vartheta(v)]$
\end{tabbing}

\begin{picture}(70,10)(10,30)
\put(45,45){\makebox(0,0)[b]{\scriptsize $a$}} \put(65,45){\makebox(0,0)[b]{\scriptsize $e$}}
\put(45,25){\makebox(0,0)[b]{\scriptsize $i$}} \put(65,25){\makebox(0,0)[b]{\scriptsize $o$}}

\put(50,41.5){\line(1,-1){10}}
\put(60,41.5){\line(-1,-1){10}}
\put(45,41){\vector(0,-1){10}} \put(65,41){\vector(0,-1){10}}

\put(43,37){\line(1,0){4}}
\put(63,37){\line(1,0){4}}
\put(51,46.5){\line(1,0){8}}
\put(55,44.5){\line(0,1){4}}
\end{picture}

\begin{tabbing}
\noindent$\mathbf{\llbracket new\, page\; ii\rrbracket}$\\*[1ex]
$32'$\hspace{3em}\=$\sim (\exists x)\varphi(x)\supset \varphi (x)\supset_x \psi(x)$\\*[.5ex]
$32''$\>$(x)\psi(x)\supset \varphi(x)\supset_x \psi(x)$\\[.5ex]
$\times$ 32. \> $\sim [\varphi(x)\supset_x \psi(x)]\equiv (\exists x)[\varphi(x)\:. \sim \psi(x)]$ \\[.5ex]
33.\>$\varphi(x)\supset_x \psi(x)\supset (\exists x)\varphi(x)\:.\: \chi(x) \supset (\exists x) \varphi(x)\:.\: \chi(x)$\\[.5ex]
!34\>$\varphi(x) \supset_x \varphi(x) \vee \chi(x) \supset \varphi(x)\supset_x \psi(x) \vee (\exists x)\varphi(x)\:.\: \chi(x)$\\[.5ex]
35\> $\varphi(x) \supset_x (p \supset \psi(x))\equiv p \supset (\varphi(x)\supset_x \psi(x))$\\[.5ex]
36.\> $(x\lfloor ,\rfloor y)\varphi(xy)\equiv (yx)\varphi(xy)$\\[.5ex]
\> $(\exists x\lfloor ,\rfloor y)\varphi (xy) \equiv (\exists yx)\varphi(xy)$\\[.5ex]
$\times$ 37\> $\sim (x)(\exists y)\varphi(xy) \equiv (\exists x)(y)\sim \varphi(xy)$\\[.5ex]
38.\> $(x\lfloor ,\rfloor y)\varphi(x)\vee \psi(y) \equiv (x)\varphi(x) \vee (y)\psi(y)$\\[.5ex]
\> $(x\lfloor ,\rfloor y)\varphi(x)\supset \psi(y) \equiv (\exists x)\varphi(x) \supset (y)\psi(y)$\\[.5ex]
39\> $(\exists x\lfloor ,\rfloor y)\varphi(x)\:.\: \psi(y) \equiv (\exists x) \varphi(x)\:.\: (\exists y)\psi(y)$\\[.5ex]
40\> $(\exists x\lfloor ,\rfloor y)\varphi(x)\:.\: \psi(xy) \equiv (\exists x)[\varphi(x)\:.\:(\exists y)\psi(xy)]$
\end{tabbing}

\begin{tabbing}
\noindent$\mathbf{\llbracket new\, page\; iii\rrbracket}$\\*[1ex]
$\times$ 24.\hspace{3em}\=$(x)\sim \varphi(x)\supset \: \sim (x)\varphi(x)$\\*[.5ex]
25.\>$(z)[\varphi(z)\supset \psi(z)]\:.\: \varphi(x)\supset \psi(x)$\\[.5ex]
$\times$ 32. \> $\sim [\varphi(x)\supset_x \psi(x)]\equiv (\exists x)[\varphi(x)\:. \sim \psi(x)]$ \\[.5ex]
26.\>$(x)[\varphi(x)\supset \psi(x)]\:.\: (x)[\varphi(x)\supset \chi(x)]\supset (x)[\varphi(x)\supset \psi(x)\:.\: \chi(x)]$\\[.5ex]
$\times$ 27.\>$\varphi(x) \supset_x \psi(x)\:.\: \psi(x)\supset_x \chi(x) \supset \varphi(x)\supset \chi(x)$\\[.5ex]
28\> $\varphi(x) \supset_x \psi(x) \supset \varphi(x)\:.\: \chi(x)\supset_x \psi(x)\:.\: \chi(x)$\\[.5ex]
\> \sout{26}\;\; $27'$\zl ,\zd\ $28'$ analog f\" ur\zl German: analogous for\zd\ $\equiv$ \\[.5ex]
29\> $\varphi(x)\supset_x \psi(x) \:.\: \chi(x)\supset_x \vartheta(x) \supset \varphi(x)\:.\: \chi(x) \supset_x \psi(x) \:.\: \vartheta(x)$\\[.5ex]
\> $29'$ for \zl A\zd equiv\zl alence\zd \\[.5ex]
\> 30 \quad $30'$ \quad 26 \qquad $\vee$ \zl unreadable symbols\zd\ \\[.5ex]
31.\> $\varphi(x) \supset_x \psi(x) \:.\: \chi(x) \supset_x \psi(x) \supset \varphi(x) \vee \chi(x) \supset_x \psi(x)$
\end{tabbing}

\noindent $\mathbf{\llbracket new\, page\; iv\rrbracket}$ \zl The beginning of this page is in shorthand in the manuscript except for the following:
\begin{tabbing}
\hspace{5em}$p \:. \sim q$ \quad $p \circ q$ \quad $p \:. \sim q$\\
\hspace{17em} $F\:. \sim q$\\*[.5ex]
$\times$ \sout{1.}\hspace{3em}\=\sout{$\sim$ \zl ``In'' or $\mid$\zd\ $\vee \; . \; \equiv \; \supset$}\\[.5ex]
$\times$ 2.\> $\sim p\:. \sim q$ \\[.5ex]
$\times$ 3.\> $0 \supset$ \zl unreadable text\zd\ \quad $p \supset p+p \equiv \: \sim p$\zd\
\end{tabbing}

\noindent \zl The following column of formulae is crossed out in the manuscript:\zd
\begin{tabbing}
\hspace{3em}\=$[p \supset (p \supset q)]\:.\: p$ \=$\supset \lfloor (\rfloor p\supset q)$\\[.5ex]
\>\> $\supset p$\\[.5ex]
\> $\text{Vor} \supset (p \supset q)\:.\: p)$\\[.5ex]
\> $\lfloor (\rfloor p \supset q)\:.\: p$\\[.5ex]
\> $\text{Vor} \supset q$\\[2ex]
\> $[p \supset (p \supset q)]\supset (p \supset q)$\\[.5ex]
\> \zl unreadable formula\zd \\[.5ex]
\> $\lfloor (\rfloor r \supset p) \supset (\sim q \supset p \supset (q \supset$ \zl The formula breaks at this point.\zd \\[.5ex]
\> $(q \supset p) \:.\: (\sim q \supset p) \supset p$
\end{tabbing}

\noindent \zl Here the crossed out column of formulae ends, and the following column of formulae, which is not crossed out, is in the manuscript on the right of it on the same page:\zd

\begin{tabbing}
\hspace{3em}\=\underline{$(p \equiv q)\vee (p \equiv r) \vee (q \equiv r)$}\\[.5ex]
\underline{1.!}\> \underline{$(p \:.\ q \supset r)\:{\supset \atop \equiv}\:(p \:. \sim r \supset\: \sim q)$}\\[.5ex]
\>\underline{$\sim p \supset (\sim q \supset \: \sim (p \vee q))$}\\[.5ex]
\>\underline{$\sim (p \supset q) \equiv p\:. \sim q$}\\[.5ex]
\sout{2.}\> \sout{Red.\zl uctio ad\zd\ abs.\zl urdum\zd\ $(\sim p \supset p)\supset p$}\\[.5ex]
$^3$!\> \underline{$(\sim p \supset p) \supset p$}
\end{tabbing}

\noindent$\mathbf{\llbracket new\, page\; v\rrbracket}$\hspace{3em}$p\; \vee \sim p$

\vspace{2ex}

$
\hspace{1.5em}
\begin{cases}
1. \vee \;\; \lfloor \text{unreadable text}\rfloor\\
2. \; A \supset B \;\;\; B \supset A \;\;\; p\supset p \vee q\\
3. \; A \equiv B \;\;\; B \equiv A
\end{cases}
$

\vspace{1ex}

$
 \begin{array}{r|l}
 \text{\zl unreadable symbol\zd }\, \times \,2' & \text{Dualit\"at \zl German: duality\zd }\\
 1' & \lfloor \text{unreadable symbols with} \equiv \rfloor
 \end{array}
$

\begin{tabbing}
\` \underline{sec.\zl ond\zd\ law of distr\zl ibution\zd }
\end{tabbing}

$
1.
\begin{cases}
\lfloor \text{unreadable symbols with}\, \vee \,\text{and}\, \subset \rfloor \;\; \times ! \quad (p\supset q) \supset q \equiv p \vee q \quad \times \\
\lfloor \text{\underline{unreadable text}}\rfloor \;\; \lfloor \text{unreadable formula}\rfloor \quad \times
\end{cases}
$

\vspace{1ex}

\zl The following two columns of formulae are separated by a sinuous vertical line in the manuscript:\zd

\vspace{-2ex}

\begin{tabbing}
\hspace{22.5em} $p\:\equiv$\\
\hspace{2em}\= ! $\boxed{q\supset [(q \supset p) \,{\supset \atop \equiv}\, p]} $ \quad $\times$ \hspace{5.5em}\= ! $\boxed{p\supset (p \:.\: q\equiv q)}$ $\times$ viell.\zl per- \\
\` haps ``vielleicht'', German: perhaps\zd\ \\[.5ex]
\> \underline{$(p \supset q \vee r)\equiv q \vee (p \supset r)$} assoc.\zl iativity\zd\ !$^{\times}$\\
 \>\> $\sim p \supset (p \vee q \equiv q)$ $\times$\\[.5ex]
\> T, F, $p,q,\sim p,\sim q, p\equiv q, \sim p \equiv q$ \> \sout{$\sim(p \supset q) \equiv p \:. \sim q$}\\[.5ex]
\> $\mid$\underline{$[p \equiv (p \equiv q)]\equiv q$}$\mid$ \> $p \equiv p \vee p\:.\: q$\\[.5ex]
\> $\mid$\underline{\zl unreadable formula\zd }$\mid$ \> $p \equiv p\:.\: (p \vee q)$\\[.5ex]
\> $\mid$\underline{$[p \supset (q \supset r)] \supset [(p \supset q)\supset (p \supset r)]$}$\mid$ \> $\times$ $[(p \supset q) \supset p]\equiv p$ $\times$\\[.5ex]
\> $[(\sim p \supset q)\supset q]\supset (p \supset q)$ \> \sout{$(\sim p \supset p)\equiv p$}\\[.5ex]
\> $\mid$\underline{$(\sim q \supset r)\supset [(q \supset p)\:.\:(r \supset p)\supset p]$}$\mid$ \> $\sim (\sim p \equiv p)$\\[.5ex]
\underline{4$^2$!}\> $\mid$\underline{$p \supset (p \supset q) \equiv [p \supset q]$}$\mid$ $\times$ \> $\mid$\underline{$(q \supset p) \vee (r \supset p) \equiv (q\:.\:r \supset p)$}$\mid$\\[.5ex]
\>\> $\mid$\underline{$(q \supset p)\:.\:(r \supset p) \equiv (q \vee r \supset p)$}$\mid$ $\times$\\[.5ex]
\>\> \sout{\zl unreadable formula\zd }\\[.5ex]
\>\> ! $\boxed{p \supset q .\equiv. \: p\lfloor .\rfloor q \equiv p}$ $\times$\\*[.5ex]
\>\> \qquad \quad ~$.\equiv. \: q \equiv p \vee q$ $\times$\\[.5ex]
\>\> \sout{\zl unreadable formula\zd }
\end{tabbing}

\noindent$\mathbf{\llbracket new\, page\; vi\rrbracket}$ \zl This page is in shorthand in the manuscript except for the following: perspicuous, implicans, shorten the proof?, degenerated, formida\-ble, internal, \underline{manage}, \underline{I claim}, \underline{$(p_1 \vee p_2)\ldots \vee p_n$}, prove with their help, designated role, Moore, \zl unreadable word\zd , schlechthin\zl German: absolutely\zd \zd

\begin{tabbing}
\noindent$\mathbf{\llbracket new\, page\; vii\rrbracket}$\\*[1ex]
$\times$ 0.\hspace{2em}\=$(x)\varphi(x) \vee (\exists x) \sim \varphi(x)$ \qquad \= $\sim (\exists x)[\varphi(x)\:. \sim \varphi(x)]$\\*[.5ex]
$\times$ 1.\> $\sim (x) \equiv (\exists x) \sim$ \> $\times 1\cdot 1$\;\; $\sim (\exists x) \equiv (x) \sim$ \\[.5ex]
$\times$ 2. \> Verschieb.\zl perhaps ``verschieben'', German: move or postpone\zd\ \\
\` $(x)\;(\exists x)$ \zl unreadable text with $\vee\: .\:$\zd \\[.5ex]
?$\times^{?}_{\times}$ 3.\>$(x)[\varphi(x)\supset p] \equiv (\exists x)\varphi(x) \supset p$\\[.5ex]
\zl ?\zd $\times^{?}_{\times}$ 4.\>\;$(\exists x)[\varphi(x)\supset p] \equiv (x)\varphi(x)\supset p$\\[.5ex]
$\times$ 5.\> $(x)[p\supset \varphi(x)]\equiv p \supset (x)\varphi(x)$ \quad \zl unreadable text with $\exists$ perhaps\zd\ \\[.5ex]
6. \> $(x)[\varphi(x)\equiv p]\equiv (p\equiv(x)\varphi(x))\equiv$ \\
\` $[p\equiv(\exists x)\varphi(x)]\supset (x)\varphi(x)\equiv (\exists x)\varphi(x)$ \\[.5ex]
$\times$ 7.\> $(x)\lfloor [\rfloor \varphi(x)\:.\:\psi(x)\lfloor ]\rfloor\equiv(x)\varphi(x)\:.\:(x)\psi(x)$\\[.5ex]
$\times$ 8.\> $(\exists x)[\varphi(x)\vee \psi(x)]\equiv(\exists x)\varphi(x)\vee (\exists x)\psi(x)$\\[.5ex]
\sout{$\times$}$\times$ 9.\> $(x)\varphi(x)\vee (x)\psi(x)\supset(x)\lfloor [\rfloor\varphi(x)\vee \psi(x)\lfloor ]\rfloor$\\[.5ex]
\zl perhaps ``share'' and ``sicher''\zl German: sure\zd \zd\ 9.\\[.5ex]
$\times \times$ 10.\> $(\exists x)\lfloor [\rfloor\varphi(x)\:.\: \psi(x)\lfloor ]\rfloor\supset (\exists x)\varphi(x)\:.\:(\exists x)\psi(x)$\\[.5ex]
$\times$ 11.\> $(x)[\varphi(x)\supset \psi(x)]$ \= $\supset (x)\varphi(x)\supset (x)\psi(x)$\\[.5ex]
12.\>\> $\supset (\exists x)\varphi(x)\supset (\exists x)\psi(x)$ \\[1ex]
$\mathbf{\llbracket new\, page\; viii\rrbracket}$\\*[1ex]
13. \> $(x)[\varphi(x)\equiv \psi(x)]\supset$ \= $(x)\varphi(x)\equiv (x)\psi(x)$ \\[.5ex]
14. \> \zl unreadable word beginning perhaps with ``eben''; ``ebenfalls'' is\\
\` German for ``also''\zd\ \quad $\exists$ \\[.5ex]
15 \> $(\exists x)[\varphi(x)\supset \psi(x)]\equiv (x)\varphi(x)\supset (\exists x)\psi(x)$\\[.5ex]
16.\> Vert.\zl perhaps ``Vertauschung'', German: exchange\zd\ in der Reihen\\ \` \zl perhaps ``Reihenfolge'',
German: in the order\zd\ \\
\` $(\exists x)$\zl $\varphi(x)\supset$\zd $(x)$\zl $\psi(x)$\zd\ \\
17\> $(x\lfloor ,\rfloor y)\varphi(xy)\supset (x)\varphi(xx)$\\[.5ex]
18\> $(\exists x)\varphi(xx)\supset (\exists x\lfloor ,\rfloor y)\varphi(xy)$\\[.5ex]
$\times$ 19\> $(x)\varphi(x)\:.\:(\exists x)\psi(x) \supset (\exists x)\lfloor [\rfloor \varphi(x)\:.\: \psi(x)\lfloor ]\rfloor$
\zl An arrow points from\\
\` this formula to: Umkehrung\zl German: reversal\zd\ v.\ 10.\zd \\[.5ex]
$18'$\> $(x)[\varphi(xx)\supset (\exists u\lfloor ,\rfloor v)\varphi(uv)]$\\[.5ex]
$17'$\> $(x)[(u\lfloor ,\rfloor v)\varphi(uv)\supset \varphi(xx)]$\\[.5ex]
$\times$!? 20.!\> $(x)[\varphi(x) \vee \psi(x)]\:.\: (x)\sim \varphi(x) \supset (x)\psi (x)$\\[.5ex]
$\times$ 21\> $(x)\lfloor [\rfloor\varphi(x) \vee \sim \varphi (x)\lfloor ]\rfloor$\\[.5ex]
22\> $(x)\varphi(x) \supset (\exists x)\varphi(x)$\\[.5ex]
$\times$ 23\; \ul not inverse \ud $(\exists x)(y)\varphi(xy)\supset (y)(\exists x)\varphi(xy)$
\end{tabbing}

\noindent$ \mathbf{\llbracket new\, page\; ix\rrbracket}$ \zl This page is in shorthand in the manuscript except for the follow\-ing:\zd\

\vspace{1ex}

\noindent 1. Tauto1.\zl ogy\zd\ ? \\
2. Tauto1.\zl ogy\zd\ Taut.\zl ology\zd\ \\
4. Theorie \zl \&\zd\ Df.\\
5. demonstrable\\
6. impl.\zl ication\zd\ \\
primit\zl ive\zd\ rules of inf\zl erence\zd \\
7. fundamental conj.\zl unction\zd \\
Syll\zl ogism\zd \\
$p \supset p$\\
$\sim p \vee p$\\
$p \supset \:\sim \sim p$

\vspace{1.5ex}

\noindent 1. Wajsb.\zl erg\zd \\
2. Post Sep.\zl aratabdruck, German: offprint\zd\ Am\zl erican\zd\ Jour\zl nal\zd\ 43 \zl Emil Post's paper ``Introduction to a general theory of elementary propositions'', with his completeness proof for the propositional calculus is in the \emph{American Journal of Mathematics} vol.\ 43 (1921), pp.\ 163-185\zd \\
3. Zentralbl.\zl att f\" ur Mathematik und ihre Grenzgebiete\zd \\
5. \zl unreadable symbol\zd\ Father O\zl '\zd Hara \zl President of the University of Notre Dame from 1934 until 1939\zd\

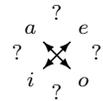
\begin{picture}(70,50)(0,15)

\put(45,45){\makebox(0,0)[b]{\scriptsize $a$}} \put(65,45){\makebox(0,0)[b]{\scriptsize $e$}}
\put(45,25){\makebox(0,0)[b]{\scriptsize $i$}} \put(65,25){\makebox(0,0)[b]{\scriptsize $o$}}

\put(50,41.5){\vector(1,-1){10}}
\put(60,41.5){\vector(-1,-1){10}}

\put(50,31.5){\vector(1,1){10}}
\put(60,31.5){\vector(-1,1){10}}

\put(40,35){\makebox(0,0)[b]{\scriptsize $?$}} \put(55,50){\makebox(0,0)[b]{\scriptsize $?$}}
\put(55,20){\makebox(0,0)[b]{\scriptsize $?$}} \put(70,35){\makebox(0,0)[b]{\scriptsize $?$}}

\end{picture}

\noindent \zl On the right of this picture one finds a question mark, the symbol $\leq$ rotated counter-clockwise for approximately 45 degrees and an unreadable symbol.\zd\

\vspace{1ex}

\noindent$\mathbf{\llbracket new\, page\; x\rrbracket}$ \zl This and the following three pages in the manuscript, \textbf{new page x-xiii}, are loose, not bound to the notebook with a spiral and without holes for the spiral. In all of the notebooks the only other loose leafs are to be found towards the end of Notebook~V and at the end of Notebook VII. In the upper half of the present page in the manuscript one finds the following, turned counter-clockwise for 90 degrees and crossed out:
\begin{tabbing}
\hspace{1.7em}\=$\mathfrak{A}\rightarrow \mathfrak{B}$ \hspace{2em} \=$\rightarrow\; \mathfrak{A}\supset \mathfrak{B}$\\[.5ex]
\>$\mathfrak{B}\rightarrow C$\\[.5ex]
\>$\mathfrak{A}\rightarrow C$\\[2ex]
\>$\sim (A\:. \sim B)$ \hspace{2em} $A\rightarrow B$\\[.5ex]
\>$A\:. \sim B \rightarrow A$\\[.5ex]
\>$A\:. \sim B \rightarrow \: \sim B$\\[2ex]
\>$A,B\:\&\: C$\\[.5ex]
\>$A,B, C \rightarrow A,B\:\&\: C$\\[2ex]
\>$A\:. \sim A$ \>$\sim B, \mathfrak{A}, A \rightarrow B$\\[.5ex]
\>\>$\sim B, \mathfrak{A} \rightarrow \: \sim A$\\[.5ex]
\>$\mathfrak{A}, A \rightarrow B$\\[.5ex]
\>$\mathfrak{A}\lfloor ,\rfloor \sim B \rightarrow \: \sim A$\\[.5ex]
\>$\mathfrak{A}\lfloor ,\rfloor \sim B\lfloor ,\rfloor A \rightarrow B$\\[.5ex]
\>$\mathfrak{A}\lfloor ,\rfloor \sim B\lfloor ,\rfloor A \rightarrow \: \sim B$
\end{tabbing}
\zl The lower half of this of the page is in shorthand in the manuscript except for the following:\zd\

\vspace{1ex}

\noindent 1. 1$\cdot$2\quad P, G, lie, = \\
3. 3$\cdot$1, 3$\cdot$2\quad $<$, Z \\
4. \zl two ditto marks referring to ``$<$, Z'' followed by =\zd\

\begin{tabbing}
\hspace{1.7em}\=$\mathfrak{A}\rightarrow \mathfrak{B}$ \hspace{2em} \=$\rightarrow\; \mathfrak{A}\supset \mathfrak{B}$\kill

$\mathbf{\llbracket new\, page\; xi\rrbracket}$\\*[1ex]
\> $\sim R, \mathfrak{A},p$ \= $\rightarrow R$\\*[.5ex]
\>\>$\rightarrow \: \sim R$\\[.5ex]
\>$\sim R, \mathfrak{A} \rightarrow \: \sim p$
\end{tabbing}

\noindent$\mathbf{\llbracket new\, page\; xii\rrbracket}$ \zl The following list of formulae is crossed out in the manu\-script:\zd\
\begin{tabbing}
\hspace{1.7em}\=\underline{$(p \supset q)\:.\:(r \supset q)\supset (p\vee r \supset q)$}\\[.5ex]
\> $\sim [(\sim p \vee q)\:.\:(\sim r \vee q)]\vee [\sim (p \vee r) \vee q]$\\[.5ex]
\> $\sim (\sim p \vee q) \vee \sim (\sim r \vee q)\vee (\sim p \:. \sim r)\vee q$\\[.5ex]
\> \underline{$(p\:. \sim q)\vee (r\:. \sim q)\vee (\sim p\:. \sim r)\vee q$}\\[.5ex]
\>$(p \vee r\: \vee \sim p \vee q)$.\\[.5ex]
\>$($\underline{$p$}\:$\vee$\:\underline{$r$}\:$\vee \sim \:$\underline{$r$}\:$\vee\: q)$.\\[.5ex]
\> \underline{$(p\: \vee \sim q\: \vee \sim p \vee q)$}\\[.5ex]
\>$(p\: \vee \sim q\: \vee \sim r \vee q)$\\[.5ex]
\>$(\sim q \vee r\: \vee \sim p \vee q)$\\[.5ex]
\>$(\sim q \vee r\: \vee \sim r \vee q)$\\[.5ex]
\>$(\sim q\: \vee \sim q\: \vee \sim p \vee q)$\\[.5ex]
\>$(\sim q\: \vee \sim q\: \vee \sim r \vee q)$\\[2ex]
\> $\lfloor (\rfloor p \supset \sim p)\supset\: \sim p$\\[.5ex]
\> $\sim (\sim p \lfloor \vee \:\sim\rfloor p)\: \vee \sim p$\\[2ex]
\> $[ \sim (\sim p \vee q) \vee (p \:. \sim q)]$\\[.5ex]
\> $[ \sim (p \:. \sim q)\: \vee \sim p \vee q]$\\[.5ex]
\> $[ (p \:. \sim q) \vee (p \:. \sim q)]$\\[.5ex]
\>$[\sim p \vee q\: \vee \sim p \vee q]$\\[.5ex]
\>\underline{$\sim p \:.\:p\:. \sim q \vee q\:.\: p \:. \sim q$}
\end{tabbing}

\noindent \zl Here the list of formulae crossed out in the manuscript ends, and the following not crossed out list is given:\zd
\begin{tabbing}
\hspace{1.7em}\=\underline{$p \supset q .\supset . \sim q \supset\: \sim p$}\\[.5ex]
\> $\sim (\sim p \vee q) \vee (q\: \vee \sim p)$\\[.5ex]
\>$(p \:. \sim q)\vee q\: \vee \sim p$\\[.5ex]
\>$(p \vee q\: \vee \sim p)\:.\:(\sim q \vee q\: \vee \sim p)$\\[2ex]
$\times$\>\underline{$(p \supset q)\supset (r \vee p \supset r \vee q)$}\\[.5ex]
\>$\sim (\sim p \vee q)\vee (\sim (r \vee p) \vee r \vee q)$\\[.5ex]
\>$(p \:. \sim r)\vee (\sim r \:. \sim p)\vee r \vee q$\\[.5ex]
\>$(p\: \vee \sim r)\:.\:(p\: \vee \sim p)\:.\:(\sim q\: \vee \sim r)\:.\:(\sim q\: \vee \sim p) \vee (r \vee q)$
\end{tabbing}

\noindent $\mathbf{\llbracket new\, page\; xiii\rrbracket}$ \zl In the left margin turned counter-clockwise for 90 degrees one finds first on this page of the manuscript:
\begin{tabbing}
\hspace{1.7em}\=\zl unreadable text with: Arist.\zl otelian\zd\ Syll.\zl ogisms\zd \zd\ \\[.5ex]
\>$\sim [a\lfloor\cdot\rfloor b=0\:.\:c\lfloor\cdot\rfloor\bar{b}=0\:.\:a\cdot c\neq 0]$
\end{tabbing}
Next one finds in the left half of the page a column of propositional formulae, partly effaced, partly crossed out and mostly unreadable, which is not given here. In the rest of the page one finds
\begin{tabbing}
\hspace{1.7em}$\aleph^{\aleph_1}_{\omega_1} \cdot \aleph^{\aleph_0}_{\omega_1 + \omega} \geq \aleph^{\aleph_1}_{\omega_1 + \omega}$
\end{tabbing}
followed by an unreadable inequality with $\aleph_2$. One finds also the following, turned counter-clockwise for 90 degrees:\zd\
\begin{tabbing}
\hspace{1.7em}\= $\mid$\underline{$\aleph^{\aleph_0}_2 > \aleph^{\aleph_0}_1$}$\mid$ \qquad \qquad \= \underline{$\aleph^{\aleph_0}_{\alpha +1} > \aleph^{\aleph_0}_1$}\\[.5ex]
\> $\aleph^{\aleph_0}_1 \geq \aleph_2$ \> $\aleph^{\aleph_0}_{\alpha} = \aleph^{\aleph_0}_1$\\[.5ex]
\> $\aleph^{\aleph_0}_1 = \aleph_{\alpha}$
\end{tabbing}
\begin{tabbing}
\hspace{1.7em}\= \kill
$\times$\>\underline{$\sim [(\sim p \vee q)\:.\: p]\vee q$}\\[.5ex]
\> $[\sim (\sim p \vee q)\: \vee \sim p] \vee q$\\[.5ex]
\>$(p \:. \sim q)\:\vee \sim p \vee q$ \qquad disj\zl unctive\zd \\[.5ex]
\> $(p\: \vee \sim p \vee q)\:\lfloor .\rfloor (\sim q\: \vee \sim p \vee q)$\qquad \zl conjunctive\zd\
\end{tabbing}
\zl and at the end the following, turned clockwise for 90 degrees:\zd\
\begin{tabbing}
\hspace{1.7em}\=$xR^Sy \equiv (z)\lfloor [\rfloor zSy \supset xRz\lfloor ]\rfloor$\\[.5ex]
\> $(R^S)^T=R^{(S\mid T)}$\\[.5ex]
\> $R^{S+T}=R^S\lfloor \cdot \rfloor R^T$ \\[.5ex]
\hspace{2.5em}\= $(z)\lfloor [\rfloor zSy \vee zTy \supset xRz\lfloor ]\rfloor$\\[.5ex]
\sout{$(\exists u)$}\> $(z)\lfloor [\rfloor zSy \supset xRz\lfloor ]\rfloor \; .\; (z)\lfloor [\rfloor zTy \supset xRz\lfloor ]\rfloor$
\end{tabbing}

\section{Notebook IV}\label{0IV}
\pagestyle{myheadings}\markboth{SOURCE TEXT}{NOTEBOOK IV}
\zl Folder 62, on the front cover of the notebook ``Log.\zl ik\zd\ Vorl.\zl esungen\zd\ \zl German: Logic Lectures\zd\ N.D.\ \zl Notre Dame\zd\ IV''\zd\

\vspace{1ex}

\zl Before p.\ {\bf 7}., the first numbered page in this notebook, there are in the manuscript four not numbered pages with formulae. These pages are numbered here with the prefix $\mathbf{new\, page}$. The formulae with $R$, $S$ and $T$ on $\mathbf{new\, page\; i}$ are in boxes on the right of this page.\zd\

\begin{tabbing}
$\mathbf{\llbracket new\, page\; i\rrbracket}$\\*[1ex]
\zl unreadable text\zd \quad\= $(p\supset q)\supset [(r \supset p)\supset(r \supset q)]$\quad\= (1)\zl 1.\zd \\*[.5ex]
\zl unreadable text\zd \> $p\supset \;\sim\sim p$ \>(2)\zl 2.\zd\ \\
\hspace{26em}$R$ $\,$..\zl :\zd\ $\sim p$\\[.5ex]
\hspace{26em}$S$ $\,$..\zl :\zd\ $\sim\sim\sim p$\\[.5ex]
\hspace{26em}$T$ $\,$..\zl :\zd\ $p$\\[.5ex]
Su (2)\zl 2.\zd\ \> $\sim p\supset \;\sim\sim\sim p$\` \zl one line below\zd\ $R\supset S$ \\[.5ex]
\zl unreadable text\zd \> $\sim p\vee p\supset \;\sim\sim\sim p\vee p$ \` \zl one line below implication\\ \` with unreadable left-hand side and $R\vee T$ or $R\vee S$ on\\
\` the right; the implication in this line is $R\vee T\supset S\vee T$\zd \\[.5ex]
Su\zl .\zd\ III\zl (3)\zd \> $\sim p\vee p\supset p\;\vee\sim p$ \>(3)\zl 3.\zd \` $R\vee T\supset T\vee R$ \\[.5ex]
Su IV\zl (4)\zd \> $\lfloor (\rfloor\sim p\supset\;\sim\sim\sim p \lfloor )\rfloor\supset [p\;\vee \sim p\supset p\;\vee\sim\sim\sim p]$\quad (4)\zl 4.\zd \\[.5ex]
Imp 2\zl .\zd , 4\zl .\zd \> $p\;\vee \sim p\supset p\;\vee\sim\sim\sim p$\>(5)\zl 5.\zd \` \zl one line above\zd \\
\` $T\vee R\supset T\vee S$ \\[.5ex]
Su III\zl (3)\zd \> $p\;\vee\sim\sim\sim p\supset\;\sim\sim\sim p\vee p$\>(6)\zl 6.\zd \` \zl one line above\zd \\ \` $T\vee S\supset S\vee T$ \\[.5ex]
Su (1)\zl 1.\zd \>$(p\;\vee\sim p\supset p\;\vee\sim\sim\sim p)\supset$\\ \` $[(\sim p\vee p \supset p\;\vee\sim p)\supset(\sim p\vee p \supset p\;\vee\sim\sim\sim p)]$\quad (7)\zl 7.\zd \\[.5ex]
Imp 2mal\zl zweimal, German: twice\zd\ 5\zl .\zd , 7 \zl . \sout{(unreadable word 3)}; 3.\zd \\
\> $\sim p\vee p \supset p\;\vee\sim\sim\sim p$\> (8)\zl 8.\zd \\[.5ex]
Su III\zl (3)\zd \> $p\;\vee\sim\sim\sim p\supset\;\sim\sim\sim p\vee\;p$\> (9)\zl 9.\zd \\ \` \zl occurs already as (6)\zd \\[.5ex]
Su (1)\zl 1.\zd \>$(p\;\vee\sim\sim\sim p\supset \;\sim\sim\sim p\vee p)\supset$\\ \` $[(\sim p\vee p \supset p\;\vee\sim\sim\sim p)\supset(\sim p\vee p \supset \;\sim\sim\sim p\vee p)]$\quad (10)\zl 10.\zd \\[.5ex]
Imp 2mal\zl zweimal, German: twice\zd\ 9\zl 6.\zd , 10\zl .\zd ; 8\zl .\zd \` $\sim p\vee p \supset \;\sim\sim\sim p\vee p$
\end{tabbing}

\begin{tabbing}
$\mathbf{\llbracket new\, page\; ii\rrbracket}$
\hspace{2,5em}\=\underline{$p\supset q\vee p$}\\*[.5ex]
\> $p\supset p\vee q$\hspace{3em} \=I\zl (1)\zd \\[.5ex]
\> \underline{$p\vee q\supset q\vee p$}\> III\zl (3)\zd \\[1ex]
\ul Su (1)\zl 1.\zd\ \ud \hspace{5em}\> $(p\vee q\supset q\vee p)\supset [(p \supset p\vee q)\supset(p \supset q\vee p)]$\hspace{1.3em}\=(2)\zd 2.\zd \\
\` Su\quad $\f{p\vee q}{p}$\quad $\f{q\vee p}{q}$\quad$\f{p}{r}$\\[.5ex]
Imp (2\zl .\zd , III\zl (3)\zd )\> $(p \supset p\vee q)\supset(p \supset q\vee p)$\> (3)\zl 3.\zd \\[-2ex]
\underline{\hspace{32em}}\\
Imp (3\zl .\zd , I\zl (1)\zd )\> $p \supset q\vee p$\> (4)\zl 4.\zd
\end{tabbing}

\begin{tabbing}
$\mathbf{\llbracket new\, page\; iii\rrbracket}$\\*[1ex]
\underline{1.} \hspace{1em}\=$(\sim p\supset p)\supset p$\hspace{3em}\= $(\sim\sim p\vee p)\supset p$\\*[1ex]
\> A. \hspace{.27em} $p\supset p$\\[.5ex]
\> \underline{$\sim\sim p\supset p$}\\[.5ex]
\> $\sim\sim p\vee p\supset p$\> \zl D\zd ilemma\\[1ex]
\underline{2.}\>$(p\: .\: q\supset r)\supset (p\: .\sim r\supset
\;\sim q)$\\[1ex]
1.\> $(p\: .\: q\supset r)\supset[p\supset(q\supset r)]$\hspace{4.5em}\=Exp.\zl ortation\zd \\[.5ex]
\>$(q\supset r)\supset(\sim r\supset \:\sim q)$\> Transpos.\zl ition\zd \\[.5ex]
2.\> $[p\supset(q\supset r)]\supset[p\supset(\sim r\supset \:\sim q)]$\> Add.\zl ition\zd\ from \zl the\zd\ left\\[.5ex]
3\zl .\zd \> $[p\supset(\sim r\supset \:\sim q)]\supset [p\: .\sim r\supset \:\sim q]$\> Imp\zl ortation\zd \\[-2ex]
\underline{\hspace{17em}}\\
\> $(p\: .\: q\supset r)\supset (p\: .\sim r\supset\;\sim q)$\> 1\zl .\zd , 2\zl .\zd , 3\zl .\zd\ Syll.\zl ogism\zd \\[1ex]
\underline{3.1}\> $(p\supset q)\supset(p\supset(p\supset q))$\\[.5ex]
\> $r\supset(p\supset r)$\hspace{9em}\=$\f{p\supset q}{r}$\\[.5ex]
\underline{3.2}\> $[p\supset(p\supset q)]\supset(p\supset q)$\>$\sim p\vee(\sim p\vee q)\supset \;\sim p\vee q$\\[2ex]

$\mathbf{\llbracket new\, page\; iv\rrbracket}$\\*[1ex]

1.\> $\sim p\vee(\sim p\vee q)\supset(\sim p\;\vee\sim p)\vee q$\\[.5ex]
\> $\sim p\;\vee\sim p\supset\; \sim p$\\[.5ex]
2.\> $(\sim p\;\vee\sim p)\vee q\supset\;\sim p \vee q$\>Add.\zl ition\zd\ from \zl the\zd\ right\\[.5ex]
\> $\sim p\vee(\sim p\vee q)\supset\;\sim p \vee q$\> Syll.\zl ogism\zd\ 1.\zl ,\zd\ 2.\\[.5ex]
\>$[p\supset(p\supset q)]\supset(p\supset q)$\>Rule of def.\zl ined\zd\ symb.\zl ol\zd
\end{tabbing}

\noindent $\mathbf{\llbracket 7. \rrbracket}$ \zl This page starts with the ending of the sentence started as follows at the end of p.\ \textbf{6}.\ towards the end of Notebook III: I want to refer\zl ence\zd\ briefly about this system\zd\ or to be more exact on a system which is based on Gentzen's idea, but simpler than his. The idea consists in \ul \ud introducing another kind of implication (denoted by an arrow $\rightarrow$). \zl The remainder of p.\ \textbf{7}. is crossed out in the manuscript.\zd\ $\lceil$such that $P\rightarrow Q$ means $Q$ is true under the assumption $P$. The diff.\zl erence\zd\ of this implication as opposed to our former one is

1. There can be any number of premis\zl \sout{s}\zd es\zl ,\zd\ e.g.\ $P,Q\rightarrow R$ means $R$ holds under the ass.\zl umptions\zd\ $P,Q$ (i.e\zl .\zd\ the same thing which would be \zl unreadable text, could be: den\zl oted\zd \zd\ by $P\: .\: Q\supset R$. \sout{e.g} In particular the number of premis\zl \sout{s}\zd es \zl Here p.\ \textbf{7}. ends and pp.\ \textbf{8}. and \textbf{9}. are missing, while p.\ \textbf{10}. begins with the second part of a broken sentence.\zd\

$\mathbf{\llbracket 10. \rrbracket}$ system with altogether three prim\zl itive\zd\ terms $\rightarrow$, $\sim$, $\supset$\zl .\zd\ We have now to distinguish between expressions in the former sense\zl , i\zl .\zd e.\ containing only $\sim$\zl ,\zd\ $\supset$ and var\zl iables,\zd\ e\zl .\zd g\zl .\zd\ $p\supset q$, $\sim p\supset q$, $q\supset p\vee r$, etc\zl .,\zd\ and sec.\zl ondary\zd\ formulas containing the arrow\zl ,\zd\ e\zl .\zd g.\ $p,p\supset q\rightarrow q$\zl .\zd\ I shall use capital Latin letters $P,Q$ only to denote expr\zl essions\zd\ of the first kind\zl ,\zd\ i\zl .\zd e\zl .\zd\ expressions in our former sense\zl ,\zd\ and I use cap.\zl ital\zd\ Greek letters \ul $\Delta,\Gamma$ \ud to denote sequences of an arb.\zl itrary\zd\ \sout{nu.} number of ass.\zl umptions\zd\ $\underbrace{P,Q,R\ldots}_\text{$\Delta$}$ \quad \sout{$\underbrace{p,p\supset q,\sim q}_\text{$\Delta$}$}\quad \sout{may be denoted by} $\mathbf{\llbracket 11. \rrbracket}$ \sout{$\Delta$}. \sout{So the cap\zl ital\zd\ Greek letters denote possible premis\zl \sout{s}\zd es to the formulas of the}

Hence a formula of G.\zl entzen's\zd\ system will \ul always \ud have the form $\Delta\rightarrow S$\zl ,\zd\ a cert.\zl ain\zd\ sequence of expr.\zl essions\zd\ of the first kind implies an expr.\zl ession\zd\ of the first kind. \sout{And} Now to the axioms and rules of inference.

\vspace{1ex}

I\quad Any form\zl ula\zd\ $P\rightarrow P$ where $P$ is an arb.\zl itrary\zd\ expr\zl ession\zd\ of the first kind is an ax.\zl iom\zd\ and only those form\zl ulas\zd\ are ax\zl ioms\zd . \sout{(So that is the law of identity)}

\vspace{1ex}

\noindent $\mathbf{\llbracket 12. \rrbracket}$ \sout{$P$ may be} So that is the law of identity which appears here as an axiom and as the only axiom.

\vspace{1ex}

As to the rules of inference we have 4\zl four,\zd\ namely

\zl crossed out: 1. $\Delta\rightarrow A$ \quad\quad $\f{\Delta, P\rightarrow A}{P,\Delta\rightarrow A}$\zd\

\vspace{1ex}

\noindent 1.\quad The rule of addition of premis\zl \sout{s}\zd es\zl ,\zd\ i.e.\ from $\Delta\rightarrow A$ one can conclude $\Delta,P \rightarrow A$ and $P,\Delta\rightarrow A$\zl ,\zd\ i\zl .\zd e.\ if $A$ is true under the assumptions $\Delta$ then it is \ul a fortiori \ud true under the assumptions $\Delta$ and the further ass.\zl umption\zd\ $P$\zl .\zd\

\vspace{2ex}

\noindent $\mathbf{\llbracket 13. \rrbracket}$
\begin{tabbing}
2.\quad \= The \zl R\zd ule of exportation:\\[1ex]
\> $\Delta,P\rightarrow Q$\hspace{2em}\= : \hspace{1.2em} $\Delta\rightarrow (P\supset Q)$
\end{tabbing}
If the prop\zl ositions\zd\ $\Delta$ and $P$ imply $Q$ then the prop\zl ositions\zd\ $\Delta$ imply that $P$ implies $Q$.
\begin{tabbing}
3.\quad The Rule of implication:\\[1ex]
\hspace{1.2em}\begin{tabular}{ l|l }
$\Delta\rightarrow P$ & \\[-1ex]
& \hspace{1.5em}$\Delta\rightarrow Q$ \\[-1ex]
$\Delta\rightarrow(P\supset Q)$ &
\end{tabular}
\end{tabbing}
So that is so to speak the rule of implication under some assumptions: If $A$ and $A\supset B$ both hold under the ass.\zl umptions\zd\ $\Delta$ then $B$ also holds under the ass.\zl umptions\zd\ $\Delta$\zl .\zd\
\begin{tabbing}
4.\quad Rule of Reductio ad abs\zl urdum\zd\ or \ul rule of \ud indirect proof\zl :\zd \\[1ex]
\hspace{1.2em}\begin{tabular}{ l|l }
$\Delta,\sim P\rightarrow Q$ & \\[-1ex]
& \hspace{1.2em}$\Delta\rightarrow P$ \\[-1ex]
$\Delta,\sim P\rightarrow\;\sim Q$ &
\end{tabular}
\end{tabbing}
Here the prem.\zl ises\zd\ mean that from the ass\zl umptions\zd\ $\Delta$ and $\sim P$ a contradiction follows\zl ,\zd\ i\zl .\zd e.\ $\sim P$ is incompatible $\mathbf{\llbracket 14. \rrbracket}$ with the ass.\zl umptions\zd\ $\Delta$\zl ,\zd\ i\zl .\zd e.\ from $\Delta$ follows $P$.

Again it can be proved that every tautology follows from the ax.\zl ioms\zd\ and rules of inf\zl erence\zd . Of course only the tautologies which can be expressed in terms of the symbols introd.\zl uced,\zd\ i\zl .\zd e.\ $\sim$\zl ,\zd\ $\supset$ \zl and\zd\ $\rightarrow$\zl .\zd\ If we want to introduce also $\vee$\zl ,\zd\ $.$ etc.\ we have to add the rule of the defined symbol $.$ or other rules concerning $\vee$\zl ,\zd\ $.$ etc.

\zl new paragraph\zd\ Now you see that in this system the aforementioned disadvantages have been avoided\zl .\zd\ All the axioms are really very simple and $\mathbf{\llbracket 15. \rrbracket}$ evident. It is particularly interesting that also the pseudo-paradoxical prop.\zl ositions\zd\ about the impl.\zl ication\zd\ follow from our system of axioms although nobody will have any objections against the axioms themselves\zl ,\zd\ i\zl .\zd e.\ everybody would admit them if we interpret both the $\rightarrow$ and the $\supset$ to mean \zl ``\zd if\ldots\ then\zl ''\zd . Perhaps I shall derive these two \sout{prop.} \ul pseudo\zl -\zd paradoxes \ud as \sout{an} examples for \sout{a} derivations from this system. The first reads:
\begin{tabbing}
By \= I\hspace{1.5em}\= $q\rightarrow q$\kill

\>\>$q\rightarrow p\supset q$\hspace{2em} Proof:\\[1ex]

$\mathbf{\llbracket 16. \rrbracket}$\\[1ex]

By \> I\> $q\rightarrow q$\\[.5ex]
$\; ''$\> 1\> $q,p\rightarrow q$\\[.5ex]
$\; ''$\> 2\> $q\rightarrow (p\supset q)$
\end{tabbing}
Incidentally\zl ,\zd\ again app.\zl lying\zd\ 2 we get $\rightarrow q\supset(p\supset q)$ which is another form for the same theorem. The sec.\zl ond\zd\ paradox reads like this:
\begin{tabbing}
By \= I\hspace{1.5em}\= $q\rightarrow q$\kill

\>\>$\sim p\rightarrow p\supset q$\hspace{2em} Proof\zl :\zd \\[1ex]
\> I\> $p\rightarrow p$\\[.5ex]
\> 1\> $\sim p,p,\sim q\rightarrow p$\\[.5ex]
\> I\> $\sim p\rightarrow \:\sim p$\\[.5ex]
\> 1\> $\sim p,p,\sim q\rightarrow \:\sim p$\\[.5ex]
\> 4\> $\sim p,p\rightarrow q$\\[.5ex]
\> 2\> $\sim p\rightarrow (p\supset q)$
\end{tabbing}

$\mathbf{\llbracket 17. \rrbracket}$ Incidentally this form\zl ula\zd\ $\sim p,p\rightarrow q$ which we derived as an intermediate step of the proof is interesting also on its own account\zl ;\zd\ it says: From a contrad.\zl ictory\zd\ assumption everything follows since the formula is true whatever the prop.\zl osition\zd\ $q$ may be. I am sorry I have no time left to go into more details about this Gent.\zl zen\zd\ system. I want to conclude now this chapter about the calc.\zl ulus\zd\ of prop\zl osition\zd . \zl Here p.\ \textbf{17}. ends and pp.\ \textbf{18}.-\textbf{23}. are missing.\zd\

$\mathbf{\llbracket 24. \rrbracket}$ I am concl.\zl uding\zd\ now the chapt.\zl er\zd\ about the calc.\zl ulus\zd\ of prop\zl osi\-tions\zd\ and begin with the next chapt.\zl er\zd\ which is to deal with the so called calc.\zl ulus\zd\ of functions \ul or predicates \ud . As I explained formerly the calc\zl ulus\zd\ of prop.\zl ositions\zd\ is c\zl h\zd aracteri\zl z\zd ed by this that only prop.\zl ositions\zd as a whole occur in it\zl .\zd\ \sout{You know} The letters $p,q,r$ etc\zl .\zd\ denoted arbitrary propositions and all the formulas and rules \ul which we proved \ud are valid whatever \sout{the} propos.\zl itions\zd\ $p,q,r$ may be\zl ,\zd\ i\zl .\zd e.\ they are independent of the structure of the prop.\zl ositions\zd\ involved. Therefore we could use \sout{a} single letters \ul $p,q\ldots$ \ud to denote \sout{a} whole propositions.

$\mathbf{\llbracket 25. \rrbracket}$ But now we shall be concerned with inferences which depend on the structure of the prop.\zl ositions\zd\ involved and therefore we shall have to study at first how prop.\zl ositions\zd\ are built up of their constituents. To this end we ask at first what do the simplest prop.\zl ositions\zd\ which one can imagine look like. Now \sout{\zl unreadable text\zd }\, evidently the simplest kind of prop.\zl ositions\zd\ are those in which simply some predicate is asserted of some subject\zl ,\zd\ e.g.\ Socrates is mortal\zl .\zd\ Here the predicate mortal is asserted to belong to the subject Socrates. Thus far we are in agree- $\mathbf{\llbracket 26. \rrbracket}$ ment with classical logic.

\zl new paragraph\zd\ But there is another type of simple prop.\ osition\zd\ which was very \ul much \ud neglected in classical logic, although this second type is \sout{even} more important for the applications of logic in mathem\zl atics\zd\ and other sciences\zl .\zd\ This second type \ul of simple prop.\zl osition\zd\ \ud consists in this that a predicate is asserted of several subjects\zl ,\zd\ e.g.\ New York is larger than Washington. \sout{or Socrates is the teacher of Plato} Here you have two subj.\zl ects,\zd\ New Y\zl ork\zd\ and W.\zl ashington,\zd\ and the pre\zl dicate\zd\ \sout{greater} larger says that a certain relation subsists between those two subj\zl ects\zd . Another ex.\zl ample is\zd \zl ``\zd Socrates is the teacher of Plato\zl ''\zd \zl ``again'' is superfluous after the first occurrence of this sentence having been crossed out above\zd . So you see there are two different kinds $\mathbf{\llbracket 27. \rrbracket}$ of predicates\zl ,\zd\ namely pred.\zl icates\zd\ with one subj\zl ect\zd\ as e.g.\ \underline{mortal} and predicates with several subj.\zl ects\zd\ as e.g.\ \underline{greater}.

\zl new paragraph\zd\ The pred.\zl icates\zd\ of the first kind may be called properties \sout{or qualities}, \sout{and} those of the sec.\zl ond\zd\ kind \ul are called \ud relations. So e.g.\
,,\zl ``\zd mortal'' is a property\zl ,\zd\ ,,\zl ``\zd greater'' is a relation. \zl M\zd ost of the pred.\zl icates\zd\ of everyday lang\zl uage\zd\ are relations and not properties\zl .\zd\ The relation ,,\zl ``\zd greater'' as you see requires two subjects and therefore is called a dyadic relation. There are also relations which require three or more subjects\zl ,\zd\ e.g.\ \underline{betweenness} is a relation with three subj.\zl ects,\zd\ i.e.\ triadic relation. If I say e.g.\ New York $\mathbf{\llbracket 28. \rrbracket}$ lies between Wash\zl ington\zd\ and Boston.\zl ,\zd\ \zl t\zd he relation of betweenness is asserted to subsist for the three subjects N.\zl ew York,\zd\ W\zl ashington\zd\ and B.\zl oston,\zd\ and always if I form a meaningful prop.\ osition\zd\ involving the word between I must mention three objects of which one is to be in between the others. \sout{So} \ul Theref\zl ore\zd\ \ud ,,\zl ``\zd betweenness'' is \ul called \ud a triadic rel.\zl ation\zd\ and similarly there are tetradic, pentadic rel.\zl ations\zd\ etc. Properties may \sout{also} be called monadic \sout{rel} \ul pred.\zl icates\zd\ \ud in this order of ideas.

I don't want to go into any discussions of what \zl \ul \ud to be deleted\zd\ predicates are (that could lead $\mathbf{\llbracket 29. \rrbracket}$ to a discussion of nominalism and realism\zl ).\zd\ \sout{\zl unreadable text, perhaps: But\zd }\, I want to say about the essence of a predicate only this. In order that a predicate be well\zl -\zd defined it must be (uniquely and) unambiguously determined of any objects (whatsoever) whether the predicate belongs to them or not. So e.g.\ a property is given if it is uniquely determined of any object whether or not the pred.\zl icate\zd\ bel.\zl ongs\zd\ to it and a dyadic rel\zl ation\zd\ is given if it is \dots\zl uniquely\zd\ det.\zl ermined\zd\ of any two obj.\zl ects\zd\ whether or not the rel.\zl ation\zd\ subsists betw.\zl een\zd\ them\zl .\zd\ \sout{is the only essential property to be required of a predicate} I shall use capital letters \sout{greek letters $\varphi,\psi,\chi$} $M,P,$ to denote individual predicates---as e\zl .\zd g.\ mortal\zl ,\zd\ greater etc.\ \sout{\zl unreadable text\zd\ $p,q,r$ \zl unreadable text\zd\ to denote arbitrary prop. and I shall use} $\mathbf{\llbracket 30. \rrbracket}$ \zl and\zd\ small letters $a,b,c$ to denote \sout{arbitrary} \ul individual \ud objects \ul as e\zl .\zd g.\ Socr\zl ates\zd , New Yor\zl k\zd\ etc\zl .\zd\ \ud (of which the pred\zl icates\zd\ \sout{$\varphi,\psi$} $M,P\ldots$ are \sout{to be} asserted). Those objects are usually called \underline{individuals} in math.\zl ematical\zd\ logic\zl .\zd\ \zl The following sentence is crossed out in the text: \ul So the individuals are the domain of things for which the pred.\zl icates\zd\ are defined so that it is uniquely det.\zl ermined\zd\ for any ind.\zl ividual\zd\ whether or not a cert\zl ain\zd\ pred.\zl icate\zd\ bel.\zl ongs\zd\ to them. \ud \zd\

Now let $M$ be a monadic pred\zl icate\zd\ \sout{(i.e.\ a quality)} \ul \zl ,\zd\ e.g\zl .\zd\ ,,\zl ``\zd man $\!$\zl mor\-tal\zd '', \ud and $a$ an indiv.\zl idual\zd\ \ul \zl ,\zd\ e.g\zl .\zd\ Socr\zl ates\zd . \ud Then the prop.\zl osi\-tion\zd\ that $M$ belongs to $a$ is denoted by $M(a)$\zl .\zd\ So $M(a)$ means ,,\zl ``\zd Socrates is mortal\zl ''\zd\ and similarly if $G$ is a di\zl y\zd adic relation \ul\zl ,\zd\ e.g\zl .\zd\ larger\zl ,\zd\ \ud and $b,c$ two ind.\zl ividuals\zd\ \ul\zl ,\zd\ e.g.\ New \zl York and\zd\ Wash\zl ington,\zd\ \ud then $G(b,c)$ means ,,\zl ``\zd The rel.\zl ation\zd\ $G$ subsists between $b$ and $c$.''\zl $c$'',\zd\ i\zl .\zd e.\ in our case ,,\zl ``\zd New York \zl is\zd\ larger than Wash.\zl ington\zd ''. So in this notation there is no copula\zl ,\zd\ but e.g\zl .\zd\ the prop.\zl osition\zd \zl ``\zd Soc\-ra\-tes is mortal\zl ''\zd\ $\mathbf{\llbracket 31. \rrbracket}$ has to be expr.\zl essed\zd\ like this Mortality(Socrates)\zl ,\zd\ and that New York is greater than W.\zl ashington\zd\ by Larger(New York, Wash.\zl ington\zd )\zl .\zd\

That much I have to say about the simplest type of prop.\zl ositions\zd\ which \sout{simply} say that some \ul def.\zl inite\zd\ \ud pred.\zl icate\zd\ belongs to some \ul def.\zl inite\zd\ subject or subjects. These prop\zl ositions\zd\ are sometimes called atomic prop.\zl ositions\zd\ \sout{in math\zl ematical\zd\ logic} bec.\zl ause\zd\ they constitute so to speak the atoms of which the more compl.\zl ex\zd\ propositions are built up. But now how are they built up? We know already one way of forming $\mathbf{\llbracket 32. \rrbracket}$ compound propositions namely by means of the operations of the propos.\zl itional\zd\ calculus $.\, \lfloor ,\rfloor \vee\lfloor ,\rfloor \supset$ etc.\zl ,\zd\ e.g.\ from the two atomic prop.\zl ositions ``\zd Socr.\zl ates\zd\ is a man\zl ''\zd\ and \zl ``\zd Socr.\zl ates\zd\ is mortal\zl ''\zd\ we can form the composit prop.\zl osition ``\zd If Socr.\zl ates\zd\ is a man Socr.\zl ates\zd\ is mortal\zl '';\zd\ \sout{\zl unreadable symbol\zd }\, i\zl written over I\zd n symb.\zl ols,\zd\ if $T$ denotes \sout{\zl unreadable text\zd }\, \ul the pred.\zl icate\zd\ of \ud mortality \sout{\zl unreadable text\zd\ \ul the indiv.\zl idual\zd\ \ud Socrates} it would read $M(a)\supset T(a)$\zl ,\zd\ or e.g\zl .\zd\ $M(a)\:\vee\sim M(a)$ would mean ,,\zl ``\zd Either Socr.\zl ates\zd\ is a man or Socr.\zl ates\zd\ is not a man''. $M(a)\: .\: T(a)$ would mean ,,\zl ``\zd Socr.\zl ates\zd\ is a man and Socrates is mortal'', and so on. The prop.\zl ositions\zd which we can obtain in this way\zl ,\zd\ i\zl .\zd e.\ by combining atomic prop.\zl ositions\zd\ by means \zl The next two pages are again numbered \textbf{31}.\ and \textbf{32}.\ in the manuscript; they are numbered \textbf{31}.\textbf{a} and \textbf{32}.\textbf{a} here.\zd\ $\mathbf{\llbracket 31.a \rrbracket}$ of the truth functions $\vee\lfloor ,\rfloor .$ etc\zl .\zd\ are sometimes called molecular prop\zl ositions\zd .

\zl new paragraph\zd\ But there is still another way of forming compound prop.\zl osi\-tions\zd\ which we have not yet taken account of in our symbolism\zl ,\zd\ namely by means of the particles ,,\zl ``\zd every'' and ,,\zl ``\zd some''. These are expressed in logistics by the use of variables as follows: Take e.g.\ the prop\zl osi\-tion\zd\ ,,\zl ``\zd Every man is mortal''\zl .\zd\ We can express \ul that in other words like this\zl :\zd\ \ud ,,\zl ``\zd Every object which is a man is mortal'' or ,,\zl ``\zd For every \sout{individual} \ul object \ud $x$ it is true that $M(x)\supset T(x)$''\zl .\zd\ Now in order to indicate \zl comma from the manuscript deleted\zd\ that this implication $\mathbf{\llbracket 32.a \rrbracket}$ is asserted of any object $x$ one puts $x$ in brackets in front of the prop.\zl osition\zd\ and includes the whole prop\zl osition\zd\ in bracket\zl s\zd\ to indicate that the whole prop.\zl osition\zd\ is asserted to be true for every $x$. And generally if we have an arb.\zl itrary\zd\ exp\zl ression,\zd\ say $\Phi(x)$ which involves a variable $x$\zl ,\zd\ then $(x)[\Phi(x)]$ means ,,\zl ``\zd For every object $x$, $\Phi(x)$ is true''\zl ,\zd\ i\zl .\zd e.\ if you take an arbitrary individual $a$ and substitute it for $x$ then the resulting prop.\zl osition\zd\ $\Phi(a)$ is true. As in our example $(x)[M(x)\supset T(x)]$, \sout{where $M$ means man} $\mathbf{\llbracket 33. \rrbracket}$ \sout{and $\Psi$ means mortal} if you subst.\zl itute\zd\ Socrates for $x$ you get the true prop\zl osition\zd . And gen.\zl erally\zd\ if you subst\zl itute\zd\ for $x$ something which is a man you get a true prop.\zl osition\zd\ bec.\zl ause\zd\ then the first and sec.\zl ond\zd\ term of the impl.\zl ication\zd\ are true. If however you subst\zl itute\zd\ someth.\zl ing\zd\ which is not a man you also get a true prop\zl osition\zd\ \ul bec.\zl ause\zd\ \ldots\ \ud So for \ul any \ud arb.\zl itrary\zd\ obj.\zl ect\zd\ which you subst\zl itute\zd\ for $x$ \zl \sout{$\Big\lceil$}\zd\ you get a true prop\zl osition\zd\ and this is indicated by writing $(x)$ in front of the prop\zl osition\zd . $(x)$ is called the universal quantifier\zl .\zd\

\zl The following text in square brackets is crossed out in the manuscript.\zd\ [I \ul wish to \ud \zl unreadable word\zd\ that exactly as formerly I used \zl unreadable word\zd\ to denote arb\zl itrary\zd\ expressions I denote now by $\Phi$ etc expr.\zl essions\zd\ which may involve variable $x$ which I indicate by writing them after the $\Phi$. An expression which involves variables and which becomes a prop.\zl osition\zd\ if you replace the var\zl iable\zd\ by $\mathbf{\llbracket 34. \rrbracket}$ individual objects is called a prop\zl ositional\zd\ funct\zl ion\zd . So e.g.\ $\varphi(x)$ is a prop\zl ositional\zd\ funct\zl ion\zd\ or $\varphi(x)\supset \psi(x)$ because\ldots]

As to the particle ,,\zl ``\zd some'' or ,,\zl ``\zd there exists'' \sout{(which is the same thing)} it is expr\zl essed\zd\ by a reversed $\exists$ put in brackets together with a var\zl iable\zd\ $(\exists x)$. So that means: there is an object $x$\zl ;\zd\ e.g.\ if we want to express that some men are not mortal we have to write $(\exists x)[M(x)\: .\: \sim T(x)]$ and generally if $\Phi(x)$ is a prop.\zl ositional\zd\ funct\zl ion\zd\ with the var.\zl iable\zd\ $x$\zl ,\zd\ $(\exists x)[\Phi(x)]$ means $\mathbf{\llbracket 35. \rrbracket}$ ,,\zl ``\zd There exits some object $a$ such that $\Phi(a)$ is true''. Nothing is said about the nu.\zl mber\zd\ of obj.\ ects\zd\ \ul for which $\Phi(a)$ is true \ud \sout{that exist}; there may be one or several\zl .\zd\ $(\exists x)\Phi(x)$ only means there is at least one obj\zl ect\zd\ $x$ such that $\Phi(x)$. $(\exists x)$ is called the existential quantifier\zl .\zd\ From this def\zl inition\zd\ you see at once that we have the following equivalences:
\begin{tabbing}
\hspace{10em}$(\exists x)\Phi(x)\;$\=$\equiv\; \sim(x)[\sim\Phi(x)]$\\[.5ex]
\hspace{10em}$(x)\Phi(x)$\>$\equiv\; \sim(\exists x)[\sim\Phi(x)]$
\end{tabbing}
\zl After these displayed formulae the page is divided in the manuscript by a horizontal line.\zd\

Generally $(x)[\sim\Phi(x)]$ means $\Phi(x)$ holds for no obj.\zl ect\zd\ and $\sim (\exists x)[\Phi(x)]$ $''$\zl means\zd \zl t\zd here is no \ul object \ud $x$ such that $\Phi(x)$\zl .\zd\ Again you see that these two \ul statements \ud are equivalent \ul with each other \ud . It is \sout{now} easy e.g.\ to express the traditional \ul four \ud $\mathbf{\llbracket 36. \rrbracket}$ types of prop.\zl ositions\zd\ a, e, i, o in our notation. In each case we have two predicates\zl ,\zd\ say $P$, $S$ and
\begin{tabbing}
\hspace{4em}\=$S$a$P$ \hspace{.3em}\=means\hspace{.4em} \=every $S$ is a $P$ \hspace{1.3em}\=i\zl .\zd e.\ \hspace{.2em}\=$(x)[S(x)\supset P(x)]$\\[.5ex]
\>$S$i$P$ \>means \>some $S$ are $P$ \>i\zl .\zd e.\ \>$(\exists x)[S(x)\: .\: P(x)]$\\[.5ex]
\>$S$e$P$ \>means \>no $S$ is a $P$ \>i\zl .\zd e.\ \>$(x)[S(x)\supset\;\sim P(x)]$\\[.5ex]
\>$S$o$P$ \>means \>some $S$ are $\sim P$ \>i\zl .\zd e.\ \>$(\exists x)[S(x)\: .\: \sim P(x)]$
\end{tabbing}
You see the universal prop.\zl ositions\zd\ have the universal quantifier in front of them and the part.\zl icular\zd\ prop.\zl ositions\zd\ the exist.\zl ential\zd\ quantifier. I want to mention that in classical logic two entirely different types of prop\zl osi\-tions\zd\ are counted as univ.\zl ersal\zd\ affirm.\zl ative,\zd\ namely prop\zl ositions\zd\ of the type \zl ``\zd Socrates is mortal\zl ''\zd\ expressed by $P(a)$ and ,,\zl ``\zd Every man is mortal\zl ''\zd\ $(x)[S(x)\supset P(x)]$\zl .\zd\

$\mathbf{\llbracket 37. \rrbracket}$ Now the existential and univ\zl ersal\zd\ quantifier can be combined with each other and with the truth\zl $\,$\zd f\zl unctions\zd\ \ul $\sim$,$\,$\ldots\ \ud in many ways so as to express more complicated prop\zl ositions\zd . \zl Here one finds in the manuscript a page numbered \textbf{37}.\textbf{1} inserted within p.\ \textbf{37}.\zd\

$\mathbf{\llbracket 37.1 \rrbracket}$ Thereby one uses some abbrev.\zl iations,\zd\ namely: Let $\Phi(xy)$ be an expr.\zl ession\zd\ cont\zl aining\zd\ 2\zl two\zd\ var.\zl iables;\zd\ then we may form: $(x)[(y)$ $[\Phi(xy)]]$ \zl .\zd\ That means ,,\zl ``\zd For any obj\zl ect\zd\ $x$ it is true that for any obj.\zl ect\zd\ $y$ $\Phi(xy)$'' that evidently means ,,\zl ``\zd $\Phi(xy)$ is true whatever objects you take for $x,y$'' and this is den\zl oted\zd\ by $(x,y)\Phi(xy)$. Evidently the order of the var\zl iables\zd\ is arb.\zl itrary\zd\ here\zl ,\zd\ i.e.\ $(x\lfloor ,\rfloor y)\Phi(xy)\equiv (y\lfloor ,\rfloor x)\Phi(xy)$. Similarly $(\exists x)[(\exists y)[\Phi(xy)]]$ means ,,\zl ``\zd\ There are some obj\zl ects\zd\ $x,y$ such that $\Phi(xy)$'' and this is abbr.\zl eviated\zd\ by $(\exists x,y)\Phi(xy)$ \ul and means: \zl text missing\zd\ \ud \ul But \ud it has to be noted \zl comma from the manuscript deleted\zd\ \sout{here} that this does not mean that there are really two diff.\zl erent\zd\ obj.\zl ects\zd\ $x,y$ satisfying $\Phi(xy)$\zl .\zd\ This formula is also \zl \sout{be}\zd\ true if there is one obj\zl ect\zd\ $a$ such that $\Phi(a,a)$\zl $\Phi(aa)$\zd\ bec.\zl ause\zd\ then there exists an $x$\zl ,\zd\ namely $a$, such that there exists a $y$\zl ,\zd\ namely again $a$\zl ,\zd\ such that etc. \zl At this place \ul \underline{Expl.} \ud is inserted, and the following text from the end of p.\ \textbf{37}.\textbf{1} seems to refer, by having at its end ``p 37'' and a sign for insertion, to this spot:
 \[
 {\rm Ex}
 \begin{cases}
  \text{Throug\zl h\zd\ any two points there exists a straight line\zl .\zd } \\
  \text{In any plane there exist to\zl two\zd\ $\mid\mid$ \zl parallel\zd\ lines\zl .\zd }
 \end{cases}
\]
These may be examples of universal and existential quantification that, unlike $(\exists x,y)\Phi(xy)$, involve variables standing for different objects, but the first is related to an example for notation on p.\ \textbf{39}.\ below.\zd\ Again $(\exists x,y)\Phi(xy)\equiv (\exists y,x)\Phi(xy)$\zl .\zd\

But it is to be noted that this interchangeability holds $\mathbf{\llbracket new\, page \rrbracket}$ only for two univ.\zl ersal\zd\ or two exist.\zl ential\zd\ quant\zl ifiers\zd . It does not hold for an univ.\zl ersal\zd\ and an exist\zl ential\zd\ quant.\zl ifier,\zd\ i\zl .\zd e.\ $(x)[(\exists y)[\Phi(yx)]]\not\equiv (\exists y)[(x)[\Phi(yx)]]$\zl .\zd\ Take e.g.\ for $\Phi(yx)$ the prop\zl osition\zd \zl ``\zd $y$ greater than $x$\zl '';\zd\ then the first means ,,\zl ``\zd For any obj.\zl ect\zd\ $x$ it is true that there exists \sout{exists} an obj\zl ect\zd\ $y$ greater than $x$''\zl ;\zd\ in other words \zl ``F\zd or any object there exists something greater''. The right\zl -\zd hand side however means ,,t\zl ``T\zd here exists an obj\zl ect\zd\ $y$ such that for any $x$ $y$ is greater than $x$''\zl ,\zd\ there exists a greatest obj\zl ect\zd . So that \sout{means} in our case \ul the right side \ud \zl says\zd\ just the oppos.\zl ite\zd\ of what the left \sout{hand} side says. \sout{As to the brackets} \zl T\zd he above abbrev.\zl iation is\zd\ also used for more than two var.\zl iables,\zd\ i\zl .\zd e.\ $(x\lfloor ,\rfloor y\lfloor ,\rfloor z)[\Phi(xyz)]$\zl ,\zd\ $(\exists x\lfloor ,\rfloor y\lfloor ,\rfloor z)[\Phi(xyz)]$\zl .\zd\

\zl Here one returns to p.\ \textbf{37}.\zd I want now to give some examples for the notation introduced. Take e.g\zl .\zd\ the prop.\zl osition\zd\ ,,\zl ``\zd For any integer there exists a greater one''\zl .\zd\ The pred\zl icates\zd\ occurring in this prop\zl osition\zd\ are: 1. integer and 2. greater\zl .\zd\ Let us denote them by $I$ and $>$ \ul so $I(x)$ is to be read\ldots\ \zl ``$x$ is an integer'' and\zd\ $>\!(x,y)$ $''$ $''$ $''$ $''$\zl $>\!(xy)$ is to be read\zd \zl ``\zd $x$ greater $y$\zl ''\zd\ or \zl ``\zd $y$ smaller $x$\zl ''.\zd\ Then the prop\zl osition\zd\ \ul \zl unreadable text\zd\ \ud is expressed in log.\zl istic\zd\ symb\zl olism\zd\ as follows:
\begin{tabbing}
\hspace{1.7em}$(x)[I(x)\supset (\exists y)[I(y)\: .\: >\!(y,x)\lfloor >\!(yx)\rfloor ]]$\zl .\zd\
\end{tabbing}
We can express the same fact by saying $\mathbf{\llbracket 38. \rrbracket}$ there is no greatest integer:\zl .\zd\ What would that \ul look like in logist.\zl ic\zd\ symb.\zl olism:\zd\ \ud
\begin{tabbing}
\hspace{1.7em}$\sim\!(\exists x)[I(x)\: .\:$ \ul such that no int\zl eger\zd\ is greater i\zl .\zd e.\ \ud \\
\` $(y)[I(y)\supset \;\sim \;>\!(yx)]$\zl ].\zd\
\end{tabbing}
As another ex.\zl ample\zd\ take the prop\zl osition\zd\ ,,\zl ``\zd There is a smallest int.\zl e\-ger\zd '' that would read\zl :\zd\
\begin{tabbing}
\hspace{1.7em}$(\exists x)[I(x)\: .\:$ \ul such that no int.\zl eger\zd\ is smaller i\zl .\zd e.\ \ud \\
\` $(y)[I(y)\supset \;\sim \;>\!(x,y)\lfloor >\!(xy)\rfloor ]]$.
\end{tabbing}
I wish to call your attention to a near at hand mistake. It would be wrong to express this \ul last \ud prop.\zl osition\zd\ like this:
\begin{tabbing}
\hspace{1.7em}$(\exists x)[I(x)\: .\:(y)[I(y)\supset\; >\!(yx)]]$
\end{tabbing}
bec.\zl ause\zd\ that would \ul mean \ud there is an int.\zl eger\zd\ smaller than every int\zl eger\zd . But such an int.\zl eger\zd\ does not exist $\mathbf{\llbracket 39. \rrbracket}$ since it would have to be smaller than itself. An integer smaller than every int.\zl eger\zd\ would have to be smaller than \ul itself\zl ---\zd that is clear\zl .\zd\ \ud So the sec.\zl ond\zd\ prop.\zl osition\zd\ is false whereas the first is true, bec.\zl ause\zd\ it says only there exists an int.\zl eger\zd\ $x$ \zl full stop deleted\zd\ which is not greater than any int\zl eger\zd . \sout{and that is true, because \zl unreadable text\zd\ has this prop. that it is greater \ul than \ud \zl unreadable text\zd\ itself (not greater than itself either).}

Another ex.\zl ample\zd\ for our not.\zl ation\zd\ may be taken from Geom\zl etry\zd . Consider the prop.\zl osition\zd\ ,,\zl ``\zd Through any two different points there is exactly one straight line''. The pred.\zl icates\zd\ which occur in this prop.\zl osition\zd\ are 1. \zl p\zd oint $P(x)$\zl ,\zd\ $\mathbf{\llbracket 40. \rrbracket}$ 2. straight line $L(x)$\zl ,\zd\ 3. \underline{different} that is the neg\zl ation\zd\ of identity\zl .\zd\ Identity is den.\zl oted\zd\ by $=$ and diff\zl erence\zd\ \sout{sometimes} by $\neq$\zl .\zd\ =$(xy)$ means $x$ and $y$ are the same thing\zl ,\zd\ e.g.\ =$\,$(Shake\-speare, author of Hamlet)\zl ,\zd\ and $\neq(xy)$ \zl means\zd\ $x$ and $y$ are different from each other\zl .\zd\ There is \ul still \ud another relation \zl comma from the manuscript deleted\zd\ that occurs in \ul our geom.\zl etric\zd\ prop.\zl osition,\zd\ namely the one \ud expressed by \ul the word \ud ,,\zl ``\zd through'' \sout{w}\zl .\zd\ That is the rel.\zl ation\zd\ which holds betw\zl een\zd\ a point \ul $x$ \ud and a line \ul $y$ \ud if ,,\zl ``\zd $y$ passes through $x$'' or in other words \zl ,, deleted\zd\ if \zl ``\zd $x$ lies on $y$''. Let us den\zl ote\zd\ that \ul relation \ud by $J(x,y)$ \zl $J(xy)$\zd . Then the \ul geom.\zl etric\zd\ \ud prop.\zl osition\zd\ ment.\zl ioned\zd , in order to be expressed in \sout{our} \ul log\zl istic\zd\ \ud symb\zl olism\zd , has to be splitted into to\zl two\zd\ parts\zl ,\zd\ namely there is at least one line and there is at most one line. The first reads: $(x,y)[P(x)\: .\: P(y)\: .\: \neq(xy)\supset$ $\mathbf{\llbracket 41. \rrbracket}$ $(\exists u)[L(u)\: .\: J(xu)\: .\: J(yu)]]$\zl .\zd\ So that means that through any two diff\zl erent\zd\ points there is\ldots\ But it is not excl.\zl uded\zd\ \ul by that statement \ud that there are two or three diff.\zl erent\zd\ lines passing through two points. \sout{To express} That there are no \ul two \ud diff\zl erent\zd\ lines could be expr.\zl essed\zd\ like this
\begin{tabbing}
\hspace{1.7em}$(x,y)[P(x)\: .\: P(y)\: .\: \neq(xy)\supset\;\sim(\exists u,v)[L(u)\: .\: L(v)\: .\: \neq(u,v)\lfloor \neq(uv)\: .\:\rfloor$ \\
\` $J(xu)\: .\: J(yu)\: .\: J(xv)\: .\: J(yv)]]$
\end{tabbing}

I hope these ex.\zl amples\zd\ will suffice to make \ul clear how the quantifiers are to be used. \ud
For any quantifier occurring in an expr\zl ession\zd\ there is a definite portion of the expr\zl ession\zd\ to which it relates (called the scope of the expression)\zl ,\zd\ e.g.\ scope of $x$ whole expr.\zl ession,\zd\ of $y$ only this portion\ldots\ So the scope it\zl is\zd\ the prop.\zl osition\zd\ of which it is asserted that it holds for all or every obj\zl ect\zd . \zl I\zd t is indicated by the brackets which begin\zl \sout{s}\zd\ immediately behind the quantifier. There are some conv\zl entions\zd\ about leaving out this\zl these\zd\ brack.\zl ets,\zd\ namely they may be left out 1. \zl i\zd f \ul the \ud scope is atomic\zl ,\zd\ e.g.\ $(x)\varphi(x)\vee\lfloor \supset\rfloor p$\zl :\zd\ $(x)[\varphi(x)]\supset p$\zl ,\zd\ not $(x)[\varphi(x)\supset p]$\zl ,\zd\ 2. if the scope begins with $\sim$ or a quant\zl ifier,\zd\ e.g\zl .\zd\
\begin{tabbing}
\hspace{1.7em}\=$(x)\sim[\varphi(x)\: .\: \psi(x)]\vee p\;$ \=: $\; (x)[\sim[\varphi(x)\: .\: \psi(x)]]\vee p$\\[.5ex]
\` (\zl perhaps proof correction mark for delete,\\ \` indicating that $\varphi,\psi$ are to be replaced by\zd\ $Q,R$)\\[.5ex]
\>$(x)(\exists y)\varphi(x)\vee p$\>: $\; (x)[(\exists y)[\psi\lfloor \varphi\rfloor (x)]]\vee p$
\end{tabbing}
But these rules are only facultative\zl ,\zd\ i\zl .\zd e.\ we may also write all the brackets if \ul it is \ud expedient for the sake of clarity\zl .\zd\

A variable to which a quantifier $(x)$\zl ,\zd\ $(y)$\zl ,\zd\ $(\exists x)$\zl ,\zd\ $(\exists y)$ refers is called a ,,\zl ``\zd bound variable''. In the examples which I gave, all variables $\mathbf{\llbracket 42. \rrbracket}$ are bound (e.g.\ to this $x$ relates this quant.\zl ifier\zd\ etc\zl .\zd ) and similarly to any var.\zl iable\zd\ occurring in those expr.\zl essions\zd\ you can associate a quantifier which refers to it. If however you take e.g.\ the exp\zl ression\zd : $I(y)\: .\: (\exists x)[I(x)\: .\: >(yx)]$.\zl ,\zd\ which means: there is an int.\zl eger\zd\ $x$ smaller than $y$.\zl , t\zd hen here $x$ is a bound var.\zl iable\zd\ bec\zl ause\zd\ the quantifier $(\exists x)$ refers to it. But $y$ is not bound bec\zl ause\zd\ the expr\zl ession\zd\ contains no quantifier referring to it\zl .\zd\ Therefore $y$ is called a free variable of this expression. An expr.\zl ession\zd\ containing free variables is not a propos.\zl ition\zd , but it only becomes a prop.\zl osition\zd\ if the free variables are replaced by individual objects, e.g.\ this expression here means $\mathbf{\llbracket 43. \rrbracket}$ ,,\zl ``\zd There is an int.\zl eger\zd\ smaller than \ul the int.\zl eger\zd\ \ud $y$''. That evidently is not a \ul definite \ud assertion which is either true or wrong. But if you subst.\zl itute\zd\ for the free var.\zl iable\zd\ $y$ a definite obj.\zl ect,\zd\ e.g.\ 7\zl ,\zd\ then you obtain a definite prop.\zl osition,\zd\ namely:\zl ``\zd There is an int.\zl eger\zd\ \ul smaller \zl than\zd\ 7\zl ''.\zd\ \ud

\zl The paragraph that starts here and the next, which are entirely crossed out, are on p.\ \textbf{43}.\textbf{1}, inserted within p.\ \textbf{43} of the manuscript.\zd\ Expressions containing free var.\zl iables\zd\ \ul and such that they become prop.\zl ositions\zd\ if\ldots\ \ud are called prop.\zl ositional\zd\ funct\zl ions\zd . Here we have a prop.\zl ositional\zd\ funct\zl ion\zd\ with one free var\zl iable\zd . There are also such \zl functions\zd\ with two or more free variables. Any prop.\zl ositional\zd\ f\zl u\zd unct\zl ion\zd\ with one var.\zl iable\zd\ def.\zl ines\zd\ a cert.\zl ain\zd\ prop\zl erty\zd\ and one with two variables a cert.\zl ain\zd\ dyadic rel\zl ation\zd .

\zl The scope of quantifiers mentioned in the crossed out paragraph that starts here is considered in a text inserted on p.\ \textbf{41}.\zd\ To any quantif.\zl ier\zd\ occurring in an expr\zl ession\zd\ there is a definite portion of the expr\zl ession\zd\ to which it relates, which is called the scope of the quantif.\zl ier\zd ,\zl ;\zd\ it is indicated by the brackets\zl ,\zd\ which opens immediately after the quantifier\zl ,\zd\ e.g.\ the scope of $(x)$ in\ldots\ is the whole expr.\zl ession:\zd it says for any $x$ the whole \ul subsequ.\zl ent\zd\ \ud prop.\zl osition\zd\ is true\zl ;\zd\ the scope of $y$ \ul here \ud is the rest of this exp.\zl ression\zd\ bec\zl ause\zd\ it says there is a $y$ for which\ldots\ You see also that this bracket closes up here and this bracket here\zl .\zd\

The bound variables have the property that it is entirely irrelevant by which letters they are denoted\zl ;\zd\ e.g.\ $(x)(\exists y)[\Phi(xy)]$ means exactly the same thing as $(u)(\exists v)[\Phi(uv)]$\zl . T\zd he only requirement is that you must use different letters for different bound variables\zl .\zd\ But even that is only necessary for variables $\mathbf{\llbracket 44. \rrbracket}$ one of whom is \sout{one} contained in \ul the scope of the \ud \sout{each} other as e\zl .\zd g\zl .\ in\zd\ $(x)[(\exists y)\Phi(xy)]$\zl , w\zd here $y$ is in the scope of $x$ which is the whole expr\zl ession, and\zd\ therefor it has to be den\zl oted\zd\ by a letter diff.\zl erent\zd\ from $x$\zl ;\zd\ $(x)[(\exists x)\Phi(xx)]$ would be ambiguous. Bound variables whose scopes lie outside of each other \ul however can \ul be denoted by the same letter without any ambiguity\zl ,\zd\ e.g.\ $(x)\varphi(x)\supset(x)\psi(x)$. For the sake of clarity we also require that the free variables in a prop.\zl ositional\zd\ f\zl u\zd nct\zl ion\zd\ should always be denoted by letters different from the bound var\zl iables;\zd\ so e\zl .\zd g\zl .\zd\ $\varphi(x)\: .\: (x)\psi(x)$ is no\zl t a\zd\ correctly formed prop\zl ositional\zd\ \ul f\zl u\zd nct\zl ion,\zd\ \ud but $\varphi(x)\: .\: (y)\psi(y)$ is one\zl .\zd\

The examples \ul of formulas \ud which I gave \ul last time and also the problems to be \ul solved \zl unreadable word, perhaps ``for''\zd\ \ud \ud \sout{so far} were propositions concerning cert.\zl ain\zd\ definite \ul \sout{\zl unreadable text\zd }\, \ud predicates $I$, $<$, =, etc. They are true only for those part\zl icular\zd\ pred.\zl icates\zd\ occurring in them. But now exactly as \sout{we had} in the calc\zl ulus\zd\ of prop\zl ositions\zd\ \ul there are \ud cert.\zl ain\zd\ formulas which are true whatever prop.\zl ositions\zd\ the letters $p,q,r$ may be so also in the \sout{enlarged} calculus of pred.\zl icates\zd\ $\mathbf{\llbracket 45. \rrbracket}$ there will be certain formulas which are true for any \ul arbitrary \ud predicates. I denote arb\zl itrary\zd\ pred.\zl icates\zd\ by small Greek letters $\varphi,\psi$\zl .\zd\ So these are supposed to be variables for predicates exactly as $p,q\ldots$ are variables for prop.\zl ositions\zd\ and $x,y,z$ are variables for obj.\zl ects.\zd\ \zl \sout{\ul individual pred\zl i\-cates\zd . \ud }\zd\

Now take e.g.\ \ul the prop\zl osition\zd\ $(x)\varphi(x)\vee(\exists x)\sim\varphi(x)$\zl ,\zd\ i\zl .\zd e.\ ,,\zl ``\zd Ei\-ther every ind.\zl ividual\zd\ has the prop.\zl erty $\varphi$ or ther is an indiv\zl idual\zd\ which has not the prop\zl erty\zd\ $\varphi$''\zl .\zd\ That will be true for any arbitrary \ul monadic \ud pred.\zl icate\zd\ $\varphi$\zl .\zd\ We \zl had\zd\ other examples before\zl ,\zd\ e\zl .\zd g\zl .\zd\ ${(x)\varphi(x)\equiv\:\sim(\exists x)\sim\varphi(x)}$ that again is true for \zl text omitted in the manu\-script, should be: any arbitrary monadic predicate $\varphi$.\zd\ Now exactly as in the calc.\zl u\-lus\zd\ of prop.\zl ositions\zd\ such expr.\zl essions\zd\ which are true for all pred.\zl icates\zd are called tautologies \sout{or logically true} or universally true. Among them are e.g.\ all the form.\zl ulas\zd\ which express the Arist.\zl otelian\zd\ $\mathbf{\llbracket 46. \rrbracket}$ moods of inf.\zl erence,\zd\ e.g.\ \ul the \ud mood b\zl B\zd arb.\zl ara\zd\ is expr.\zl essed\zd\ like this:
\begin{tabbing}
\hspace{1.7em}$(x)[\varphi(x)\supset\psi(x)]\: .\: (x)[\psi(x)\supset\chi(x)]\supset (x)[\varphi(x)\supset\chi(x)]$
\end{tabbing}
The mood d\zl D\zd arii \zl \sout{I}\zd\ like this
\begin{tabbing}
\hspace{11em}\=$\varphi$\hspace{1em}\= $M$a$P$\hspace{1em}\=$\psi$\\[.5ex]
\>$\chi$\>\uline{$S$i$M$}\>$\varphi$\\[.5ex]
\>\>$S$i$P$\\[1ex]
\hspace{1.7em}$(x)[\varphi(x)\supset\psi(x)]\: .\: (\exists x)[\chi(x)\: .\:\varphi(x)]\supset (\exists x)[\chi(x)\: .\:\psi(x)]$
\end{tabbing}

It is of course the chief aim of logic to investigate the\zl written over something else\zd\ \sout{\zl unreadable symbol\zd }\, tautologies and exactly as in the calc.\zl ulus\zd\ of prop.\zl ositions\zd\ there are \ul again \ud two chief problems which arise. Namely \zl :\zd\ 1.\ To develop methods for finding out about a given expr.\zl ession\zd\ whether or not it is a tautology\zl ,\zd\ 2.\ To reduce all taut.\zl ologies\zd\ to a finite nu.\zl mber\zd\ of logical axioms and rules of inf.\zl erence\zd\ from which they can be derived. I wish to mention right now that only $\mathbf{\llbracket 47. \rrbracket}$ the second problem can be solved \sout{satisfactorily} for the calc.\zl ulus\zd\ of pred\zl icates\zd . One has actually succeeded in setting up a system of ax.\zl ioms\zd\ for it and in proving its completeness (i\zl .\zd e.\ that every taut.\zl ology\zd\ can be derived from it)\zl .\zd\

\zl new paragraph\zd\ As to the first problem\zl ,\zd\ \ul the so called decision probl.\zl em,\zd\ \ud it has also been solved \ul in a sense \ud but \ul in the \ud negative\zl ,\zd\ i\zl .\zd e.\ one has succeeded in \ul proving \ud that there does not \ul exist any \ud mechanical proced.\zl ure\zd\ to decide of any given expression whether or not it is a tautology \ul of the calc.\zl ulus\zd\ of pred\zl icates\zd . \ud That does not mean mean that there are \ul any individual \ud formulas of which one could not decide whether or not they are $\mathbf{\llbracket 48. \rrbracket}$ taut\zl ologies\zd . It only means that it is not poss.\zl ible\zd\ to decide that by a \ul purely \ud mech.\zl anical\zd\ procedure. For the calc.\zl ulus\zd\ of prop.\zl ositions\zd\ this was possible\zl ,\zd\ e.g.\ the truth\zl -\zd table method is a purely mec.\zl hanical\zd\ proc.\zl edure\zd\ which allows to decide of any given expr.\zl ession\zd\ whether or not it is a taut\zl ology\zd . So what has been proved is only that a similar thing cannot exist for the calc\zl ulus\zd\ of pred\zl icates\zd . However for certain \sout{particular} \ul special \ud kinds of formulas such methods of decision have been developed\zl ,\zd\ e.g.\ for all form.\zl ulas\zd\ with only monadic pred.\zl icates\zd\ (i\zl .\zd e.\ formulas without relations in it)\zl ;\zd\ $\mathbf{\llbracket 49. \rrbracket}$ e.g.\ all form.\zl ulas\zd\ expressing the Arist.\zl otelian\zd\ moods are of this type \zl full stop deleted\zd\ bec.\zl ause\zd\ no relations occur in the Arist.\zl otelian\zd\ moods.

Before going into more detail about that I must say a few more words about the notion of a taut.\zl ology\zd\ of the calc.\zl ulus\zd\ of pred\zl icates\zd .

There are also taut\zl ologies\zd\ which involve variables both for propositions and for pred.\zl icates,\zd\ e.g.\
\begin{tabbing}
\hspace{1.7em}$p\: .\:(x)\varphi(x)\equiv (x)[p\: .\: \varphi(x)]$
\end{tabbing}
i\zl .\zd e.\ if $p$ is an arb.\zl itrary\zd\ prop\zl osition\zd\ and $\varphi$ an arb.\zl itrary\zd\ pred.\zl icate\zd\ then the assertion on the left\zl ,\zd\ i.e\zl .\zd\ ,,\zl ``\zd $p$ is true and for every $x$\zl ,\zd\ $\varphi(x)$ is true'' is equivalent with the assertion on the right\zl ,\zd\ i\zl .\zd e.\ ,,\zl ``\zd for every obj.\zl ect\zd\ $\mathbf{\llbracket 50. \rrbracket}$ $x$ the conjunction $p\: .\: \varphi(x)$ is true''. Let us prove that\zl ,\zd\ i\zl .\zd e.\ let us prove that the left side implies the right side and vice versa the right side implies the left side \sout{I}. If the left side is true that means: $p$ is true and for every $x$\zl ,\zd\ $\varphi(x)$ is true\zl ,\zd\ but then the right side is also true bec.\zl ause\zd\ then for every $x$\zl ,\zd\ $p\: .\: \varphi(x)$ is evidently true \sout{[So the left side implies the right side]}. But also vice versa\zl :\zd\ If for every $x$\zl ,\zd\ \sout{$(x)[$} $p\: .\: \varphi(x)$ is true then 1. $p$ must be true bec.\zl ause\zd\ otherwise $p\: .\: \varphi(x)$ would be true for no $x$ and 2. $\varphi(x)$ must be true for every $x$ since by ass.\zl umption\zd\ even $p\: .\: \varphi(x)$ is true for every $x$. So you see this equiv.\zl alence\zd\ holds for any pred.\zl icate\zd\ $\varphi$\zl ,\zd\ $\mathbf{\llbracket 51. \rrbracket}$ i\zl .\zd e.\ it is a tautology.

\zl new paragraph\zd\ There are four analogous taut.\zl ologies\zd\ obtained by repl.\zl ac\-ing\zd\ $.$ by $\vee$ and the un.\zl iversal\zd\ qu.\zl antifier\zd\ by the exist.\zl ential\zd\ qu\zl antifier,\zd\ namely
\begin{tabbing}
\hspace{1.7em}\=2.\quad\=$p\vee (x)\varphi(x)\equiv (x)[p\vee\varphi(x)]$\\[.5ex]
\>3.\>$p\: .\:(\exists x)\varphi(x)\equiv (\exists x)[p\: .\: \varphi(x)]$\\[.5ex]
\>4.\>$p\vee(\exists x)\varphi(x)\equiv (\exists x)[p\vee\varphi(x)]$
\end{tabbing}
I shall give the proof for them later on\zl .\zd\ These 4\zl four\zd\ \ul formulas \ud are of a great importance because they allow to shift a quantifier over a symb\zl ol\zd\ of conj\zl unction\zd\ or disj\zl unction\zd . If you write $\sim p$ inst\zl ead\zd\ of $p$ in the first you get $[p\supset(x)\varphi(x)]\equiv (x)[p\supset\varphi(x)]$. This law of logic is used particularly frequently in proofs as you will see later\zl .\zd\ Other ex.\zl amples\zd\ of tautologies are e.g\zl .\zd\
\begin{tabbing}
\hspace{1.7em}\=$(x)\varphi(x)\: .\: (x)\psi(x)\equiv (x)[\varphi(x)\: .\:\psi(x)]$\\[.5ex]
\>$(\exists x)\varphi(x)\vee (\exists x)\psi(x)\equiv (\exists x)[\varphi(x)\vee\psi(x)]$
\end{tabbing}
or e.g.\
\begin{tabbing}
\hspace{1.7em}$\sim(x)(\exists y)\varphi(xy)\equiv (\exists x)(y)\sim\varphi(xy)$
\end{tabbing}
$\mathbf{\llbracket 52. \rrbracket}$ That means:

Proof\zl .\zd\ $\sim(x)(\exists y)\varphi(xy)\,{\equiv \atop \mbox{\scriptsize\rm means}}\,(\exists x)\sim(\exists y)\varphi(xy)$, but $\sim(\exists y)\varphi(xy)\equiv(y)\sim\varphi(xy)$ as we saw before. Hence the whole expr.\zl ession\zd\ is equiv.\zl alent\zd\ with $\equiv(\exists x)(y)\sim\varphi(xy)$ which was to be proved.
\vspace{-3ex}
\begin{tabbing}
\underline{\hspace{33em}}
\end{tabbing}
\vspace{-1ex}

Another ex\zl ample\zd : $(x)\varphi(x)\supset(\exists x)\varphi(x)$\zl ,\zd\ i\zl .\zd e.\ If \sout{$\varphi$ bel. to} every ind.\zl ividual\zd\ \ud has the prop\zl erty\zd\ $\varphi$ \ud then a fort.\zl iori\zd\ there are ind.\zl ividuals\zd\ which have the prop.\zl erty\zd\ $\varphi$. The inverse of this prop.\zl osition\zd\ \sout{no} is no taut.\zl ology,\zd\ i\zl .\zd e\zl .\zd\
\begin{tabbing}
\hspace{1.7em}$(\exists x)\varphi(x)\supset(x)\varphi(x)$ is not a taut.\zl ology\zd\
\end{tabbing}
bec.\zl ause\zd\ if there is an obj.\zl ect\zd\ $x$ which has the prop\zl erty\zd\ $\varphi$ that does not imply that every ind.\zl ividual\zd\ has the prop.\zl erty\zd\ $\varphi$.

\zl new paragraph\zd\ But here there is an \sout{\zl unreadable text, perhaps: ast.\zd }\, \ul impor\-tant \ud remark $\mathbf{\llbracket 53. \rrbracket}$ to be made. Namely: In order to prove that this form\zl ula\zd\ here is not a taut\zl ology\zd\ we must know that there exists more than one obj.\zl ect\zd\ in the world. For if we assume that there exists only one obj.\zl ect\zd\ in the world then this form\zl ula\zd\ would be true for every pred\zl icate\zd\ $\varphi$\zl ,\zd\ hence would be \sout{a taut.} \ul universally true \ud bec\zl ause\zd\ if there is only one obj\zl ect,\zd\ \ul say $a$\zl ,\zd\ \ud in the world then if there is an obj\zl ect\zd\ $x$ for which $\varphi(x)$ is true this obj\zl ect\zd\ must be $a$ (since by ass\zl umption\zd\ there is no other obj.\zl ect\zd )\zl ,\zd\ hence $\varphi(a)$ is true\zl ;\zd\ but then $\varphi$ is true for every obj.\zl ect\zd\ bec.\zl ause\zd\ by ass.\zl umption\zd\ there exists only this obj.\zl ect\zd\ $a$. \zl I.\zd e.\ in a world with only one $\mathbf{\llbracket 54. \rrbracket}$ obj.\zl ect\zd\ $(\exists x)\varphi(x)\supset(x)\varphi(x)$ is a taut\zl ology\zd . It is easy to \zl find\zd\ some expressions which are \sout{tautol.} \ul universally true \ud if there are only two ind\zl ividuals\zd\ in the world etc\zl .,\zd\ e.g.\
\begin{tabbing}
\hspace{1.7em}$(\exists x,y)[\psi(x)\: .\:\psi(y)\: .\:\varphi(x)\: .\:\sim\varphi(y)]\supset (x)[\psi(x)]$
\end{tabbing}

At present \ul I only wanted to point out that \ud the notion of a taut.\zl ology\zd\ of the calc.\zl ulus\zd\ of pred\zl icates\zd\ needs a further specific\zl ation\zd\ in order to be precise\zl .\zd\ This specif\zl ication\zd\ consists in this that an expr.\zl ession\zd\ is called a taut.\zl ology\zd\ only if it \ul \zl is\zd\ universally \ud true \sout{for \ul every pred. \ud} no matter how many ind.\zl ividuals\zd\ are in the world assuming only that there is at least one (otherwise the meaning \ul of the quantifiers is not \zl unreadable text, perhaps ``definite''\zd\ \ud \zl ).\zd\ So e\zl .\zd g.\ \ul $(x)\varphi(x)\supset(\exists y)\varphi(y)$\zl ;\zd\ \ud this is \zl a\zd\ taut.\zl ology\zd\ bec.\zl ause\zd\ it is true\ldots\ \zl but this \ul inverse \ud is not bec.\zl ause\zd\ \ldots\ It can be proved that this means the same thing as if I said: An expr\zl ession is\zd\ a taut\zl ology\zd\ if \zl it\zd\ is true in a world with infinitely many ind.\zl ividuals,\zd\ i.e.\ one can prove that \ul whenever an expr\zl ession\zd\ is univ.\zl ersally\zd\ \ud true in a world

\vspace{2ex}

\zl This text is continued on p.\ \textbf{55}., the first page of Notebook~V. On a new page after p.\ \textbf{54}., the last page of the present notebook, one finds the following jottings:\zd

\vspace{1ex}

\noindent interest lies in this, choice fortunate \underline{Ideenrealismus}, lie betw, \sout{greater} essence, (predicate is asserted of), individuals \sout{property (quality)} \underline{copula}, (built up of), (every), \zl unreadable text, presumably in shorthand\zd , (reversed $\exists$) \zl sign pointing to (every) above\zd\ \sout{property} \zl unreadable symbol\zd\ Hamlet. \sout{property belongs to} \zl underlined unreadable text, presumably in shorthand, pointing to Hamlet above\zd\ author
\[
(x,y)[P(x)\: .\: P(y)\: .\:\neq(xy)\supset (\exists u)(v)[L(v)\; J(xu)\: .\: J(yv)\equiv .\: v=u]]
\]
strict. implic.

\section{Notebook V}\label{0V}
\pagestyle{myheadings}\markboth{SOURCE TEXT}{NOTEBOOK V}
\zl Folder 63, on the front cover of the notebook ``Log.\zl ik\zd\ Vorl.\zl esungen\zd\ \zl German: Logic Lectures\zd\ N.D.\ \zl Notre Dame\zd\ V''\zd\

\vspace{1ex}

\zl The first page of this notebook, p.\ \textbf{55}., begins with the second part of a sentence interrupted at the end of p.\ \textbf{54}. of Notebook IV.\zd\

$\mathbf{\llbracket 55. \rrbracket}$ with infinitely many obj\zl ects\zd\ it is true in any world no matter how many ind.\zl ividuals\zd\ there may be and of course also vice versa. I shall not prove this equiv.\zl alence\zd\ but shall stick to the first definition.

The formulas by which we expressed the taut.\zl ologies\zd\ contain free var.\zl i\-ables\zd\ (not for individuals) but for predicates and for prop.\zl ositions,\zd\ e.g.\ $\varphi$ here is a free var\zl iable\zd\ in this expr.\zl ession\zd\ (no quant\zl ifier\zd\ related to it\zl ,\zd\ i\zl .\zd e.\ no $(\varphi)$ $(\exists \varphi)$ occurs)\zl ;\zd\ similarly here\zl ,\zd\ \ul \zl s\zd o these form\zl ulas\zd\ are really prop\zl ositional\zd\ f\zl u\zd nct\zl ions\zd\ since they contain free var\zl iables.\zd\ \zl \sout{and bec\zl ause\zd\ prop.\zl ositions\zd\ if etc.}\zd\ \ud [And the def\zl inition\zd\ of a taut\zl ology\zd\ was that whatever part.\zl icular\zd\ prop.\zl osition\zd\ or pred.\zl icate\zd\ you subst.\zl itute\zd\ for those free var\zl iables\zd\ of pred\zl i\-cates\zd\ or prop\zl ositions\zd\ you get a true prop\zl osition.\zd\ The var\zl iables\zd\ for ind\zl ividu\-als\zd\ were all bound\zl .\zd ] We can extend the notion of a $\mathbf{\llbracket 56. \rrbracket}$ taut.\zl ology\zd\ also to such expr.\zl essions\zd\ as contain free variables for indiv\zl iduals,\zd\ e.g\zl .\zd\

\begin{tabbing}
\hspace{1.7em} $\varphi(x) \:\vee \sim \varphi(x)$
\end{tabbing}
This is a prop.\zl ositional\zd\ f\zl u\zd nct\zl ion\zd\ containing one free funct\zl ional\zd\ var\zl i\-able\zd\ and one free indiv\zl idual\zd\ variable $x$ and whatever obj\zl ect\zd\ and pred.\zl i\-cate\zd\ you subst\zl itute\zd\ for $\varphi, x$ you get a true prop\zl osition.\zd\ For\zl mula\zd\

\begin{tabbing}
\hspace{1.7em}$(x)\varphi(x) \supset \varphi(y)$
\end{tabbing}
contains $\varphi ,y$ and \ul is \ud univ.\zl ersally\zd\ true bec.\zl ause\zd\ if $M$ is an \ul arb.\zl itrary\zd\ \ud pred.\zl icate and\zd\ $a$ \zl an\zd\ \ul arb.\zl itrary\zd\ \ud ind.\zl ividual\zd\ then

\begin{tabbing}
\hspace{1.7em}$(x)M(x)\supset M(a)$
\end{tabbing}
So in gen.\zl eral\zd\ a tauto\zl logical\zd\ \ul logical formula \ud of the calc.\zl ulus\zd\ of funct.\zl ions\zd\ is a \sout{expr.\zl ession\zd }\, \ul prop.\zl ositional\zd\ f\zl u\zd nct\zl ion\zd\ \ud composed of the above mentioned symbols and which is true whatever part.\zl icular\zd\ $\mathbf{\llbracket 57. \rrbracket}$ objects and predi\-c.\zl ates\zd\ and prop.\zl ositions\zd\ you subst\zl itute\zd\ for free var.\zl iables\zd\ \ul no matter how many ind\zl ividuals\zd\ there exist\zl .\zd\ \ud We can of course express this \ul fact\zl ,\zd\ namely \ud that a cert.\zl ain\zd\ formula is a universally\zl ,\zd\ true by writing quantifiers in front\zl ,\zd\ e.g\zl .\zd\

\begin{tabbing}
\hspace{1.7em}$(\varphi ,x)[\varphi(x) \:\vee \sim \varphi(x)]$
\end{tabbing}
or

\begin{tabbing}
	\hspace{1.7em}$(\varphi, y)[(x)\varphi(x) \supset \varphi(y)]$
\end{tabbing}
\zl unreadable text\zd\ \zl F\zd or the taut\zl ology\zd\ of the calc\zl ulus\zd\ of prop.\zl ositions\zd\

\begin{tabbing}
	\hspace{1.7em}$(p,q)[p \supset p \vee q]$
\end{tabbing}
But it is more convenient to make the convention that univ.\zl ersal\zd\ quantifiers whose scope is the whole expr.\zl ession\zd\ may be left out\zl .\zd\ So if a formula cont.\zl ain\-ing\zd\ free var.\zl iables\zd\ is written down as an assertion\zl ,\zd\ \ul e.g.\ as an axiom or theorem\zl ,\zd\ \ud it means that it holds for everything subst.\zl ituted\zd\ for the \ul free \ud var.\zl iables,\zd\ i.e.\ it means the same thing as if all var.\zl iables\zd\ were bound by quantifiers whose scope is the whole expr\zl ession.\zd\ \ul This \ul convention \ud is in agreement with the way in which \sout{the} theorems are expressed in math.\zl ematics,\zd\ e.g.\ the law of raising a sum to the square is written $(x+y)^2=x^2+2xy+y^2$\zl ,\zd\ i\zl .\zd e.\ with free var.\zl iables\zd\ $x,y$ which express that this holds for any num\zl bers.\zd\ \ud $\mathbf{\llbracket 57.1 \rrbracket}$ \zl This page begins with a crossed out part of a sentence.\zd\ It is also in agreement with our use of \sout{the} variables for propositions in the calc.\zl ulus\zd\ of prop\zl ositions\zd . The axioms and theorems of the prop.\zl ositional\zd\ calc.\zl ulus\zd\ were written with free var.\zl iables,\zd\ \sout{for prop.\zl ositions\zd }\, e.g\zl .\zd\ $p \supset p \vee q$\zl ,\zd\ and \sout{such} a formula like this was understood to mean that it holds for any prop.\zl ositions\zd\ $p,q$\zl .\zd\ \zl The remainder of this page, until the line near the top of p.\ \textbf{58}.\ beginning with ``I hope that'', is crossed out in the manuscript: (So it means what we would have to express by the use of quantifiers by $(p,q)[q \supset p \vee q]$. And in a similar sense we shall also use free variables for pred.\zl icates\zd\ to express that something holds for any arb.\zl itrary\zd\ pred.\zl icate.\zd\ So it is quite\zd\ $\mathbf{\llbracket 58. \rrbracket}$ \sout{natural that we make the same convention.}

I hope that these examples will be sufficient and that I can \ul now \ud begin with setting up the axiomatic system for the calc.\zl ulus\zd\ of pred\zl icates\zd\ \ul which allows to derive all taut.\zl ologies\zd\ of the calc.\zl ulus\zd\ of pred\zl icates\zd . \ud The primit.\zl ive\zd\ notions will be $1.$ the former $\sim ,\vee$ $2.$ the univ\zl ersal\zd\ quant.\zl ifier\zd\ $(x), (y)$\zl .\zd\ The exist\zl ential\zd\ quant\zl ifier\zd\ need not be taken as a primit.\zl ive\zd\ notion because it can be def\zl ined\zd\ in terms of $\sim$ and $(x)$ by $(\exists x)\varphi(x)\equiv \sim (x) \sim \varphi (x)$\zl .\zd\ The form\zl ulas\zd\ of the calc.\zl ulus\zd\ of pred.\zl icates\zd\ will be composed of three kinds of letters\zl :\zd\ $p,q,\ldots$ prop\zl ositional\zd\ var.\zl ia\-bles,\zd\ $\varphi, \psi,\ldots$ \ul functional \ud var\zl iables\zd\ for pred.\zl icates,\zd\ $x,y,\ldots$ var.\zl iables\zd\ for individuals. Furthermore they will contain $\mathbf{\llbracket 59. \rrbracket}$ $(x)\lfloor ,\rfloor (y)\lfloor ,\rfloor \sim \lfloor ,\rfloor \vee$ and the notions defined by those 3\zl three,\zd\ i\zl .\zd e.\ $(\exists x), (\exists y), \supset, \:.\:, \equiv, \mid$ etc. \zl The following text written on the right of p.\ \textbf{59}.\ in the manuscript is numbered \textbf{59}.\textbf{1}, but since the whole of that text is marked in the manuscript for insertion on p.\ \textbf{59}., the number of the page \textbf{59}.\textbf{1}.\ is deleted.\zd\ \ul So the quantifiers apply only to ind.\zl ividual\zd\ var.\zl iables,\zd\ prop\zl ositional\zd\ and funct.\zl ional\zd\ var.\zl iables\zd\ are free\zl ,\zd\ \ul i.e.\ that something holds for all $p,\varphi$ is to be expressed by free var\zl iables\zd\ according to the conv.\zl ention\zd\ mentioned before\zl .\zd\ \ud

So all formulas given as ex.\zl amples\zd\ \ul before \ud are examples for expr.\zl es\-sions\zd\ of the calc\zl ulus\zd\ of funct\zl ions\zd\ but also e\zl .\zd g\zl .\zd\ $(\exists x)\psi(xy)$ \zl and\zd\ $[p\:.\:(\exists x)\psi(xy)]\vee \varphi(y)$ \ul would be ex\zl amples\zd\ \ud etc. I am using the letters $\Phi,\Psi,\Pi$\zl comma from the manuscript deleted\zd\ to denote arbitrary expressions of the calc.\zl ulus\zd\ of pred\zl icates\zd\ and if I wish to ind.\zl icate\zd\ that some var\zl iable\zd\ say $x$ occurs in a form\zl ula\zd\ as a free var\zl iable\zd\ denote the form.\zl ula\zd\ by $\Phi(x) {\text{or}\atop\vee} \Psi(xy)$ \ul if $x,y$ occur both free\zl ,\zd\ \ud which does not exclude that there may be other free var.\zl iables\zd\ bes.\zl ides\zd\ $x$, or $x$ and $y$\zl ,\zd\ in the form\zl ula\zd . \ud

The axioms are like this:
\begin{tabbing}
\hspace{1.7em}\=I.\hspace{1.5em}\= The four ax.\zl ioms\zd\ of the calc.\zl ulus\zd\ of prop.\zl ositions\zd \\[1ex]
\>\>\hspace{2em}\=$p \supset p \vee q$\hspace{2em}\=$p \vee q \supset q \vee p$\\[.5ex]
\>\>\>$p \vee p \supset p$\>$(p \supset q)\supset (r \vee p \supset r\vee q)$\\[1ex]
\>II\zl .\zd \>One specific ax.\zl iom\zd\ for the univ.\zl ersal\zd\ quantifier \\[1ex]
\>\>\> \zl Ax.\ 5\zd \qquad $(x)\varphi(x)\supset \varphi(y)$
\end{tabbing}
This is the formula mentioned before which says: ,,\zl ``\zd For any $y$\zl ,\zd\ \ul $\varphi$ \ud it is true that if $\varphi$ holds for every $x$ then it holds for $y$''\zl .\zd\

These are all ax.\zl ioms\zd\ which we need. [\sout{They are expressed by using free var.\zl iables\zd\ $p, \varphi, y$ in the sense just discussed.}] The rules of inf\zl erence\zd\ are the following $4$\zl four:\zd\

\vspace{1ex}

\noindent$\mathbf{\llbracket 60. \rrbracket}$
\begin{itemize}
\item[1.\zl 1\zd] The rule of impl\zl ication which reads exactly as for the calc.\zl ulus\zd\ of prop\zl o\-sitions:\zd\ If $\Phi,\Psi$ are any expr.\zl essions\zd\ then from $\Phi,\Phi \supset \Psi$ you can conclude \zl $\Psi$\zd \zl .\zd\
\end{itemize}
The only diff\zl erence\zd\ is that now $\Phi,\Psi$ are expr.\zl essions\zd\ which may involve quantifiers and funct\zl ional\zd\ var.\zl iables\zd\ and individual var.\zl iables\zd\ in add.\zl ition\zd\ to the symb\zl ols\zd\ occuring in the calc.\zl ulus\zd\ of prop\zl ositions\zd .
\begin{tabbing}
\ul So e\zl .\zd g\zl .\zd\ \= from $[p \vee (x)[\varphi(x)\supset \varphi(x)]]\supset \varphi(y) \: \vee \sim \varphi (y)$\\[.5ex]
\>\underline{and $[p \vee (x)[\varphi (x)\supset \varphi (x)]]$}\\[.5ex]
\>concl.\zl ude\zd\ $\varphi (y)\: \vee \sim \varphi (y)$ \ud
\end{tabbing}
\begin{itemize}
\item[2.\zl 2\zd] The rule of Subst.\zl itution\zd\ which has now 3\zl three\zd\ parts (accord.\zl ing\zd\ to the 3\zl three\zd\ kinds of var.\zl iables\zd )\zl :\zd\
\begin{itemize}
\item[] \sout{1. For prop.\zl ositional\zd\ var.\zl iables\zd\ $p,q$ any expr.\zl ession may be subst.\zl i\-tuted\zd }
\item[1.\zl a)\zd ] For ind.\zl ividual\zd\ var\zl iables\zd\ $x,y$ \ul bound or free \ud any other ind.\zl ividu\-al\zd\ var\zl iable\zd\ may be subst\zl ituted\zd\ as long as our conventions \ul about the not.\zl ion\zd\ of free var.\zl iables\zd\ \ud are observed\zl ,\zd\ i.e.\ bound variable \sout{are} whose scopes do not ly\zl ie\zd\ outside of each other must be denoted by diff.\zl erent\zd\ letters and \sout{that} all free variables must be denoted by letters different from all bound var.\zl iables\zd\ -- [Rule \ul of \ud renaming the ind\zl ividual\zd\ variables.]\zl ].\zd\
\end{itemize}
\end{itemize}
\vspace{-2ex}
$\mathbf{\llbracket 61. \rrbracket}$
\begin{itemize}
\item[]
\begin{itemize}
\item[2.\zl b)\zd ] For a prop.\zl ositional\zd\ var\zl iable\zd\ any expre\zl ession\zd\ may be subst\-\zl ituted\zd\ \ul with a cert\zl ain\zd\ restriction form\zl ulated\zd\ later\zl .\zd\
\item[3.\zl c)\zd ] If you have \sout{the} an expr\zl ession\zd\ $\Pi$ \zl \sout{(e.g\ldots)}\zd\ and $\varphi$ a prop.\zl osition\-al\zd\ \zl functional\zd\ variable occurring \zl in\zd\ $\Pi$ perhaps on sev.\zl eral\zd\ places and with diff.\zl erent\zd\ arg\zl uments\zd\ $\varphi (x)$\zl ,\zd\ $\varphi(y)$\zl ,\zd\ \ldots\ and if $\Phi(x)$ is an expr.\zl ession\zd\ containing $x$ free then you may subs.\zl ti\-tute\zd\ $\Phi(x)$ for $\varphi(x)$\zl ,\zd\ $\Phi(y)$ for $\varphi(y)$ etc\zl .\zd\ simultaneously in all places wher\zl e\zd\ $\varphi$ occurs. Similarly for $\varphi(xy)$ \zl and\zd\ $\Phi(xy)$\zl \sout{\ul it concerns the letters by which the\ldots\ \ud} \zd\
\end{itemize}
\end{itemize}
\zl The following text on the rest of this page is crossed out in the manuscript:

Take e.g.\ $(x)\varphi(x) \supset \varphi(y)$ and consider the expr.\zl ession\zd\ $(\exists z)\psi(zx)$ which is a prop.\zl ositional\zd\ f\zl u\zd nct\zl tion\zd\ with one free ind.\zl ividual\zd\ var\zl iable\zd . If we subst\zl itute\zd\ this expr\zl ession\zd\ for $\varphi$ of the first expr.\zl ession\zd

In all those three rules of subst\zl itution\zd\ we have only to be careful about one thing which may be expr.\zl essed\zd\ roughly speaking by saying\zl :\zd\ The bound variables must not get mixed up. But\zd

$\mathbf{\llbracket 61.1 \rrbracket}$ It is clear that this is a correct inf\zl erence,\zd\ i\zl.\zd e\zl.\zd\ gives a taut\zl olo\-gy\zd\ if the formula in which we subst\zl itute\zd\ is a taut\zl ology,\zd\ bec.\zl ause\zd\ if a form.\zl ula\zd\ \ul is \ud \zl a\zd\ taut\zl ology\zd\ that means that it holds for any propert\zl y\zd\ or rel.\zl ation\zd\ $\varphi, \psi$\zl ,\zd\ but \ul any \ud prop\zl ositional\zd\ f\zl u\zd nct\zl ion\zd\ with one or several free var.\zl iables\zd\ defines a cert\zl ain\zd\ prop\zl erty\zd\ or rel.\zl ation;\zd\ therefore the form\zl ula\zd\ must hold for them. \ul Take e.g.\ the taut.\zl ology\zd\ \ud $(x)\varphi(x)\supset \varphi(y)$ and subst\zl itute\zd\ for $\varphi$ the expr\zl ession\zd\ $(\exists z)\psi(zx)$ \ul which has one free ind.\zl ividual\zd\ variable \ud . Now the last form.\zl ula\zd\ says that for every prop.\zl erty\zd\ $\varphi$ and any ind\zl ividual\zd\ $y$ we have: ,,\zl ``\zd If for any $x$ $\varphi(x)$ then $\varphi(y)$''\zl .\zd\ \sout{Since this holds for any prop.\zl erty\zd\ \ul $\varphi$ \ud }. But if $\psi$ is an arb.\zl itrary\zd\ rel.\zl ation\zd\ then $(\exists z)\psi(zx)$ defines a cert\zl ain\zd\ prop.\zl erty\zd\ bec.\zl ause\zd\ it is a prop\zl ositional\zd\ f\zl u\zd nct\zl ion\zd\ with one free var\zl iable\zd\ $x$. Hence the ab.\zl ove\zd\ form.\zl ula\zd\ must hold also for this prop\zl erty,\zd\ i.e.\ we have: If for every object $(x)[(\exists z)\psi(zx)]$ then also for $y$ \zl \sout{$\supset$}\zd $(\exists z)\psi(zy)$ and that will be true whatever the rel.\zl ation\zd\ $\psi$ \ul and the object $y$ \ud may be\zl ,\zd\ i\zl .\zd e.\ it is again a taut\zl ology\zd .

$\mathbf{\llbracket 62. \rrbracket}$ $\Big\lceil$You see in this process of subst.\zl itution\zd\ we have sometimes to change the free variables\zl,\zd\ like\zl as\zd\ here we have to change $x$ into $y$ bec.\zl ause\zd\ the $\varphi$ occurs with the var\zl iable\zd\ $y$ here\zl ;\zd\ if the $\varphi$ occurred with the var.\zl iable\zd\ $u$ $\varphi(u)$ we would have to subst.\zl itute\zd\ $(\exists z)\psi(zu)$ in this place.$\Big\rfloor$ In this ex.\zl ample\zd\ we subst.\zl ituted\zd\ an expr.\zl ession\zd\ cont.\zl aining\zd\ $x$ as \ul \zl the\zd\ only free var.\zl iable,\zd\ but \ud we can subst\zl itute\zd\ for $\varphi(x)$ here also an expr\zl ession\zd\ which contains other free \ul ind.\zl ividual\zd\ variables besides $x$\zl ,\zd\ \sout{and} i\zl .\zd e.\ \ul also in this case we shall obtain a taut\zl ology\zd . Take e.g.\ the expr.\zl ession\zd\ $(\exists z)\chi(zxu)$. This is a prop.\zl ositional\zd\ funct\zl ion\zd\ with the free ind.\zl ividual\zd\ var.\zl iable\zd\ $x$ but it has the free ind.\zl ividual\zd\ var\zl iable\zd\ $u$ in addition. Now if we replace $\chi$ by a spec.\zl ial\zd\ triadic rel.\zl ation\zd\ $R$ and $u$ by a spec.\zl ial\zd\ obj\zl ect\zd\ \underline{a} then $(\exists z)R(zxa)$ is a prop.\zl ositional\zd\ f\zl u\zd nct\zl ion\zd\ with one free var.\zl iable\zd\ $x$\zl ;\zd\ hence

\zl As indicated by ``63.1'' at the bottom on the right of this page, the sentence interrupted here is continued on p.\ \textbf{63}.\textbf{1}, after the last sentence on this page which is crossed out, and the entirely crossed out p.\ \textbf{63}, which together make the following text: Therefore \zl unreadable text\zd\ rel\zl ation\zd\ between $x,u$ but if we replace $u$ by $\mathbf{\llbracket 63. \rrbracket}$ an individ.\zl ual\zd\ obj.\zl ect\zd\ say $a$ then $(\exists z)\chi(zxa)$ is now a prop.\zl ositional\zd\ f\zl u\zd nct\zl ion\zd\ with one free var.\zl iable\zd\ $x$\zl ,\zd\ i.e.\ defines a cert.\zl ain\zd\ property of $x$. Therefore we can substitute it for $\varphi$ in the above taut\zl ology\zd\ and obtain
$$(x)[(\exists z)\chi(zxa)]\supset (\exists z)\chi(zya)$$
But now this will be correct whatever the obj\zl ect\zd\ $a$ may be\zl ,\zd\ i.e.\ we can replace $a$ by a variable $u$ and obt.\zl ain\zd\
$$(x)[(\exists z)\chi(zxu)]\supset (\exists z)\chi(zyu)$$
and this will be a taut.\zl ology,\zd\ i\zl .\zd e.\ true whatever $u,y,\chi$ may be. So the rule of subst.\zl itution\zd\ is to be understood to mean for $\varphi(x)$ one can subst.\zl itute\zd\ an expr\zl ession\zd\ containing at least the free var\zl iable\zd\ $x$ but\zd\

$\mathbf{\llbracket 63.1 \rrbracket}$ \zl it\zd\ defines a cert.\zl ain\zd\ prop.\zl erty,\zd\ hence the above form\zl ula\zd\ holds\zl ,\zd\ i\zl .\zd e\zl .\zd\
\begin{tabbing}
\hspace{1.7em}$(x)(\exists z)R(zxa)\supset (\exists z)R(zya)$
\end{tabbing}
whatever $y$ may be\zl ,\zd\ but this will be true whatever $R,a$ may be\zl ;\zd\ therefore if we replace them by var\zl iables\zd\ \ul $\chi$, $u$ \ud the form\zl ula\zd\ obtained:
\begin{tabbing}
\hspace{1.7em}$(x)(\exists z)\chi(zxu)\supset (\exists z)\chi(zyu)$
\end{tabbing}
\sout{and this} will be true for any $\chi,u,y$\zl ,\zd\ i\zl .\zd e.\ it is a taut\zl ology\zd . So the rule of subst.\zl itution\zd\ is also correct for expr.\zl essions\zd\ containing add.\zl itional\zd\ free var.\zl iables\zd\ $u$, and therefore this $\Phi(x)$ is to mean an expr.\zl ession\zd\ containing \ul the free var\zl iable\zd\ \ud  $x$ but perhaps some other free var.\zl iables\zd\ in addition.

$\mathbf{\llbracket 64.\rrbracket}$ Examples for the other two rules of subst.\zl itution:\zd\
\begin{tabbing}
\hspace{1.7em}\=\zl \sout{1.}\zd\ \zl F\zd or prop.\zl ositional\zd\ var.\zl iable\zd\ \\[1ex]
\>$p\:.\:(x)\varphi(x) \equiv (x)[p\:.\:\varphi(x)]$
\end{tabbing}
subst\zl itute\zd\ $(\exists z)\psi(z)$. Since this holds for every prop\zl osition\zd\ it holds also for $(\exists z)\psi(z)$ which is a prop\zl osition\zd\ if $\psi$ is any arb.\zl itrary\zd\ pred.\zl icate.\zd\ \sout{\zl unreadable word\zd }\, Hence we have for any pred.\zl icates\zd\ $\psi,\varphi$
\begin{tabbing}
\hspace{1.7em}$(\exists z)\psi(z)\:.\:(x)\varphi(x)\equiv (x)[(\exists z)\psi(z)\:.\: \varphi(x)]$
\end{tabbing}
But we are also allowed to subst\zl itute\zd\ expr.\zl essions\zd\ containing free var.\zl ia\-bles\zd\ and prop.\zl ositional\zd\ var.\zl iables\zd\ e.g\zl .\zd\ \sout{$p \supset$} $(z)\chi(zu)$ (free var\zl iable\zd\ $u$) bec.\zl ause\zd\ if \ul you \ud take \ul for \ud $u$ \sout{be} any ind.\zl ividual\zd\ obj.\zl ect\zd\ \ul $a$ \ud [and $p$ any indiv.\zl idual\zd\ prop\zl osition\zd\ \ul $\pi$ \ud] and $\chi$ any rel.\zl ation\zd\ \ul $R$ \ud then $\mathbf{\llbracket 65.\rrbracket}$ this will be a prop\zl osition\zd . \sout{hence} And $p\:.\:(x)\varphi(x) \lfloor \equiv\rfloor (x)[p\:.\:\varphi(x)]$ holds for any prop\zl osition\zd . So it will also hold for this\zl ,\zd\ i\zl .\zd e.\
\begin{tabbing}
\hspace{1.7em}$[$\sout{$p \supset$}$(z)\chi(zu)]\:.\:(x)\varphi(x) \equiv (x)[$\sout{$p \supset$}$(z)\chi(zu)\:.\:\varphi(x)]$
\end{tabbing}
will be true whatever $p,\chi,\varphi,u$ may be\zl ,\zd\ i\zl .\zd e\zl .\zd\ a tautology.

Finally an example for subst.\zl itution\zd\ of ind.\zl ividual\zd\ var\zl iables:\zd
\begin{itemize}
\item[\zl \sout{1.}\zd ] For a bound $(x)\varphi(x)\supset \varphi(y)$\; \zl :\zd\ \; $(z)\varphi(z)\supset \varphi(y)$. So this inf.\zl erence\zd\ merely brings out the fact that the notation of bound variables is arb\zl itrary\zd .
\item[\zl \sout{2.}\zd ] The rule of subst.\zl itution\zd\ applied for free var.\zl iables\zd\ is more essen\-tial\zl ;\zd\ e.g.\ \zl f\zd rom $(x\lfloor ,\rfloor y)\varphi(xy)\supset \varphi(uv)$ we can conclude $(x\lfloor ,\rfloor y)\varphi(xy)$ $\supset\varphi(uu)$ \ul by subst.\zl ituting\zd\ $u$ for $v$. This is an all.\zl owable\zd\ subst.\zl itu\-tion\zd\ because the variable which you subst.\zl itute,\zd\ $u$\zl ,\zd\ does\zl \,\zd not occur as a bound var\zl iable\zd . It occurs as a free var\zl iable\zd\ but that does\zl \,\zd not matter\zl .\zd
\end{itemize}

Of course if a var.\zl iable\zd\ occurs in sev.\zl eral\zd\ places it has to be replaced by the same other var\zl iable\zd\ $\mathbf{\llbracket 66.\rrbracket}$ in all places where it occurs. In the rule of subst\zl itution\zd\ for prop.\zl ositional\zd\ and functional \ul variable there is one restriction to be made as I mentioned before, namely one has \ud to be careful about the letters which \sout{we} \ul one \ud uses for the bound variables\zl ,\zd\ e.g.\

\begin{tabbing}
\hspace{1.7em}$(\exists x)[p\:.\:\varphi(x)]\:.\:(x)\varphi(x)\supset (x)[p\:.\:\varphi(x)]$
\end{tabbing}
\ul is a tautol\zl ogy\zd . \ud Here we cannot subst\zl itute\zd\ $\psi(x)$ for $p$ bec\zl ause\zd\ \zl \sout{ie.}\zd\

\begin{tabbing}
\hspace{1.7em}$(\exists x)[\psi(x)\:.\:\varphi(x)]\:.\:(x)\varphi(x)\supset (x)[\psi(x)\:.\:\varphi(x)]$
\end{tabbing}
is not a tautology\zl ,\zd\ \sout{e.g.\ we cannot subst.\zl itute\zd\ here for $p$ the expr.\zl ession\zd\ $\psi(x)$ i\zl .\zd e.\ $\psi(x)\:.\:\varphi(x)\equiv (x)[\psi(x)\:.\:\varphi(x)]$} \zl \sout{is not a tautology}\zd\ bec.\zl ause\zd\ here the expr.\zl ession\zd\ which we subst.\zl ituted\zd\ contains a var\zl iable\zd\ $x$ which is bound in the expr\zl ession\zd\ in which we substitute\zl .\zd\ \ul Reason\zl :\zd\ This form\zl ula\zd\ holds for any prop.\zl osition\zd\ $p$ but not for any prop.\zl ositional\zd\ f\zl u\zd nct.\zl ion\zd\ with the free var.\zl iable\zd\ $x$\zl .\zd\ \zl Before the next sentence a horizontal line is drawn in the manu\-script.\zd\ Now if we subst\zl itute\zd\ for $p$ an expr.\zl ession\zd\ $\Phi$ containing perhaps free var\zl iables\zd\ $y,z,\ldots$ (but not the free var\zl iable\zd\ $x$) then $y,z$ will be free in the whole expr\zl ession\zd . Therefore if $y,z,\ldots$ are replaced by definite things then $\Phi$ will bec.\zl ome\zd\ a prop.\zl osition\zd\ bec\zl ause\zd\ then all free var\zl iables\zd\ con\zl tained\zd\ in it are repl.\zl aced\zd\ by def.\zl i\-nite\zd\ obj\zl ects\zd . \zl After the preceding sentence a horizontal line is drawn in the manuscript.\zd\

Therefore the expr\zl ession\zd\ to be subst\zl ituted\zd\ must not contain $x$ as a free var.\zl iable\zd\ because it would play the role of a prop.\zl ositional\zd\ f\zl u\zd nc\-t\zl ion\zd\ and not of a prop\zl osition\zd . In order to avoid such \ud \ul occurrences \ud we have to make in the rule of subst\zl itution\zd\ the \sout{further} stipulation that the expr.\zl ession\zd\ to be subst\zl ituted\zd\ should contain no variable $\mathbf{\llbracket 67.\rrbracket}$ (bound or free) which occurs in the expr.\zl ession\zd\ in which we substitute bound or free\zl ,\zd\ exc.\zl luding\zd\ \zl \sout{[}\zd of course the variable $x$ here\zl \sout{]}\zd \zl .\zd\ If you add this restriction you obtain the formulation of the rule of subst\zl itution\zd\ which you have in your notes that were distributed.

\zl The following text is crossed out in the manuscript: which are identified with $x$ \zl unreadable word\zd\ $\varphi(x)$. But besides these the expr.\zl ession\zd\ should contain no var.\zl iable\zd\ which occurs in the expr.\zl ession\zd\ in which we subst.\zl itute\zd\ \ul So this restriction has to be added to the rule of subst\zl itution\zd . \ud So the final form of the rule of subst\zl itution\zd\ is as follows:\zd\

So far I formulated two rules of inf\zl erence\zd\ (impl\zl ication,\zd\ subst.\zl itution\zd ). The third is \zl displayed with number 3\zd\ the rule of defined symb\zl ol\zd\ which reads:
\begin{itemize}
\item [1.] For any expre\zl essions\zd\ $\Phi,\Psi$\zl ,\zd\ $\Phi \supset \Psi$ may be repl.\zl aced\zd\ by $\sim \Phi \vee \Psi$ and similarly for $.$ \zl and\zd\ $\equiv$.
\end{itemize}

\vspace{-2ex}

\noindent $\mathbf{\llbracket 68.\rrbracket}$

\vspace{-2ex}

\begin{itemize}
\item[2.] $(\exists x)\Phi(x)$ may be repl.\zl aced\zd\ by $\sim (x) \sim \Phi(x)$ \ul and vice versa \ud where $\Phi(x)$ is any expr.\zl ession\zd\ containing the free var.\zl iable\zd\ $x$\zl .\zd\ (So that means that the exist.\zl ential\zd\ quantifier is def.\zl ined\zd\ by means of the univ.\zl ersal\zd\ quant.\zl ifier\zd\ in our syst\zl em\zd .)
\end{itemize}

\sout{\zl unreadable word\zd }\, \zl T\zd he three rules of inf\zl erence\zd\ ment.\zl ioned\zd\ so far (impl.\zl i\-cation\zd , subst\zl itution,\zd\ def.\zl ined\zd\ symb\zl ol\zd ) corresp\zl ond\zd\ exactly to the three rules of inf.\zl erence\zd\ which we had in the calc.\zl ulus\zd\ of prop\zl ositions\zd . Now we set up a fourth one which is specific for the univ\zl ersal\zd\ quantifier\zl ,\zd\ namely:

\begin{itemize}
\item [4.\zl 4\zd ] \ul Rule of the universal quantifier\zl :\zd\ \ud From $\Pi\supset \Phi(x)$\zl ,\zd\ if $\Pi$ does not contain $x$ as a free var.\zl iable\zd\ we can conclude $\mathbf{\llbracket 69.\rrbracket}$ $\Pi \supset (x)\Phi(x)$.
\end{itemize}

That this inf.\zl erence\zd\ is correct can be seen like this: Assume $\pi$ is a definite propos\zl ition\zd\ and $M(x)$ a \sout{\zl unreadable word\zd }\, \ul definite \ud prop.\zl ositional\zd\ f\zl u\zd nct\zl ion\zd\ with \ul exactly one free var.\zl iable\zd\ $x$ and let us assume we know: \underline{$\pi \supset M(x)$ \zl colon deleted\zd\ holds for every $x$}\zl .\zd\ Then I say we can conclude: $\pi\supset (x)M(x)$\zl .\zd\ For 1. \zl i\zd f $\pi$ is false the concl\zl usion\zd\ holds\zl ,\zd\ 2. if $\pi$ is true then by ass\zl umption\zd\ $M(x)$ is true for every $x$\zl ,\zd\ i\zl .\zd e\zl .\zd\ $(x)M(x)$ is true\zl ;\zd\ hence the conclusion again holds bec.\zl ause\zd\ it is an impl.\zl ication\zd\ both terms of which are true\zl .\zd

\zl The following text is crossed out in the manuscript: $\pi \supset M(x)$ reason: For every obj\zl ect\zd\ $x$ it is true that: If $\pi$ then $x$ has the prop.\zl erty\zd\ \ul def\zl ined\zd\ by $M$ \ud \sout{\zl unreadable symbol\zd }\, But then it follows: If $\pi$ is true then every obj.\zl ect\zd\ has the prop.\zl erty\zd\ $M$ i\zl .\zd e.\ $\pi\supset (x)M(x)$ bec\zl ause\zd\ assume $\pi$ is true then owing to this $M(x)$ is true whatever $x$ may \zl be\zd\ bec\zl ause\zd\ \ul this impl.\zl ication\zd\ holds \ud i\zl .\zd e.\ $(x)M(x)$ is true\zd\

So we have proved that in any case $\pi\supset (x)M(x)$ \ul is true if $\pi\supset M(x)$ is true for every $x$ \ud . But from this consid.\zl eration\zd\ about a part.\zl icular\zd\ prop.\zl osition\zd\ $\pi$ and a part.\zl icular\zd\ prop.\zl ositional\zd\ $\mathbf{\llbracket 70.\rrbracket}$ \ul f\zl u\zd nct\zl ion\zd\ with one free var\zl iable\zd\ \ud $M(x)$ it follows that the above rule of inf.\zl erence\zd\ yields tautologies if applied to tautologies. Bec.\zl ause\zd \zl a\zd ssume $\Pi\supset \Phi(x)$ is a taut\zl ology.\zd\ Now then $\Pi$ will cont.\zl ain\zd\ some free var.\zl iables\zd\ for prop\zl osi\-tions\zd\ $p, q,\ldots$ for fu\zl nctions\zd\ $\varphi, \psi,\ldots$ and for ind.\zl ividuals\zd\ $y, z,\ldots$ ($x$ does not occur among them) and $\Phi(x)$ will also contain \ul free \ud var.\zl iables\zd\ $p, q,\ldots$\zl ,\zd\ $\varphi, \psi, \ldots$ and \ul free \ud var.\zl iables\zd\ for ind\zl ividuals\zd\ \underline{$x$}, $y, z$ ($x$ among them). Now if you subst\zl itute\zd\ def.\zl inite\zd\ prop\zl osi\-tions\zd\ for $p,q$\zl ,\zd\ def.\zl inite\zd\ pred\zl icates\zd\ for $\varphi, \psi$ and def.\zl inite\zd\ obj.\zl ects\zd\ for $y,z,\ldots$ but leave $x$ w.\zl here\zd\ it stands then $\mathbf{\llbracket 71.\rrbracket}$ by this subst.\zl itution\zd\ all free var.\zl iables\zd\ of $\Pi$ are replaced by indiv.\zl idual\zd\ objects, hence $\Pi$ becomes a definite \sout{assertion} prop.\zl osition\zd\ $\pi$ and all free var.\zl iables\zd\ of $\Phi$ exc.\zl luding\zd\ $x$ are repl.\zl aced\zd\ by obj.\zl ects;\zd\ hence $\Phi(x)$ becomes a prop.\zl ositional\zd\ f\zl u\zd nct\zl ion\zd\ with one free var.\zl iable\zd\ $M(x)$ \sout{which defines a cert.\zl ain\zd\ monadic predicate $M$} and we know $\pi \supset M(x)$ is true for any obj.\zl ect\zd\ $x$ bec.\zl ause\zd\ \sout{the} it is obt.\zl ained\zd\ by subst\zl itution\zd\ of indiv.\zl idual\zd pred\zl icates\zd , prop\zl ositions and\zd\ obj\zl ects\zd\ in a taut\zl ology\zd . But then \ul as we have just seen under this ass\zl umption\zd\ $\pi \supset (x)M(x)$ is true. But this argum\zl ent\zd\ applies whatever part.\zl icular\zd\ pred.\zl icate,\zd\ $\mathbf{\llbracket 72.\rrbracket}$ prop.\zl osition\zd\ etc\zl .\zd\ we subst.\zl itute;\zd\ always the result $\pi \supset (x)M(x)$ is true\zl ,\zd\ i\zl .\zd e.\ $\Pi\supset (x)\Phi(x)$ is a taut\zl ology\zd . \ul \zl T\zd his rule of course is meant \ul to apply \ud to any other ind.\zl ividual\zd\ var.\zl iable\zd\ $y,z$ instead of $x$ \zl .\zd\ \ud So these are the axioms and rules of inf.\zl erence\zd\ of which one can prove that they are complete: i.e\zl .\zd\ every taut.\zl ology\zd\ of the cal\zl culus\zd\ of f\zl u\zd nct\zl ions\zd\ can be derived\zl .\zd\ \zl Here one finds in the manuscript an insertion sign to which no text to be inserted corresponds, and the page is divided by a sinuous horizontal line.\zd\

Now I want to give some examples \ul for derivations from these ax\zl ioms\zd . Again an expression will be called demonstrable or derivable if it can be obtained from Ax\zl ioms\zd\ $1\ldots 5$ \zl (1)\ldots(4) and Ax.\ 5\zd\ by rules 1 -- 4. \ud First of all I wish to remark that, since among our ax.\zl ioms\zd\ and rules all ax\zl ioms\zd\ and rules of the calc\zl ulus\zd\ of prop.\zl ositions\zd\ occur, we can derive from our ax.\zl ioms\zd\ and rules all formulas and rules which we formerly derived in the calc.\zl ulus\zd\ of prop\zl ositions\zd . \ul But \ul the rules \zl are\zd\ now \ud formulated \zl \sout{now}\zd\ for \sout{the} all expr.\zl essions\zd\ of the calc.\zl ulus\zd\ of pred.\zl icates,\zd\ e.g.\ \ul if $\Phi,\Psi$ \ud \zl are such expressions\zd\
\begin{tabbing}
\hspace{15em}\=$\Phi \supset \Psi$\\[.5ex]
\> \underline{$\Psi \supset \Pi$}\\[.5ex]
\> $\Phi \supset \Pi$ \ud
\end{tabbing}

So we are justified to use them in the subsequ.\zl ent\zd\ $\mathbf{\llbracket 73.\rrbracket}$ derivations. At first I mention some further rules of \zl the\zd\ calc.\zl ulus\zd\ \ul of prop.\zl ositions\zd\ \ud which I shall need:
\begin{tabbing}
\hspace{1.7em}\=1. \hspace{2em}\=$P \equiv Q \qquad \, : \qquad$ \=$P \supset Q,\; Q \supset P$ \;\; and vice versa \\[.5ex]
\>2.\>$P \equiv Q \qquad \lfloor :\rfloor \qquad$ \>$\sim P \equiv \: \sim Q$ \\[.5ex]
\>$1'$.\> $p \equiv\: \sim \sim p$ \> ($2'$. $p \equiv p$)\\[.5ex]
\>$3'$.\> $(p \supset q)\:.\: p \supset q$ \> $(p \supset q)\supset (p \supset q)$ \;\; Import.\zl ation\zd
\end{tabbing}
\begin{tabbing}
\hspace{1.7em}\=1. \hspace{2em}\=\underline{$\varphi(y) \supset (\exists x)\varphi(x)$} \\[.5ex]
\>\>$(x)[\sim \varphi(x)] \supset\: \sim \varphi(y)$ \qquad \= Subst.\zl itution,\zd\ Ax\zl . \zd 5 \\[.5ex]
\>\>$\varphi(y)\supset\: \sim (x)[\sim \varphi(x)]$ \>Transp.\zl osition\zd\ \quad $\f{\sim \varphi(x)}{\varphi(x)}$ \\
\>\>$\varphi(y) \supset(\exists x)\varphi(x)$ \> \underline{def.\zl ined\zd\ symb.\zl ol\zd }\\[2ex]
\>2. \> \underline{$(x)\varphi(x) \supset (\exists x)\varphi(x)$}\\[.5ex]
\>\>$(x)\varphi(x) \supset \varphi(y)$ \> Ax.\ 5 \\[.5ex]
\>\>$\varphi(y)\supset (\exists x)\varphi(x)$ \> 1\zl .\zd\
\end{tabbing}
\zl The next page of the manuscript is not numbered and contains only the following heading:
\begin{center}
\underline{Log.\zl ik\zd\ Vorl.\zl esungen\zd\ \zl German: Logic Lectures\zd\ Notre Dame}\\
\underline{1939}
\end{center}
This page and the pages following it up to p.\ \textbf{73}.\textbf{7}, which makes nine pages, are on loose, torn out, leafs, with holes for a spiral, but not bound with the spiral to the rest of the notebook, as the other pages in this Notebook~V are. In all of the notebooks the only other loose leafs are to be found at the end of Notebooks III and VII.\zd\

$\mathbf{\llbracket 73.1\rrbracket}$ Last time I set up a system of axioms and rules of inf.\zl erence\zd\ from which it is possible to derive all tautologies of the calc.\zl ulus\zd\ of predicates. Incidentally I wish to mention that \zl the\zd\ technical term tautology is somewhat out of \sout{use} \ul fashion \ud at present\zl ,\zd\ the word analytical (which goes back to Kant) is used in it's \zl its\zd\ place, and that has certain advantages because analytical is an indifferent term whereas the term tautological suggests a certain philosophy of logic\zl ,\zd\ namely \ul the theory \ud that the propositions \ul of logic \ud are in some sense void of content\zl ,\zd\ that they say nothing\zl .\zd\ Of course it is by no means necessary for a $\mathbf{\llbracket 73.2\rrbracket}$ mathematical logician to adopt this theory, bec.\zl ause\zd\ math.\zl ematical\zd\ logic is a purely math.\zl ematical\zd\ theory which is wholly indiff.\zl erent\zd\ towards \zl any\zd\ phil.\zl osophical\zd\ question. So if I use this term tautological I don't want to imply \ul by that \ud any definite standpoint as to the essence of logic\zl ,\zd\ but \ul the term taut.\zl ological\zd\ \ud is only to be understood as a shorter expr.\zl ession\zd\ for universally true. Now as to our axiomatic syst.\zl em\zd\ the Axioms were as follows 1. 2. Rules of inf\zl erence\zd
\begin{tabbing}
\hspace{1.7em}\= 1.\zl 1\zd\ \hspace{.5em}\= Implic\zl ation\zd\ \quad\= $\Phi, \Phi \supset \Psi \quad : \quad \Psi$\\[1ex]
\>2.\zl 2\zd\ \> Subst.\zl itution\zd\hspace{1em}\> \zl a\zd )\hspace{1.3em}\= indiv.\zl idual\zd\ var\zl iables\zd\ \\[.5ex]
\>\>\>b)\> prop.\zl ositional\zd\ var.\zl iables\zd \\[.5ex]
\>\>\> d.\zl c\zd )\> funct.\zl ional\zd\ \zl variables\zd\ \\[1ex]
\>3.\zl 3\zd\ \> Rule of def.\zl ined\zd\ symb.\zl ol\zd \\[.5ex]

\>\hspace{3em}1. For $.$\zl ,\zd\ $\supset$\zl ,\zd\ $\equiv$\,\, as formerly\\[.5ex]
\>\hspace{3em}2. $(\exists x)\Phi(x)$ may be repl.\zl aced\zd\ by $\sim (x) \sim \Phi(x)$ and vice versa\\[1ex]

\>4.\zl 4\zd\ \>Rule of the univ\zl ersal\zd\ quantifier\\[.5ex]
\>\>$\Phi\supset \Psi(x) \quad \:.\:.\:\lfloor :\rfloor \quad \Phi \supset (x)\Psi(x)$
\end{tabbing}

$\mathbf{\llbracket 73.3\rrbracket}$ It may seem superfluous to formulate \ul so carefully \ud the stipulations about the letters which we have to use for the bound var.\zl iables\zd\ here in rule 3.\zl 2\zd\ because if you take account of the meaning of the expr.\zl essions\zd\ involved you will observe these rules automatically\zl ,\zd\ because otherwise they would either be ambiguous or not have the intended \ul meaning\zl .\zd\ \ud To this it is to be answered that it is exactly \sout{one of} the \ul chief \ud purpose of the axiomatization of logic \ul to avoid this reference to the meaning of the formulas\zl ,\zd\ i\zl .\zd e.\ we want to \ud \zl \sout{to}\zd\ set up a calculus which can be handled purely mechanically (i\zl .\zd e.\ \ul a calculus \ud which makes thinking superfluous $\mathbf{\llbracket 73.4\rrbracket}$ and which can replace thinking for cert\zl ain\zd\ quest\zl ions\zd )\zl .\zd\

\zl new paragraph\zd\ In other words we want to put into effect as far as possible Leibnitz\zl 's; or perhaps ``Leibnitzian''\zd\ program of a ,,\zl ``\zd calculus ratiocinator'' which he c\zl h\zd aracter\zl izes\zd\ by saying:\zl colon from the manuscript deleted\zd\ \zl that\zd\ \zl h\zd e expects there will be a time in the future when there will be no discussion \ul or reasoning \ud necessary for deciding logical questions but when one will be able \sout{to} simply to say ,,\zl ``\zd calculemus''\zl ,\zd\ \ul let us reckon \ud exactly as in questions of elementary arith\zl metic\zd . This program has been partly carried out by this axiomatic syst\zl em\zd\ \ul for logic \ud . For you \sout{will} see that the rules of inference can be applied $\mathbf{\llbracket 73.5 \rrbracket}$ purely mechanical\zl ly,\zd\ e.g\zl .\zd\ in order to apply the rule of syll.\zl ogism\zd\ \zl comma from the manuscript deleted\zd\ \underline{$\Phi$},\underline{$\Phi$}$\supset \Psi$ you don't have to know what $\Phi$ or $\Psi$ or the sign of impl.\zl ication\zd\ means\zl ,\zd\ but you have only to look at the outward structure of the two prem\zl ises\zd . \zl The following insertion is found in the scanned manuscript on a not numbered page after p.\ \textbf{73}.\textbf{6}.\zd\ \ul All you have to know in order to apply this rule to two premises is that the sec.\zl ond\zd\ premise contains the $\supset$ and that the part preceding the $\supset$ is conform with the first premise. And similar remarks apply to the other axioms\zl .\zd\ \ud

\zl new paragraph\zd\ Therefore \ul as I men\zl tioned\zd\ already \ud it would \ul actual\-ly \ud be possible to construct a machine which would do the following thing: The \ul supposed \ud machine is to have a crank and whenever you turn the crank once around the machine would write \ul down \ud a tautology of the calc\zl ulus\zd\ of predicates and it would write down every \ul existing \ud taut.\zl ology\zd\ of the calc.\zl ulus\zd\ of pred\zl icates\zd\ $\mathbf{\llbracket 73.6 \rrbracket}$ if you turn the crank sufficiently often. So this machine would really replace thinking completely as far as deriving of form\zl ulas\zd\ of the calc.\zl ulus\zd\ of pred\zl icates\zd\ \ul is concerned. \ud It would be a thinking machine in the literal sense of the word.

\zl new paragraph\zd\ For the calculus of prop.\zl ositions\zd\ you can do even more\zl .\zd\ You could construct a machine in \zl the\zd\ form of a typewriter such that if you type down a formula of the calc.\zl ulus\zd\ of prop.\zl ositions\zd\ then the machine would ring a bell \zl if it is a tautology\zd\ and if it is not it would not. You could do the same thing for the calculus \zl The next page of the scanned manuscript, which is not numbered,
contains just an insertion for the text on p.\ \textbf{73}.\textbf{5}, to be found at the appropriate place there.\zd\ $\mathbf{\llbracket 73.7 \rrbracket}$ of monadic pred\zl icates\zd . But one can prove that it is impossible to construct a machine which would do the same thing for the \ul whole \ud calculus of pred\zl icates\zd . So here already one can prove that Leibnitz\zl 's; or perhaps ``Leibnitzian''\zd\ program of the ,,\zl ``\zd calculemus'' cannot be carried through\zl ,\zd\ i\zl .\zd e.\ one knows that the human mind will never be able to be replaced by a machine already for this comparatively simple quest.\zl ion\zd\ to decide whether a form\zl ula\zd\ is a taut.\zl ology\zd\ or not.

\zl The next page of the manuscript, which is not numbered, but is not on a loose leaf as the preceding nine pages in the scanned manuscript are, contains only the following two lines:\zd
\begin{tabbing}
\hspace{3em}\= \sout{$(x)\varphi(x)\supset (\exists x)\varphi(x)$ \qquad Syll.\zl ogism\zd }\\[.5ex]
?$4'$ \> $(p \vee q)\supset (\sim p \supset q) \quad \mid \quad (\sim p \supset q)\supset (p \vee q)$
\end{tabbing}
\begin{tabbing}
$\mathbf{\llbracket 74\rrbracket}$
\hspace{0.9em}\= $(x)\varphi(x) \supset (\exists x)\varphi(x)$ \quad Syll.\zl ogism\zd \\[2ex]
3.\> \underline{$\sim (\exists x)\varphi(x)\equiv (x)\sim \varphi(x)$} \\
\> $\sim \sim (x)\sim \varphi(x) \equiv (x) \sim \varphi(x)$ \quad  $p \equiv \:\sim \sim p$ \quad
$\f{p}{(x) \sim \varphi(x)}$ \zl $\f{(x) \sim \varphi(x)}{p}$ \\
\` fraction bar omitted in the manuscript\\*
\` and arrow pointing to line under 3.\zd\ \\[.5ex]
\>$\sim (\exists x)\varphi(x) \equiv (x)\sim \varphi(x)$ \quad def.\zl ined\zd\ symb\zl ol\zd \\[2ex]
4.\> \underline{$p\:.\: (x)\varphi(x) \equiv (x)[p\:.\: \varphi(x)]$}\\[.5ex]
\> $(x)\varphi(x) \supset \varphi(x)$ \\[.5ex]
\> $p\:.\: (x)\varphi(x) \supset p\:.\: \varphi(y)$ \quad Mult\zl iplication\zd\ from left\\[.5ex]
\> \underline{$p$}$\: .\: (x)\varphi(x)\supset (y)[p\:.\: \varphi(y)]$ \quad \underline{Rule} 4 \quad $\Phi:p\:.\: (x)\varphi(x)$ \quad $\Psi(y):p\:.\: \varphi(y)$\\[.5ex]
\> $(x)[p\:.\: \varphi(x)]\supset p\:.\: \varphi(y)$ \quad Ax.\ 5 \quad Subst.\zl itution\zd\ $\f{p\:.\: \varphi(x)}{\varphi(x)}$ \\[.5ex]
\> $p\:.\: \varphi(y) \supset \varphi(y)$ \quad \=$p\:.\: q \supset q$ \quad $\f{q}{\varphi(y)}$ \zl $\f{\varphi(y)}{q}$ fraction bar omitted in the\\*[-2ex]
\` manuscript\zd\ \\
\> $p\:.\: \varphi(y) \supset p$ \>$p\:.\: q \supset p$\\[.5ex]
\> $(x)[p\:.\: \varphi(x)]\supset \varphi(y)$ \quad \= Syll\zl ogism\zd \\[.5ex]
\> $(x)[p\:.\: \varphi(x)]\supset p$ \> Syll\zl ogism\zd\ \\[.5ex]
$\mathbf{\llbracket 75\rrbracket}$\> $(x)[p\:.\: \varphi(x)]\supset (y)\varphi(y)$ \quad Rule 4\\[.5ex]
\> $(x)[p\:.\: \varphi(x)]\supset p\:.\: (y)\varphi(y)$ \quad Compos.\zl ition\zd\ \\[2ex]
5.?\> \underline{$p \vee (x)\varphi(x) \equiv (x)[p \vee \varphi(x)]$}\\[.5ex]
\> $(x)\varphi(x) \supset \varphi(y)$ \quad Ax\zl .\zd\ 5\\[.5ex]
\> $p \vee (x)\varphi(x)\supset p \vee \varphi(y)$ \quad Add\zl ition\zd\ from left\\[.5ex]
\> \underline{$p$}\,$\vee (x)\varphi(x)\supset (y)[p \vee \varphi(y)]$ \quad Rule 4\\[.5ex]
\> $(x)[p \vee \varphi(x)]\supset p \vee \varphi(y)$ \quad Ax\zl . \zd 5\\[.5ex]
\> $p \vee \varphi(y) \supset (\sim p \supset \varphi(y))$ \quad $p \vee q \supset (\sim p \supset q)$\\[.5ex]
\> $(x)[p \vee \varphi(x)]\supset (\sim p \supset \varphi(y))$ \quad Syll\zl ogism\zd \\[.5ex]
\> $(x)[p \vee \varphi(x)]\:.\: \sim p \supset \varphi(y)$ \quad Imp\zl ortation\zd\ \\[.5ex]
\> $(x)[p \vee \varphi(x)]\:.\: \sim p \supset (y)\varphi(y)$ \quad Rule 4\\[.5ex]
\> $(x)[p \vee \varphi(x)]$\=$\supset [\sim p \supset (y)\varphi(y)]$ \quad Exp.\zl ortation\zd \\[.5ex]
\>\> $\supset [p \vee (y)\varphi(y)]$\\

\zl \sout{\ul 6. \ud }\zd\ \\[1ex]

$\mathbf{\llbracket 76\rrbracket}$\\*[.5ex]

6.\> \underline{$(x)[\varphi(x) \supset \psi(x)]\supset [(x)\varphi(x) \supset (x)\psi(x)]$}
\end{tabbing}

\vspace{-4ex}

\begin{equation*}
\begin{rcases}
(x)[\varphi(x)\supset \psi(x)]\supset [\varphi(y)\supset \psi(y)]\text{\hspace{1.3em}}\\
(x)\varphi(x)\supset \varphi(y)\text{\hspace{5.2em}}
\end{rcases}
\text{\hspace{.2em} Ax}\lfloor . \:\rfloor 5 \quad \f{\varphi(x)\supset \psi(x)}{\varphi(x)}
\end{equation*}

\vspace{-4ex}

\begin{tabbing}
$\mathbf{\llbracket 74\rrbracket}$ \hspace{0.9em}\= \kill

\>$(x)[\varphi(x)\supset \psi(x)]\:.\: (x)\varphi(x) \supset [\varphi(y)\supset \psi(y)]\:.\: \varphi(y)$ \quad Mult.\zl iplication\zd\ \\
\> $[\varphi(y)\supset \psi(y)]\:.\: \varphi(y) \supset \psi(y)$ \quad $(p \supset q)\:.\: p \supset q$ \quad $\f{p}{\varphi(y)}$\;\; $\f{q}{\psi(y)}$ \\*[-1ex]
\` \zl $\f{\varphi(y)}{p}$\;\; $\f{\psi(y)}{q}$\;\; fraction bars omitted in the manuscript\zd\ \\
\> $(x)[\varphi(x)\supset \psi(x)]\:.\: (x)\varphi(x)$\:\=$\supset \psi(y)$ \hspace{3.2em}\= Syll.\zl ogism\zd \\[.5ex]
\>\> $\supset (y)\psi(y)$ \> Rule 4\\[.5ex]
\> $(x)[\varphi(x)\supset \psi(x)]\supset [(x)\varphi(x)\supset (y)\psi(y)]$ \>\> Exp.\zl ortation\zd\ \\[2ex]

7.\> \underline{Derived Rule I}\\*[.5ex]
\> $\Phi(x) \quad \:.\:.\:\lfloor :\rfloor \quad (x)\Phi(x)$\\
\hspace{21.85em}$P \supset Q \quad \lfloor :\rfloor \quad P\:.\: R \supset Q$\\*
\> $p \:\vee\sim p \supset \Phi(x)$ \quad by add.\zl ition\zd\ of premises \quad $Q \quad \lfloor :\rfloor \quad P\supset Q$\\[.5ex]

$\mathbf{\llbracket 77\rrbracket}$\> $p \:\vee\sim p \supset (x)\Phi(x)$ \quad Rule 4\\[.5ex]
\>\underline{$p \:\vee \sim p$}\\[.5ex]
\> $(x)\Phi(x)$ \quad Rule of impl.\zl ication\zd \\[2ex]

8\zl .\zd \> \underline{Derived rule II}\\*[.5ex]
\> $\Phi(x) \supset \Psi(x)\quad : \quad (x)\Phi(x)\supset (x)\Psi(x)$\\[.5ex]
~~~1.\> $(x)[\Phi(x)\supset \Psi(x)]$\\[.5ex]
~~~2.\> Subst\zl itution\zd : $(x)[\Phi(x)\supset \Psi(x)] \supset (x)\Phi(x)\supset (x)\Psi(x)$\\[.5ex]
~~~3\zl .\zd \> Impl\zl ication\zd\ \\[2ex]
?9.\> \underline{Derived rule III}\\*[.5ex]
\>$\Phi(x) \equiv \Psi(x)\quad : \quad$\=$(x)\Phi(x)\equiv (x)\Psi(x)$\\[.5ex]
\>$\Phi(x) \supset \Psi(x)$ \> $(x)\Phi(x)\supset (x)\Psi(x)$\\[.5ex]
\>$\Psi(x) \supset \Phi(x)$ \> \underline{$(x)\Psi(x)\supset (x)\Phi(x)$}\\[.5ex]
\>\> \hspace{3em} $\:\backprime\:\backprime\:\backprime\:$ \\[.5ex]

\noindent$\mathbf{\llbracket 78\rrbracket}$\\*[.5ex]

?10.\> \underline{$\sim (x)\varphi(x) \equiv (\exists x)\sim \varphi(x)$}\\[.5ex]
\> $\varphi(x) \equiv\: \sim \sim \varphi(x)$ \qquad double neg\zl ation\zd\ \\[.5ex]
\> $(x)\varphi(x) \equiv (x)\sim \sim \varphi(x)$ \qquad Rule II\\[.5ex]
\> $\sim (x) \varphi(x)$\:\=$\equiv\: \sim (x) \sim \sim \varphi(x)$ \qquad Transp.\zl osition\zd\ \\[.5ex]
\>\> $\equiv (\exists x)\sim \varphi(x)$ \qquad def.\zl ined\zd\ symb.\zl ol\zd \\[2ex]

\ul ?$10'$.\> \underline{$(x)\varphi(x) \vee (\exists x) \sim \varphi(x)$}\\[.5ex]
\> $(x)\varphi(x) \:\vee \sim (x)\varphi(x)$ \qquad Excl.\zl uded\zd\ middle\\[.5ex]
\> $\sim (x)\varphi(x)\supset (\exists x) \sim \varphi(x)$ \qquad \zl ?\zd 10. \\[.5ex]
\> $\mid$\underline{$(x)\varphi(x) \:\vee \sim (x) \varphi(x)$}$\mid$ $\supset (x)\varphi(x) \vee (\exists x)\sim \varphi(x)$ \qquad Implic.\zl ation\zd\ \ud \\[2ex]

?11.\> \underline{$(x)[\varphi(x)\:.\:\psi(x)]\equiv (x)\varphi(x)\:.\:(x)\psi(x)$}\\[.5ex]
\> $\varphi(x)\:.\:\psi(x)\supset \varphi(x)$\\[.5ex]
\> $(x)[\varphi(x)\:.\:\psi(x)]\supset (x)\varphi(x)$ \qquad Rule II\\*[.5ex]
\> $(x)[\varphi(x)\:.\:\psi(x)]\supset (x)\psi(x)$ \qquad \quad $''$\\[.5ex]
\> $(x)[\varphi(x)\:.\:\psi(x)]\supset (x)\varphi(x)\:.\:(x)\psi(x)$ \qquad Comp.\zl osition\zd\
\end{tabbing}

\vspace{-6ex}

\begin{equation*}
\begin{rcases}
(x)\varphi(x)\supset \varphi(y)\;\;\\
(x)\psi(x)\supset \psi(y)\;\;
\end{rcases}
\text{\hspace{.2em} Ax}\lfloor . \:\rfloor 5 \text{\hspace{15.1em}}
\end{equation*}

\vspace{-4ex}

\begin{tabbing}
$\mathbf{\llbracket 74\rrbracket}$ \hspace{0.9em}\= \kill

\> $(x)\varphi(x)\:.\:(x)\psi(x) \supset \varphi(x)\:.\:\psi(x)$ \qquad Comp\zl osition\zd\ \\[.5ex]

$\mathbf{\llbracket 79\rrbracket}$\> $(x)\varphi(x)\:.\:(x)\psi(x) \supset (x)[\varphi(x)\:.\:\psi(x)]$ \qquad Rule 4.\zl 4\zd\ \\[2ex]

?12.\>\underline{$(x)[\varphi(x)\supset \psi(x)]\:.\:(x)[\psi(x)\supset \chi(x)]\supset (x)[\varphi(x)\supset \chi(x)]$}\quad \zl \sout{\ul g''s \ud }\zd\ \\[.5ex]
*\> $(x)[\varphi(x)\supset \psi(x)]\:.\:(x)[\psi(x)\supset \chi(x)]\supset (x)\{ [\varphi(x)\supset \psi(x)]\:.$\\[.5ex]
\` $[\psi(x)\supset \chi(x)]\}$\quad Subst\zl itution ?\zd 11.\\[.5ex]
\> $[\varphi(x)\supset \psi(x)]\:.\:[\psi(x)\supset \chi(x)]\supset [\varphi(x)\supset \chi(x)]$\quad  Subst.\zl itution\zd\ \\*
\` Syll.\zl ogism\zd \\[.5ex]
**\>$(x)\{ [\varphi(x)\supset \psi(x)]\:.\:[\psi(x)\supset \chi(x)]\} \supset (x)[\varphi(x)\supset \chi(x)]$\` Rule 4.\zl II\zd\ \\[1ex]
\> * and ** \zl with\zd\ Syll\zl ogism\zd\ give the result.\\[.5ex]

$\mathbf{\llbracket 80.\rrbracket}$\\*[.5ex]

13\zl .\zd \> \underline{Rule \quad $\Psi(x)\supset \Phi \quad \lfloor:\rfloor \quad (\exists x)\Psi(x) \supset \Psi$}\\[.5ex]
\> $\sim \Phi \supset\: \sim \Psi(x)$\\[.5ex]
\> $\sim \Phi \supset (x)\sim \Psi(x)$\\[.5ex]
\> $\sim (x)\sim \Psi (x) \supset \Phi$\\[.5ex]
\> $(\exists x)\Psi(x)\supset \Phi$\\[2ex]

\ul $13'$\zl .\zd \> \underline{$\varphi(y)\supset (\exists x)\varphi(x)$}\\[.5ex]
\> $(x)\sim \varphi(x)\supset\: \sim \varphi(y)$\\[.5ex]
\> $\varphi(y)\supset\: \sim (x)\sim \varphi(x)$ \qquad def.\zl ined\zd\ symb.\zl ol\zd\ \ud \\[2ex]

14.\> \underline{$(x)[\varphi(x) \supset \psi(x)]\supset [(\exists x)\varphi(x)\supset (\exists x)\psi(x)]$}\\[.5ex]
\> \sout{$(x)$} $[\varphi(x) \supset \psi(x)]\supset$\:\=$[\sim \psi(x) \supset\: \sim \varphi (x)]$\\*[.5ex]
$\times$\> $(x)$ \hspace{1.7em}..\zl $''$\zd \>$(x)$ \hspace{1.7em}\zl $''$\zd\\[.5ex]
$\times$\> $(x)[\sim \psi(x) \supset\: \sim \varphi (x)]\supset (x)\sim \psi(x)\supset (x) \sim \varphi (x) $ \\[.5ex]
$\times$\> $[(x)\sim \psi(x)\supset (x) \sim \varphi (x)] \supset \:\sim (x) \sim \varphi (x) \supset \:\sim (x) \sim \psi(x)$\\
\` $(p \supset q)\supset (\sim q \supset\: \sim p)$ \quad $\f{p}{(x)\sim \psi(x)}$ \;\; $\f{q}{(x)\sim \varphi(x)}$\\*
\` \zl  $\f{(x)\sim \psi(x)}{p}$ \; $\f{(x)\sim \varphi(x)}{q}$ \; fraction bars omitted in the manuscript\zd \\
\> $(x)[\varphi(x)\supset \psi(x)]\supset[\sim (x)\sim \varphi(x) \supset \:\sim (x) \sim \psi(x)]$ Rule of def\zl ined\zd\ \\*[.5ex]
\`  symb.\zl ol\zd\ \\[.5ex]

$\mathbf{\llbracket 81.\rrbracket}$\\*[.5ex]

15\zl .\zd \> Rule corresp.\zl onding\zd\ to $14.$\\[2ex]
16. \> \underline{$(\exists x)[\varphi(x) \vee \psi (x)] \equiv (\exists x)\varphi(x) \vee (\exists x)\psi(x)$}\\[.5ex]
\> $\varphi(x) \supset \varphi(x) \vee \psi(x)$\\[.5ex]
\> $(\exists x)\varphi(x) \supset (\exists x)[\varphi(x) \vee \psi(x)]$\\
\> \ldots\\[.5ex]
\> Dilemma\\[-1.5ex]
\> \underline{\hspace{12.5em}}\\[.5ex]
\> $\varphi(y) \supset (\exists x)\varphi(x)$\\[.5ex]
\> $\psi(y) \supset (\exists x)\psi(x)$\\[.5ex]
\> $\varphi(y) \vee \psi(y)$\:\=$\supset (\exists x)\varphi(x) \vee (\exists x) \psi(x)$\\*[.5ex]
\> $(\exists y)[$\quad $''$\quad $]$ \>$\supset$ \qquad $''$ \qquad \qquad $''$
\end{tabbing}

An example where we have to subst.\zl itute\zd\ for $\varphi(x)$ something containing other free var.\zl iables\zd\ besides $x$\zl :\zd
\begin{tabbing}
\hspace{3em}\=\underline{$(y)(x)\psi(xy)\equiv (x)(y)\psi(xy)$}\\[.5ex]
\> \underline{$(x)\varphi(x) \supset \varphi(y)$} \qquad \sout{$(z)\psi(xz)$}\\[.5ex]
\> $(x)\varphi(x) \supset \varphi(u)$ \qquad $\f{\psi(xy)}{\varphi(x)}$\\[.5ex]
\hspace{1.7em}*\> $(x)\psi(xy) \supset \psi(uy)$\\[.5ex]
\> $(z)\varphi(z) \supset \psi(y) \lfloor \varphi(y)\rfloor$ \qquad\= $\f{(x)\psi(xz)}{\varphi(z)}$\\[.5ex]
\hspace{1.7em}*\> $(z)(x)\psi(xz) \supset (x)\psi(xy)$ \> *\:* Syllog.\zl ism\zd \\[.5ex]
\> $(z)(x)\psi(xz) \supset \psi(uy)$ \> Rule 4 \quad\= $y$\\*[.5ex]
\> $(z)(x)\psi(xz) \supset (y)\psi(uy)$ \>\hspace{.9em} $''$ \> $u$\\[.5ex]
\> $($\=$z)(x)\psi(xz) \supset ($\=$u)(y)\psi(uy)$ \\[.5ex]
\>\> $y$ \> $x$\\[.5ex]
\> $(y)(x)\psi(xz) \supset (x)(y)\psi(uy)$
\end{tabbing}

$\mathbf{\llbracket 82.\rrbracket}$I have mentioned already that among the taut.\zl ological\zd\ form.\zl ulas\zd\ of the calc\zl ulus\zd\ of pred.\zl icates\zd\ are in part.\zl icular\zd\ those which express the Aristotelian moods of inf\zl erence\zd , \zl b\zd ut \ul \zl unreadable symbol\zd\ that \ud not all of the 19 Arist\zl otelian\zd\ moods are really valid in the calc.\zl ulus\zd\ of prop\zl ositions.\zd\ \zl \sout{but}\zd\ \sout{only 15 of them \zl unreadable word\zd\ the remaining 4} \ul \zl S\zd ome of them \ud require an add.\zl itional\zd\ third premise in order to be valid\zl ,\zd\ \ul namely that the predicates involved be not vacuous\zl ;\zd\ \ud  e.g.\ the mood Darapti is one of those not valid\zl ,\zd\ it says
\begin{tabbing}
\hspace{1.7em}\=$M$a$S$, $M$a$P$ \; : \; $S$i$P$, \quad in symbols:\\[.5ex]
\>$(x)[M(x)\supset S(x)]\:.\:(x)[M(x)\supset P(x)] \supset (\exists x)[S(x)\:.\:P(x)]$
\end{tabbing}
But this is not a tautological formula because that would mean it holds for any monadic pred.\zl icates\zd\ $M,S,P$ whatsoever. But $\mathbf{\llbracket 83.\rrbracket}$ we can easily name pred.\zl icates\zd\ for which it is wrong namely\zl ;\zd\ if you take for $M$ a \ul vacuous \ud pred.\zl icate\zd\ which belongs to no object\zl ,\zd\ say e.g.\ the pred\zl icate\zd\ president of A.\zl merica\zd\ born in South Bend \zl The following text is crossed out in the manuscript: that is a perfectly meaningful \ul correctly formed \ud pred\zl icate,\zd\ only by a historical accident there exists no object to which it belongs [or water snake is another ex.\zl ample\zd\ when a water snake is defined to be a snake living in the water.] Now I say if you take for $M$ such a vacuous pred.\zl icate\zd \zd\ and \sout{take} for $S$ and $P$ any two mutually exclusive pred.\zl icates,\zd\ i\zl .\zd e\zl .\zd\ such that no $S$ is $P$\zl ,\zd\ then the above formula will be wrong because \zl \sout{1.}\zd\ the two premises are both true\zl .\zd\ \zl S\zd ince $\mathbf{\llbracket 84.\rrbracket}$ $M(x)$ is false for every $x$ we have $M(x)\supset S(x)$ is true for every $x$ (bec.\zl ause it is an\zd\ impl.\zl ication\zd\ with false first term)\zl ;\zd\ likewise $M(x)\supset P(x)$ is true for every $x$\zl .\zd\ i\zl I.\zd e.\ the premises are both true but the conclusion is false bec.\zl ause\zd\ $S,P$ are supposed to be two predicates such that there is no $S$ which is a $P$\zl .\zd\ Hence for the part.\zl icular\zd\ pred.\zl icate\zd\ we chose the first term of this whole impl.\zl ication\zd\ is true \zl and\zd\ the sec\zl ond\zd\ is false\zl ,\zd\ i\zl .\zd e\zl .\zd\ the whole form.\zl ula\zd\ is false. So there are pred\zl icates\zd\ which substituted in this form\zl ula\zd\ yield a false prop\zl osition,\zd\ hence this form\zl ula\zd\ is not a taut\zl ology\zd . If we want to transform \ul that expr.\zl ession\zd\ into \ud \zl a\zd\ real taut.\zl ology\zd\ we have to add the further premise that $M$ is not $\mathbf{\llbracket 85.\rrbracket}$ vacuous\zl ,\zd\ i\zl .\zd e.\

\begin{tabbing}
\hspace{1.7em}$(\exists x)M(x)\:.\:(x)[M(x)\supset S(x)]\:.\:(x)[M(x)\supset P(x)]\supset (\exists x)[S(x)\:.\:P(x)]$
\end{tabbing}
would really be a tautology. Altogether there are 4\zl four\zd\ \sout{some} of the \ul $19$ \ud Arist.\zl otelian\zd\ moods which require this additional premise. \ul Furthermore $S\text{a}P \supset S\text{i}P$\zl ,\zd\ \ul $P$i$S$ (conversion) \ud as I mentioned last time also requires that $S$ \zl is\zd\ n\zl written over t\zd on\zl -\zd vac\zl uous\zd . Also $S\text{a}P\supset\: \sim (S\text{e}P)$\zl ,\zd\ i\zl .\zd e.\ $S$a$P$ and $S$e$P$ cannot both be true\zl ,\zd\ does not hold in the log.\zl ical\zd\ calc.\zl ulus\zd\ bec.\zl ause\zd\ if $S$ \zl is\zd\ vacuous both $S$a$P$ and $S$e$P$ are true $(x)[S(x)\supset P(x)]\:.\:(x)[S(x)\supset \sim P(x)]$\zl ;\zd\ $S(x)=$ $x$ \ul is a \ud pres.\zl ident\zd\ of the States born in Southb.\zl South Bend,\zd\ $P(x)\lfloor =\rfloor$ $x$ is bald, then both
\begin{tabbing}
\hspace{1.7em}\= Every \=presid.\zl ent\zd\ \ldots  \;\;\= is bald \\[.5ex]
\> No \> $\prime \prime$\zl president \ldots \zd\ \> is bald \ud
\end{tabbing}

So \ul we see \ud Arist.\zl otle\zd\ makes the implicit assumption that all pred.\zl i\-cates\zd\ which he speaks of are non-vacuous; in the logistic calc.\zl ulus\zd\ \sout{of pred.\zl icates\zd }\; however we do\zl $\,$\zd not make this assumption\zl ,\zd\ i\zl .\zd e.\ all tautologies and all formulas derivable from our axioms hold for any pred.\zl icates\zd\ whatsoever they may be\zl ,\zd\ vacuous or not. $\mathbf{\llbracket 86.\rrbracket}$ Now one may ask: which \sout{is the more expedient} procedure is preferable, to form\zl ulate\zd\ the laws of logic in such a way that they hold for all pred.\zl icates\zd\ \ul vacuous and non\zl -\zd vacuous \ud or in such a way that they hold only for non\zl -\zd vacuous. I think there can be no doubt that the logistic way is preferable for many reasons:

1. As we saw it may depend on purely empirical facts whether or not a pred.\zl icate\zd\ is vacuous (as we saw in the ex.\ul ample \ud of a presid.\zl ent\zd\ of America born in South\zl $\,$B\zd end). Therefore if we don't admit vacuous predic\zl ates\zd\ \sout{at all} it will depend on empirical facts which pred.\zl icates\zd\ \ul are \ud \sout{have} to be admitted in logical reasonings \ul or which inferences are valid, \ud but that $\mathbf{\llbracket 87.\rrbracket}$ is very undesirable. Whether a pred.\zl icate\zd\ can be used in reasoning (drawing inf.\zl erences\zd ) should depend only on mere logical considerations and not on empir.\zl ical\zd\ facts.

But a second and still more important argument is this\zl :\zd\ that to exclude vacuous predicates would be a very serious hampering\zl ,\zd\ e.g.\ in mathematical reasoning, because it happens frequently that we have to form pred\zl icates\zd\ of which we don't know in the beginning of \sout{the} \ul an \ud argument whether or not they are vacuous\zl ,\zd\ e.g.\ in indirect proofs\zl .\zd\ If we want to prove that there does not exist an alg.\zl ebraic\zd\ equ.\zl ation\zd\ whose root is $\pi$ we operate $\mathbf{\llbracket 88.\rrbracket}$ with the pred\zl icate\zd\ ,,\zl ``\zd algeb\zl raic\zd\ equ.\zl ation\zd\ with root $\pi$'' and use it in conclusions\zl ,\zd\ and later on it turns out that this pred.\zl icate\zd\ is vacuous. \ul But also in everyday life it happens frequently that we \sout{want} \ul have \ud to make general assertions about predicates of which we don't know whether they are vacuous [e\zl . E\zd .g.\ \sout{if} \ul A\zl a\zd ssume that \ud in a \zl u\zd niversity \zl deleted from the manuscript: \ul in Muham\zl medan\zd countries we have the true prop.\zl osition\zd\ \ud \zd there is the rule that examinations may be repeated arbitrarily often\zl ;\zd\ \ul then \ud we can make the statement\zl :\zd\ A student which has\ldots\ ten times is allowed to\ldots\ for an eleventh time\zl \sout{]}\zd . But if we want to exclude vacuous pred\zl icates\zd\ we cannot express this true prop\zl osition\zd\ \zl deleted from the manuscript: about \ul Turkey \ud (the univ.\zl ersity\zd\ under cons\zl isderation\zd )\zd\ if we don't know whether there exists \sout{such} a student who has\ldots\ But of course this (rule) \ul prop\zl osition\zd\ \ud has nothing to do with the exist\zl ence\zd\ of a student\ldots\ \ul \zl O\zd r e.g.\ excluding vac.\zl uous\zd\ pred.\zl icates\zd\ has the consequ\zl ence\zd\ that we cannot always form the conj.\zl unction\zd\ of two pred.\zl icates,\zd\ e.g.\ presid.\zl ent\zd\ of U.S.A. is \zl an\zd\ adm\zl issible\zd\ \ul pred.\zl icate,\zd\ \ud born in South Bend is adm.\zl issible,\zd\ but presid\zl ent\zd\ of Am\zl erica\zd\ born in South Bend is not admissible. \ud \ud So if we want to avoid absolutely unnecessary complications we \ul must not exclude the vacuous pred\zl icates\zd\ and \ud have to form.\zl ulate\zd\ the laws of logic in such a way that they apply both to vacuous and non-vacuous pred\zl icates\zd . I don't say that it is false to exclude them\zl ,\zd\ but it leads to abs.\zl olutely\zd\ unnec.\zl essary\zd\ complic\zl ations\zd . \zl After this paragraph the page is divided in the manuscript by a horizontal line.\zd\

As to the 15 valid moods of Arist.\zl otle\zd\ they can all be expressed by one logistic formula\zl .\zd\ \ul However \ud \zl i\zd n order to do that I have first to embody the calc.\zl ulus\zd\ of monadic pred\zl icates\zd\ in a different form\zl ,\zd\ namely in the form of the calc.\zl ulus\zd\ of classes. \sout{This \ul transformation however applies only to the \ud however applies only to the formulas containing only monadic pred\zl icates\zd }\, $\mathbf{\llbracket 89.\rrbracket}$ \sout{i.e.\ such that no var.\zl iables\zd\ for rel\zl ations\zd\ $\varphi(xy)$ occur)}. The calc.\zl ulus\zd\ of classes also yields \sout{also} the \sout{decision} solution of the decision problem for formulas with only monadic predicates.

If we have an arb.\zl itrary\zd\ \ul monadic \ud predicate\zl ,\zd\ say $P$\zl ,\zd\ then we can consider the extension of this pred.\zl icate,\zd\ i\zl .\zd e\zl .\zd\ the totality of all obj.\zl ects\zd\ satisfying $P$\zl ;\zd\ it is denoted by $\hat{x}[P(x)]$. These ext.\zl ensions\zd\ of monad.\zl ic\zd\ predicates \zl are\zd\ all called classes. So this \ul symb\zl ol\zd\ $\hat{x}$ \ud means: the class of obj.\zl ects\zd\ $x$ such that the subsequ.\zl ent\zd\ is true. It is applied also to prop.\zl ositional\zd\ f\zl u\zd nct\zl ions,\zd\ e.g.\ $\hat{x}[I(x)\:.\: x>7]$ means ,,\zl ``\zd the class of integers greater \zl than\zd\ seven''\zl .\zd\ $\mathbf{\llbracket 90.\rrbracket}$ \sout{$\hat{x}[T(x)]$ the class of most beings}. So to any monadic predicate belongs a uniquely det.\zl ermined\zd\ class of obj\zl ects\zd\ as its ,,\zl ``\zd extension''\zl ,\zd\ but of course there may be different predicates with the same extension\zl ,\zd\ as e\zl .\zd g.\ the two pred\zl icates\zd : \sout{good} heat conducting, elasticity conducting\zl . T\zd hese are two entirely diff.\zl erent\zd\ pred.\zl icates,\zd\ but every obj.\zl ect\zd\ which has the first prop\zl erty\zd\ also has the sec.\zl ond\zd\ one and vice versa\zl ;\zd\ therefore their ext\zl ension\zd\ is the same\zl ,\zd\ i\zl .\zd e.\ if $H,E$ denotes them, $\hat{x}[H(x)]=\hat{x}[E(x)]$ although $H \neq E$ \ul I am writing the symbol of identity \ul and distinctness \ud in between the two ident.\zl ical\zd\ obj.\zl ects\zd\ as is usual in math\zl ematics\zd . I shall speak \ul about \ud this way \ul of writing \ud in more detail later\zl .\zd\ \ud \\
In gen.\zl eral\zd\ we have if $\varphi,\psi$ are two mon.\zl adic\zd\ pred.\zl icates\zd\ then
\begin{tabbing}
\hspace{1.7em}$\hat{x}[\varphi(x)]=\hat{x}[\psi(x)] \equiv (x)[\varphi(x)\equiv \psi(x)]$
\end{tabbing}
\ul This equivalence expresses the essential property of extensions of pred\zl i\-cates\zd . It is to be noted \ul that \ud we have not defined what classes are bec\zl ause\zd\ we \ul explained it by the term extension\zl ,\zd\ and extensions we explained by the term totality\zl ,\zd\ and a totality is the same thing as a class. So this def.\zl inition\zd\ would be circular. The real state of affairs is this\zl :\zd\ that we consider $\widehat{x}$ as a \sout{primitive term} \zl \sout{a}\zd\ new primit.\zl ive\zd\ (undefined) term, which satisfies this axiom here. Russell \ul however \ud has shown that one can dispense with this \ul $\hat{x}$ as a \ud primit.\zl ive\zd\ term by introducing it by a kind of implicit def.\zl inition,\zd\ but that would take too much time \ul to explain it\zl ;\zd\ \ud so we simply can consider it as a primit\zl ive.\zd\ \ud \ud

The letters $\alpha, \beta, \gamma,\ldots$ are used as variables for classes and the statement that \zl The text interrupted here is continued on p.\ \textbf{91}., the first page of Notebook VI.\zd\ \sout{an obj\zl ect\zd\ $a$ bel\zl ongs\zd\ to $\alpha$ is den\zl oted\zd\ by $a \varepsilon \alpha$}.

\vspace{1ex}

\zl On the remaining not numbered, last page, of Notebook~V, one finds many lines in shorthand or crossed out, and one finds also: individual variables, \sout{faculta\-tive} \ul \underline{optional} \ud , convention, \sout{The interest lies in this that} propriety, \sout{choice is fortunate}, specific(individual, definite), Def~1. \sout{Expression $(P,\Phi(x) \rightarrow$}, 2.~\sout{\underline{Conv.}}, 4.~\sout{Taut.}, \sout{embody}.\zd\

\section{Notebook VI}\label{0VI}
\pagestyle{myheadings}\markboth{SOURCE TEXT}{NOTEBOOK VI}
\zl Folder 64, on the front cover of the notebook ``Log.\zl ik\zd\ Vorl.\zl esungen\zd\ \zl German: Logic Lectures\zd\ N.D.\ \zl Notre Dame\zd\ VI''\zd\

\vspace{1ex}

\zl The first page of this notebook, p.\ \textbf{91}., begins with the second part of a sentence interrupted at the end of p.\ \textbf{90}. of Notebook~V.\zd\

$\mathbf{\llbracket 91. \rrbracket}$ an obj\zl ect\zd\ $a$ bel.\zl ongs\zd\ to $\alpha$ (or is an el\zl ement\zd\ of $\alpha$) by $a\,\varepsilon\,\alpha$. Hence
\[
\fbox{\begin{minipage}{7.8em}
$y\,\varepsilon\,\hat{x}[\varphi(x)]\equiv\varphi(y)$
\end{minipage}} \quad\quad {\rm Furthermore}
 \begin{cases}
  \alpha = \hat{x}[x\,\varepsilon\,\alpha]\\
  (x)[x\,\varepsilon\,\alpha\equiv x\,\varepsilon\,\beta] \supset \alpha =\beta
 \end{cases}
\]
So far we spoke only of extensions of monadic predicates\zl ;\zd\ we can also introduce extensions of dyadic (and polyadic\zl )\zd\ pred\zl icates.\zd\ If e\zl .\zd g\zl .\zd\ $Q$ is a dyadic pred\zl icate\zd\ then $\hat{x}\hat{y}[Q(xy)]$ (called the ext.\zl ension\zd\ of $Q$) will be something that satisfies the condition:
\begin{tabbing}
\hspace{1.7em}$\hat{x}\hat{y}[\psi(xy)]=\hat{x}\hat{y}[\chi(xy)\lfloor ]\rfloor \lfloor .\rfloor\equiv \lfloor .\rfloor (x\lfloor ,\rfloor y)[\psi(xy)\equiv \chi(xy)]$
\end{tabbing}
e.g.\ the class of pairs \ul $(x,y)$ \ud such that $Q(xy)$ would $\mathbf{\llbracket 92. \rrbracket}$ be something which satisfies this cond.\zl ition,\zd\ but the ext.\zl ension\zd\ of \zl a\zd\ rel.\zl ation\zd\ is not defined as the class of ord\zl ered\zd\ pairs\zl ,\zd\ but is consid.\zl ered\zd\ as an und.\zl efined\zd\ term bec.\zl ause\zd\ ordered pair is defined in terms of ext.\zl ension\zd\ of relations. An example for this \ul formula\zl ,\zd\ \ud i\zl .\zd e.\ an \ul example \ud of two different dyadic \ul pred\zl icates\zd\ \ud which have the same extension would be $x<y$, $x>y\vee x=y$\zl ,\zd\ $x$ exerts an electrostatic a\zl t\zd traction on $y$\zl ,\zd\ $x$ and $y$ are loaded by electricities of different sign\zl .\zd\

\zl new paragraph\zd\ Ext.\zl ensions\zd\ of monadic pred.\zl icates\zd\ are called classes, \sout{\zl unreadable symbol\zd }\, extensions of polyadic pred.\zl icates\zd\ are called relations in logistic. So in log.\zl istic\zd\ the term rel\zl ation\zd\ is \sout{reserved} \ul used \ud not for the polyadic pred\zl icates\zd\ themselves but for their extensions, that \sout{\zl unreadable text\zd }\, conflicts with the meaning of the term rel.\zl ation\zd\ in everyday \zl life\zd\ \ul and also with the meaning in which I introduced this term a few lectures ago, \ud but since it is usual to use this term rel\zl ation\zd\ in \ul this ext\zl ensional\zd\ sense \ud I shall stick to this use \ul and the trouble is that \zl there is\zd\ no better term \ud . If $R$ is a rel\zl ation\zd , the statement that $x$ bears $R$ $\mathbf{\llbracket 93. \rrbracket}$ to $y$ is den\zl oted\zd\ by $xRy$. This way of writing\zl ,\zd\ \ul namely to write the symb\zl ol\zd\ denoting the rel.\zl ation\zd\ between the symbols denoting the obj.\zl ects\zd\ for which the rel\zl ation\zd\ is asserted to hold\zl ,\zd\ \ud is adapted to the notation of math\zl ematics,\zd\ e.g.\ $<$\zl ,\zd\ $x<y$, $=$, $x=y$\zl . O\zd f course we have:
\begin{tabbing}
\hspace{1.7em}$(x\lfloor ,\rfloor y)[xRy\equiv xSy]\supset R=S$
\end{tabbing}
for any two rel\zl ations\zd\ $R,S$\zl ,\zd\ \zl The text that follows, until the end of the paragraph, is inserted in the manuscript.\zd\ exactly as before $\lfloor (x)[\:\rfloor x\,\varepsilon\,\alpha\equiv x\,\varepsilon\,\beta\lfloor\: ] \supset \alpha =\beta\rfloor$. So a relation is uniquely det.\zl ermined\zd\ if you know all the pairs which have this relation bec.\zl ause\zd\ \ul by this form\zl ula\zd\ \ud there cannot exist two different rel.\zl ations\zd\ which subsist between the same pairs (although there can exist many different dyadic pred.\zl icates\zd )\zl .\zd

\zl The text that follows, until the end of the paragraph, is in big square brackets in the manuscript.\zd\ Therefore a relation can be represented e.g\zl .\zd\ by a figure of arrows
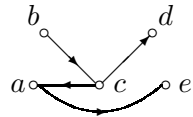
\begin{center}
\begin{picture}(100,40)(-5,10)

\put(45,20){\circle{3}}
\put(20,20){\circle{3}}
\put(24,40){\circle{3}}
\put(65,40){\circle{3}}
\put(70,20){\circle{3}}

\put(43.5,20){\line(-1,0){22}}
\put(43.5,21){\line(-1,1){18}}
\put(46,21){\line(1,1){18}}

\put(17,20){\small\makebox(0,0)[r]{$a$}}
\put(20,43){\small\makebox(0,0)[b]{$b$}}
\put(55,20){\small\makebox(0,0)[r]{$c$}}
\put(70,43){\small\makebox(0,0)[b]{$d$}}
\put(80,20){\small\makebox(0,0)[r]{$e$}}

\put(32,20){\vector(-1,0){2}}
\put(35,29.3){\vector(1,-1){2}}
\put(62,37){\vector(1,1){2}}
\put(46,10){\vector(1,0){2}}

\qbezier(22,20)(45,0)(68.5,20)

\end{picture}
\end{center}
or by a quadratic scheme e\zl .\zd g\zl .\zd\
\begin{center}
\begin{tabular}{ c|c c c c c }
& $a$ & $b$ & $c$ & $d$ & $e$\\
\hline
$a$ &&&&& $\bullet$\\
$b$ &&& $\bullet$ \\
$c$ & $\bullet$ &&& $\bullet$\\
$d$ & \\
$e$ &
\end{tabular}
\end{center}
Such a figure determines a unique rel.\zl ation;\zd\ in general it will be infinite\zl .\zd\

\zl The l\zd etters $R,S,T$ are mostly used as var.\zl iables\zd\ for rel\zl ations\zd . But now let us return to the ext.\zl ensions\zd\ of mon.\zl adic\zd\ pred.\zl icates,\zd\ i\zl .\zd e.\ the classes for which we want to set up a calculus.

First we have two part.\zl icular\zd\ classes $\bigwedge$ \zl written over 0\zd\ (vacuous class) \zl ,\zd\ $\bigvee$ (the universal class) which are defined as the ext.\zl ension\zd\ $\mathbf{\llbracket 94. \rrbracket}$ of a vacuous pred\zl icate\zd\ and of a pred\zl icate\zd\ that bel.\zl ongs\zd\ to everything. So
\begin{tabbing}
\hspace{1.7em}\=$\bigwedge\;$\=$=\; \hat{x}[\varphi(x)\: .\, \sim\varphi(x)]$\\[.5ex]
\>$\bigvee$\>$=\; \hat{x}[\varphi(x)\:\vee \sim\varphi(x)]$
\end{tabbing}
\sout{It is clear that} It makes no difference which vacuous pred\zl icate\zd\ I take for defining $\bigwedge$. If $A$\zl ,\zd\ $B$ are two diff\zl erent\zd\ vacuous pred.\zl icates\zd\ then $\hat{x}(A(x))=\hat{x}(B(x))$ \zl $\hat{x}[A(x)]=\hat{x}[B(x)]$\zd\  bec\zl ause\zd\ $(x)[A(x)\equiv B(x)]$. And similarly if $C,D$ are two diff.\zl erent\zd\ pred.\zl icates\zd\ belonging to everything $\hat{x}[C(x)]=\hat{x}[D(x)]$ bec.\zl ause\zd\ $(x)[C(x)\equiv D(x)]$\zl ,\zd\ i\zl .\zd e.\ there exists exactly one 0-class and exactly one $\mathbf{\llbracket 95. \rrbracket}$ universal class\zl ,\zd\ \ul although of course there exist many different vacuous pred\zl i\-cates\zd . But they all have the same extension\zl ,\zd\ namely nothing which is denoted by $\bigwedge$\zl .\zd\ So the zero class is the class with no el.\zl ements\zd\ $(x)[\sim x\,\varepsilon\, \bigwedge]$\zl $\bigwedge$ written over 0\zd \zl ,\zd\ the universal class is the class of which every obj.\zl ect\zd\ is an el.\zl ement\zd \zl unreadable text\zd\ $(x) \lfloor (\rfloor x\,\varepsilon\, \bigvee )$\zl ;\zd\ $\bigwedge$ and $\bigvee$ are sometimes denoted by 0 and 1 because of \zl \sout{[}\zd\ cert.\zl ain\zd\ analogies with arithm\zl etic\zd . \ud

\zl new paragraph\zd\ Next we can introduce cert.\zl ain\zd\ operations for classes which are analogous to the arithm\zl etical\zd\ operations: namely
\begin{tabbing}
\hspace{.5em}\= Add.\zl ition\zd\ or sum\hspace{7.5em}\= $\alpha +\beta$ \= $=\;\; \hat{x}[x\,\varepsilon\,\alpha\vee x\,\varepsilon\,\beta]$\\[.5ex]
\` $y\,\varepsilon\, \alpha +\beta\equiv y\,\varepsilon\,\hat{x}[x\,\varepsilon\,\alpha\vee x\,\varepsilon\,\beta]\equiv y\,\varepsilon\,\alpha\vee y\,\varepsilon\,\beta$\\[.5ex]
\hspace{15.8em} mathem.\zl atician\zd\ or dem.\zl ocrat\zd\ \\[.5ex]
\> Mult.\zl iplication\zd\ or inters.\zl ection\zd\ \> $\alpha\cdot\beta$ \> $=\;\; \hat{x}[x\,\varepsilon\,\alpha\: .\: x\,\varepsilon\,\beta]$\\[.5ex]
\hspace{15.8em} mathem\zl atician\zd\ democr.\zl at\zd\ \\[.5ex]
\> Op.\zl posite\zd\ or compl.\zl ement\zd\ \> $-\alpha$ \> $=\;\; \hat{x}[\sim x\,\varepsilon\,\alpha]$\quad or \quad$\overline{\alpha}$\\[.5ex]
\hspace{15.8em} non mathem.\zl atician\zd\ \\[.5ex]
\>Difference\> $\alpha -\beta$\> $=\;\;\alpha\cdot(-\beta)\;
=\;\hat{x}[x\,\varepsilon\,\alpha\: .\, \sim x\,\varepsilon\,\beta]$\\[.5ex]
\hspace{15.8em} mathem\zl atician\zd\ not democr\zl at\zd\ \\
\hspace{15.8em} (New Yorke\zl r\zd\ not sick)
\end{tabbing}
\zl On the right of the table above, two intersecting circles, as in Euler or Venn diagrams, are drawn in the manuscript.\zd\

Furthermore we have a rel.\zl ation\zd\ classes which corresponds to the arith\-m\zl e\-tic\zd\ rel\zl ation\zd\ of $<$\zl ,\zd\ namely the relation of subclass
\begin{tabbing}
\hspace{1.7em}$\alpha\subseteq\beta\equiv(x)[x\,\varepsilon\,\alpha\supset x\,\varepsilon\,\beta]$\qquad \ul Man $\subseteq$ Mortal \ud
\end{tabbing}
All these op.\zl erations\zd\ obey laws very similar $\mathbf{\llbracket 96. \rrbracket}$ to the corresponding arithmetical laws: e.g.\
\begin{tabbing}
\hspace{1.7em}\= $\alpha +\beta = \beta +\alpha$\hspace{2em}$\alpha\cdot\beta = \beta\cdot\alpha$\\[.5ex]
\>$(\alpha +\beta)+\gamma = \alpha +(\beta +\gamma)$\hspace{2em}$(\alpha\cdot\beta)\cdot\gamma = \alpha\cdot(\beta\cdot\gamma)$\\[.5ex]
\>$(\alpha +\beta)\cdot\gamma = \alpha\cdot\gamma +\alpha\cdot\gamma \lfloor \beta\cdot\gamma\rfloor$\\[.5ex]
\>$(\alpha\cdot\beta)+\gamma = (\alpha +\gamma)\cdot(\alpha +\gamma) \lfloor \beta +\gamma\rfloor$
\end{tabbing}
\ul \zl T\zd hey follow from the corresponding laws of the calculus of prop.\zl ositions:\zd\ e.g.\
\begin{tabbing}
\hspace{1.7em}\= $x\,\varepsilon(\alpha +\beta)\equiv x\,\varepsilon\,\alpha\vee x\,\varepsilon\,\beta\equiv x\,\varepsilon\,\beta\vee x\,\varepsilon\,\alpha\equiv x\,\varepsilon(\beta +\alpha)$\\[1ex]
\>$x\,\varepsilon(\alpha +\beta)\cdot\gamma\equiv x\,\varepsilon(\alpha +\beta)\: .\: x\,\varepsilon\,\gamma\equiv(x\,\varepsilon\,\alpha\vee x\,\varepsilon\,\beta)\: .\: x\,\varepsilon\,\gamma$\\[.5ex]
\>\hspace{4em}\=$\equiv(x\,\varepsilon\,\alpha\: .\: x\,\varepsilon\,\gamma)\vee (x\,\varepsilon\,\alpha
\lfloor x\,\varepsilon\,\beta\rfloor .\: x\,\varepsilon\,\gamma)$\\[.5ex]
\>\> $\equiv x\,\varepsilon\,\alpha\cdot\beta\vee x\,\varepsilon\,\alpha\gamma\lfloor x\,\varepsilon\,\alpha\cdot\gamma\vee x\,\varepsilon\,\beta\cdot\gamma\rfloor$\\[.5ex]
\>\>$\equiv x\,\varepsilon(\alpha\cdot\beta +\alpha\cdot\gamma)\lfloor x\,\varepsilon(\alpha\cdot\gamma +\beta\cdot\gamma)\rfloor$\quad \zl $(\alpha +\beta)\cdot\gamma$ deleted\zd\ \ud\\[1ex]
\>$\alpha + 0\;$\=$=\alpha$\hspace{2em}\=$\alpha\cdot 0$\hspace{.7em}\=$=0$\\[.5ex]
\>$\alpha\cdot 1$\>$=\alpha$\>$\alpha +1$\>$=1$\\[1ex]
\;\;\;\ul \>$(x)\sim(x\,\varepsilon\, 0)$\hspace{3em}$x\,\varepsilon(\alpha +0)\equiv x\,\varepsilon\,\alpha\vee x\,\varepsilon\, 0\equiv x\,\varepsilon\,\alpha$\\[.5ex]
\>$(x)\lfloor (\rfloor x\,\varepsilon\, 1)$ \ud
\end{tabbing}
\zl On the right of the table above, three intersecting circles, as in Euler or Venn diagrams, with $\alpha$, $\beta$ and perhaps $\gamma$ marked in them, and some areas shaded, are drawn in the manuscript.\zd\
\begin{tabbing}
\hspace{3.5em}$\alpha\subseteq\beta$\hspace{7em}\=$\alpha\subseteq\beta\lfloor .\rfloor\beta\subseteq\gamma\supset\alpha\subseteq\gamma$\\[.5ex]
\hspace{3.5em}$\gamma\subseteq\delta$\hspace{8em} Law of transitivity\\[-2ex]
\hspace{1.7em}\underline{\hspace{6em}}\\
\hspace{1.7em}$\alpha +\gamma\subseteq\beta +\delta$\\
\hspace{2.2em}$\alpha\cdot\gamma\subseteq\beta\cdot\delta$\>$\alpha\subseteq\beta\: .\:\beta\subseteq\alpha\supset\alpha =\beta$.
\end{tabbing}
Laws different from arithm\zl etical:\zd\
\begin{tabbing}
\hspace{1.7em}\=$\alpha +\alpha =\alpha\cdot\alpha=\alpha$\hspace{3.5em}$x\,\varepsilon\,\alpha \cdot\lfloor +\rfloor\alpha\equiv x\,\varepsilon\,\alpha\vee x\,\varepsilon\,\alpha\equiv x\,\varepsilon\,\alpha$\\[.5ex]
\>$\alpha\subseteq\beta\supset[\alpha +\beta =\beta\;\; .\;\; \alpha\cdot\beta=\alpha]$\hspace{2em}$\beta\subseteq\alpha +\beta$\hspace{2em}\=$\alpha\subseteq\beta$\\
\>\>$\beta\subseteq\beta$\\[-2ex]
\hspace{22.7em}\underline{\hspace{8em}}\\
\hspace{22.7em}$\alpha +\beta\subseteq\beta +\beta=\beta$
\end{tabbing}
$\mathbf{\llbracket 97. \rrbracket}$
\begin{tabbing}
\hspace{1.7em}\=$-(\alpha +\beta)=(-\alpha)\cdot(-\beta)$\hspace{2em} De Morgan\\[.5ex]
\` $x\,\varepsilon\, -(\alpha +\beta)\equiv\;\sim x\,\varepsilon\, (\alpha +\beta)\equiv\;\sim(x\,\varepsilon\,\alpha\vee x\,\varepsilon\,\beta)\equiv\; \sim(x\,\varepsilon\,\alpha)\: .\:\sim(x\,\varepsilon\,\beta)\equiv$\\
\` $x\,\varepsilon\, -\alpha\: .\: x\,\varepsilon\, -\beta\equiv\; x\,\varepsilon\, (-\alpha)\cdot(-\beta)$\\
\>$-(\alpha\cdot\beta)=(-\alpha)+(-\beta)$\\[.5ex]
\>$\alpha\cdot(-\alpha)=0$\hspace{2em}$\alpha +(-\alpha)=1$\\[.5ex]
\>$-(-\alpha)=\alpha$
\end{tabbing}
\ul The compl.\zl ement\zd\ of $\alpha$ is sometimes also denoted by $\overline{\alpha}$ (so that $\overline{\alpha}= -\alpha$)\zl .\zd\;\ud

\zl The exercise that follows, with three displayed formulae, is in big square brackets in the manuscript.\zd\ Exercise \zl unreadable text\zd\ Law for diff.\zl erence:\zd\
\begin{tabbing}
\hspace{1.7em}\=$\alpha\lfloor\cdot\rfloor(\beta -\gamma)=\alpha\cdot\beta -\alpha\cdot\gamma$\\[.5ex]
\>$\alpha\cdot\beta = \alpha -(\alpha -\beta)$\\[.5ex]
\>$\alpha\subseteq\beta \supset \overline{\beta}\subseteq\overline{\alpha}$
\end{tabbing}

If $\alpha\cdot\beta=0$\zl ,\zd\ that means the classes $\alpha$ and $\beta$ have no common element\zl ,\zd\ then $\alpha$ \zl and\zd\ $\beta$ are called mutually exclusive. We can now formulate the four Aristotelian types of judgement a, e, i, o also in the symbolism of the calc.\zl ulus\zd\ of classes as follows\zl :\zd\
\begin{tabbing}
\hspace{1.7em}\=$\alpha\, {\rm a}\, \beta$\=$\;\equiv\alpha\subseteq\beta\equiv\;$$|$\underline{$\alpha\cdot\overline{\beta}=0$}$|$\\[.5ex]
$\mathbf{\llbracket 98. \rrbracket}$\\[.5ex]
\>$\alpha\, {\rm e}\, \beta$\>$\;\equiv\;$ \underline{$\alpha\cdot\beta=0$} \=$\;\equiv$ \hspace{1.2em}$\alpha\subseteq\overline{\beta}$\hspace{1.7em}\=$\equiv\;\alpha\lfloor\cdot\rfloor\beta=0$\\[.5ex]
\>$\alpha\, {\rm i}\, \beta$\>$\;\equiv\;$ \underline{$\alpha\cdot\beta\neq 0$} \>$\;\equiv$ $\;\sim(\alpha\subseteq\overline{\beta})$\>$\equiv\;\alpha\cdot\beta\neq 0$\\[.5ex]
\>$\alpha\, {\rm o}\, \beta$\>$\;\equiv\;$ \underline{$\alpha\cdot\overline{\beta}\neq 0$} \>$\;\equiv$ $\;\sim(\alpha\subseteq\beta)\,$.\zl \>$\equiv\rfloor\;\alpha\cdot\overline{\beta}\neq 0$
\end{tabbing}
\zl In the last three lines, the underlined formulae and the $\equiv$ symbol that follows them are to be deleted, since they are repeated at the end of the lines.\zd\ So all of these 4\zl four\zd\ types of judgements can be expressed by the vanishing\zl ,\zd\ resp.\zl ectively\zd\ not vanishing\zl ,\zd\ of cert.\zl ain\zd\ intersections.

Now the formula which compresses all of the 15 valid Aristotelian inferences reads like this
\[
\sim(\alpha\cdot\beta =0\;\; .\;\; \overline{\alpha}\cdot\gamma =0\;\; .\;\; \beta\cdot\gamma\neq 0)
\]
So this is a universally true formula bec\zl ause\zd\ $\alpha\cdot\beta\lfloor =0\rfloor$ means $\beta$ outside of $\alpha$\zl ,\zd\ $\overline{\alpha}\cdot\gamma =0$ means $\gamma$ inside of $\alpha$\zl .\zd\ If $\beta$ outside $\gamma$ inside they can have no element in $\mathbf{\llbracket 99. \rrbracket}$ common\zl ,\zd\ i\zl .\zd e.\ the two first prop\zl ositions\zd\ imply $\beta\cdot\gamma =0$\zl ,\zd\ i\zl .\zd e.\ it cannot be that all three of them are true\zl .\zd\ Now since this says that all \zl written over ``All''\zd\ three of them cannot be true you can always conclude the negation of the third from the two others\zl ;\zd\ e.g\zl .\zd\
\begin{tabbing}
\hspace{1.7em}\=$\alpha\cdot\beta =0\;\; .\;\; \alpha\cdot\overline{\gamma}=0\;\lfloor \overline{\alpha}\cdot\gamma =0\rfloor\;\supset\;\beta\cdot\gamma =0$\\[.5ex]
\>$\alpha\cdot\beta =0\;\; .\;\; \beta\cdot\gamma\neq 0\;\supset\;\overline{\alpha}\cdot\gamma\neq 0$\quad etc\zl .\zd\
\end{tabbing}
and in this way you obtain all valid 15 moods if you substitute for $\alpha,\beta,\gamma$ \sout{the} in an appropriate way the minor term\zl ,\zd\ the major term and the middle term or their neg\zl ation,\zd\ e\zl .\zd g\zl .\zd\
\begin{tabbing}
$\mathbf{\llbracket 100. \rrbracket}$\\[1ex]
\hspace{1.7em} I \quad \underline{\zl B\zd arbara}
\hspace{2em}\begin{tabular}{ l|l }
$M$a$P$ & \\[-1ex]
& \hspace{.5em}$S$a$P$ \\[-1ex]
$S$a$M$&
\end{tabular}\\[1ex]
\hspace{2.7em}$\lfloor M\cdot\overline{P} =0\;\; .\;\; S\cdot\overline{M} =0\rfloor\supset S\cdot \overline{P}=0$\\[.5ex]
\hspace{1.7em}$\sim(M\cdot\overline{P} =0\;\; .\;\; S\cdot\overline{M} =0\,\;\; .\,\;\; S\cdot\overline{P}\neq 0)$\\[.5ex]
\hspace{1.7em}$\alpha=M$\hspace{2em}$\beta=\overline{P}$\hspace{2em}$\gamma =S$\\[2ex]

\hspace{1.7em} III \quad \underline{\zl F\zd eriso}
\hspace{2em}\begin{tabular}{ l|l }
$M$e$P$ & \\[-1ex]
& \hspace{.5em}$S$o$P$ \\[-1ex]
$M$i$S$&
\end{tabular}\\[1ex]
\hspace{3.6em}$M\cdot \underline{P}\lfloor P\rfloor =0\;\; .\;\; M\cdot S \neq 0\;\supset\; S\cdot\overline{P}\neq 0$\\[.5ex]
\hspace{1.7em}$\sim(\lfloor M\cdot P =0\;\; .\;\; M\cdot S \neq 0\rfloor \;\; .\;\; S\cdot\underline{\overline{P}}\lfloor \overline{P}\rfloor= 0)$\\[.5ex]
\hspace{1.7em}$\alpha=P$\hspace{2em}$\beta=M$\hspace{2em}$\gamma =P\lfloor \gamma =S\rfloor$.
\end{tabbing}

The \ul 4\zl four\zd\ \ud moods which require an additional premise can also be expressed by one formula\zl ,\zd\ namely:
\[
\sim(\alpha\neq 0\;\; .\;\;\alpha\cdot\beta =0\;\; .\;\; \alpha\cdot\gamma =0\;\; .\;\; \overline{\beta}\cdot\overline{\gamma}= 0)
\]

$\mathbf{\llbracket 101. \rrbracket}$ Darapti
\begin{tabbing}
\hspace{9.8em}\=$M$a$P$\\[.5ex]
\>$M$a$S$\\[-2ex]
\>\underline{\hspace{2.8em}}\\
\>$\;S$i$P$
\end{tabbing}
e.g.\ is obtained by taking
\begin{tabbing}
\hspace{1.7em}\=$M=\alpha\lfloor \alpha =M\rfloor$\hspace{2em}$\beta=\overline{P}$\hspace{2em}$\gamma =\overline{S}$\\[.5ex]
\>$M$a$P\: .\: M$a$S\supset S$i$P$\\[.5ex]
\>$M\lfloor\cdot\rfloor\overline{P}=0\lfloor .\rfloor M\lfloor\cdot\rfloor\overline{S}=0\lfloor \supset\rfloor S\cdot P\lfloor\neq\rfloor 0$\\[.5ex]
\` \zl $\beta=\overline{P}$ and $\gamma =\overline{S}$, which are written already above, are deleted\zd\
\end{tabbing}
However\zl ,\zd\ this sec.\zl ond\zd\ formula is an easy consequence of the first\zl ,\zd\ i\zl .\zd e.\ we can derive it by two applications of the first. To this end we have only to note that $\alpha\neq 0$ can be expressed by $\alpha\,{\rm i}\,\alpha$ bec.\zl ause\zd\ \sout{\zl unreadable symbol\zd }
\begin{tabbing}
\hspace{1.7em}\=$\varphi\,{\rm i}\,\psi\equiv(\exists x)[\varphi(x)\lfloor .\rfloor\psi(x)]$\\[.5ex]
\>$\varphi\,{\rm i}\,\varphi\equiv(\exists x)[\varphi(x)\: .\:\varphi(x)]\equiv(\exists x)\varphi(x)$\\[1ex]
\>$\sim(\alpha\lfloor\cdot\rfloor\beta =0\;\; .\;\; \overline{\alpha}\lfloor\cdot\rfloor\gamma =0\;\; .\;\; \beta\lfloor\cdot\rfloor\gamma\neq 0)$\\[.5ex]
\>$\alpha\lfloor\cdot\rfloor\alpha\neq 0$\hspace{2em}$\alpha\underline{\beta}=0\:\lfloor\alpha\cdot\beta=0\rfloor$\hspace{2em}
$\alpha\underline{\overline{\beta}}=0\:\lfloor\alpha\cdot\overline{\beta}=0\rfloor$\\[.5ex]
\>$\alpha:\beta$\hspace{2em}$\beta :\gamma\alpha$ \zl perhaps $\beta,\gamma :\alpha$, which should mean: $\beta :\alpha$, \hspace{.1em} $\gamma :\alpha$\zd\ \\[1ex]
\>$\alpha\lfloor\cdot\rfloor\overline{\beta}\neq 0$\hspace{2em}$\alpha\lfloor\cdot\rfloor\gamma=0$\hspace{2em}
$\overline{\beta}\overline{\gamma}\neq 0\:\lfloor\overline{\beta}\cdot\overline{\gamma}= 0\rfloor$\\[.5ex]
\>$\alpha:\gamma$\hspace{2em}$\beta :\alpha$\hspace{2em}$\gamma :\overline{\beta}$\\[1ex]

III \;\zl F\zd eriso\quad $\alpha\lfloor\cdot\rfloor\alpha\neq 0\;\; .\;\;\alpha\lfloor\cdot\rfloor\beta=0\;\; .\;\;\alpha\lfloor\cdot\rfloor\gamma=0 \;\supset\; \overline{\beta}\lfloor\cdot\rfloor\overline{\gamma}\neq 0$\\[.5ex]

\hspace{5.5em}$\f{\f{\lfloor\alpha\cdot\alpha\neq 0\hspace{1.6em}\alpha\cdot\beta=0\rfloor}{\alpha\lfloor\cdot\rfloor\overline{\beta}\neq 0}\hspace{1.2em}\afrac{\alpha\lfloor\cdot\rfloor\gamma=0}}
{\overline{\beta}\lfloor\cdot\rfloor\overline{\gamma}\neq 0}$
\end{tabbing}

$\mathbf{\llbracket 102. \rrbracket}$ In general it can be shown that every correct formula express.\zl ed\zd\ by the Arist\zl otelian\zd\ terms a, e, i, o and op\zl erations\zd\ of \zl the\zd calc.\zl ulus\zd\ of prop.\zl osi\-tions\zd\ can be derived from this principle\zl ;\zd\ to be more exact\zl ,\zd\ fund\zl amental\zd\ notions a, i
\begin{tabbing}
def\hspace{2em}\=$\alpha\,{\rm e}\,\beta\equiv\;\sim(\alpha\,{\rm i}\,\beta)$\\[.5ex]
\>$\alpha\,{\rm o}\,\beta\equiv\;\sim(\alpha\,{\rm a}\,\beta)$\\[1ex]
\>1. \hspace{.5em}\=$\alpha\,{\rm a}\,\alpha$\hspace{2em}Identity\\[.5ex]
\>2. \>$\alpha\,{\rm a}\,\beta\: .\: \beta\,{\rm a}\,\gamma\supset\alpha\,{\rm a}\,\gamma$\hspace{2em}I\hspace{.5em}Barbar\zl a\zd\ \\[.5ex]
\>3. \> $\alpha\,{\rm i}\,\underline{\beta}\lfloor\beta\rfloor\: .\: \underline{\beta}\lfloor\beta\rfloor\,{\rm a}\,\gamma\supset\gamma\,{\rm i}\,\underline{\alpha}\lfloor\alpha\rfloor$\hspace{1em}\sout{Darii}\hspace{1em} IV\hspace{.5em}\zl D\zd imatis
\end{tabbing}
and all axioms of the prop.\zl ositional\zd\ calculus\zl ;\zd\ then if we have a form\zl ula\zd\ composed only of such expr.\zl essions\zd\ $\alpha\,{\rm a}\,\beta$\zl ,\zd\ $\alpha\,{\rm i}\,\gamma$ and $\sim,\vee\ldots$ and which is universally true\zl ,\zd\ i\zl .\zd e.\ holds for all classes $\alpha\lfloor ,\rfloor\beta,\gamma$ involved\zl ,\zd\ then it is derivable from these ax.\zl ioms\ by rule of subst\zl itution\zd\ and impl.\zl ication\zd\ and def.\zl ined\zd\ symb\zl ol\zd . $\mathbf{\llbracket 103. \rrbracket}$ I am sorry I have no time to give the proof.

\zl new paragraph\zd\ So we can say that the Aristotelian theory of syllogisms for expressions of this part.\zl icular\zd\ type a, e, i, o is complete\zl ,\zd\ i\zl .\zd e.\ every true formula follows from the Aristotelian moods. \zl The following inserted jottings from the manuscript are deleted: $\mu$\zl or $u$\zd $\cdot\beta=0$, $\overline{\nu}\gamma=0$, $\overline{\mu}$\zl or $\overline{u}$\zd $\nu=0$.\zd\ \zl B\zd ut those Arist\zl otelian\zd\ moods are even abundant because those two moods alone are already sufficient to obtain everything else. \sout{But} The incompleteness of the Aristot\zl elian\zd\ theory lies in this that there are many
$\mathbf{\llbracket 104. \rrbracket}$ propositions which cannot be expressed in terms of the Arist.\zl otelian\zd\ primit.\zl ive\zd\ terms. E.g.\ all form.\zl ulas\zd\ which I wrote down for $+\lfloor ,\rfloor\cdot\lfloor ,\rfloor -$ (distrib\zl utive\zd\ law, De Morgan law etc.) bec\zl ause\zd\ those symb.\zl ols\zd\ $+\lfloor ,\rfloor\cdot\lfloor ,\rfloor -$ do\zl $\,$\zd not occur in Arist\zl otle\zd . But there are even simpler things not expr.\zl essible\zd\ in Arist\zl otelian terms;\zd \zl left square bracket deleted\zd\ e\zl .\zd g.\ $\overline{a}\cdot\overline{c}=0$ \zl full stop deleted\zd\ (some not $a$ are not $c$)\zl ,\zd\ e.g.\ ${\alpha\, {\rm e}\,\beta \atop \underline{\beta\, {\rm o}\,\gamma}}$ according to Arist\zl otle\zd\ there is no concl.\zl usion\zd\ from that (there is a \zl principle\zd\ that from two neg.\zl ative\zd\ premises no conclusion can be drawn)

\vspace{1ex}

\zl On the right of p.\ \textbf{104}.\ one finds in the manuscript the following jottings:\zd\
\begin{tabbing}
\hspace{1.7em}\= $\alpha$\zl written over $\beta$\zd\ = Comm.\zl unist?\zd\ \\[.5ex]
\>$\beta$ = Dem\zl ocrat\zd\ \hspace{15em}\= $\beta\,{\rm a}\, \overline{\alpha}$ \\[.5ex]
\>$\gamma$ = Math\zl ematician\zd\ \> $\beta\,{\rm o}\, \overline{\gamma}$ \\[-2ex]
\>\underline{\hspace{9.3em}} \> \underline{\hspace{2em}}\\
\> \zl the conclusion is presumably in shorthand\zd\ \> $\overline{\alpha}\,{\rm o}\,\gamma$
\end{tabbing}

\noindent $\mathbf{\llbracket 105. \rrbracket}$ and that is true if we take account only of propositions expressible by the a, e, i, o\zl .\zd\ But there is a concl.\zl usion\zd\ to be drawn from that\zl ,\zd\ namely \zl ``\zd Some not $\alpha$ are not $\gamma$\zl ''\zd\ $\overline{\alpha}\cdot\overline{\gamma}\neq 0$\zl .\zd\ Since some $\beta$ are not $\gamma$ and every $\beta$ is not $\alpha$ we have some not $\alpha$ (namely the $\beta$) are not $\gamma$\zl .\zd \zl right parenthesis deleted\zd \zl The relation\zd\ which \zl holds\zd\ between two classes $\alpha,\gamma$ if $\overline{\alpha}\cdot\overline{\gamma}\neq 0$ cannot be expressed by a, e, i, o\zl ,\zd\ but it is arb.\zl itrary\zd\ to exclude that rel\zl ation\zd . \zl $\Big\rfloor$ deleted\zd\ Another ex.\zl ample\zd
\begin{tabbing}
\hspace{1.7em}\= $\alpha\, {\rm i}\,\beta$\\[.5ex]
\>$\alpha\, {\rm o}\,\beta$\\[-2ex]
\>\underline{\hspace{2.5em}}\\
\> $\beta$\zl $\alpha$\zd\ contains at least t\zl w\zd o elements
\end{tabbing}

\zl On the right of p.\ \textbf{105}.\ one finds in the manuscript:
\begin{tabbing}
\hspace{1.7em}\= $M$e$P$\hspace{5em}\= $M$a$\overline{P}$\\[.5ex]
\>$S$a$M$\>$S$a$M$\\[-2ex]
\>\underline{\hspace{2.5em}}\>\underline{\hspace{2.5em}}\\
\> $S$e$P$\>$S$a$\overline{P}$
\end{tabbing}
which show that the mood Celarent of the first figure is really Barbara.\zd

\vspace{1ex}

\noindent $\mathbf{\llbracket 106. \rrbracket}$ Such prop.\zl ositions: ``\zd There are two diff.\zl erent\zd\ objects \sout{$a,b$} to which the pred\zl icate\zd\ $\alpha$ belongs\zl ''\zd\ can of course not be expr.\zl essed\zd\ by a, e, i, o\zl ,\zd\ but they can in the logistic calc.\zl ulus\zd\ by
\begin{tabbing}
\hspace{1.7em}$(\exists x\lfloor ,\rfloor y)[x\neq y\: .\: x\,\varepsilon\,\alpha \: .\: y\,\varepsilon\,\beta\lfloor y\,\varepsilon\,\alpha\rfloor]$\zl .\zd\
\end{tabbing}
\zl Here, after ``Another ex:'' the text is interrupted in the manuscript.\zd\

$\mathbf{\llbracket 107. \rrbracket}$ \zl The following paragraph is crossed out in the manuscript:

\vspace{1ex}

We have seen already in the theory of the monadic pred\zl icates\zd\ for classes \sout{that many} that many concepts \zl unreadable text\zd\ laws of \zl unreadable text\zd\ are missing in the Arist\zl otelian\zd\ treatment\zl .\zd\ But the proper domain of logic where the incompleteness of Arist\zl otelian\zd \zl unreadable text\zd\ in terms of diff\zl erent\zd \zl unreadable text\zd is the theory of relations. \sout{\zl unreadable text\zd\ we are going to deal with in more detail \zl unreadable text\zd\ relations}\zd\

\vspace{1ex}

\ul Last time I developed in outline the calc.\zl ulus\zd\ of classes in which we introduced certain operations +\zl ,\zd\ $\cdot$\zl ,\zd\ $-$ which obey laws similar \sout{laws} \ul to those \ud of arithmetic\zl .\zd\ \ud One can develop a \ul similar \ud calc.\zl ulus\zd\ \ul for relations\zl .\zd\ \ud First of all we can introduce for relations operations +\zl ,\zd\ $\cdot$\zl ,\zd\ $-$ in a manner perfectly analogous to the calc.\zl ulus\zd\ of classes.

$\mathbf{\llbracket 108. \rrbracket}$ If $R$ \zl and\zd\ $S$ are any two dyad.\zl ic\zd\ rel.\zl ations\zd\ I put
\begin{tabbing}
\hspace{1.7em}\= $R +S\;$ \= $=\;\; \hat{x}\hat{y}[xRy\vee xSy]$\\[.5ex]
\>$R\lfloor\cdot\rfloor S$ \> $=\;\; \hat{x}\hat{y}[xRy\; .\; xSy]$\\[.5ex]
\>$-R$\> $=\;\; \hat{x}\hat{y}[\sim xRy]$ \zl unreadable word\zd\ p\zl .\zd\ 110\\[.5ex]
\;\;\;\ul \>$R-S$\> $=\;\; \hat{x}\hat{y}[xRy\; .\: \sim xSy]$ \ud
\end{tabbing}

$\Big\lceil$So e.g.\ if $R$ is the rel\zl ation\zd\ of father, $S$ the rel\zl ation\zd\ of mother \zl unreadable text; should be: one has for the relation\zd\ of parent\zl :\zd\
\begin{tabbing}
\hspace{1.7em}\= parent \= = \hspace{.2em} \= father + mother\\[.5ex]
\> $x$ \zl is a\zd\ parent of $y$ $\equiv$ $x$ is a father of $y$ $\vee$ $x$ is \zl a\zd\ mother of $y$\\[.5ex]
\> $\leq$\> = \>$(< + =)$\\[.5ex]
\> child\> = \>son + daughter$\Big\rfloor$
\end{tabbing}

\zl The following unfinished paragraph at the end of p.\ \textbf{108}.\ is crossed out:

\vspace{1ex}

subrel\zl ation\zd . $R$ is called a subrelation of $S$
\begin{tabbing}
\hspace{1.7em}$R\subseteq S$ if $(x\lfloor ,\rfloor y)[xRy \supset xSy]$
\end{tabbing}
e.g.\ father $\subseteq$ ancestor\zl ,\zd\ but not\zd\

\vspace{1ex}

$\mathbf{\llbracket 109. \rrbracket}$ Or consider similarity for polygons and \zl the\zd rel\zl ation\zd\ of \zl unread\-able text, perhaps in shorthand, maybe: same size\zd\ and the rel\zl ation\zd\ of congr.\zl uence\zd , then Congr\zl uence\zd\ = Simil\zl arity\zd\ $\cdot$ \zl unreadable text, perhaps in shorthand, same as the preceding one, maybe: Same size\zd , or consider\zl the text until ``then we have'' is partly crossed out\zd\ the 4\zl four, written over 3\zd\ rel\zl ations\zd\ $\|$ \zl parallelism\zd , without com\zl mon\zd\ points, co\zl \sout{m}\zd planar, and \zl unreadable text, perhaps in shorthand, maybe: skew,\zd\ then we have
\begin{tabbing}
\hspace{1.7em}\=\zl \sout{or ,,}\zd\ Parallelism = without com\zl mon\zd\ point \zl $\cdot$\zd\ co\zl \sout{m}\zd planar,\\[.5ex]
\>or Parallelism \zl =\zd\ without com\zl mon\zd\ point ,\zl $\cdot$\zd\ $-$ \zl unreadable text,\\
\`  perhaps in shorthand, same as the preceding one, maybe: skew\zd\
\end{tabbing}
or $-$brother will subsist\zl unreadable letter\zd\ between two obj\zl ects\zd\ $x,y$ if 1. $x,y$ are two human beings and $x$ is not a brother of $y$ or 2\zl .\zd\ if $x$ or $y$ is not a human being bec.\zl ause\zd\ $x$ brother $y$ is true only if $x$ and $y$ are human beings and in addition $x$ is a brother of $y$. So if $x$ or $y$ are not human beings the relation eo ipso will not $\mathbf{\llbracket 110. \rrbracket}$ \zl \sout{will not}\zd\ hold\zl ,\zd\ i\zl .\zd e.\ \ul the rel\zl ation \ud $-$brother will hold. \ul \ul Exactly \ud as for classes there will exist also a vacuous and a universal relation denoted by $\dot{\Lambda}$ and $\dot{\rm{V}}$. $\dot{\Lambda}$ is the rel\zl ation\zd\ which subsists between no obj\zl ects\zd\ $(x\lfloor ,\rfloor y)\sim x \dot{\Lambda} y$\zl , and\zd\ $(x\lfloor ,\rfloor y)x \dot{\rm{V}} y$\zl ,\zd\ e\zl .\zd g\zl .\zd\
\begin{tabbing}
\hspace{1.7em}\=greater $\cdot$ smaller = $\dot{\Lambda}$\\[.5ex]
\>greater + (not greater) = $\dot{\rm{V}}$
\end{tabbing}
Also there exists an analogon to the notion of subcl.\zl ass,\zd\ namely $R\subseteq S$ if $xRy \supset xSy$\zl ,\zd\ e.g.\
\begin{tabbing}
\hspace{1.7em}\=father \hspace{1.2ex}\=$\subseteq$ ancestor\\[.5ex]
\> brother \>$\subseteq$ relative\\[.5ex]
\> smaller \>$\subseteq$ not greater \ud
\end{tabbing}

$\Big\lceil$These \ul operations \ud for rel.\zl ations\zd\ \sout{considered so far} (i\zl .\zd e.\ +\zl ,\zd\ $\cdot$\zl ,\zd\ $-$) are exactly analogous to the corresp.\zl onding\zd\ for classes and therefore will obey the same laws, e.g\zl .\zd\ $(R+S)\lfloor\cdot\rfloor T = R\cdot T + S\cdot T$. But in addition to them there are cert\zl ain\zd\ operations specific for relations and therefore more interesting\zl ,\zd\ e.g.\ for any \ul rel.\zl ation\zd\ \ud $R$ we can form what is called the \underline{inverse of $R$} (denoted by $\breve{R}$ \ul or $R^{-1}$ \ud ) where $\breve{R}=\hat{x}\hat{y}[yRx]$\zl ,\zd\ hence $x\breve{R}y \equiv yRx$\zl ,\zd\ i.e.\ if $y$\zl written over $x$\zd\ has the rel\zl ation\zd\ $R$ to $x$ then $x$ has the rel\zl ation\zd\ $\breve{R}$ $\mathbf{\llbracket 111. \rrbracket}$ to $y$.\zl ,\zd\ e.g.\
\begin{tabbing}
\hspace{1.7em}\=child = (parent)$^{-1}$\\[.5ex]
\> $x$ child $y$ $\equiv$ $y$ parent $x$\\[.5ex]
\> $< \;= (>)^{-1}$ \\[.5ex]
\> smaller = (greater)$^{-1}$\\[.5ex]
\> (nephew + niece) = (uncle + aunt)$^{-1}$
\end{tabbing}
There are also relations which are identical with their inverse \zl the following text until $I=I^{-1}$ is crossed out: e.g.\ identity \zl unreadable word, perhaps: to\zd\ $(=) =(=)^{-1}$\zl ,\zd\ bec\zl ause\zd\ $(x=y)\equiv (y=x)$ (in order to make the form more conspicuous one writes \ul also \ud $I$ for identity such that $I=I^{-1}$\zd \zl ,\zd\ i\zl .\zd e.\ $xRy\equiv yRx$\zl .\zd\ Such relations are called symmetric. \zl O\zd ther ex.\zl ample\zd\ (brother + sister) is sym.\zl metric\zd\ because - \zl \ldots;\zd\ brother is not sym.\zl metric,\zd\ sister is \ul n't \ud either. $\mathbf{\llbracket 112. \rrbracket}$ \zl The following text until $\Big\rfloor$ is crossed out: The op.\zl eration\zd\ of inverse obeys the law $(R^{-1})^{-1} =R$ \ul and is connected by laws of distr\zl ibution\zd\ with the former oper\zl ation\zd\ + .\zl ,\zd\ e.g.\ $(R+S)^{-1} =R^{-1}+S^{-1}$ \ud \sout{$(R+S)^{-1} =R^{-1}+S^{-1}$} .$\Big\rfloor$\zd\

Another oper.\zl ation\zd\ specific for rel.\zl ations\zd\ \ul and particularly important \ud is the so called relative prod.\zl uct\zd\ of two rel.\zl ations\zd\ ren.\zl dered\zd\ by $R|S$ and defined by
\[
R|S=\hat{x}\hat{y}[(\exists z)(xRz\; .\; zSy)]
\]
i\zl .\zd e\zl .\zd\ $R|S$ subsists between $x$ and $y$ if there is some obj.\zl ect\zd\ $z$ to which $x$ has the \zl r\zd\ el.\zl ation\zd\ $R$ and which has the rel.\zl ation\zd\ $S$ to $y$\zl ,\zd\ e.g.\
\begin{tabbing}
\hspace{1.7em}\=nephew = son$|$(brother or sister)
\end{tabbing}
$\mathbf{\llbracket 113. \rrbracket}$ $x$ is a nephew to $y$ if \sout{there is a person $z$ such that} $x$ \ul is \ud son of \ul some person \ud $z$ \sout{and th} which is brother or sister of $y$. In everyday langu\zl age\zd\ the prop\zl osition\zd\ $xRy$ is usually expressed by $x$ is an $R$ of $y$ \ul or $x$ is the $R$ of $y$ \ud \zl (e.g\zl .\zd\ \zl missing text\zd )\zd . Using this \zl unreadable text, perhaps in shorthand\zd\ we can say $xR|Sy$ means $x$ is an $R$ of an $S$ of $y$\zl ,\zd\ e.g.\ $x$ is a nephew of $y$ means $x$ is a son of a brother or sister of $y$\zl . O\zd ther example:
\begin{tabbing}
\hspace{1.7em}\=paternal uncle = brother$|$father\\[.5ex]
\` \ul Forts.\zl German: continued\zd \zl unreadable word, should be: p.\zd\ \textbf{119}. \ud
\end{tabbing}
The relative prod\zl uct\zd\ can also be applied to a relation and the same rel\zl ation\zd\ again \zl , i.e.\zd\ we can form $R|R$ (by def= $R^2$) square of a rel\zl ation,\zd\ $\mathbf{\llbracket 114. \rrbracket}$ e\zl .\zd g.\
\begin{tabbing}
\hspace{1.7em} \= \zl p\zd aternal grandfather = (father)$^2$\\[.5ex]
\> grandchild = (child)$^2$
\end{tabbing}
\zl S\zd imilarly we can form $(R|R)|R= R^3$\zl ,\zd\ e\zl .\zd g.\
\begin{tabbing}
\hspace{1.7em} \= great grandchild = (child)$^3$\\[.5ex]
\` \ul Forts.\zl German: continued\zd\ p\zl .\zd\ \textbf{117}. \ud
\end{tabbing}

$\Big\lceil$The relative product again follows laws very similar to the arithmetic one\zl '\zd s\zl ,\zd\ e.g.\
\begin{tabbing}
\hspace{1.7em} \= Associat\zl ivity\zd : \hspace{.5em} \=$(R|S)|T = R|(S|T)$\\[.5ex]
\> Distrib\zl utivity:\zd \>$R|(S+T)=R|S +R|T$\\[.5ex]
\> also\> $R|(S\lfloor\cdot\rfloor T)\subseteq R|S \cdot R|T$
\end{tabbing}
\zl on the right of the formulae just displayed, there is a pale, unreadable and crossed out text with formulae, probably a derivation of some of the displayed formulae\zd\ but not commutativity
\begin{tabbing}
\hspace{1.7em} \= $R|S=S|R$ is false\\[.5ex]
\> brother$|$father $\neq$ father$|$brother
\end{tabbing}
since paternal uncle \zl unreadable text\zd\ but $\neq$\zl is not\zd\ father\zl .\zd $\Big\rfloor$

\vspace{1ex}

\zl The whole of pages \textbf{115}.\ and \textbf{116}.\ are crossed out.\zd\

\vspace{1ex}

$\mathbf{\llbracket 115. \rrbracket}$ Identity $I$ is a unity for this prod.\zl uct,\zd\ i\zl .\zd e.\ $R|I=I|R= R$ bec.\zl ause
\begin{tabbing}
\hspace{1.7em} \= $xR|Iy\;$\=$\equiv xIz\; .\; zRy$ for some $z$\\[.5ex]
\>\>$\equiv xRy$\\[1ex]
\> Monotonicity: $R\subseteq S, P\subseteq R\supset R|P\subseteq S|Q$
\end{tabbing}

$\Big\lceil$Furthermore the class of all \sout{\zl unreadable text\zd }\, obj\zl ects\zd\ which have the rel\zl ation\zd\ $R$ to some \ul obj\zl ect\zd\ \ud $y$ is called domain \zl unreadable text\zd\ $\mathbf{D}`R= \hat{x}[(\exists y)xRy]$ and the class of all obj.\zl ects\zd\ to which some obj\zl ect\zd\ has the rel.\zl ation\zd\ is called converse domain $\mathbf{C}`R= \hat{x}[(\exists y)yRx]$ so that $\mathbf{C}`R= \mathbf{D}`R^{-1}$\zl ,\zd\ e.g.\ $\mathbf{D}`$(father) = men that have children\zl .\zd\

$\mathbf{\llbracket 116. \rrbracket}$ $\Big\lceil$In ord.\zl inary\zd\ language this class is also denoted by ,,\zl ``\zd father''. So you see in everyday lang\zl uage\zd\ the same word is used for two diff\zl erent\zd\ things\zl ,\zd\ a rel.\zl ation\zd\ and its domain\zl :\zd\
\begin{tabbing}
\hspace{1.7em} \= $\mathbf{C}`$father = class of \ul all \ud men (except Adam and Eve)\\[.5ex]
\> $\mathbf{D}`$(brother or sister) = $\mathbf{C}`$(brother or sister) \\
\` = class of men which have a brother or sister\zl ,\zd\ \\
\` \sout{\zl unreadable symbol\zd }\, hence\\[.5ex]
\> Man$-\mathbf{D}`$(brother or sister) = unique children\\[.5ex]
\> $\mathbf{D}`R+\mathbf{C}`R = \mathbf{C}`R$\zl ,\zd\ \\[.5ex]
\> $\mathbf{C}`$father = class of all men
\end{tabbing}

An important property which belongs to many relations is ,,\zl ``\zd Transitivi\-ty''\zl .\zd $\Big\rfloor$

$\mathbf{\llbracket 117. \rrbracket}$ A rel.\zl ation\zd\ $R$ is called transitive if
\begin{tabbing}
\hspace{1.7em} $(x\lfloor ,\rfloor y\lfloor ,\rfloor z)[xRy\; .\; yRz\supset xRz]:\!\!\lfloor\equiv\rfloor \; R$ is transitive
\end{tabbing}
In other words if an $R$ of an $R$ of $z$ is an $R$ of $z$\zl ;\zd\ e.g.\ brother is transitive\zl ,\zd\ a brother of a brother of a person is a brother of this person\zl ,\zd\ in other words
\begin{tabbing}
\hspace{1.7em} $x$ brother $y\; .\; y$ brother $z \supset x$ brother $z$
\end{tabbing}
Smaller is also transitive\zl ,\zd\ i.e.\
\begin{tabbing}
\hspace{1.7em} $x<y\; .\; y<z \supset x<z$
\end{tabbing}
\ul Very many rel\zl ations\zd\ in math.\zl ematics\zd\ are transit\zl ive:\zd \zl unreadable word\zd\ \zl ,\zd\ congr.\zl uence,\zd\ $\mid\mid$ \zl parallelism,\zd\ isom\zl orphism,\zd\ ancestor\zl .\zd\ \ud Son is not transitive, a son of a son of a person is not a son of a person\zl .\zd\ \zl The following sentence, under a line drawn in the text, is crossed out: The relation of son even has the opposite prop\zl erty\zd\
\begin{tabbing}
\hspace{1.7em} $x$ son $y\; .\; y$ son $z \supset\;\sim(x$ son $y)$ \zl $\sim(x$ son $z)$\zd $\Big\rfloor$\zd\
\end{tabbing}
$\mathbf{\llbracket 118. \rrbracket}$ $\Big\lceil$t\zl T\zd herefore called intransitive\zl ;\zd\ friend is an ex\zl ample\zd\ of a relation which is neither transitive nor intransitive. A friend of a friend of $x$ is not always a friend of $x$\zl ,\zd\ but is sometimes a friend of $x$. By means of the previously introduced op.\zl eration\zd\ transitiv\zl ity\zd\ \ul can be \ud expressed by
\begin{tabbing}
\hspace{1.7em}\= $R^2 \subseteq R$ \quad bec.\zl ause\zd\ \\[.5ex]
\> $ xR^2 y \:.\!\supset (\exists z)\lfloor (\rfloor xRz\; .\; zRy) \supset xRy$
\end{tabbing}
if $R$ is transitive\zl ,\zd\ but also vice versa if $R$ satisfies the cond.\zl ition\zd\ $R^2 \subseteq R$ then $R$ is trans.\zl itive\zd\
\begin{tabbing}
\hspace{1.7em}\= $xRy\; .\; yRz \supset xR^2z \supset xRz$$\Big\rfloor$
\end{tabbing}
\zl The following inserted sentence is crossed out: Ex.\zl amples\zd\ of trans\zl itive\zd\ rel.\zl a\-tions:\zd\ sim\zl symmetry\zd , congr\zl uence\zd , \zl unreadable word, presumably in shorthand\zd , =\zl equality\zd , $\mid\mid$ \zl parallelism\zd , ancestor, \zl unreadable word, presumably in short\-hand\zd ,\zd\

$\mathbf{\llbracket 119. \rrbracket}$ A very important prop.\zl erty\zd of relations is the following one: A binary rel.\zl ation\zd\ $R$ is called one-many if for any obj.\zl ect\zd\ $y$ there exists at most one obj.\zl ect\zd\ $x$ such that $xRy$\zl :\zd\
\begin{tabbing}
\hspace{1.7em}$(x\lfloor ,\rfloor y\lfloor ,\rfloor z)[xRy\; .\; zRy\supset x=z]\equiv \; R$ is one\zl -\zd many
\end{tabbing}
and many\zl -\zd one if $R^{-1}$ is one\zl -\zd many\zl ;\zd\ e.g.\ father is one\zl -\zd many\zl ,\zd\ every obj.\zl ect\zd\ $x$ can have at most one father\zl ,\zd\ it can have no father if it is no man\zl ,\zd\ but it never has two \zl unreadable text in parentheses\zd\ \ul or more \ud fathers. The rel\zl ation\zd\ $<$ is not one\zl -\zd many\zl :\zd\ for any nu.\zl mber\zd\ there are many diff.\zl erent\zd\ nu.\zl mbers\zd\ $<$\zl smaller than it\zd .

\zl The following text at the end of p.\ \textbf{119}.\ is crossed out, though its continuation on p.\ \textbf{120}.\ is not: $\Big\lceil$\zl deleted: or e.g.\zd\ t\zl T\zd he rel\zl ation\zd\ $x$ is the reciproc.\zl al\zd\ of n.\zl umber\zd $y$ is one\zl -\zd many. Every nu.\zl mber\zd\ has at most$\Big\rfloor$\zd\ $\mathbf{\llbracket 120. \rrbracket}$ $\Big\lceil$one reciprocal. Some numbers have no reciprocal\zl ,\zd\ namely 0 (but that makes no difference). The rel.\zl ation\zd\ of reciprocal is at the same time many\zl -\zd one\zl ;\zd\ such relations are called one\zl -\zd one\zl .\zd $\Big\rfloor$

\zl The following inserted text is crossed out: The inverse rel\zl ation\zd\ ,,\zl ``\zd son'' is not one\zl -\zd many\zl ;\zd\ there \sout{\zl unreadable word\zd }\, can be several persons having this relation of son to one person. A rel.\zl ation\zd\ which is one\zl -\zd many and many\zl -\zd one is called one-one\zl .\zd \zd\

The relation of husband \ul in Christian \zl coun\zd tries \ud e.g.\ is an \sout{other} ex.\zl ample\zd\ of a one-one relation. The rel.\zl ation\zd\ smaller is neither one-many nor many\zl -\zd one\zl ; for\zd\ any nu.\zl mber\zd\ there exist many different nu.\zl mbers\zd\ smaller than it and many diff\zl erent\zd\ numbers greater than it. \ul \sout{One-many doesnot mean that} \ud

One-many-ness can also be defined for polyadic relations $\mathbf{\llbracket 121. \rrbracket}$ \zl \sout{namely}\zd . \ul A triadic rel.\zl ation\zd\ \ud $M$ is called one\zl -\zd many if
\begin{tabbing}
\hspace{1.7em}$(x\lfloor ,\rfloor y\lfloor ,\rfloor z\lfloor ,\rfloor u)[xM(zu)\; .\; yM(zu)\supset x=y]$
\end{tabbing}
e.g.\ $\hat{x}\hat{y}\hat{z}(x=y+z)$\zl ,\zd\ $\hat{x}\hat{y}\hat{z}[x -\lfloor =\rfloor \frac{y}{z}]$ have this prop\zl erty\zd . For any two nu.\zl mbers $y$ and $z$\zd there exists at most one $x$ which is the sum or difference\zl .\zd\ $\hat{x}\hat{y}(x$\zl = deleted\zd\ \sout{$\surd y)$} is a square root of $y$) is not one\zl -\zd many because there are in gen.\zl eral\zd\ two different nu.\zl mbers\zd\ which are square roots of $y$. \sout{but $\hat{x}\hat{y}[x=y^2]$ is onemany} You see the one\zl -\zd many \sout{dyadic} relations are exactly the same thing which is called ,,\zl ``\zd functions'' in math\zl ematics\zd . The dyadic one\zl -\zd many relations are the f\zl u\zd nct\zl ions\zd\ with one argument \ul as e\zl .\zd g.\ $x^2$\zl ,\zd\ \ud the $\mathbf{\llbracket 122. \rrbracket}$ triadic one\zl -\zd many relations are the funct\zl ions\zd\ with two arg.\zl uments\zd\ as e.g.\ $x+y$\zl .\zd\ \zl The inserted text that follows from $\Big\lceil$ to $\Big\rfloor$ is crossed out.\zd\ \ul $\Big\lceil$Relations which are not one\zl -\zd many may also be thought of as f\zl u\zd nct.\zl ions,\zd\ but as many\zl -\zd valued f\zl u\zd nct\zl ions,\zd\ e.g\zl .\zd\ the log\zl arithm\zd\ for complex\zl full stop deleted\zd\ nu.\zl mbers\zd\ $\log x$ has inf.\zl initely\zd\ many values for a given $x$\zl .\zd\ There this symb\zl ol\zd\ $\log$ \zl full stop deleted\zd\ from the log.\zl ical\zd\ standpoint denotes a \zl \sout{not}\zd\ \sout{one many} dyadic relation \ul which is not one\zl -\zd many\zl .\zd\ \ud This rel\zl ation\zd\ subsists between two nu.\zl mbers\zd\ $y$\zl ,\zd\ $x$ if $y$ is one of the values of the log.\zl arithm\zd\ for the argument $x$. But if the word f\zl u\zd nct\zl ion\zd\ is used without further specification then always single\zl -\zd valued f\zl u\zd nct\zl ion\zd\ \sout{are} \ul is \ud meant in math.\zl ematics\zd ; \sout{and a} \ul the term \ud ,,\zl ``\zd single-valued f\zl u\zd nct\zl ion\zd '' \sout{is} \ul denotes exactly \ud the same \ul thing \ud as \ul the term \ud \sout{a} ,,\zl ``\zd one\zl -\zd many relation''.$\Big\rfloor$ \ud

In \ul order to \ud make statements about f\zl u\zd nct\zl ions,\zd\ \ul i\zl .\zd e.\ one\zl -\zd many rel.\zl ations\zd\ \ud it is very convenient to introduce a notation usual in mathematics and also in everyday lang.\zl uage;\zd\ namely $R`x$ \sout{means} \ul denotes \ud the $y$ which has the rel.\zl ation\zd\ $R$ to $x$\zl ,\zd\ i\zl .\zd e.\ the $y$ such that $yRx$ provided that this $y$ exists and is unique. Similarly for a triadic rel.\zl ation\zd\ $M`(yz)$ means\zl denotes\zd\ the $x$ such that\zl \ldots\zd\ Inst\zl ead\zd\ of this also $yMz$ is written\zl ,\zd\ e.g.\ + denotes a triadic rel.\zl ation\zd\ between $\mathbf{\llbracket 123. \rrbracket}$ numbers \ul (sum) \ud and $y+z$ denotes the number which has this triadic rel.\zl ation\zd\ to $y$ and $z$ \zl \ul provided that it exists \ud . \zl The following inserted sentence from $\Big\lceil$ to $\Big\rfloor$ is crossed out: $\Big\lceil$the statement that it exists is \zl unreadable text, perhaps: seen\zd\ by $E!R`x$ (e.g\zl .\zd\ $E!\frac{1}{2}$\zl ,\zd\ $\sim E!\frac{1}{0}$$\Big\rfloor$\zd\ \sout{(This notation is not ambiguous \ul \zl unreadable text\zd\ \ud )} In everyday language the $`$ is expressed by the words The\ldots\ of\zl ,\zd\ e.g\zl .\zd\ t\zl T\zd he sum of $x$ and $y$\zl ,\zd\ The father of $y$.

\zl The text that follows until the end of p.\ \textbf{125}.\ is crossed out in the manuscript.\zd\ $\Big\lceil$There is only one tricky point \zl in this\zd\ notation. Namely \zl w\zd hat meaning are we to \zl \sout{attribute}\zd\ \ul assign \ud to propositions containing this symbol $R`x$ if there does\zl $\,$\zd not exist \zl a\zd\ unique $y$ such that $yRx$ (i\zl .\zd e.\ none or several)\zl ,\zd\ e.g.\ The present king of$\Big\rfloor$ $\mathbf{\llbracket 124. \rrbracket}$ $\Big\lceil$France is bald. We may convene that such propositions are meaningless \zl full stop deleted\zd\ (neither true nor false). But that has certain undesirable consequences, namely whether or not the present king of France exists or not is an empirical question\zl .\zd\ Therefore it would depend on an empirical fact whether or not this sequence of words is a meaningful statement or nonsense whereas one should expect that it can depend only on the grammar of the language concerned whether something makes sense\zl .\zd $\Big\rfloor$ $\mathbf{\llbracket 125. \rrbracket}$ $\Big\lceil$\sout{Therefore eg.} Russell makes the convention\zl \sout{s}\zd\ that such statements are false and not meaningless. The conv.\zl ention\zd\ is: That every \underline{atomic prop.\zl osition\zd} in which such an $R`x$ (describing something nonexistent) occurs is false\zl ,\zd\ i\zl .\zd e.\
\begin{tabbing}
\hspace{1.7em}$\varphi(R`x)\equiv (\exists y)[(z)[zRx\equiv\ z=y]\; .\; \varphi(y)]$
\end{tabbing}
\zl \sout{e.g.}\zd $\Big\rfloor$

$\mathbf{\llbracket 126. \rrbracket}$ All afore\ul mentioned \ud notions \sout{defined} of the calc.\zl ulus\zd\ of classes and relations are themselves relations; e.g.\ $\alpha\subseteq\beta$ is a binary rel.\zl ation\zd\ between classes\zl ,\zd\ $\alpha + \beta$ is a dyadic f\zl u\zd nct.\zl ion,\zd\ i\zl .\zd e.\ a triadic rel\zl ation\zd\ between classes (which subsists between $\alpha,\beta,\gamma$ if $\gamma=\alpha+\beta$)\zl .\zd\ The op.\zl eration\zd\ of inverse is a rel.\zl ation\zd\ between relations subsisting between $R$ and $S$ if $R=S^{-1}$ \zl or\zd\ the rel.\zl ative\zd\ prod.\zl uct\zd\ is a triadic rel.\zl ation\zd\ between relations subsisting between $R,S,T$ if $R=S|T$. Symmetry defines a cert\zl ain\zd\ class of rel.\zl ations\zd\ (the \ul class of \ud sym.\zl metric\zd\ relations)\zl .\zd\ So we see that we have obtained a $\mathbf{\llbracket 127. \rrbracket}$ new kind of \sout{objects} concepts (called concepts of second type or sec.\zl ond\zd\ order) which refer to the concepts of first order\zl ,\zd\ \ul i\zl .\zd e.\ which expresses properties of conc.\zl epts\zd\ of first order or rel.\zl ations\zd\ between conc.\zl epts\zd\ of first order \zl full stop deleted\zd\ or to be more exact prop\zl erties\zd\ and rel\zl ations\zd\ of extensions of concepts of first order\zl .\zd\ But this is not very essential since we can define corresponding conc.\zl epts\zd\ which express prop.\zl erties\zd\ and rel\zl ations\zd\ of the pred.\zl icates\zd\ themselves\zl ,\zd\ e.g\zl .\zd\ $\chi$\zl written over $\psi$\zd\ sum of $\varphi,\psi$ \zl if\zd\ $\chi(x)\equiv \varphi(x)\vee\psi(x)$ etc\zl .\zd\ \ud

And it is possible to (go on) continue in this way\zl ,\zd\ i\zl .\zd e.\ we can define concepts of third \sout{order or} \ul type or \ud order, which refer to the concepts of sec.\zl ond\zd\ order. \sout{as eg} \ul An example would be: \ud ,,\zl ``\zd mutually exclusive''\zl ;\zd\ a class of classes $\cal U$\zl ,\zd\ \ul i\zl .\zd e\zl .\zd\ a class whose el\zl ements\zd\ are themselves classes\zl ,\zd\ \ud is called a mut.\zl ually\zd\ excl.\zl usive\zd class of classes if $\alpha,\beta\:\varepsilon\:{\cal U} \supset \alpha\cdot\beta = \Lambda$. This concept of ,,\zl ``\zd mut.\zl ually\zd\ excl.\zl usive\zd\ class of classes'' expresses a prop.\zl erty\zd\ of classes \ul of classes\zl ,\zd\ i\zl .\zd e\zl .\zd\ of an obj\zl ect\zd\ of 3\zl third\zd\ order\zl ,\zd\ therefore is \ud of third order. \zl On the right of p.\ \textbf{127}.\ one finds the following text to be inserted it is not clear where: e\zl .\zd g.\ the word community of Am.\zl erica\zd\ or army the present states \sout{of} on the earth\zd\ So \ul you see\zd\ in this way we get a whole hierarchy of concepts $\mathbf{\llbracket 128. \rrbracket}$ which is called the hierarchy of types\zl .\zd\ In fact there are two diff\zl erent\zd\ hierarchies of types \sout{\zl unreadable symbol\zd }\zl ,\zd\ namely the hierarchy of ext.\zl ensions\zd\ and the hierarchy of predicates. \zl The following sentence is crossed out: So far I have spoken only of the former\zl ;\zd\ the latter would begin with predicates\zl ,\zd\ then \ul we have \ud predicates of predicates (i\zl .\zd e.\ prop.\zl erties\zd\ of pred\zl icates\zd\ or relations between pred.\zl icates\zd ) \zl \dots\zd \zd

\zl Following an unreadable symbol, perhaps in shorthand, there is a vertical line on the left margin for the remaining text on p.\ \textbf{128}.\ and the whole text on p.\ \textbf{129}.\zd\ An interesting ex.\zl ample\zd\ of predicates of highe\zl r\zd\ \ul type are \ud the \sout{nat} natural numbers. According to Russell and Frege the nat\zl ural\zd\ nu.\zl mbers\zd\ are properties of pred\zl icates\zd . \sout{\zl unreadable text\zd }\, If I say e.g\zl .\zd : There are eight planet\zl \sout{e}\zd s \zl full stop deleted\zd \zl ,\zd\ this expresses a property of the predicate $\mathbf{\llbracket 129. \rrbracket}$ ,,\zl ``\zd planet''. So the nu.\zl mber\zd\ 8 can be defined to be a property \ul of predicates \ud which belongs to a pred.\zl icate\zd\ $\varphi$ if there are exactly 8 obj\zl ects\zd\ falling under this pred\zl icate\zd . If this definition is followed up it turns out that all notions of arithm\zl etic\zd\ can be defined in terms of logical notions and that the laws of arithm.\zl etic\zd\ can be derived from the laws of logic except for one thing\zl ,\zd\ namely \ul for building up arithmetic\zd\ one needs the prop\zl osition\zd\ that there are infinitely many obj.\zl ects,\zd\ which cannot be proved from the ax.\zl ioms\zd\ of logic.

$\mathbf{\llbracket 130. \rrbracket}$ The lowest layer in the hierarchy of types described are the individuals or obj.\zl ects\zd\ of the world\zl ;\zd\ what these ind.\zl ividuals\zd\ are is an extralog.\zl ical\zd\ question which depends on the theory of the world which we assume\zl ;\zd\ in a material.\zl ist\zd\ theory it would be the atoms or the points of space and time\zl ,\zd\ i\zl written over I\zd n a spiritualist theory it would be the spirits and so on. As to the higher types (classes\zl ,\zd\ classes of classes\zl ,\zd\ \ul predicates of pred.\zl icates\zd\ \ud etc\zl .\zd ) each \ul type \ud must be distinguished very carefully \ul from any other \ud as can be shown e.g\zl .\zd\ by the foll\zl owing\zd\ $\mathbf{\llbracket 131. \rrbracket}$ example. If $a$ is an obj\zl ect\zd\ one can form the class whose only element is $a$ (denoted by $\riota\, ` a$)\zl .\zd\ So this $\riota\, ` a$ would be the extension of a predicate, which belongs to $a$ and only to $a$. \sout{Now} It \sout{should be} \ul is \ud near at hand to identify this $a$ and $\riota\, ` a$\zl ,\zd\ i.e.\ to assume that the obj\zl ect\zd\ $a$ and the class whose only element is $a$ are the same. However it can be shown that this is not admissible\zl ,\zd\ i.e.\ it would lead to contradictions to $\mathbf{\llbracket 132. \rrbracket}$ assume this identity $\riota\, ` a=a$ \zl comma from the manuscript deleted\zd\ to be generally true because \zl comma from the manuscript deleted\zd\ if we take for $x$ a class (which has several elements) then certainly $\riota\, ` \alpha$ and $\alpha$ are distinct from each other; since $\riota\, ` \alpha$ is a class which has only one element\zl ,\zd\ \ud namely $\alpha$\zl ,\zd\ \ud whereas $\alpha$ is a class which has several elements\zl ,\zd\ so they are certainly distinct from each other. But \zl \sout{on the other hand}\zd\ although we have to distinguish very carefully between the different type\zl s\zd\ there is \ul on the other hand \ud a very close analogy between the diff.\zl erent\zd\ type\zl s\zd . \sout{Sup} E.g.\ classes of individuals $\mathbf{\llbracket 133. \rrbracket}$ and classes of classes \ul of individuals \ud will obey exactly the same laws. For both of them we can define an \sout{\zl unreadable letter\zd }\, and a multiplication and the same laws of calculus \ul will hold \ud for them. Therefore it is desirable not to formulate these laws separately for classes of classes and classes of individuals, but to introduce a general notion of a class comprising \ul in it \ud all those particular cases\zl :\zd\ classes of ind\zl ividuals\zd , classes of rel\zl ations\zd , classes of classes etc. And it was actually in $\mathbf{\llbracket 134. \rrbracket}$ this way that the logistic calculus was first set up (with such a general notion of a class \ul and \ud of a predicate \ul and \ud of a relation and so on \sout{comprising all} embracing under it all types) and this way also corresponds \sout{certainly} more to \sout{the} natural thinking. In ordinary \ul language e.g.\ \ud we have such a general notion of a class without a distinction of \sout{the} different types\zl .\zd\

\zl new paragraph\zd\ The more detailed working out of logic on this \ul typeless \ud base \sout{(in natural thinking)} has led to \ul the discovery \zl \sout{of}\zd\ \ud of the most interesting $\mathbf{\llbracket 135. \rrbracket}$ facts in modern logic. Namely to the fact that the evidences of natural thinking are not consistent with themselves\zl ,\zd\ i\zl .\zd e.\ lead to contradictions which are called ,,\zl ``\zd logical paradoxes''\zl .\zd\ The first of these contradictions was \sout{found} discovered by the mathematician Burali-Forti \zl in\zd\ 1897. A few years later Russell produced a similar contradiction which however \sout{was cleaned} \ul avoided \ud \zl \sout{the}\zd\ the \ul un \ud essential mathematical by\zl -\zd work \ul of Burali-Forti\zl '\zd s contrad.\zl iction\zd\ \ud and showed the real logical structure of the contradiction \ul much clearer \ud . This \sout{Ru} so\zl hyphen deleted\zd\ $\mathbf{\llbracket 136. \rrbracket}$ called \underline{Russell} \zl hyphen deleted\zd\ \underline{paradox} has remained \ul up to now \ud the classical example of a logical paradox and I want to explain it now in \sout{all} detail\zl .\zd\ I shall first \ul enumerate \ud some apparently evident propositions from which the paradox follows in a few steps\zl .\zd\

The paradox under consideration involves \ul only \ud the following notions\zl :\zd\
\begin{itemize}
\item[1.] object in the most general sense\zl ,\zd\
which embraces everything that can be made an object of thinking\zl ;\zd\ in part.\zl icular\zd\ it embraces the indiv.\zl iduals\zd , classes, pred\zl icates\zd\ of all types
\end{itemize}

\zl at the bottom of this page: Forts.\ Heft\zl German: continued in Notebook\zd\ VII\zd\

$\mathbf{\llbracket 137. \rrbracket}$ \zl at the top of this page: Heft\zl German: Notebook\zd\ VII.\zd\

\begin{itemize}
\item[2.] monadic predicate (briefly pred\zl icate\zd )\zl ,\zd\ also in the most general sense comprising \ul predicates of ind\zl ividuals\zd\ as well as \ud predicates of predicates etc. And this term \zl dash deleted\zd\ predicate is to be so understood that it is an essential requirement of a predic\zl ate\zd \zl comma from the manuscript deleted\zd\ that it is well\zl -\zd defined for any object \ul whatsoever \ud whether the given predicate belongs to it or \zl not\zd\
\end{itemize}

Now of these two notions ,,\zl ``\zd object'' and ,,\zl ``\zd predicate'' we have the following apparently evident propositions\zl :\zd\
\begin{itemize}
\item[1.] \uline{If $\varphi$ is a pred\zl icate\zd\ and $x$ an obj\zl ect\zd\ then it is uniquely det\zl er\-mined\zd\ } \underline{whether $\varphi$ belongs to $x$ or not}.
\end{itemize}
Let us denote the prop.\zl osition\zd\ $\varphi$ bel.\zl ongs\zd\ to $x$ by $\varphi(x)$\zl .\zd\ So we have \zl \sout{)}\zd\ if $\varphi$ is a well\zl -\zd def.\zl ined\zd\ pred.\zl icate\zd\ and $x$ \zl \sout{)}\zd\ an obj\zl ect\zd\ then $\varphi(x)$ is always a meaningful\zl \sout{l}\zd\ prop.\zl osition\zd\ $\mathbf{\llbracket 138. \rrbracket}$ which is either true or false\zl .\zd\
\begin{itemize}
\item[2.] Vice versa\zl :\zd\ If we have a combination of words or symbols \ul $A(x)$ \ud which contains the letter $x$ and is such that it becomes a \ul meaningful \ud prop\zl osi\-tion\zd\ for any arbitrary object which you substit.\zl ute\zd\ for $x$ then $A(x)$ defines a cert\zl ain\zd\ predicate $\varphi$ which belongs to an obj.\zl ect\zd\ $x$ if and only if $A(x)$ is true\zl .\zd\
\end{itemize}
\ul So the assumption means that \zl \sout{ie.}\zd\ if you subst.\zl itute\zd\ for $x$ the name of an arb.\zl itrary\zd\ object then it is always uniquely determined whether the resulting propos\zl ition\zd\ is true or false\zl .\zd\ \ud

\zl The first item numbered 3 and the text which follows it until the page ends with ``whatever $x$'' has a big square bracket on its left margin.\zd\
\begin{itemize}
\item[3.] It is uniquely determined of any obj\zl ect\zd\ whether or not it is a pred\zl i\-cate\zd .
\end{itemize}
Let us denote by \sout{P and unreadable symbol} $P(x)$ the prop\zl osition\zd\ ,,\zl ``\zd $x$ is a predicate'' so that $P$(red)\zl ,\zd\ $\sim P$(smaller)\zl ,\zd\ $\sim P$(New York)\zl ;\zd\ then by 3 \sout{for} $P(x)$ is always a meaningful prop.\zl osition\zd\ whatever $x$ $\mathbf{\llbracket 139. \rrbracket}$ may be\zl \sout{]}\zd \zl .\zd\
\begin{itemize}
\item[3.\zl 4.\zd ] \underline{Any predicate is an obj\zl ect\zd }.
\end{itemize}

I think these 3\zl written over 4; in the next paragraph four assumptions are mentioned (see also the corresponding four assumptions on pp.\ \textbf{138}.-\textbf{140}. of Notebook VII), so it should be: four\zd\ prop.\zl ositions\zd\ are all evident to natural thinking\zl .\zd\ \zl \sout{[}\zd\ 1 and 2 can be considered as a def\zl inition\zd\ of the term predicate and 3 says that the notion of pred.\zl icate\zd\ thus defined is well\zl -\zd defined.\zl \sout{]}\zd\

And now let us consider the following statement ${P(x)\: .\: \sim x(x)}$ that means $x$ is a predicate and it belongs to $x$ (i\zl .\zd e.\ to itself). According to our \ul four ass\zl umptions\zd\ that is a meaningf\zl ul\zd\ pro\zl position\zd\ which is either true or false whatever you subst\zl itute\zd\ for $x$. \zl N\zd amely\zl ,\zd\ \ul at first by 3 it is uniquely det\zl er\-mined:\zd\ \ud if you $\mathbf{\llbracket 140. \rrbracket}$ subst.\zl itute\zd\ for $x$ something which is not a pred.\zl icate\zd\ it becomes false\zl ,\zd\ if you subst.\zl itute\zd\ for $x$ a pred.\zl icate\zd\ then \ul $P(x)$ is true but \ud $x(y)$ is either true or false for any obj.\zl ect\zd\ $y$ \zl written over $x$\zd\ by 1. \zl B\zd ut $x$ is a pred.\zl icate,\zd\ hence an obj.\zl ect\zd\ by ass\zl umption\zd\ 3\zl 4,\zd\ hence $x(x)$ is either true or false\zl ,\zd\ hence the whole statement is always meaningful\zl ,\zd\ i\zl .\zd e.\ either true or false\zl .\zd\ Therefore by 2 it defines a cert\zl ain\zd\ pred.\zl icate\zd\ $\Phi$ \sout{which} such that $\Phi(x)\,{\equiv \atop \mbox{\scriptsize\rm means}}\, P(x)\: .\: \sim x(x)$\zl .\zd\ \ul \sout{$\equiv x$ is pro impredicable} \ud

\zl Next comes a page again numbered \textbf{140}. with a crossed out text.\zd\
\[
\Phi(\Phi)\equiv P(\Phi)\: .\sim \Phi(\Phi)
\]
But this leads immediately to a contradiction since this equ.\zl ation\zd\ means two implications
 \[
\Phi(\Phi)\supset P(\Phi)\: . \sim \Phi(\Phi)
\]
\[
P(\Phi)\: \lfloor .\rfloor \sim \Phi(\Phi) \supset \Phi(\Phi)
\]

\zl The last two pages of Notebook VI are not numbered. These two pages will not be entirely reproduced here, since they contain only rather unconnected notes and jottings, presumably for exercises, written without much order and care. These notes will however be described here up to a point.

On the first of these pages is first an exercise involving reduction to normal form, in which one finds the following (the unsystematically written $.$ is here deleted, as well as the unreadable crossed out beginning of the third line):
\begin{tabbing}
\hspace{1.7em} \=$c(b+y\bar{a})+d(\bar{b}(\bar{y}+a))\neq 0$\\[.5ex]
\>$d\bar{b}a+cb+yc\bar{a}+d\bar{b}\bar{y}\neq 0$\\[.5ex]
\>$cb+c\bar{a}+d\bar{b}\neq 0$\hspace{1.7em} $ab\neq 0$\\[.5ex]
\>$b+b\bar{a}+a\bar{b}\neq 0$\\[.5ex]
\>$d=0$\hspace{1.7em} $c=a$
\end{tabbing}
In the first line, and above it, one finds in the margin:
\begin{tabbing}
\hspace{1.7em} \=$a$ \hspace{1.7em}\=$b$\hspace{2.7em} $b\neq 0$\\[.5ex]
\hspace{1.65em}$(c$\>\>$d)$\\[.5ex]
\>$b$\>$a$
\end{tabbing}
In the remainder of this page one finds ``$x$ is a parent of $y$'' and ``\underline{child} = son or daughter''. The rest is either in shorthand, or it is unreadable, or it is crossed out.

On the remaining, last page, of Notebook VI, one finds first a few lines, mostly in shorthand (once crossed out), in which one finds also: $(x)\varphi(x)$, \sout{equally shaped}, tautological entailment, out of fashion, unfeasible, fail, permitted to take, unpracticable. Next, at the end of the notebook, one finds notes, rather difficult to read, written without much order and care, and partly erased, which involve some equations, Boolean expressions, perhaps a syllogism, a Venn-Euler diagram, ``$2^4-1(=15)$'', ``$2^{15}=32000$'' ($2^{15}$ is 32768) and ``$2^{15}-2$''. At the top of this last page, one finds the caption ``illegible text'', presumably put by the archive where the manuscript is preserved.\zd\

\section{Notebook VII}\label{0VII}
\pagestyle{myheadings}\markboth{SOURCE TEXT}{NOTEBOOK VII}
\zl Folder 65, on the front cover of the notebook ``Logik Vorl.\zl esungen\zd\ \zl German: Logic Lectures\zd\ N.D.\ \zl Notre Dame\zd\ VII''\zd\

\zl Notebook VII starts with nine, not numbered, pages of numbered remarks and questions, more than eighty of them, partly unreadable, partly in shorthand, and all seemingly not closely related to the remaining notes for the course. They will be reproduced here up to a point only.\zd\

\vspace{1ex}

$\mathbf{\llbracket new\, page \rrbracket}$ $-$1. Every \zl unreadable abbreviated word\zd\ prop\zl osition\zd\ is true\zl .\zd\

$-$2. Everyone \sout{\zl not\zd }\, (Christ\zl ian\zd , cathol\zl ic\zd ) bel.\zl ieving\zd\ \ul the neg\zl ation\zd\ of \ud \sout{in} dogm\zl a\zd\ commits a mortal sin\zl .\zd\

$-$3. Everyone not bel\zl ieving\zd\ a dogma although he knows that it is \zl dogma\zd\ commits a mortal sin\zl .\zd\

$-$4. Everyone teaching \ul publicly \ud the neg\zl ation\zd\ doctr.\zl inal\zd\ prop.\zl osi\-tion\zd\ as the truth (although) commits a\ldots

$-$5\zl .\zd\ Everyone asserting privately\ldots

$\circ$ 6. The world was existed appr.\zl oximately\zd\ 6000 years \zl 9\zd vid\zl .\zd\ 25\zl ).\zd\

$\circ$ 7. The sky is \ul made \ud of solid material\zl .\zd\

$\circ$ 8. There exist angels and evil spirits\zl .\zd\

$\circ$ 9. Some of the \zl unreadable text\zd\ are caused by evil spirits\zl .\zd\

$\circ$ 10\zl .\zd\ \sout{Hypnot} The phen\zl omena\zd\ of hypnotism \ul (telepathy \ul telekinesis\zl ,\zd\ prophecy \ud ) \ud are caused by evil spirits (spirits?)

$-1'$. If $A$ is a dogma at some time it is a dogma at any later time.

$-2'$. If $A$ is a \zl unreadable abbreviated word, same as in 1.\zd\ prop\zl osition\zd\ at some time it is\ldots

$\mathbf{\llbracket new\, page \rrbracket}$ $\blacksquare$ 11. Will logic and mathematics be the same in the after the end of this world\zl ?\zd\

$\blacksquare$ 12. \sout{Woul \zl or ``Word''\zd\ The death of Christ was} It was in the power of Christ (inqu. R homo) not to dy for manhood\zl .\zd\ \zl Sentence 12 is in big square brackets.\zd\

$\blacksquare\circ$ 13. It would have been no sin of Christ if he had not \sout{died for} sacrificed himself for manhood\zl .\zd\

$\circ$ 14. Can an infidel \ul cath\zl olic\zd\ \ud priest \sout{deal out} administer sacraments if he keeps the outward form\zl ?\zd\

$\circ$ 15. Can an infidel make a valid baptism if he keeps the form\zl ?\zd\

$\blacksquare$ 16. \sout{Does} Is everyone not baptis\zl z\zd ed and living after Christ's death \sout{go to hell} damned\zl ?\zd\

$\blacksquare$ 17. Does everybody baptis\zl z\zd ed which has committed a mort\zl al\zd\ sin without being $\circ$ absolved by a cath\zl olic\zd\ priest \sout{go to hell} damned\zl ?\zd\

$\mathbf{\llbracket new\, page \rrbracket}$ $\blacksquare$ 17\zl .\zd\ \zl D\zd oes it make sense to speak of a mortal sin of \zl two unreadable symbols, should be: a\zd\ Christian (\zl perhaps: not\zd\ cathol\zl ic\zd \zl In the remainder, one recognizes ``ipso facto'' and ``etc'', which seem to be mixed with shorthand.\zd

$\blacksquare$ 17.1\zl .\zd\ Is any action of which one does\zl $\,$\zd not know that it is a mortal sin \zl text interupted\zd\

\zl Here an arrow originating in a $\blacksquare$ put before 23 $\circ$ and 24 $\circ$ points to the space between 17.1. and 18.\zd\

sp.\ 18\zl .\zd\ Is \zl unreadable word, perhaps: avarice\zd\ a mortal sin in any case\zl ?\zd\

sp.\ 19. Is every lie {intended for deceiving (maybe \zl unreadable text\zd ) a mortal sin] \zl ?\zd\
sp.\ 20\zl .\zd\ Is every action whose final aim is to damage anybody a mort\zl al\zd\ sin\zl ?\zd\

sp.\ 21\zl .\zd\ Is it a mort\zl al\zd\ sin to kill the enemy \zl in\zd\ a war waged by the \sout{\zl unreadable word\zd }\, secular power\zl ?\zd\

sp.\ 22. Is it a mortal sin \zl unreadable text\zd\

\zl see the remark after 17.1\zd\

Th 23 $\circ$ \sout{Are the mort\zl al\zd\ sins for a noncath.\zl olic\zd\ Christian} the same \sout{(even without \zl unreadable word\zd\ teaching)}

Th 24 $\circ$ Are they the same to \sout{someone} Christian who has \zl unreadable text\zd\ relig.\zl ious\zd\ teaching\zl ?\zd\

$\blacksquare$ 25\zl .\zd\ vide 74 Are all \zl unreadable word\zd\ made by God or also by other spirits\zl ?\zd\

\zl The following item, sp.\ 20, is set under a line at the bottom of the page, like a footnote.\zd\

sp.\ 20\zl unreadable symbol\zd\ to procure someth\zl ing\zd\ good for oneself by dam\zl ag\-ing\zd\ another.

$\mathbf{\llbracket new\, page \rrbracket}$ $\blacksquare$ 26. Is every \ul act\zl ual\zd\ \ud suff\zl ering\zd\ a punishment for a preceding (succeeding) \zl unreadable text\zd\ sin (of the parents)\zl ?\zd\ Animals?

sp.\ 27. Is it a mortal sin \sout{to ask a \zl unreadable word, perhaps shorthand\zd } \zl unreadable text\zd\ or to ask them

$\blacksquare$ 21\zl .\zd\ Is it in the power of anybody to make the world better by his acts or is it all the same\zl ?\zd\

sp.\ 22. Is the use of \zl unreadable word\zd\ means to make mon\zl e\zd y a mort\zl al\zd\ sin\zl ?\zd\

$\times\circ$ \zl preceded by symbols in shorthand, ``c/c.'' and ``bib''\zd\ $\blacksquare$ 23. Is it possible that any\zl body\zd\ who goes to heaven has a worse caract\zl er\zd\ than anyb\zl ody who\zd\ goes to hell (bec.\zl ause\zd\ \ud \zl unreadable word\zd\ \ud by \zl unreadable text\zd\ he was \zl unreadable word\zd\ from sins)\zl ?\zd\

$\blacksquare\times$ 24. Are some of the physi\zl cal\zd\ laws caused by \sout{evil} regular action of evil spirits\zl ?\zd\

$\mathbf{\llbracket new\, page \rrbracket}$ $\circ$ 25\zl .\zd\ Are the \sout{geo} fossils a work of the devil\zl ?\zd\

$\blacksquare\times$ 26\zl .\zd\ Do there exist exist any animals by natural reasons \zl (\zd without action of demons)\zl ?\zd\

$\blacksquare\circ$ 27. \zl written over ``vide 29''\zd\ Did the \zl unreadable word\zd\ really \zl unread\-able text\zd\ or was it a deception\zl ?\zd\

\zl The questions numbered $\blacksquare\circ$ 28., $\circ$ 29\zl .\zd , $-$30, $\circ$ \zl 31.\zd\ and $\circ$ 32\zl .\zd and Z.\ 33., which is on a new page, are mostly unreadable.\zd\ In $\circ$ 29\zl .\zd medicine seems to be mentioned together with spirits, in $-$30.\ the propositions of the Bible as dogmas, while in Z.\ 33.\ one finds ``Satisfact.\ for div.\ mortal sins.'' and ``Protest.'', in $-$34.\ one finds ``Def'', ``Dogma'' and ``Doctrine'' besides shorthand, and $-$35.\ is entirely in shorthand.\zd\

$\blacksquare$? 36\zl .\zd\ Is praying only caused by sin\zl ?\zd\

$\blacksquare$ 37\zl .\zd\ Are the saints in heaven at present have conscience and are praying\zl ?\zd\

$\blacksquare$ \zl preceding the next two numbers\zd\ 38. Is heaven where they are a \sout{spa} place in space\zl ?\zd\

39.\ Similarly (hell)

$\blacksquare$ 40.\ Has the body of J.\zl esus\zd\ Chr.\zl ist\zd\ moved to heaven\zl ?\zd\

\zl The text in the remaining remarks, numbered until 80, is mostly in shorthand, or too fragmentary to be understandable. We single out some readable words and fragments that seem important: in Lit 46.\ ``Martyrology'', in Th 57.\ ``\sout{Confessio} Excom.'' and ``Status'' and a word that seems to begin with ``amor'', in $\blacksquare$ 60\zl .\zd\ ``renasci de spiritu'', in Lit 64.\ ``Synchron.'', in $\blacksquare$ 69.\ ``omnes qui filii diaboli vocantur'', in ? $\blacksquare$ 70\zl .\zd\ ``filii diaboli'' and ``diabolus'', in $\blacksquare$ 72.\ ``Christus'', in 74.\ ``Lex iis qui sub lege sunt loquitur'', in $\blacksquare$ 77.\ ``corpus Christi'', in $-$78.\ ``theor\zl .\zd\ Physik'', and here are some complete numbered remarks:

$\blacksquare$ 75.\ Trinit. dogma\zl :\zd\ una natura\zl ,\zd\ tres persones

$\blacksquare$ 76.\ Homousi\zl os, homoousios:\zd\ una persona\zl ,\zd\ duae naturae

80.\ Ja[mes]cob 1,5, \zl 1\zd ,8 ?

New pages start within 44.\ and after 58. and 72.\zd\

\zl At the bottom of p.\ \textbf{136}.\ of Notebook VI, for subsequent pages, \textbf{137}.\ and later, the reader was directed to Notebook VII, and at the top of p.\ \textbf{137}.\ of Notebook VI it is written: ``Notebook VII''. It seems one should assume that pp.\ \textbf{137}.-\textbf{140}. of Notebook VI are to be superseded by pages which follow here, starting with p.\ \textbf{137}.\zd\

\vspace{2ex}

\noindent $\mathbf{\llbracket 137. \rrbracket}$
\begin{itemize}
\item[2.] The notion of a ,,\zl ``\zd well\zl -\zd defined \ul monadic \ud predicate''.
\end{itemize}
That is \sout{\zl unreadable word\zd }\, a monadic predicate $\varphi$ such that for any obj\zl ect\zd\ $x$ whatsoever it is uniquely det.\zl ermined\zd\ by the def.\zl inition\zd\ of $\varphi$ whether or not $\varphi$ belongs to $x$, so that \ul for any arb.\zl itrary\zd\ obj\zl ect\zd\ $x$ \ud $\varphi(x)$ is \sout{always} a meaningful prop.\zl osition\zd\ which is either true or false\zl .\zd\ Since I need no other kind of pred\zl icate\zd\ in the subsequ.\zl ent\zd\ considerations but only well\zl -\zd defined monadic pred\zl icates\zd , I shall use the term ,,\zl ``\zd \underline{pred.\zl icate\zd }'' in the sense of monadic well\zl -\zd def.\zl ined\zd\ pred\zl icate\zd .
\begin{itemize}
\item[\ul 3.] The concept which is expressed by \ul the word \ud ,,\zl ``\zd is'' or ,,\zl ``\zd belongs'' in ord\zl inary\zd\ langu\zl age\zd\ and which we expressed by $\varphi(x)$\zl ,\zd\ which means the pred\zl icate\zd\ $\varphi$ belongs to $x$\zl .\zd\ \ud
\end{itemize}

Now for these notions (of obj.\zl ect\zd\ and pred.\zl icate\zd ) we have the foll.\zl ow\-ing\zd\ apparently evident prop.\zl ositions:\zd\

\vspace{2ex}

\noindent $\mathbf{\llbracket 138. \rrbracket}$
\begin{itemize}
\item[1.] For any obj.\zl ect\zd\ $x$ it is \sout{well} uniquely det\zl ermined\zd\ whether or not it is a \sout{welldef.} pred.\zl icate;\zd\ \ul in other word\zl s\zd\ \ul well\zl -\zd def\zl ined\zd\ \ud predicate is itself a well\zl -\zd defined predicate\zl .\zd\ \ud
\item[2.] If $y$ is a pred.\zl icate\zd\ and $x$ an obj\zl ect\zd\ then it is well\zl -\zd defined whether the pred.\zl icate\zd\ $y$ belongs to $x$. \ul This is an immed.\zl iate\zd\ consequence of the def.\zl inition\zd\ of a well\zl -\zd defined pred\zl icate\zd . \ud
\end{itemize}

Let us denote \ul for any two obj\zl ects\zd\ $y$\zl , $x$\zd\ \ud by $y(x)$ the prop.\zl osition $y$ is a pred.\zl icate\zd\ and belongs to $x$\zl .\zd\ So for any two obj\zl ects\zd\ \ul $y,x$ \ud $y(x)$ will be a meaningful prop.\zl osition\zd\ \ul of \ud which it is uniquely determ\zl ined\zd\ whether it is true or false\zl ,\zd\ namely if $y$ is no pred\zl icate\zd\ it is false \ul whatever $x$ may be\zl ,\zd\ \ud if it is a pred\zl icate\zd\ then it is true or false according as the pred.\zl icate\zd\ $y$ bel.\zl ongs\zd\ to $x$ or does not belong to $x$\zl ,\zd\ which is uniquely det\zl ermined\zd .

\vspace{2ex}

\noindent $\mathbf{\llbracket 139. \rrbracket}$
\begin{itemize}
\item[3.] If we have a combination of symbols or words \sout{$A(x)$} contain\zl ing\zd\ the letter $x$ (denote it by $A(x)$) and if this comb.\zl ination\zd\ is such that it becomes a mean.\zl ingful\zd\ prop.\zl osition\zd\ whatever \sout{you} obj\zl ect\zd\ you subst.\zl itute\zd\ for $x$ then $A(x)$ defines a cert\zl ain\zd\ \ul well\zl -\zd def\zl ined\zd\ \ud predicate $\varphi$ which belongs to an obj\zl ect\zd\ $x$ if and only if $\varphi(x)$ \zl $A(x)$\zd\ is true.
\end{itemize}
(I repeat the hypothesis of this statement: It \sout{means} is as follows\zl ,\zd\ that if you subst\zl itute\zd\ for $x$ \ul the name of \ud an arb.\zl itrary\zd\ obj.\zl ect\zd\ then the resulting expr.\zl ession\zd\ is always a meaningf.\zl ul\zd\ prop.\zl osition\zd\ of which it is uniquely det.\zl ermined\zd\ whether it is true or false.) \ul Now this statement too could be consid.\zl ered\zd\ as a consequence of the def.\zl inition\zd\ of a well\zl -\zd def.\zl ined\zd\ pred\zl icate\zd . \ud
\begin{itemize}
\item[4.] Any pred\zl icate\zd\ is an obj\zl ect\zd . That $\mathbf{\llbracket 140. \rrbracket}$ follows bec.\zl ause\zd\ we took the term obj\zl ect\zd\ in the most general sense according to which anything one can think of is an object.
\end{itemize}

I think these 4\zl four, written over 3\zd\ prop.\zl ositions\zd\ are all evident to natural thinking. But nevertheless they lead to contradictions\zl ,\zd\ namely in the following way\zl .\zd\ Consider the expr\zl ession\zd\ $\sim x(x)$ that is an expr.\zl ession\zd\ involving\zd\ the var\zl iable\zd\ $x$ and such that for any obj\zl ect\zd\ \sout{\zl unreadable symbol\zd }\, substituted for this var.\zl iable\zd\ \ul $x$ \ud you do \ul obtain \ud a \ul mean.\zl ing\zd ful\zl \sout{l}\zd\ propos.\zl ition\zd\ of which it is uniquely det\zl ermined\zd\ whether it is true or \zl missing from the manuscript: false.\zd\ $\mathbf{\llbracket 141. \rrbracket}$ \zl N\zd amely if $x$ is not a pred.\zl i\-cate\zd\ this bec.\zl omes\zd\ false by the above definition of $y(x)$\zl ;\zd\ if $x$ is a pred\zl i\-cate\zd\ then \ul by 1 \ud for any obj\zl ect\zd\ $y$ it is uniquely det.\zl ermined\zd\ whether $x$ bel\zl ongs to\zd\ $y$\zl ,\zd\ hence also for $x$ it is uniquely det\zl ermined\zd\ bec\zl ause\zd\ $x$ is a pred\zl icate,\zd\ hence an object (by 4)\zl .\zd\ \sout{\zl unreadable word\zd }\, $\sim x(x)$ means $x$ is a pred.\zl icate\zd\ not belonging to itself. It is easy to name pred\zl icates\zd\ which do belong to themselves\zl ,\zd\ e.g.\ the pred\zl icate\zd\ ,,\zl ``\zd predicate''\zl ;\zd\ we have \ul the concept \ud ,,\zl ``\zd predicate'' is a predicate. Most of the pred.\zl icates\zd\ of course do\zl $\,$\zd not belong to thems\zl elves.\zd\ \zl S\zd ay e\zl .\zd g.\ \zl t\zd he predicate man is not a man\zl ,\zd\ $\mathbf{\llbracket 142. \rrbracket}$ so it does\zl $\,$\zd not belong to itself\zl .\zd\ But e\zl .\zd g.\ the pred\zl icate\zd\ not man \zl hyphen between these two words deleted, since it is omitted in the text later\zd\ does belong to itself since the pred\zl icate\zd\ not man is certainly not a man\zl ,\zd\ so it is a not man\zl ,\zd\ i\zl .\zd e.\ belongs to itself\zl .\zd\

Now since $\sim x(x)$ is either true or false for any obj\zl ect\zd\ $x$ it defines a cert\zl ain\zd\ pred\zl icate\zd\ by 3. Call this \ul well\zl -\zd def\zl ined\zd\ \ud pred.\zl icate\zd\ $\Phi$\zl ,\zd\ so that $\Phi(x)\equiv\;\sim x(x)$\zl .\zd\ For $\Phi$ even a term in ord.\zl inary\zd\ lang.\zl uage\zd\ was introduced\zl ,\zd\ namely the word ,,\zl ``\zd impredicable''\zl ,\zd\ and \ul for \ud the neg\zl ation\zd\ of it \ul the word \ud ,,\zl ``\zd predicable''\zl ;\zd\ so \ul an obj\zl ect\zd\ is called \ud predicable if it $\mathbf{\llbracket 143. \rrbracket}$ is a pred.\zl icate\zd\ belonging to itself and impredicable in the opposite case\zl ,\zd\ \ul i\zl .\zd e.\ if it is either not a pred\zl icate\zd\ or is a pred\zl icate\zd\ and does\zl $\,$\zd not belong to itself. \ud \zl S\zd o predicate is predicable\zl ,\zd\ not man is pred\zl icable,\zd\ man is impred\zl icable,\zd\ Socr\zl ates\zd\ is impred\zl icable\zd . \zl A line at the end of this paragraph separates it from the text below it.\zd\

And now we ask is \sout{predicable} the pred\zl icate\zd\ ,,\zl ``\zd impred\zl icable\zd '' \ul predi\-cable or \ud impredicable\zl .\zd\ Now we know this equiv.\zl alence\zd\ holds for any obj\zl ect\zd\ $x$ (it is the def\zl inition\zd\ of impred\zl icable\zd ).\zl ;\zd\ $\Phi$ is a pred\zl icate,\zd\ hence an obj.\zl ect,\zd\ hence this equiv.\zl alence\zd\ holds \ul for $\Phi$\zl ,\zd\ \ud i\zl .\zd e.\zd\ $\Phi(\Phi)\equiv\;\sim \Phi(\Phi)$. \sout {And} \ul What does \ud $\Phi(\Phi)$ say\zl ?\zd\ Since $\Phi$ means impred.\zl icable\zd\ it says \ul \sout{the pre} \ud impred\zl icable\zd\ is impredicable. \sout{and} So we see that this prop.\zl osition\zd\ is equivalent \sout{\zl unreadable symbol\zd }\, with its \ul own \ud negation\zl .\zd\

$\mathbf{\llbracket 144. \rrbracket}$ But from that it follows that it must be both true and false\zl ,\zd\ bec\zl ause\zd\ we can conclude from this equiv\zl alence\zd :
\begin{tabbing}
\hspace{1.7em}\=$\Phi(\Phi)\supset\; \sim \Phi(\Phi)$\\[.5ex]
\>$\sim \Phi(\Phi) \supset \Phi(\Phi)$
\end{tabbing}
By the first impl.\zl ication,\zd\ $\Phi(\Phi)$ cannot be true\zl ,\zd\ bec\zl ause\zd\ the ass\zl ump\-tion\zd\ that it is true leads to the concl.\zl usion\zd\ that it is false\zl ,\zd\ i\zl .\zd e.\ \ul it leads \ud to a contradiction\zl ;\zd\ but \sout{\zl unreadable symbol\zd }\, $\Phi(\Phi)$ cannot be false either because \sout{the ass} by the sec.\zl ond\zd\ impl.\zl ication\zd\ the ass.\zl umption\zd\ that it is false leads to the concl.\zl usion\zd\ that it is true.\zl ,\zd\ i\zl .\zd e.\ \ul again \ud to a contrad\zl iction\zd . So this $\Phi(\Phi)$ would be a prop.\zl osition\zd\ which is neither true nor false\zl ,\zd\ hence it would be both true and false $\mathbf{\llbracket 145. \rrbracket}$ bec\zl ause\zd\ that it is not true implies that it is false and that it is not false implies that it is true. So we apparently have discovered a prop\zl osition\zd\ which is both true and false\zl ,\zd\ which is impossible by the law of contradiction\zl .\zd\

\zl The text in the following paragraph is inserted in the manuscript on the right of p.\ \textbf{145}., which is numbered \textbf{145}.\textbf{1}., and at the top of the not numbered page on the right of p.\ \textbf{146}.\zd\ \ul The same argument can be given without log.\zl ical\zd\ symb\zl ols\zd\ in the following form\zl .\zd\ The quest\zl ion\zd\ is: \ul Is the pred.\zl icate \ud ,,\zl ``\zd impredicable'' pred.\zl icable\zd\ or impred\zl icable\zd . 1.\ If \sout{it} \ul impred\zl icable\zd\ \ud were pred.\zl icable\zd\ that would mean that it belongs to itself\zl ,\zd\ i\zl .\zd e.\ then impred\zl icable\zd\ is impred\zl icable\zd . So from the ass\zl umption\zd\ that \ul impred.\zl icable\zd\ \ud is pred.\zl i\-cable\zd\ we derived that it is impred\zl icable;\zd\ so it is not \sout{im} predic\zl able\zd . 2.\zl \sout{)}\zd\ On the other hand assume impred\zl icable\zd\ is impred\zl icable;\zd\ then it belongs to itself\zl ,\zd\ hence \ul \sout{\zl unreadable word\zd }\, \ud is predicable. So from the ass.\zl umption\zd\ that it is impred\zl icable\zd\ we derived that it is pred\zl icable\zd . So it is cert\zl ainly\zd\ not impred\zl icable\zd . So it is neither pred\zl icable\zd\ nor impred\zl icable\zd . But then it must be both pred.\zl icable\zd\ and impred.\zl icable\zd\ because since it is not pred.\zl icable\zd\ it is impr.\zl edicable\zd\ and since it is not impred\zl icable\zd\ it is pred\zl icable\zd . So again we have a prop.\zl osition\zd\ which is both true and false\zl .\zd\ \ud

Now what are we to \zl unreadable word, should be: do\zd\ about this situation? One may first try to say\zl :\zd\ Well\zl ,\zd\ the law of contradiction is an error. There do exist such strange things as prop.\zl ositions\zd\ which are both true and false. But this \ul way out of the diff\zl iculty\zd\ \ud is \ul evidently \ud not possible $\mathbf{\llbracket 146. \rrbracket}$ because that would imply that every prop.\zl osition\zd\ \ul whatsoever \ud is both true and false\zl .\zd\ We had \sout{the form. of} in the calc.\zl ulus\zd\ of prop\zl ositions\zd\ the form\zl ula\zd\ $p\; .\sim p\supset q$ \ul for any $p,q$\zl ,\zd\ \ud hence also $p\; .\sim p\supset \; \sim q$ where $p$ and $q$ are arb\zl itrary\zd\ prop\zl ositions\zd . So if we have one prop\zl osition\zd\ $p$ which is both true and false then any prop\zl osition\zd\ $q$ has the undesirable prop\zl erty\zd\ of being both true and false\zl ,\zd\ which would make any thinking completely meaningless. So we have to conclude that we arrived at this contradictory concl.\zl usion\zd\
\[
\Phi(\Phi) \quad \rm{and} \quad \sim \Phi(\Phi)
\]
$\mathbf{\llbracket 147. \rrbracket}$ by some error or fallacy\zl ,\zd\ and the question is what does this error consist in [i\zl .\zd e.\ which one of our evident prop\zl ositions\zd\ is wrong]\zl .\zd\

\zl new paragraph\zd\ The nearest at hand conjecture about this error is that there is some circular fallacy hidden in this argument.\zl ,\zd\ because we are speaking of pred.\zl icates\zd\ belonging to themselves or not belonging to themselves. One may say that it is meaningless \sout{from the beginning} to apply a predicate to itself\zl .\zd\ \zl ' deleted\zd\ I don't think that this is the correct solution. For the following reasons\zl :\zd\
\begin{itemize}
\item[1.] It is \ul not possible to \ud except for any pred.\zl icate\zd\ $P$ $\mathbf{\llbracket 148. \rrbracket}$ just this pred\zl icate\zd\ $P$ itself from the things to which it can be applied
\end{itemize}
i\zl .\zd e\zl .\zd\ \sout{\zl unreadable word\zd }\, we \ul cannot \ud modify the assumption 1.\ by \ul saying \ud the \zl written over another unreadable word\zd\ prop.\zl property, or perhaps: proposi\-tion\zd\ $\varphi(x)$ is well\zl -\zd def.\zl ined\zd\ for any $x$ except $\varphi$ itself because if you define \ul e.g.\ \ud a pred\zl icate\zd\ $\mu$ \sout{say} by two pred\zl icates\zd\ $\varphi,\psi$ by $\mu(x)\equiv \varphi(x)\; .\; \psi(x)$ then we would have already three \zl written over another unreadable word\zd\ pred.\zl icates\zd\ $\mu$, $\varphi$ and $\psi$ to which $\mu$ cannot be applied\zl :\zd\
\begin{tabbing}
\hspace{1.7em} $\mu(\varphi)\;{\equiv\atop\rm{Df}}\;\varphi(\varphi)\; .\; \psi(\varphi)$\quad where this makes no sense\zl .\zd\
\end{tabbing}
$\mathbf{\llbracket 149. \rrbracket}$ So it is certainly not sufficient to exclude just self\zl -\zd reflexivity \ul of a pred.\zl icate\zd\ \ud \ul because that entails automatically that we have to exclude also other thing\zl s\zd\ and it is very difficult and leads to \ul very \ud undesirable results if one tries to formulate what is to be excluded \ul \sout{\zl unreadable text\zd }\, \ud on the basis of this idea to avoid self\zl -\zd reflexivities. That was done by Russel\zl l\zd\ in his so called ramified theory of types which since has been abandoned by practically all logicians. \ud On the other hand \ul it is not even justified to exclude self\zl -\zd reflexivities of every formula \ud \ul bec.\zl ause\zd\ \ud self\zl -\zd reflexivity does\zl $\,$\zd not always lead to contradiction but is perfectly legitimate in many cases\zl .\zd\ If \ul e.g\zl .\zd\ \ud I say \zl \sout{e.g.}:\zd\ ,,\zl ``\zd Any sent\zl ence\zd\ of the English language contains a verb\zl ''\zd\ then it is perfectly alright to apply this proposition to itself and to conclude from it that also this prop.\zl osition\zd\ under consideration contains a verb.

\zl new paragraph\zd\ \sout{The} \ul Therefore the \ud real fallacy seems to ly\zl lie\zd\ $\mathbf{\llbracket 150. \rrbracket}$ in something else \ul tha\zl \sout{t}\zd n the self\zl -\zd reflexivity\zl ,\zd\ \ud namely in these \sout{\zl unreadable symbol\zd }\, notions of object and predicate in the most general sense \ul embracing obj\zl ects\zd\ of all logical types \ud . The Russell paradox seems to show that there does not exist such \zl a\zd\ concept of everything \sout{because}\zl .\zd\ A\zl written over a\zd s we saw the logical objects form a \sout{certain} hierarchy of types and however far you may proceed in the\zl ``e'' written over ``is''\zd\ construction of these types you will always be able to continue the process \sout{\zl unreadable symbol\zd }\, still farther and therefore it is illegitimate and makes no sense to speak of the totality of all obj\zl ects\zd .

$\mathbf{\llbracket 151. \rrbracket}$ One might think that one could obtain the totality of all obj.\zl ects\zd\ in the following way: take first the indiv.\zl iduals\zd\ and call them obj.\zl ects\zd\ of type 0\zl ,\zd\ then take the concepts of type 1\zl ,\zd\ then the conce\zl pts\zd\ of type 2\zl ,\zd\ 3 etc\zl .\zd\ for any natural nu\zl mber\zd . But it is by no means true that we obtain in this manner the totality of all concepts.\zl ,\zd\ \sout{But that isnt true} because \ul e.g\zl .\zd\ \ud the concept of the\zl ``e'' written over ``is''\zd\ totality of concepts thus obtained \ul for all int\zl egers\zd\ $n$ as types is itself a \ud concept not occurring in this totality\zl ,\zd\ i\zl .\zd e.\ it is a concept of a ty\zl ``y'' written over another letter\zd pe higher than $\mathbf{\llbracket 152. \rrbracket}$ any finite nu.\zl mber,\zd\ i\zl .\zd e.\
of an infinite type. It is denoted as \ul a concept of \ud type $\omega$. But even with this type \ul $\omega$ \ud we are by \ul no means \ud at an end, \sout{either} because we can \ul \sout{e.g.} \ud define \sout{concepts which are} e.g\zl .\zd\ relations between conc\zl epts\zd\ of
type $\omega$ and they would be of \ul a still higher \ud type $\omega + 1$\zl .\zd\ So we see there are \ul in a sense \ud much more than infinitely many log\zl ical\zd\ type\zl s\zd ; \sout{and} there are so many that it is not possible to form a concept of the totality of all of them\zl ,\zd\ because whichever concept we form we can define a concept of a higher type\zl ,\zd\ hence not falling under $\mathbf{\llbracket 153. \rrbracket}$ the given concept.

\zl new paragraph\zd\ So if we want to take account of this fundamental fact of logic \ul that there does\zl $\,$\zd not exist a concept of the totality of all objects whatsoever\zl ,\zd\ \ud we must drop the words ,,\zl ``\zd object''\zl , ``\zd predicate''\zl , ,,\zl ``\zd everything'' from our language and replace them by the words: object of a given type\zl ,\zd\ predicate of a given type\zl ,\zd\ everything which belongs to a given type. \ul In part.\zl icular,\zd\ prop\zl osition\zd\ 4 has now \zl to\zd\ be formul\zl ated\zd\ like this. \zl If\zd\ $A(x)$ is an expr.\zl ession\zd\ which becomes a
meaningf.\zl ul\zd\ prop.\zl osition\zd\ for any obj\zl ect\zd\ $x$ of a given type $\alpha$ then it defines a concept of type $\alpha +1$\zl .\zd\ \sout{Now} We cannot even
formulate the prop.\zl osition\zd\ in its previous form.\zl ,\zd\ because we don't have such words as obj\zl ect\zd , pred\zl icate\zd\ etc\zl .\zd\ in our
lang\zl uage\zd . \ud Then the Russell paradox disappears immediately because
we can form the concept $\Phi$ defined by $\Phi(x)\equiv\;\sim x(x)$ only for $x$'s of a given type $\alpha$\zl ,\zd\ i\zl .\zd e.\ $\mathbf{\llbracket 154. \rrbracket}$ we can define a concept $\Phi$ such that this equivalence holds for every $x$ of
type $\alpha$\zl .\zd\ (We cannot even formulate that it holds for every obj.\zl ect\zd\ because we have dropped these words from our langu\zl age\zd ).\zl .)\zd\ But then $\Phi$ will be \zl a\zd\ concept of next higher type because it is a property of objects of type $\alpha$. Therefore we
cannot substitute $\Phi$ here for $x$ because this equiv\zl alence\zd\ holds only for obj.\zl ects\zd\ of type $\alpha$.

\zl new paragraph\zd\ So this seems to me to be the (satisfactory) true solution of the $\mathbf{\llbracket 155. \rrbracket}$ Russell paradox\zl \sout{e}\zd \sout{s}.
I only wish to mention that the hierarchy of types as I sketched it here is
considerably more general than it was when it was first presented by it's\zl its\zd\ inventor B.\
Russell. Russell's theory of types was given \zl in\zd\ two different forms\zl ,\zd\ the so called
simplified and the ramified theory of types\zl ,\zd\ both of which are much more
restrictive then the one I explained here\zl ; e.\zd g\zl .\zd\ in both of them it would be
imp.\zl ossible\zd\ to form concepts of type $\omega$, \ul also the statement $x(x)$ would always be
meaningless\zl .\zd\ \ud Russell\zl '\zd s theory\zl ``he'' written over an unreadable word\zd\ of $\mathbf{\llbracket 156. \rrbracket}$ types is more based on the first idea of
s\zl writen over: ex\zd olving the paradoxes (namely to exclude self\zl -\zd reflexivities) and the tot.\zl ality\zd\ of all obj\zl ects\zd\ is only excl.\zl uded\zd\ because it would be self\zl -\zd reflexive (since it would itself be
an object\zl )\zd . However the develop\zl ment\zd\ of ax\zl ioms\zd\ of set theory has shown that
Russell\zl '\zd s syst\zl em\zd\ is too restrictive\zl ,\zd\ i\zl .\zd e.\ it excludes many arguments \ul which (as far as one can see) do\zl $\,$\zd not lead to contradictions and which are necessary for building up
abstract set\zl $\,$t\zd heory\zl .\zd\ \ud

There are other logical paradoxes which are solved by the theory of types\zl ,\zd\ i\zl .\zd e.\ by
excluding the terms obj\zl ect\zd , every etc\zl .\zd\ But there are others in which the fallacy is of
an \ul entirely \ud different nature. They are the so called epistemological paradoxes. $\mathbf{\llbracket 157. \rrbracket}$ The oldest of them is the Epimenides\zl .\zd\ In the form it is \sout{\zl unreadable symbols\zd }\, usually
presented, it is no paradox. But if a man says ,,\zl ``\zd I am lying now'' \ul and says nothing else\zl ,\zd\ \ul or if he says: The prop.\zl osition\zd\ which I am \sout{jus} pronouncing right now is false\zl ,\zd\ \ud then th\zl written over something unreadable\zd is
statement can be proved to be both true and false, because this prop\zl osition\zd\ $p$ says
that $p$ is false\zl ;\zd\ so we have $p\equiv (p$ is false)\zl ,\zd\ $p\equiv\;\sim p$\zl ,\zd\ from which it follows that $p$ is both
true and false as we saw before. The same para\zl dox\zd\ can be brought to a much more
conclusive form as follows:

\vspace{2ex}

\zl Here, at the end of p.\ \textbf{157}., the text in the manuscript is interrupted, and subsequent pages are not numbered until p.\ \textbf{1}. below. In between are four pages of jottings given here, presumably for exercises.\zd

\vspace{2ex}

$\mathbf{\llbracket new\, page\; i \rrbracket}$ Ableitung d.\zl der\zd\ paradoxen\zl ,,ox'' in this word written over something else\zd\ Aussagen \" uber Impl.\ \zl Implikation\zd\ aus den unten angeg.\zl an\-gegeben\zd\ 5 Axiomen \zl German: Derivation of paradoxical propositions about implication from the five axioms given below:\zd\

\zl The next three lines, before 1., are crossed out:\zd\
\begin{tabbing}
\hspace{1.7em}$p\supset p$\hspace{3.4em}\=$r\supset(p\supset p\lfloor .\rfloor r)$\\*[.5ex]
\> \underline{$p\;.\; r\supset r$}\\*[.5ex]
\> $p\supset r$
\end{tabbing}
\begin{tabbing}
\hspace{1.7em}\=1.) \=$p\supset q$ \hspace{3.4em}\=2.) \=$p\supset q$\\*[.5ex]
\>\>\underline{$q\supset r$}\>\>\underline{$r$\hspace{3em}}\\*[.5ex]
\>\>$p\supset r$\>\>$p\supset q\;.\; r$\\[1ex]
\>3.)\>$q\;.\; r\supset r$\>4.)\>$p\supset p$
\end{tabbing}
Ableitung \zl derivation\zd
$
\begin{cases}
\begin{tabular}{ l|l }
$r$ & \\[-1ex]
$p\supset p\; .\; r$ & \\[-.5ex]
\underline{$p\; .\; r \supset r$}& \\
$p\supset r$ & \rm{\zl unreadable sign\zd\ 5.) tollendo tollens}\\[-.5ex]
$\sim r\supset \;\sim p$ & \hspace{8em} $p\supset r\: . \sim r : \supset\; \sim p$
\end{tabular}
\end{cases}
$

\begin{tabbing}
$\mathbf{\llbracket new\, page\; ii \rrbracket}$\\*[1ex]
\hspace{1.7em}\=$\mu (x)$ \zl $=$ or $\equiv$\zd\ $\varphi (x) \; .\; \psi (x)$\\*[.5ex]
\>$\varphi(\mu)$\\[.5ex]
\>strike out, drop\\[.5ex]
\>something else but (than)\\[.5ex]
\>falling under a concept
\end{tabbing}

\noindent \zl The following pages \textbf{new page iii-iv} and pp.\ \textbf{1}.-\textbf{7}.\ following them until the end of the scanned manuscript, which makes nine pages, are on loose, torn out, leafs,
with holes for a spiral, but not bound with the spiral to the rest of the notebook, as the other pages in this Notebook VII are. In all of the notebooks the only other loose leafs are to be found at the end of Notebooks III and towards the end of Notebook~V.\zd\

\begin{tabbing}
$\mathbf{\llbracket new\, page\; iii \rrbracket}$\\*[1ex]
\hspace{1.7em}\=\textcircled{1.} $p\rightarrow p$ \hspace{3em} \= \sout{\textcircled{1}. $q\rightarrow (\Delta\rightarrow q)$}\\*[.5ex]
\>\textcircled{2.}$\; p,q\rightarrow p$\>\sout{\textcircled{2}. $\sim q\rightarrow (\lfloor {\rm crossed\; out\; symbol}\rfloor q,\Delta\rightarrow p)$}\\[.5ex]
\>\> \sout{\textcircled{3} $\sim\sim q\rightarrow q$}\\[.5ex]
\>\> \sout{\textcircled{4.} $\Delta, q\rightarrow \;\sim(\Delta\rightarrow q)$}\\[.5ex]
\> \sout{\textcircled{5.} $\Delta\rightarrow p$}\\[.5ex]
\> \hspace{1em} \underline{\sout{$\Delta\rightarrow q$}\hspace{3em} \sout{$p, q\rightarrow r$}}\\[.5ex]
\hspace{6.5em}\sout{$\Delta\rightarrow r$}\\[1ex]
\> \sout{\textcircled{6.} $\Delta,p\rightarrow q$}\hspace{3em}\zl unreadable formula\zd\ \\[.5ex]
\> \hspace{1em} \underline{\sout{$\Delta,\sim p\rightarrow q$}\hspace{1.7em} \sout{\zl unreadable formula\zd }}\\[.5ex]
\hspace{6.5em} \sout{$\Delta\rightarrow q$}\\[1ex]
\> \textcircled{5} ${\mathfrak A}\rightarrow p_1$\\[.5ex]
\>\hspace{1em} ${\mathfrak A}\rightarrow p_n$\hspace{5em}\= ${\mathfrak A}\rightarrow (p\rightarrow r)$ \\[.5ex]
\>\hspace{1em} \underline{$p_1\ldots p_n\rightarrow q$}\> \underline{${\mathfrak A}\rightarrow p$} \\[.5ex]
\>\hspace{1em} ${\mathfrak A}\rightarrow q$\>${\mathfrak A}\rightarrow r$\\[1ex]
\>\textcircled{4.} Ind.\zl uktive\zd\ Bew.\zl eis\zd \zl German: Inductive proof\zd\ \\[.5ex]
\>\textcircled{3} Export.\ \zl gen\zd\ a Import.\ \hspace{2em}\= $(p,p\rightarrow r)\rightarrow r$\\[.5ex]
\>\> $\lfloor (\rfloor p\rightarrow r)\rightarrow p\rightarrow r$
\end{tabbing}

\begin{tabbing}
$\mathbf{\llbracket new\, page\; iv \rrbracket}$\\*[1ex]
\zl 1.\zd\ \hspace{1em}\= $R+S$\zl ,\zd\ $R\cdot S$\zl ,\zd\ $R\subset S$\zl ,\zd\ $-R$\zl ,\zd\ $R-S$\\*[.5ex]
\>\sout{$\dot{\rm V}$\zl ,\zd\ $\dot\Lambda$\zl ,\zd\ $R`x$\zl ,\zd\ $E!R`x$}\\[.5ex]
\>[$\overrightarrow{R}$\zl ,\zd\ $\overleftarrow{R}$\zl ,\zd\ $R``\beta$\zl ,\zd\ $R_{\varepsilon}`\beta$]\\[.5ex]
\>$\breve{R}$\zl ,\zd\ $D`R$\zl , $C`R$, $C`R$,\zd\ $R|S$
\end{tabbing}
\zl inside two incomplete boxes and crossed out: sym\zl metry\zd , as, 1, [, I\zd\
trans\zl itivity,\zd\ one many\zl ,\zd\ father, $\bot\atop {x}$,
\begin{tabbing}
\hspace{1.7em}\= $xM(y,z)$\zl ,\zd\ $M`(y\lfloor ,\rfloor z)$\zl ,\zd\ $yMz$\\[.5ex]
\> $i`x$\zl ,\zd\ $\{x\}$\zl ,\zd\ $0,1,2,\ldots$\zl ,\zd\ $1\rightarrow 1$
\end{tabbing}
Abstractions \zl $\;$\zd\ prinz.\zl ip, perhaps German: principle,\zd\ aeq.\zl perhaps: equal, or something of the same root,\zd\ Ind\zl uction\zd .

\vspace{1ex}

$\mathbf{\llbracket 1. \rrbracket}$ \zl Here the numbering of pages in this notebook starts anew.\zd\ All four rules are purely formal\zl ,\zd\ i\zl .\zd e.\ for applying them it is \zl apostrophe deleted\zd\ not necessary to know the meaning of the expressions. Examples of derivations from the axioms. Since all axioms and rules of the calculus of propositions are also axioms and rules of the calculus of functions we are justified \sout{to} in as\zl sum\zd ing all formulas and rules formerly derived \zl the order of the last two words corrected in the manuscript\zd\ \sout{for} in the calculus of propositions.

\vspace{1ex}

1. Example $\varphi(y)\supset (\exists x)\varphi(x)$
\begin{tabbing}
Derivation:\\[.5ex]
\hspace{1.7em}\= (1)\hspace{1em}\= $(x)[\sim\varphi(x)]\supset\;\sim\varphi(y)$\quad obtained by substituting $\sim\varphi(x)$ for $\varphi(x)$\\
\` in Ax5.\zl Ax.\ 5\zd\ \\
\>(2)\>$\varphi(y)\supset\;\sim(x)[\sim\varphi(x)]$\quad by rule of transposition applied to (1)\\[.5ex]
\>(3)\>$\varphi(y)\supset (\exists x)\varphi(x)$\quad by rule of defined symbol from (2)
\end{tabbing}
$\mathbf{\llbracket 2. \rrbracket}$ 2. Example $(x)[\varphi(x)\supset\psi(x)]\supset[(x)\varphi(x)\supset(x)\psi(x)]$
\begin{tabbing}
\hspace{1.7em}\= (1)\hspace{1em}\= $(x)[\varphi(x)\supset\psi(x)]\supset[\varphi(y)\supset\psi(y)]$\quad by substituting $\varphi(x)\supset\psi(x)$\\
\` for $\varphi(x)$ in Ax.\ 5\\
\>(2)\>$(x)\varphi(x)\supset\varphi(y)$\quad Ax\zl . \zd 5\\[.5ex]
\>(3)\>$(x)[\varphi(x)\supset\psi(x)]\: .\: (x)\varphi(x)\supset[\varphi(y)\supset\psi(y)]\: .\: \varphi(y)$\quad by rule of\\
\` multiplication of implications applied to (1) and (2)\\[.5ex]
\>(4)\> $[\varphi(y)\supset\psi(y)]\: .\: \varphi(y)\supset\psi(y)$\quad by substituting $\varphi(y)$ for $p$ and $\psi(y)$\\
\` for $q$ in the demonstrable formula $(p\supset q)\: .\: p\supset q$\\[.5ex]
$\mathbf{\llbracket 3. \rrbracket}$\>(5)\>$(x)[\varphi(x)\supset\psi(x)]\: .\: (x)\varphi(x)\supset\psi(y)$\quad by rule of syllogism applied to\\
\` (3) \zl and\zd\ (4)\\[.5ex]
\>(6)\>$(x)[\varphi(x)\supset\psi(x)]\: .\: (x)\varphi(x)\supset(y)\psi(y)$\quad by rule of quantifier from (5)\\[.5ex]
\>(7)\>$(x)[\varphi(x)\supset\psi(x)]\supset[(x)\varphi(x)\supset (y)\psi(y)]$\quad by rule of exportation\\
\` from \zl (\zd 6\zl )\zd\ \\[.5ex]
\>(8)\>$(x)[\varphi(x)\supset\psi(x)]\supset[(x)\varphi(x)\supset (x)\psi(x)]$\quad by rule of substitution for\\
\` individual variables
\end{tabbing}

Predicates which belong to no object are called vacuous (e.g.\ president of U.S.A. born in South Bend). $S$a$P$ and $S$e$P$ are both true if $S$ is vacuous whatever $P$ may be. $\mathbf{\llbracket 4. \rrbracket}$ All tautologies are true also for vacuous predicates but some of the Aristotelian inferences are not\zl ,\zd\ e.g.\
\begin{tabbing}
\hspace{1.7em}\= $S$a$P\supset S$i$P$\hspace{3em}\= (false if $S$ is vacuous)\\[.5ex]
\> $S$a$P\supset \;\sim(S$e$P)$\> (false $\, ''$ $\, ''$ $\,''$ \hspace{1.55em}$''$\hspace{1.55em}),
\end{tabbing}
the mood Darapti $M$a$P\: .\: M$a$S\supset S$i$P$ is false if $M$ is vacuous and if $S,P$ are any two predicates such that $\sim(S$i$P)$.

The totality of all objects to which a monadic predicate $P$ belongs is called the extension of $P$ and denoted by $\hat{x}[P(x)]$, so that the characteristic $\mathbf{\llbracket 5. \rrbracket}$ property of the symbol $\hat{x}$ is:
\[
\hat{x}\varphi(x)=\hat{x}\psi(x)\equiv (x)[\varphi(x)\equiv\psi(x)]
\]

Extensions of monadic predicates are called classes (denoted by $\alpha,\beta,\gamma\ldots$) \zl .\zd\ That $y$ belongs to the class $\alpha$ is expressed by $y\varepsilon\alpha$ so that $y\varepsilon\hat{x}\varphi(x)\equiv \varphi(y)$\zl .\zd\ $\hat{x}$ is \sout{also} applied to arbitrary propositional functions $\Phi(x)$\zl ,\zd\ i\zl .\zd e.\ $\hat{x}\Phi(x)$ means the class of objects satisfying $\Phi(x)$\zl ,\zd\ e.g\zl .\zd\ $\hat{x}[I(x)\: .\: x>7]=$ class of integers greater \zl than\zd\ seven\zl .\zd\

Also for dyadic pred\zl icates\zd\ \ul $Q(xy)$ \ud extensions \ul denoted by $\hat{x}\hat{y}[Q(xy)]$ \ud are introduced\zl ,\zd\ which satisfy the equivalence
\[
\hat{x}\hat{y}[\psi(xy)]=\hat{x}\hat{y}[\chi(xy)]\equiv (x\lfloor ,\rfloor y)[\psi(xy)\equiv\chi(xy)]
\]

$\mathbf{\llbracket 6. \rrbracket}$ It is usual to call these extensions (not the dyadic predicates themselves) relations. If $\Phi(xy)$ is a propositional function with two variables $\hat{x}\hat{y}\Phi(xy)$ denotes the relation which is defined by $\Phi(xy)$\zl .\zd\ If $R$ is a relation $xRy$ means that $x$ bears the relation $R$ to $y$ so that
\[
u\{\hat{x}\hat{y}[\varphi(xy)]\}v\equiv\varphi(uv)
\]

The extension of a vacuous predicate is called 0\zl zero\zd\ class and denoted by 0\ul (or $\Lambda$) \ud\zl ;\zd\ the extension of a pred.\zl icate\zd\ belonging to every object is called universal class and denoted by 1 (or V)\zl .\zd\

$\mathbf{\llbracket 7. \rrbracket}$ For classes operation of $+$, $\cdot$~, $-$ which obey laws similar to the arithmetic laws are introduced by the following definitions:
\begin{tabbing}
\hspace{1.7em}\= $\alpha +\beta$ \= $=\;\; \hat{x}[x\,\varepsilon\,\alpha\vee x\,\varepsilon\,\beta]$\quad\=(sum)\\[.5ex]
\> $\alpha\cdot\beta$ \> $=\;\; \hat{x}[x\,\varepsilon\,\alpha\: .\: x\,\varepsilon\,\beta]$\> (intersection)\\[.5ex]
\> $-\alpha$ \> $=\;\; \hat{x}[\sim x\,\varepsilon\,\alpha]$\> (complement)\\[.5ex]
\>$\alpha -\beta$\> $=\;\;\alpha\cdot(-\beta)$\> (difference)
\end{tabbing}

\end{document}